%% file: pBMS1master.tex
\documentclass{amsart}
\usepackage{enumerate,color,amsmath,amssymb,amsfonts,bm}
\usepackage{frcursive,fullpage,graphicx,psfrag}

\newtheorem{Lemma}{Lemma}
\newtheorem{Remark}{Remark}
\newtheorem{Theorem}{Theorem}
\newtheorem{Proposition}{Proposition}
\newtheorem{Corollary}{Corollary}

\def\bbone{{\mathchoice {\rm 1\mskip-4mu l} {\rm 1\mskip-4mu l}
{\rm 1\mskip-4.5mu l} {\rm 1\mskip-5mu l}}}

\begin{document}

\author{Abdelmalek Abdesselam}
\address{Abdelmalek Abdesselam, Department of Mathematics,
P. O. Box 400137,
University of Virginia,
Charlottesville, VA 22904-4137, USA}
\email{malek@virginia.edu}

\author{Ajay Chandra}
\address{Ajay Chandra, Department of Mathematics,
P. O. Box 400137,
University of Virginia,
Charlottesville, VA 22904-4137, USA}
\email{ac2yx@virginia.edu}

\author{Gianluca Guadagni}
\address{Gianluca Guadagni, Department of Mathematics,
P. O. Box 400137,
University of Virginia,
Charlottesville, VA 22904-4137, USA}
\email{gg5d@virginia.edu}

\title{Rigorous quantum field theory functional integrals over the $p$-adics I: anomalous dimensions}

\begin{abstract}
In this article we provide the complete proof of the result announced in~\cite{ACG0} about
the construction of scale invariant non-Gaussian generalized stochastic
processes over three dimensional $p$-adic space. The construction includes that of the associated squared field and
our result shows this squared field has a dynamically generated anomalous dimension which rigorously
confirms a prediction made more than forty years ago, in an essentially identical situation, by K. G. Wilson. 
We also prove a mild form of universality for the model under consideration.
Our main innovation is that our rigourous renormalization group formalism allows for space dependent couplings.
We derive the relationship between mixed correlations and the dynamical systems features of our extended renormalization group
transformation at a nontrivial fixed point. The key to our control of the composite field is a partial linearization theorem
which is an infinite-dimensional version of the K{\oe}nigs Theorem in holomorphic dynamics.
This is akin to a nonperturbative construction of a nonlinear scaling field in the sense of F. J. Wegner
infinitesimally near the critical surface.
Our presentation is essentially self-contained and geared towards a wider audience. While primarily concerning the areas
of probability and mathematical physics we believe this article will be of interest to researchers in
dynamical systems theory, harmonic analysis and number theory. It can also be profitably read by graduate students in theoretical
physics with a craving for mathematical precision while struggling to learn the renormalization group.
\end{abstract}

\maketitle

\tableofcontents 

\input{Intro} 
\input{Prelim} 
\input{Formalres} 
\input{DefextRG} 
\input{Prelimest} 
\input{Mainest} 
\input{BulkRG} 
\input{DynsysI} 
\input{DynsysII} 
\input{Deviations} 
\input{Constructmeas} 
\input{Nontriv} 

\bigskip
\noindent{\bf Acknowledgements:}
{\small
We thank P. K. Mitter for generously sharing his notes~\cite{Mitter} containing
a nonperturbative proof of $\eta_{\phi}=0$ for the real BMS model in finite volume as well
as a formal perturbative calculation supporting the conjecture $\eta_{\phi^2}>0$, also in the real case.
We thank P. Cartier and M. Emerton for helping clarify
some of the material in \S\ref{padintrosec} and \S\ref{padmeassec}. We thank I. Herbst for the reference to Nelson's book~\cite{Nelson}.
The work of A. A. was supported in part by the National Science Foundation 
under grant DMS \# 0907198.}

\input{Bibli} 

\ 

\end{document}

%% file: Intro.tex
\section{Introduction}

Although the constructive approach to quantum field theory (QFT) and the rigorous
methods based on the renormalization group (RG) already have a long history,
the work from that area concerning the construction of composite fields and the
operator product expansion (OPE) is rather scarce.
All we managed to find after a review of the literature are: the work of
Feldman and R\c{a}czka~\cite{FeldmanR} followed by Constantinescu~\cite{Constantinescu}
(see also~\cite{EckmannE})
on composite fields (up to $\phi^3$) in the three-dimensional
massive $\phi^4$ model, 
and that of Iagolnitzer and Magnen~\cite{IagoM1,IagoM2} on the OPE for the two-dimensional
massive Gross-Neveu model.
Note that we do not count $\mathcal{P}(\phi)_2$ theories since the renormalization
of composite fields in that case requires nothing beyond what is already needed for the free field, namely,
Wick ordering~\cite{GlimmJcomp}.
Neither are we concerned here with perturbative results on composite field renormalization and the OPE.
For the latter the reader is referred to the excellent presentations in~\cite{Zimmermann,KellerK1,KellerK2,HollandsK}
and references therein.
Yet, for $\phi_{3}^{4}$ and $GN_2$ which respectively are superrenormalizable and asymptotically free in the
ultraviolet, the short distance behaviour is governed by a Gaussian RG fixed point. Thus the concerned
composite fields do not exhibit anomalous scaling dimensions.
In this article we construct a composite field with {\it dynamically} generated anomalous dimension
governed by a nontrivial RG fixed point.
Regarding similar anomalous dimensions for elementary rather
than composite fields we should mention the previous works~\cite{daVeiga} and~\cite{BenfattoGPS}.
The first concerns a model believed to have anomalous scaling (see, e.g., the review~\cite[\S2.7]{RosensteinWP}). However, after
a more than heroic effort, the author stopped at the construction of correlation functions. To get a hold on the
anomalous dimension would have required the extra work consisting of more precise estimates on the
quantity $\delta^{\ast}(N)$ in~\cite[p. 189]{daVeiga} together with the short distance asymptotics of some correlation function.
The second goes the full distance and proves the existence of anomalous dimension. However, the latter is governed by a line
of fixed points and depends on the coupling, unlike the situation for the paradigmatic Wilson-Fisher fixed point~\cite{WilsonF}.
We should from the onset warn the reader that our main theorem, given in \S\ref{formalstatsec},
seems in contradiction with a statement made in~\cite[p. 277]{GawedzkiK2}. Whether this contradiction is real or only apparent
remains to be seen. By apparent we mean a contradiction which could be explained away, e.g., by the difference between the models
or the objects to which the conflicting statements apply
as often happens when one takes limits in different orders.

The particular model studied in this article is what we call the $p$-adic BMS model
in honor of Brydges, Mitter and Scoppola who initiated the rigorous study of its real counterpart~\cite{BMS}.
This followed the earlier study of a similar model by Brydges, Dimock and Hurd~\cite{BDHeps}.
The real BMS model formally is a Radon-Nikodym perturbation of a massless Gaussian measure ${\rm d}\mu_{C_{-\infty}}(\phi)$
for a random scalar field $\phi$ on $\mathbb{R}^3$ with covariance
\[
C_{-\infty}=(-\Delta)^{-\left(\frac{3+\epsilon}{4}\right)}\ .
\]
Here $\epsilon>0$ is a small bifurcation parameter. The Radon-Nikodym weight is heuristically given by
a constant times
\begin{equation}
\exp\left(
-\int_{\mathbb{R}^3}\ 
\left\{g\phi(x)^4+\mu\phi(x)^2\right\}{\rm d}^{3}x
\right)\ .
\label{RadonNeq}
\end{equation}
Compared with the model studied in~\cite{BDHeps}, the BMS model has several advantages. The first is that it is technically
simpler since it does not require wave function renormalization when $\epsilon<1$.
It is also more physical. Indeed, the BMS family of models indexed by
$\epsilon$ includes the massless $\phi_{3}^{4}$ model at $\epsilon=1$.
Even for $\epsilon<1$ the BMS model is also, in all likelihood, Osterwalder-Schrader
positive and should therefore lead to the construction of a unitary QFT in Minkowski space.
Note that this kind of models with fractional powers of the Laplacian represent the best one can
presently tackle as far as the rigorous study of the phenomena associated to the Wilson-Fisher fixed point.
Indeed, the so-called $\epsilon$-expansion in the physics literature is usually based on dimensional regularization which,
as far as we know, has not been defined rigorously and nonperturbatively. See however~\cite{GurauMR} for an intriguing conjecture
which may lead to progress on this issue, although it needs some amendments as shown in~\cite{ScottS}.

In~\cite{BMS} the authors considered an {\it ad hoc} RG transformation for this model in the formal infinite
volume limit and they proved the existence of a nontrivial fixed point. They also constructed its local stable manifold.
This was followed by~\cite{Abdesselam} where connecting orbits joining the trivial fixed point to the nontrivial one
were constructed.
In this article we study the natural $p$-adic analogue of the real BMS model considered in~\cite{BMS,Abdesselam}.
Here the random field $\phi$ lives on $\mathbb{Q}_{p}^{3}$ instead of $\mathbb{R}^3$ but it still is a real-valued field.
There already is a rather large body of literature on $p$-adic QFT models (see~\cite{Missarov2012} and references therein).
However, we have not seen in this literature an explicit rigorous nonperturbative construction of the self-similar
scalar field studied in this article (or the very similar version on $\mathbb{Q}_p$ instead of $\mathbb{Q}_{p}^{3}$),
although the technology for doing that has been available for a long time, i.e., since the ground-breaking work of
Bleher and Sinai~\cite{BleherS1,BleherS2}. Indeed, the RG transformation naturally associated to a random field over the $p$-adics
falls under the umbrella of hierarchical RGs. This connection was briefly pointed out in~\cite{Bleher82}.
The only $p$-adic QFT model for which we have seen an explicit construction is that in~\cite{KochubeiS}.
Yet, the propagator in that model has a very mild singularity in the ultraviolet so that Wick ordering and
hypercontractivity techniques are enough to do the construction.
Furthermore, the issue of infinite volume limit~\cite{GlimmJinfvol,GlimmJS} is not settled.

In the next few paragraphs we will try to give a rough idea of our result and methods. The precise statement of our theorem
which was already announced in~\cite{ACG0} will be recalled in \S\ref{formalstatsec} after the necessary definitions, in particular
regarding $p$-adic analysis, are presented in \S\ref{basicprelsec}. For the sake of pedagogy and efficiency we will pretend to
be working over $\mathbb{R}^3$ in what follows. This will delay having to come to terms with some oddities of the $p$-adic
world such as: the size of $L$ is $L^{-1}$, the lattice is $\mathbb{Q}_{p}^{3}/\mathbb{Z}_{p}^{3}$ while $\mathbb{Z}_{p}^{3}$
is a lattice cell, and so on. The reader should bear in mind that it is the $p$-adic analogues of the following statements
which are addressed in this article.
First introduce a number $L>1$ which serves as a yardstick for measuring changes of scale. This is the analogue of $L=2$ used in
the dyadic decomposition methods in harmonic analysis.
We then define a cut-off covariance $C_{r}$, for $r\in\mathbb{Z}$, by suppressing Fourier modes or momenta $k$ with size greater than $L^{-r}$
starting from the non-cut-off covariance 
\[
\widehat{C}_{-\infty}(k)=\frac{1}{|k|^{\left(\frac{3+\epsilon}{2}\right)}}=\frac{1}{|k|^{3-2[\phi]}}
\]
where the symbol $[\phi]$ stands for the quantity $\frac{3-\epsilon}{4}$, i.e., the scaling dimension of the free massless Gaussian.
We also change the integration set in (\ref{RadonNeq}) from $\mathbb{R}^3$ to a finite volume $\Lambda_{s}$
of linear size $L^s$, $s\in\mathbb{Z}$.
As renormalization theory tells us to, we also replace the couplings $g$, $\mu$ in (\ref{RadonNeq}) by $r$-dependent quantities
$\tilde{g}_r$ and $\tilde{\mu}_r$. This dependence is also called a bare ansatz. Finally, we also replace the monomials
$\phi(x)^4$, $\phi(x)^2$ by their Wick ordered analogues $:\phi^4:_{C_r}(x)$, $:\phi^2:_{C_r}(x)$. This corresponds to a triangular
change of coordinates, namely, switching from the monomial basis to that of Hermite polynomials which is more convenient.
These modifications result in a well defined probability measure ${\rm d}\nu_{r,s}(\phi)$
with moments or correlators
\[
\langle\phi(x_1)\cdots\phi(x_n)\rangle_{r,s}\ .
\]
We construct the wanted measure ${\rm d}\nu$ as the limit of the ${\rm d}\nu_{r,s}$ when the ultraviolet cut-off
$r$ is taken to $-\infty$ and the infrared cut-off $s$ is taken to $\infty$.
The particular bare ansatz we use is of the form $\tilde{g}_r=L^{-(3-4[\phi])r} g$,
$\tilde{\mu}_r=L^{-(3-2[\phi])r} \mu$ where $g$, $\mu$ are fixed quantitites. In fact the mass $\mu$ has to be fine-tuned to
its so-called critical value $\mu_{\rm c}(g)$.
The correlators $\langle\phi(x_1)\cdots\phi(x_n)\rangle_{r,s}$ are obtained as derivatives of
a moment generating function $\mathcal{S}_{r,s}(\tilde{f})$ in terms of a test function $\tilde{f}$.
This generating function is of the form
\[
\mathcal{S}_{r,s}(\tilde{f})=\frac{\int {\rm d}\mu_{C_r}(\phi)\ e^{-V_{r,s}(\phi)+\phi(\tilde{f})}}
{\int {\rm d}\mu_{C_r}(\phi)\ e^{-V_{r,s}(\phi)}}
\]
where $\phi(\tilde{f})$ denotes the distributional pairing.

The first step of the analysis involves a `rescaling to unit lattice' which is a change of variables in the field $\phi(x)$.
Indeed, the latter sampled according to the Gaussian measure with covariance $C_r$ has the same law as
the field $L^{-r[\phi]}\phi(L^{-r} x)$ where the new $\phi$ is sampled according to the Gaussian measure with covariance $C_{0}$.
One then repeatedly applies an RG transformation to both the numerator and denominator
of the resulting expression for $\mathcal{S}_{r,s}(\tilde{f})$.
In its naive form this transformation is based on the identity
\begin{equation}
\int {\rm d}\mu_{C_0}(\phi)\ e^{-V(\phi)}=
\int {\rm d}\mu_{C_0}(\phi)\ e^{-V'(\phi)}
\label{naiveRGeq}
\end{equation}
with a new potential $V'$ given by
\[
V'(\phi)=-\log\left(
\int {\rm d}\mu_{\Gamma}(\zeta)\ e^{-V(\zeta+L^{-[\phi]}\phi(L^{-1}\bullet))}
\right)
\]
where $\Gamma=C_0-C_1$ is the fluctuation covariance involving Fourier modes in the shell $L^{-1}< |k|\le 1$.
In fact we need to extract a field independent quantity $\delta b$ in (\ref{naiveRGeq})
which becomes
\[
\int {\rm d}\mu_{C_0}(\phi)\ e^{-V(\phi)}=e^{\delta b(V)}
\int {\rm d}\mu_{C_0}(\phi)\ e^{-V'(\phi)}
\]
while the new potential rather is
\[
V'(\phi)=\delta b(V)-\log\left(
\int {\rm d}\mu_{\Gamma}(\zeta)\ e^{-V(\zeta+L^{-[\phi]}\phi(L^{-1}\bullet))}
\right)\ .
\]
While the above transformation $V\rightarrow V'$ is well defined (in finite volume), it is difficult to exploit it.
Indeed, this transformation is most useful where it is most singular, i.e., in infinite volume. In their seminal article~\cite{BrydgesY},
Brydges and Yau, based on earlier work in constructive QFT, found a solution to this difficulty
which is to introduce a suitable lift:
\[
\begin{array}{ccc}
\vec{V} & \longrightarrow & \vec{V}'\\
\downarrow & \ & \downarrow \\
V & \longrightarrow & V'
\end{array}
\]
for the naive RG transformation $V\rightarrow V'$. 
Such a Brydges-Yau lift is highly nonunique and quite complicated but it has the advantage of providing
a dynamical system which can be analyzed in the infinite volume limit using rigorous estimates.
We define the Brydges-Yau lift
relevant to the present model in \S\ref{defBYliftsec}.
After the initial rescaling to unit lattice, the numerator of the ratio expressing $\mathcal{S}_{r,s}(\tilde{f})$ gives rise
to an initial vector $\vec{V}^{(r,r)}(\tilde{f})$ while the denominator produces a similar one $\vec{V}^{(r,r)}(0)$.
The moment generating function is obtained from the log-moment generating function $\mathcal{S}^{\rm T}(\tilde{f})$
which itself is the limit of quantities roughly given by
\[
\mathcal{S}_{r,s}^{\rm T}(\tilde{f})=\sum_{r\le q<s} \left\{\delta b[\vec{V}^{(r,q)}(\tilde{f})]
-\delta b[\vec{V}^{(r,q)}(0)]\right\}
\]
where $\vec{V}^{(r,r)}\rightarrow \vec{V}^{(r,r+1)}\rightarrow  \vec{V}^{(r,r+2)}\rightarrow \cdots$ denote
the RG iterates of the initial vectors.

We remove the cut-offs by
showing the convergence of the series over (logarithmic) scales $q$
\begin{equation}
\mathcal{S}^{\rm T}(\tilde{f})=
\sum_{q\in\mathbb{Z}} \left\{\delta b[\vec{V}^{(-\infty,q)}(\tilde{f})]
-\delta b[\vec{V}^{(-\infty,q)}(0)]\right\}\ .
\label{qseriesintroeq}
\end{equation}
This hinges on controlling the deviations from the bulk $\vec{V}^{(r,q)}(\tilde{f})-\vec{V}^{(r,q)}(0)$
corresponding to the effect of a local perturbation due to the test function.

Most of the previous work in rigorous RG theory relies on a translation trick (completing the square)
which reduces the RG for the numerator to that of the denominator (see, e.g.,~\cite[Ch. 4]{BenfattoG}).
However this creates technical difficulties under RG iteration since the effect of this trick is to
add an $\tilde{f}$ dependent term to the so-called background field (the $\phi$ on the right-hand side of (\ref{naiveRGeq}))
which becomes more and more singular. Since we would like to treat similar perturbations by higher smeared powers of the field,
with a view towards composite field renormalization, a more systematic approach would be to define a more general RG map
which handles the data $\vec{V}^{(r,q)}(\tilde{f})$ produced by the numerator.
This is one of the main technical innovations in this article: we define in \S\ref{defBYliftsec} what we call the {\it extended}
RG transformation
which handles potentials with couplings $g,\mu$ which are allowed to vary with the position $x$ of the $\phi(x)^m$ monomials.
To first order of approximation this extended RG amounts to taking local averages of the previous couplings
and multiplying them by the appropriate power-counting factor. This is somewhat dual to the usual
block-spinning approach which acts on the field rather than the couplings.
The convergence of (\ref{qseriesintroeq}) is proved by exhibiting suitable decay in both infinite directions.
Very roughly one can introduce
$L^{q_{+}}$ which corresponds to the size of the support of the test function $\tilde{f}$ in direct space,
while $L^{-q_{-}}$ corresponds to size of the support of the Fourier transform of $\tilde{f}$.
Over the $p$-adics this is strictly true while only approximately so over the reals.
The uncertainty principle can be seen as the relation $q_{-}\le q_{+}$.
Most of the contribution to the series (\ref{qseriesintroeq}) comes from the range of scales where the test function lives, i.e., for $q$
between $q_{-}$ and $q_{+}$.
For the $\phi$ perturbations, because of our specific choice of cut-offs, the contribution of scales $q<q_{-}$ is identically zero.
However, when considering $\phi^2$ perturbations smeared by another test function $\tilde{j}$ this is no longer the case.
In fact, a large part of the analysis is devoted to proving the needed decay in $q$ after subtraction of a suitable linear term
in $\tilde{j}$.
Two sectors of the extended RG space are important for the analysis. The bulk RG corresponds to spatially uniform couplings
(as for the denominator of $\mathcal{S}_{r,s}(\tilde{f})$) and is what is needed for the control of the ultraviolet region.
Indeed, after the initial rescaling, the test functions become diluted over a very large volume and thus appear spatially constant, i.e.,
bulk-like. On the other hand, in the infrared region the deviations from the bulk appear point-like. There is no more dangerous $L^3$
factor responsible for the expanding or relevant RG directions and this fact is key to the needed decay in $q$ for $q>q_{+}$.

In this article we construct the mixed correlations
\[
\langle
\phi(x_1)\cdots\phi(x_n)\ \phi^{2}(y_1)\cdots \phi^{2}(y_m)
\rangle
\]
as distributions even at coinciding points~\cite{EckmannE}.
This is done by controlling a more complicated generating function which very roughly is given by
\[
\mathcal{S}_{r,s}(\tilde{f},\tilde{j})=\frac{\int {\rm d}\mu_{C_r}(\phi)\ e^{-V_{r,s}(\phi)+\phi(\tilde{f})+\phi^2(\tilde{j})}}
{\int {\rm d}\mu_{C_r}(\phi)\ e^{-V_{r,s}(\phi)}}
\]
instead of $\mathcal{S}_{r,s}(\tilde{f})$.
In the initial rescaling or dilution, the test function $\tilde{f}$ is weakened by a factor $L^{(3-[\phi])r}$
while the test function $\tilde{j}$ should normally be weakened by a factor $L^{(3-2[\phi])r}$. We say normally because this would be
the case in the absence of anomalous dimension. It turns out, however, the correct weakening factor one needs to use is $\alpha_{\rm u}^r$
where $\alpha_{\rm u}$ is the eigenvalue in the expanding direction at the nontrivial RG fixed point.
The anomalous dimension is due to the fact $\alpha_{\rm u}$ is strictly less than the
analogous eigenvalue $L^{3-2[\phi]}=L^{\frac{3+\epsilon}{2}}$ at the Gaussian fixed point.
This strict inequality was known for a long time. Indeed, for a very similar hierarchical model,
Bleher and Sinai obtained small $\epsilon$ asymptotic expansions for the eigenvalues at the fixed point~\cite[Thm. 3.1]{BleherS2}.
The previous fact follows from such results. The new and nontrivial part of our proof concerns the translation of this fact
about the eigenvalue $\alpha_{\rm u}$ into information about the $\phi^2$ correlation functions.
This we accomplish thanks to a, possibly new, infinite-dimensional version of the Theorem of K{\oe}nigs in holomorphic dynamics.

In~\cite[\S8]{Koenigs} K{\oe}nigs proves the following result.
Let $F(z)$ be an analytic function defined near zero such that $F(0)=0$ and $F'(0)=\alpha$ with $0<|\alpha|<1$.
Then the limit
\[
\Psi(z)=\lim_{n\rightarrow\infty} \alpha^{-n} F^n(z)
\]
exists and is analytic near zero. It also satisfies $\Psi'(0)=1$ and therefore provides a conjugation
of $F$ to its linearization at zero. Indeed, one has the intertwining relation
\[
\alpha\Psi=\Psi\circ F\ .
\]
A two-line proof (which however was a great source of inspiration to us) was given in~\cite[p. 6]{Abate}.
It amounts to showing the uniform absolute convergence of the telescopic sum
with general term 
\[
\alpha^{-(n+1)} F^{n+1}(z)-\alpha^{-n} F^n(z).
\]
In our situation we have an expanding eignenvalue and we do not know if our map, i.e., the RG is invertible
(common wisdom says not since
the RG is an irreversible process of erasing degrees of freedom).
So we are in a situation where we rather have to construct the limit
\[
\Psi(z)=\lim_{n\rightarrow\infty} F^n(\alpha^{-n}z)
\]
instead.
In fact, the result we establish is more general. It is the existence and analyticity
of
\[
\lim_{n\rightarrow\infty} RG^n(v+\alpha_{\rm u}^{-n}w)
\]
for a transformation $RG$ in an infinite-dimensional Banach space. This map has a fixed point $v_{\ast}$
with a dimension one unstable manifold $W^{\rm u}$ with eigenvalue $\alpha_{\rm u}$. The previous statement
applies to points $v$ on the codimension one stable manifold $W^{\rm s}$, while $w$ can be any vector which is not too large.
This realizes a partial linearization of the RG map and is tantamount to constructing a nonlinear scaling field
in the sense of Wegner~\cite{Wegner} infinitesimally close to the critical surface $W^{\rm s}$.
Our result holds regardless of nonresonance conditions because the unstable manifold in the mass direction is one-dimensional.
One should also note the similarity of the above construction of a conjugating function
to that of the classical scattering analogue of M{\o}ller wave operators~\cite{Hunziker,Nelson,Simonscat}.
Indeed, Nelson in particular in his book popularized the scattering idea as a way to obtain such conjugations. He also
attempted to derive Sternberg's Linearization Theorem in this way.
Linearization results in infinitely many dimensions are rather rare. See however~\cite{Grunbaum} which contains an interesting
discussion of the relation between linearization and classical scattering.
In essence, this is also related to the age old method of variation of constants~\cite{CarlettiMV}.
Other authors attempted to obtain a partial linearization result (in the $C^{\infty}$ category)
analogous to ours~\cite{ColletE1}. However, their
proof is incorrect as they later acknowledged in~\cite[\S12]{ColletE2}.
Our partial linearization theorem derived in \S\ref{Koenigssec} is in the analytic category.
This is essential for the construction of our generalized stochastic processes.
Indeed, this analyticity is inherited by the series (\ref{qseriesintroeq}), in the presence of the test function $\tilde{j}$.
As a result we get $n!$ bounds on the moments of both the elementary field $\phi$ and composite field $\phi^2$.
These are needed for reconstructing the measures from the moments.

This allows us to produce generalized stochastic process which are translation, rotation and scale invariant by the subgroup $L^{\mathbb{Z}}$
of the full group of scale transformations.
In a companion paper~\cite{ACG2} we will show that
our model has full scale invariance. This would open the door to the investigation of conformal
invariance along the lines of~\cite{LernerM,Lerner}. If one also makes progress on the OPE one could even contemplate
the possibility of an exact solution (see~\cite[p. 70]{PordtW} and~\cite{HarlowSSS} for related work).
Our article can be seen as taking some steps towards the {\it systematic} elaboration of a very general theory
of `local fields on local fields', to borrow a pun attributed to S. Evans.
A brief sketch of such a theory was given in~\cite{AbdesselamOW}. 
It is a natural continuation of a line of thought pursued by the Soviet School of probability theory and mathematical physics,
see in particular~\cite{Dobrushin1,Dobrushin2,Sinai}.
Note that the presence of anomalous dimension (for $\phi^2$)
should make it clear that our field $\phi$ is not subordinated
to a Gaussian field in the sense of~\cite{Dobrushin1,Major}. However proving this would require the beginning of the OPE
(namely generating the $\phi^2$ field by collapsing two $\phi$'s) which we do not address in this article and leave for a future publication.
We believe our article opens the way to a rigorous investigation of the OPE.
The latter is of fundamental importance. Indeed, vertex operator algebras (see, e.g.,~\cite{FrenkelB}), chiral algebras~\cite{BeilinsonD}
and factorization algebras~\cite{CostelloG} can be seen
as mathematical constructions which try to capture the OPE structure in QFT.

Our primary motivation for considering the $p$-adic version of the BMS model is that it is a good
toy model for the original version over the reals. 
In fact the model over $\mathbb{R}$  (for $\epsilon$ small)
is expected to exhibit similar features: absence of anomalous dimension ($\eta_{\phi}=0$) for the elementary field
$\phi$ and presence of anomalous dimension ($\eta_{\phi^2}>0$) of order $\epsilon$ for the composite field $\phi^2$.
Regarding this real model, P. K. Mitter~\cite{Mitter} proved nonperturbatively but in finite volume
that the elementary field has no anomalous scaling. He also did a formal perturbative calculation for
the anomalous dimension of $\phi^2$ (see also~\cite{FisherMN}).
Regarding the $p$-adic model, one can argue that the prediction of the properties $\eta_{\phi}=0$
and $\eta_{\phi^2}>0$ was made forty years ago by Wilson himself in~\cite{Wilson72}. Indeed, the discussion
in that article was in the framework of Wilson's approximate RG recursion which is what we now call the hierarchical RG.
In other words, the situation considered in~\cite{Wilson72} is essentially identical to that of the $p$-adic BMS model.
Using hierarchical models such as the $p$-adic model in order to shed light on Euclidean ones such as the real
BMS model has been and will continue to be a fruitful approach. Indeed, some of Wilson's ideas on the RG were already
present in the article~\cite{Wilson65}. Nevertheless, the first systematic exposition of what is now known as the Wilson
RG philosophy most likely is the article~\cite{WilsonII} which is about the approximate recursion.
Only later came the adaptation of his methods to the model over $\mathbb{R}$ which found its definitive presentation
in the famous lectures~\cite{WilsonK}. Our methodology is to try to follow a similar path.
In doing so we took great care in choosing, for our treatment of the $p$-adic case, methods which are known to work
over the reals. We simply transposed such methods, in as natural a way as possible, to the $p$-adic
setting. Indeed, the definition of our RG transformation in \S\ref{defBYliftsec}, the corresponding estimates in \S\ref{longlemsec}
as well as some of the dynamical systems techniques in \S\ref{dyn1sec} were directly adapted from~\cite{BMS,Abdesselam}.
We believe the main nontrivial task which remains in order to extend the results of the present paper to the real case
is to devise a proper analogue of the extended RG given in \S\ref{defBYliftsec}.
Such a transformation should essentially reduce, in the special case of spatially uniform potentials, to the RG transformation
in~\cite{BMS,Abdesselam}. This problem is a matter of harmonic analysis and is thematically similar to
extending a result about Walsh series to Fourier series (see, e.g.,~\cite{DemeterLTT}).

Our secondary motivation is to help facilitate the investigation of the connection
between QFT and number theory. This is perhaps still at the speculative stage. Nevertheless,
see~\cite{Burnol1,Burnol2,GerasimovLO1,GerasimovLO2,GerasimovLO3,Leichtnam}
for interesting work in this direction.
We hope our article will help number theorists unterstand how the {\it rigorous} RG works when used for the construction
and study of QFT functional integrals.

We will end this introduction by commenting on the length of this article.
There are two reasons for this: the choice of methods and the high level of detail.
As we said earlier we did not try to prove our result about the $p$-adic BMS model in the quickest and most direct manner.
The potential for adaptability to the real situation was the overarching principle that guided our investigation.
Also note that some of our results are stated in more generality than needed for the sole purpose of proving our main result which is
Theorem \ref{themainthm}.
An example is Theorem~\ref{mainestthm}. The benefit reaped from this methodological choice is that
we will be able to reuse Theorem~\ref{mainestthm}, exactly as stated,
and in combination with the techniques from~\cite{Abdesselam},
in order to construct the generalized processes or QFTs corresponding to the RG orbits connecting the Gaussian
and infrared fixed points.
As for the amount of detail, we note the following.
Interest in QFT by mathematicians is high while the community of people with a working
knowledge of constructive QFT and rigorous RG theory is very small. This gives us a strong incentive to write
this article so it is understandable to a wider audience.
Our presentation is essentially self-contained and only uses very modest prerequisites
to be mentioned in the next section.
We believe our article can also be read profitably by students in theoretical physics who would like to see, on a
simple example,
what is the precise relation between the RG dynamical system near a fixed point and the behaviour of the correlation functions.
Such a reader may skip the sections containing estimates such as \S\ref{longlemsec}.
Last but not least, our article concerns matters on which much was published that contained errors.
We needed the high amount of detail to be absolutely sure that our proof is correct.

%% file: Prelim.tex
\section{Preliminaries}\label{basicprelsec}

\subsection{Generalities about $p$-adics}\label{padintrosec}

Here we breifly review the basic notions about $p$-adics which are needed in this article.
More details can be found in~\cite{VladimirovVZ,AlbeverioKS,GoldfeldH}. See also~\cite{Zabrodin}
for a quick introduction to the $p$-adics which includes very helpful pictures.
Let $p$ be a prime number and consider the $p$-adic absolute value $|\cdot|_p$ on $\mathbb{Q}$
defined by $|x|_p=0$ if $x=0$ and $|x|_p=p^{-k}$ if $x=\frac{a}{b}\times p^k$ where $a,k\in\mathbb{Z}$ and $b$, a positive integer, 
are such that $a,b$ are coprime and neither are divisible by $p$.
The field $\mathbb{Q}_p$ of $p$-adic numbers is the completion of $\mathbb{Q}$ with respect to this absolute value.
Every $x$ in $\mathbb{Q}_p$
has a unique convergent series representation
\[
x=\sum_{n\in\mathbb{Z}} a_n p^n
\]
where the digits $a_n$ belong to $\{0,1,\ldots,p-1\}$ and at most finitely many of them are nonzero for negative $n$.
The absolute value of $x\neq 0$ can be recovered from this representation as $|x|_p=p^{-v_p(x)}$ where
\[
v_p(x)=\min\{n\in\mathbb{Z}\ |\ a_n\neq 0\}\ .
\]
Using the same representation one can define the fractional (or polar) part of $x$ which is
$\{x\}_p=\sum_{n<0} a_n p^n$.
The closed unit ball $\mathbb{Z}_p=\{x\in\mathbb{Q}_p|\ |x|\le 1\}$ is a compact subring of $\mathbb{Q}_p$.
From now on we will drop the $p$ subscript from the absolute value.
The additive Haar measure on $\mathbb{Q}_p$ normalized so that $\mathbb{Z}_p$ has measure one
will simply be denoted by ${\rm d}x$.
In $d$ dimensions,
the $p$-adic norm of a vector $x=(x_1,\ldots,x_d)\in\mathbb{Q}_p^d$ is defined as
$|x|=\max\{|x_1|,\ldots,|x_d|\}$. The product measure ${\rm d}^d x$ obtained from the previous one-dimensional measure
is invariant by the subgroup $GL_d(\mathbb{Z}_p)$ of $GL_d(\mathbb{Q}_p)$.
The subgroup $GL_d(\mathbb{Z}_p)$ is defined as the set of $d\times d$ matrices which
together with their inverses have entries in $\mathbb{Z}_p$.
This subgroup is the maximal compact subgroup of $GL_d(\mathbb{Q}_p)$ (unique up to conjugacy) and is the
natural analogue of the orthogonal group $O(d)$ acting on $\mathbb{R}^d$.
The use of the maximum in the definition of the norm is motivated by the resulting invariance with respect to $GL_d(\mathbb{Z}_p)$.

The space of real (resp. complex) test functions $S(\mathbb{Q}_p^d,\mathbb{R})$ (resp. $S(\mathbb{Q}_p^d,\mathbb{C})$)
is the Schwartz-Bruhat space of compactly supported locally constant real-valued (resp. complex-valued) functions on $\mathbb{Q}_p^d$.
If we do not specify the target, then we mean $\mathbb{R}$.
Recall that a seminorm on $S(\mathbb{Q}_p^d)$ is a function $\mathcal{N}:S(\mathbb{Q}_p^d)\rightarrow [0,\infty)$
which satisfies the usual norm axioms except the requirement that $\mathcal{N}(f)=0$ implies $f=0$.
The coarsest topology on $S(\mathbb{Q}_p^d)$ which makes all possible seminorms continuous is called the finest locally convex
topology and it is the one we use. The space of distributions $S'(\mathbb{Q}_p^d)$
simply is the topological dual of $S(\mathbb{Q}_p^d)$ which turns out to be the algebraic dual.
Note that $S(\mathbb{Q}_p^d)$ is a nonmetrizable topological vector space. Therefore
the theory of denumerably Hilbert nuclear spaces does not apply to it.

The Fourier transform of a complex valued test function $f$ is defined by
\[
\widehat{f}(k)=\int_{\mathbb{Q}_p^d} f(x)\exp(-2i\pi\{k\cdot x\}_p)\ {\rm d}^d x
\]
where $k\cdot x=k_1 x_1+\cdots+k_d x_d$ and the rational $\{k\cdot x\}_p$ is seen as a real number.
One has that the characteristic function of $\mathbb{Z}_p^d$ is fixed by the Fourier transform, that is
$\widehat{\bbone}_{\mathbb{Z}_p^d}=\bbone_{\mathbb{Z}_p^d}$.
From this it easily follows that the space $S(\mathbb{Q}_p^d,\mathbb{C})$ is stable by Fourier transform.
One can also define the Fourier transform of distributions by duality.

One has an analogue of the nuclear theorem in this setting which allows one to identify an
$n$-linear form $W:S(\mathbb{Q}_p^d)\times\cdots\times S(\mathbb{Q}_p^d)\rightarrow \mathbb{R}$
with a distribution in $S'(\mathbb{Q}_p^{nd})$. We believe the most expedient way
of proving such results as well as the ones in \S\ref{padmeassec} 
is by following B. Simon's philosophy of exploiting a topological vector space isomorphism
with a very concrete space of sequences (see~\cite{SimonJMP} and~\cite[\S{I.2}]{SimonFI}).
Indeed, it is easy to see that $S(\mathbb{Q}_p^d)$ can be written as a countable union
of an increasing sequence of finite-dimensional vector spaces $V_1\subset V_2\subset\cdots$.
One can construct a basis $(f_n)_{n\in\mathbb{N}}$ by taking a basis of $V_1$, then appending
vectors needed to complete it into a basis of $V_2$, etc. If one takes $(e_n)_{n\in\mathbb{N}}$
to be the Gram-Schmidt orthonormalization for the $L^2$ inner product $\langle\bullet,\bullet\rangle$ then the map
\[
f\longmapsto \left(\langle e_n, f\rangle\right)_{n\in\mathbb{N}}
\]
realizes such an isomorphism from $S(\mathbb{Q}_p^d)$ to $\mathsf{s}=\oplus_{n\in\mathbb{N}}\mathbb{R}$.
Namely, $\mathsf{s}$ is the space of almost finite sequences $(x_n)_{n\ge 0}$ of real numbers,
also equipped with the finest locally convex topology. The dual is the space $\mathsf{s}'=\mathbb{R}^{\mathbb{N}}$
of all sequences $(y_n)_{n\ge 0}$ with duality pairing given by $\sum_{n\ge 0} x_n y_n$.
In this new setting it is very easy to prove statements such as: a bilinear form on $\mathsf{s}$ is automatically continuous.
This follows from the, somewhat counterintuitive, remark that if $(a_{i,j})_{(i,j)\in\mathbb{N}^2}$
is an array of nonnegative numbers, there exists a sequence $b_n\ge 0$ such that $a_{(i,j)}\le b_i b_j$ for all $i$ and $j$.
The same property also holds for multilinear maps.

Next we need to define some transformations which will allow us to give a precise formulation for the notions
of translation, rotation and scale invariance.
If one views a point $x$ in $\mathbb{Q}_p^d$ as a column vector
then one has a left-action of $GL_d(\mathbb{Z}_p)$ on points simply by matrix multiplication.
It results in left-actions on test functions $f$, distributions 
$\phi$ and more generally $n$-linear forms $W$
on $S(\mathbb{Q}_p^d)$, using
\[
(M\cdot f)(x)=f(M^{-1}x)\ ,
\]
\[
(M\cdot\phi)(f)=\phi(M^{-1}\cdot f)\ ,
\]
\[
(M\cdot W)(f_1,\ldots,f_n)=W(M^{-1}\cdot f_1,\ldots,M^{-1}\cdot f_n)\ .
\]
Such objects are called rotation invariant if they are preserved by all $M\in GL_d(\mathbb{Z}_p)$.
If one formally thinks of a distribution $\phi$ as a `function' via the $L^2$ pairing
\[
\phi(f)=\int_{\mathbb{Q}_p^d} \phi(x) f(x)\ {\rm d}^d x
\]
then the choice of definition means ``$(M\cdot\phi)(x)=\phi(M^{-1}x)$''.
Thus, a distribution $\phi$ is rotation invariant if ``$\phi(M^{-1}x)=\phi(x)$'' for all $M$ and $x$. 

Likewise regarding translations, one can define for $y\in\mathbb{Q}_p^d$ the transformations
\[
\tau_y(x)=x+y\ ,
\]
\[
\tau_y(f)(x)=f(x-y)\ ,
\]
\[
\tau_y(\phi)(f)=\phi(\tau_{-y}(f))\ ,
\]
\[
\tau_y(W)(f_1,\ldots,f_n)=W(\tau_{-y}(f_1),\ldots,\tau_{-y}(f_n))\ .
\]
One then defines the notion of invariance by translation for such objects
in the same way as before.

We now consider scaling transformations. Given $\lambda\in\mathbb{Q}_p^\ast=\mathbb{Q}_p\backslash\{0\}$, we write:
\[
(\lambda\cdot f)(x)=f(\lambda^{-1}x)\ ,
\]
\[
(\lambda\cdot\phi)(f)=|\lambda|^d\ \phi(\lambda^{-1}\cdot f)\ .
\]
This corresponds to the formal equation ``$(\lambda\cdot \phi)(x)=\phi(\lambda^{-1}x)$''.
A distribution $\phi$ is called partially scale invariant with homogeneity $\alpha\in\mathbb{R}$ with respect to a subgroup
$H$ of the full scaling group $p^{\mathbb{Z}}\subset \mathbb{Q}_p^\ast$ if $\lambda\cdot\phi=|\lambda|^{-\alpha}\phi$
for all $\lambda\in H$.
This formally means ``$|\lambda|^{\alpha}\phi(\lambda^{-1} x)=\phi(x)$''.

\subsection{Probability measures on the space of distributions}\label{padmeassec}

The generalized stochastic processes we will be interested in are probability measures
on the space of distributions $S'(\mathbb{Q}_{p}^{d})$. The $\sigma$-algebra
is the cylindrical one $\mathcal{C}$ which is the smallest that makes the maps $\phi\rightarrow \phi(f)$
measurable, for all test functions $f$.
Such a probability measure $\nu$ is called rotation invariant if for any $M\in GL_d(\mathbb{Z}_p)$
the push-forward (or direct image) of $\nu$ by the map $\phi\mapsto M\cdot\phi$ is $\nu$ itself.
Invariance by translation is defined in the same way.
A probability measure on $S'(\mathbb{Q}_p^d)$ is called partially scale invariant with homogeneity $\alpha$
with respect to the subgroup $H$ if the push-forward of $\nu$ by
the map $\phi\mapsto |\lambda|^{\alpha}(\lambda\cdot\phi)$
is $\nu$ itself, for all $\lambda\in H$.

We now need two theorems which allow us to construct and identify probability measures from their characteristic
functions or their moments. The first one is the analogue of the Bochner-Minlos Theorem (see~\cite[\S{I.2}]{SimonFI})
in the $p$-adic setting. The second is a reconstructions theorem with $n!$ bounds for the Hamburger Moment Problem
on $S'(\mathbb{Q}_p^d)$.

\begin{Theorem}\label{BochnerMthm}
Let $\Phi$ be a function $S(\mathbb{Q}_p^d)\rightarrow\mathbb{C}$
which satisfies
\begin{enumerate}
\item
$\Phi(0)=1$,
\item
$\Phi$ is continuous,
\item
For all $n\ge 1$, all test functions $f_1,\ldots,f_n$ in $S(\mathbb{Q}_p^d)$ and all
complex numbers $z_1,\ldots,z_n$,
\[
\sum_{a,b=1}^n \bar{z}_a z_b \ \Phi(f_b-f_a)\ \in\ [0,\infty)\ ;
\]
\end{enumerate}
then there exists a unique probability measure $\nu$ on the measurable space $(S'(\mathbb{Q}_p^d),\mathcal{C})$
such that for all $f\in S(\mathbb{Q}_p^d)$ we have
\[
\Phi(f)=\int_{S'(\mathbb{Q}_p^d)}\ {\rm d}\nu(\phi)\ e^{i\phi(f)}\ .
\]
\end{Theorem}

\begin{Theorem}\label{Hamburgerthm}
Let $(S_n)_{n\ge 0}$ be a sequence of distributions with $S_n\in S'(\mathbb{Q}_{p}^{nd})$
which satisfies
\begin{enumerate}
\item
$S_0=1$,
\item
for any $n$, $S_n$ is invariant by the permutation group $\mathfrak{S}_n$,
\item
for all almost finite sequence of test functions $(h_n)_{n\ge 0}$ with $h_n\in S(\mathbb{Q}_p^{nd},\mathbb{C})$
one has
\[
\sum_{n,m\ge 0} S_{n+m}(\overline{h_n}\otimes h_m)\ \in\ [0,\infty)\ ,
\]
\item
For all finite dimensional complex subspace $V$ of $S(\mathbb{Q}_p^{d},\mathbb{C})$
there exists a semi-norm $\mathcal{N}_V$ on $S(\mathbb{Q}_p^{d},\mathbb{C})$ such that for all $n\ge 0$ and all
$f_1,\ldots,f_n$ in $V$ one has
\[
|S_n(f_1\otimes\cdots\otimes f_n)|\le n!\times \mathcal{N}_V(f_1)\times\cdots\times \mathcal{N}_V(f_n)\ ;
\]
then there exists a unique probability measure with finite moments $\nu$ on the measurable space $(S'(\mathbb{Q}_p^d),\mathcal{C})$
such that for all $f_1,\ldots,f_n\in S(\mathbb{Q}_p^d,\mathbb{C})$ we have
\[
S_n(f_1\otimes\cdots\otimes f_n)=\int_{S'(\mathbb{Q}_p^d)}\ {\rm d}\nu(\phi)\ \phi(f_1)\cdots\phi(f_n)\ .
\]
\end{enumerate}
\end{Theorem}
Note that we used in the statement of the second theorem some obvious functorial properties of the nuclear theorem
with respect to complexification.
We do not give the proofs of these two theorems. We simply note that, via the isomorphism with $\mathsf{s}$ and $\mathsf{s}'$
and Kolmogorov's Extension Theorem, they reduce to their finite dimensional versions which are classical results in analysis.
Finally, we will use some basic formulas for Gaussian integration and manipulation of Wick monomials. These are left
as easy exercises for the reader. This material is also covered in~\cite[Ch. 9]{GlimmJ}, \cite[\S{I.1}]{Simon},
\cite[\S2.2]{Salmhofer}.

\subsection{Banach space analyticity}

A tool which we use a lot in this article is the theory of analytic maps in the complex Banach space context.
See for instance~\cite{Chae} for an introduction or reminder.
We also need to consider the analyticity of maps on $S(\mathbb{Q}_p^d,\mathbb{C})$. However when we do so
we always restrict to a finite-dimensional subspace so we only need the usual notion of analyticity on $\mathbb{C}^n$.
This allows us to avoid using the more involved theory of analyticity in locally convex spaces~\cite{Dineen}.
We only state here a lemma that is used many times in this article and which allows us to get Lipshitz estimates in an effortless manner.
We use the notation $B(x_0,r)$ for the open ball of radius $r$ centered at $x_0$.
We likewise use $\bar{B}(x_0,r)$ to denote the corresponding closed ball.

\begin{Lemma}\label{Lipschitzlem}
Let $X$ and $Y$ be two complex Banach spaces. Suppose $r_1>0$ and $r_2\ge 0$.
Let $x_0\in X$ and $y_0\in Y$, and let $f$ be an analytic map
\[
\begin{array}{cccc}
f: & B(x_0,r_1) & \longrightarrow & \bar{B}(y_0,r_2) 
\end{array}\ .
\]
Let $\nu\in (0,\frac{1}{2})$, then for any $x_1,x_2\in \bar{B}(x_0,\nu r_1)$
\[
||f(x_1)-f(x_2)||\le \frac{r_2(1-\nu)}{r_1(1-2\nu)}\ ||x_1-x_2||\ .
\]
\end{Lemma}

\noindent{\bf Proof:}
Suppose $x_1\neq x_2$ satisfy the hypothesis of the proposition.
For $z\in\mathbb{C}$ define
\[
g(z)=f\left(\frac{x_1+x_2}{2}+z\frac{x_1-x_2}{2}\right)-y_0\ .
\]
We first find a bound on $|z|$ which garantees that the argument of $f$ is in the ball
$B(x_0,r_1)$.
Since $\nu<\frac{1}{2}$, we have
\[
2r_1(1-\nu)>2\nu r_1\ge ||x_1-x_0||+||x_2-x_0||\ge||x_1-x_2||\ .
\]
Therefore 
\[
R_{\rm max}=\frac{2r_1(1-\nu)}{||x_1-x_2||}>1\ .
\]
Now the open interval $(1,R_{\rm max})$ is nonempty, and for any $R$ in this interval
as well as for any $z$ with $|z|\le R$ we have
\[
||\frac{x_1+x_2}{2}+z\frac{x_1-x_2}{2}-x_0||\le
\nu r_1+\frac{R}{2}||x_1-x_2||<r_1\ .
\]
Let $\gamma$ be the circle of radius $R$ around the origin in the complex plane.
For such an $R\in (1,R_{\rm max})$
we have by Cauchy's Theorem
\[
f(x_1)-f(x_2)= g(1)-g(-1)=
\frac{1}{\pi i}\oint_{\gamma} \frac{g(z)}{z^2-1}\ dz \ .
\]
Hence
\[
||f(x_1)-f(x_2)||  \le \frac{1}{\pi}\times 2\pi R r_2\times
\max_{|z|=R} \frac{1}{|z^2-1|}
= \frac{2R r_2}{R^2-1}\ .
\]
We now minimize this bound with respect to $R\in (1,R_{\rm max})$.
Since $R\mapsto \frac{2R}{R^2-1}$ is decreasing on $(1,\infty)$,
\[
\inf_{R\in (1,R_{\rm max})}\ \frac{2R}{R^2-1}=\frac{2R_{\rm max}}{R_{\rm max}^2-1}\ .
\]
Inserting the formula for $R_{\rm max}$ in the upper bound for $||f(x_1)-f(x_2)||$ and
simplifying the resulting expression gives the desired Lipschitz estimate.
\qed

%% file: Formalres.tex
\section{Formal statement of the results}\label{formalstatsec}

Now let us pick $d=3$ and for $0<\epsilon<1$ let us denote the quantity $\frac{3-\epsilon}{4}$ by the symbol $[\phi]$.
Let $L=p^l$ for some integer $l\ge 1$. For $r\in\mathbb{Z}$ (typically negative), we consider
the bilinear form on $S(\mathbb{Q}_p^3)$ given by
\[
C_r(f,g)=\int_{\mathbb{Q}_p^3} 
\frac{\widehat{f}(-k)\widehat{g}(k)\bbone\{|k|\le L^{-r}\}}{|k|^{3-2[\phi]}}\ {\rm d}^3 k
\]
where we use $\bbone\{\cdots\}$ for the characteristic function of the condition between braces.
By Theorem \ref{BochnerMthm}, there is a unique probability measure $\mu_{C_r}$
on $S'(\mathbb{Q}_p^3)$ such that for any $f\in S(\mathbb{Q}_p^3)$
\[
\left\langle e^{i\phi(f)}\right\rangle_{\mu_{C_r}}=\exp\left(-\frac{1}{2}C_r(f,f)\right)
\]
where we used the statistical mechanics notation for the expectation with respect to $\phi$
sampled according to the measure $\mu_{C_r}$.
Note that one can write, with a slight abuse of notation
\[
C_r(f,g)=\int_{\mathbb{Q}_p^{3\times 2}}
\ C_r(x-y)f(x)g(y)\ {\rm d}^3 x\ {\rm d}^3 y
\]
where the function $C_r$ is explicitly given by
\[
C_r(x)=\sum_{n=lr}^{\infty}
p^{-2n[\phi]}\left[
\bbone_{\mathbb{Z}_p^3}(p^n x)-p^{-3}\bbone_{\mathbb{Z}_p^3}(p^{n+1} x)
\right]\ .
\]

The measure $\mu_{C_r}$ is supported on distributions given by bonafide functions which are locally constant at scale $L^r$, namely,
constant on each coset in $\mathbb{Q}_p^3/(L^{-r}\mathbb{Z}_p)^3$, the latter quotient playing the role of the lattice of mesh $L^r$.
Note that since $|p|=p^{-1}$ where the $p$ on the left wears its $p$-adic hat while the one on the right is viewed as a real number,
the volume of $(L^{-r}\mathbb{Z}_p)^3$ is $L^{3r}$ in accordance with the intuitive image of a three-dimensional
box with linear dimension $L^r$.
For $s\in\mathbb{Z}$ (typically positive), we use the notation $\Lambda_s=\{x\in\mathbb{Q}_p^3|\ |x|\le L^s\}$ and we also define the Wick
powers
\[
:\phi^2:_{C_r}(x)=\phi(x)^2-C_r(0)\ ,
\]
\[
:\phi^4:_{C_r}(x)=\phi(x)^4-6\ C_r(0)\ \phi(x)^2+3\ C_r(0)^2
\]
and, given $g>0$ as well as $\mu\in\mathbb{R}$, the potential
\[
\tilde{V}_{r,s}(\phi)=\int_{\Lambda_s}
\left\{
L^{-(3-4[\phi])r}\ g :\phi^4:_{C_r}(x)+L^{-(3-2[\phi])r}\ \mu:\phi^2:_{C_r}(x)
\right\}\ {\rm d}^3 x\ .
\]
By the previous remarks, the measure 
\[
{\rm d}\nu_{r,s}(\phi)=\frac{1}{\mathcal{Z}_{r,s}}e^{-\tilde{V}_{r,s}(\phi)}{\rm d}\mu_{C_r}(\phi)
\]
is a well defined probability measure on $S'(\mathbb{Q}_p^3)$ with finite moments.
The normalization factor $\mathcal{Z}_{r,s}$ is at least equal to one as can be seen from Jensen's inequality.
We will denote expectations with respect to $\nu_{r,s}$ by $\langle\cdots\rangle_{r,s}$.
Finally, given a locally constant $\phi$ at scale $L^r$,
we define an element $N_{r}[\phi^2]$ of $S'(\mathbb{Q}_p^3)$
by letting it act on $j\in S(\mathbb{Q}_p^3)$ via
\[
N_{r}[\phi^2](j)=Z_2^{\ r} \int_{\mathbb{Q}_p^3}
\left(
Y_2 :\phi^2:_{C_r}(x)-Y_0\ L^{-2r[\phi]}
\right)\ j(x)\ {\rm d}^3 x
\]
where $Z_2,Y_0, Y_2$ are parameters used in the construction.

We will also use the notation
\[
\bar{g}_\ast=\frac{(p^\epsilon-1)}{36\ L^{\epsilon}(1-p^{-3})}\ .
\]

The main result of this article is the following theorem.

\begin{Theorem}
\label{themainthm}
\ 

$\exists \rho>0$, $\exists L_0$, $\forall L\ge L_0$,
$\exists \epsilon_0>0$, $\forall\epsilon\in (0,\epsilon_0]$,
one can find $\eta_{\phi^2}>0$ and functions $\mu(g),Y_0(g),Y_2(g)$
of $g$ in the interval $(\bar{g}_\ast-\rho\epsilon^{\frac{3}{2}},\bar{g}_\ast+\rho\epsilon^{\frac{3}{2}})$,
such that if one sets $\mu=\mu(g)$, $Z_2=L^{-\frac{1}{2}\eta_{\phi^2}}$, $Y_0=Y_0(g)$
and $Y_2=Y_2(g)$ in the previous definitions, then for all collections of test functions
$f_1,\ldots,f_n,j_1,\ldots,j_m$, the limits
\[
\lim_{\substack{r\rightarrow -\infty\\ s\rightarrow\infty}}
\left\langle
\phi(f_1)\cdots\phi(f_n)
N_r[\phi^2](j_1)\cdots N_r[\phi^2](j_m)
\right\rangle_{r,s}
\] 
exist and do not depend on the
order in which the $r\rightarrow -\infty$ and $s\rightarrow\infty$ limits are taken.
Moreover, the resulting quantities or correlators henceforth similarly and formally denoted by dropping the $r$ and $s$ subscripts
(and using squares, 4-th powers, etc., for repeats)
satisfy the following properties:

1) They are invariant by translation and rotation.

2) They satisfy the partial scale invariance property
\[
\left\langle
\phi(\lambda\cdot f_1)\cdots\phi(\lambda\cdot f_n)
\ N[\phi^2](\lambda\cdot j_1)\cdots N[\phi^2](\lambda\cdot j_m)
\right\rangle=\qquad\qquad\qquad\qquad
\]
\[
\qquad\qquad\qquad|\lambda|^{(3-[\phi])n+(3-2[\phi]-\frac{1}{2}\eta_{\phi^2})m}
\left\langle
\phi(f_1)\cdots\phi(f_n)
\ N[\phi^2](j_1)\cdots N[\phi^2](j_m)
\right\rangle
\]
for all $\lambda\in L^{\mathbb{Z}}$.

3) They satisfy the nontriviality conditions
\[
\langle \phi(\bbone_{\mathbb{Z}_p^3})^4
\rangle - 3\langle \phi(\bbone_{\mathbb{Z}_p^3})^2
\rangle <0\ ,
\]
\[
\langle N[\phi^2](\bbone_{\mathbb{Z}_p^3})^2
\rangle=1\ .
\]

4) The pure $\phi$ correlators are the moments of a unique probability measure $\nu_{\phi}$ on $S'(\mathbb{Q}_p^2)$
with finite moments. This measure
is translation and rotation invariant. It is also partially scale invariant with homogeneity $-[\phi]$ with respect to
the scaling subgroup $L^{\mathbb{Z}}$.

5) The pure $N[\phi^2]$ correlators are the moments of a unique probability measure $\nu_{\phi^2}$ on $S'(\mathbb{Q}_p^2)$
with finite moments. This measure
is translation and rotation invariant. It is also partially scale invariant with homogeneity $-2[\phi]-\frac{1}{2}\eta_{\phi^2}$
with respect to the scaling subgroup $L^{\mathbb{Z}}$.

6) The measures $\nu_\phi$ and $\nu_{\phi^2}$ satisfy a mild form of universality: they do not depend on $g$ in the above-mentioned interval.
\end{Theorem}

%% file: DefextRG.tex
\section{Definition of the extended RG}\label{defBYliftsec}

\subsection{Functional spaces}\label{funcspsec}

In this section we will introduce the different spaces on which the RG transformation will act.
The basic space we will use is $C_{{\rm bd}}^9(\mathbb{R},\mathbb{C})$, namely,
the space of nine times continuously differentiable functions from $\mathbb{R}$ to $\mathbb{C}$
which, together with their derivatives up to order nine, are bounded.
On this space we will use seminorms $||\cdot||_{\partial\phi,\psi,\theta}$ defined for
$K\in C_{{\rm bd}}^9(\mathbb{R},\mathbb{C})$ by
\[
||K(\phi)||_{\partial\phi,\psi,\theta}=\sum_{j=0}^{9}\frac{\theta^j}{j!}\left|
\frac{{\rm d}^j K}{{\rm d}\phi^j}(\psi)\right|\ .
\]
Here $\partial\phi$ is merely a symbol which indicates the variable with respect to which the derivatives are taken.
This will be especially useful when the function may depend on several such variables.
By contrast, $\psi$ is an argument of the seminorm. The derivatives are evaluated at $\phi=\psi$ and therefore the result
depends on $\psi$. Finally $\theta\in [0,\infty)$ is a parameter used to properly calibrate this seminorm.
We will mainly use two values for this parameter denoted by $h$ and $h_\ast$ to be specified later.
As an example of use of the previous notation, we have $||\phi^2||_{\partial\phi,\psi,\theta}=|\psi|^2+2\theta|\psi|+\theta^2$.
In the important special case where $\psi=0$, we will abbreviate the notation into
\[
|K(\phi)|_{\partial\phi,\theta}=||K(\phi)||_{\partial\phi,0,\theta}\ .
\]
A nice property of these seminorms is multiplicativity.
Indeed for any two functions $K_1$, $K_2$ in 
$C_{{\rm bd}}^9(\mathbb{R},\mathbb{C})$ 
we have
\[
||K_1(\phi)K_2(\phi)||_{\partial\phi,\psi,\theta}\le
||K_1(\phi)||_{\partial\phi,\psi,\theta}
\times||K_2(\phi)||_{\partial\phi,\psi,\theta}
\]
which is an easy consequence of the Leibniz rule and the choice of $\frac{1}{j!}$ weights.

To a parameter $\bar{g}>0$ called a calibrator we associate a norm $|||\cdot|||_{\bar{g}}$ on the complex Banach
space $C_{{\rm bd}}^9(\mathbb{R},\mathbb{C})$
defined by
\[
|||K|||_{\bar{g}}=\max\left\{
|K(\phi)|_{\partial\phi,h_\ast},
\bar{g}^2 \sup_{\phi\in\mathbb{R}} 
||K(\phi)||_{\partial\phi,\phi,h}
\right\}\ .
\] 
The parameter $h>0$ will be a function of $\bar{g}$ while $h_\ast>0$ only depends on $L$ and $\epsilon$.
We also introduce the notation $C_{{\rm bd,ev}}^9(\mathbb{R},\mathbb{C})$
for the closed subspace of $C_{{\rm bd}}^9(\mathbb{R},\mathbb{C})$ made of functions $K$  with the even
symmetry $K(-\phi)=K(\phi)$.

We use the notation $\mathbb{L}_q=\mathbb{Q}_p^3 / (L^{-q}\mathbb{Z}_p)^3$ for the lattice with
mesh $L^q$. The unit lattice $\mathbb{L}_0$ will simply be denoted by $\mathbb{L}$.
We will typically denote an element of the latter
by $\Delta$. We will call such an element a unit cube, a unit block or simply a box.
Elements of $\mathbb{L}_1$ will be called $L$-blocks.
These are all of the form $L^{-1}\Delta$ for some $\Delta\in\mathbb{L}$.
We will denote by $[L^{-1}\Delta]$ the set of unit blocks contained in the $L$-block
$L^{-1}\Delta$. For $x\in\mathbb{Q}_p^3$ we denote by $\Delta(x)$ the unique box in $\mathbb{L}$ which contains $x$.

We will need to work with complex valued test functions on $\mathbb{Q}_{p}^{3}$ belonging to suitable
finite-dimensional subspaces of $S(\mathbb{Q}_{p}^{3},\mathbb{C})$.
For any $q_{-},q_{+}\in\mathbb{Z}$ such that $q_{-}\le q_{+}$
we define the space $S_{q_{-},q_{+}}(\mathbb{Q}_{p}^{3},\mathbb{C})$
of test functions with support in $\Lambda_{q_{+}}$ and which are locally constant at scale $L^{q_{-}}$.
Namely, a complex-valued test function $f$ belongs to $S_{q_{-},q_{+}}(\mathbb{Q}_{p}^{3},\mathbb{C})$
if and only if
\[
\forall x\in\mathbb{Q}_{p}^{3},\ f(x)\neq 0\ \Longrightarrow\ |x|\le L^{q_{+}}
\]
and
\[
\forall x,y\in\mathbb{Q}_{p}^{3},\ |x-y|\le L^{q_{-}}\ \Longrightarrow\ f(x)=f(y)\ .
\]
By a trivial compacity argument, it is clear that $S(\mathbb{Q}_{p}^{3},\mathbb{C})$ is the union of all the
subspaces $S_{q_{-},q_{+}}(\mathbb{Q}_{p}^{3},\mathbb{C})$. The latter have complex dimension $L^{3(q_{+}-q_{-})}$. 

The most general RG transformation considered in this article is denoted by $RG_{\rm ex}$ and is called the extended RG map.
It will be defined as a transformation $\vec{V}\mapsto\vec{V}'$ on the Banach space $\mathcal{E}_{\rm ex}$
defined as follows.
An element of that space is an indexed family
\[
\vec{V}=(V_\Delta)_{\Delta\in\mathbb{L}}
\]
where
\[
V_\Delta=(\beta_{4,\Delta},\beta_{3,\Delta},\beta_{2,\Delta},\beta_{1,\Delta},
W_{5,\Delta},W_{6,\Delta}, f_\Delta, R_\Delta)\in \mathbb{C}^7\times C_{{\rm bd}}^9(\mathbb{R},\mathbb{C})\ .
\]
Given suitable positive exponents $e_1,e_2,e_3,e_4,e_W,e_R$ we define the norm
\[
||V_\Delta||=\max\left\{
|\beta_{4,\Delta}|\bar{g}^{-e_4},
|\beta_{3,\Delta}|\bar{g}^{-e_3},
|\beta_{2,\Delta}|\bar{g}^{-e_2},
|\beta_{1,\Delta}|\bar{g}^{-e_1},
\right.\qquad\qquad\qquad
\]
\[
\qquad\qquad\qquad\left.
|W_{5,\Delta}|\bar{g}^{-e_W},
|W_{6,\Delta}|\bar{g}^{-e_W},
|f_\Delta| L^{(3-[\phi])}, |||R_\Delta|||_{\bar{g}}\ \bar{g}^{-e_R}
\right\}
\]
for the components living in a box $\Delta$, and
\[
||\vec{V}||=\sup_{\Delta\in\mathbb{L}} ||V_\Delta||
\]
for the whole infinite vector.
Now $\mathcal{E}_{\rm ex}$ is by definition the Banach space
\[
\mathcal{E}_{\rm ex}=\left\{
\vec{V}\in\prod_{\Delta\in\mathbb{L}}\left(\mathbb{C}^7\times C_{{\rm bd}}^9(\mathbb{R},\mathbb{C}) \right)
\ |\ ||\vec{V}||<\infty
\right\}\ .
\]

We let $\mathcal{E}_{\rm bk}$ denote the closed subspace of $\mathcal{E}_{\rm ex}$ made of vectors $\vec{V}=(V_\Delta)_{\Delta\in\mathbb{L}}$
such that $V_\Delta$ is constant with respect to $\Delta$. This corresponds to uniform potentials which are relevent for the bulk RG evolution.
We have a canonical isometric identification between $\mathcal{E}_{\rm bk}$ and the one-box space $\mathcal{E}_{\rm 1B}=\mathbb{C}^7\times
C_{{\rm bd}}^9(\mathbb{R},\mathbb{C})$. 
By definition $\mathcal{E}=\mathbb{C}^2\times C_{{\rm bd,ev}}^9(\mathbb{R},\mathbb{C})$ and consists of elements of the form $(g,\mu,R)$, with $R$ even,
to which one canonically associates
\[
V=(g,0,\mu,0,0,0,0,R)
\]
in $\mathcal{E}_{\rm 1B}$.
Finally we define $\mathcal{E}_{\rm pt}$ to be the closed subspace
of $\mathcal{E}_{\rm ex}$ given by vectors $(V_\Delta)_{\Delta\in\mathbb{L}}$
such that $V_\Delta=0$ if $\Delta\neq \Delta(0)$.
This corresponds to point-like perturbations living at the origin.
It is easy to see that $\mathcal{E}_{\rm bk}$ and $\mathcal{E}_{\rm pt}$ are in direct sum inside $\mathcal{E}_{\rm ex}$.
We will later see the that $\mathcal{E}_{\rm bk}\oplus\mathcal{E}_{\rm pt}$ is stable by the extended RG.

\subsection{Algebraic definition of $RG_{\rm ex}$}\label{algdefsec}

We now proceed with the presentation of the formulas which express the map $\vec{V}\mapsto \vec{V}'=RG_{\rm ex}[\vec{V}]$.
We will also define a collection $\delta b[\vec{V}]=(\delta b_\Delta[\vec{V}])_{\Delta\in\mathbb{L}}\in \mathbb{C}^{\mathbb{L}}$
of field independent quantitites also called vacuum contributions.

We will need the following notations for scaling shifts.
For a field $\phi$ living on $\mathbb{Q}_p^3$ and for any $q\in\mathbb{Z}$
we let
\[
\phi_{\leadsto q}(x) =
L^{-[\phi]q} \phi(L^q x)\ .
\]
For a test function $f$ which is typically paired with a $\phi$ field
we write
\[
f_{\rightarrow q}(x) =
L^{-(3-[\phi])q} f(L^q x)\ .
\]
Finally for a test function $j$ which is typically paired with a $\phi^2$ field
we write
\[
j_{\Rightarrow q}(x) =
L^{-(3-2[\phi])q} j(L^q x)\ .
\]

At the beginning of the RG analysis one needs to do a rescaling which transforms
the ultraviolet scale $L^r$ into $L^0=1$, i.e., unit scale.
Quantities resulting from this initial rescaling will typically be denoted without tildes
whereas the original or native quantities will usually be denoted by tildes.
Given $g$ and $\mu$ as in \S\ref{formalstatsec} we use the notation
\[
\tilde{g}_r=L^{-(3-4[\phi])r}\ g\qquad{\rm and}\qquad
\tilde{\mu}_r=L^{-(3-2[\phi])r}\ \mu\ .
\]
Then for a field $\tilde{\phi}$ which is locally constant at scale $L^r$ one has
\[
\tilde{V}_{r,s}(\tilde{\phi}) = \int_{\Lambda_{s}} {\rm d}^3 x \left[
  \tilde{g}_r :\tilde{\phi}^4(x):_{C_r} +
  \tilde{\mu}_r :\tilde{\phi}^2(x):_{C_r} \right]
\]
according to the definition in \S\ref{formalstatsec}.

For test functions $\tilde{f}$ and $\tilde{j}$ in $S_{q_-,q_+}(\mathbb{Q}_p^3,\mathbb{C})$,
our main quantity of interest will be
\[
\mathcal{Z}_{r,s}(\tilde{f},\tilde{j})
=\int_{S'(\mathbb{Q}_p^3)} {\rm d}\mu_{C_r}(\tilde{\phi})
\ \exp\left(
-\tilde{V}_{r,s}(\tilde{\phi})+\tilde{\phi}(\tilde{f})+Y_2 Z_2^r :\tilde{\phi}^2:_{C_r}(\tilde{j})
-Y_0 Z_0^r\int_{\mathbb{Q}_p^3} \tilde{j}(x)\ {\rm d}^3x
\right)
\]
where we used the pairing
\[
:\tilde{\phi}^2:_{C_r}(\tilde{j})=\int_{\mathbb{Q}_p^3} 
:\tilde{\phi}^2:_{C_r}(x)
\ \tilde{j}(x)\ {\rm d}^3x
\]
and where $Z_2,Z_0>0$ and $Y_2,Y_0\in\mathbb{R}$ are yet to be defined.
Indeed, one can obtain the correlators
\[
\left\langle
\tilde{\phi}(\tilde{f}_1)\cdots\tilde{\phi}(\tilde{f}_n)
N_r[\tilde{\phi}^2](\tilde{j}_1)\cdots N_r[\tilde{\phi}^2](\tilde{j}_m)
\right\rangle_{r,s}
=\int_{S'(\mathbb{Q}_p^3)} {\rm d}\nu_{r,s}(\tilde{\phi})\ 
\tilde{\phi}(\tilde{f}_1)\cdots\tilde{\phi}(\tilde{f}_n)
N_r[\tilde{\phi}^2](\tilde{j}_1)\cdots N_r[\tilde{\phi}^2](\tilde{j}_m)
\]
as multiple derivatives at $\tilde{f}=\tilde{j}=0$
of the moment generating function
\[
\mathcal{S}_{r,s}(\tilde{f},\tilde{j})=\frac{\mathcal{Z}_{r,s}(\tilde{f},\tilde{j})}{\mathcal{Z}_{r,s}(0,0)}\ .
\]

Note that if $\tilde{\phi}$ is distributed according to ${\rm d}\mu_{C_{r}}$ then if we define $\phi$ so that $\tilde{\phi} =
\phi_{\leadsto r}$ then we have that $\phi$ is distributed
according to $d \mu_{C_{0}}$. We also have
\[
\tilde{V}_{r,s}(\phi_{\leadsto r}) = \int_{\Lambda_{s-r}} {\rm d}^3 x \left[
  L^{(3-4[\phi])r} \tilde{g}_r :\phi^4(x):_{C_0} +
  L^{(3-2[\phi])r}\tilde{\mu}_r :\phi^2(x):_{C_0} \right]\ .
\]

Therefore by  a simple change of variables from $\tilde{\phi}$ to $\phi$ we have
\[
\mathcal{Z}_{r,s}(\tilde{f},\tilde{j})
=
\exp\left(-Y_0 Z_0^rL^{2[\phi]r}\int_{\mathbb{Q}_p^3} j(x)\ {\rm d}^3x
\right)
\]
\[
\times
\int_{S'(\mathbb{Q}_p^3)} {\rm d}\mu_{C_0}(\phi)
\exp\left(-V_{r,s}(\phi)
+\phi(f)+Y_2 Z_2^r :\phi^2:_{C_0}(j)\right)
\]
where
$f=\tilde{f}_{\rightarrow -r}$ and $j=\tilde{j}_{\Rightarrow -r}$
and
\[
V_{r,s}(\phi)=\int_{\Lambda_{s-r}} {\rm d}^3 x \left[
  g :\phi^4(x):_{C_0} + \mu :\phi^2(x):_{C_0} \right]\ .
\]
When $r\le q_-\le q_+\le s$, the latter functional integral can be written in the form
\[
\int_{S'(\mathbb{Q}_p^3)} {\rm d} \mu_{C_0} (\phi)\ \mathcal{I}_{s-r}[\vec{V}](\phi)
\]
for a suitable vector $\vec{V}$ in $\mathcal{E}_{\rm ex}$.
Here the integrand associated to such a vector is defined as follows.
Note that such a vector can be written in a compact way as
\[
\vec{V} = (\beta_4,...,\beta_1,W_5,W_6,f,R)
\]
where each entry in $\vec{V}$ is an indexed collection of unit box dependent quantities, for
example:
\[
\beta_4 = (\beta_{4,\Delta})_{\Delta \in \mathbb{L}}\ .
\]
Note that we will also make the natural identification between functions
on $\mathbb{Q}_p^3$ which are constant over unit
blocks and $\mathbb{L}$-indexed vectors. For example if $f \in
S(\mathbb{Q}_p^3,\mathbb{C})$ is constant over unit blocks then we can
just as well think of 
$f$ as the collection $(f_{\Delta})_{\Delta \in\mathbb{L}}$ where
$f_{\Delta(x)}=f(x)$ for all $x\in\mathbb{Q}_p^3$. We will also use this device for the $\beta$'s and $W$'s.

The correspondance between $\vec{V}$'s and integrands
is given for any integer $t\ge 0$ as follows:
\[
\mathcal{I}_t[\vec{V}](\phi) = \prod_{\substack{\Delta \in \mathbb{L}
\\ \Delta\subset\Lambda_t}}
\mathcal{I}_{\Delta}[\vec{V}](\phi)
\]
with
\[
\mathcal{I}_{\Delta}[\vec{V}](\phi) =
e^{f_{\Delta}\phi_{\Delta}}\times\{
    \exp \left[ -\beta_{4,\Delta} :\phi_{\Delta}^4:_{C_0} -
      \beta_{3,\Delta} :\phi_{\Delta}^3:_{C_0} - \beta_{2,\Delta}
      :\phi_{\Delta}^2:_{C_0} - \beta_{1,\Delta} :\phi_{\Delta}:_{C_0}
    \right] \times 
\]
\[
\left( 1+ W_{5,\Delta} :\phi_{\Delta}^5:_{C_0} +
      W_{6,\Delta} :\phi_{\Delta}^6:_{C_0} \right) +
    R_{\Delta}(\phi_{\Delta})\}\ .
\]

The RG evolution of $\vec{V}$ and the definition of the field independent
quantities $\delta b$ will gives us, for $t\ge 1$, the following identity:
\[
\int_{S'(\mathbb{Q}_p^3)} {\rm d}\mu_{C_0} (\phi)\ \mathcal{I}_t[\vec{V}](\phi) = \exp \left[
  \frac{1}{2} (f,\Gamma f)_{\Lambda_t} + \sum_{\substack{
\Delta \in \mathbb{L} \\ \Delta\subset\Lambda_{t-1}}} \delta
    b_{\Delta}[\vec{V}] \right]\times \int_{S'(\mathbb{Q}_p^3)} {\rm d}\mu_{C_0} (\phi)
 \ \mathcal{I}_{t-1}\left[RG_{\rm ex}[\vec{V}]\right](\phi)
\]
where we used the notation
\[
(f,\Gamma f)_X=\int_{X^2} {\rm d}^3 x\ {\rm d}^3 y\ 
f(x)\ \Gamma(x-y)\ f(y) 
\]
for any measurable subset $X$ of $\mathbb{Q}_p^3$.

Each choice of $r,\tilde{f},\tilde{j}$ gives us a sequence of vectors
$(\vec{V}^{(r,q)}(\tilde{f},\tilde{j}))_{q \in \mathbb{Z}, r\le q \le s}$. The first such
$\vec{V}$ in the sequence is given by:
\[\vec{V}^{(r,r)}(\tilde{f},\tilde{j}) =
(\beta_4,\beta_3,\beta_2,\beta_1,W_5,W_6,f,R)
\]
where
\begin{eqnarray*}
\beta_3 & = &0\\
\beta_1 & = &0\\
W_5 &=& 0\\
W_6&= &0\\
R &=&0\\
f& = &\tilde{f}_{\rightarrow (-r)}\\
\beta_{4}(x) &=& g\ {\rm for\ all\ }x\\
\beta_{2} (x)& =& \mu-Y_2 Z_2^r\ L^{(3-2[\phi])r}\ \tilde{j}(L^{-r}x)\ {\rm for\ all\ }x\ .
\end{eqnarray*}

With these definitions we have that:
\[
\mathcal{Z}_{r,s}(\tilde{f},\tilde{j}) = 
\exp\left(-Y_0 Z_0^rL^{2[\phi]r}\int_{\mathbb{Q}_p^3} j(x)\ {\rm d}^3x
\right)\times
\int {\rm d}\mu_{C_0}(\phi)
\ \mathcal{I}_{s-r}[\vec{V}^{(r,r)}(\tilde{f},\tilde{j})](\phi)
\]
where from now on we will drop the domain $S'(\mathbb{Q}_p^3)$ of the functional integrals.
We define $\vec{V}^{(r,q)}(\tilde{f},\tilde{j})$ for $r<q\le s$ by iterating the extended RG transformation
$q-r$ times, namely:
\[
\vec{V}^{(r,q)}(\tilde{f},\tilde{j}) = RG_{\rm ex}^{q-r}[\vec{V}^{(r,r)}(\tilde{f},\tilde{j})]\ .
\]
Iterating the RG transformation gives us the following equation for
$\mathcal{Z}_{r,s}(\tilde{f},\tilde{j})$:

\[
\mathcal{Z}_{r,s}(\tilde{f},\tilde{j}) =
\exp\left(-Y_0 Z_0^rL^{2[\phi]r}\int_{\mathbb{Q}_p^3} j(x)\ {\rm d}^3x
\right)
\]
\[
\times
 \exp \left[ \sum_{r \le q \le t}
  \left( \frac{1}{2} (f^{(r,q)}, \Gamma f^{(r,q)})_{\Lambda_{s-q}} +
  \sum_{\substack{ \Delta \in\mathbb{L} \\ \Delta\subset\Lambda_{s-q-1}}} 
\delta b_{\Delta}[\vec{V}^{r,q}(\tilde{f},\tilde{j})] \right)
\right]\times 
\int {\rm d}\mu_{C_0}(\phi)
\ \mathcal{I}_{s-t-1}[\vec{V}^{(r,t+1)}(\tilde{f},\tilde{j})](\phi)
\]
which holds for every scale $t$ with $r\le t<s$.
The notation $f^{(r,q)}$ stands for the $f$ component of
$\vec{V}^{(r,q)}(\tilde{f},\tilde{j})$.

We stop the RG iterations at $t=s-1$, i.e., when the finite volume $\Lambda$ becomes a unit box.
This gives
\[
\mathcal{Z}_{r,s}(\tilde{f},\tilde{j}) = 
\exp\left(-Y_0 Z_0^rL^{2[\phi]r}\int_{\mathbb{Q}_p^3} j(x)\ {\rm d}^3x
\right)\times
\]
\[
\exp \left[ \sum_{r \le q \le s-1}
  \left( \frac{1}{2} (f^{(r,q)}, \Gamma f^{(r,q)})_{\Lambda_{s-q}} + \sum_{
\substack{ \Delta \in \mathbb{L} \\ \Delta\subset\Lambda_{s-q-1}}} 
\delta b_{\Delta}[\vec{V}^{r,q}(\tilde{f})] \right)
\right]\times 
\partial\mathcal{Z}_{r,s}(\tilde{f},\tilde{j})
\]
where the boundary factor is:
\[
\partial\mathcal{Z}_{r,s}(\tilde{f},\tilde{j})=
\int {\rm d}\mu_{C_0}(\phi)
\ \mathcal{I}_{0}[\vec{V}^{(r,s)}(\tilde{f},\tilde{j})](\phi)
\]
which reduces to an integral over a single real variable.

The full RG transformation $\vec{V}\rightarrow \vec{V}'=RG_{\rm ex}[\vec{V}]$
will be defined by specifying
\[
\vec{V}' = (\beta'_4,...,\beta'_1,W'_5,W'_6,f',R')
\]
starting from the analogous unprimed quantities.
We will also define the corresponding $\delta b=\delta b[\vec{V}]$ at the same time.

The easiest terms are given by:

\begin{equation}
f'_{\Delta'}=L^{3-[\phi]}
\begin{array}{c}
\ \\
{\rm avg}\\
{\scriptstyle \Delta\in[L^{-1}\Delta']}
\end{array}
f_{\Delta}
\label{newfdefeq}
\end{equation}
where ``${\rm avg}$'' means the average.

We also have
\[
W'_{6,\Delta'}=
L^{3-6[\phi]}
\begin{array}{c}
\ \\
{\rm avg}\\
{\scriptstyle \Delta\in[L^{-1}\Delta']}
\end{array}
W_{6,\Delta}
+8 L^{-6[\phi]}
\parbox{1.5cm}{
\psfrag{b}{$\beta_4$}
\raisebox{-5ex}{
\includegraphics[width=1.6cm]{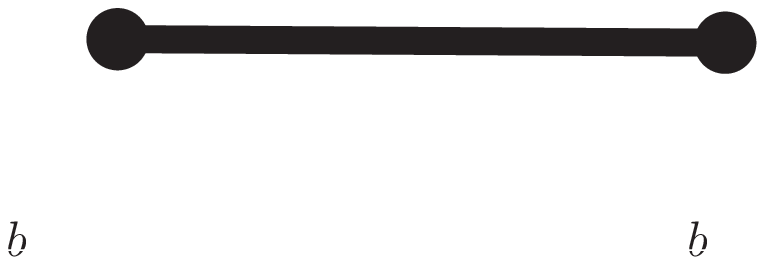}}
}
\]
and
\[
W'_{5,\Delta'}=
L^{3-5[\phi]}
\begin{array}{c}
\ \\
{\rm avg}\\
{\scriptstyle \Delta\in[L^{-1}\Delta']}
\end{array}
W_{5,\Delta}
+ 6 L^{-5[\phi]}
\parbox{1cm}{
\psfrag{w}{$W_6$}\psfrag{f}{$f$}
\raisebox{1ex}{
\includegraphics[width=1cm]{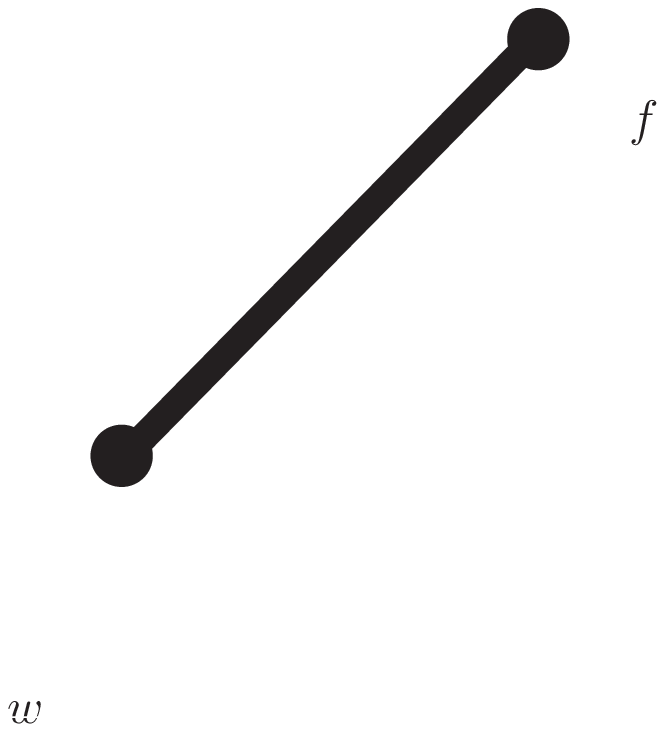}}
}
+ 12 L^{-5[\phi]}
\parbox{1.5cm}{
\psfrag{a}{$\beta_4$}\psfrag{b}{$\beta_3$}
\raisebox{-5ex}{
\includegraphics[width=1.6cm]{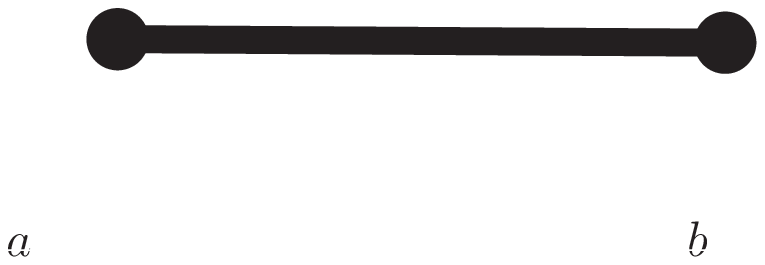}}
}\ \ 
+ 48 L^{-5[\phi]}
\parbox{2.1cm}{
\psfrag{b}{$\beta_4$}\psfrag{f}{$f$}
\raisebox{1ex}{
\includegraphics[width=2.1cm]{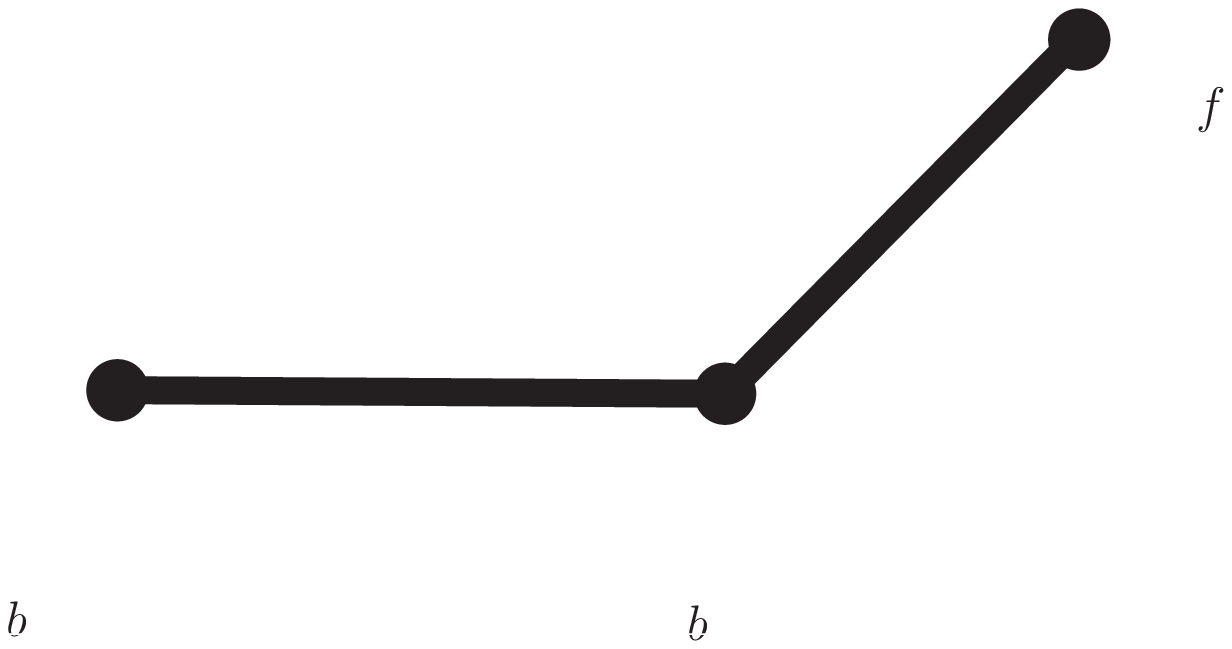}}
}\ \ 
\]
The Feynman diagrams are given explicitly by
\[
\parbox{1.5cm}{
\psfrag{b}{$\beta_4$}
\raisebox{-5ex}{
\includegraphics[width=1.6cm]{Fig1.eps}}
}\ \ 
=\int_{(L^{-1}\Delta')^2} {\rm d}^3 x\ {\rm d}^3 y\ 
\beta_4(x)\ \Gamma(x-y)\ \beta_4(y)
\]
\[
\parbox{1cm}{
\psfrag{w}{$W_6$}\psfrag{f}{$f$}
\raisebox{1ex}{
\includegraphics[width=1cm]{Fig2.eps}}
}\ \ 
=\int_{(L^{-1}\Delta')^2} {\rm d}^3 x\ {\rm d}^3 y\ 
W_6(x)\ \Gamma(x-y)\ f(y)
\]
\[
\parbox{1.5cm}{
\psfrag{a}{$\beta_4$}\psfrag{b}{$\beta_3$}
\raisebox{-5ex}{
\includegraphics[width=1.6cm]{Fig3.eps}}
}\ \ 
=\int_{(L^{-1}\Delta')^2} {\rm d}^3 x\ {\rm d}^3 y\ 
\beta_4(x)\ \Gamma(x-y)\ \beta_3(y)
\]
and
\[
\parbox{2.1cm}{
\psfrag{b}{$\beta_4$}\psfrag{f}{$f$}
\raisebox{1ex}{
\includegraphics[width=2.1cm]{Fig4.eps}}
}\ \ 
=\int_{(L^{-1}\Delta')^3} {\rm d}^3 x\ {\rm d}^3 y\ {\rm d}^3 z\ 
\beta_4(x)\ \Gamma(x-y)\ \beta_4(y)\ \Gamma(y-z)\ f(z)
\]
where we again used the correspondence between vectors indexed
by unit cubes and functions on $\mathbb{Q}_p^3$, i.e., we used $\beta_4(x)=\beta_{4,\Delta(x)}$, etc.

We will need the following intermediate quantities
\[
\hat{\beta}_{k,\Delta'}=L^{3-k[\phi]}
\begin{array}{c}
\ \\
{\rm avg}\\
{\scriptstyle \Delta\in[L^{-1}\Delta']}
\end{array}
\beta_{k,\Delta}
\]
for $1\le k\le 4$.

We will also write
\[
V_\Delta(\phi)=\sum_{k=1}^4 \beta_{k,\Delta} :\phi^k:_{C_0}
\]
\[
Q_\Delta(\phi)=W_{5,\Delta} :\phi^5:_{C_0}+
W_{6,\Delta} :\phi^6:_{C_0}
\]
and
\[
K_\Delta(\phi)=Q_\Delta(\phi) e^{-V_\Delta(\phi)}+ R_\Delta(\phi)\ .
\]
These are functions of a single variable $\phi=\phi_\Delta$.
To lighten the notations we drop the reference to $\Delta$ for the field $\phi$,
when this causes no ambiguity.
We have, using the decomposition of Gaussian measures
\begin{eqnarray*}
\int {\rm d} \mu_{C_0} (\phi)\ \mathcal{I}_t[\vec{V}](\phi) 
& = & \int {\rm d}\mu_{C_0} (\phi) \prod_{\substack{\Delta\in\mathbb{L} \\ \Delta\subset\Lambda_t}}
\left\{ e^{f_\Delta \phi_\Delta}\times \left[
e^{-V_\Delta(\phi_\Delta)}+K_\Delta(\phi_\Delta)
\right]
\right\}\\
 & =& \int {\rm d}\mu_{C_0} (\phi) \int {\rm d}\mu_{\Gamma} (\zeta)
 \prod_{\substack{\Delta\in\mathbb{L} \\ \Delta\subset\Lambda_t}}
\left\{ e^{f_\Delta \phi_{1,\Delta}+f_\Delta\zeta_\Delta}\times \left[
e^{-V_\Delta(\phi_{1,\Delta}+\zeta_{\Delta})}+K_\Delta(\phi_{1,\Delta}+\zeta_\Delta)
\right]
\right\}
\end{eqnarray*}
where $\phi_1=\phi_{\leadsto 1}$.

We then organize the product according to the $L$-blocks containing $\Delta$
and use the independence of the $\zeta$ random variables living in different
$L$-blocks to obtain
\begin{eqnarray*}
\int {\rm d} \mu_{C_0} (\phi)\ \mathcal{I}_t[\vec{V}](\phi) 
& = & \int {\rm d}\mu_{C_0} (\phi) 
\prod_{\substack{\Delta'\in\mathbb{L} \\ \Delta'\subset\Lambda_{t-1}}}\left(
\int {\rm d}\mu_{\Gamma} (\zeta)
\begin{array}{c}
\ \\
\ \\
\ 
\end{array}
\right. \\
 & & \left. 
\begin{array}{c}
\ \\
\ \\
\ 
\end{array}
\prod_{\Delta\in[L^{-1}\Delta']}
\left\{ e^{f_\Delta \phi_{1,\Delta}+f_\Delta\zeta_\Delta}\times \left[
e^{-V_\Delta(\phi_{1,\Delta}+\zeta_{\Delta})}+K_\Delta(\phi_{1,\Delta}+\zeta_\Delta)
\right]
\right\}\right) \\
 & = & \int {\rm d}\mu_{C_0} (\phi) 
\prod_{\substack{\Delta'\in\mathbb{L} \\ \Delta'\subset\Lambda_{t-1}}}\left(
e^{f'_{\Delta'}\phi_{\Delta'}} \times \mathcal{B}_{\Delta'}\right)
\end{eqnarray*}
where
\[
\mathcal{B}_{\Delta'}=
\int {\rm d}\mu_{\Gamma} (\zeta)
 \prod_{\Delta\in[L^{-1}\Delta']}
\left\{ e^{f_\Delta\zeta_\Delta} \times\left[
e^{-V_\Delta(\phi_{1,\Delta}+\zeta_{\Delta})}+K_\Delta(\phi_{1,\Delta}+\zeta_\Delta)
\right]\right\}\ .
\]
With a slight abuse of notation we define
\[
\tilde{V}_\Delta(\phi_1)=\sum_{k=1}^4
\beta_{k,\Delta} :\phi_1^k:_{C_1}
\]
We also let
\[
\hat{V}_{\Delta'}(\phi)=\sum_{k=1}^4
\hat{\beta}_{k,\Delta'} :\phi^k:_{C_0}\ .
\]
Note that $\sum_{\Delta\in[L^{-1}\Delta']}\tilde{V}_\Delta(\phi_1)= \hat{V}_{\Delta'}(\phi)$
where $\phi$ is in fact the component $\phi_{\Delta'}$ of the field but we suppressed this from the notation.
Now define
\[
p_{\Delta}=p_{\Delta}(\phi_1,\zeta)=V_\Delta(\phi_1+\zeta)-\tilde{V}_\Delta(\phi_1)
\]
namely
\[
p_{\Delta}=\sum_{a,b}
\bbone\left\{
\begin{array}{c}
a+b\le 4\\
a\ge 0\ ,\ b\ge 1
\end{array}
\right\}
\frac{(a+b)!}{a!\ b!}\ \beta_{a+b,\Delta}
:\phi_1^a:_{C_1}\ \times
\ :\zeta^b:_{\Gamma}\ .
\]
Now let
\[
P_\Delta(\phi_1,\zeta)=e^{-V_\Delta(\phi_1+\zeta)}-e^{-\tilde{V}_\Delta(\phi_1)}\ .
\]
We expand $\mathcal{B}_{\Delta'}$ by writing the factors as
\[
e^{-V_\Delta(\phi_{1,\Delta}+\zeta_{\Delta})}+K_\Delta(\phi_{1,\Delta}+\zeta_\Delta)=
e^{-\tilde{V}_\Delta(\phi_1)}+P_\Delta(\phi_1,\zeta)+K_\Delta(\phi_{1,\Delta}+\zeta_\Delta)\ .
\]
This results in
\[
\mathcal{B}_{\Delta'}=
e^{\frac{1}{2}(f,\Gamma f)_{L^{-1}\Delta'}-\hat{V}_{\Delta'}(\phi)}
+\hat{K}_{\Delta'}(\phi)
\]
where
\[
\hat{K}_{\Delta'}(\phi)=
\sum_{Y_P,Y_K}
\int {\rm d}\mu_{\Gamma} (\zeta)
\ e^{\int_{L^{-1}\Delta'} f\zeta}
\times
\prod_{\substack{\Delta\in[L^{-1}\Delta'] \\ \Delta\notin Y_P\cup Y_K}}\left[
e^{-\tilde{V}_\Delta(\phi_1)}
\right]
\times
\prod_{\Delta\in Y_P}\left[
P_\Delta(\phi_1,\zeta)
\right]
\times
\prod_{\Delta\in Y_K}\left[
K_\Delta(\phi_1+\zeta)
\right]
\]
where the sum is over pairs of disjoint subsets $Y_P$, $Y_K$ of $[L^{-1}\Delta']$
such that at least one of them is nonempty.
We now assume that we are given collections of numbers $\delta\beta_{k,\Delta'}$
for $0\le k\le 4$ and $\Delta'\in\mathbb{L}$.
We will write $\delta b_{\Delta'}=\delta\beta_{0,\Delta'}$
We therefore have
\begin{eqnarray*}
\int {\rm d}\mu_{C_0} (\phi)\ \mathcal{I}_t[\vec{V}](\phi) 
& = & \exp\left(\frac{1}{2}(f,\Gamma f)_{\Lambda_t}+\sum_{\substack{\Delta'\in\mathbb{L} \\
\Delta'\subset\Lambda_{t-1}}} \delta b_{\Delta'}
\right)\times \\
 & & \int {\rm d}\mu_{C_0} (\phi) \prod_{\Delta'\in\mathbb{L}}
\left\{
e^{f'_{\Delta'}\phi_{\Delta'}}
\times\left[
e^{-\hat{V}_{\Delta'}(\phi_{\Delta'})-\delta b_{\Delta'}}+
\hat{K}_{\Delta'}(\phi_{\Delta'})
e^{-\delta b_{\Delta'}-\frac{1}{2}(f,\Gamma f)_{L^{-1}\Delta'}}
\right]
\right\}
\end{eqnarray*}
Define
\[
\delta V_{\Delta'}(\phi)=\sum_{k=0}^4 \delta\beta_{k,\Delta'} :\phi^k:_{C_0} 
\]
and
\[
V'_{\Delta'}(\phi)=\sum_{k=1}^4 (\hat{\beta}_{k,\Delta'}-\delta\beta_{k,\Delta'}) :\phi^k:_{C_0} 
\]
so that
\[
V'_{\Delta'}(\phi)=\hat{V}_{\Delta'}(\phi)-\delta V_{\Delta'}(\phi)+\delta b_{\Delta'}\ .
\]
One can check that
\begin{eqnarray*}
\int {\rm d}\mu_{C_0} (\phi)\ \mathcal{I}_t[\vec{V}](\phi) 
& = & \exp\left(\frac{1}{2}(f,\Gamma f)_{\Lambda_t}+\sum_{\substack{\Delta'\in\mathbb{L} \\ \Delta'\subset\Lambda_{t-1}}} \delta b_{\Delta'}
\right)\times \\
 & & \int {\rm d}\mu_{C_0} (\phi) \prod_{\substack{\Delta'\in\mathbb{L} \\ \Delta'\subset\Lambda_{t-1}}}
\left\{
e^{f'_{\Delta'}\phi_{\Delta'}}
\times\left[
e^{-V'_{\Delta'}(\phi_{\Delta'})}+K'_{\Delta'}(\phi_{\Delta'})
\right]\right\}
\end{eqnarray*}
where
\[
K'_{\Delta'}(\phi)=
e^{-\delta b_{\Delta'}-\frac{1}{2}(f,\Gamma f)_{L^{-1}\Delta'}}
\times\left\{
\hat{K}_{\Delta'}(\phi)-e^{-\hat{V}_{\Delta'}(\phi)
+\frac{1}{2}(f,\Gamma f)_{L^{-1}\Delta'}}
\left(e^{\delta V_{\Delta'}(\phi)}-1\right)
\right\}\ .
\]

We now do the $\lambda$ expansion which introduces a new complex parameter $\lambda$.
Define
\[
r_{1,\Delta}=r_{1,\Delta}(\phi_1,\zeta)
=e^{-\tilde{V}_\Delta(\phi_1)}\left[
e^{-p_\Delta}-1+p_\Delta-\frac{1}{2} p_\Delta^2
\right]
\]
and let
\[
P_\Delta(\lambda,\phi_1,\zeta)=
e^{-\tilde{V}_\Delta(\phi_1)}\left[
-\lambda p_\Delta+\frac{\lambda^2}{2} p_\Delta^2
\right]+\lambda^3 r_{1,\Delta}(\phi_1,\zeta)
\]
so that
\[
\left. P_\Delta(\lambda,\phi_1,\zeta)\right|_{\lambda=1}=
P_\Delta(\phi_1,\zeta)\ .
\]
We also define
\[
K_\Delta(\lambda,\phi_1,\zeta)=
\lambda^2 Q_\Delta(\phi_1+\zeta)\ e^{-\tilde{V}_\Delta(\phi_1)}
+\lambda^3\left[
Q_\Delta(\phi_1+\zeta)\left(e^{-p_\Delta}-1\right) e^{-\tilde{V}_\Delta(\phi_1)}+R_\Delta(\phi_1+\zeta)
\right]
\]
so that
\[
\left. K_\Delta(\lambda,\phi_1,\zeta)\right|_{\lambda=1}=K_\Delta(\phi_1+\zeta)
\]
We use the same expansion formula as before in order to define the $\lambda$-deformation
\[
\hat{K}_{\Delta'}(\lambda,\phi)=
\sum_{Y_P,Y_K}
\int {\rm d}\mu_{\Gamma} (\zeta)
\ e^{\int_{L^{-1}\Delta'} f\zeta}
\times
\prod_{\substack{\Delta\in[L^{-1}\Delta'] \\ \Delta\notin Y_P\cup Y_K}}\left[
e^{-\tilde{V}_\Delta(\phi_1)}
\right]
\times
\prod_{\Delta\in Y_P}\left[
P_\Delta(\lambda,\phi_1,\zeta)
\right]
\times
\prod_{\Delta\in Y_K}\left[
K_\Delta(\lambda,\phi_1,\zeta)
\right]\ .
\]
This is a polynomial expression in $\lambda$ with no constant term. We can write it as
\[
\hat{K}_{\Delta'}(\lambda,\phi)=
\mathsf{A}\lambda
+\mathsf{B}\lambda^2
+\mathsf{C}\lambda^3
+\hat{K}_{\Delta'}^{\ge 4}(\lambda,\phi)
\]
where $\hat{K}_{\Delta'}^{\ge 4}(\lambda,\phi)$ contains the terms of order 4 or more.

We now assume that we are given collections of numbers $\delta\beta_{k,j,\Delta'}$
for $0\le k\le 4$, $1\le j\le 3$ and $\Delta'\in\mathbb{L}$ such that
\[
\delta\beta_{k,\Delta'}=\delta\beta_{k,1,\Delta'}+
\delta\beta_{k,2,\Delta'}+\delta\beta_{k,3,\Delta'}\ .
\]
Define
\[
\delta\beta_{k,\Delta'}(\lambda)=\lambda\ \delta\beta_{k,1,\Delta'}+
\lambda^2\delta\beta_{k,2,\Delta'}+\lambda^3\delta\beta_{k,3,\Delta'}\ .
\]
In particular this defines $\delta b_{\Delta'}(\lambda)=\delta\beta_{0,\Delta'}(\lambda)$.
We also let
\[
\delta V_{\Delta'}(\lambda,\phi)=\sum_{k=0}^4 \delta\beta_{k,\Delta'}(\lambda) :\phi^k:_{C_0}\ . 
\]
Using the same formula as before for $K'_{\Delta'}$, we define
the corresponding $\lambda$-deformation:
\[
K'_{\Delta'}(\lambda,\phi)=
e^{-\delta b_{\Delta'}(\lambda)-\frac{1}{2}(f,\Gamma f)_{L^{-1}\Delta'}}
\times\left\{
\hat{K}_{\Delta'}(\lambda,\phi)-e^{-\hat{V}_{\Delta'}(\phi)
+\frac{1}{2}(f,\Gamma f)_{L^{-1}\Delta'}}
\left(e^{\delta V_{\Delta'}(\lambda,\phi)}-1\right)
\right\}\ .
\]
We again expand this in $\lambda$ up to order 3:
\[
K'_{\Delta'}(\lambda,\phi)=
\mathsf{A}'\lambda
+\mathsf{B}'\lambda^2
+\mathsf{C}'\lambda^3
+O(\lambda^4)\ .
\]
We now choose the order 1 counterterms $\delta\beta_{k,1,\Delta'}$ so that $\mathsf{A}'=0$.
Namely, for any $k$, $0\le k\le 4$, we let
\begin{equation}
\delta\beta_{k,1,\Delta'}=-\sum_b
\bbone\left\{
\begin{array}{c}
k+b\le 4 \\
b\ge 1
\end{array}
\right\}
\frac{(k+b)!}{k!\ b!}\ L^{-k[\phi]}\ 
\parbox{2.1cm}{
\psfrag{a}{$\beta_{k+b}$}\psfrag{f}{$f$}\psfrag{b}{$b$}
\raisebox{-1ex}{
\includegraphics[width=2.1cm]{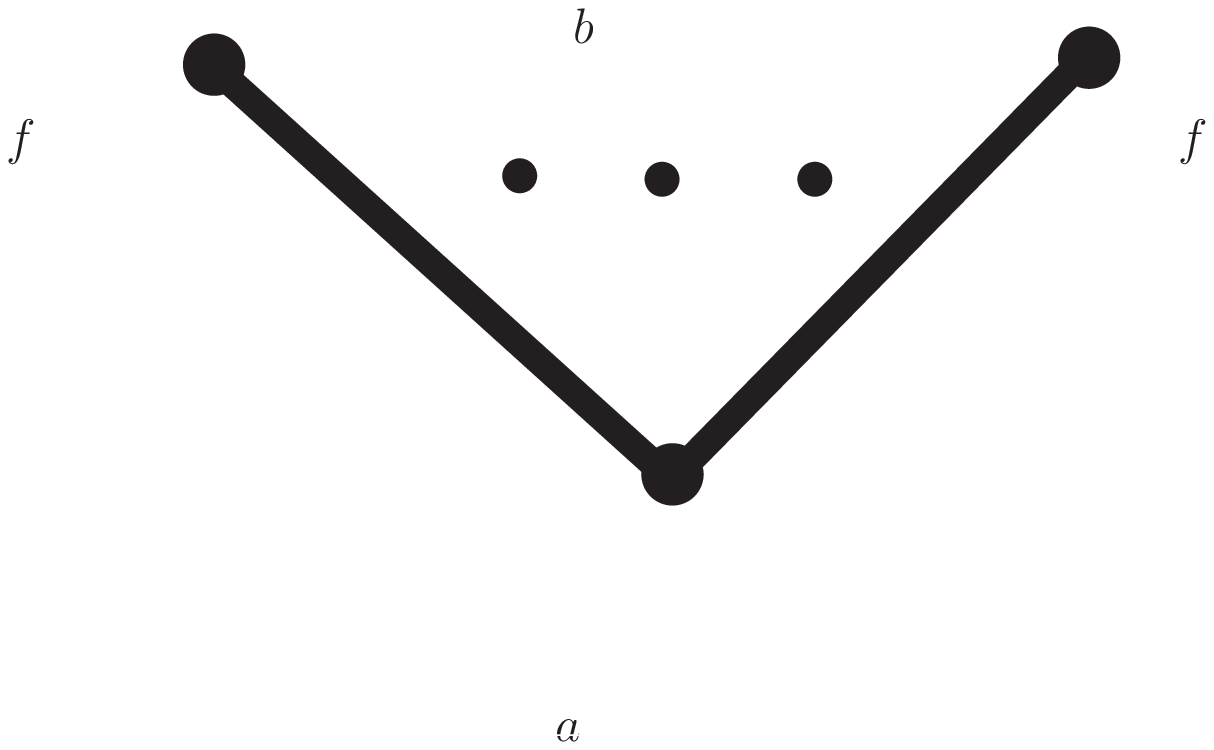}}
}\ \ 
\label{deltb1eq}
\end{equation}
where
\[
\parbox{2.1cm}{
\psfrag{a}{$\beta_{k+b}$}\psfrag{f}{$f$}\psfrag{b}{$b$}
\raisebox{-1ex}{
\includegraphics[width=2.1cm]{Fig5.eps}}
}\ \ 
=\int_{(L^{-1}\Delta')^{b+1}} {\rm d}^3 x\ {\rm d}^3 y_1\cdots {\rm d}^3 y_b\ 
\beta_{k+b}(x)\times\prod_{i=1}^{b}\left[\Gamma(x-y_i)\ f(y_i)\right]\ .
\]

We define the order 2 counterterms $\delta\beta_{k,2,\Delta'}$ so that
\[
\mathsf{B}'=e^{-\hat{V}_{\Delta'}(\phi)} Q'_{\Delta'}(\phi)
\]
where
\[
Q'_{\Delta'}(\phi)=W'_{5,\Delta'} :\phi_{\Delta'}^5:_{C_0}+
W'_{6,\Delta'} :\phi_{\Delta'}^6:_{C_0}
\]
and $W'_{5,\Delta'}$, $W'_{6,\Delta'}$ are new coefficients. The later will turn out to be
the output quantities defined previously.
This hinges on imposing the choice:
\[
\delta\beta_{k,2,\Delta'}=\sum_{a_1,a_2,b_1,b_2,m}
\bbone\left\{
\begin{array}{c}
a_i+b_i\le 4 \\
a_i\ge 0\ ,\ b_i\ge 1\\
1\le m\le \min(b_1,b_2)
\end{array}
\right\}
\frac{(a_1+b_1)!\ (a_2+b_2)!}{a_1!\ a_2!\ m!\ (b_1-m)!\ (b_2-m)!}
\]
\[
\times \frac{1}{2} C(a_1,a_2|k)\times
L^{-(a_1+a_2)[\phi]}\times C_0(0)^{\frac{a_1+a_2-k}{2}}\times
\ \ 
\parbox{4cm}{
\psfrag{a}{$\beta_{a_1+b_1}$}
\psfrag{b}{$\beta_{a_2+b_2}$}
\psfrag{c}{$\scriptstyle{b_1-m}$}
\psfrag{d}{$\scriptstyle{b_2-m}$}
\psfrag{m}{$m$}
\psfrag{f}{$f$}
\raisebox{-1ex}{
\includegraphics[width=4cm]{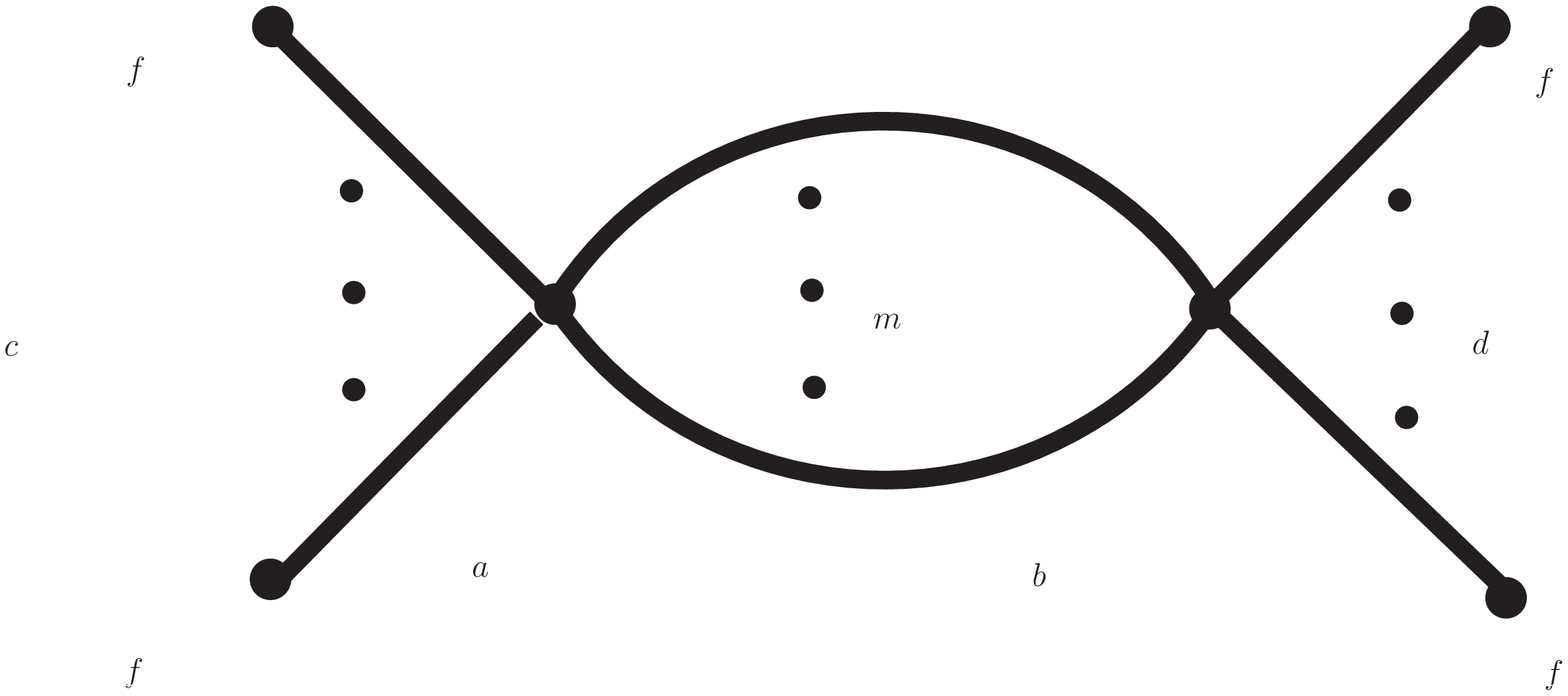}}
}\ \ 
\]
\begin{equation}
+\sum_b
\bbone\left\{
\begin{array}{c}
k+b=5\ {\rm or}\ 6 \\
b\ge 0
\end{array}
\right\}
\frac{(k+b)!}{k!\ b!}\  L^{-k[\phi]}
\ \parbox{2.1cm}{
\psfrag{a}{$W_{k+b}$}\psfrag{f}{$f$}\psfrag{b}{$b$}
\raisebox{-1ex}{
\includegraphics[width=2.1cm]{Fig5.eps}}
}\ \ 
\label{deltb2eq}
\end{equation}
where
\[
\parbox{4cm}{
\psfrag{a}{$\beta_{a_1+b_1}$}
\psfrag{b}{$\beta_{a_2+b_2}$}
\psfrag{c}{$\scriptstyle{b_1-m}$}
\psfrag{d}{$\scriptstyle{b_2-m}$}
\psfrag{m}{$m$}
\psfrag{f}{$f$}
\raisebox{-1ex}{
\includegraphics[width=4cm]{Fig6.eps}}
}\qquad\ 
=
\int_{(L^{-1}\Delta')^{b_1+b_2-2m+2}} {\rm d}^3 x_1\ {\rm d}^3 x_2
\ {\rm d}^3 y_1\cdots {\rm d}^3 y_{b_1-m}\ \ {\rm d}^3 z_1\cdots {\rm d}^3 z_{b_2-m}
\]
\[
\beta_{a_1+b_1}(x_1)\ \beta_{a_2+b_2}(x_2)
\ \Gamma(x_1-x_2)^m\times
\prod_{i=1}^{b_1-m}\left[\Gamma(x_1-y_i)\ f(y_i)\right]\times
\prod_{i=1}^{b_2-m}\left[\Gamma(x_2-z_i)\ f(z_i)\right]
\]
\[
\parbox{2.1cm}{
\psfrag{a}{$W_{k+b}$}\psfrag{f}{$f$}\psfrag{b}{$b$}
\raisebox{-1ex}{
\includegraphics[width=2.1cm]{Fig5.eps}}
}\ \ 
=
\int_{(L^{-1}\Delta')^{b+1}} {\rm d}^3 x\ {\rm d}^3 y_1\cdots {\rm d}^3 y_b\ 
W_{k+b}(x)\times\prod_{i=1}^{b}\left[\Gamma(x-y_i)\ f(y_i)\right]
\]
and where $C(a_1,a_2|k)$
are connection coefficients for Hermite polynomials.
More precisely
\[
C(a_1,a_2|k)=
\bbone\left\{
\begin{array}{c}
|a_1-a_2|\le k\le a_1+a_2 \\
a_1+a_2+k\in 2\mathbb{Z}
\end{array}
\right\}
\times\frac{a_1!\ a_2!}{\left(\frac{a_1+a_2-k}{2}\right)!
\ \left(\frac{a_1+k-a_2}{2}\right)!
\ \left(\frac{a_2+k-a_1}{2}\right)!}\ .
\]
These satisfy the property
\[
:\phi_{\Delta'}^{a_1}:_{C_0}\ \times\ 
:\phi_{\Delta'}^{a_2}:_{C_0}=
\sum_k C(a_1,a_2|k)\ C_0(0)^{\frac{a_1+a_2-k}{2}}
:\phi_{\Delta'}^{k}:_{C_0}\ .
\]

The order 3 counterterms $\delta\beta_{k,3,\Delta'}$ 
will be defined as $(\vec{\beta},f)$-dependent linear functions of $R$.
This is a bit lengthy so we need a few preparatory steps before we can give the explicit
formulas for these counterterms.

First notice that the quantity $\mathsf{C}$ splits as $\mathsf{C}=\mathsf{C}_0+\mathsf{C}_1$
where
\begin{eqnarray*}
\mathsf{C}_0 &=& -\frac{1}{6}
\sum_{\substack{\Delta_1,\Delta_2,\Delta_3\in[L^{-1}\Delta'] \\ {\rm distinct}}}
e^{-\hat{V}_{\Delta'}(\phi)}
\int {\rm d} \mu_{\Gamma} (\zeta)
\ e^{\int_{L^{-1}\Delta'} f\zeta}
\ p_{\Delta_1}\ p_{\Delta_2}\ p_{\Delta_3}\\
 &  & -\frac{1}{2}
\sum_{\substack{\Delta_1,\Delta_2\in[L^{-1}\Delta'] \\ {\rm distinct}}}
e^{-\hat{V}_{\Delta'}(\phi)}
\int {\rm d} \mu_{\Gamma} (\zeta)
\ e^{\int_{L^{-1}\Delta'} f\zeta}
\ p_{\Delta_1}\ p_{\Delta_2}^2\\
 & & -\sum_{\substack{\Delta_1,\Delta_2\in[L^{-1}\Delta'] \\ {\rm distinct}}}
e^{-\hat{V}_{\Delta'}(\phi)}
\int {\rm d} \mu_{\Gamma} (\zeta)
\ e^{\int_{L^{-1}\Delta'} f\zeta}
\ p_{\Delta_1}\ Q_{\Delta_2}(\phi_1+\zeta)\\
 & & +\sum_{\Delta_1\in[L^{-1}\Delta']}
e^{-\hat{V}_{\Delta'}(\phi)}
\int {\rm d} \mu_{\Gamma} (\zeta)
\ e^{\int_{L^{-1}\Delta'} f\zeta}
\ Q_{\Delta_1}(\phi_1+\zeta)\left(e^{-p_{\Delta_1}}-1\right)\\
 & & +\sum_{\Delta_1\in[L^{-1}\Delta']}
\left(\prod_{\substack{\Delta\in[L^{-1}\Delta'] \\ \Delta\neq\Delta_1}}
e^{-\tilde{V}_\Delta(\phi_1)}\right)\times
\int {\rm d} \mu_{\Gamma} (\zeta)
\ e^{\int_{L^{-1}\Delta'} f\zeta}\ r_{1,\Delta_1}
\end{eqnarray*}
and
\[
\mathsf{C}_1=\sum_{\Delta_1\in[L^{-1}\Delta']}
\left(\prod_{\substack{\Delta\in[L^{-1}\Delta'] \\ \Delta\neq\Delta_1}}
e^{-\tilde{V}_\Delta(\phi_1)}\right)\times
\int {\rm d} \mu_{\Gamma} (\zeta)
\ e^{\int_{L^{-1}\Delta'} f\zeta}\ R_{\Delta_1}(\phi_1+\zeta)\ .
\]
Note that we will not need the detailed evaluation of $\mathsf{C}_0$, but we simply need the remark
that it is $R$-independent.

Define
\[
\delta V_{j,\Delta'}(\phi)=\sum_{k=0}^4 \delta\beta_{k,j,\Delta'} :\phi^k:_{C_0} 
\]
for $1\le j\le 3$.
Then the $\lambda^3$ coefficient of $K'_{\Delta'}(\lambda,\phi)$ is given by $\mathsf{C}'=
\mathsf{C}'_0+\mathsf{C}'_1$
where
\[
\mathsf{C}'_0=
e^{-\frac{1}{2}(f,\Gamma f)_{L^{-1}\Delta'}}\ \mathsf{C}_0
-e^{-\hat{V}_{\Delta'}(\phi)}\left(\frac{1}{6}\delta V_{1,\Delta'}(\phi)^3
+\delta V_{1,\Delta'}(\phi)\ \delta V_{2,\Delta'}(\phi) \right)
-e^{-\hat{V}_{\Delta'}(\phi)}\ Q'_{\Delta'}(\phi)\ \delta\beta_{0,1,\Delta'}
\]
and
\[
\mathsf{C}'_1=
e^{-\frac{1}{2}(f,\Gamma f)_{L^{-1}\Delta'}}\ \mathsf{C}_1
-e^{-\hat{V}_{\Delta'}(\phi)} \ \delta V_{3,\Delta'}(\phi)\ .
\]
We now suppose that we are given collections of numbers
$\delta\beta_{k,3,\Delta',\Delta_1}$ for $0\le k\le 4$, $\Delta'\in\mathbb{L}$
and $\Delta_1\in[L^{-1}\Delta']$
such that
\[
\delta\beta_{k,3,\Delta'}=\sum_{\Delta_1\in[L^{-1}\Delta']}
\delta\beta_{k,3,\Delta',\Delta_1}\ .
\]
We then have
\[
\mathsf{C}'_1=\sum_{\Delta_1\in[L^{-1}\Delta']}
\left(\prod_{\substack{\Delta\in[L^{-1}\Delta'] \\ \Delta\neq\Delta_1}}
e^{-\tilde{V}_\Delta(\phi_1)}\right)\times
J_{\Delta',\Delta_1}(\phi)
\]
where
\[
J_{\Delta',\Delta_1}(\phi)=
e^{-\frac{1}{2}(f,\Gamma f)_{L^{-1}\Delta'}}\times
\int {\rm d}\mu_{\Gamma} (\zeta)
\ e^{\int_{L^{-1}\Delta'} f\zeta}\ R_{\Delta_1}(\phi_1+\zeta)
\]
\begin{equation}
-\left(
\sum_{k=0}^4 \delta\beta_{k,3,\Delta',\Delta_1} :\phi^k:_{C_0} 
\right)\times e^{-\tilde{V}_{\Delta_1}(\phi_1)}\ .
\label{jddeq}
\end{equation}
The quantities $\delta\beta_{k,3,\Delta',\Delta_1}$
are uniquely determined by imposing the following normalization conditions on the derivatives up to order 4:
\[
J_{\Delta',\Delta_1}^{(\nu)}(0)=0
\]
for all $\Delta'\in\mathbb{L}$, $\Delta_1\in[L^{-1}\Delta']$
and $\nu$ such that $0\le \nu\le 4$.

Write $J_{\Delta',\Delta_1}(\phi)=J_{+}(\phi)-J_{-}(\phi)$
where
\[
J_{+}(\phi)=e^{-\frac{1}{2}(f,\Gamma f)_{L^{-1}\Delta'}}\times
\int {\rm d} \mu_{\Gamma} (\zeta)
\ e^{\int_{L^{-1}\Delta'} f\zeta}\ R_{\Delta_1}(\phi_1+\zeta)
\]
and
\[
J_{-}(\phi)=
\left(
\sum_{k=0}^4 \delta\beta_{k,3,\Delta',\Delta_1} :\phi^k:_{C_0} 
\right)\times e^{-\tilde{V}_{\Delta_1}(\phi_1)}\ .
\]
For any $\nu$, $0\le \nu\le 4$, we have
\[
J_{+}^{(\nu)}(0)=L^{-\nu[\phi]}\ e^{-\frac{1}{2}(f,\Gamma f)_{L^{-1}\Delta'}}\times
\int {\rm d} \mu_{\Gamma} (\zeta)
\ e^{\int_{L^{-1}\Delta'} f\zeta}\ R_{\Delta_1}^{(\nu)}(\zeta)\ .
\]
Whereas
\[
J_{-}(\phi)=u(\phi)\ e^{v(\phi)}
\]
with
\[
u(\phi)=u_4 \phi^4+u_3 \phi^3+u_2 \phi^2+u_1 \phi+u_0
\]
and
\[
v(\phi)=v_4 \phi^4+v_3 \phi^3+v_2 \phi^2+v_1 \phi+v_0
\]
with coefficients explicitly given by
\begin{eqnarray*}
u_4 & = & \delta\beta_4\\
u_3 & = & \delta\beta_3\\
u_2 & = & \delta\beta_2-6C\delta\beta_4\\
u_1 & = & \delta\beta_1-3C\delta\beta_3\\
u_0 & = & \delta\beta_0-C\delta\beta_2+3C^2\delta\beta_4
\end{eqnarray*}
and
\begin{eqnarray*}
v_4 & = & -L^{-4[\phi]}\beta_4\\
v_3 & = & -L^{-3[\phi]}\beta_3\\
v_2 & = & -L^{-2[\phi]}\beta_2+6CL^{-4[\phi]}\beta_4\\
v_1 & = & -L^{-[\phi]}\beta_1+3CL^{-3[\phi]}\beta_3\\
v_0 & = & CL^{-2[\phi]}\beta_2-3C^2L^{-4[\phi]}\beta_4\ .
\end{eqnarray*}
Note that we used the abbreviated notation $\delta\beta_k=\delta\beta_{k,3,\Delta',\Delta_1}$,
$\beta_k=\beta_{k,\Delta_1}$ and $C=C_0(0)$.
Using Maple we found for the Taylor expansion of $J_{-}(\phi)$ up to order 4:
\begin{eqnarray*}
J_{-}(\phi) & = &
e^{v_0}\times
\left\{
{u_{0}}
+ ({u_{0}}\,{v_{1}} + {u_{1}})\,\phi
+ \left({u_{1}}\,{v_{1}} + {u_{0}}\,{v_{2}}
+ {\displaystyle \frac {1}{2}}\,{u_{0}}\,{v_{1}}^{2} + {u_{2}}\right)\,\phi^{2}
\right. \\
 & & + \left({u_{1}}\,{v_{2}}
+ {\displaystyle \frac {1}{2}}\,{u_{1}}\,{v_{1}}^{2}
+ {u_{2}}\,{v_{1}} + {u_{0}}\,{v_{3}} + {u_{0}}\,{v_{1}}\,{v_{2}}
+ {\displaystyle \frac {1}{6}}\,{u_{0}}\,{v_{1}}^{3} + {u_{3}}\right)\,\phi^{3} \\
 & & + \left({u_{4}} + {u_{1}}\,{v_{3}} + {u_{1}}\,{v_{1}}\,{v_{2}}
+ {\displaystyle \frac {1}{6}}\,{u_{1}}\,{v_{1}}^{3} 
+ {u_{0}}\,{v_{4}} + {u_{0}}\,{v_{1}}\,{v_{3}}\right. \\
 & & \left. \left.
+ {\displaystyle \frac {1}{2}}\,{u_{0}}\,{v_{2}}^{2}
+ {\displaystyle \frac {1}{2}}\,{u_{0}}\,{v_{2}}\,{v_{1}}^{2}
+ {\displaystyle \frac {1}{24}}\,{u_{0}}\,{v_{1}}^{4}
+ {u_{2}}\,{v_{2}} + {\displaystyle \frac {1}{2}}\,{u_{2}}\,{v_{1}}^{2}
+ {u_{3}}\,{v_{1}}\right)\phi^{4} \right\}
+ O(\phi^{5})\ .
\end{eqnarray*}
Write $a_\nu=e^{-v_0} J_{+}^{(\nu)}(0)$.
We therefore have to solve for $u_0,\ldots,u_4$ in the triangular polynomial system
\begin{eqnarray*}
a_0 & = & u_0 \\
a_1 & = & {u_{1}} + {u_{0}}\,{v_{1}} \\
\frac{1}{2}\,a_2 & = & {u_{2}} + {u_{1}}\,{v_{1}} + {u_{0}}\,{v_{2}}
+ {\displaystyle \frac {1}{2}}\,{u_{0}}\,{v_{1}}^{2} \\
\frac{1}{6}\,a_3 & = & {u_{3}} + {u_{1}}\,{v_{2}}
+ {\displaystyle \frac {1}{2}}\,{u_{1}}\,{v_{1}}^{2}
+ {u_{2}}\,{v_{1}} + {u_{0}}\,{v_{3}} + {u_{0}}\,{v_{1}}\,{v_{2}}
+ {\displaystyle \frac {1}{6}}\,{u_{0}}\,{v_{1}}^{3} \\
\frac{1}{24}\,a_4 & = & {u_{4}} + {u_{1}}\,{v_{3}} + {u_{1}}\,{v_{1}}\,{v_{2}}
+ {\displaystyle \frac {1}{6}}\,{u_{1}}\,{v_{1}}^{3} 
+ {u_{0}}\,{v_{4}} + {u_{0}}\,{v_{1}}\,{v_{3}} \\
 & & + {\displaystyle \frac {1}{2}}\,{u_{0}}\,{v_{2}}^{2}
+ {\displaystyle \frac {1}{2}}\,{u_{0}}\,{v_{2}}\,{v_{1}}^{2}
+ {\displaystyle \frac {1}{24}}\,{u_{0}}\,{v_{1}}^{4}
+ {u_{2}}\,{v_{2}} + {\displaystyle \frac {1}{2}}\,{u_{2}}\,{v_{1}}^{2}
+ {u_{3}}\,{v_{1}}\ .
\end{eqnarray*}

This is straightforward but leads to complicated intermediate formulas which we skip.
We then replace the $v's$ by their expressions in terms of the $\beta$'s.
Finally we use the obtained formulas for the $u$'s in order to get
\begin{eqnarray*}
\delta\beta_4 & = & u_4 \\
\delta\beta_3 & = & u_3 \\
\delta\beta_2 & = & u_2+6C u_4 \\
\delta\beta_1 & = & u_1+3C u_3 \\
\delta\beta_0 & = & u_0 + Cu_2 +3C^2 u_4\ .
\end{eqnarray*}

The final result, obtained with the help of Maple, and using the notation $d_k=L^{-k[\phi]}\beta_k$ is:
\begin{eqnarray*}
\delta\beta_4 & = &
{\displaystyle \frac {1}{24}} \,{a_{4}} 
+ \left(
{\displaystyle \frac {1}{6}} \,{d_{1}}
- {\displaystyle \frac {1}{2}} \,C\,{d_{3}}\right)\,{a_{3}}
+ \left(
{\displaystyle \frac {1}{4}} \,{d_{1}}^{2} 
- {\displaystyle \frac {3}{2}} \,C\,{d_{1}}\,{d_{3}} 
+ {\displaystyle \frac {9}{4}} \,C^{2}\,{d_{3}}^{2}
- 3\,C\,{d_{4}}
+ {\displaystyle \frac {1}{2}} \,{d_{2}}\right)\,{a_{2}} \\
 & & + \left(
{\displaystyle \frac {9}{2}} \,C^{2}\,{d_{1}}\,{d_{3}}^{2} 
- {\displaystyle \frac {3}{2}} \,C\,{d_{1}}^{2}\,{d_{3}}
- 6\,C\,{d_{1}}\,{d_{4}}
- 3\,C\,{d_{3}}\,{d_{2}} \right. \\
 & & \left. 
+ 18\,C^{2}\,{d_{3}}\,{d_{4}}
+ {\displaystyle \frac {1}{6}} \,{d_{1}}^{3}
+ {d_{1}}\,{d_{2}}
- {\displaystyle \frac {9}{2}} \,C^{3}\,{d_{3}}^{3}
+ {d_{3}}\right){a_{1}} \\
 & & + \left(
{d_{4}}
- 6\,C\,{d_{2}}\,{d_{4}}
- {\displaystyle \frac {1}{2}} \,C\,{d_{1}}^{3}\,{d_{3}}
+ {\displaystyle \frac {9}{4}} \,C^{2}\,{d_{1}}^{2}\,{d_{3}}^{2}
- {\displaystyle \frac {9}{2}} \,C^{3}\,{d_{1}}\,{d_{3}}^{3}
- 3\,C\,{d_{1}}^{2}\,{d_{4}} \right. \\
 & &  + {\displaystyle \frac {9}{2}} \,C^{2}\,{d_{3}}^{2}\,{d_{2}}
- 3\,C\,{d_{3}}^{2}
+ {\displaystyle \frac {1}{2}} \,{d_{1}}^{2}\,{d_{2}}
+ 18\,C^{2}\,{d_{1}}\,{d_{3}}\,{d_{4}}
- 3\,C\,{d_{1}}\,{d_{3}}\,{d_{2}} \\
 & & \left.
- 27\,C^{3}\,{d_{3}}^{2}\,{d_{4}}
+ 18\,C^{2}\,{d_{4}}^{2}
+ {d_{3}}\,{d_{1}}
+ {\displaystyle \frac {1}{24}} \,{d_{1}}^{4}
+ {\displaystyle \frac {1}{2}} \,{d_{2}}^{2}
+ {\displaystyle \frac {27}{8}} \,C^{4}\,{d_{3}}^{4}\right){a_{0}} \ ,
\end{eqnarray*}
\vskip1cm

\begin{eqnarray*}
\delta\beta_3 & = &
{\displaystyle \frac {1}{6}} \,{a_{3}}
+ \left({\displaystyle \frac {1}{2}} \,{d_{1}}
- {\displaystyle \frac {3}{2}} \,C\,{d_{3}}\right)\,{a_{2}}
+ \left({d_{2}}
- 6\,C\,{d_{4}}
+ {\displaystyle \frac {1}{2}} \,{d_{1}}^{2}
- 3\,C\,{d_{1}}\,{d_{3}}
+ {\displaystyle \frac {9}{2}} \,C^{2}\,{d_{3}}^{2}\right)\,{a_{1}} \\
 & & + \left(
{\displaystyle \frac {9}{2}} \,C^{2}\,{d_{1}}\,{d_{3}}^{2}
- {\displaystyle \frac {3}{2}} \,C\,{d_{1}}^{2}\,{d_{3}}
- 6\,C\,{d_{1}}\,{d_{4}}
- 3\,C\,{d_{3}}\,{d_{2}} \right. \\
 & & \left.
+ 18\,C^{2}\,{d_{3}}\,{d_{4}}
+ {\displaystyle \frac {1}{6}} \,{d_{1}}^{3}
+ {d_{1}}\,{d_{2}}
- {\displaystyle \frac {9}{2}} \,C^{3}\,{d_{3}}^{3}
+ {d_{3}}\right){a_{0}} \ , 
\end{eqnarray*}

\vskip1cm
\begin{eqnarray*}
\delta\beta_2 & = & {\displaystyle \frac {1}{4}} \,C\,{a_{4}}
+ \left( - 3\,C^{2}\,{d_{3}}
+ C\,{d_{1}}\right)\,{a_{3}}
+ \left( - 9\,C^{2}\,{d_{1}}\,{d_{3}}
+ 3\,C\,{d_{2}}
+ {\displaystyle \frac {3}{2}} \,C\,{d_{1}}^{2}
+ {\displaystyle \frac {1}{2}} 
+ {\displaystyle \frac {27}{2}} \,C^{3}\,{d_{3}}^{2}
- 18\,C^{2}\,{d_{4}}\right)\,{a_{2}} \\
 & & + \left(
108\,C^{3}\,{d_{3}}\,{d_{4}}
+ C\,{d_{1}}^{3} - 36\,C^{2}\,{d_{1}}\,{d_{4}}
- 18\,C^{2}\,{d_{3}}\,{d_{2}}
+ 6\,C\,{d_{1}}\,{d_{2}} \right. \\
 & & \left.
+ {d_{1}}
+ 3\,C\,{d_{3}} 
- 9\,C^{2}\,{d_{1}}^{2}\,{d_{3}}
+ 27\,C^{3}\,{d_{1}}\,{d_{3}}^{2}
- 27\,C^{4}\,{d_{3}}^{3}\right){a_{1}} \\
 & & + \left(
{\displaystyle \frac {27}{2}} \,C^{3}\,{d_{1}}^{2}\,{d_{3}}^{2}
- 27\,C^{4}\,{d_{1}}\,{d_{3}}^{3} 
- 36\,C^{2}\,{d_{2}}\,{d_{4}}
+ {\displaystyle \frac {1}{2}} \,{d_{1}}^{2}
+ {d_{2}}
+ {\displaystyle \frac {1}{4}} \,C\,{d_{1}}^{4}
- 18\,C^{2}\,{d_{1}}^{2}\,{d_{4}}
+ 27\,C^{3}\,{d_{3}}^{2}\,{d_{2}}\right. \\
 & & - 162\,C^{4}\,{d_{3}}^{2}\,{d_{4}} 
+ 108\,C^{3}\,{d_{1}}\,{d_{3}}\,{d_{4}}
+ 3\,C\,{d_{1}}^{2}\,{d_{2}}
+ 3\,C\,{d_{1}}\,{d_{3}}
- {\displaystyle \frac {27}{2}} \,C^{2}\,{d_{3}}^{2}
+ 3\,C\,{d_{2}}^{2}
- 3\,C^{2}\,{d_{1}}^{3}\,{d_{3}} \\
 & & \left.
+ 108\,C^{3}\,{d_{4}}^{2}
+ {\displaystyle \frac {81}{4}} \,C^{5}\,{d_{3}}^{4}
- 18\,C^{2}\,{d_{1}}\,{d_{3}}\,{d_{2}}\right){a_{0}} \ ,
\end{eqnarray*}

\vskip1cm
\begin{eqnarray*}
\delta\beta_1 & = & {\displaystyle \frac {1}{2}} \,C\,{a_{3}}
+ \left({\displaystyle \frac {3}{2}} \,C\,{d_{1}}
- {\displaystyle \frac {9}{2}} \,C^{2}\,{d_{3}}\right)\,{a_{2}}
+ \left(3\,C\,{d_{2}}
+ {\displaystyle \frac {3}{2}} \,C\,{d_{1}}^{2}
+ {\displaystyle \frac {27}{2}} \,C^{3}\,{d_{3}}^{2}
- 18\,C^{2}\,{d_{4}}
+ 1
- 9\,C^{2}\,{d_{1}}\,{d_{3}}\right)\,{a_{1}} \\
 & & + \left(
- {\displaystyle \frac {9}{2}} \,C^{2}\,{d_{1}}^{2}\,{d_{3}}
+ {d_{1}}
- {\displaystyle \frac {27}{2}} \,C^{4}\,{d_{3}}^{3}
+ 54\,C^{3}\,{d_{3}}\,{d_{4}}
+ 3\,C\,{d_{1}}\,{d_{2}}\right. \\
 & & \left.
- 18\,C^{2}\,{d_{1}}\,{d_{4}}
- 9\,C^{2}\,{d_{3}}\,{d_{2}} 
+ {\displaystyle \frac {1}{2}} \,C\,{d_{1}}^{3}
+ {\displaystyle \frac {27}{2}} \,C^{3}\,{d_{1}}\,{d_{3}}^{2}\right){a_{0}} \ ,
\end{eqnarray*}

\vskip1cm
\begin{eqnarray*}
\delta\beta_0 & = & {\displaystyle \frac {1}{8}} \,C^{2}\,{a_{4}}
+ \left( - {\displaystyle \frac {3}{2}} \,C^{3}\,{d_{3}}
+ {\displaystyle \frac {1}{2}} \,C^{2}\,{d_{1}}\right)\,{a_{3}} \\
 & & + \left({\displaystyle \frac {3}{4}} \,C^{2}\,{d_{1}}^{2}
+ {\displaystyle \frac {27}{4}} \,C^{4}\,{d_{3}}^{2}
+ {\displaystyle \frac {3}{2}} \,C^{2}\,{d_{2}}
+ {\displaystyle \frac {1}{2}} \,C
- {\displaystyle \frac {9}{2}} \,C^{3}\,{d_{1}}\,{d_{3}}
- 9\,C^{3}\,{d_{4}}\right)\,{a_{2}} \\
 & & + \left(
{\displaystyle \frac {27}{2}} \,C^{4}\,{d_{1}}\,{d_{3}}^{2}
- 9\,C^{3}\,{d_{3}}\,{d_{2}}
+ 54\,C^{4}\,{d_{3}}\,{d_{4}}
+ C\,{d_{1}}
- 18\,C^{3}\,{d_{1}}\,{d_{4}}
+ {\displaystyle \frac {1}{2}} \,C^{2}\,{d_{1}}^{3}\right. \\
 & & \left.
+ 3\,C^{2}\,{d_{1}}\,{d_{2}}
- {\displaystyle \frac {9}{2}} \,C^{3}\,{d_{1}}^{2}\,{d_{3}} 
- {\displaystyle \frac {27}{2}} \,C^{5}\,{d_{3}}^{3}\right){a_{1}} \\
 & & + \left(
{\displaystyle \frac {3}{2}} \,C^{2}\,{d_{1}}^{2}\,{d_{2}}
- 9\,C^{3}\,{d_{1}}\,{d_{3}}\,{d_{2}}
- {\displaystyle \frac {3}{2}} \,C^{3}\,{d_{1}}^{3}\,{d_{3}}
+ {\displaystyle \frac {27}{4}} \,C^{4}\,{d_{1}}^{2}\,{d_{3}}^{2} 
- {\displaystyle \frac {27}{2}} \,C^{5}\,{d_{1}}\,{d_{3}}^{3}
- 18\,C^{3}\,{d_{2}}\,{d_{4}}
+ {\displaystyle \frac {1}{8}} \,C^{2}\,{d_{1}}^{4}\right. \\
 & & + {\displaystyle \frac {81}{8}} \,C^{6}\,{d_{3}}^{4}
+ C\,{d_{2}}
- 9\,C^{3}\,{d_{1}}^{2}\,{d_{4}}
+ {\displaystyle \frac {27}{2}} \,C^{4}\,{d_{3}}^{2}\,{d_{2}} 
- 81\,C^{5}\,{d_{3}}^{2}\,{d_{4}}
+ 1
+ 54\,C^{4}\,{d_{1}}\,{d_{3}}\,{d_{4}} \\
 & & \left.
- {\displaystyle \frac {9}{2}} \,C^{3}\,{d_{3}}^{2}
+ 54\,C^{4}\,{d_{4}}^{2}
- 3\,C^{2}\,{d_{4}}
+ {\displaystyle \frac {1}{2}} \,C\,{d_{1}}^{2} 
+ {\displaystyle \frac {3}{2}} \,C^{2}\,{d_{2}}^{2}\right){a_{0}} \ .
\end{eqnarray*}

\vskip1cm
Also recall that
\begin{eqnarray}
a_i & = & \exp\left[
-C L^{-2[\phi]}\beta_2+3C^2 L^{-4[\phi]}\beta_4
-\frac{1}{2} (f,\Gamma f)_{L^{-1}\Delta'}
\right] \nonumber \\
 & & \times L^{-i[\phi]}\times \int {\rm d}\mu_{\Gamma}(\zeta)
\ e^{\int_{L^{-1}\Delta'} f\zeta}\ R_{\Delta_1}^{(i)}(\zeta)\ .
\label{aieq}
\end{eqnarray}

Note that we get formulas of the form
\[
\delta\beta_k=\sum_{i=0}^{4} M_{k,i}\ a_i
\]
where
the matrix elements $M_{k,i}$ are given by finite sums of the form
\begin{equation}
M_{k,i}=\sum\ \#\ C^j L^{-(l_1+\cdots+l_n)[\phi]}\ \beta_{l_1}\cdots\beta_{l_n}
\label{mkieq}
\end{equation}
with $j\ge 0$, $n\ge 0$, and $1\le l_m\le 4$ for every $m$, $1\le m\le n$.
Here the symbol $\#$ stands for some purely numerical constants.
Furthermore, the terms which appear satisfy the homogeneity constraint
\begin{equation}
l_1+\cdots+l_n-2j=k-i\ .
\label{mkiconsteq}
\end{equation}
We also have a limitation on the range of allowed $n$'s:
\[
n\le (k-i)+2\left\lfloor\frac{4-k}{2}\right\rfloor\ .
\]

This completes the definition of the $\delta\beta_{k,3,\Delta',\Delta_1}$
and therefore of the order 3 counterterms
\[
\delta\beta_{k,3,\Delta'}=\sum_{\Delta_1\in[L^{-1}\Delta']}
\delta\beta_{k,3,\Delta',\Delta_1}\ .
\]
We also have a complete definition of
\[
K'_{\Delta'}(\lambda,\phi)=\lambda^2 e^{-\hat{V}_{\Delta'}(\phi)} Q'_{\Delta'}(\phi)
+\lambda^3 \mathsf{C}'_0
+\lambda^3 \mathsf{C}'_1
+O(\lambda^4)\ .
\]
We now define
\[
\mathcal{L}_{\Delta'}^{(\vec{\beta},f)}(R)=\mathsf{C}'_1
\]
with the previous choices for the $\delta\beta_{k,3,\Delta',\Delta_1}$.
This makes $\mathcal{L}^{(\vec{\beta},f)}$ a $(\vec{\beta},f)$-dependent linear operator on the
space where $R$ lives.

Now the new couplings $\beta'_{k,\Delta'}$ as well as the quantities $\delta b_{\Delta'}$
are fully defined.
We just need the new $R$.
It is given by
\[
R'_{\Delta'}=\mathcal{L}_{\Delta'}^{(\vec{\beta},f)}(R)
+\xi_{R,\Delta'}(\vec{V})
\]
where the formula for remainder term is
\begin{eqnarray*}
\xi_{R,\Delta'}(\vec{V})(\phi) & = &
\frac{1}{2\pi i}\oint_{\gamma_0} \frac{{\rm d}\lambda}{\lambda^4}
\left. K'_{\Delta'}(\lambda,\phi) \right|_{R=0} \\
 & & + \frac{1}{2\pi i}\oint_{\gamma_{01}} \frac{{\rm d}\lambda}{\lambda^4(\lambda-1)}
K'_{\Delta'}(\lambda,\phi) \\
 & & +\left(e^{-\hat{V}_{\Delta'}(\phi)}-e^{-V'_{\Delta'}(\phi)}\right) Q'_{\Delta'}(\phi)
\end{eqnarray*}
where $\gamma_0$ is any positively oriented contour around $\lambda=0$,
and $\gamma_{01}$ is any positively oriented contour which encloses both $\lambda=0$
and $\lambda=1$.
In~\cite{BMS} the three terms for the remainder are respectively denoted
by $R_{\rm main}$, $R_3$, and $R_4$.
It is important to note that in the first term we set $R=0$, which means that all $\delta\beta_{k,3,\Delta'}$
are set equal to zero.
Also note that $\mathcal{L}^{(\vec{\beta},f)}(R)$ corresponds to the $R_{\rm linear}$ notation in~\cite{BMS}.

To finish setting up the notation
we write for $1\le k\le 4$
\[
\xi_{k,\Delta'}(\vec{V})=-\delta\beta_{k,3,\Delta'}
\]
whereas
\[
\xi_{0,\Delta'}(\vec{V})=\delta\beta_{0,3,\Delta'}
\]
In this way the RG evolution for the couplings is
\[
\beta'_{k,\Delta'}=\hat{\beta}_{k,\Delta'}-\delta\beta_{k,1,\Delta'}
-\delta\beta_{k,2,\Delta'} + \xi_{k,\Delta'}(\vec{V})
\]
for $1\le k\le 4$.
Likewise the collected terms which contribute to the progressive calculation of
$\mathcal{Z}_{r,s}(\tilde{f},\tilde{j})$
are
\[
\delta b_{\Delta'}[\vec{V}]=\delta\beta_{0,1,\Delta'}
+\delta\beta_{0,2,\Delta'} + \xi_{0,\Delta'}(\vec{V})\ .
\]

%% file: Prelimest.tex
\section{Preliminary estimates}

\subsection{Properties of covariances}

In this section we collect some of the properties satisfied by the covariances and
needed in the sequel. Recall that $L=p^l$ where $l$ is an integer $l>0$.

\begin{Lemma}\label{L1lem}
The covariance $\Gamma$ can be expressed pointwise as follows.
\begin{enumerate}
\item
If $|x|\le 1$ then
\[
\Gamma(x)=\frac{1-p^{-3}}{1-p^{-2[\phi]}} (1-L^{-2[\phi]})\ .
\]
\item
If $|x|=p^i$ with $1\le i\le l$, then
\[
\Gamma(x)=-p^{-3+2[\phi]}p^{-2l[\phi]}+
\frac{1-p^{-3+2[\phi]}}{1-p^{-2[\phi]}}(p^{-2i[\phi]}-p^{-2l[\phi]})\ . 
\]
\item
If $|x|>L$ then $\Gamma(x)=0$.
\end{enumerate}
\end{Lemma}

\noindent{\bf Proof:}
Recall that
\[
\Gamma(x)=\sum_{j=0}^{l-1} p^{-2j[\phi]}\left(
\bbone_{\mathbb{Z}_p^3}(p^j x)-p^{-3}\bbone_{\mathbb{Z}_p^3}(p^{j+1} x)
\right)\ .
\]
By Abel summation, or discrete integration by parts,
this can be rewritten as
\begin{equation}
\Gamma(x)=\bbone_{\mathbb{Z}_p^3}(x)-p^{-3-2(l-1)[\phi]}\bbone_{\mathbb{Z}_p^3}(p^l x)
+\sum_{j=1}^{l-1}p^{-2j[\phi]}(1-p^{-3+2[\phi]})\bbone_{\mathbb{Z}_p^3}(p^j x)\ .
\label{Abelgam}
\end{equation}
Now we also have
\[
\bbone_{\mathbb{Z}_p^3}(p^j x)=\bbone\{|p^jx|\le 1\}=\bbone\{|x|\le p^j\}=\sum_{i\le j}\bbone\{|x|=p^i\}\ .
\]
We insert the last expression into the sum in (\ref{Abelgam}) and get after commuting the sums over $i$ and $j$ that
\[
\Gamma(x)=\bbone_{\mathbb{Z}_p^3}(x)-p^{-3-2(l-1)[\phi]}\bbone_{\mathbb{Z}_p^3}(p^l x)
+\sum_{i\in\mathbb{Z}} U_i \bbone\{|x|=p^i\}
\]
where
\[
U_i=\sum_{j\in\mathbb{Z}}
\bbone\left\{
\begin{array}{c}
1\le j\le l-1\\
i\le j
\end{array}
\right\} p^{-2j[\phi]}\ .
\]
Now note that if $i\ge l$ then $U_i=0$. Also, if $i\le 0$ then
\[
U_i=\frac{p^{-2[\phi]}-p^{-2l[\phi]}}{1-p^{-2[\phi]}}\ .
\]
Finally, if $1\le i\le l-1$ then
\[
U_i=\frac{p^{-2i[\phi]}-p^{-2l[\phi]}}{1-p^{-2[\phi]}}\ .
\]
As a result we have
\[
\Gamma(x)=\bbone\{|x|\le 1\}-p^{-3+2[\phi]} p^{-2l[\phi]}
\bbone\{|x|\le p^l\}
+\frac{1-p^{-3+2[\phi]}}{1-p^{-2[\phi]}}
\sum_{i\le l-1}\bbone\{|x|=p^i\}
\left(p^{-2[\phi]\max(i,1)}-p^{-2l[\phi]}
\right)
\]
from which the result follows by specialization to the different cases mentioned.
\qed

As a result of the previous lemma we have a precise control over the sign of the function $\Gamma$.

\begin{Lemma}\label{L2lem}

\ 
\begin{enumerate}
\item
If $|x|<p^l$ then $\Gamma(x)>0$.
\item
If $|x|=p^l$ then $\Gamma(x)<0$.
\item
If $|x|>p^l$ then $\Gamma(x)=0$.
\end{enumerate}
\end{Lemma}

\noindent{\bf Proof:}
Recall that $\epsilon\in (0,1]$ and therefore $[\phi]=\frac{3-\epsilon}{4}\in \left[
\frac{1}{2},\frac{3}{4}\right)$.
We also have $l\ge 1$ and of course the prime number $p$ is at least $2$.
From Lemma \ref{L1lem} 1), we then readily get that $\Gamma(x)>0$ if $|x|\le 1$. The case $|x|>p^l$ has already been
considered.
For $|x|=p^l$ the formula in Lemma \ref{L1lem} 2) reduces to $\Gamma(x)=-p^{-3+2[\phi]}p^{-2l[\phi]}<0$.
Finally when $|x|=p^i$, $2\le i\le l-1$
then the formula in Lemma \ref{L1lem} 2) shows that $\Gamma(x)$ decreases with $i$ in that range.
We only need look at the case $i=l-1$ where one has
\[
\Gamma(x)=p^{-2(l-1)[\phi]}\left[
1-p^{-3}-p^{-3+2[\phi]}
\right]\ .
\]
Simply using $p^{-3}\le \frac{1}{8}$ and $3-2[\phi]>\frac{3}{2}$, which implies $p^{-3+2[\phi]}<2^{-\frac{3}{2}}$,
we get $1-p^{-3}-p^{-3+2[\phi]}>0$ and thus $\Gamma(x)>0$.
\qed

\begin{Corollary}\label{gamL1cor}
The fluctuation covariance satisfies the $L^1$ bound
\[
||\Gamma||_{L^1}<\frac{1}{\sqrt{2}} L^{3-2[\phi]}\ .
\]
\end{Corollary}

\noindent{\bf Proof:}
Indeed, by $\Gamma=C_0-C_1$ and the definitions of the $C_r$ covariances in \S\ref{formalstatsec}
we have that $\int_{\mathbb{Q}_p^3} {\rm d}^3 x\ \Gamma(x)\ =\widehat{\Gamma}(0)=0$.
In other words the positive part exactly cancels the negative part which is easy to compute since it only involves
$x$'s with $|x|=p^l$.
Therefore
\begin{eqnarray*}
||\Gamma||_{L^1} & = & -2 \int_{\mathbb{Q}_p^3} {\rm d}^3 x\ \Gamma(x)\ \bbone\{|x|=p^l\}\\
 & = & 2(1-p^{-3})p^{-3+2[\phi]} L^{3-2[\phi]}\ .
\end{eqnarray*}
We use $1-p^{-3}<1$ and again $p^{-3+2[\phi]}<2^{-\frac{3}{2}}$ to conclude.
\qed

As for the unit cut-off covariance $C_0$, the following easy property will useful in the sequel.

\begin{Lemma}\label{czerobdlem}
When $\epsilon\in(0,1]$, we have $1<C_0(0)<2$.
\end{Lemma}

\noindent{\bf Proof:}
Recall that
\[
C_0(0)=\frac{1-p^{-3}}{1-p^{-2[\phi]}}=\frac{1-p^{-3}}{1-p^{-\left(\frac{3-\epsilon}{2}\right)}}\ .
\]
Only using $p\ge 1$ and the given range for $\epsilon$ we get
\[
p^{-\frac{3}{2}}\le p^{-\left(\frac{3-\epsilon}{2}\right)}\le p^{-1}\le \frac{1}{2}\ .
\]
Hence
\[
1<\frac{1-p^{-3}}{1-p^{-1}}\le C_0(0)\le \frac{1-p^{-3}}{1-p^{-\frac{3}{2}}}=1+p^{-\frac{3}{2}}< 2\ .
\]
\qed

We will also need some information on the $L^\infty$ and $L^2$ norms of $\Gamma$ which are provided by the following two easy lemmas.

\begin{Lemma}\label{gamL00lem}
We have the simple estimate
\[
||\Gamma||_{L^{\infty}}\le 2\ .
\]
\end{Lemma}

\noindent{\bf Proof:}
If $|x|\le 1$, it follows from Lemmas \ref{L1lem} and \ref{czerobdlem} that $0<\Gamma(x)<2$.
If $|x|>L$, then $\Gamma(x)=0$. If $|x|=L$, then 
\[
|\Gamma(x)|=|-p^{-(3-2[\phi])}L^{-2[\phi]}|\le 1\ .
\]
Finally if $|x|=p^i$ with $1\le i\le l-1$, then by Lemma \ref{L2lem}
\[
|\Gamma(x)|=\Gamma(x)=-p^{-3+2[\phi]}p^{-2l[\phi]}+
\frac{1-p^{-3+2[\phi]}}{1-p^{-2[\phi]}}(p^{-2i[\phi]}-p^{-2l[\phi]})
\]
\[
\le \frac{1-p^{-3+2[\phi]}}{1-p^{-2[\phi]}}(p^{-2i[\phi]}-p^{-2l[\phi]})
\le \frac{1-p^{-3+2[\phi]}}{1-p^{-2[\phi]}}\le \frac{1-p^{-3}}{1-p^{-2[\phi]}}=C_0(0)<2\ .
\]
This shows $|\Gamma(x)|\le 2$ in all cases.
\qed

\begin{Lemma}\label{gamL2lem}
We have
\[
\int_{\mathbb{Q}_p^3} |\Gamma(x)|^2\ {\rm d}^3 x=\frac{(1-p^{-3})(L^\epsilon -1)}{p^\epsilon -1}
\longrightarrow (1-p^{-3})\times l 
\]
when $\epsilon\rightarrow 0$,
with $l$ defined by $L=p^l$ and the limit taken with $L$ fixed.
\end{Lemma}

\noindent{\bf Proof:}
By the Plancherel formula over the $p$-adics
\[
\int_{\mathbb{Q}_p^3} |\Gamma(x)|^2\ {\rm d}^3 x=
\int_{\mathbb{Q}_p^3} |\widehat{\Gamma}(k)|^2\ {\rm d}^3 k\ .
\]
But 
\[
\widehat{\Gamma}(k)=\widehat{C}_0(k)-\widehat{C}_1(k)=\frac{\bbone\{L^{-1}<|k|\le 1\}}{|k|^{3-2[\phi]}}
\]
and therefore
\begin{eqnarray*}
\int_{\mathbb{Q}_p^3} |\Gamma(x)|^2\ {\rm d}^3 x & = & \int_{\mathbb{Q}_p^3}
\ \frac{\bbone\{L^{-1}<|k|\le 1\}}{|k|^{6-4[\phi]}}{\rm d}^3 k \\
 & = & \sum_{j=0}^{l-1}  \int_{\mathbb{Q}_p^3}
\ \frac{\bbone\{|k|=p^{-j}\}}{(p^{-j})^{6-4[\phi]}}{\rm d}^3 k \\
 & = & \sum_{j=0}^{l-1} (1-p^{-3}) p^{-3j} p^{j(6-4[\phi])}\ .
\end{eqnarray*}
The result follows since $3-4[\phi]=\epsilon$ and of course the $\epsilon\rightarrow 0$ limit is trivial.
\qed

\subsection{Gaussian integration bound}

\begin{Lemma}\label{L3lem}
Let $\Delta'$ be a block in $\mathbb{L}$. 
Let the real parameter $\alpha$ satisfy $0\le\alpha\le\frac{\sqrt{2}}{4} L^{-(3-2[\phi])}$.
If $f$ is a real-valued function on $L^{-1}\Delta'$ which is constant on unit cubes and
such that $||f||_{L^\infty}\le\frac{1}{2}  L^{-\frac{1}{2}(3-2[\phi])}$, then
for any finite set $X\subset [L^{-1}\Delta']$ we have the bound
\[
\int {\rm d}\mu_\Gamma(\zeta)\ e^{\int_{L^{-1}\Delta'} f\zeta} \prod_{\Delta\in X} e^{\alpha \zeta_{\Delta}^2}
\le 2^{|X|} e^{\frac{1}{2} (f,\Gamma f)_{L^{-1}\Delta'}}\ .
\]
\end{Lemma}

\noindent{\bf Proof:}

First 
note that one can view the integral we would like to bound, $I$,  as an expectation with respect to the centered Gaussian vector
$(\zeta_\Delta)_{\Delta\in[L^{-1}\Delta']}$ in
$\mathbb{R}^{L^3}$ with covariance
$\mathbf{E}(\zeta_{\Delta_1}\zeta_{\Delta_2})=\Gamma_{\Delta_1,\Delta_2}=\Gamma(x_1-x_2)$
where $x_1$ is any point in $\Delta_1$ and likewise for $x_2$ in $\Delta_2$.
Let $u_1,\ldots,u_{L^3}$ be an orthonormal basis which diagonalizes $\Gamma$ (seen as an $L^3\times L^3$ matrix).
Let $\lambda_1,\ldots,\lambda_{L^3}$ be the corresponding eigenvalues and suppose we arranged the numbering so that
$\lambda_1\ge \lambda_2\ge \cdots$.
Note that the matrix $\Gamma$ is singular and therefore only positive semi-definite, because of the property that
$\int_{L^{-1}\Delta'} \zeta=0$ almost surely.
We therefore introduce $m=\max\{i|\lambda_i>0\}$. We now have that $\zeta$ has the same law
as $\sum_{i=1}^{m} a_i u_i$ where the $a_i$'s are independent centered Gaussian random variables with
variance $\lambda_i$.
Thus
\[
I=\prod_{i=1}^m (2\pi\lambda_i)^{-\frac{1}{2}} \times
\int_{\mathbb{R}^m} {\rm d}a_1\ldots{\rm d}a_m
\ \exp\left[
-\frac{1}{2}\sum_{i=1}^m \frac{a_i^2}{\lambda_i}+
\sum_{\substack{\Delta\in[L^{-1}\Delta'] \\ 1\le i\le m}} f_\Delta a_i u_{i,\Delta}
+\alpha \sum_{\Delta\in X}\left(
\sum_{i=1}^{m} a_i u_{i,\Delta}
\right)^2
\right]\ .
\]
Since $X\subset[L^{-1}\Delta']$
\[
\sum_{\Delta\in X}\left(
\sum_{i=1}^{m} a_i u_{i,\Delta}
\right)^2
\le
\sum_{\Delta\in [L^{-1}\Delta']}\left(
\sum_{i=1}^{m} a_i u_{i,\Delta}
\right)^2
=\sum_{i=1}^m a_i^2
\]
because of the orthonormality of the $u$'s.
Therefore a sufficient condition for the convergence of the integral is that $2\alpha\lambda_i<1$
for all $i$, $1\le i\le m$.
Granting this condition for now, we define $\tilde{f}_i=\sum_{\Delta\in[L^{-1}\Delta']} f_\Delta u_{i,\Delta}$
and use the standard `completing the square' trick by writing
\[
-\frac{1}{2}\sum_{i=1}^m \frac{a_i^2}{\lambda_i} +\sum_{i=1}^m a_i \tilde{f}_i
=-\frac{1}{2}\sum_{i=1}^m \frac{1}{\lambda_i}(a_i-\lambda_i \tilde{f}_i)^2
+\frac{1}{2}\sum_{i=1}^m \lambda_i \tilde{f}_{i}^{2}
\]
and changing variables to $a_i-\lambda_i \tilde{f}_i$. Hence
\[
I=\prod_{i=1}^m (2\pi\lambda_i)^{-\frac{1}{2}} \times
\int_{\mathbb{R}^m} {\rm d}a_1\ldots{\rm d}a_m
\ \exp\left[
-\frac{1}{2}\sum_{i=1}^m \frac{a_i^2}{\lambda_i}
+\frac{1}{2}\sum_{i=1}^m \lambda_i \tilde{f}_{i}^{2}
+\alpha \sum_{\Delta\in X}\left(
\sum_{i=1}^{m} (a_i+\lambda_i\tilde{f}_i) u_{i,\Delta}
\right)^2
\right]\ .
\]
Note that
\begin{eqnarray*}
\sum_{i=1}^m \lambda_i \tilde{f}_i^2 & = & \sum_{i=1}^m \sum_{\Delta_1,\Delta_2\in[L^{-1}\Delta']} \lambda_i
f_{\Delta_1} f_{\Delta_2} u_{i,\Delta_1} u_{i,\Delta_2}\\
 & = &  \sum_{\Delta_1,\Delta_2\in[L^{-1}\Delta']} f_{\Delta_1} f_{\Delta_2} \Gamma_{\Delta_1,\Delta_2}\\
 & = & (f,\Gamma f)_{L^{-1}\Delta'}
\end{eqnarray*}
by construction of the $u$'s.
We also have
\begin{eqnarray*}
\sum_{i=1}^m (a_i+\lambda_i \tilde{f}_i)u_{i,\Delta} & = & \zeta_\Delta+ \sum_{i=1}^m \sum_{\Delta_1\in[L^{-1}\Delta']}
\lambda_i f_{\Delta_1} u_{i,\Delta_1} u_{i,\Delta}\\
 & = & \zeta_\Delta +\sum_{\Delta_1\in[L^{-1}\Delta']} \Gamma_{\Delta,\Delta_1} f_{\Delta_1}\\
 & = & \zeta_\Delta+(\Gamma f)_{\Delta} 
\end{eqnarray*}
where we reverted to the use of the $\zeta_\Delta$ variables of integration which have the same law as the quantities
$\sum_{i=1}^m a_i u_{i\Delta}$,
and where $(\Gamma f)(x)$ denotes $\int_{\mathbb{Q}_p^3} {\rm d}^3 y\ \Gamma(x-y) f(y)$. By the finite range property of $\Gamma$
we have, for $x\in\Delta\in[L^{-1}\Delta']$, $(\Gamma f)(x)=(\Gamma f)_\Delta= \sum_{\Delta_1\in[L^{-1}\Delta']}
\Gamma_{\Delta,\Delta_1} f_{\Delta_1}$.
As a result of the previous calculations
\[
I= e^{\frac{1}{2}(f,\Gamma f)_{L^{-1}\Delta'}}\times
\int {\rm d}\mu_\Gamma(\zeta)\ e^{\alpha \sum_{\Delta\in X}((\Gamma f)_\Delta+\zeta_\Delta)^2}\ .
\]
We now expand the square in the last exponential and we also introduce the covariance matrix $\Gamma_X$ for the marginal
random vector $\zeta|_{X}=(\zeta_\Delta)_{\Delta\in X}$ in order to write
\[
I=e^{\frac{1}{2}(f,\Gamma f)_{L^{-1}\Delta'}+\alpha(\Gamma f,\Gamma f)_X}\times
\int {\rm d}\mu_{\Gamma_X}(\zeta|_{X})\ e^{\alpha\langle \zeta|_X,\zeta|_X\rangle +2\alpha
\langle \Gamma f|_X,\zeta|_X\rangle}
\]
where the inner products are the ones of $l^2(X)$, namely $\langle w,w'\rangle=\sum_{\Delta\in X}w_\Delta w'_\Delta$ for vectors 
in $l^2(X)$ which are indexed by boxes in the finite set $X$.

Let $(v_i)_{1\le i\le |X|}$ be an orthonormal basis diagonalizing the symmetric positive semi-definite
matrix $\Gamma|_X$, with eigenvalues $\mu_i$ arranged so that $\mu_1\ge \mu_2\ge\cdots$ and let $n=\max\{i|\mu_i>0\}$.
As before, we have that the random vector $\zeta|_X$ has the same law as $\sum_{i=1}^n b_i v_i$
where the $b_i$ are independent centered Gaussian random variables  with variance $\mu_i$.
Following this change of variables of integration
$\langle \zeta|_X,\zeta|_X\rangle$ becomes $\sum_{i=1}^{n}b_i^2$
whereas $\langle \Gamma f|_X,\zeta|_X\rangle$ becomes $\sum_{i=1}^{n} g_i b_i$
with $g_i=\sum_{\Delta\in X} (\Gamma f)_\Delta v_{i,\Delta}$.
Hence
\begin{eqnarray*}
\int {\rm d}\mu_{\Gamma_X}(\zeta|_{X})\ e^{\alpha\langle \zeta|_X,\zeta|_X\rangle +2\alpha
\langle \Gamma f|_X,\zeta|_X\rangle} & = & \prod_{i=1}^n\left[
\frac{1}{\sqrt{2\pi\mu_i}}\int_{\mathbb{R}} {\rm d}b_i\ 
e^{-\frac{b_i^2}{2\mu_i}+\alpha b_i^2+2\alpha g_i b_i}
\right] \\
 & = & \prod_{i=1}^n\left[
\frac{1}{\sqrt{2\pi\mu_i}}
\times\sqrt{2\pi\left(\frac{1}{\mu_i}-2\alpha\right)^{-1}}
\times e^{\frac{1}{2}\left(\frac{1}{\mu_i}-2\alpha\right)^{-1} (2\alpha g_i)^2}
\right] \\
 & = & \prod_{i=1}^n\left[\frac{1}{\sqrt{1-2\alpha\mu_i}}
e^{2\alpha^2\frac{\mu_i}{1-2\alpha\mu_i} g_i^2}
\right]
\end{eqnarray*}
provided $2\alpha\mu_i<1$ for all $i$, $1\le i\le n$.

Now $\mu_i\le ||\Gamma_X||$ where the latter quantity is the operator norm of $\Gamma_X$
induced by the norm on $l^2(X)$ coming from the inner product $\langle\cdot,\cdot\rangle$.
For $v$ a real vector in $l^2(X)$, we have $||\Gamma_X v||^2=\sum_{\Delta\in X}(\Gamma_X v)_{\Delta}^{2}
=\sum_{\Delta\in X}(\Gamma w)_{\Delta}^{2}$
where $w\in l^2([L^{-1}\Delta'])$ is the extension of $v$ by zero outside $X$.
Thus $||\Gamma_X v||^2\le \sum_{\Delta\in [L^{-1}\Delta']}(\Gamma_X w)_{\Delta}^{2}=||\Gamma w||^2\le ||\Gamma||^2 ||w||^2=
||\Gamma||^2 ||v||^2$. As a result $||\Gamma_X||\le ||\Gamma||$ where the latter is the operator norm of the matrix
$\Gamma$ coming from the inner product norm of $l^2([L^{-1}\Delta'])$.
However we have the bound $||\Gamma||\le ||\Gamma||_{L^1}=\int_{\mathbb{Q}_p^3} |\Gamma(x)| {\rm d}^3 x$.
Indeed, given $w\in l^2([L^{-1}\Delta'])$
which we can identify with a function $w(x)$ on $\mathbb{Q}_p^3$ with support in $L^{-1}\Delta'$
and which is constant on unit blocks, we have
\begin{eqnarray*}
||\Gamma w||^2 & = & \int_{\mathbb{Q}_p^3} [(\Gamma w)(x)]^2\ {\rm d}^3x \\
 & = & 
\int_{\mathbb{Q}_p^{3\times 3}}
\Gamma(x-y)\Gamma(x-z)w(y)w(z)\ {\rm d}^3 x\ {\rm d}^3 y\ {\rm d}^3 z \\
 & \le & \int_{\mathbb{Q}_p^{3\times 3}}
|\Gamma(x-y)|\ |\Gamma(x-z)|\ |w(y)|\ |w(z)|\ {\rm d}^3 x\ {\rm d}^3 y\ {\rm d}^3 z \\
 & \le & \int_{\mathbb{Q}_p^{3\times 3}}
|\Gamma(x-y)|\ |\Gamma(x-z)|\ \left(\frac{1}{2}|w(y)|^2+\frac{1}{2}|w(z)|^2\right)\ {\rm d}^3 x\ {\rm d}^3 y\ {\rm d}^3 z \\
 & = & 2\times\frac{1}{2}\times ||\Gamma||_{L^1}^2 ||w||_{L^2}^2\ .
\end{eqnarray*}
Therefore from Corollary \ref{gamL1cor} we get $||\Gamma||\le ||\Gamma||_{L^1}<\frac{1}{\sqrt{2}}L^{3-2[\phi]}$.
Since the $\lambda_i$ are bounded by $||\Gamma||$ (the case where $X=[L^{-1}\Delta']$),
the hypothesis $\alpha \le \frac{\sqrt{2}}{4}L^{-3+2[\phi]}$
implies that the previous convergence requirement $2\alpha\lambda_i<1$ is satisfied
and also that not only $2\alpha\mu_i<1$ holds but so does the stronger inequality $2\alpha\mu_i\le \frac{1}{2}$.
From the latter we have $\frac{\mu_i}{1-2\alpha\mu_i}\le 2\mu_i$
and thus
\begin{eqnarray*}
\int {\rm d}\mu_{\Gamma_X}(\zeta|_{X})\ e^{\alpha\langle \zeta|_X,\zeta|_X\rangle +2\alpha
\langle \Gamma f|_X,\zeta|_X\rangle} & \le & \prod_{i=1}^{n} \left(\sqrt{2}
e^{4\alpha^2\mu_i g_i^2} \right) \\
 & \le & 2^{\frac{|X|}{2}}\exp\left(\frac{\sqrt{2}}{4}L^{-(3-2[\phi])}\sum_{i=1}^{n} g_i^2\right)
\end{eqnarray*}
where we used $n\le |X|$, $\alpha\le \frac{\sqrt{2}}{4}L^{-(3-2[\phi])}$ and $\mu_i<\frac{1}{\sqrt{2}}L^{3-2[\phi]}$.
Besides, $g_i=\sum_{\Delta\in X} (\Gamma f)_\Delta v_{i,\Delta}=\langle v_i,(\Gamma f)|_X\rangle$
and therefore
\[
\sum_{i=1}^{n} g_i^2\le \sum_{i=1}^{|X|}\langle v_i,(\Gamma f)|_X\rangle^2=\langle (\Gamma f)|_X,(\Gamma f)|_X\rangle
=(\Gamma f,\Gamma f)_X\ .
\]
But $(\Gamma f,\Gamma f)_X=\sum_{\Delta\in X} (\Gamma f)_{\Delta}^2$
and clearly $|(\Gamma f)_{\Delta}|\le ||\Gamma||_{L^1} ||f||_{L^\infty}$
so $(\Gamma f,\Gamma f)_X\le |X|\ ||\Gamma||_{L^1}^2 ||f||_{L^\infty}^2$.

Putting all the previous bounds together we see that the desired inequality holds provided
\[
\exp\left[\frac{\sqrt{2}}{4}L^{3-2[\phi]} ||f||_{L^\infty}^2\right]\le \sqrt{2}
\]
which is true since, by hypothesis, $||f||_{L^\infty}\le\frac{1}{2}L^{-\frac{1}{2}(3-2[\phi])}$ and 
$\frac{4}{\sqrt{2}}\times\frac{1}{2}\log 2\simeq 0.980\ldots>\frac{1}{4}$.
\qed

\subsection{Two easy lemmas}

The following are simple bounds which will however be used many times in order to bound individual field factors
using the exponential of a quartic or a quadratic expression. The quartic case will typically apply to background fields $\phi$
whereas quadratic bounds will typically apply to fluctuation fields $\zeta$.

Note that the possibly complex $\phi^4$ couplings $\beta_4$ will sit in an open ball of the form $|\beta_4-\bar{g}|<\frac{1}{2}\bar{g}$
with $\bar{g}>0$.
By elementary trigonometry it easily follows that $\frac{\Re \beta_4}{|\beta_4|}\ge\frac{\sqrt{3}}{2}$.
We of course also have $\frac{1}{2}<\frac{\Re \beta_4}{\bar{g}}<\frac{3}{2}$.

\begin{Lemma}\label{L4lem}
$\forall j\in\mathbb{N}$, $\forall \bar{g}>0$, $\forall \gamma>0$, $\forall \beta_4\in\mathbb{C}$
such that $|\beta_4-\bar{g}|<\frac{1}{2}\bar{g}$, $\forall \phi\in\mathbb{R}$ we have
\[
|\phi|^j\le \left(\frac{j}{4e}\right)^{\frac{j}{4}} (\gamma \Re \beta_4)^{-\frac{j}{4}} e^{\gamma (\Re \beta_4)\phi^4} 
\le \left(\frac{j}{2e}\right)^{\frac{j}{4}} (\gamma \bar{g})^{-\frac{j}{4}} e^{\gamma (\Re \beta_4)\phi^4} 
\]
with the convention $j^j=1$ if $j=0$.
\end{Lemma}

\noindent{\bf Proof:}
The function $u^{\frac{j}{4}} e^{-u}$ for $u\ge 0$ is maximized when $u=\frac{j}{4}$. Simply apply this to $u=\gamma(\Re \beta_4)\phi^4$
and use $\frac{1}{2}<\frac{\Re \beta_4}{\bar{g}}$ for the second inequality.
\qed

\begin{Lemma}\label{L5lem}
$\forall j\in\mathbb{N}$, $\forall \kappa>0$, $\forall \zeta\in\mathbb{R}$ we have
\[
|\zeta|^j\le \left(\frac{j}{2e}\right)^{\frac{j}{2}} \kappa^{-\frac{j}{2}} e^{\kappa \zeta^2} 
\]
again with the convention $j^j=1$ if $j=0$.
\end{Lemma}

The proof is similar.

\subsection{The key stability bound}

The following lemma is essential to our estimates since it provides bounds on the seminorms introduced in \S\ref{funcspsec}
for functions given by the exponential of a degree four polynomial in the real-valued field, with complex coefficients.

\begin{Lemma}\label{L6lem}
Let $U(\phi)=a_4\phi^4+a_3\phi^3+a_2\phi^2+a_1\phi+a_0$ where the possibly complex coefficients $a_0,\ldots,a_4$
satisfy $|a_4|>0$, $\Re a_4\ge\frac{\sqrt{3}}{2} |a_4|$, $|a_k|\le\frac{1}{3}\log\left(\frac{1+\sqrt{2}}{2}\right)|a_k|^{\frac{k}{4}}$
for $k=1,2,3$, and $|a_0|\le \log 2$.
Then
\begin{enumerate}
\item
the condition
\[
0\le\theta\le\frac{\sqrt{2}-1}{4} e^{-918785}\times |a_4|^{-\frac{1}{4}}
\]
implies
\[
||e^{-U(\phi)}||_{\partial\phi,\phi,\theta}\le 2 e^{-\frac{1}{2}(\Re a_4)\phi^4}
\]
for all $\phi\in\mathbb{R}$;
\item
the condition
\[
0\le\theta\le\frac{(\sqrt{2}-1)^2}{e}\times |a_4|^{-\frac{1}{4}}
\]
implies
\[
|e^{-U(\phi)}|_{\partial\phi,\theta}\le 2 \ .
\]
\end{enumerate}
\end{Lemma}

\noindent{\bf Proof:}
It follows from the definition of our seminorms that
\[
||e^{-U(\phi)}||_{\partial\phi,\phi,\theta}=
e^{-\Re U(\phi)}+\sum_{n=1}^9 \frac{\theta^n}{n!}\left|D^n e^{-U(\phi)}\right|
\]
where $D$ denotes the differentiation operator $\frac{{\rm d}}{{\rm d}\phi}$.
An easy induction provides the following explicit formula of Faa di Bruno type for the derivatives of functions of the form $e^{f(\phi)}$:
\begin{equation}
D^n e^{f(\phi)}=\sum_{k\ge 0}\frac{1}{k!}
\sum_{\substack{m_1,\ldots,m_k\ge 1 \\ \Sigma m_i=n}}
\frac{n!}{m_1!\cdots m_k!}\left(\prod_{i=1}^{k} D^{m_i} f(\phi)\right) e^{f(\phi)}\ .
\label{FaadiBeq}
\end{equation}
This will be used in order to bound the quantities $\left|D^n e^{-U(\phi)}\right|$.
First, let us introduce the notation $\alpha=\frac{\sqrt{3}}{2}$ and $r=\frac{1}{3}\log\left(\frac{1+\sqrt{2}}{2}\right)$.
We have
\begin{eqnarray*}
-\Re U(\phi) & = & -\sum_{k=0}^{4} (\Re a_k)\phi^k \\
 & \le & -\frac{1}{2}(\Re a_4)\phi^4-\frac{\alpha}{2}|a_4|\phi^4+\left(\sum_{k=1}^{3}|a_k|\ |\phi|^k\right)+|a_0|
\end{eqnarray*}
from the hypothesis $\Re a_4\ge\alpha |a_4|$.
Using the assumption $|a_k|\le r|a_4|^{\frac{k}{4}}$
we then obtain
\[
-\Re U(\phi)\le -\frac{1}{2}(\Re a_4)\phi^4 +\Omega_1(|a_4|^{\frac{1}{4}}|\phi|)+|a_0|
\]
where $\Omega_1(x)=-\frac{\alpha}{2}x^4+r(x^3+x^2+x)$.
We first write a convenient upper bound on $\sup_{x\ge 0} \Omega_1(x)$.
For $0\le x\le 1$, we simply use $\Omega_1(x)\le r(x^3+x^2+x)\le 3r$.
For $x\ge 1$, we write $\Omega_1(x)\le -\frac{\alpha}{2}x^4+3r x^3$ and maximize the right-hand side over $[0,\infty)$. The maximum occurs
at $x=\frac{9r}{2\alpha}$ and is equal to $\frac{3r}{4}\left(\frac{9r}{2\alpha}\right)^3<3r$. The last inequality used the fact $9r<2\alpha$
which can be checked from the chosen numerical values of $r\simeq 0.0627\ldots$ and $\alpha\simeq 0.866\ldots$
As a result
\[
e^{-\Re U(\phi)}\le e^{-\frac{1}{2}(\Re a_4)\phi^4+3r+|a_0|} \ .
\]
We now use the formula (\ref{FaadiBeq}) and write, for $1\le n\le 9$,
\[
D^n e^{-U(\phi)}=e^{-U(\phi)}\times\sum_{k=1}^n\frac{1}{k!}
\sum_{\substack{1\le m_1,\ldots,m_k\le 4 \\ \Sigma m_i=n}}
\frac{n!}{m_1!\cdots m_k!}\times\prod_{i=1}^{k}\left(- D^{m_i} U(\phi)\right)\ .
\]
Using the condition $\Sigma m_i=n$ for handling the $\theta$ exponents we get the bound
\begin{equation}
\frac{\theta^n}{n!}\left|D^n e^{-U(\phi)}\right|\le
e^{-\Re U(\phi)}\times
\sum_{k=1}^n\frac{1}{k!}
\sum_{\substack{1\le m_1,\ldots,m_k\le 4 \\ \Sigma m_i=n}}
\ \prod_{i=1}^{k}\left[\frac{\theta^{m_i}\left|D^{m_i} U(\phi)\right|}{m_i!}\right]\ .
\label{dnUbdeq}
\end{equation}
We now assume $\theta\le {\gamma_1}|a_4|^{-\frac{1}{4}}$
for some suitable ${\gamma_1}\ge 0$ to be specified later. We insert this inequality in (\ref{dnUbdeq}) and pull out ${\gamma_1}^{\Sigma m_i}={\gamma_1}^n$
before throwing away the constraint $\Sigma m_i=n$ which results in
\begin{eqnarray*}
\frac{\theta^n}{n!}\left|D^n e^{-U(\phi)}\right| & \le & e^{-\Re U(\phi)} {\gamma_1}^n\times
\sum_{k=1}^n\frac{1}{k!}\left(
\sum_{m=1}^4 \frac{|a_4|^{-\frac{m}{4}}}{m!}\left|D^{m} U(\phi)\right|
\right)^k \\
 & \le & {\gamma_1}^n\ \exp\left[
 -\Re U(\phi)+\sum_{m=1}^4 \frac{|a_4|^{-\frac{m}{4}}}{m!}\left|D^{m} U(\phi)\right|
\right]\ .
\end{eqnarray*}
The individual quantities in the last exponential are bounded in terms of $x=|a_4|^{\frac{1}{4}}|\phi|$ as follows:
\begin{eqnarray*}
|a_4|^{-\frac{1}{4}} \left|D U(\phi)\right| & = & |a_4|^{-\frac{1}{4}}\times
|4 a_4\phi^3+3 a_3\phi^2+2 a_2\phi+a_1| \\
 & \le & 4x^3+3rx^2+2rx+r\ ,
\end{eqnarray*}
\begin{eqnarray*}
\frac{|a_4|^{-\frac{2}{4}}}{2} \left|D^2 U(\phi)\right| & = & |a_4|^{-\frac{1}{2}}\times
|6 a_4\phi^2+3 a_3\phi+a_2| \\
 & \le & 6x^2+3rx+r\ ,
\end{eqnarray*}
\begin{eqnarray*}
\frac{|a_4|^{-\frac{3}{4}}}{3!} \left|D^3 U(\phi)\right| & = & |a_4|^{-\frac{3}{4}}\times
|4 a_4\phi+a_3| \\
 & \le & 4x+r\ ,
\end{eqnarray*}
whereas
\[
\frac{|a_4|^{-\frac{4}{4}}}{4!} \left|D^4 U(\phi)\right|=1 \ .
\]
Therefore
\[
\sum_{m=1}^4 \frac{|a_4|^{-\frac{m}{4}}}{m!}\left|D^{m} U(\phi)\right|
\le 4x^3+(3r+6)x^2+(5r+4)x+(3r+1)
\]
and
\[
-\Re U(\phi)+\sum_{m=1}^4 \frac{|a_4|^{-\frac{m}{4}}}{m!}\left|D^{m} U(\phi)\right|
\le -\frac{1}{2}(\Re a_4)\phi^4 +\Omega_2(|a_4|^{\frac{1}{4}}|\phi|)+|a_0| 
\]
where
\[
\Omega_2(x)=-\frac{\alpha}{2}x^4
+(r+4)x^3+(4r+6)x^2+(6r+4)x+(3r+1)\ .
\]
We again find a convenient bound on $\sup_{x\ge 0} \Omega_2(x)$.
Simply using $r<1$ and dropping the $x^4$ term
we have, for $0\le x\le 1$, the bound $\Omega_2(x)\le 5+10+10+4=29$.
For $x\ge 1$, we have the crude bound $\Omega_2(x)\le -\frac{\alpha}{2}x^4+29 x^3$ and we proceed with maximizing the right-hand side
over $x\in [0,\infty)$. When $\alpha=\frac{\sqrt{3}}{2}$ the
maximum occurs at $x=\frac{87}{\sqrt{3}}$ and is equal to
$\frac{29^4\times 3\sqrt{3}}{4}\simeq 918784.97\ldots<918785$. We denote the latter numerical constant by $M$.
The previous considerations now give
\begin{eqnarray*}
||e^{-U(\phi)}||_{\partial\phi,\phi,\theta} & \le &
e^{-\frac{1}{2}(\Re a_4)\phi^4+3r+|a_0|}\\
 & & +\sum_{n=1}^9 {\gamma_1}^n
\ \exp\left[-\frac{1}{2}(\Re a_4)\phi^4+M+|a_0|\right]\\
 & \le & e^{-\frac{1}{2}(\Re a_4)\phi^4}\times
 e^{|a_0|}\times
\left[
e^{3r}+e^M\times\frac{{\gamma_1}}{1-{\gamma_1}} 
\right]
\end{eqnarray*}
provided ${\gamma_1}<1$. If one requires the stronger condition  ${\gamma_1}\le \frac{1}{2}$
then $e^{3r}+e^M\times\frac{{\gamma_1}}{1-{\gamma_1}}\le e^{3r}+2 e^M {\gamma_1}$. From our choice for $r$ we have $e^{3r}=\frac{1+\sqrt{2}}{2}$.
If we now set ${\gamma_1}=\frac{\sqrt{2}-1}{4} e^{-M}$ which clearly is less than $\frac{1}{2}$
then $e^{3r}+2 e^M {\gamma_1}=\sqrt{2}$. On the other hand, by assumption on $a_0$ we have $e^{|a_0|}\le \sqrt{2}$.
The statement in 1) is therefore proved.

For the statement in 2) concerning the bound on $|e^{-U(\phi)}|_{\partial\phi,\theta}=||e^{-U(\phi)}||_{\partial\phi,0,\theta}$,
with derivatives taken at zero, we follow the same steps. However, the situation
simplifies considerably. Indeed,
\[
|e^{-U(\phi)}|_{\partial\phi,\theta}=e^{-\Re U(0)}+
\sum_{n=1}^9 \frac{\theta^n}{n!}\left|
\left. D^n e^{-U(\phi)}\right|_{\phi=0}
\right|
\]
can be bounded as we did before, under the new hypothesis $\theta\le {\gamma_2} |a_4|^{-\frac{1}{4}}$ for suitable ${\gamma_2}\in[0,1)$,
by the estimate
\[
|e^{-U(\phi)}|_{\partial\phi,\theta}
\le e^{-\Re a_0}+\frac{{\gamma_2}}{1-{\gamma_2}}\times
\exp\left[
-\Re U(0)+\sum_{m=1}^4 \frac{|a_4|^{-\frac{m}{4}}}{m!} \left|D^m U(0)\right|
\right]\ .
\]
Now
\[
\sum_{m=1}^4 \frac{|a_4|^{-\frac{m}{4}}}{m!} \left|D^m U(0)\right|=
\sum_{m=1}^4 |a_4|^{-\frac{m}{4}} |a_m|\le 3r+1\ .
\]
If one imposes the condition ${\gamma_2}\le \frac{1}{2}$, then
\[
|e^{-U(\phi)}|_{\partial\phi,\theta}\le e^{|a_0|}\times
\left[1+2{\gamma_2} e^{3r+1}\right]\ .
\]
Because of the chosen value of $r$, one  will have $1+2{\gamma_2} e^{3r+1}=\sqrt{2}$
if one now sets ${\gamma_2}=\frac{(\sqrt{2}-1)^2}{e}\simeq 0.0631\ldots$ which is less than $\frac{1}{2}$.
The statement in 2) then follows easily.
\qed

%% file: Mainest.tex
\section{The main estimates on a single extended RG step}

\subsection{Statement of the theorem}

Recall that $\epsilon\in (0,1]$, $L=p^l$ with $p\ge 2$ a prime number and where $l\ge 1$ is an integer.
The symbol $[\phi]$ denotes the quantity $\frac{3-\epsilon}{4}$.
We now introduce the numerical constants
\[
c_1=2^{-\frac{9}{4}}(\sqrt{2}-1) e^{-918785}\qquad{\rm and}\qquad
c_2=2^{\frac{3}{4}}
\]
which are used to calibrate the parameters
\[
h=c_1 \bar{g}^{-\frac{1}{4}}\qquad{\rm and}\qquad
h_\ast=c_2 L^{\frac{3+\epsilon}{4}}
\]
for the seminorms we use. With these choices the norm
\[
|||R|||_{\bar{g}}=\max\left\{
|R(\phi)|_{\partial\phi,h_\ast},
\bar{g}^2 \sup_{\phi\in\mathbb{R}} 
||R(\phi)||_{\partial\phi,\phi,h}
\right\}
\]
is now unambiguously defined in terms of the calibrator $\bar{g}$.
In the present article we will need to take $\bar{g}$
of order $\epsilon$, however the next theorem will be stated in greater generality as far as
allowed values for this calibrator. Indeed, we intend to reuse the rather expensive theorem that follows
for the construction of another quantum field theory on $\mathbb{Q}_p^3$
which corresponds to an RG trajectory joining the Gaussian and the infrared fixed points as in \cite{Abdesselam}.

\begin{Theorem}\label{mainestthm}

$\exists B_{R\mathcal{L}}\ge 0$, $\forall l\ge 1$, $\exists B_0,\ldots,B_4,B_{R\xi}\ge 0$,

$\forall \eta\in\left[0,\frac{1}{4}\right)$, $\forall \eta_R\in\left[3\eta,\frac{3}{16}+\frac{9}{4}\eta\right]$,

$\forall A_{\bar{g}}>0$, $\exists \epsilon_0>0$, $\forall \epsilon\in (0,\epsilon_0]$,
$\forall \bar{g}\in \left(0,A_{\bar{g}}\epsilon\right]$ and $\Delta'\in\mathbb{L}$, then on the domain
\[
\forall \Delta\in [L^{-1}\Delta'],\ 
\left\{
\begin{array}{ccccc}
|\beta_{4,\Delta}-\bar{g}| & < & \frac{1}{2}\bar{g} & &  \\
|\beta_{k,\Delta}| & < & \bar{g}^{1-\eta} & {\rm for} & k=1,2,3 \\
|W_{k,\Delta}| & < & \bar{g}^{2-2\eta} & {\rm for} & k=5,6 \\
|f_\Delta| & < & L^{-(3-[\phi])} & & \\
|||R_\Delta|||_{\bar{g}} & < & \bar{g}^{\frac{11}{4}-\eta_R}
\end{array}
\right.
\]
the maps $\xi_{0,\Delta'},\ldots,\xi_{4,\Delta'}$, $\mathcal{L}_{\Delta'}$ and $\xi_{R,\Delta'}$
are well-defined, analytic, send real data to real data and satisfy the bounds
\[
|\xi_{k,\Delta'}(\vec{V})|\le B_k \max_{\Delta\in[L^{-1}\Delta']} |||R_\Delta|||_{\bar{g}}\qquad
{\rm for}\ k=0,\ldots,4\ ,
\]
\[
|||\mathcal{L}_{\Delta'}^{\vec{\beta},f}(R)|||_{\bar{g}}\le
B_{R\mathcal{L}} L^{3-5[\phi]}
\max_{\Delta\in[L^{-1}\Delta']} |||R_\Delta|||_{\bar{g}}\ ,
\]
and
\[
|||\xi_{R,\Delta'}(\vec{V})|||_{\bar{g}}\le B_{R\xi} \bar{g}^{\frac{11}{4}-3\eta}\ .
\]
\end{Theorem}

\begin{Remark}
In fact, except for the map $\xi_{R,\Delta'}$, the conclusions
of the theorem are valid without restriction on the size of $R$. This is because
$\xi_{0,\Delta'},\ldots,\xi_{4,\Delta'}$, $\mathcal{L}_{\Delta'}$ are linear in $R$.
\end{Remark}

\subsection{The standard hypotheses}\label{stdhypsec}

In this section we collect a set of conditions which we call the standard hypotheses
and which will be assumed throughout \S\ref{longlemsec} and therefore will not be repeated in the local statements of the lemmas
in that section. Finally in \S\ref{concthmsec} we will show that it is possible to satisfy all these conditions and
thus wrap up the proof of Theorem \ref{mainestthm}.
The list of conditions repeats some statements made before and it is also redundant.
We meant this list to collect in one place the various requirements for the validity of the lemmas in \S\ref{longlemsec}.
We assume:
\[
h=c_1 \bar{g}^{-\frac{1}{4}}\qquad{\rm with}\qquad c_1=2^{-\frac{9}{4}}(\sqrt{2}-1) e^{-918785}\ ,
\]
\begin{equation}
h_\ast=c_2 L^{\frac{3+\epsilon}{4}}\qquad{\rm with}\qquad c_2=2^{\frac{3}{4}}\ ,
\label{hstarstdeq}
\end{equation}
\begin{equation}
0<\bar{g}\le 1\ ,\ \eta\ge 0\ ,\ \eta_R\ge 0\ ,\ \eta<\frac{1}{4}\ ,
\label{etastdeq}
\end{equation}
\begin{equation}
|\beta_{4,\Delta}-\bar{g}|<\frac{1}{2}\bar{g}\ ,\ |\beta_{k,\Delta}|<\bar{g}^{1-\eta}\ {\rm for}\ k=1,2,3\ ,
\label{betastdeq}
\end{equation}
\begin{equation}
|f_\Delta|<L^{-(3-[\phi])}\ ,\ |W_{k,\Delta}|<\bar{g}^{2-2\eta}\ {\rm for}\ k=5,6\ ,
\label{fwstdeq}
\end{equation}
\begin{equation}
|||R_\Delta|||_{\bar{g}}<\bar{g}^{\frac{11}{4}-\eta_R}\ {\rm where}\ 
|||R_\Delta|||_{\bar{g}}=\max\left\{
|R_\Delta(\phi)|_{\partial\phi,h_\ast},
\bar{g}^2 \sup_{\phi\in\mathbb{R}} 
||R_\Delta(\phi)||_{\partial\phi,\phi,h}
\right\}\ ,
\label{normstdeq}
\end{equation}
\begin{equation}
\left.
\begin{array}{c}
2^{\frac{3}{4}}\times L^{\frac{3}{4}}\times\bar{g}^{\frac{1}{4}-\eta} \\
19\sqrt{2}\times L^{\frac{3}{2}}\times\bar{g}^{\frac{1}{2}-\eta} \\
2^{\frac{1}{4}}\times 7\times L^{\frac{9}{4}}\times\bar{g}^{\frac{3}{4}-\eta}
\end{array}
\right\}\le\frac{1}{3}\log\left(\frac{1+\sqrt{2}}{2}\right)\ ,
\label{logstdeq}
\end{equation}
\begin{equation}
L^{\epsilon}\le\frac{4}{3}\ ,\ 20\times L^2\times \bar{g}^{1-\eta}\le \log 2\ ,\ 
\bar{g}^{\frac{1}{4}}\le \frac{(\sqrt{2}-1)^2}{2e} \times L^{-1}\ ,
\label{g14stdeq}
\end{equation}
\begin{equation}
\bar{g}\le c_1^4 c_2^{-4} L^{-1}\ ,\  \bar{g}\le c_1^{36} c_2^{-36} L^{-36}\ ,
\label{c1c2stdeq}
\end{equation}
\begin{equation}
\eta_R\ge 3\eta\ ,\ \eta_R\le 1+3\eta\ ,
\label{etarstdeq}
\end{equation}
\begin{equation}
\exp\left(2\bar{g}^{\frac{1}{4}}+18\bar{g}\right)\le 2\ ,\ 8L^{4}\bar{g}^{1-\eta}\le 1\ ,
\label{expgstdeq}
\end{equation}
\begin{equation}
\eta_R\le \frac{3}{16}+\frac{9}{4}\eta\ {\rm and}\ \eta_R<\frac{3}{4}\ .
\label{etar2stdeq}
\end{equation}
We also assume the additional conditions
\begin{equation}
\left.
\begin{array}{c}
2\mathcal{O}_{32}\times L^{15}\times \bar{g}^{\frac{1}{4}-\frac{\eta_R}{3}} \\
2\mathcal{O}_{35}\times L^{15}\times \bar{g}^{\frac{11}{12}-\frac{\eta_R}{3}} \\
3\mathcal{O}_{27}\times L^{5}\times \bar{g}^{\frac{1}{4}-\eta} \\
\mathcal{O}_{28}\times L^{9}\times \bar{g}^{\frac{11}{12}-\frac{\eta_R}{3}}
\end{array}
\right\}\le 1
\label{extrastdeq}
\end{equation}
which involve purely numerical constants $\mathcal{O}$ to be specified in \S\ref{longlemsec}.

\subsection{A long series of lemmas}\label{longlemsec}

Assuming the conditions stated in \S\ref{stdhypsec} we now embark on the following series of lemmas which
will lead to the proof of Theorem \ref{mainestthm}. The estimates will involve a collection of numerical constants
which are given explicitly and are numbered as $\mathcal{O}_1,\mathcal{O}_2$, etc.
Since we are not aiming for optimal estimates, our motivation for keeping such constants explicit
is to serve as an Ariadne thread for the reader in her/his journey through the following estimates
which involve many interdependent parameters.
The notations will continue those of \S\ref{algdefsec}.

\begin{Lemma}\label{L7lem}
For all $t\in (0,1]$ and all unit cube $\Delta$ we have
\[
\forall \phi\in\mathbb{R}\ ,\ \ ||e^{-tV_\Delta(\phi)}||_{\partial\phi,\phi,h}\le 2^{-\frac{t}{2}(\Re\beta_{4,\Delta})\phi^4}
\]
as well as
\[
|e^{-tV_\Delta(\phi)}|_{\partial\phi,h_\ast}\le 2\ .
\]
\end{Lemma}

\noindent{\bf Proof:}
From the definition and by undoing the Wick ordering, we have
\begin{eqnarray*}
V_\Delta(\phi) & = & \sum_{k=1}^4 \beta_{k,\Delta} :\phi^k:_{C_0} \\
 & = & \beta_{4,\Delta}\left(\phi^4-6C_0(0)\phi^2+3C_0(0)^2\right) \\
 & & + \beta_{3,\Delta}\left(\phi^3-3C_0(0)\phi\right) \\
 & & + \beta_{2,\Delta}\left(\phi^2-C_0(0)\right) \\
 & & + \beta_{1,\Delta}\phi
\end{eqnarray*}
and therefore $tV_\Delta(\phi)=\sum_{k=0}^{4} a_k\phi^k$
with
\begin{eqnarray*}
a_4 & = & t\beta_{4,\Delta} \ ,\\
a_3 & = & t\beta_{3,\Delta} \ ,\\
a_2 & = & t\left(\beta_{2,\Delta}-6C_0(0)\beta_{4,\Delta}\right) \ ,\\
a_1 & = & t\left(\beta_{1,\Delta}-3C_0(0)\beta_{3,\Delta}\right) \ ,\\
a_0 & = & t\left(-C_0(0)\beta_{2,\Delta}+3 C_0(0)^2\beta_{4,\Delta}\right) \ . \\
\end{eqnarray*}

We now simply check that the requirements in Lemma \ref{L6lem} are satisfied.
From the standard hypothesis (\ref{betastdeq}) we have $|a_4|=t|\beta_{4,\Delta}|>0$ as well as
$(\Re a_4)\times|a_4|^{-1}=(\Re \beta_{4,\Delta})\times |\beta_{4,\Delta}|\ge \frac{\sqrt{3}}{2}$.
We have, since $t\in (0,1]$,
\begin{eqnarray*}
|a_3|\times|a_4|^{-\frac{3}{4}} & = & t^{\frac{1}{4}} |\beta_{3,\Delta}|\times |\beta_{4,\Delta}|^{-\frac{3}{4}}\\
 & \le & \bar{g}^{1-\eta} \times \left(\frac{1}{2}\bar{g}\right)^{-\frac{3}{4}} \\
 & \le & \frac{1}{3}\log\left(\frac{1+\sqrt{2}}{2}\right)
\end{eqnarray*}
by  (\ref{logstdeq}).
Likewise, using Lemma \ref{czerobdlem}, we have
\begin{eqnarray*}
|a_2|\times|a_4|^{-\frac{2}{4}} & \le & t^{\frac{1}{2}}
\left[ |\beta_{2,\Delta}|+12 |\beta_{4,\Delta}|\right]\times |\beta_{4,\Delta}|^{-\frac{1}{2}}\\
 & \le & \left(\bar{g}^{1-\eta} +12\times\frac{3}{2}\bar{g}\right) \times \left(\frac{1}{2}\bar{g}\right)^{-\frac{1}{2}} \\
 & \le & 19\sqrt{2}\times \bar{g}^{\frac{1}{2}-\eta}
\end{eqnarray*}
since $\eta\ge 0$ and $\bar{g}\le 1$. The standard hypothesis (\ref{logstdeq}) then gives us the desired $\frac{1}{3}\log\left(\frac{1+\sqrt{2}}{2}\right)$ upper bound.
In the same way, we have
\begin{eqnarray*}
|a_1|\times|a_4|^{-\frac{1}{4}} & \le & t^{\frac{3}{4}}
\left[ |\beta_{1,\Delta}|+6 |\beta_{3,\Delta}|\right]\times |\beta_{4,\Delta}|^{-\frac{1}{4}}\\
 & \le & 7\bar{g}^{1-\eta} \times \left(\frac{1}{2}\bar{g}\right)^{-\frac{1}{4}} \\
 & \le & \frac{1}{3}\log\left(\frac{1+\sqrt{2}}{2}\right)
\end{eqnarray*}
by (\ref{logstdeq}).
Now 
\[
|a_0|\le 2|\beta_{2,\Delta}|+12|\beta_{4,\Delta}|\le
2\bar{g}^{1-\eta}+18\bar{g}\le\log 2
\]
by (\ref{expgstdeq}).
If one takes $\theta$ in Lemma \ref{L6lem} to be $h$, then the needed condition from part 1) of that lemma
is equivalent to $|a_4|\le 2\bar{g}$ which holds since $|a_4|=t|\beta_{4,\Delta}|\le \frac{3}{2}\bar{g}$.
Likewise if one takes $\theta=h_\ast$
and wants to use part 2) of the lemma then the required condition
is $2^{\frac{3}{4}} L^{\left(\frac{3+\epsilon}{4}\right)}\le \frac{(\sqrt{2}-1)^2}{e} |a_4|^{-\frac{1}{4}}$.
Since $L\ge 2$ and $\epsilon\in(0,1]$, it is enough to have
\[
|a_4|^{\frac{1}{4}}
\le L^{-1}2^{-\frac{3}{4}}\times  \frac{(\sqrt{2}-1)^2}{e}\ .
\]
Using $|a_4|\le 2\bar{g}$, we see that the requirement follows from the standard hypothesis (\ref{g14stdeq}).
The desired inequalities now follow from Lemma \ref{L6lem}.
\qed

\begin{Lemma}\label{L8lem}
For all $t\in (0,1]$ and all unit cube $\Delta$ we have
\[
\forall \psi\in\mathbb{R}\ ,\ \ ||e^{-t\tilde{V}_\Delta(\psi)}||_{\partial\psi,\psi,h}
\le 2^{-\frac{t}{2}(\Re\beta_{4,\Delta})\psi^4}
\]
as well as
\[
|e^{-t\tilde{V}_\Delta(\psi)}|_{\partial\psi,h_\ast}\le 2\ .
\]
\end{Lemma}

\noindent{\bf Proof:}
Recall that $\tilde{V}_\Delta$ is defined in the same way as $V_\Delta$ except that the Wick ordering is with respect
to $C_1$ instead of $C_0$. We again use Lemma \ref{L6lem} in order to prove the wanted result.
The formulas for the $a_k$'s are exactly the same as in the previous lemma apart for changing $C_0(0)$
to $C_1(0)=L^{-2[\phi]}C_0(0)$ which is also bounded by 2. Since the latter property is the only thing we used about $C_0$, the
present lemma then follows in the same manner as the previous one.
\qed

\begin{Lemma}\label{L9lem}
For all unit cube $\Delta'$ and for all subset $Y_0\subset[L^{-1}\Delta']$
we have
\[
\forall \phi\in\mathbb{R}\ ,\ \ 
\left|\left|\prod_{\Delta\in Y_0} e^{-\tilde{V}_\Delta(\phi_1)}\right|\right|_{\partial\phi,\phi,h}
\le 2
\]
as well as
\[
\left|\prod_{\Delta\in Y_0} e^{-\tilde{V}_\Delta(\phi_1)}\right|_{\partial\phi,h_\ast}\le 2\ .
\]
If $|Y_0|\ge\frac{L^3}{2}$ (which holds if $|Y_0|=L^3$ or $L^3-1$ because $L\ge 2$) then we have the improved bound
\[
\forall \phi\in\mathbb{R}\ ,\ \ 
\left|\left|\prod_{\Delta\in Y_0} e^{-\tilde{V}_\Delta(\phi_1)}\right|\right|_{\partial\phi,\phi,h}
\le 2 e^{-\frac{\bar{g}}{16}\phi^4}\ .
\]
Here $\phi_1$ denotes the rescaled field $L^{-[\phi]}\phi$.
\end{Lemma}

\noindent{\bf Proof:}
The argument is similar to the previous lemmas provided one keeps in mind that $\phi_1=L^{-[\phi]}\phi$
but derivatives are with respect to $\phi$. In the degenerate case where $Y_0$ is empty there is nothing to prove
so we will assume that $|Y_0|>0$.
By definition
\[
\prod_{\Delta\in Y_0} e^{-\tilde{V}_\Delta(\phi_1)}=e^{-U(\phi)}
\]
where $U(\phi)=\sum_{k=0}^{4} a_k \phi^k$ with
\begin{eqnarray*}
a_4 & = & \sum_{\Delta\in Y_0} L^{-4[\phi]}\beta_{4,\Delta} \\
a_3 & = & \sum_{\Delta\in Y_0} L^{-3[\phi]}\beta_{3,\Delta} \\
a_2 & = & \sum_{\Delta\in Y_0} 
\left( L^{-2[\phi]}\beta_{2,\Delta}-6L^{-4[\phi]} C_0(0)\beta_{4,\Delta}\right) \\
a_1 & = & \sum_{\Delta\in Y_0} 
\left( L^{-[\phi]}\beta_{1,\Delta}-3L^{-3[\phi]} C_0(0)\beta_{3,\Delta}\right) \\
a_0 & = & \sum_{\Delta\in Y_0} 
\left(- L^{-2[\phi]}C_0(0)\beta_{2,\Delta}+3L^{-4[\phi]} C_0(0)^2\beta_{4,\Delta}\right) \ .\\
\end{eqnarray*}
From
\[
a_4=|Y_0|L^{-4[\phi]}\bar{g}+\sum_{\Delta\in Y_0} L^{-4[\phi]}(\beta_{4,\Delta}-\bar{g})
\]
we get
\[
\left|a_4-|Y_0|L^{-4[\phi]}\bar{g}\right|\le \frac{1}{2}\bar{g}|Y_0|L^{-4[\phi]}
\]
and therefore
\[
|a_4|\ge  \frac{1}{2}\bar{g}|Y_0|L^{-4[\phi]}>0\ .
\]
Since $a_4$ is $|Y_0|L^{-4[\phi]}$ times a barycenter of elements in the convex set
$\{\beta\in\mathbb{C}|\ |\beta-\bar{g}|<\frac{1}{2}\bar{g}\}$, we easily see that $\Re a_4\ge \frac{\sqrt{3}}{2} |a_4|$ holds.
We also have
\[
|a_3|\times|a_4|^{-\frac{3}{4}}\le |Y_0|L^{-3[\phi]}\bar{g}^{1-\eta
\left(\frac{1}{2}|Y_0|\bar{g}L^{-4[\phi]}\right)^{-\frac{3}{4}}}
\]
\[
\le 2^{\frac{3}{4}}|Y_0|^{\frac{1}{4}}\bar{g}^{\frac{1}{4}-\eta}\le
2^{\frac{3}{4}}L^{\frac{3}{4}}\bar{g}^{\frac{1}{4}-\eta}\le\frac{1}{3}\log\left(\frac{1+\sqrt{2}}{2}\right)
\]
by (\ref{logstdeq}).
Likewise
\[
|a_2|\times|a_4|^{-\frac{2}{4}}\le |Y_0|
\left(L^{-2[\phi]}\bar{g}^{1-\eta}+L^{-4[\phi]}\times 12\times \frac{3}{2} \bar{g}\right)
\times \left(\frac{1}{2}|Y_0|\bar{g} L^{-4[\phi]}\right)^{-\frac{1}{2}}
\]
\[
\le 2^{\frac{1}{2}}|Y_0|^{\frac{1}{2}}\left(\bar{g}^{\frac{1}{2}-\eta}+18L^{-2[\phi]}\bar{g}^{\frac{1}{2}}\right)
\le 19\sqrt{2} L^{\frac{3}{2}} \bar{g}^{\frac{1}{2}-\eta}
\]
which is bounded using (\ref{logstdeq}).
Also,
\[
|a_1|\times|a_4|^{-\frac{1}{4}}\le |Y_0|
\left(L^{-[\phi]}\bar{g}^{1-\eta}+6L^{-3[\phi]} \bar{g}^{1-\eta}\right)
\times \left(\frac{1}{2}|Y_0|\bar{g}L^{-4[\phi]}\right)^{-\frac{1}{4}}
\]
\[
\le 2^{\frac{1}{4}}\times 7\times |Y_0|^{\frac{3}{4}}\bar{g}^{\frac{3}{4}-\eta}
\le 2^{\frac{1}{4}}\times 7\times L^{\frac{9}{4}} \bar{g}^{\frac{3}{4}-\eta}\le 
\frac{1}{3}\log\left(\frac{1+\sqrt{2}}{2}\right)
\]
by (\ref{logstdeq}).
Finally, 
\[
|a_0|\le |Y_0|\left(
2 L^{-2[\phi]}\bar{g}^{1-\eta}+12 L^{-4[\phi]}\times\frac{1}{2}\bar{g}
\right)
\]
\[
\le 20 L^{3-2[\phi]} \bar{g}^{1-\eta}\le 20 L^2 \bar{g}^{1-\eta}\le \log 2
\]
by (\ref{g14stdeq}) and using $L^{\frac{3+\epsilon}{2}}\le L^2$.

In order to apply Lemma \ref{L6lem} 1) it is enough to have $|a_4|\le 2\bar{g}$.
Since
\[
|a_4|\le L^{3-4[\phi]}\times \frac{3}{2}\bar{g}=\frac{3}{2} L^{\epsilon}\bar{g}\ ,
\]
all we need is $L^{\epsilon}\le\frac{4}{3}$ which is a condition in (\ref{g14stdeq}).
As for the use of Lemma \ref{L6lem} 2), a condition in (\ref{g14stdeq}), namely,
\[
\bar{g}^{\frac{1}{4}}\le\frac{(\sqrt{2}-1)^2}{2e} L^{-1}
\]
implies
\[
|a_4|^{\frac{1}{4}}\le (2\bar{g})^{\frac{1}{4}}\le 
\frac{(\sqrt{2}-1)^2}{e} L^{-1}\times 2^{-\frac{3}{4}}
\le \frac{(\sqrt{2}-1)^2}{e} L^{-\left(\frac{3+\epsilon}{4}\right)}\times 2^{-\frac{3}{4}}
\]
which is the required relationship between $|a_4|$ and $h_\ast$.

As a result of the previous considerations and following the application of Lemma \ref{L6lem}
we arrive at
\begin{equation}
\left|\left|\prod_{\Delta\in Y_0} e^{-\tilde{V}_\Delta(\phi_1)}\right|\right|_{\partial\phi,\phi,h}
\le 2 e^{-\frac{1}{2}(\Re a_4)\phi^4}\le 2
\label{crudeL9eq}
\end{equation}
since $\Re a_4\ge\frac{\sqrt{3}}{2}|a_4|>0$, as well as 
\[
\left|\prod_{\Delta\in Y_0} e^{-\tilde{V}_\Delta(\phi_1)}\right|_{\partial\phi,,h_\ast}
\le 2\ .
\]

Now we impose the stronger hypothesis $|Y_0|\ge \frac{L^3}{2}$.
Then
\[
\Re a_4=\sum_{\Delta\in Y_0} L^{-4[\phi]}\Re \beta_{4,\Delta}
\ge |Y_0|L^{-4[\phi]}\times \frac{1}{2}\bar{g}
\]
\[
\ge \frac{1}{4}L^{3-4[\phi]} \bar{g}\ge \frac{1}{4}\bar{g}
\]
because $L^\epsilon> 1$.
As a result, the crude bound by $2$ in (\ref{crudeL9eq}) can be amended to
\[
\left|\left|\prod_{\Delta\in Y_0} e^{-\tilde{V}_\Delta(\phi_1)}\right|\right|_{\partial\phi,\phi,h}
\le 2 e^{-\frac{1}{16}\bar{g}\phi^4}\ .
\]
\qed

The special case $Y_0=[L^{-1}\Delta']$ of the previous lemma, which is important for bounding $e^{-\hat{V}}$ expressions,
will be stated as a separate lemma. It simply follows from the observation
\[
\prod_{\Delta\in Y_0} e^{-\tilde{V}_\Delta(\phi_1)}=e^{-\hat{V}_{\Delta'}(\phi)}
\]
when $Y_0=[L^{-1}\Delta']$.

\begin{Lemma}\label{L10lem}
For all unit boxes $\Delta'$, we have
\[
\forall \phi\in\mathbb{R},\ \   ||e^{-\hat{V}_{\Delta'}(\phi)}||_{\partial\phi,\phi,h}
\le 2 e^{-\frac{1}{16}\bar{g}\phi^4}\ .
\]
and
\[
|e^{-\hat{V}_{\Delta'}(\phi)}|_{\partial\phi,h_\ast}
\le 2\ .
\]
\end{Lemma}

\begin{Lemma}\label{L11lem}
For all unit cube $\Delta$, for all $\kappa,\gamma\in(0,1]$
and for all $\phi,\zeta\in\mathbb{R}$ we have
\[
||p_{\Delta}(\phi_1,\zeta)||_{\partial\phi,\phi,h}\le
\mathcal{O}_1 \kappa^{-2}\gamma^{-\frac{3}{4}}\bar{g}^{\frac{1}{4}-\eta}
e^{\kappa\zeta_{\Delta}^2} e^{\gamma(\Re\beta_{4,\Delta})\phi_1^4}
\]
where $\mathcal{O}_1=54600$.
\end{Lemma}

\noindent{\bf Proof:}
Recall that
\[
p_{\Delta}(\phi_1,\zeta)=
\sum_{a,b}
\bbone\left\{
\begin{array}{c}
a+b\le 4\\
a\ge 0\ ,\ b\ge 1
\end{array}
\right\}
\frac{(a+b)!}{a!\ b!}\ \beta_{a+b,\Delta}
:\phi_1^a:_{C_1}\ \times
\ :\zeta^b:_{\Gamma}
\]
where in fact $:\phi_1^a:_{C_1}$ means $\left. :\psi^a:_{C_1}\right|_{\psi=\phi_1}$ with $\phi_1=L^{-[\phi]}\phi$.
We use
\begin{equation}
||p_{\Delta}(\phi_1,\zeta)||_{\partial\phi,\phi,h}\le
||p_{\Delta}(\psi,\zeta)||_{\partial\psi,\phi_1,h}
\label{chainruleeq}
\end{equation}
because by the chain rule each derivative with respect to $\phi$ brings a factor $L^{-[\phi]}\le 1$ times the corresponding
derivative with respect to $\psi$.
We therefore have 
\[
||p_{\Delta}(\phi_1,\zeta)||_{\partial\phi,\phi,h}\le
\sum_{a,b}
\bbone\left\{
\begin{array}{c}
a+b\le 4\\
a\ge 0\ ,\ b\ge 1
\end{array}
\right\}
\frac{(a+b)!}{a!\ b!}\ |\beta_{a+b,\Delta}|
\times
||:\psi^a:_{C_1}||_{\partial\psi,\phi_1,h}\ \times
\ |:\zeta^b:_{\Gamma}|\ .
\]
We undo the Wick ordering noting that $1\le b\le 4$ and $0\le a\le 3$.
From
\begin{eqnarray*}
:\psi^3:_{C_1} & = & \psi^3-3L^{-2[\phi]}C_0(0)\psi \\
:\psi^2:_{C_1} & = & \psi^2-L^{-2[\phi]}C_0(0) \\
:\psi^1:_{C_1} & = & \psi \\
:\psi^0:_{C_1} & = & 1
\end{eqnarray*}
and the bound $C_0(0)<2$
we get
\[
||:\psi^a:_{C_1}||_{\partial\psi,\phi_1,h}\le 7\max_{0\le k\le a}
||\psi^k||_{\partial\psi,\phi_1,h}
\]
for $0\le a\le 3$.
Since, by Lemma \ref{L1lem}, $\Gamma(0)=C_0(0)(1-L^{-2[\phi]})<C_0(0)<2$, the explicit formulas
\begin{eqnarray*}
:\zeta^4:_{\Gamma} & = & \zeta^4-6\Gamma(0)\zeta^2+3\Gamma(0)^2 \\
:\zeta^3:_{\Gamma} & = & \zeta^3-3\Gamma(0)\zeta \\
:\zeta^2:_{\Gamma} & = & \zeta^2-\Gamma(0) \\
:\zeta^1:_{\Gamma} & = & \zeta
\end{eqnarray*}
similarly imply
\[
|:\zeta^b:_{\Gamma}|\le 25 \max_{0\le j\le b} |\zeta|^j
\]
for $1\le b\le 4$.
The $|\beta|$'s are bounded by $\frac{3}{2}\bar{g}$ or $\bar{g}^{1-\eta}$. Since $\eta\ge 0$ and $\bar{g}\le 1$
we use the uniform bound by the worst case scenario $\frac{3}{2}\bar{g}^{1-\eta}$.
Thus
\[
||p_{\Delta}(\phi_1,\zeta)||_{\partial\phi,\phi,h}\le
7\times 25\times
\sum_{a,b}
\bbone\left\{
\begin{array}{c}
a+b\le 4\\
a\ge 0\ ,\ b\ge 1
\end{array}
\right\}
\frac{(a+b)!}{a!\ b!}\ |\beta_{a+b,\Delta}|
\]
\[
\times \frac{3}{2}\bar{g}^{1-\eta}
\times \max_{0\le k\le a}
||\psi^k||_{\partial\psi,\phi_1,h}
\times \max_{0\le j\le b} |\zeta|^j\ .
\]
The binomial coefficients appearing in the sum add up to 26 and therefore using a maximum over $a$, $b$ instead of a sum
gives rise to the numerical coefficient $7\times 25\times \frac{3}{2}\times 26=6825$, namely,
\[
||p_{\Delta}(\phi_1,\zeta)||_{\partial\phi,\phi,h}\le 6825 \bar{g}^{1-\eta}
\times
\max_{a,b}\left[\max_{0\le k\le a}
||\psi^k||_{\partial\psi,\phi_1,h}
\times \max_{0\le j\le b} |\zeta|^j\right]
\]
where the maximum is over pairs of integers which satisfy $a\ge 0$, $b\ge 1$ and $a+b\le 4$.
By Lemma \ref{L5lem} and for $\kappa\in(0,1]$ and $b\le 4$ we have
\[
\max_{0\le j\le b} |\zeta|^j\le \kappa^{-2} e^{\kappa\zeta^2} \max_{0\le j\le 4} \left(\frac{j}{2e}\right)^{\frac{j}{2}}
\le \kappa^{-2} e^{\kappa\zeta^2}\ .
\]
For $0\le k\le 3<9$
we have
\[
||\psi^k||_{\partial\psi,\psi,h}=(h+|\psi|)^k=\sum_{n=0}^k \left(\begin{array}{c}
k \\ n
\end{array}\right)
(c_1 \bar{g}^{-\frac{1}{4}})^{k-n} |\psi|^n\ .
\]
We use Lemma \ref{L4lem} to write
\[
|\psi|^n\le \left(\frac{n}{2e}\right)^{\frac{n}{4}} [\gamma\bar{g}]^{-\frac{n}{4}}
e^{\gamma(\Re\beta_{4,\Delta})\psi^4}
\]
We drop the numerical factor since $n\le 3<2e$ and use $\gamma\in(0,1]$ to arrive at
\[
||\psi^k||_{\partial\psi,\psi,h}\le \gamma^{-\frac{3}{4}} e^{\gamma(\Re\beta_{4,\Delta})\psi^4}
\bar{g}^{-\frac{k}{4}}\times \sum_{n=0}^k \left(\begin{array}{c}
k \\ n
\end{array}\right) c_1^{k-n}\ .
\]
The last sum reduces to $(1+c_1)^k<2^k\le 8$.
Therefore
\[
\max_{0\le k\le 3}
||\psi^k||_{\partial\psi,\phi_1,h}
\le 8 \gamma^{-\frac{3}{4}} \bar{g}^{-\frac{3}{4}} e^{\gamma(\Re\beta_{4,\Delta})\phi_1^4}
\]
and the result follows.
\qed

\begin{Lemma}\label{L12lem}
For all $\kappa$ such that $0<\kappa\le 2^{-\frac{3}{2}}L^{-(3-2[\phi])}$ and for all $\zeta\in\mathbb{R}$
we have
\[
|p_{\Delta}(\phi_1,\zeta)|_{\partial\phi,h_\ast}\le
\mathcal{O}_2 \kappa^{-2}\bar{g}^{1-\eta}
e^{\kappa\zeta_{\Delta}^2}
\]
where $\mathcal{O}_2=6825$.
\end{Lemma}

\noindent{\bf Proof:}
We proceed as in the previous lemma, thus arriving
at the bound
\[
|p_{\Delta}(\phi_1,\zeta)|_{\partial\phi,h_\ast}\le
|p_{\Delta}(\psi,\zeta)|_{\partial\psi,h_\ast}
\]
\[
\le
6825\bar{g}^{1-\eta}\times
\max_{a,b}\left[\max_{0\le k\le a}
|\psi^k|_{\partial\psi,h_\ast}
\times \max_{0\le j\le b} |\zeta|^j\right]
\]
with the same conditions on $a$ and $b$.
However, now $|\psi^k|_{\partial\psi,h_\ast}=h_{\ast}^k$. By again using Lemma \ref{L5lem}
we get
\[
|p_{\Delta}(\phi_1,\zeta)|_{\partial\phi,h_\ast}\le
6825\bar{g}^{1-\eta}\times
\max_{\substack{k,j\ge 0 \\ k+j\le 4}}
\left[
h_{\ast}^k \kappa^{-\frac{j}{2}}\left(\frac{j}{2e}\right)^{\frac{j}{2}} e^{\kappa\zeta^2}
\right]\ .
\]
We drop the cumbersome factor $\left(\frac{j}{2e}\right)^{\frac{j}{2}}\le 1$ and note that the hypothesis
on $\kappa$ ensures that $h_{\ast}\le \kappa^{-\frac{1}{2}}$. This implies
\[
h_{\ast}^k \kappa^{-\frac{j}{2}}\le \kappa^{-\frac{k+j}{2}}\le \kappa^{-2}
\]
because of the condition $k+j\le 4$ and $\kappa\le 1$ which also follows from the hypothesis.
\qed

\begin{Lemma}\label{L13lem}
For all $\kappa\in(0,3]$ we have
\[
||r_{1,\Delta}(\phi_1,\zeta)||_{\partial\phi,\phi,h}\le 
\mathcal{O}_3 \kappa^{-6} e^{\kappa\zeta_\Delta^2} \bar{g}^{\frac{3}{4}-3\eta}
\]
where $\mathcal{O}_3=2^{\frac{21}{4}}\times 3^{\frac{29}{4}}\times\mathcal{O}_1^3$.
\end{Lemma}

\noindent{\bf Proof:}
By definition
\[
r_{1,\Delta}(\phi_1,\zeta)=e^{-\tilde{V}_\Delta(\phi_1)}\left[
e^{-p\Delta}-1+p_\Delta-\frac{1}{2}p_\Delta^2
\right]
\]
where $p_\Delta$ is shorthand for $p_\Delta(\phi_1,\zeta)$.
This is the third order Taylor remainder
when expanding $e^{-\tilde{V}_\Delta(\phi_1)-s p_\Delta}$ at $s=1$ about 0.
We can also write
\[
r_{1,\Delta}(\phi_1,\zeta)=\int_0^1 {\rm d}t\ \frac{(1-t)^2}{2} (-p_\Delta)^3 
e^{-\tilde{V}_\Delta(\phi_1)-t p_\Delta}\ .
\]
But
\begin{eqnarray*}
\tilde{V}_\Delta(\phi_1)+t p_\Delta & = & \tilde{V}_\Delta(\phi_1) +t\left(
V_\Delta(\phi_1+\zeta)-\tilde{V}_\Delta(\phi_1)
\right) \\
 & = & t V_\Delta(\phi_1+\zeta) +(1-t)\tilde{V}_\Delta(\phi_1)\ .
\end{eqnarray*}
By the multiplicative property of the seminorms
\begin{eqnarray*}
||e^{-\tilde{V}_\Delta(\phi_1)-t p_\Delta}||_{\partial\phi,\phi,h} & = & 
||e^{ t V_\Delta(\phi_1+\zeta) } \times e^{(1-t)\tilde{V}_\Delta(\phi_1)}||_{\partial\phi,\phi,h} \\
 & \le & ||e^{ t V_\Delta(\phi_1+\zeta) } ||_{\partial\phi,\phi,h} \times
|| e^{(1-t)\tilde{V}_\Delta(\phi_1)}||_{\partial\phi,\phi,h}
\end{eqnarray*}
and thus
\[
||r_{1,\Delta}(\phi_1,\zeta)||_{\partial\phi,\phi,h}
\le \int_0^1 {\rm d}t\ \frac{(1-t)^2}{2} ||(-p_\Delta)^3 
e^{-\tilde{V}_\Delta(\phi_1)-t p_\Delta}||_{\partial\phi,\phi,h}
\]
\[
\le \frac{1}{2}\int_0^1 {\rm d}t\ \frac{(1-t)^2}{2} 
|| p_\Delta||_{\partial\phi,\phi,h}^3\times
||e^{ t V_\Delta(\phi_1+\zeta) } ||_{\partial\phi,\phi,h} \times
|| e^{(1-t)\tilde{V}_\Delta(\phi_1)}||_{\partial\phi,\phi,h}\ .
\]
Using the same inequality comparing derivatives with respect to $\phi$ versus $\psi=\phi_1$
as in (\ref{chainruleeq}) we obtain
\[
||e^{ t V_\Delta(\phi_1+\zeta) } ||_{\partial\phi,\phi,h}\le
||e^{ t V_\Delta(\psi+\zeta) } ||_{\partial\psi,\phi_1,h}=
||e^{ t V_\Delta(\psi) } ||_{\partial\psi,\phi_1+\zeta,h}
\]
\begin{equation}
\le 2 e^{-\frac{t}{2}(\Re\beta_{4,\Delta})(\phi_1+\zeta)^4}\le 2
\label{crudeL13eq}
\end{equation}
thanks to Lemma \ref{L7lem}.
Likewise
\[
|| e^{(1-t)\tilde{V}_\Delta(\phi_1)}||_{\partial\phi,\phi,h} \le 
|| e^{(1-t)\tilde{V}_\Delta(\psi)}||_{\partial\psi,\phi_1,h} 
\le  2 e^{-\frac{(1-t)}{2}(\Re\beta_{4,\Delta})\phi_1^4}
\]
by Lemma \ref{L8lem}.
Although 
$||p_\Delta||_{\partial\phi,\phi,h}$ does not depend on $t$, we bound it in a $t$ dependent way
using Lemma \ref{L11lem} with
$\gamma=\frac{1-t}{6}$ and $\frac{\kappa}{3}$ instead of $\kappa$.
Namely, we write
\[
||p_\Delta||_{\partial\phi,\phi,h}
\le \mathcal{O}_1\times 9\times \kappa^{-2}\times (1-t)^{-\frac{3}{4}}\times
6^{\frac{3}{4}}
\times e^{\frac{\kappa}{3}\zeta_\Delta^2}
e^{\frac{(1-t)}{6}(\Re\beta_{4,\Delta})\phi_1^4}
\times \bar{g}^{\frac{1}{4}-\eta}\ .
\]
Altogether this results in the bound
\[
||r_{1,\Delta}(\phi_1,\zeta)||_{\partial\phi,\phi,h}\le
\frac{1}{2}\mathcal{O}_1^3\times 9^3\times 6^{\frac{9}{4}} \kappa^{-6} e^{\kappa\zeta_\Delta^2}
\times 4\times \bar{g}^{\frac{3}{4}-3\eta}\times \int_0^1 (1-t)^{-\frac{1}{4}}\ {\rm d}t
\]
which is the desired result.
\qed

\begin{Lemma}\label{L14lem}
For all $\kappa$ such that $0<\kappa\le 2^{-\frac{3}{2}}\times 3\times L^{-(3-2[\phi])}$ and for all $\zeta\in\mathbb{R}$ we have
\[
|r_{1,\Delta}(\phi_1,\zeta)|_{\partial\phi,h_\ast}\le
\mathcal{O}_4 \kappa^{-6} \bar{g}^{3-3\eta} e^{\kappa \zeta_\Delta^2}
\]
where
$\mathcal{O}_4=2\times 3^5\times \mathcal{O}_2^3$.
\end{Lemma}

\noindent{\bf Proof:}
Proceeding as in the previous lemma we arrive at
\[
|r_{1,\Delta}(\phi_1,\zeta)|_{\partial\phi,h_\ast}
\le \int_0^1 {\rm d}t\ \frac{(1-t)^2}{2} |p_\Delta|_{\partial\phi,h_\ast}^3\times
|e^{ t V_\Delta(\phi_1+\zeta) }|_{\partial\phi,h_\ast} \times
| e^{(1-t)\tilde{V}_\Delta(\phi_1)}|_{\partial\phi,h_\ast}\ .
\]
We use
\[
|e^{ t V_\Delta(\phi_1+\zeta) }|_{\partial\phi,h_\ast}=
||e^{ t V_\Delta(\phi_1+\zeta) }||_{\partial\phi,0,h_\ast}
\le ||e^{ t V_\Delta(\phi_1+\zeta) }||_{\partial\phi,0,h}
\]
since $h_\ast\le h$. Indeed this follows from $c_2L^{\frac{3+\epsilon}{4}}\le c_2 L\le c_1\bar{g}^{-\frac{1}{4}}$, i.e., from (\ref{c1c2stdeq})
We can then reuse the bound (\ref{crudeL13eq}) at $\phi=0$,
namely,
\[
|e^{ t V_\Delta(\phi_1+\zeta) }|_{\partial\phi,h_\ast}\le 2\ .
\]
Now
\[
| e^{(1-t)\tilde{V}_\Delta(\phi_1)}|_{\partial\phi,h_\ast}\le
| e^{(1-t)\tilde{V}_\Delta(\psi)}|_{\partial\psi,h_\ast}\le 2
\]
by Lemma \ref{L8lem}.
Finally, we use Lemma \ref{L12lem} with $\frac{\kappa}{3}$ instead of $\kappa$ in order to get
\[
|p_\Delta|_{\partial\phi,h_\ast}\le \mathcal{O}_2 \times 9\times \kappa^{-2} e^{\frac{\kappa}{3}\zeta_\Delta^2} \bar{g}^{1-\eta}\ .
\]
Altogether this results in
\[
|r_{1,\Delta}(\phi_1,\zeta)|_{\partial\phi,h_\ast}
\le \frac{1}{2}\left( \int_0^1 {\rm d}t\ \frac{(1-t)^2}{2}\right)\times 4\times \mathcal{O}_2^3
\times 3^6 \kappa^{-6} e^{\kappa\zeta_\Delta^2} \bar{g}^{3-3\eta}
\]
which is the wanted bound.
\qed

\begin{Lemma}\label{L15lem}
For all $\kappa\in (0,1]$, and all $\lambda\in\mathbb{C}$ which satisfies $|\lambda|\bar{g}^{\frac{1}{4}-\eta}\le 1$,
we have $\forall \phi, \zeta\in\mathbb{R}$
\[
||P_\Delta(\lambda,\phi_1,\zeta)||_{\partial\phi,\phi,h}\le
\mathcal{O}_5 \kappa^{-6} e^{\kappa\zeta_\Delta^2}\times |\lambda|\bar{g}^{\frac{1}{4}-\eta}
\]
where
$\mathcal{O}_5=2^{\frac{9}{2}}\mathcal{O}_1+2^7 \mathcal{O}_1^2+\mathcal{O}_3$.
\end{Lemma}

\noindent{\bf Proof:}
By definition
\[
P_\Delta(\lambda,\phi_1,\zeta)=e^{-\tilde{V}_\Delta(\phi_1)}\left[
-\lambda p_\Delta+\frac{\lambda^2}{2} p_{\Delta}^2
\right] +\lambda^3 r_{1,\Delta}(\phi_1,\zeta)
\]
and thus from the properties of the seminorm we get
\[
||P_\Delta(\lambda,\phi_1,\zeta)||_{\partial\phi,\phi,h}\le
||e^{-\tilde{V}_\Delta(\phi_1)}||_{\partial\phi,\phi,h}\times
\left[
|\lambda|\ ||p_{\Delta}(\phi_1,\zeta)||_{\partial\phi,\phi,h}
+\frac{|\lambda|^2}{2}\ ||p_{\Delta}(\phi_1,\zeta)||_{\partial\phi,\phi,h}^2
\right]
\]
\[
+ |\lambda|^3\ ||r_{1,\Delta}(\phi_1,\zeta)||_{\partial\phi,\phi,h}\ .
\]
We bound $||p_{\Delta}(\phi_1,\zeta)||_{\partial\phi,\phi,h}$ using Lemma \ref{L11lem} with $\gamma=\frac{1}{4}$ and with
$\frac{\kappa}{2}$ instead of $\kappa$. We bound $||e^{-\tilde{V}_\Delta(\phi_1)}||_{\partial\phi,\phi,h}$ using Lemma \ref{L7lem}.
Finally, we bound $||r_{1,\Delta}(\phi_1,\zeta)||_{\partial\phi,\phi,h}$ using Lemma \ref{L13lem}.
Put together, this results in
\[
||P_\Delta(\lambda,\phi_1,\zeta)||_{\partial\phi,\phi,h}\le
2e^{-\frac{1}{2}(\Re\beta_{4,\Delta})\phi_1^4}
\times\left[
|\lambda|\mathcal{O}_1 2^{2+\frac{3}{2}}\bar{g}^{\frac{1}{4}-\eta}\kappa^{-2} e^{\frac{\kappa}{2}\zeta_{\Delta}^2}
e^{\frac{1}{4}(\Re\beta_{4,\Delta})\phi_1^4}
\right.
\]
\[
\left.
\frac{1}{2}|\lambda|^2 \mathcal{O}_1^2 2^{7}\bar{g}^{\frac{1}{2}-2\eta}\kappa^{-4} e^{\kappa\zeta_{\Delta}^2}
e^{\frac{1}{2}(\Re\beta_{4,\Delta})\phi_1^4}
\right]
+|\lambda|^3 \mathcal{O}_3 \bar{g}^{\frac{3}{4}-3\eta}\kappa^{-6} e^{\kappa\zeta_{\Delta}^2}\ .
\]
Since $0<\kappa\le 1$, $\kappa^{-6}$ is the worst power of that kind. Since also $\Re\beta_{4,\Delta}>0$,
the worst exponential factor left is $e^{\kappa\zeta_\Delta^2}$. Hence
\[
||P_\Delta(\lambda,\phi_1,\zeta)||_{\partial\phi,\phi,h}\le
\kappa^{-6} e^{\kappa\zeta_\Delta^2}\times
\left[
2^{\frac{9}{2}}\mathcal{O}_1 |\lambda| \bar{g}^{\frac{1}{4}-\eta}
+2^7 \mathcal{O}_1^2 |\lambda|^2 \bar{g}^{\frac{1}{2}-2\eta}
+\mathcal{O}_3 |\lambda|^3 \bar{g}^{\frac{3}{4}-3\eta}
\right]\ .
\]
We then conclude using the hypothesis $|\lambda| \bar{g}^{\frac{1}{4}-\eta}\le 1$.
\qed

\begin{Lemma}\label{L16lem}
For all $\kappa\in (0,2^{-\frac{1}{2}}L^{-(3-2[\phi])}]$,
and all $\lambda\in\mathbb{C}$ which satisfies $|\lambda|\bar{g}^{1-\eta}\le 1$,
we have $\forall \zeta\in\mathbb{R}$
\[
|P_\Delta(\lambda,\phi_1,\zeta)|_{\partial\phi,h_\ast}\le
\mathcal{O}_6 \kappa^{-6} e^{\kappa\zeta_\Delta^2}\times |\lambda|\bar{g}^{1-\eta}
\]
where
$\mathcal{O}_6=8\mathcal{O}_2+16\mathcal{O}_2^2+\mathcal{O}_4$.
\end{Lemma}

\noindent{\bf Proof:}
Similarly to the proof of the previous lemma we have
\[
|P_\Delta(\lambda,\phi_1,\zeta)|_{\partial\phi,h_\ast}\le
|e^{-\tilde{V}_\Delta(\phi_1)}|_{\partial\phi,h_\ast}\times
\left[
|\lambda|\ |p_{\Delta}(\phi_1,\zeta)|_{\partial\phi,h_\ast}
+\frac{|\lambda|^2}{2}\ |p_{\Delta}(\phi_1,\zeta)|_{\partial\phi,h_\ast}^2
\right]
\]
\[
+ |\lambda|^3\ |r_{1,\Delta}(\phi_1,\zeta)|_{\partial\phi,h_\ast}\ .
\]
We bound $|p_{\Delta}(\phi_1,\zeta)|_{\partial\phi,h_\ast}$ by Lemma \ref{L12lem}
with $\frac{\kappa}{2}$ instead of $\kappa$.
We have
\[
|e^{-\tilde{V}_\Delta(\phi_1)}|_{\partial\phi,h_\ast}\le
|e^{-\tilde{V}_\Delta(\psi)}|_{\partial\psi,h_\ast}\le 2
\]
thanks to Lemma \ref{L8lem}.
We also use Lemma \ref{L14lem} to bound $|r_{1,\Delta}(\phi_1,\zeta)|_{\partial\phi,h_\ast}$.
As a result we get
\[
|P_\Delta(\lambda,\phi_1,\zeta)|_{\partial\phi,h_\ast}\le
2\left[
|\lambda|\mathcal{O}_2\times 4\times \bar{g}^{1-\eta}\kappa^{-2} e^{\frac{\kappa}{2}\zeta_{\Delta}^2}
\right.
\]
\[
\left.
\frac{1}{2}|\lambda|^2 \mathcal{O}_2^2 2^{4}\bar{g}^{2-2\eta}\kappa^{-4} e^{\kappa\zeta_{\Delta}^2}
\right]
+|\lambda|^3 \mathcal{O}_4 \bar{g}^{3-3\eta}\kappa^{-6} e^{\kappa\zeta_{\Delta}^2}\ .
\]
The hypothesis on $\kappa$ implies that $\kappa\le 1$ and therefore the worst negative power
which appears is $\kappa^{-6}$. We also use the hypothesis $|\lambda|\bar{g}^{1-\eta} \le 1$
to bound the square and cube of that quantity as in the previous lemma and the result follows.
\qed

\begin{Lemma}\label{L17lem}
For all $\kappa\in(0,1]$, $\gamma\in(0,1]$,and $\phi,\zeta\in\mathbb{R}$ we have
\[
||Q_\Delta(\phi_1+\zeta)||_{\partial\phi,\phi,h}\le
\mathcal{O}_7 \bar{g}^{\frac{1}{2}-2\eta}\gamma^{-\frac{3}{2}} e^{\gamma(\Re\beta_{4,\Delta})\phi_1^4}
\kappa^{-3} e^{\kappa\zeta^2}
\]
where
\[
\mathcal{O}_7
= 331^2\times 96\times 2^6\times
\max_{0\le j\le 6} \left(\frac{j}{2e}\right)^{\frac{j}{2}}
\times \max_{0\le n\le 6} \left(\frac{n}{2e}\right)^{\frac{n}{4}}\ .
\]
\end{Lemma}

\noindent{\bf Proof:}
We proceed as in the proof of Lemma \ref{L10lem}.
By definition and by the elementary properties of Wick monomials
\[
Q_\Delta(\phi_1+\zeta)=
\sum_{a,b}\bbone\left\{
\begin{array}{c}
5\le a+b\le 6 \\
a,b\ge 0\end{array}
\right\}
\frac{(a+b)!}{a!b!} W_{a+b,\Delta}
\left.:\psi^a:_{C_1}\right|_{\psi=\phi_1}
:\zeta^b:_{\Gamma}\ .
\]
Therefore, again dominating $\phi$ derivatives by $\psi$ derivatives and using (\ref{fwstdeq}), we get
\[
||Q_\Delta(\phi_1+\zeta)||_{\partial\phi,\phi,h}\le
\bar{g}^{2-2\eta}
\sum_{a,b}\bbone\left\{
\begin{array}{c}
5\le a+b\le 6 \\
a,b\ge 0\end{array}
\right\}
\frac{(a+b)!}{a!b!}
||:\psi^a:_{C_1}||_{\partial\psi,\phi_1,h}
\times|:\zeta^b:_{\Gamma}|\ .
\]
In addition to the Wick ordering formulas in the proof of Lemma \ref{L11lem} we have
\begin{eqnarray*}
:\psi^4:_{C_1} & = & \psi^4-6L^{-2[\phi]}C_0(0)\psi^2+3L^{-4[\phi]}C_0(0)^2 \\
:\psi^5:_{C_1} & = & \psi^5-10L^{-2[\phi]}C_0(0)\psi^3+15L^{-4[\phi]}C_0(0)^2\psi \\
:\psi^6:_{C_1} & = & \psi^6-15L^{-2[\phi]}C_0(0)\psi^4+45  L^{-4[\phi]}C_0(0)^2\psi^2-15 L^{-6[\phi]}C_0(0)^3
\end{eqnarray*}
as well as
\begin{eqnarray*}
:\zeta^5:_{\Gamma} & = & \zeta^5-10\Gamma(0)\zeta^3+15 \Gamma(0)^2\zeta \\
:\zeta^6:_{\Gamma} & = & \zeta^6-15\Gamma(0)\zeta^4+45 \Gamma(0)^2\zeta^2-15 \Gamma(0)^3\ .
\end{eqnarray*}
Therefore when bounding these expressions using $\Gamma(0),C_0(0)\in [0,2]$, the worst numerical
factor coming from the sixth power case is $1+15\times 2+45\times 2^2+15\times 2^3=331$.
We therefore have
\begin{eqnarray*}
||:\psi^a:_{C_1}||_{\partial\psi,\phi_1,h} & \le & 331 \max_{0\le k\le a} ||\psi^k||_{\partial\psi,\phi_1,h} \\
|:\zeta^b:_{\Gamma}| & \le & 331 \max_{0\le j\le b} |\zeta|^j\ .
\end{eqnarray*}
This result in the rather coarse bound
\[
||Q_\Delta(\phi_1+\zeta)||_{\partial\phi,\phi,h}\le
331^2\bar{g}^{2-2\eta}
\sum_{a,b}\bbone\left\{
\begin{array}{c}
5\le a+b\le 6 \\
a,b\ge 0\end{array}
\right\}
\frac{(a+b)!}{a!b!}
\left(\max_{0\le k\le a} ||\psi^k||_{\partial\psi,\phi_1,h} \right)
\left(\max_{0\le j\le b} |\zeta|^j\right)
\]
\[
\le 331^2\times(2^5+2^6)\bar{g}^{2-2\eta}\max_{a,b}
\left[
\left(\max_{0\le k\le a} ||\psi^k||_{\partial\psi,\phi_1,h} \right)
\left(\max_{0\le j\le b} |\zeta|^j\right)
\right]
\]
where the new numerical factor $2^5+2^6=96$ comes from the sum of binomial coefficients and the maximum is over pairs
of nonnegative integers $a,b$ such that $a+b=5$ or $6$.
By Lemma \ref{L5lem}, and given that $\kappa\in (0,1]$, we have
\[
\max_{0\le j\le b} |\zeta|^j\le
\max_{0\le j\le 6} |\zeta|^j\le
\kappa^{-3} e^{\kappa\zeta^2}\times \max_{0\le j\le 6} \left(\frac{j}{2e}\right)^{\frac{j}{2}} \ . 
\]
For $0\le k\le 6<9$ we still have
\begin{eqnarray*}
||\psi^k||_{\partial\psi,\psi,h}  & = & (h+|\psi|)^k \\
 & \le & \sum_{n=0}^{k} \left(
\begin{array}{c}
k \\ n
\end{array} 
\right) (c_1 \bar{g}^{-\frac{1}{4}})^{k-n} \left(\frac{n}{2e}\right)^{\frac{n}{4}}
(\gamma\bar{g})^{-\frac{n}{4}} e^{\gamma(\Re\beta_{4,\Delta})\psi^4} \\
 & \le & \left(
\max_{0\le n\le 6} \left(\frac{n}{2e}\right)^{\frac{n}{4}} 
\right)\times (1+c_1)^k \bar{g}^{-\frac{k}{4}}\gamma^{-\frac{k}{4}} e^{\gamma(\Re\beta_{4,\Delta})\psi^4}
\end{eqnarray*}
again by Lemma \ref{L4lem}. We collect all these estimates and in the final result we
bound $(1+c_1)^k$ by $2^k\le 2^6$ and the powers of $\bar{g}$, $\gamma$ and $\kappa$
by their worst case values, i.e., respectively $\bar{g}^{\frac{1}{2}-2\eta}$, $\gamma^{-\frac{3}{2}}$ and $\kappa^{-3}$.
\qed

\begin{Lemma}\label{L18lem}
For all $\kappa\in(0,2^{-\frac{3}{2}}L^{-(3-2[\phi])}]$, and $\zeta\in\mathbb{R}$ we have
\[
|Q_\Delta(\phi_1+\zeta)|_{\partial\phi,h_\ast}\le
\mathcal{O}_8 \bar{g}^{2-2\eta}
\kappa^{-3} e^{\kappa\zeta^2}
\]
where
\[
\mathcal{O}_8
= 331^2\times 96\times
\max_{0\le j\le 6} \left(\frac{j}{2e}\right)^{\frac{j}{2}}\ .
\]
\end{Lemma}

\noindent{\bf Proof:}
As in the previous lemma we have
\[
|Q_\Delta(\phi_1+\zeta)|_{\partial\phi,h_\ast}\le
\bar{g}^{2-2\eta}
\sum_{a,b}\bbone\left\{
\begin{array}{c}
5\le a+b\le 6 \\
a,b\ge 0\end{array}
\right\}
\frac{(a+b)!}{a!b!}
|:\psi^a:_{C_1}|_{\partial\psi,h_\ast}
\times|:\zeta^b:_{\Gamma}|\ .
\]
We again have $|:\zeta^b:_{\Gamma}|\le 331 \max_{0\le j\le b} |\zeta|^j$
as well as $|:\psi^a:_{C_1}|_{\partial\psi,h_\ast}\le 331 \max_{0\le k\le a} |\psi^k|_{\partial\psi,h_\ast}$,
but now $|\psi^k|_{\partial\psi,h_\ast}=h_{\ast}^k$.
As a result
\[
|Q_\Delta(\phi_1+\zeta)|_{\partial\phi,h_\ast}\le
331^2\times 96\times
\bar{g}^{2-2\eta}
\max_{a,b}\left[
\left(\max_{0\le k\le a} h_{\ast}^k\right)\times
\left(\max_{0\le j\le b} |\zeta|^j\right)
\right]
\]
with the maximum again over of nonnegative integers $a,b$ such that $a+b=5$ or $6$.
Therefore
\[
|Q_\Delta(\phi_1+\zeta)|_{\partial\phi,h_\ast}\le
331^2\times 96\times
\bar{g}^{2-2\eta}
\max_{\substack{j,k\ge 0 \\ j+k\le 6}} h_{\ast}^{k} |\zeta|^j
\]
\[
\le
331^2\times 96\times
\bar{g}^{2-2\eta}\times e^{\kappa\zeta^2}
\max_{\substack{j,k\ge 0 \\ j+k\le 6}}
h_{\ast}^{k} \kappa^{-\frac{j}{2}}\left(\frac{j}{2e}\right)^{\frac{j}{2}}
\]
after applying Lemma \ref{L5lem}.
From the hypothesis on $\kappa$ and (\ref{hstarstdeq}) we have that $h_\ast\le \kappa^{-\frac{1}{2}}$
and therefore the quantities $h_{\ast}^{k} \kappa^{-\frac{j}{2}}$ are bounded by $\kappa^{-3}$
and the result follows.
\qed

\begin{Lemma}\label{L19lem}
For all $\phi,\zeta\in\mathbb{R}$ we have
\[
||Q_\Delta(\phi_1+\zeta) e^{-V_\Delta(\phi_1+\zeta)}||_{\partial\phi,\phi,h}\le
\mathcal{O}_9 \bar{g}^{\frac{1}{2}-2\eta}
\]
where
\[
\mathcal{O}_9=331\times 2^{\frac{19}{2}}\times\max_{0\le n\le 6} \left(\frac{n}{2e}\right)^{\frac{n}{4}}\ .
\]
\end{Lemma}

\noindent{\bf Proof:}
By definition
\[
Q_\Delta(\phi_1+\zeta) e^{-V_\Delta(\phi_1+\zeta)}=\sum_{k=5,6} W_{k,\Delta}
\left.
\left(:\psi^k:_{C_0} e^{-V_\Delta(\psi)}\right)
\right|_{\psi=\phi_1+\zeta}
\]
and therefore
\begin{eqnarray*}
||Q_\Delta(\phi_1+\zeta) e^{-V_\Delta(\phi_1+\zeta)}||_{\partial\phi,\phi,h}
 & \le &  \bar{g}^{2-2\eta} \sum_{k=5,6}
||:\psi^k:_{C_0} e^{-V_\Delta(\psi)}||_{\partial\psi,\phi_1+\zeta,h} \\
 & \le & 2\times 331\times \max_{0\le k\le 6}
||\psi^k||_{\partial\psi,\phi_1+\zeta,h}\ ||e^{-V_\Delta(\psi)}||_{\partial\psi,\phi_1+\zeta,h}
\end{eqnarray*}
by undoing the Wick ordering as in the proof of Lemma \ref{L17lem}.
By proceeding as in the latter, but with the specific choice $\gamma=\frac{1}{2}$ in the step involving
the application of Lemma \ref{L4lem}, we get
\[
||\psi^k||_{\partial\psi,\psi,h}
\le \left(\max_{0\le n\le 6} \left(\frac{n}{2e}\right)^{\frac{n}{4}} \right)
\times \left(\frac{\bar{g}}{2}\right)^{-\frac{k}{4}} e^{\frac{1}{2}(\Re\beta_{4,\Delta})\psi^4}\times 2^6 \ .
\]
Hence, by Lemma \ref{L7lem} with $t=1$, we get
\[
||\psi^k||_{\partial\psi,\phi_1+\zeta,h}\ ||e^{-V_\Delta(\psi)}||_{\partial\psi,\phi_1+\zeta,h}
\le \left(\max_{0\le n\le 6} \left(\frac{n}{2e}\right)^{\frac{n}{4}} \right)
\times \bar{g}^{-\frac{k}{4}}\times 2^{\frac{k}{4}+7}\ .
\]
Taking the worst case $k=6$ gives the desired bound.
\qed

\begin{Lemma}\label{L20lem}
If $0<\kappa\le 2^{-\frac{3}{2}} L^{-(3-2[\phi])}$ then for all $\zeta\in\mathbb{R}$ we have
\[
|Q_\Delta(\phi_1+\zeta) e^{-V_\Delta(\phi_1+\zeta)}|_{\partial\phi,h_\ast}\le
\mathcal{O}_{10} \kappa^{-3} e^{\kappa \zeta^2} \bar{g}^{2-2\eta}
\]
where $\mathcal{O}_{10}=2\times\mathcal{O}_9$.
\end{Lemma}

\noindent{\bf Proof:}
We have
\begin{eqnarray*}
|Q_\Delta(\phi_1+\zeta) e^{-V_\Delta(\phi_1+\zeta)}|_{\partial\phi,h_\ast}
 & \le & |Q_\Delta(\phi_1+\zeta)|_{\partial\phi,h_\ast}
\times |e^{-V_\Delta(\phi_1+\zeta)}|_{\partial\phi,h_\ast} \\
 & \le & \mathcal{O}_8 \kappa^{-3} e^{\kappa \zeta^2} \bar{g}^{2-2\eta}\times
 |e^{-V_\Delta(\phi_1+\zeta)}|_{\partial\phi,h_\ast} 
\end{eqnarray*} 
by Lemma \ref{L18lem}. The last factor is bounded in a coarse way using
\[
|e^{-V_\Delta(\phi_1+\zeta)}|_{\partial\phi,h_\ast} = ||e^{-V_\Delta(\phi_1+\zeta)}||_{\partial\phi,0,h_\ast} 
\]
\[
\le ||e^{-V_\Delta(\psi)}||_{\partial\psi,\zeta,h_\ast}\le
 ||e^{-V_\Delta(\psi)}||_{\partial\psi,\zeta,h}
\]
since $h_\ast\le h$ by (\ref{c1c2stdeq}).
We finally use Lemma \ref{L7lem} with $t=1$ in order to write
\[
|e^{-V_\Delta(\phi_1+\zeta)}|_{\partial\phi,h_\ast} \le 2 e^{-\frac{1}{2}(\Re\beta_{4,\Delta})\zeta^4}\le 2
\]
and the result follows.
\qed

\begin{Lemma}\label{L21lem}
For all $\kappa\in(0,1]$ and all $\phi,\zeta\in\mathbb{R}$
we have
\[
\left|\left|Q_\Delta(\phi_1+\zeta)
\left(e^{-V_\Delta(\phi_1+\zeta)}-e^{-\tilde{V}_\Delta(\phi_1)}
\right)\right|\right|_{\partial\phi,\phi,h}\le
\mathcal{O}_{11} \kappa^{-5} e^{\kappa\zeta^2} \bar{g}^{\frac{3}{4}-3\eta}
\]
with
\[
\mathcal{O}_{11}=
2^6\mathcal{O}_1\times\mathcal{O}_9+2^{\frac{43}{4}} \mathcal{O}_1\times\mathcal{O}_7\ .
\]
\end{Lemma}

\noindent{\bf Proof:}
Define
\[
j(s)=e^{-(1-s)V_\Delta(\phi_1+\zeta)-s\tilde{V}_\Delta(\phi_1)}=e^{-V_\Delta(\phi_1+\zeta)+s p_\Delta}
\]
for $s\in [0,1]$.
Thus
\begin{eqnarray*}
Q_\Delta(\phi_1+\zeta)
\left(e^{-V_\Delta(\phi_1+\zeta)}-e^{-\tilde{V}_\Delta(\phi_1)}
\right) & = & Q_\Delta(\phi_1+\zeta) (j(0)-j(1)) \\
 & = & -Q_\Delta(\phi_1+\zeta) \int_0^1 {\rm d}s\ j'(s) \\
 & = & -Q_\Delta(\phi_1+\zeta) \int_0^1 {\rm d}s\ 
 p_\Delta(\phi_1,\zeta) e^{-(1-s)V_\Delta(\phi_1+\zeta)-s\tilde{V}_\Delta(\phi_1)} \\
  & = & A+B
\end{eqnarray*}
with
\[
A= -Q_\Delta(\phi_1+\zeta) \int_0^{\frac{1}{2}} {\rm d}s\ 
 p_\Delta(\phi_1,\zeta) e^{-(1-s)V_\Delta(\phi_1+\zeta)-s\tilde{V}_\Delta(\phi_1)} 
\]
and
\[
B= -Q_\Delta(\phi_1+\zeta) \int_{\frac{1}{2}}^1 {\rm d}s\ 
 p_\Delta(\phi_1,\zeta) e^{-(1-s)V_\Delta(\phi_1+\zeta)-s\tilde{V}_\Delta(\phi_1)} \ .
\]
We bound these two pieces separately.
Note that
\[
A= -Q_\Delta(\phi_1+\zeta) e^{-\frac{1}{2}V_\Delta(\phi_1+\zeta)}
\int_0^{\frac{1}{2}} {\rm d}s\ 
 p_\Delta(\phi_1,\zeta) e^{-(\frac{1}{2}-s)V_\Delta(\phi_1+\zeta)-s\tilde{V}_\Delta(\phi_1)} 
\]
which implies the estimate
\[
||A||_{\partial\phi,\phi,h}\le
||Q_\Delta(\phi_1+\zeta) e^{-\frac{1}{2}V_\Delta(\phi_1+\zeta)}||_{\partial\phi,\phi,h}\times
\int_0^{\frac{1}{2}} {\rm d}s\ 
||p_\Delta(\phi_1,\zeta)||_{\partial\phi,\phi,h}
\]
\[
\times ||e^{-(\frac{1}{2}-s)V_\Delta(\phi_1+\zeta)}||_{\partial\phi,\phi,h}\times
||e^{s\tilde{V}_\Delta(\phi_1)}||_{\partial\phi,\phi,h}\ .
\]
Repeating the proof of Lemma \ref{L19lem}, but this time taking $\gamma=\frac{1}{4}$ instead of $\frac{1}{2}$
which results in an extra factor $\left(\frac{1}{2}\right)^{-\frac{6}{4}}=2^{\frac{3}{2}}$,
we obtain
\[
||Q_\Delta(\phi_1+\zeta) e^{-\frac{1}{2}V_\Delta(\phi_1+\zeta)}||_{\partial\phi,\phi,h}
\le \mathcal{O}_9\times 2^{\frac{3}{2}}\times\bar{g}^{\frac{1}{2}-2\eta}\ .
\]
For $0<s<\frac{1}{2}$, we get from Lemma \ref{L7lem} with $t=\frac{1}{2}-s$
\[
||e^{-(\frac{1}{2}-s)V_\Delta(\phi_1+\zeta)}||_{\partial\phi,\phi,h}
\le ||e^{-(\frac{1}{2}-s)V_\Delta(\psi)}||_{\partial\psi,\phi_1+\zeta,h}
\le 2 e^{-\frac{1}{2}\left(\frac{1}{2}-s\right)(\Re\beta_{4,\Delta})(\phi_1+\zeta)^4}\le 2\ .
\]
The above steps result in the bound
\[
||A||_{\partial\phi,\phi,h}\le
2^{\frac{5}{2}}\mathcal{O}_9\bar{g}^{\frac{1}{2}-2\eta}\times
\int_0^{\frac{1}{2}} {\rm d}s\ 
||p_\Delta(\phi_1,\zeta)||_{\partial\phi,\phi,h}
||e^{s\tilde{V}_\Delta(\phi_1)}||_{\partial\phi,\phi,h}\ .
\]
By Lemma \ref{L8lem}
\[
||e^{s\tilde{V}_\Delta(\phi_1)}||_{\partial\phi,\phi,h}\le
||e^{s\tilde{V}_\Delta(\psi)}||_{\partial\psi,\phi_1,h}\le 2 e^{-\frac{s}{2}(\Re\beta_{4,\Delta})\phi_1^4}\ .
\]
Now we use Lemma \ref{L11lem} with $\gamma=\frac{s}{2}$ and with the present $\kappa$
in order to derive
\[
||p_\Delta(\phi_1,\zeta)||_{\partial\phi,\phi,h}\le\mathcal{O}_1
\kappa^{-2} \left(\frac{s}{2}\right)^{-\frac{3}{4}}\bar{g}^{\frac{1}{4}-\eta}
e^{\kappa\zeta^2} e^{\frac{s}{2}(\Re\beta_{4,\Delta})\phi_1^4}\ .
\]
This produces the bound
\[
||A||_{\partial\phi,\phi,h}\le 2^{\frac{5}{2}}
\mathcal{O}_9\times 2\times \mathcal{O}_1\times 2^{\frac{3}{4}}\times
\bar{g}^{\frac{3}{4}-3\eta}\times \kappa^{-2}e^{\kappa\zeta^2}\times
\int_0^{\frac{1}{2}} {\rm d}s\ s^{-\frac{3}{4}}\ ,
\]
namely,
\[
||A||_{\partial\phi,\phi,h}\le
2^6\times \mathcal{O}_1\times\mathcal{O}_9\times \bar{g}^{\frac{3}{4}-3\eta}\times \kappa^{-2}e^{\kappa\zeta^2}\ .
\]
We now take care of $B$. From the definition we readily obtain
\[
||B||_{\partial\phi,\phi,h}\le
||Q_\Delta(\phi_1+\zeta)||_{\partial\phi,\phi,h}
\int_{\frac{1}{2}}^1 {\rm d}s\ 
||p_\Delta(\phi_1,\zeta)||_{\partial\phi,\phi,h}
|| e^{-(1-s)V_\Delta(\phi_1+\zeta)}||_{\partial\phi,\phi,h}
||e^{-s\tilde{V}_\Delta(\phi_1)}||_{\partial\phi,\phi,h} \ .
\]
Since $\kappa\in(0,1]$, we have $\frac{\kappa}{2}\in(0,2]$
and therefore Lemma \ref{L17lem} with $\frac{\kappa}{2}$ instead of $\kappa$ and with $\gamma=\frac{1}{8}$
gives us the estimate
\[
||Q_\Delta(\phi_1+\zeta)||_{\partial\phi,\phi,h}\le
\mathcal{O}_7 \bar{g}^{\frac{1}{2}-2\eta}\times 8^{\frac{3}{2}} e^{\frac{1}{8}(\Re\beta_{4,\Delta})\phi_1^4}
\times 8 \kappa^{-3}e^{\frac{\kappa}{2}\zeta^2}\ .
\]
We also use Lemma \ref{L11lem} with $\frac{\kappa}{2}$ instead of $\kappa$ and with $\gamma=\frac{1}{8}$
and get
\[
||p_\Delta(\phi_1,\zeta)||_{\partial\phi,\phi,h}
\le \mathcal{O}_1\times 4\times \kappa^{-2}\times 8^{\frac{3}{4}}\times 
\bar{g}^{\frac{1}{4}-\eta} e^{\frac{\kappa}{2}\zeta^2}\times e^{\frac{1}{8}(\Re\beta_{4,\Delta})\phi_1^4}\ .
\]
From Lemma \ref{L7lem} with $t=1-s$ we obtain
\[
|| e^{-(1-s)V_\Delta(\phi_1+\zeta)}||_{\partial\phi,\phi,h}
\le || e^{-(1-s)V_\Delta(\psi)}||_{\partial\psi,\phi_1+\zeta,h}
\le 2 e^{-\frac{(1-s)}{2}(\Re\beta_{4,\Delta})(\phi_1+\zeta)^4}\le 2\ .
\]
Finally, the last ingredient is the use of Lemma \ref{L8lem} with $t=s$
which results in
\[
||e^{-s\tilde{V}_\Delta(\phi_1)}||_{\partial\phi,\phi,h}\le
||e^{-s\tilde{V}_\Delta(\psi)}||_{\partial\psi,\phi_1,h}
\le e^{-\frac{s}{2}(\Re\beta_{4,\Delta})\phi_1^4}\le 2
e^{-\frac{1}{4}(\Re\beta_{4,\Delta})\phi_1^4}
\]
since $s\ge\frac{1}{2}$.
Altogether the previous bounds imply
\begin{eqnarray*}
||B||_{\partial\phi,\phi,h} & \le & \int_{\frac{1}{2}}^1{\rm d}s\ 
\mathcal{O}_7\times \bar{g}^{\frac{1}{2}-2\eta}\times 8^{\frac{3}{2}}\times 8\times\kappa^{-3} \\
 & & \times \mathcal{O}_1 \times 4\times \kappa^{-2}\times 8^{\frac{3}{4}}\times  \bar{g}^{\frac{1}{4}-\eta}
 \times 2\times e^{\kappa\zeta^2} \\
 & \le & 2^{\frac{43}{4}} \mathcal{O}_1\times \mathcal{O}_7 
 \kappa^{-5}\bar{g}^{\frac{3}{4}-3\eta} e^{\kappa\zeta^2}\ .
\end{eqnarray*}
Combining the bounds for $A$ and $B$ we obtain the desired estimate.
\qed

\begin{Lemma}\label{L22lem}
If $0<\kappa\le 2^{-\frac{3}{2}} L^{-(3-2[\phi])}$ then for all $\zeta\in\mathbb{R}$ we have
\[
\left|Q_\Delta(\phi_1+\zeta)
\left(e^{-V_\Delta(\phi_1+\zeta)}-e^{-\tilde{V}_\Delta(\phi_1)}
\right)\right|_{\partial\phi,h_\ast}\le
\mathcal{O}_{12} \kappa^{-5} e^{\kappa\zeta^2} \bar{g}^{3-3\eta}
\]
where $\mathcal{O}_{12}=2^{7}\mathcal{O}_2\times \mathcal{O}_8 $.
\end{Lemma}

\noindent{\bf Proof:}
The proof is simpler than that of the previous lemma because we do not need to split the quantity at hand.
We directly bound the latter, namely,
\[
Q_\Delta(\phi_1+\zeta)
\left(e^{-V_\Delta(\phi_1+\zeta)}-e^{-\tilde{V}_\Delta(\phi_1)}
\right)=
-Q_\Delta(\phi_1+\zeta) \int_0^1 {\rm d}s\ 
 p_\Delta(\phi_1,\zeta) e^{-(1-s)V_\Delta(\phi_1+\zeta)-s\tilde{V}_\Delta(\phi_1)}
\]
by
\[
\left|Q_\Delta(\phi_1+\zeta)
\left(e^{-V_\Delta(\phi_1+\zeta)}-e^{-\tilde{V}_\Delta(\phi_1)}
\right)\right|_{\partial\phi,h_\ast}\le
|Q_\Delta(\phi_1+\zeta)|_{\partial\phi,h_\ast}\times
\int_0^1 {\rm d}s\ 
 |p_\Delta(\phi_1,\zeta)|_{\partial\phi,h_\ast}
\]
\[
\times  |e^{-(1-s)V_\Delta(\phi_1+\zeta)}|_{\partial\phi,h_\ast}
\times
|e^{-s\tilde{V}_\Delta(\phi_1)}|_{\partial\phi,h_\ast}\ .
\]
From Lemma \ref{L18lem} with $\frac{\kappa}{2}$ instead of $\kappa$ we have
\[
|Q_\Delta(\phi_1+\zeta)|_{\partial\phi,h_\ast}\le 8\times\mathcal{O}_8
\kappa^{-3} e^{\frac{\kappa}{2}\zeta^2} \bar{g}^{2-2\eta}\ .
\]
Likewise, from Lemma \ref{L12lem} with $\frac{\kappa}{2}$ instead of $\kappa$ we get
\[
|p_\Delta(\phi_1,\zeta)|_{\partial\phi,h_\ast}\le 4\times\mathcal{O}_2
\kappa^{-2} e^{\frac{\kappa}{2}\zeta^2} \bar{g}^{1-\eta}\ .
\]
We also have
\[
|e^{-(1-s)V_\Delta(\phi_1+\zeta)}|_{\partial\phi,h_\ast}=||e^{-(1-s)V_\Delta(\phi_1+\zeta)}||_{\partial\phi,0,h_\ast}
\le ||e^{-(1-s)V_\Delta(\psi)}|_{\partial\psi,\zeta,h_\ast}
\]
\[
\le ||e^{-(1-s)V_\Delta(\psi)}|_{\partial\psi,\zeta,h}\le 2 e^{-\frac{(1-s)}{2}(\Re\beta_{4,\Delta})\zeta^4}\le 2
\]
by Lemma \ref{L7lem}.
Finally Lemma \ref{L8lem} provides the estimate
\[
|e^{-s\tilde{V}_\Delta(\phi_1)}|_{\partial\phi,h_\ast}\le |e^{-s\tilde{V}_\Delta(\psi)}|_{\partial\psi,h_\ast}\le 2\ .
\]
Altogether the previous bounds imply
\[
\left|Q_\Delta(\phi_1+\zeta)
\left(e^{-V_\Delta(\phi_1+\zeta)}-e^{-\tilde{V}_\Delta(\phi_1)}
\right)\right|_{\partial\phi,h_\ast}\le
8\mathcal{O}_8 \kappa^{-5} e^{\kappa\zeta^2}\bar{g}^{3-3\eta}\times 16\times\mathcal{O}_2
\]
from which the result follows.
\qed

\begin{Lemma}\label{L23lem}
For all $K\in C_{{\rm bd}}^{9}(\mathbb{R},\mathbb{C})$ and for all $\sigma\in\mathbb{R}$ we have
\[
||K(\psi)||_{\partial\psi,\sigma,h_\ast}\le \mathcal{O}_{13}
e^{h_{\ast}^{-2}\sigma^2}\times\left[
|K(\psi)|_{\partial\psi,h_\ast}+h_{\ast}^{9}h^{-9}\sup_{\psi\in\mathbb{R}}
||K(\psi)||_{\partial\psi,\psi,h_\ast}
\right]
\]
where
\[
\mathcal{O}_{13}=1+511\times\max_{0\le j\le 9}\left(\frac{j}{2e}\right)^{\frac{j}{2}}\ .
\]
\end{Lemma}

\noindent{\bf Proof:}
Recall that by definition
\[
||K(\psi)||_{\partial\psi,\sigma,h_\ast}=\sum_{n=0}^9\frac{h_{\ast}^n}{n!}|K^{(n)}(\sigma)|\ .
\]
The term with $n=9$ is bounded by writing
\[
\frac{h_{\ast}^9}{9!}|K^{(9)}(\sigma)|=h_{\ast}^{9}h^{-9}\times
\frac{h^9}{9!}|K^{(9)}(\sigma)|\le h_{\ast}^{9}h^{-9}\times
\sup_{\psi\in\mathbb{R}}
||K(\psi)||_{\partial\psi,\psi,h_\ast}\ .
\]
For terms with $0\le n\le 8$ we use a Taylor expansion around zero of order
$8-n$ so that the integral remainder involves $(9-n)$-th derivatives of $K^{(n)}$, i.e., $9$-th derivatives
of the original function $K$.
Indeed, one can write
\[
K^{(n)}(\sigma)=\sum_{m=0}^{8-n} \frac{\sigma^m}{m!} K^{(n+m)}(0)
+\frac{1}{(8-n)!}\int_0^1 (1-s)^{8-n} \sigma^{9-n} K^{(9)}(s\sigma)\ {\rm d}s
\]
and therefore
\[
|K^{(n)}(\sigma)|\le \sum_{m=0}^{8-n} \frac{|\sigma|^m}{m!}
(n+m)!\ h_{\ast}^{-(n+m)} |K(\psi)|_{\partial\psi,h_\ast}
\]
\[
+\frac{1}{(8-n)!} |\sigma|^{9-n}\times 9!\ h^{-9}\left(
\sup_{\psi\in\mathbb{R}} ||K(\psi)||_{\partial\psi,\psi,h}
\right)
\int_0^1 (1-s)^{8-n}\ {\rm d}s\ .
\]
We use Lemma \ref{L5lem} with $\kappa=h_{\ast}^{-2}$ in order to bound powers of $|\sigma|$
by
\[
|\sigma|^m\le \left(\frac{m}{2e}\right)^{\frac{m}{2}}\times
h_{\ast}^{m} e^{h_{\ast}^{-2}\sigma^2} 
\]
which inserted in the previous inequality gives
\[
\frac{h_{\ast}^{n}}{n!}|K^{(n)}(\sigma)|\le \left(
\max_{0\le j\le 9} \left(\frac{j}{2e}\right)^{\frac{j}{2}}
\right)\times e^{h_{\ast}^{-2}\sigma^2} 
\]
\[
\times\left[
\sum_{m=0}^{8-n} \frac{(n+m)!}{n!m!} |K(\psi)|_{\partial\psi,h_\ast}
+ \frac{9!\ h_{\ast}^{9-n}}{n!(9-n)!} h_{\ast}^{n}h^{-9} 
\sup_{\psi\in\mathbb{R}} ||K(\psi)||_{\partial\psi,\psi,h}
\right]\ .
\]
Putting together the bounds for the different values of $n$ we obtain
\[
||K(\psi)||_{\partial\psi,\sigma,h_\ast}\le h_{\ast}^{9}h^{-9}  \sup_{\psi\in\mathbb{R}} ||K(\psi)||_{\partial\psi,\psi,h}
\]
\[
+e^{h_{\ast}^{-2}\sigma^2} \left(
\max_{0\le j\le 9} \left(\frac{j}{2e}\right)^{\frac{j}{2}}
\right)\sum_{n=0}^{8} \left[
\left(\sum_{m=0}^{8-n} \frac{(n+m)!}{n!m!}\right) |K(\psi)|_{\partial\psi,h_\ast}
+\frac{9!}{n!(9-n!)} h_{\ast}^{9}h^{-9}  \sup_{\psi\in\mathbb{R}} ||K(\psi)||_{\partial\psi,\psi,h}
\right]\ .
\]
The result as well as the given value of $\mathcal{O}_{13}$ then follow since
\[
\sum_{n=0}^{8}\sum_{m=0}^{8-n} \frac{(n+m)!}{n!m!} =\sum_{n=0}^{8} \frac{9!}{n!(9-n!)}=2^9-1=511.
\]
\qed

\begin{Lemma}\label{L24lem}
For all $K\in C_{{\rm bd}}^{9}(\mathbb{R},\mathbb{C})$, $\beta_4\in\mathbb{C}$ such that
$|\beta_4-\bar{g}|<\frac{1}{2}\bar{g}$,
$\gamma\in (0,1]$ and $\phi\in\mathbb{R}$ we have
\[
||K(\phi)||_{\partial\phi,\phi,h}\le \mathcal{O}_{14}
\gamma^{-\frac{9}{4}} e^{\gamma(\Re\beta_{4})\phi^4}
\left[
|K(\psi)|_{\partial\psi,h}+L^{-9[\phi]} \sup_{\psi\in\mathbb{R}} ||K(\psi)||_{\partial\psi,\psi,L^{[\phi]}h}
\right]
\]
with
\[
\mathcal{O}_{14}=1+\left((1+c_{1}^{-1})^9-1\right)\times
\max_{0\le j\le 9} \left(\frac{j}{2e}\right)^{\frac{j}{4}}\ .
\] 
\end{Lemma}

\noindent{\bf Proof:}
We proceed as in the proof of the previous lemma and
write
\[
\frac{h^9}{9!}|K^{(9)}(\phi)|=L^{-9[\phi]}\times \frac{(L^{[\phi]}h)^9}{9!} |K^{(9)}(\phi)|
\le L^{-9[\phi]}\times \sup_{\psi\in\mathbb{R}} ||K(\psi)||_{\partial\psi,\psi,L^{[\phi]}h}
\]
in order to handle the $n=9$ term in the sum defining $||K(\phi)||_{\partial\phi,\phi,h}$.
For the other terms with $0\le n\le 8$ one has, as before,
\[
|K^{(n)}(\phi)|\le \sum_{m=0}^{8-n} \frac{|\phi|^m}{m!}
(n+m)!\ h^{-(n+m)} |K(\psi)|_{\partial\psi,h}
\]
\[
+\frac{1}{(9-n)!} |\phi|^{9-n}\times 9! (L^{[\phi]}h)^{-9}
\sup_{\psi\in\mathbb{R}} ||K(\psi)||_{\partial\psi,\psi,L^{[\phi]}h}\ .
\]
We this time use Lemma \ref{L4lem} in order to bound powers of $|\phi|$ by
\[
|\phi|^m\le \left(\frac{m}{2e}\right)^{\frac{m}{4}} \gamma^{-\frac{m}{4}}\bar{g}^{-\frac{m}{4}}
e^{\gamma(\Re\beta_{4})\phi^4}\ .
\]
Note that $\gamma^{-\frac{m}{4}}\le \gamma^{-\frac{9}{4}}$ since $0<\gamma\le 1$, $0\le n\le 8$ and $0\le m\le 9-n$.
Besides $\bar{g}^{-\frac{m}{4}}=(c_{1}^{-1} h)^{m}$
and therefore
\[
\frac{h^n}{n!}|K^{(n)}(\phi)|\le
\left(\max_{0\le j\le 9} \left(\frac{j}{2e}\right)^{\frac{j}{4}}\right)\times \gamma^{-\frac{9}{4}}
\times e^{\gamma(\Re\beta_{4})\phi^4}
\]
\[
\times \left[
\sum_{m=0}^{8-n} \frac{h^m c_1^{-m}}{m!} \frac{h^n}{n!} (n+m)! h^{-(n+m)} |K(\psi)|_{\partial\psi,h}
+\frac{9!}{n!(9-n)!} h^n h^{9-n} c_1^{-(9-n)} (L^{[\phi]}h)^{-9}\sup_{\psi\in\mathbb{R}} ||K(\psi)||_{\partial\psi,\psi,L^{[\phi]}h}
\right]\ .
\]
Altogether this gives the estimate
\[
||K(\phi)||_{\partial\phi,\phi,h}\le L^{-9[\phi]}  \sup_{\psi\in\mathbb{R}} ||K(\psi)||_{\partial\psi,\psi,L^{[\phi]}h}
\]
\[
+\left(\max_{0\le j\le 9} \left(\frac{j}{2e}\right)^{\frac{j}{4}}\right)\times \gamma^{-\frac{9}{4}}
\times e^{\gamma(\Re\beta_{4})\phi^4}\times
\left\{
\left(
\sum_{m=0}^{8-n} \left(\begin{array}{c}
n+m \\ m \end{array}\right) c_1^{-m}
\right)|K(\psi)|_{\partial\psi,h}
\right.
\]
\[
\left.
+\left(\begin{array}{c}
9 \\ n \end{array}\right) c_1^{-(9-n)} L^{-9[\phi]}  \sup_{\psi\in\mathbb{R}} ||K(\psi)||_{\partial\psi,\psi,L^{[\phi]}h}
\right\}\ .
\]
The result with the given value for $\mathcal{O}_{14}$ follows from this last inequality since
\[
\sum_{n=0}^8 \sum_{m=0}^{8-n} \left(\begin{array}{c}
n+m \\ m \end{array}\right) c_1^{-m}
=c_1\left[(1+c_1^{-1})^9-1\right]<(1+c_1^{-1})^9-1=
\sum_{n=0}^8 \left(\begin{array}{c}
9 \\ n \end{array}\right) c_1^{-(9-n)}\ .
\]
\qed

\begin{Lemma}\label{L25lem}
Let $\kappa\in(0,1]$.
\begin{enumerate}
\item
If $|\lambda|\bar{g}^{\frac{1}{4}-\frac{1}{3}\eta_R}\le 1$ then $\forall \phi,\zeta\in\mathbb{R}$,
\[
||K_\Delta(\lambda,\phi_1,\zeta)||_{\partial\phi,\phi,h}
\le \mathcal{O}_{15}\kappa^{-5} e^{\kappa\zeta^2}\left(|\lambda|\bar{g}^{\frac{1}{4}-\frac{1}{3}\eta_R}\right)^2
\]
where $\mathcal{O}_{15}=2^{\frac{5}{2}}\mathcal{O}_7+\mathcal{O}_{11}+1$.
\item
If $R=0$ and $|\lambda|\bar{g}^{\frac{1}{4}-\eta}\le 1$ then we have the improvement $\forall \phi,\zeta\in\mathbb{R}$,
\[
||K_\Delta(\lambda,\phi_1,\zeta)||_{\partial\phi,\phi,h}
\le \mathcal{O}_{16}\kappa^{-5} e^{\kappa\zeta^2}\left(|\lambda|\bar{g}^{\frac{1}{4}-\eta}\right)^2
\]
where $\mathcal{O}_{16}=2^{\frac{5}{2}}\mathcal{O}_7+\mathcal{O}_{11}$.
\end{enumerate}
\end{Lemma}

\noindent{\bf Proof:}
By definition
\[
K_\Delta(\lambda,\phi_1,\zeta)=\lambda^2 Q_\Delta(\phi_1+\zeta) e^{-\tilde{V}_\Delta(\phi_1)}
+\lambda^3 Q_\Delta(\phi_1+\zeta)\left(e^{-V_\Delta(\phi_1+\zeta)}- e^{-\tilde{V}_\Delta(\phi_1)}\right)
+\lambda^3 R_\Delta(\phi_1+\zeta)\ .
\]
Thus
\[
||K_\Delta(\lambda,\phi_1,\zeta)||_{\partial\phi,\phi,h}
\le |\lambda|^2 ||Q_\Delta(\phi_1+\zeta)||_{\partial\phi,\phi,h} ||e^{-\tilde{V}_\Delta(\phi_1)}||_{\partial\phi,\phi,h}
\]
\[
+|\lambda|^3 \left|\left|Q_\Delta(\phi_1+\zeta)\left(e^{-V_\Delta(\phi_1+\zeta)}
- e^{-\tilde{V}_\Delta(\phi_1)}\right)\right|\right|_{\partial\phi,\phi,h}
+|\lambda|^3 ||R_\Delta(\phi_1+\zeta)||_{\partial\phi,\phi,h}\ .
\]
From Lemma \ref{L8lem} we have
\[
||e^{-\tilde{V}_\Delta(\phi_1)}||_{\partial\phi,\phi,h}\le ||e^{-\tilde{V}_\Delta(\psi)}||_{\partial\psi,\phi_1,h}
\le 2 e^{-\frac{1}{2}(\Re\beta_{4,\Delta})\phi_1^4}\ .
\]
We use Lemma \ref{L17lem} with $\gamma=\frac{1}{2}$ and get
\[
||Q_\Delta(\phi_1+\zeta)||_{\partial\phi,\phi,h}\le \mathcal{O}_7
\bar{g}^{\frac{1}{2}-2\eta}\times 2^{\frac{3}{2}}\times\kappa^{-3} 
e^{\kappa\zeta^2} e^{\frac{1}{2}(\Re\beta_{4,\Delta})\phi_1^4}\ .
\]
As a result
\[
|\lambda|^2 ||Q_\Delta(\phi_1+\zeta)||_{\partial\phi,\phi,h} ||e^{-\tilde{V}_\Delta(\phi_1)}||_{\partial\phi,\phi,h}
\le
2^{\frac{5}{2}}  \mathcal{O}_7 \bar{g}^{\frac{1}{2}-2\eta}
\kappa^{-3} e^{\kappa\zeta^2} |\lambda|^2\ .
\]
From Lemma \ref{L21lem} we get
\[
|\lambda|^3 \left|\left|Q_\Delta(\phi_1+\zeta)\left(e^{-V_\Delta(\phi_1+\zeta)}
- e^{-\tilde{V}_\Delta(\phi_1)}\right)\right|\right|_{\partial\phi,\phi,h}
\le
\mathcal{O}_{11} \bar{g}^{\frac{3}{4}-3\eta}
\kappa^{-5} e^{\kappa\zeta^2} |\lambda|^3\ .
\]
Finally the last term is bounded using
\[
||R_\Delta(\phi_1+\zeta)||_{\partial\phi,\phi,h}\le
||R_\Delta(\psi)||_{\partial\psi,\phi_1+\zeta,h}\le \sup_{\psi\in\mathbb{R}}  ||R_\Delta(\psi)||_{\partial\psi,\psi,h}
\]
\[
\le \bar{g}^{-2} |||R_\Delta|||_{\bar{g}}
\le \bar{g}^{-2}\times \bar{g}^{\frac{11}{4}-\eta_R}=\bar{g}^{\frac{3}{4}-\eta_R}
\]
from (\ref{normstdeq}).
Collecting the previous estimates we arrive at
\[
||K_\Delta(\lambda,\phi_1,\zeta)||_{\partial\phi,\phi,h}
\le \kappa^{-5} e^{\kappa\zeta^2} 
\left[
2^{\frac{5}{2}}\mathcal{O}_{7} |\lambda|^2 \bar{g}^{\frac{1}{2}-2\eta}
+ \mathcal{O}_{11} |\lambda|^3 \bar{g}^{\frac{3}{4}-3\eta}
+  |\lambda|^3 \bar{g}^{\frac{3}{4}-\eta_R}
\right]\ .
\]
By the standard hypothesis (\ref{etarstdeq}), $\eta_R\ge 3\eta$ and since $0<\bar{g}\le 1$ we have $\bar{g}^{\frac{1}{4}-\eta}\le
\bar{g}^{\frac{1}{4}-\frac{1}{3}\eta_R}$ and therefore
\[
||K_\Delta(\lambda,\phi_1,\zeta)||_{\partial\phi,\phi,h}
\le \kappa^{-5} e^{\kappa\zeta^2} 
\left[
2^{\frac{5}{2}}\mathcal{O}_{7} \left(|\lambda|\bar{g}^{\frac{1}{4}-\frac{1}{3}\eta_R}\right)^2
+ \mathcal{O}_{11} \left(|\lambda|\bar{g}^{\frac{1}{4}-\frac{1}{3}\eta_R}\right)^3
+  \left(|\lambda|\bar{g}^{\frac{1}{4}-\frac{1}{3}\eta_R}\right)^3
\right]\ .
\]
from which part 1) follows.
As for part 2), the $R$ term being absent from the start, the bound on $K$ reduces to
\[
||K_\Delta(\lambda,\phi_1,\zeta)||_{\partial\phi,\phi,h}
\le \kappa^{-5} e^{\kappa\zeta^2} 
\left[
2^{\frac{5}{2}}\mathcal{O}_{7} |\lambda|^2 \bar{g}^{\frac{1}{2}-2\eta}
+ \mathcal{O}_{11} |\lambda|^3 \bar{g}^{\frac{3}{4}-3\eta}
\right]\ .
\]
which immediately yealds the desired result.
\qed

\begin{Lemma}\label{L26lem}
\begin{enumerate}
\item
If $|\lambda|\bar{g}^{\frac{11}{12}-\frac{1}{3}\eta_R}\le 1$ then $\forall \zeta\in\mathbb{R}$,
\[
|K_\Delta(\lambda,\phi_1,\zeta)|_{\partial\phi,h_\ast}
\le \mathcal{O}_{17} h_{\ast}^{10} e^{h_{\ast}^{-2}\zeta^2}
\left(|\lambda|\bar{g}^{\frac{11}{12}-\frac{1}{3}\eta_R}\right)^2
\]
where $\mathcal{O}_{17}=2\mathcal{O}_8+\mathcal{O}_{12}+2\mathcal{O}_{13}$.
\item
If $R=0$ and $|\lambda|\bar{g}^{1-\eta}\le 1$ then we have the improvement $\forall \zeta\in\mathbb{R}$,
\[
|K_\Delta(\lambda,\phi_1,\zeta)|_{\partial\phi,h_\ast}
\le \mathcal{O}_{18} h_{\ast}^{10} e^{h_{\ast}^{-2}\zeta^2}
\left(|\lambda|\bar{g}^{1-\eta}\right)^2
\]
where $\mathcal{O}_{18}=2\mathcal{O}_8+\mathcal{O}_{12}$.
\end{enumerate}
\end{Lemma}

\noindent{\bf Proof:}
As before we start with
\[
|K_\Delta(\lambda,\phi_1,\zeta)|_{\partial\phi,h_\ast}
\le 
|\lambda|^2 |Q_\Delta(\phi_1+\zeta)|_{\partial\phi,h_\ast} |e^{-\tilde{V}_\Delta(\phi_1)}|_{\partial\phi,h_\ast}
\]
\[
+|\lambda|^3 \left|Q_\Delta(\phi_1+\zeta)\left(e^{-V_\Delta(\phi_1+\zeta)}- e^{-\tilde{V}_\Delta(\phi_1)}\right)
\right|_{\partial\phi,h_\ast}
+|\lambda|^3 |R_\Delta(\phi_1+\zeta)|_{\partial\phi,h_\ast}\ .
\]
Then by Lemma \ref{L18lem} with $\kappa=h_{\ast}^{-2}$
\[
|Q_\Delta(\phi_1+\zeta)|_{\partial\phi,h_\ast} \le \mathcal{O}_8 h_{\ast}^{6} e^{h_{\ast}^{-2}\zeta^2}\bar{g}^{2-2\eta}\ .
\]
From Lemma \ref{L8lem} we have
\[
|e^{-\tilde{V}_\Delta(\phi_1)}|_{\partial\phi,h_\ast}=||e^{-\tilde{V}_\Delta(\phi_1)}||_{\partial\phi,0,h_\ast}
\le ||e^{-\tilde{V}_\Delta(\psi)}||_{\partial\psi,0,h_\ast}=|e^{-\tilde{V}_\Delta(\psi)}|_{\partial\psi,h_\ast}\le 2\ .
\]
By Lemma \ref{L22lem} with $\kappa=h_{\ast}^{-2}$ we have
\[
\left|Q_\Delta(\phi_1+\zeta)\left(e^{-V_\Delta(\phi_1+\zeta)}- e^{-\tilde{V}_\Delta(\phi_1)}\right)
\right|_{\partial\phi,h_\ast}
\le \mathcal{O}_{12} h_{\ast}^{10} e^{h_{\ast}^{-2}\zeta^2}\bar{g}^{3-3\eta}\ .
\]
As a result of the estimates we have so far
\[
|K_\Delta(\lambda,\phi_1,\zeta)|_{\partial\phi,h_\ast}
\le 
2|\lambda|^2 \mathcal{O}_8 h_{\ast}^{6} e^{h_{\ast}^{-2}\zeta^2}\bar{g}^{2-2\eta}
+|\lambda|^3 \mathcal{O}_{12} h_{\ast}^{10} e^{h_{\ast}^{-2}\zeta^2}\bar{g}^{3-3\eta}
+|\lambda|^3 |R_\Delta(\phi_1+\zeta)|_{\partial\phi,h_\ast}\ .
\]
The last term will be estimated as follows.
Note that
\[
|R_\Delta(\phi_1+\zeta)|_{\partial\phi,h_\ast}=||R_\Delta(\phi_1+\zeta)||_{\partial\phi,0,h_\ast}
\le ||R_\Delta(\psi+\zeta)|_{\partial\psi,0,h_\ast}=||R_\Delta(\psi)||_{\partial\psi,\zeta,h_\ast}
\]
\[
\le \mathcal{O}_{13}  e^{h_{\ast}^{-2}\zeta^2}
\left[|R_\Delta(\psi)|_{\partial\psi,h_\ast}+h_{\ast}^{-9} h^{9}
\sup_{\psi\in\mathbb{R}} ||R_\Delta(\psi)||_{\partial\psi,\psi,h} \right]
\]
by Lemma \ref{L23lem}.
Hence
\[
|R_\Delta(\phi_1+\zeta)|_{\partial\phi,h_\ast}\le \mathcal{O}_{13}  e^{h_{\ast}^{-2}\zeta^2}
|||R_\Delta|||_{\bar{g}}\left(1+\bar{g}^{-2} h_{\ast}^{-9} h^{9}\right)\ .
\]
Now
\[
\bar{g}^{-2} h_{\ast}^{-9} h^{9}=c_2^{9} c_1^{-9} L^{\frac{9}{4}(3+\epsilon)} \bar{g}^{\frac{1}{4}}
\le c_2^{9} c_1^{-9} L^{9} \bar{g}^{\frac{1}{4}}\le 1
\]
by the standard hypothesis (\ref{c1c2stdeq}).
Also using (\ref{normstdeq}) we now arrive at
\[
|K_\Delta(\lambda,\phi_1,\zeta)|_{\partial\phi,h_\ast}
\le h_{\ast}^{10} e^{h_{\ast}^{-2}\zeta^2}
\times\left[
2\mathcal{O}_8 |\lambda|^2 \bar{g}^{2-2\eta}
+\mathcal{O}_{12}|\lambda|^3 \bar{g}^{3-3\eta}
+2  \mathcal{O}_{13} |\lambda|^3 \bar{g}^{\frac{11}{4}-\eta_R}
\right]\ .
\]
Since $\eta_R\ge 3\eta$ we have $\bar{g}^{1-\eta}\le  \bar{g}^{\frac{11}{12}-\frac{1}{3}\eta_R}$
and part 1) follows.
When the $R$ term is absent, the previous estimate on $K$ reduces to
\[
|K_\Delta(\lambda,\phi_1,\zeta)|_{\partial\phi,h_\ast}
\le h_{\ast}^{10} e^{h_{\ast}^{-2}\zeta^2}
\times\left[
2\mathcal{O}_8 |\lambda|^2 \bar{g}^{2-2\eta}
+\mathcal{O}_{12}|\lambda|^3 \bar{g}^{3-3\eta}\right]
\]
from which part 2) follows.
\qed

\begin{Lemma}\label{L27lem}
If $|\lambda|\bar{g}^{\frac{1}{4}-\frac{1}{3}\eta_R}\le 1$
then for all unit cube $\Delta'$ and $\phi\in\mathbb{R}$ we have
\[
||\hat{K}_{\Delta'}(\lambda,\phi)||_{\partial\phi,\phi,h}\le 2 e^{\frac{1}{2}(\Re f,\Gamma \Re f)_{L^{-1}\Delta'}}
\sum_{n=1}^{\infty} \left(
\mathcal{O}_{19}L^{15}|\lambda|\bar{g}^{\frac{1}{4}-\frac{1}{3}\eta_R}
\right)^n
\]
where $\mathcal{O}_{19}=2^{11}\max(\mathcal{O}_{5},\mathcal{O}_{15})$.
\end{Lemma}

\noindent{\bf Proof:}
Recall that by definition
\[
\hat{K}_{\Delta'}(\lambda,\phi)=
\sum_{Y_P,Y_K}
\int {\rm d}\mu_{\Gamma} (\zeta)
\ e^{\int_{L^{-1}\Delta'} f\zeta}
\times
\prod_{\substack{\Delta\in[L^{-1}\Delta'] \\ \Delta\notin Y_P\cup Y_K}}\left[
e^{-\tilde{V}_\Delta(\phi_1)}
\right]
\times
\prod_{\Delta\in Y_P}\left[
P_\Delta(\lambda,\phi_1,\zeta)
\right]
\times
\prod_{\Delta\in Y_K}\left[
K_\Delta(\lambda,\phi_1,\zeta)
\right]
\]
where $(Y_P,Y_K)$ ranges over pairs of disjoint subsets of $[L^{-1}\Delta']$ such that not both are empty.
It is easy to see that one therefore has the following bound on $\hat{K}$:
\[
||\hat{K}_{\Delta'}(\lambda,\phi)||_{\partial\phi,\phi,h}\le
\sum_{Y_P,Y_K}
\int {\rm d}\mu_{\Gamma} (\zeta)
\ e^{\int_{L^{-1}\Delta'} (\Re f)\zeta}
\times
\left|\left|
\prod_{\substack{\Delta\in[L^{-1}\Delta'] \\ \Delta\notin Y_P\cup Y_K}}\left[
e^{-\tilde{V}_\Delta(\phi_1)}
\right]\right|\right|_{\partial\phi,\phi,h}
\]
\[
\times
\prod_{\Delta\in Y_P}
||P_\Delta(\lambda,\phi_1,\zeta)||_{\partial\phi,\phi,h}
\times
\prod_{\Delta\in Y_K}
||K_\Delta(\lambda,\phi_1,\zeta)||_{\partial\phi,\phi,h}\ .
\]
By Lemma\ref{L9lem} we have
\[
\left|\left|
\prod_{\substack{\Delta\in[L^{-1}\Delta'] \\ \Delta\notin Y_P\cup Y_K}}\left[
e^{-\tilde{V}_\Delta(\phi_1)}
\right]\right|\right|_{\partial\phi,\phi,h}\le 2\ .
\]
By Lemma \ref{L25lem} 1) with $\kappa=h_{\ast}^{-2}$ we have
\[
\prod_{\Delta\in Y_K}
||K_\Delta(\lambda,\phi_1,\zeta)||_{\partial\phi,\phi,h}
\le
\prod_{\Delta\in Y_K}\left[
\mathcal{O}_{15} h_{\ast}^{10} e^{h_{\ast}^{-2}\zeta^2} \left(|\lambda|\bar{g}^{\frac{1}{4}-\frac{1}{3}\eta_R}\right)^2
\right]\ .
\]
Since $\eta_R\ge 3\eta$, it follows from the hypotheses that
\[
|\lambda|\bar{g}^{\frac{1}{4}-\eta}\le
|\lambda|\bar{g}^{\frac{1}{4}-\frac{1}{3}\eta_R}\le 1
\]
and therefore Lemma \ref{L15lem} with $\kappa=h_{\ast}^{-2}$ 
implies
\[
\prod_{\Delta\in Y_P}
||P_\Delta(\lambda,\phi_1,\zeta)||_{\partial\phi,\phi,h}
\le
\prod_{\Delta\in Y_P}\left[
\mathcal{O}_{5} h_{\ast}^{12} e^{h_{\ast}^{-2}\zeta^2} |\lambda|\bar{g}^{\frac{1}{4}-\frac{1}{3}\eta_R}
\right]\ .
\]
Thus
\[
||\hat{K}_{\Delta'}(\lambda,\phi)||_{\partial\phi,\phi,h}\le
2\sum_{Y_P,Y_K}
\prod_{\Delta\in Y_P}\left[
\mathcal{O}_{5} h_{\ast}^{12} |\lambda|\bar{g}^{\frac{1}{4}-\frac{1}{3}\eta_R}
\right]
\times
\prod_{\Delta\in Y_K}\left[
\mathcal{O}_{15} h_{\ast}^{10} \left(|\lambda|\bar{g}^{\frac{1}{4}-\frac{1}{3}\eta_R}\right)^2
\right]
\]
\[
\times
\int {\rm d}\mu_{\Gamma} (\zeta)
\ e^{\int_{L^{-1}\Delta'} (\Re f)\zeta}
\prod_{\Delta\in Y_P\cup Y_K} e^{h_{\ast}^{-2}\zeta^2}\ .
\]
We now use Lemma \ref{L3lem} with $\alpha=h_{\ast}^{-2}=\frac{\sqrt{2}}{4} L^{-(3-2[\phi])}$ in order to bound 
the last integral. Indeed, by the standard hypothesis (\ref{fwstdeq})
\[
||(\Re f)|_{L^{-1}\Delta'}||_{L^\infty}\le
||f|_{L^{-1}\Delta'}||_{L^\infty}< L^{-(3-[\phi])}< \frac{1}{2}L^{-\frac{1}{2}(3-2[\phi])}
\]
since $L\ge 2$ implies $2L^{-\frac{3}{2}}\le 2^{-\frac{1}{2}}<1$.
Hence
\[
\int {\rm d}\mu_{\Gamma} (\zeta)
\ e^{\int_{L^{-1}\Delta'} (\Re f)\zeta}
\prod_{\Delta\in Y_P\cup Y_K} e^{h_{\ast}^{-2}\zeta^2}\le 2^{|Y_P|+|Y_K|} e^{\frac{1}{2}(\Re f,\Gamma \Re f)_{L^{-1}\Delta'}}
\]
and therefore
\[
||\hat{K}_{\Delta'}(\lambda,\phi)||_{\partial\phi,\phi,h}\le
2 e^{\frac{1}{2}(\Re f,\Gamma \Re f)_{L^{-1}\Delta'}}
\sum_{Y_P,Y_K}
\prod_{\Delta\in Y_P}\left[
2\mathcal{O}_{5} h_{\ast}^{12} |\lambda|\bar{g}^{\frac{1}{4}-\frac{1}{3}\eta_R}
\right]
\times
\prod_{\Delta\in Y_K}\left[
2\mathcal{O}_{15} h_{\ast}^{10} \left(|\lambda|\bar{g}^{\frac{1}{4}-\frac{1}{3}\eta_R}\right)^2
\right]\ .
\]
Using $h_{\ast}\le 2^{\frac{3}{4}}L$ and dropping the square in the $Y_K$ factors since
$|\lambda|\bar{g}^{\frac{1}{4}-\frac{1}{3}\eta_R}\le 1$ we arrive at
\[
||\hat{K}_{\Delta'}(\lambda,\phi)||_{\partial\phi,\phi,h}\le
2 e^{\frac{1}{2}(\Re f,\Gamma \Re f)_{L^{-1}\Delta'}}
\sum_{Y_P,Y_K}
\prod_{\Delta\in Y_P}\left[
2^{10}\mathcal{O}_{5} L^{12} |\lambda|\bar{g}^{\frac{1}{4}-\frac{1}{3}\eta_R}
\right]
\times
\prod_{\Delta\in Y_K}\left[
2^{\frac{17}{2}}\mathcal{O}_{15} L^{10} |\lambda|\bar{g}^{\frac{1}{4}-\frac{1}{3}\eta_R}
\right]
\]
\[
\le 2 e^{\frac{1}{2}(\Re f,\Gamma \Re f)_{L^{-1}\Delta'}}
\sum_{Y_P,Y_K} \rho^{|Y_P|+|Y_K|}
\]
with
\[
\rho=2^{10}\times\max(\mathcal{O}_{5},\mathcal{O}_{15})
\times L^{12}  |\lambda|\bar{g}^{\frac{1}{4}-\frac{1}{3}\eta_R}\ .
\]
Now
\[
\sum_{Y_P,Y_K} \rho^{|Y_P|+|Y_K|}=\sum_{n\ge 1}\rho^n
\sum_{\substack{i,j\ge 0 \\ i+j=n}}
\sum_{\substack{Y_P,Y_K\subset [L^{-1}\Delta'] \\ {\rm disjoint}}} \bbone\{|Y_P|=i,|Y_K|=j\}\ .
\]
Since the cardinality of $[L^{-1}\Delta']$ is $L^3$, we have from elementary combinatorics
\begin{eqnarray*}
\sum_{Y_P,Y_K} \rho^{|Y_P|+|Y_K|} & \le & \sum_{n=1}^{L^3}\rho^n
\sum_{\substack{i,j\ge 0 \\ i+j=n}} \frac{(L^3)!}{i!j!(L^3-n)!} \\
 & \le & \sum_{n=1}^{L^3}\rho^n \left(
\begin{array}{c} L^3 \\ n \end{array} 
\right) 2^n\ .
\end{eqnarray*}
We use the very coarse bound
\[
\left(
\begin{array}{c} L^3 \\ n \end{array} 
\right)=\frac{L^3(L^3-1)\cdots(L^3-n+1)}{n!}\le L^{3n}
\]
which results in
\[
\sum_{Y_P,Y_K} \rho^{|Y_P|+|Y_K|}
\le \sum_{n=1}^{\infty} (2L^3\rho)^n\ .
\]
The latter inserted in the previous estimate for $\hat{K}$
gives the desired inequality.
\qed

\begin{Lemma}\label{L28lem}
If $R=0$ and $|\lambda|\bar{g}^{\frac{1}{4}-\eta}\le 1$
then for all unit cube $\Delta'$ and $\phi\in\mathbb{R}$ we have
\[
||\hat{K}_{\Delta'}(\lambda,\phi)||_{\partial\phi,\phi,h}\le 2 e^{\frac{1}{2}(\Re f,\Gamma \Re f)_{L^{-1}\Delta'}}
\sum_{n=1}^{\infty} \left(
\mathcal{O}_{20}L^{15}|\lambda|\bar{g}^{\frac{1}{4}-\eta}
\right)^n
\]
where $\mathcal{O}_{20}=2^{11}\max(\mathcal{O}_{5},\mathcal{O}_{16})$.
\end{Lemma}

\noindent{\bf Proof:}
One can repeat the last proof verbatim except that one must use part 2) of Lemma \ref{L25lem}
instead of part 1). This accounts for $\mathcal{O}_{16}$ featuring in the new constant instead of $\mathcal{O}_{15}$.
\qed

\begin{Lemma}\label{L29lem}
If $|\lambda|\bar{g}^{\frac{11}{12}-\frac{1}{3}\eta_R}\le 1$
then for all unit cube $\Delta'$ and $\phi\in\mathbb{R}$ we have
\[
|\hat{K}_{\Delta'}(\lambda,\phi)|_{\partial\phi,h_\ast}\le 2 e^{\frac{1}{2}(\Re f,\Gamma \Re f)_{L^{-1}\Delta'}}
\sum_{n=1}^{\infty} \left(
\mathcal{O}_{21}L^{15}|\lambda|\bar{g}^{\frac{11}{12}-\frac{1}{3}\eta_R}
\right)^n
\]
where $\mathcal{O}_{21}=2^{11}\max(\mathcal{O}_{6},\mathcal{O}_{17})$.
\end{Lemma}

\noindent{\bf Proof:}
Again from the definition of $\hat{K}$
one easily deduces the estimate
\[
|\hat{K}_{\Delta'}(\lambda,\phi)|_{\partial\phi,h_\ast}\le
\sum_{Y_P,Y_K}
\int {\rm d}\mu_{\Gamma} (\zeta)
\ e^{\int_{L^{-1}\Delta'} (\Re f)\zeta}
\times
\left|
\prod_{\substack{\Delta\in[L^{-1}\Delta'] \\ \Delta\notin Y_P\cup Y_K}}\left[
e^{-\tilde{V}_\Delta(\phi_1)}
\right]\right|_{\partial\phi,h_\ast}
\]
\[
\times
\prod_{\Delta\in Y_P}
|P_\Delta(\lambda,\phi_1,\zeta)|_{\partial\phi,h_\ast}
\times
\prod_{\Delta\in Y_K}
|K_\Delta(\lambda,\phi_1,\zeta)|_{\partial\phi,h_\ast}\ .
\]
While
\[
\left|
\prod_{\substack{\Delta\in[L^{-1}\Delta'] \\ \Delta\notin Y_P\cup Y_K}}\left[
e^{-\tilde{V}_\Delta(\phi_1)}
\right]\right|_{\partial\phi,h_\ast}\le 2
\]
by Lemma \ref{L9lem}, we have
\[
|P_\Delta(\lambda,\phi_1,\zeta)|_{\partial\phi,h_\ast}\le \mathcal{O}_6 h_{\ast}^{12} e^{h_{\ast}^{-2}\zeta^2}
|\lambda|\bar{g}^{1-\eta} \le
\mathcal{O}_6 h_{\ast}^{12} e^{h_{\ast}^{-2}\zeta^2}
|\lambda|\bar{g}^{\frac{11}{12}-\frac{1}{3}\eta_R}
\]
by Lemma \ref{L16lem}. Indeed, $0<\bar{g}\le 1$ and $\eta_R\ge 3\eta$ ensure that
$|\lambda|\bar{g}^{1-\eta} \le|\lambda|\bar{g}^{\frac{11}{12}-\frac{1}{3}\eta_R}\le 1$.
Finally, Lemma \ref{L26lem} 1) provides us with the last ingredient
\[
|K_\Delta(\lambda,\phi_1,\zeta)|_{\partial\phi,h_\ast}\le \mathcal{O}_{17} h_{\ast}^{10} e^{h_{\ast}^{-2}\zeta^2}
\left(|\lambda|\bar{g}^{\frac{11}{12}-\frac{1}{3}\eta_R}\right)^2\ .
\]
The rest of the proof is exactly the same as that of Lemma \ref{L27lem}.
\qed

\begin{Lemma}\label{L30lem}
If $R=0$ and $|\lambda|\bar{g}^{1-\eta}\le 1$
then for all unit cube $\Delta'$ and $\phi\in\mathbb{R}$ we have
\[
|\hat{K}_{\Delta'}(\lambda,\phi)|_{\partial\phi,h_\ast}\le 2 e^{\frac{1}{2}(\Re f,\Gamma \Re f)_{L^{-1}\Delta'}}
\sum_{n=1}^{\infty} \left(
\mathcal{O}_{22}L^{15}|\lambda|\bar{g}^{1-\eta}
\right)^n
\]
where $\mathcal{O}_{22}=2^{11}\max(\mathcal{O}_{6},\mathcal{O}_{18})$.
\end{Lemma}

\noindent{\bf Proof:}
The argument is the same as in the last proof except for the use of Part 2) of Lemma \ref{L26lem}
instead of Part 1).
\qed

\begin{Lemma}\label{L31lem}
For all $\Delta'\in\mathbb{L}$ and $\Delta_1\in[L^{-1}\Delta']$, the quantity
$J_+(\phi)$ defined in \S\ref{algdefsec} satisfies the bound
\[
|J_+(\phi)|_{\partial\phi,L^{[\phi]}h_\ast}\le\mathcal{O}_{23} |||R_{\Delta_1}|||_{\bar{g}}
\]
where $\mathcal{O}_{23}=4\mathcal{O}_{13}\times\exp\left(2^{-\frac{3}{2}}\right)$.
\end{Lemma}

\noindent{\bf Proof:}
Recall that by definition
\[
J_+(\phi)=e^{-\frac{1}{2}(f,\Gamma f)_{L^{-1}\Delta'}}\int {\rm d}\mu_{\Gamma} (\zeta)
\ e^{\int_{L^{-1}\Delta'} f\zeta} R_{\Delta_1}(\phi_1+\zeta)
\]
and therefore one readily obtains
\[
|J_+(\phi)|_{\partial\phi,L^{[\phi]}h_\ast}\le
e^{-\frac{1}{2}\Re (f,\Gamma f)_{L^{-1}\Delta'}}\int {\rm d}\mu_{\Gamma} (\zeta)
\ e^{\int_{L^{-1}\Delta'} (\Re f)\zeta} |R_{\Delta_1}(\phi_1+\zeta)|_{\partial\phi,L^{[\phi]}h_\ast}\ .
\]
By the definitions of the seminorms and the chain rule one has
\[
|R_{\Delta_1}(\phi_1+\zeta)|_{\partial\phi,L^{[\phi]}h_\ast}=
||R_{\Delta_1}(\phi_1+\zeta)||_{\partial\phi,0,L^{[\phi]}h_\ast}=
||R_{\Delta_1}(\psi+\zeta)||_{\partial\psi,0,h_\ast}=
||R_{\Delta_1}(\psi)||_{\partial\psi,\zeta,h_\ast}\ .
\]
From Lemma \ref{L23lem} we then derive
\begin{eqnarray*}
|R_{\Delta_1}(\phi_1+\zeta)|_{\partial\phi,L^{[\phi]}h_\ast} & \le & \mathcal{O}_{13}
e^{h_{\ast}^{-2}\zeta_{\Delta_1}^2}\left[
|R_{\Delta_1}(\psi)|_{\partial\psi,h_\ast}+h_{\ast}^{9} h^{-9} \sup_{\psi\in\mathbb{R}}
||R_{\Delta_1}(\psi)||_{\partial\psi,\psi,h}
\right] \\
 & \le & \mathcal{O}_{13}
e^{h_{\ast}^{-2}\zeta_{\Delta_1}^2}
|||R_{\Delta_1}|||_{\bar{g}}\left(1+h_{\ast}^{9}h^{-9}\bar{g}^{-2}\right) \\
 & \le & 2 \mathcal{O}_{13}
e^{h_{\ast}^{-2}\zeta_{\Delta_1}^2}
|||R_{\Delta_1}|||_{\bar{g}}
\end{eqnarray*}
by the standard hypothesis (\ref{c1c2stdeq}).
As a result
\[
|J_+(\phi)|_{\partial\phi,L^{[\phi]}h_\ast}\le
e^{-\frac{1}{2}\Re (f,\Gamma f)_{L^{-1}\Delta'}}
\times 2\mathcal{O}_{13}\times |||R_{\Delta_1}|||_{\bar{g}}\times
\int {\rm d}\mu_{\Gamma} (\zeta)
\ e^{\int_{L^{-1}\Delta'} (\Re f)\zeta} e^{h_{\ast}^{-2}\zeta_{\Delta_1}^2}\ .
\]
The standard hypothesis (\ref{fwstdeq}) again allows to use Lemma \ref{L3lem} with $\alpha=h_{\ast}^{-2}$
to the effect that
\[
|J_+(\phi)|_{\partial\phi,L^{[\phi]}h_\ast}\le
4\mathcal{O}_{13}\times |||R_{\Delta_1}|||_{\bar{g}}\times
\exp\left\{
-\frac{1}{2}\Re (f,\Gamma f)_{L^{-1}\Delta'}+\frac{1}{2} (\Re f,\Gamma \Re f)_{L^{-1}\Delta'}
\right\}
\]
holds.
Note that
\[
\Re (f,\Gamma f)_{L^{-1}\Delta'}= (\Re f,\Gamma \Re f)_{L^{-1}\Delta'}-(\Im f,\Gamma \Im f)_{L^{-1}\Delta'}
\]
and thus
\[
|J_+(\phi)|_{\partial\phi,L^{[\phi]}h_\ast}\le
4\mathcal{O}_{13}\times |||R_{\Delta_1}|||_{\bar{g}}\times
\exp\left\{
\frac{1}{2} (\Im f,\Gamma \Im f)_{L^{-1}\Delta'}
\right\}\ .
\]
But
\begin{eqnarray*}
\left|(\Im f,\Gamma \Im f)_{L^{-1}\Delta'}\right| & \le & \int_{(L^{-1}\Delta')^2}
{\rm d}^3 x {\rm d}^3 y\ |\Gamma{x-y}|\ |\Im f(x)|\ |\Im f(y)| \\
 & \le & ||f|_{L^{-1}\Delta'}||_{L^\infty}^2\times L^3\times ||\Gamma||_{L^1} \\
 & \le & L^{-2(3-[\phi])}\times L^3\times \frac{L^{3-2[\phi]}}{\sqrt{2}}\\
 & \le & \frac{1}{\sqrt{2}}
\end{eqnarray*}
because of the standard hypothesis (\ref{fwstdeq}), the finite range property of $\Gamma$ and the bound in Corollary \ref{gamL1cor}.
Inserting this last inequality in the previous estimate for $J_+$
gives the wanted bound.
\qed

\begin{Lemma}\label{L32lem}
For all $\Delta'\in\mathbb{L}$, $\Delta_1\in[L^{-1}\Delta']$ and integer $k$ such that $0\le k\le 4$
the $\delta\beta$ quantities defined in \S\ref{algdefsec} satisfy
\[
|\delta\beta_{k,3,\Delta',\Delta_1}|\le \mathcal{O}_{24} \times \left(L^{[\phi]}h_\ast\right)^{-k} \times
|||R_{\Delta_1}|||_{\bar{g}}
\]
and
\[
|\delta\beta_{k,3,\Delta'}|\le \mathcal{O}_{24} \times L^{3-k[\phi]} \times\max_{\Delta_1\in[L^{-1}\Delta']}
|||R_{\Delta_1}|||_{\bar{g}}
\]
with
\[
\mathcal{O}_{24}=48\times \mathcal{O}_{23}\times
\sum_{i=0}^{4}\sum_{j,n,l} \left|\#_{k,i,j,n,l}\right| 2^j \left(\frac{3}{2}\right)^{n}
\]
where $\#_{k,i,j,n,l}$ denote the numerical coefficients in the explicit formulas produced by Maple from \S\ref{algdefsec}.
\end{Lemma}

\noindent{\bf Proof:}
Recall that
\[
\delta\beta_{k,3,\Delta',\Delta_1}=\sum_{i=0}^{4} M_{k,i} a_i
\]
where
\begin{eqnarray*}
a_i & = & \exp\left[
-C_0(0)L^{-2[\phi]}\beta_{2,\Delta_1}+3C_0(0)^2 L^{-4[\phi]}\beta_{4,\Delta_1}-\frac{1}{2}(f,\Gamma f)_{L^{-1}\Delta'}
\right] \\
 & & \times L^{-i[\phi]}\times \int {\rm d}\mu_{\Gamma} (\zeta)
\ e^{\int_{L^{-1}\Delta'} f\zeta} R_{\Delta_1}^{(i)}(\zeta) \\
 & = &  \exp\left[
-C_0(0)L^{-2[\phi]}\beta_{2,\Delta_1}+3C_0(0)^2 L^{-4[\phi]}\beta_{4,\Delta_1}
\right] \times J_{+}^{(i)}(0)
\end{eqnarray*}
and
\[
M_{k,i}=\sum_{j,n,l} \#_{k,i,j,n,l} C_0(0)^j L^{-(l_1+\cdots+l_n)[\phi]}\beta_{l_1,\Delta_1}\cdots \beta_{l_n,\Delta_1}\ .
\]
From the standard hypotheses we have $|\beta_{2,\Delta_1}|<\bar{g}^{1-\eta}\le \bar{g}^{\frac{3}{4}}$ since $\eta<\frac{1}{4}$.
We also have $|\beta_{4,\Delta_1}|<\frac{3}{2}\bar{g}$. Using $C_0(0)<2$, $L^{-[\phi]}\le 1$ and the standard hypothesis (\ref{expgstdeq})
we then deduce the bounds
\[
|a_i|\le |J_{+}^{(i)}(0)|\times \exp\left[2\bar{g}^{\frac{1}{4}}+18\bar{g}\right]\le 2|J_{+}^{(i)}(0)|\ .
\]
By definition of the seminorms
\[
|J_{+}^{(i)}(0)|\le i! (L^{[\phi]}h_\ast)^{-i}  |J_{+}(0)|_{\partial\phi,L^{[\phi]}h_\ast}\ .
\]
Since $i\le 4$ we then get from the last inequality
\[
|a_i|\le 48 (L^{[\phi]}h_\ast)^{-i}  |J_{+}(0)|_{\partial\phi,L^{[\phi]}h_\ast}\ .
\]
Now recall that the sum expressing the $M_{k,i}$ 
is quantified over $j\ge 0$, $n\ge 0$ and $l=(l_1,\ldots,l_n)\in\{1,\ldots,4\}^n$. For the numerical coefficients
$\#_{k,i,j,n,l}$ to be nonzero the constraint
\[
l_1+\cdots+l_n-2j=k-i
\]
must be satisfied.
The $\beta_{l_\nu,\Delta_1}$ are bounded by $\bar{g}^{1-\eta}$ or $\frac{3}{2}\bar{g}$ which can be replaced by a uniform
worst case scenario bound of $\frac{3}{2}\bar{g}^{1-\eta}$.
We can thus write
\[
|M_{k,i}|\le \sum_{j,n,l} |\#_{k,i,j,n,l}| 2^j L^{-(l_1+\cdots+l_n)[\phi]}\times\left(\frac{3}{2}\bar{g}^{1-\eta}\right)^n\ .
\]
We now consider two different cases in order to continue estimating the $|M_{k,i}|$.

\noindent{\bf 1st case:} Suppose $i\ge k$. Since the $l$'s are positive, we have $L^{-(l_1+\cdots+l_n)[\phi]}\le 1$. We also use the coarse
bound $\bar{g}^{1-\eta}\le 1$ which results from the standard hypothesis $\eta<\frac{1}{4}$.
We then simply write
\[
|M_{k,i}|\le \sum_{j,n,l} |\#_{k,i,j,n,l}| 2^j \times\left(\frac{3}{2}\right)^n\ .
\]

\noindent{\bf 2nd case:} Suppose $i<k$. Since $j\ge 0$, the previous constraint implies
\[
l_1+\cdots+l_n=2j+k-i\ge k-i
\]
and therefore $L^{-(l_1+\cdots+l_n)[\phi]}\le L^{-(k-i)[\phi]}$.
One can also infer that $n\ge 1$ since $l_1+\cdots+l_n\ge k-i>0$ with the consequence that $(\bar{g}^{1-\eta})^n\le \bar{g}^{1-\eta}$.
The bound on $|M_{k,i}|$ which results from these remarks can reorganized as
\[
|M_{k,i}|\le \sum_{j,n,l} |\#_{k,i,j,n,l}| 2^j \times\left(h_\ast L^{[\phi]}\right)^{-(k-i)}\times
\left(\frac{3}{2}\right)^n \bar{g}^{1-\eta} h_{\ast}^{k-i}\ .
\]
Since $0\le i<k\le 4$, $h_\ast\ge 1$ and $\epsilon\le 1$ we have
\[
h_{\ast}^{k-i}\le h_{\ast}^{4}=\left(2^{\frac{3}{4}}L^{\frac{3+\epsilon}{4}}\right)^{4}\le 8 L^4\ .
\]
The standard hypothesis (\ref{expgstdeq}) now allows us to write
\[
|M_{k,i}|\le \left(h_\ast L^{[\phi]}\right)^{-(k-i)}\times
\sum_{j,n,l} |\#_{k,i,j,n,l}| 2^j 
\left(\frac{3}{2}\right)^n
\]
which is the wanted bound for $|M_{k,i}|$ in this second case.

We now combine the previous consideration and get
\[
|\delta\beta_{k,3,\Delta',\Delta_1}|\le 
\sum_{i=k}^4 |M_{k,i}|\ |a_i| + \sum_{0\le i<k} |M_{k,i}|\ |a_i|
\]
\[
\le \sum_{i=k}^4 48 (L^{[\phi]}h_\ast)^{-i}  |J_{+}(0)|_{\partial\phi,L^{[\phi]}h_\ast}\times 
\sum_{j,n,l} |\#_{k,i,j,n,l}| 2^j \times\left(\frac{3}{2}\right)^n
\]
\[
+\sum_{0\le i<k} 48 (L^{[\phi]}h_\ast)^{-k}  |J_{+}(0)|_{\partial\phi,L^{[\phi]}h_\ast}\times 
\sum_{j,n,l} |\#_{k,i,j,n,l}| 2^j \times\left(\frac{3}{2}\right)^n
\]
Since $L^{[\phi]}$ and $h_\ast$ are greater than 1 we have $(L^{[\phi]}h_\ast)^{-i}\le (L^{[\phi]}h_\ast)^{-k}$ when $i\ge k$.
We can then more conveniently write
\[
|\delta\beta_{k,3,\Delta',\Delta_1}|\le 
48 (L^{[\phi]}h_\ast)^{-k}  |J_{+}(0)|_{\partial\phi,L^{[\phi]}h_\ast}\times 
\sum_{i=0}^{4}\sum_{j,n,l} |\#_{k,i,j,n,l}| 2^j \times\left(\frac{3}{2}\right)^n
\]
from which the desired follows thanks to Lemma \ref{L31lem}.
Finally the second bound on $|\delta\beta_{k,3,\Delta'}|$ follows simply by summing over $\Delta_1\in [L^{-1}\Delta']$
and discarding the factors $h_{\ast}^{-k}\le 1$.
\qed

\begin{Lemma}\label{L33lem}
For all unit cube $\Delta'$ and all integer $k$ such that $0\le k\le 4$, we have
\[
|\delta\beta_{k,1,\Delta'}|\le \mathcal{O}_{25}\bar{g}^{1-\eta} L^{\frac{5}{2}} \bbone\{k\le 3\}
\]
where $\mathcal{O}_{25}=\frac{27}{2}$.
\end{Lemma}

\noindent{\bf Proof:}
From the definition we get
\[
|\delta\beta_{k,1,\Delta'}|\le \sum_b
\bbone\left\{
\begin{array}{c}
k+b\le 4 \\
b\ge 1
\end{array}
\right\}
\frac{(k+b)!}{k!\ b!}\ L^{-k[\phi]}\ 
\left|\parbox{2.1cm}{
\psfrag{a}{$\beta_{k+b}$}\psfrag{f}{$f$}\psfrag{b}{$b$}
\raisebox{-1ex}{
\includegraphics[width=2.1cm]{Fig5.eps}}
}\ \ \right|
\]
where the Feynman diagram has been defined in \S\ref{algdefsec}. This already shows the vanishing when $k=4$.
We now restrict to the case $k\le 3$. We bound the $f$'s by $||f|_{L^{-1}\Delta'}||_{L^\infty}$ and perform the integration
over the corresponding points of evaluation in $\mathbb{Q}_p^3$ which give $||\Gamma||_{L^1}$ factors.
We thus get the bound
\[
\left|\parbox{2.1cm}{
\psfrag{a}{$\beta_{k+b}$}\psfrag{f}{$f$}\psfrag{b}{$b$}
\raisebox{-1ex}{
\includegraphics[width=2.1cm]{Fig5.eps}}
}\ \ \right|\le
||f|_{L^{-1}\Delta'}||_{L^\infty}^b\times ||\Gamma||_{L^1}^b\times L^3\times
\max_{\Delta\in[L^{-1}\Delta']}|\beta_{k+b,\Delta}|\ .
\]
\[
\le \frac{3}{2} \bar{g}^{1-\eta} L^3 \times\left(\frac{1}{\sqrt{2}}L^{3-2[\phi]} L^{-(3-[\phi])}\right)^b\ .
\]
We discard the $\frac{1}{\sqrt{2}}$ factors and bound the remaining power of $L$, namely $L^{3-b[\phi]}$, by
$L^{\frac{5}{2}}$ since $b\ge 1$ and $\epsilon\le 1$.
Hence
\[
|\delta\beta_{k,1,\Delta'}|\le \sum_b
\bbone\left\{
\begin{array}{c}
k+b\le 4 \\
b\ge 1
\end{array}
\right\}
\frac{(k+b)!}{k!\ b!}\ \frac{3}{2} \bar{g}^{1-\eta} L^\frac{5}{2}
\]
where we also discarded the factor $L^{-k[\phi]}$. Since
\[
\max_{0\le k\le 3} \sum_b
\bbone\left\{
\begin{array}{c}
k+b\le 4 \\
b\ge 1
\end{array}
\right\}
\frac{(k+b)!}{k!\ b!} =9
\]
the lemma is proved.
\qed

\begin{Lemma}\label{L34lem}
For all unit cube $\Delta'$ and all integer $k$ such that $0\le k\le 4$, we have
\[
|\delta\beta_{k,2,\Delta'}|\le \mathcal{O}_{26}\bar{g}^{2-2\eta} L^{5}
\]
with
\[
\mathcal{O}_{26}=\frac{9}{2}\sum_{a_1,a_2,b_1,b_2,m}
\bbone\left\{
\begin{array}{c}
a_i+b_i\le 4 \\
a_i\ge 0\ ,\ b_i\ge 1\\
1\le m\le \min(b_1,b_2)
\end{array}
\right\}
\frac{(a_1+b_1)!\ (a_2+b_2)!}{a_1!\ a_2!\ m!\ (b_1-m)!\ (b_2-m)!}
\]
\[
\times C(a_1,a_2|k)\times
2^{\frac{a_1+a_2-k}{2}}
\]
\[
+\sum_b
\bbone\left\{
\begin{array}{c}
k+b=5\ {\rm or}\ 6 \\
b\ge 0
\end{array}
\right\}
\frac{(k+b)!}{k!\ b!}
\]
where the $C(a_1,a_2|k)$ are the connection coefficients defined in S\ref{algdefsec}.
\end{Lemma}

\noindent{\bf Proof:}
From the definition we have
\[
|\delta\beta_{k,2,\Delta'}|\le\sum_{a_1,a_2,b_1,b_2,m}
\bbone\left\{
\begin{array}{c}
a_i+b_i\le 4 \\
a_i\ge 0\ ,\ b_i\ge 1\\
1\le m\le \min(b_1,b_2)
\end{array}
\right\}
\frac{(a_1+b_1)!\ (a_2+b_2)!}{a_1!\ a_2!\ m!\ (b_1-m)!\ (b_2-m)!}
\]
\[
\times \frac{1}{2} C(a_1,a_2|k)\times
L^{-(a_1+a_2)[\phi]}\times C_0(0)^{\frac{a_1+a_2-k}{2}}\times
\left|\ \ 
\parbox{4cm}{
\psfrag{a}{$\beta_{a_1+b_1}$}
\psfrag{b}{$\beta_{a_2+b_2}$}
\psfrag{c}{$\scriptstyle{b_1-m}$}
\psfrag{d}{$\scriptstyle{b_2-m}$}
\psfrag{m}{$m$}
\psfrag{f}{$f$}
\raisebox{-1ex}{
\includegraphics[width=4cm]{Fig6.eps}}
}\ \ \ \ \ \right|
\]
\[
+\sum_b
\bbone\left\{
\begin{array}{c}
k+b=5\ {\rm or}\ 6 \\
b\ge 0
\end{array}
\right\}
\frac{(k+b)!}{k!\ b!}\  L^{-k[\phi]}
\left|\ \parbox{2.1cm}{
\psfrag{a}{$W_{k+b}$}\psfrag{f}{$f$}\psfrag{b}{$b$}
\raisebox{-1ex}{
\includegraphics[width=2.1cm]{Fig5.eps}}
}\ \ \right|
\]
The $W$ diagrams are bounded in the same way as in the previous lemma by
\[
\left|\ \parbox{2.1cm}{
\psfrag{a}{$W_{k+b}$}\psfrag{f}{$f$}\psfrag{b}{$b$}
\raisebox{-1ex}{
\includegraphics[width=2.1cm]{Fig5.eps}}
}\ \ \right|
\le \bar{g}^{2-2\eta} L^{3-b[\phi]}\le \bar{g}^{2-2\eta} L^{3}\ .
\]
The diagrams with two internal $\beta$ vertices are bounded using the same method which
gives a factor of $||f|_{L^{-1}\Delta'}||_{L^\infty}\times ||\Gamma||_{L^1}\le L^{-[\phi]}$ per $f$ external vertex.
The $|\beta|$'s are bounded in a uniform manner by $\frac{3}{2}\bar{g}^{1-\eta}$ 
and this results in the estimate
\[
\left|\ \ 
\parbox{4cm}{
\psfrag{a}{$\beta_{a_1+b_1}$}
\psfrag{b}{$\beta_{a_2+b_2}$}
\psfrag{c}{$\scriptstyle{b_1-m}$}
\psfrag{d}{$\scriptstyle{b_2-m}$}
\psfrag{m}{$m$}
\psfrag{f}{$f$}
\raisebox{-1ex}{
\includegraphics[width=4cm]{Fig6.eps}}
}\ \ \ \ \ \right|\le \frac{9}{4}\bar{g}^{2-2\eta}\times
\left(L^{-[\phi]}\right)^{b_1+b_2-2m}\times
\int_{(L^{-1}\Delta')^2} {\rm d}^3 x_1\ {\rm d}^3 x_2\ |\Gamma(x_1-x_2)|^m\ .
\]
Since $b_1+b_2-2m$ can take all integer values between 0 and 6 we simply discard the factor $\left(L^{-[\phi]}\right)^{b_1+b_2-2m}$
in the bound.
By the ultrametricity and the finite range property of $\Gamma$ we have
\[
\int_{(L^{-1}\Delta')^2} {\rm d}^3 x_1\ {\rm d}^3 x_2\ |\Gamma(x_1-x_2)|^m \ =\ L^3 ||\Gamma||_{L^m}^{m}\ .
\]
For the purposes of this lemma and for the relevant values of $m$, namely 1, 2 or 3, we use the blanket estimate
\[
||\Gamma||_{L^m}^{m}\le ||\Gamma||_{L^\infty}^{m-1}\times ||\Gamma||_{L^1}\le 2^{m-1} \frac{1}{\sqrt{2}}L^{3-2[\phi]}\le 4 L^2\ .
\]
We therefore have the estimates
\[
\left|\ \ 
\parbox{4cm}{
\psfrag{a}{$\beta_{a_1+b_1}$}
\psfrag{b}{$\beta_{a_2+b_2}$}
\psfrag{c}{$\scriptstyle{b_1-m}$}
\psfrag{d}{$\scriptstyle{b_2-m}$}
\psfrag{m}{$m$}
\psfrag{f}{$f$}
\raisebox{-1ex}{
\includegraphics[width=4cm]{Fig6.eps}}
}\ \ \ \ \ \right|\le 9\bar{g}^{2-2\eta} L^5
\]
and
\[
\left|\ \parbox{2.1cm}{
\psfrag{a}{$W_{k+b}$}\psfrag{f}{$f$}\psfrag{b}{$b$}
\raisebox{-1ex}{
\includegraphics[width=2.1cm]{Fig5.eps}}
}\ \ \right|
\le \bar{g}^{2-2\eta} L^3
\]
which we insert in the previous bound on $|\delta\beta_{k,2,\Delta'}|$. We drop the
$L^{-(a_1+a_2)[\phi]}$ and $L^{-k[\phi]}$ factors and use $C_0(0)<2$ 
to arrive at the wanted result.
\qed

\begin{Lemma}\label{L35lem}
Let $\mathcal{O}_{27}=16\times 25\times[32\mathcal{O}_{25}+40\mathcal{O}_{26}+40\mathcal{O}_{24}]$.
Provided $\lambda$ satisfies $\mathcal{O}_{27}L^{5}|\lambda|\bar{g}^{\frac{1}{4}-\eta}\le 1$
we have, for all $\Delta'\in\mathbb{L}$ and $\phi\in\mathbb{R}$, 
\[
||e^{-\hat{V}_{\Delta'}(\phi)+\delta V_{\Delta'}(\lambda,\phi)}||_{\partial\phi,\phi,h}\le 3\ .
\]
\end{Lemma}

\noindent{\bf Proof:}
By the multiplicative property of the seminorms and by Lemma \ref{L10lem} we have
\begin{eqnarray*}
||e^{-\hat{V}_{\Delta'}(\phi)+\delta V_{\Delta'}(\lambda,\phi)}||_{\partial\phi,\phi,h}
 & \le & ||e^{-\hat{V}_{\Delta'}(\phi)}||_{\partial\phi,\phi,h}\times 
 ||e^{\delta V_{\Delta'}(\lambda,\phi)}||_{\partial\phi,\phi,h} \\
 & \le & 2 e^{-\frac{\bar{g}}{16}\phi^4}\times \exp\left[ ||\delta V_{\Delta'}(\lambda,\phi)||_{\partial\phi,\phi,h}
\right]\ .
\end{eqnarray*}
Now by definition
\begin{eqnarray*}
\delta V_{\Delta'}(\lambda,\phi) & = & \sum_{k=0}^{4} \delta\beta_{k,\Delta'}(\lambda)\ :\phi^k:_{C_0} \\
 & = &  \sum_{k=0}^{4}  \left(\sum_{j=1}^{3}
\lambda^j \delta\beta_{k,j,\Delta'}\right)\ :\phi^k:_{C_0}
\end{eqnarray*}
and therefore
\[
||\delta V_{\Delta'}(\lambda,\phi)||_{\partial\phi,\phi,h}\le
\sum_{\substack{0\le k\le 4 \\ 1\le j\le 3}} |\lambda|^j\ 
|\delta\beta_{k,j,\Delta'}|\ ||:\phi^k:_{C_0}||_{\partial\phi,\phi,h}\ .
\]
As in the proof of Lemma \ref{L11lem} we have, for $0\le k\le 4$,
\begin{eqnarray*}
||:\phi^k:_{C_0}||_{\partial\phi,\phi,h} & \le & 25\times\max_{0\le a\le k} ||\phi^a||_{\partial\phi,\phi,h} \\
 & \le & 25\times\max_{0\le a\le k} (h+|\phi|)^a
\end{eqnarray*}
For the definition of $h$ one can write
\begin{eqnarray*}
\bar{g}^{\frac{k}{4}}||:\phi^k:_{C_0}||_{\partial\phi,\phi,h} & \le & 25\times\max_{0\le a\le k}  \bar{g}^{\frac{a}{4}}(h+|\phi|)^a\\
 & \le & 25\times\max_{0\le a\le k} (c_1+\bar{g}^{\frac{1}{4}}|\phi|)^a
\end{eqnarray*}
since $a\le k$ and $\bar{g}\le 1$.
Since $c_1<1$ we has the more convenient bounds
\begin{eqnarray*}
\bar{g}^{\frac{k}{4}}||:\phi^k:_{C_0}||_{\partial\phi,\phi,h} & \le & 25\times\max_{0\le a\le k}  (1+\bar{g}^{\frac{1}{4}}|\phi|)^a\\
 & \le & 25 (1+\bar{g}^{\frac{1}{4}}|\phi|)^k
\end{eqnarray*}
which result in
\[
||\delta V_{\Delta'}(\lambda,\phi)||_{\partial\phi,\phi,h}\le 25
\sum_{\substack{0\le k\le 4 \\ 1\le j\le 3}} |\lambda|^j\ 
|\delta\beta_{k,j,\Delta'}|\times \bar{g}^{-\frac{k}{4}} (1+\bar{g}^{\frac{1}{4}}|\phi|)^k\ .
\]
We now bound the contributions of each $j$ separately.
For $j=1$, one has by applying Lemma \ref{L33lem}
\[
\sum_{0\le k\le 4}
|\delta\beta_{k,1,\Delta'}|\times \bar{g}^{-\frac{k}{4}} (1+\bar{g}^{\frac{1}{4}}|\phi|)^k
\le \mathcal{O}_{25}\bar{g}^{1-\eta} L^{\frac{5}{2}}
 \sum_{0\le k\le 3} \bar{g}^{-\frac{k}{4}} (1+\bar{g}^{\frac{1}{4}}|\phi|)^k\ .
\]
For nonnegative numbers $A$ and $B$ one has the classic inequality
\[
\frac{A+B}{2}\le \left(\frac{A^a+B^a}{2}\right)^{\frac{1}{a}}
\]
for all $a\ge 1$ which can be conveniently rewritten as
\begin{equation}
(A+B)^a\le 2^{a-1} (A^a+B^a)\ .
\label{L1Labdeq}
\end{equation}
For $0\le k\le 3$ we bound $\bar{g}^{-\frac{k}{4}}$ by $\bar{g}^{-\frac{3}{4}}$
and also write
\[
(1+\bar{g}^{\frac{1}{4}}|\phi|)^k\le (1+\bar{g}^{\frac{1}{4}}|\phi|)^4\le 8
(1+\bar{g}\phi^4)
\]
using (\ref{L1Labdeq}) with $a=4$.
As a result we have
\[
\sum_{0\le k\le 4}
|\delta\beta_{k,1,\Delta'}|\times \bar{g}^{-\frac{k}{4}} (1+\bar{g}^{\frac{1}{4}}|\phi|)^k
\le 4\times 8\times \mathcal{O}_{25}\bar{g}^{1-\eta} L^{\frac{5}{2}}\times (1+\bar{g}\phi^4)\ .
\]
We use a similar for the $j=2$ contribution where the sum over $k$ goes from $0$ to $4$. Namely, bounding $\bar{g}^{-\frac{k}{4}}$
by $\bar{g}$, $(1+\bar{g}^{\frac{1}{4}}|\phi|)^4$ by $8(1+\bar{g}\phi^4)$ and using Lemma \ref{L34lem}, we get
\[
\sum_{0\le k\le 4}
|\delta\beta_{k,2,\Delta'}|\times \bar{g}^{-\frac{k}{4}} (1+\bar{g}^{\frac{1}{4}}|\phi|)^k
\le 5\times 8\times \mathcal{O}_{26}\bar{g}^{1-2\eta} L^{5}\times (1+\bar{g}\phi^4)\ .
\]
The same procedure for the $j=3$ contribution, this time using Lemma \ref{L32lem} and the standard hypothesis (\ref{normstdeq}), gives
\[
\sum_{0\le k\le 4}
|\delta\beta_{k,3,\Delta'}|\times \bar{g}^{-\frac{k}{4}} (1+\bar{g}^{\frac{1}{4}}|\phi|)^k
\le 5\times 8\times \mathcal{O}_{24}\bar{g}^{\frac{7}{4}-\eta_R} L^{3}\times (1+\bar{g}\phi^4)\ .
\]
Hence one can collect the previous separate estimates into
\[
||\delta V_{\Delta'}(\lambda,\phi)||_{\partial\phi,\phi,h}\le 25 (1+\bar{g}\phi^4)
\left[32\mathcal{O}_{25}|\lambda| \bar{g}^{\frac{1}{4}-\eta} L^{\frac{5}{2}}
+40\mathcal{O}_{26}|\lambda|^2 \bar{g}^{1-2\eta} L^{5}
+40\mathcal{O}_{24}|\lambda|^3 \bar{g}^{\frac{7}{4}-\eta_R} L^{3}
\right]\ .
\]
Let $\rho=|\lambda|\bar{g}^{\frac{1}{4}-\eta}$ then clearly $|\lambda|^2\bar{g}^{1-2\eta}=\rho^2\bar{g}^{\frac{1}{2}}\le \rho^2$.
Also because of the standard hypothesis $\eta_R\le 1+3\eta$ we have $|\lambda|^3\bar{g}^{\frac{7}{4}-\eta_R}\le \rho^3$.
Notice that since for instance $\mathcal{O}_{25}=\frac{27}{2}$, we clearly have $\mathcal{O}_{27}>1$.
Thus the hypothesis of the present lemma implies in particular that $\rho\le 1$.
We therefore have the more convenient bound
\[
||\delta V_{\Delta'}(\lambda,\phi)||_{\partial\phi,\phi,h}\le 25 (1+\bar{g}\phi^4)
L^5\rho\times\left[
32\mathcal{O}_{25}+40\mathcal{O}_{26}+40\mathcal{O}_{24}
\right]
\]
and thus
\[
||e^{-\hat{V}_{\Delta'}(\phi)+\delta V_{\Delta'}(\lambda,\phi)}||_{\partial\phi,\phi,h}\le
2 e^{-\frac{\bar{g}}{16}\phi^4}
\exp\left\{25L^5\rho
\left[
32\mathcal{O}_{25}+40\mathcal{O}_{26}+40\mathcal{O}_{24}
\right](1+\bar{g}\phi^4)\right\}\ .
\]
The hypothesis and the chosen definition of $\mathcal{O}_{27}$
implies
\[
\exp\left\{25L^5\rho
\left[
32\mathcal{O}_{25}+40\mathcal{O}_{26}+40\mathcal{O}_{24}
\right](1+\bar{g}\phi^4)\right\}\le \exp\left\{
\frac{1}{16}(1+\bar{g}\phi^4)
\right\}
\]
which gives the desired bound
\[
||e^{-\hat{V}_{\Delta'}(\phi)+\delta V_{\Delta'}(\lambda,\phi)}||_{\partial\phi,\phi,h}\le
2 e^{\frac{1}{16}}<3\ .
\]
\qed

\begin{Lemma}\label{L36lem}
Let 
\[
\mathcal{O}_{28}=200\times\left\{\log\left(\frac{3}{2}\right)\right\}^{-1}\times
[4\mathcal{O}_{25}+5\mathcal{O}_{27}+5\mathcal{O}_{24}].
\]
Provided $\lambda$ satisfies $\mathcal{O}_{28}L^{9}|\lambda|\bar{g}^{\frac{11}{12}-\frac{1}{3}\eta}\le 1$
we have, for all $\Delta'\in\mathbb{L}$,
\[
|e^{-\hat{V}_{\Delta'}(\phi)+\delta V_{\Delta'}(\lambda,\phi)}|_{\partial\phi,h_\ast}\le 3\ .
\]
\end{Lemma}

\noindent{\bf Proof:}
Again by the multiplicative property of the seminorms and by Lemma \ref{L10lem} we have
\[
|e^{-\hat{V}_{\Delta'}(\phi)+\delta V_{\Delta'}(\lambda,\phi)}|_{\partial\phi,h_\ast}
\le 2 \times \exp\left[ |\delta V_{\Delta'}(\lambda,\phi)|_{\partial\phi,h_\ast}
\right]\ .
\]
We also have
\[
|\delta V_{\Delta'}(\lambda,\phi)|_{\partial\phi,h_\ast}\le
\sum_{\substack{0\le k\le 4 \\ 1\le j\le 3}} |\lambda|^j\ 
|\delta\beta_{k,j,\Delta'}|\ |:\phi^k:_{C_0}|_{\partial\phi,h_\ast}\ .
\]
For $0\le k\le 4$,
\begin{eqnarray*}
|:\phi^k:_{C_0}|_{\partial\phi,h_\ast} & \le & 25\times\max_{0\le a\le k} |\phi^a|_{\partial\phi,h_\ast} \\
 & \le & 25\times\max_{0\le a\le k} h_{\ast}^{a} \\
 & \le & 25 h_{\ast}^{k}
\end{eqnarray*}
since $h_{\ast}\ge 1$.
Using Lemmas \ref{L33lem}, \ref{L34lem} and \ref{L32lem} we then immediately get
\begin{eqnarray*}
|\delta V_{\Delta'}(\lambda,\phi)|_{\partial\phi,h_\ast} & \le & 25\left\{
|\lambda|\times\left(\sum_{k=0}^3 \mathcal{O}_{25} \bar{g}^{1-\eta} L^{\frac{5}{2}} h_{\ast}^{k}\right)
\right. \\
 & & +|\lambda|^2\times\left(\sum_{k=0}^4 \mathcal{O}_{26} \bar{g}^{2-2\eta} L^{5} h_{\ast}^{k}\right) \\
 & & \left. +|\lambda|^3\times\left(\sum_{k=0}^4 \mathcal{O}_{24} \bar{g}^{\frac{11}{4}-\eta_R} L^{3} h_{\ast}^{k}\right)
 \right\}\ .
\end{eqnarray*}
We bound powers of $h_\ast$ simply by $h_{\ast}^4=2^3 L^{\frac{3+\epsilon}{4}}\le 8L^4$ and $L^{\frac{5}{2}}$ by $L^3$.
We thus easily get
\[
|\delta V_{\Delta'}(\lambda,\phi)|_{\partial\phi,h_\ast} \le
25\times 8\times L^9 \left\{
4 |\lambda| \mathcal{O}_{25} \bar{g}^{1-\eta}
+ 5 |\lambda|^2 \mathcal{O}_{26} \bar{g}^{2-2\eta}
+ 5 |\lambda|^23 \mathcal{O}_{24} \bar{g}^{\frac{11}{4}-\eta_R}
\right\}\ .
\]
Let this time $\rho=|\lambda|\bar{g}^{\frac{11}{12}-\frac{1}{3}\eta_R}$. Since $\log\left(\frac{3}{2}\right)\simeq 0.405$...
 and, e.g., $\mathcal{O}_{25}=\frac{27}{2}$, it is clear that $\mathcal{O}_{28}>1$.
Hence the hypothesis of the lemma implies $\rho\le 1$. Besides, the standard hypothesis $\eta_R\ge 3\eta$
implies $\bar{g}^{1-\eta}\le \bar{g}^{\frac{11}{12}-\frac{1}{3}\eta_R}$. As a consequence we have
\[
\delta V_{\Delta'}(\lambda,\phi)|_{\partial\phi,h_\ast} \le
200\times L^9 \rho\left\{
4\mathcal{O}_{25} + 5\mathcal{O}_{26}+5\mathcal{O}_{24}\right\}
\]
from which the result follows easily.
\qed

\begin{Lemma}\label{L37lem}
Let 
\[
\mathcal{O}_{29}=200\times\left\{\log\left(\frac{3}{2}\right)\right\}^{-1}\times
[4\mathcal{O}_{25}+5\mathcal{O}_{27}].
\]
Under the extra assumption that $R=0$ and
provided $\lambda$ satisfies $\mathcal{O}_{29}L^{9}|\lambda|\bar{g}^{1-\eta}\le 1$
we have, for all $\Delta'\in\mathbb{L}$,
\[
|e^{-\hat{V}_{\Delta'}(\phi)+\delta V_{\Delta'}(\lambda,\phi)}|_{\partial\phi,h_\ast}\le 3\ .
\]
\end{Lemma}

\noindent{\bf Proof:}
The proof is the same as that of the previous lemma except for the absence of the $\beta_{k,3,\Delta'}$ terms.
The only modification is to let $\rho=|\lambda|\bar{g}^{1-\eta}$ instead of $|\lambda|\bar{g}^{\frac{11}{12}-\frac{1}{3}\eta_R}$.
\qed

\begin{Lemma}\label{L38lem}
Let $\mathcal{O}_{30}=\mathcal{O}_{25}+\mathcal{O}_{26}+\mathcal{O}_{24}$.
Provided $\lambda$ satisfies $\mathcal{O}_{30}L^{5}|\lambda|\bar{g}^{\frac{11}{12}-\frac{1}{3}\eta_R}\le 1$
we have, for all $\Delta'\in\mathbb{L}$,
$|\delta b_{\Delta'}(\lambda)|\le 1$.
\end{Lemma}

\noindent{\bf Proof:}
By definition
\[
\delta b_{\Delta'}(\lambda)=\delta\beta_{0,\Delta'}(\lambda)=
\lambda \delta\beta_{0,1,\Delta'}+\lambda^2 \delta\beta_{0,2,\Delta'}+
\lambda^3 \delta\beta_{0,3,\Delta'}\ .
\]
From Lemmas \ref{L33lem}, \ref{L34lem} and \ref{L32lem} we get
\begin{eqnarray*}
|\delta\beta_{0,1,\Delta'}| & \le & \mathcal{O}_{25} L^{\frac{5}{2}} \bar{g}^{1-\eta} \\
|\delta\beta_{0,2,\Delta'}| & \le & \mathcal{O}_{26} L^{5} \bar{g}^{2-2\eta} \\
|\delta\beta_{0,3,\Delta'}| & \le & \mathcal{O}_{24} L^{3} \bar{g}^{\frac{11}{4}-\eta_R} 
\end{eqnarray*}
which give
\[
|\delta b_{\Delta'}(\lambda)|\le L^5\left[
\mathcal{O}_{25}|\lambda|\bar{g}^{1-\eta}
+\mathcal{O}_{26}|\lambda|^2\bar{g}^{2-2\eta}
+\mathcal{O}_{24}|\lambda|^3 \bar{g}^{\frac{11}{4}-\eta_R} 
\right]\ .
\]
Since clearly $\mathcal{O}_{30}>1$, one can conclude as we did previously
that
\[
|\delta b_{\Delta'}(\lambda)|\le\mathcal{O}_{30}L^{5}|\lambda|\bar{g}^{\frac{11}{12}-\frac{1}{3}\eta_R}\le 1\ .
\]
\qed

\begin{Lemma}\label{L39lem}
Let $\mathcal{O}_{31}=\mathcal{O}_{25}+\mathcal{O}_{26}$.
Under the extra assumption that $R=0$ and
provided $\lambda$ satisfies $\mathcal{O}_{31}L^{5}|\lambda|\bar{g}^{1-\eta}\le 1$
we have, for all $\Delta'\in\mathbb{L}$,
$|\delta b_{\Delta'}(\lambda)|\le 1$.
\end{Lemma}

\noindent{\bf Proof:}
The proof is the same as that of the previous lemma, without the $\delta\beta_{0,3,\Delta'}$ term.
\qed

\begin{Lemma}\label{L40lem}
Let
\[
\mathcal{O}_{32}=\max(2\mathcal{O}_{19},\mathcal{O}_{27}, \mathcal{O}_{30})\ \ {\rm and}\ \ 
\mathcal{O}_{33}=7\times\exp\left(1+\frac{\sqrt{2}}{2}\right)\ .
\]
Provided $\lambda$ satisfies $\mathcal{O}_{32}L^{15}|\lambda|\bar{g}^{\frac{1}{4}-\frac{1}{3}\eta_R}\le 1$
we have, for all $\Delta'\in\mathbb{L}$ and $\phi\in\mathbb{R}$,
\[
||K'_{\Delta'}(\lambda,\phi)||_{\partial\phi,\phi,h}\le \mathcal{O}_{33}\ .
\]
\end{Lemma}

\noindent{\bf Proof:}
One can rewrite the definition of $K'_{\Delta'}(\lambda,\phi)$ as
\[
K'_{\Delta'}(\lambda,\phi)=
e^{-\delta b_{\Delta'}(\lambda)}
e^{-\frac{1}{2}(f,\Gamma f)_{L^{-1}\Delta'}}
\]
\[
\times\left\{
\hat{K}_{\Delta'}(\lambda,\phi)
-e^{\frac{1}{2}(f,\Gamma f)_{L^{-1}\Delta'}}\left(
e^{-\hat{V}_{\Delta'}(\phi)+\delta V_{\Delta'}(\lambda,\phi)}
-e^{-\hat{V}_{\Delta'}(\phi)}\right)\right\}
\]
from which one deduces
\[
||K'_{\Delta'}(\lambda,\phi)||_{\partial\phi,\phi,h}\le e^{|\delta b_{\Delta'}(\lambda)|}
\times \exp\left[2^{-\frac{3}{2}}\right]
\]
\[
\times\left\{
||\hat{K}_{\Delta'}(\lambda,\phi)||_{\partial\phi,\phi,h}+
\exp\left[2^{-\frac{3}{2}}\right]\times
\left(||e^{-\hat{V}_{\Delta'}(\phi)+\delta V_{\Delta'}(\lambda,\phi)}||_{\partial\phi,\phi,h}
+||e^{-\hat{V}_{\Delta'}(\phi)}||_{\partial\phi,\phi,h}\right)\right\}\ .
\]
Indeed, we previously showed $|(f,\Gamma f)_{L^{-1}\Delta'}|\le\frac{1}{\sqrt{2}}$.
We have $||e^{-\hat{V}_{\Delta'}(\phi)}||_{\partial\phi,\phi,h}\le 2$
by Lemma \ref{L10lem}.
Clearly $\mathcal{O}_{32}\ge \mathcal{O}_{30}\ge \mathcal{O}_{25}>1$ and therefore the hypothesis of the present lemma
implies that of Lemma \ref{L27lem}.
The latter gives the bound
\[
||\hat{K}_{\Delta'}(\lambda,\phi)||_{\partial\phi,\phi,h}\le 2 \exp\left[2^{-\frac{3}{2}}\right]\times
\sum_{n=1}^{\infty}\left(\mathcal{O}_{19}L^{15}|\lambda|\bar{g}^{\frac{1}{4}-\frac{1}{3}\eta_R}\right)^n
\]
where we used $|(\Re f,\Gamma \Re f)_{L^{-1}\Delta'}|\le\frac{1}{\sqrt{2}}$.
Since $0\le x\le \frac{1}{2}$ implies $\sum_{n=1}^{\infty} x^n\le 1$ and since the hypothesis implies
$\mathcal{O}_{19}L^{15}|\lambda|\bar{g}^{\frac{1}{4}-\frac{1}{3}\eta_R}\le\frac{1}{2}$, we have the simpler estimate
\[
||\hat{K}_{\Delta'}(\lambda,\phi)||_{\partial\phi,\phi,h}\le 2 \exp\left[2^{-\frac{3}{2}}\right]\ .
\]
From $\eta_R\ge 3\eta$ we get $\bar{g}^{\frac{1}{4}-\eta}\le\bar{g}^{\frac{1}{4}-\frac{1}{3}\eta_R}$. Since also $L^5<L^{15}$,
the hypothesis of the  present lemma implies that of Lemma \ref{L35lem}
which gives us $|\delta b_{\Delta'}(\lambda)|\le 1$.
Finally, since $\bar{g}^{\frac{11}{12}-\frac{1}{3}\eta_R}\le\bar{g}^{\frac{1}{4}-\frac{1}{3}\eta_R}$,
the hypothesis of Lemma \ref{L38lem} is satisfied. This gives us the last needed ingredient
\[
||e^{-\hat{V}_{\Delta'}(\phi)+\delta V_{\Delta'}(\lambda,\phi)}||_{\partial\phi,\phi,h}\le 3\ .
\]
Altogether we obtain
\[
||K'_{\Delta'}(\lambda,\phi)||_{\partial\phi,\phi,h}\le e\times \exp\left[2^{-\frac{3}{2}}\right]
\times\left\{2 \exp\left[2^{-\frac{3}{2}}\right] +\exp\left[2^{-\frac{3}{2}}\right](3+2)
\right\}=\mathcal{O}_{33}\ .
\]
\qed

\begin{Lemma}\label{L41lem}
Let
\[
\mathcal{O}_{34}=\max(2\mathcal{O}_{20},\mathcal{O}_{27}, \mathcal{O}_{31})\ .
\]
Under the extra assumption that $R=0$ and
provided $\lambda$ satisfies $\mathcal{O}_{34}L^{15}|\lambda|\bar{g}^{\frac{1}{4}-\eta}\le 1$
we have, for all $\Delta'\in\mathbb{L}$ and $\phi\in\mathbb{R}$,
\[
||K'_{\Delta'}(\lambda,\phi)||_{\partial\phi,\phi,h}\le \mathcal{O}_{33}\ .
\]
\end{Lemma}

\noindent{\bf Proof:}
The proof is similar to that of the last lemma.
The only modifications are as follows. We use Lemma \ref{L28lem} instead of Lemma \ref{L27lem},
noting that
$\mathcal{O}_{34}\ge \mathcal{O}_{31}\le \mathcal{O}_{25}>1$.
We use Lemma \ref{L39lem} instead of \ref{L38lem}.
\qed

\begin{Lemma}\label{L42lem}
Let
\[
\mathcal{O}_{35}=\max(2\mathcal{O}_{21},\mathcal{O}_{28}, \mathcal{O}_{30})\ .
\]
Provided $\lambda$ satisfies $\mathcal{O}_{35}L^{15}|\lambda|\bar{g}^{\frac{11}{12}-\frac{1}{3}\eta_R}\le 1$
we have, for all $\Delta'\in\mathbb{L}$ and $\phi\in\mathbb{R}$,
\[
|K'_{\Delta'}(\lambda,\phi)|_{\partial\phi,h_\ast}\le \mathcal{O}_{33}\ .
\]
\end{Lemma}

\noindent{\bf Proof:}
The proof is similar to that of Lemma \ref{L40lem}.
The only modifications are as follows. We use Lemma \ref{L29lem} instead of Lemma \ref{L27lem},
noting that
$\mathcal{O}_{35}\ge \mathcal{O}_{30}\le \mathcal{O}_{25}>1$.
We use Lemma \ref{L36lem} instead of \ref{L35lem}.
\qed

\begin{Lemma}\label{L43lem}
Let
\[
\mathcal{O}_{36}=\max(2\mathcal{O}_{22},\mathcal{O}_{29}, \mathcal{O}_{31})\ .
\]
Under the extra assumption that $R=0$ and
provided $\lambda$ satisfies $\mathcal{O}_{36}L^{15}|\lambda|\bar{g}^{1-\eta}\le 1$
we have, for all $\Delta'\in\mathbb{L}$ and $\phi\in\mathbb{R}$,
\[
|K'_{\Delta'}(\lambda,\phi)|_{\partial\phi,h_\ast}\le \mathcal{O}_{33}\ .
\]
\end{Lemma}

\noindent{\bf Proof:}
The proof is similar to that of Lemma \ref{L40lem}.
The only modifications are as follows. We use Lemma \ref{L30lem} instead of Lemma \ref{L27lem},
noting that
$\mathcal{O}_{36}\ge \mathcal{O}_{31}\le \mathcal{O}_{25}>1$.
We use Lemma \ref{L37lem} instead of \ref{L35lem}
and Lemma \ref{L39lem} instead of \ref{L38lem}.
\qed

Recall from \S\ref{algdefsec} that 
\begin{equation}
\xi_{R,\Delta'}(\vec{V})(\phi)=\xi_{R,\Delta'}^{{\rm main}}(\vec{V})(\phi)
+\xi_{R,\Delta'}^{{\rm higher}}(\vec{V})(\phi)+\xi_{R,\Delta'}^{{\rm shift}}(\vec{V})(\phi)
\label{xideceq}
\end{equation}
where
\[
\xi_{R,\Delta'}^{{\rm main}}(\vec{V})(\phi)=\frac{1}{2\pi i}\oint_{\gamma_0} \frac{{\rm d}\lambda}{\lambda^4}
\left. K'_{\Delta'}(\lambda,\phi) \right|_{R=0}\ ,
\]
\[
\xi_{R,\Delta'}^{{\rm higher}}(\vec{V})(\phi)=
\frac{1}{2\pi i}\oint_{\gamma_{01}} \frac{{\rm d}\lambda}{\lambda^4(\lambda-1)}
K'_{\Delta'}(\lambda,\phi)
\]
and
\[
\xi_{R,\Delta'}^{{\rm shift}}(\vec{V})(\phi)=
\left(e^{-\hat{V}_{\Delta'}(\phi)}-e^{-V'_{\Delta'}(\phi)}\right) Q'_{\Delta'}(\phi)\ .
\]
The next few lemmas will provide bounds for each of these terms.

\begin{Lemma}\label{L44lem}
For all unit cube $\Delta'$ we have that
\[
|||\xi_{R,\Delta'}^{{\rm main}}(\vec{V})|||_{\bar{g}}
\le \mathcal{O}_{37} L^{45} \bar{g}^{\frac{11}{4}-3\eta}
\]
where $\mathcal{O}_{37}=\mathcal{O}_{33}\times\max[\mathcal{O}_{34}^3,\mathcal{O}_{36}^3]$.
\end{Lemma}

\noindent{\bf Proof:}
We use the freedom to deform the contour of integration in order to pick for $\gamma_0$
the circle of radius $\rho$ around the origin where
\[
\rho=\left(\mathcal{O}_{34} L^{15} \bar{g}^{\frac{1}{4}-\eta}\right)^{-1}>0\ .
\]
We then use
\[
\left|\left|
\frac{1}{2\pi i}\oint_{\gamma_0} \frac{{\rm d}\lambda}{\lambda^4}
\left. K'_{\Delta'}(\lambda,\phi) \right|_{R=0}
\right|\right|_{\partial\phi,\phi,h}\le
\rho^{-3}\sup_{\lambda\in\gamma_0} \left|\left|
\left. K'_{\Delta'}(\lambda,\phi) \right|_{R=0}\right|\right|_{\partial\phi,\phi,h}
\]
\[
\le \mathcal{O}_{33} \rho^{-3}
\]
by Lemma \ref{L41lem}.
Hence
\[
||\xi_{R,\Delta'}^{{\rm main}}(\vec{V})(\phi)||_{\partial\phi,\phi,h}\le
\mathcal{O}_{33} \mathcal{O}_{34}^3 L^{45} \bar{g}^{\frac{3}{4}-3\eta}\ .
\]
The bound on $|\xi_{R,\Delta'}^{{\rm main}}(\vec{V})(\phi)||_{\partial\phi,h_\ast}$
is derived in the same manner using Lemma \ref{L43lem} and setting
\[
\rho=\left(\mathcal{O}_{36} L^{15} \bar{g}^{1-\eta}\right)^{-1}>0
\]
for the contour radius.
We get
\[
|\xi_{R,\Delta'}^{{\rm main}}(\vec{V})(\phi)||_{\partial\phi,h_\ast}\le
\mathcal{O}_{33} \mathcal{O}_{36}^3 L^{45} \bar{g}^{3-3\eta}
\]
and therefore
\[
|||\xi_{R,\Delta'}^{{\rm main}}(\vec{V})|||_{\bar{g}}
\le \mathcal{O}_{33} \times L^{45}\times
\max \left[\mathcal{O}_{36}^3 \bar{g}^{3-3\eta},
\mathcal{O}_{36}^3  \bar{g}^{\frac{11}{4}-3\eta}
\right]
\]
by definition of the $|||\cdot|||_{\bar{g}}$ norm. Since $\bar{g}\le 1$ the lemma follows.
\qed

\begin{Lemma}\label{L45lem}
For all unit cube $\Delta'$ we have that
\[
|||\xi_{R,\Delta'}^{{\rm higher}}(\vec{V})|||_{\bar{g}}
\le \mathcal{O}_{38} L^{60} \bar{g}^{\frac{11}{4}-3\eta}
\]
where $\mathcal{O}_{38}=2\mathcal{O}_{33}\times\max[\mathcal{O}_{32}^4,\mathcal{O}_{35}^4]$.
\end{Lemma}

\noindent{\bf Proof:}
We proceed as in the proof of the previous lemma. We take for the contour $\gamma_{01}$ a circle of radius $\gamma\ ge 2$
around the origin. This ensures that both $0$ and $1$ are enclosed by the contour and allows one to bound  
$\frac{1}{|\lambda -1|}$ 
by $\frac{2}{\rho}$.
We first take
\[
\rho=\left(\mathcal{O}_{32} L^{15} \bar{g}^{\frac{1}{4}-\frac{1}{3}\eta_R}\right)^{-1}\ge 2
\]
because of standard hypothesis (\ref{extrastdeq}). Lemma \ref{L40lem} then results in the bound
\[
||\xi_{R,\Delta'}^{{\rm higher}}(\vec{V})(\phi)||_{\partial\phi,\phi,h}\le
\mathcal{O}_{33} \times\frac{2}{\rho^4}
= 2\mathcal{O}_{33} \times 
\mathcal{O}_{32}^4 L^{60} \bar{g}^{1-\frac{4}{3}\eta_R}\ .
\]
Likewise, if we pick
\[
\rho=\left(\mathcal{O}_{35} L^{15} \bar{g}^{\frac{11}{12}-\frac{1}{3}\eta_R}\right)^{-1}
\]
then the latter is at least equal to $2$ by the standard hypothesis (\ref{extrastdeq})
Therefore \ref{L42lem} results in the bound
\[
|\xi_{R,\Delta'}^{{\rm higher}}(\vec{V})(\phi)|_{\partial\phi,h_\ast}\le
2\mathcal{O}_{33} \times 
\mathcal{O}_{35}^4 L^{60} \bar{g}^{\frac{11}{3}-\frac{4}{3}\eta_R}\ .
\]
Finally, we get
\[
|||\xi_{R,\Delta'}^{{\rm higher}}(\vec{V})|||_{\bar{g}}\le 2\mathcal{O}_{33} \times L^{60}\times
\max\left[
\mathcal{O}_{35}^4 \bar{g}^{\frac{11}{3}-\frac{4}{3}\eta_R},
\mathcal{O}_{32}^4 \bar{g}^{3-\frac{4}{3}\eta_R}
\right] \ .
\]
The conclusion of the proof is a matter of showing that the powers of $\bar{g}$
involved are bounded by $\bar{g}^{\frac{11}{4}-3\eta}$.
In other words, one needs to check the two inequalities
\[
\frac{11}{3}-\frac{4}{3}\eta_R\ge \frac{11}{4}-3\eta
\]
and
\[
3-\frac{4}{3}\eta_R\ge \frac{11}{4}-3\eta\ .
\]
The first inequality follows from the second since $\frac{11}{3}>3$.
Finally, the second inequality is equivalent to the standard hypothesis (\ref{etar2stdeq}) and therefore holds.
\qed

\begin{Lemma}\label{L46lem}
For all unit cube $\Delta'$ and all $\phi\in\mathbb{R}$ we have that
\[
||\xi_{R,\Delta'}^{{\rm shift}}(\vec{V})(\phi)||_{\partial\phi,\phi,h}
\le \mathcal{O}_{39} L^{\frac{15}{2}} \bar{g}^{\frac{3}{4}-3\eta}
\]
where 
\[
\mathcal{O}_{39}=2^{15}\times 3^{\frac{5}{2}}\times 11\times 103\times(5+\sqrt{2})\times e^{\frac{1}{48}}\times\mathcal{O}_{27}
\times \max_{0\le n\le 6}\left(\frac{n}{4e}\right)^{\frac{n}{4}}\ .
\]
\end{Lemma}

\noindent{\bf Proof:}
Recall that
\[
\hat{V}_{\Delta'}(\phi)=\sum_{k=1}^4
\hat{\beta}_{k,\Delta'} :\phi^k:_{C_0}
\]
while
\[
V'_{\Delta'}(\phi)=\hat{V}_{\Delta'}(\phi)-\delta V_{\Delta'}(\phi)+\delta b_{\Delta'}
\]
with
\[
\delta V_{\Delta'}(\phi)=\left.\delta V_{\Delta'}(\lambda,\phi)\right|_{\lambda=1}\qquad{\rm and}\qquad
\delta b_{\Delta'}=\left.\delta b_{\Delta'}(\lambda)\right|_{\lambda=1}\ .
\]
Let us introduce the notation
\[
U_{\Delta'}(\phi)=\hat{V}_{\Delta'}(\phi)-V'_{\Delta'}(\phi)=
\delta V_{\Delta'}(\phi)-\delta b_{\Delta'}
\]
so that
\[
U_{\Delta'}(\phi)=\sum_{k=1}^4
\delta\beta_{k,\Delta'} :\phi^k:_{C_0}=\sum_{k=1}^4\left(
\sum_{j=1}^{3}
\delta\beta_{k,j,\Delta'}\right) :\phi^k:_{C_0}\ .
\]
The latter quantity is the same as $\delta V_{\Delta'}(\lambda,\phi)$ that was estimated in Lemma \ref{L35lem}
except that the sum over $k$ goes from $1$ to $4$ instead of from $0$ to $4$, and $\lambda$ is now set equal to $1$.
By the same argument as in Lemma \ref{L35lem} we therefore get
\[
||U_{\Delta'}(\lambda,\phi)||_{\partial\phi,\phi,h}\le 25 (1+\bar{g}\phi^4)
L^5\rho
\times\left[
32\mathcal{O}_{25}+40\mathcal{O}_{26}+40\mathcal{O}_{24}
\right]
\]
with $\rho=\bar{g}^{\frac{1}{4}-\eta}\le 1$ by the standard hypothesis (\ref{etastdeq}).
Note that for simplicity we did not take advantage of the smaller range of summation
for $k$ in order to improve the coefficients $32$ and $40$ in the last bound.

We now write
\begin{eqnarray*}
\xi_{R,\Delta'}^{{\rm shift}}(\vec{V}) & = & \left(e^{-\hat{V}_{\Delta'}(\phi)}-
e^{-\hat{V}_{\Delta'}(\phi)+U_{\Delta'}(\phi)}\right) Q'_{\Delta'}(\phi) \\
 & = & -e^{-\hat{V}_{\Delta'}(\phi)}\left(e^{U_{\Delta'}(\phi)-1}\right) Q'_{\Delta'}(\phi) \\
 & = & -e^{-\hat{V}_{\Delta'}(\phi)} Q'_{\Delta'}(\phi)
 \int_0^1 {\rm d}s\ U_{\Delta'}(\phi)\ e^{s U_{\Delta'}(\phi)}
\end{eqnarray*}
which implies the bound
\[
||\xi_{R,\Delta'}^{{\rm shift}}(\vec{V})||_{\partial\phi,\phi,h}
\le ||e^{-\hat{V}_{\Delta'}(\phi)}||_{\partial\phi,\phi,h}
||Q'_{\Delta'}(\phi)||_{\partial\phi,\phi,h}
\times
\int_0^1 {\rm d}s\ ||U_{\Delta'}(\phi)||_{\partial\phi,\phi,h}
\times \exp\left\{s ||U_{\Delta'}(\phi)||_{\partial\phi,\phi,h}\right\}
\]
\begin{equation}
\le ||e^{-\hat{V}_{\Delta'}(\phi)}||_{\partial\phi,\phi,h}
\ ||Q'_{\Delta'}(\phi)||_{\partial\phi,\phi,h}
\ ||U_{\Delta'}(\phi)||_{\partial\phi,\phi,h}
\ e^{ ||U_{\Delta'}(\phi)||_{\partial\phi,\phi,h}}\ .
\label{shiftbdeq}
\end{equation}
Each of these four factors needs to be estimated separately.
By Lemma \ref{L10lem} we have
\begin{equation}
||e^{-\hat{V}_{\Delta'}(\phi)}||_{\partial\phi,\phi,h}\le
2 e^{-\frac{\bar{g}}{16}\phi^4}\ .
\label{shift1bdeq}
\end{equation}
The last exponential expression will be needed in order to control each of the other three
factors in (\ref{shiftbdeq}). Indeed, we will show that the latter can be bounded using the exponential
of $\frac{1}{3}\times\frac{\bar{g}}{16}\phi^4=\frac{\bar{g}}{48}\phi^4$.

By the previous considerations,
\begin{equation}
||U_{\Delta'}(\phi)||_{\partial\phi,\phi,h}\le
25 (1+\bar{g}\phi^4)
L^5 \bar{g}^{\frac{1}{4}-\eta}
\times\left[
32\mathcal{O}_{25}+40\mathcal{O}_{26}+40\mathcal{O}_{24}
\right]\ .
\label{shift2bdeq}
\end{equation}
Since
\[
48\times 25 
\times\left[
32\mathcal{O}_{25}+40\mathcal{O}_{26}+40\mathcal{O}_{24}
\right]=3\times \mathcal{O}_{27}
\]
the standard hypothesis (\ref{extrastdeq}) ensures that
\[
||U_{\Delta'}(\phi)||_{\partial\phi,\phi,h}\le
\frac{1}{48} (1+\bar{g}\phi^4)
\]
and therefore
\begin{equation}
e^{ ||U_{\Delta'}(\phi)||_{\partial\phi,\phi,h}}\le
e^{\frac{1}{48}}\times e^{\frac{\bar{g}}{48}\phi^4}
\label{shift3bdeq}
\end{equation}
which is the first estimate of the kind we are seeking.

Next we note that (\ref{shift2bdeq}) can be rewritten as
\begin{eqnarray}
||U_{\Delta'}(\phi)||_{\partial\phi,\phi,h} & \le &
\frac{1}{16} \mathcal{O}_{27} 
L^5 \bar{g}^{\frac{1}{4}-\eta}(1+\bar{g}\phi^4) \nonumber \\
 & \le & \frac{1}{16} \mathcal{O}_{27} 
L^5 \bar{g}^{\frac{1}{4}-\eta}(48+\bar{g}\phi^4) \nonumber \\
 & \le & 3\mathcal{O}_{27} 
L^5 \bar{g}^{\frac{1}{4}-\eta}(1+\frac{\bar{g}}{48}\phi^4) \nonumber \\
 & \le & 3\mathcal{O}_{27} 
L^5 \bar{g}^{\frac{1}{4}-\eta}e^{\frac{\bar{g}}{48}\phi^4}\ .
\label{shift4bdeq}
\end{eqnarray} 

Finally we need a similar bound on the $Q'_{\Delta'}$ factor.
Recall that
\[
Q'_{\Delta'}(\phi)=W'_{5,\Delta'}:\phi_{\Delta'}^{5}:_{C_0}+
W'_{6,\Delta'}:\phi_{\Delta'}^{6}:_{C_0}\ .
\]
Similarly to the proof of Lemma \ref{L17lem} by undoing the Wick ordering
we have
\[
||Q'_{\Delta'}(\phi)||_{\partial\phi,\phi,h}  \le
|W'_{5,\Delta'}|\times 81\times \left(\max_{0\le a\le 5} ||\phi^a||_{\partial\phi,\phi,h} \right)
\]
\[
+|W'_{5,\Delta'}|\times 331\times \left(\max_{0\le a\le 6} ||\phi^a||_{\partial\phi,\phi,h} \right)
\]
and thus
\[
||Q'_{\Delta'}(\phi)||_{\partial\phi,\phi,h}  \le 412\times
\max\left[|W'_{5,\Delta'}|,|W'_{6,\Delta'}|\right]\times \max_{0\le a\le 6} ||\phi^a||_{\partial\phi,\phi,h}\ .
\]
Now for $0\le k\le 6$ and for $\gamma>0$ we have 
\begin{eqnarray*}
||\phi^k||_{\partial\phi,\phi,h} & = & (h+|\phi|)^k \\
 & \le & \sum_{n=0}^{k} \left(
\begin{array}{c} k \\ n \end{array} 
\right) \times \left(c_1\bar{g}^{-\frac{1}{4}}\right)^{k-n}\times
\left(\frac{n}{4e}\right)^{\frac{n}{4}}\times \left(\gamma\bar{g}\right)^{-\frac{n}{4}}
e^{\gamma\bar{g}\phi^4}
\end{eqnarray*}
by the first inequaity in Lemma \ref{L4lem} with $\beta_4=\bar{g}$.
If also $\gamma\le 1$, then we bound $\gamma^{-\frac{n}{4}}$ by $\gamma^{-\frac{k}{4}}$
and get
\[
||\phi^k||_{\partial\phi,\phi,h} \le 
\left(\max_{0\le n\le 6}\left(\frac{n}{4e}\right)^{\frac{n}{4}}\right)\times
\bar{g}^{-\frac{k}{4}}\gamma^{-\frac{k}{4}} e^{\gamma\bar{g}\phi^4}\times(1+c_1)^k\ .
\]
We now pick $\gamma=\frac{1}{48}$ and simply bound $(1+c_1)^k$ by $2^k\le 2^6$.
One then easily obtains
\[
||Q'_{\Delta'}(\phi)||_{\partial\phi,\phi,h}  \le 412\times
\max\left[|W'_{5,\Delta'}|,|W'_{6,\Delta'}|\right]\times
\left(\max_{0\le n\le 6}\left(\frac{n}{4e}\right)^{\frac{n}{4}}\right)\times
\bar{g}^{-\frac{3}{2}} 48^{\frac{3}{2}}\times 2^6\times
e^{\frac{\bar{g}}{48}\phi^4}\ .
\]
In order to continue we need to now bound the $W'$ factors.
From the definition in \S\ref{algdefsec} we immediately get
\[
|W'_{6,\Delta'}|\le L^{3-6[\phi]}\bar{g}^{2-2\eta}
+8 L^{-6[\phi]} L^3 \left(\frac{3}{2}\bar{g}\right)^2\times||\Gamma||_{L^1}
\]
Since $\eta\ge 0$, we get from Corollary \ref{gamL1cor}
\[
|W'_{6,\Delta'}|\le L^{3-6[\phi]}\bar{g}^{2-2\eta}
\left(1+8\left(\frac{3}{2}\right)^2\times \frac{1}{\sqrt{2}} L^{3-2[\phi]}\right)\ ,
\]
namely,
\[
|W'_{6,\Delta'}|\le L^{3-6[\phi]}\bar{g}^{2-2\eta} (1+9\sqrt{2} L^{3-2[\phi]})
\le L^{6-8[\phi]}\bar{g}^{2-2\eta} (1+9\sqrt{2})\ .
\]
Now from the assumption $\epsilon\le 1$ we get $[\phi]=\frac{3-\epsilon}{4}\ge \frac{1}{2}$
and therefore $6-8[\phi]\le 2$.
As a result we simplify the last bound on $W'_6$ into
\[
|W'_{6,\Delta'}|\le L^2 \bar{g}^{2-2\eta} (1+9\sqrt{2})\ .
\]
Similarly, from the definition in \S\ref{algdefsec}, Corollary \ref{gamL1cor}
and the standard hypotheses (\ref{betastdeq}) and (\ref{fwstdeq}), we easily get
\[
|W'_{5,\Delta'}|\le L^{3-5[\phi]}\bar{g}^{2-2\eta}+
6L^{-5[\phi]}\bar{g}^{2-2\eta}L^3 \frac{1}{\sqrt{2}} L^{3-2[\phi]} L^{-(3-[\phi])}
\]
\[
+12  L^{-5[\phi]}\frac{3}{2}\bar{g}L^3 \frac{1}{\sqrt{2}} L^{3-2[\phi]} \bar{g}^{1-\eta}
+48 L^{-5[\phi]}\left(\frac{3}{2}\right)^2 \bar{g}^2 L^3 \frac{1}{2} L^{2(3-2[\phi])} L^{-(3-[\phi])}
\]
\[
\le L^{3-5[\phi]}\bar{g}^{2-2\eta}\left[
1+3\sqrt{2} L^{-[\phi]}+9\sqrt{2}L^{3-2[\phi]}+54 L^{3-3[\phi]}
\right]
\]
\[
\le L^{6-7[\phi]}\bar{g}^{2-2\eta}\left[
1+3\sqrt{2}+9\sqrt{2}+54
\right]
\]
since $L^{-[\phi]}$ and $L^{3-3[\phi]}$ are bounded by $L^{3-2[\phi]}$.
Since $\epsilon\le 1$ implies $6-7[\phi]\le \frac{5}{2}$ we then have
\[
|W'_{5,\Delta'}|\le L^{\frac{5}{2}}\bar{g}^{2-2\eta}
\times 11(5+\sqrt{2})
\]
which compared with the previous estimate on $W'_{6}$
results in
\[
\max\left[|W'_{5,\Delta'}|,|W'_{6,\Delta'}|\right]\le
L^{\frac{5}{2}}\bar{g}^{2-2\eta}
\times 11(5+\sqrt{2})\ .
\]
As a consequence we have
\begin{equation}
||Q'_{\Delta'}(\phi)||_{\partial\phi,\phi,h}  \le 412\times
L^{\frac{5}{2}}\bar{g}^{2-2\eta}
\times 11(5+\sqrt{2})
\times
\left(\max_{0\le n\le 6}\left(\frac{n}{4e}\right)^{\frac{n}{4}}\right)\times
\bar{g}^{-\frac{3}{2}} 48^{\frac{3}{2}}\times 2^6\times
e^{\frac{\bar{g}}{48}\phi^4}\ .
\label{shift5bdeq}
\end{equation}
Finally we use the bounds (\ref{shiftbdeq}),(\ref{shift1bdeq}),(\ref{shift3bdeq}),(\ref{shift4bdeq})
and (\ref{shift5bdeq}) in order to derive the inequality
\begin{eqnarray*}
||\xi_{R,\Delta'}^{{\rm shift}}(\vec{V})||_{\partial\phi,\phi,h}
 & \le & 2 e^{-\frac{\bar{g}}{16}\phi^4} \\
 & & \times 412\times
L^{\frac{5}{2}}\bar{g}^{\frac{1}{2}-2\eta}
\times 11(5+\sqrt{2})
\times
\left(\max_{0\le n\le 6}\left(\frac{n}{4e}\right)^{\frac{n}{4}}\right)\times
48^{\frac{3}{2}}\times 2^6\times
e^{\frac{\bar{g}}{48}\phi^4} \\
 & & \times
3\mathcal{O}_{27} 
L^5 \bar{g}^{\frac{1}{4}-\eta}e^{\frac{\bar{g}}{48}\phi^4} \\
 & & \times
e^{\frac{1}{48}}\times e^{\frac{\bar{g}}{48}\phi^4}
\end{eqnarray*}
which after cleaning up becomes the desired bound.
\qed

\begin{Lemma}\label{L47lem}
For all unit cube $\Delta'$ we have that
\[
|\xi_{R,\Delta'}^{{\rm shift}}(\vec{V})(\phi)|_{\partial\phi,h_\ast}
\le \mathcal{O}_{40} L^{\frac{35}{2}} \bar{g}^{\frac{35}{12}-2\eta-\frac{1}{3}\eta_R}
\]
where 
\[
\mathcal{O}_{40}=2^{\frac{13}{2}}\times 3\times 11\times 103\times \log\left(\frac{3}{2}\right)
\times(5+\sqrt{2})\times\mathcal{O}_{28}\ .
\]
\end{Lemma}

\noindent{\bf Proof:}
As before we have
\[
|\xi_{R,\Delta'}^{{\rm shift}}(\vec{V})(\phi)|_{\partial\phi,h_\ast}\le
|e^{-\hat{V}_{\Delta'}(\phi)}|_{\partial\phi,h_\ast}
|Q'_{\Delta'}(\phi)|_{\partial\phi,h_\ast}
|U_{\Delta'}(\phi)|_{\partial\phi,h_\ast}
e^{ |U_{\Delta'}(\phi)|_{\partial\phi,h_\ast}}\ .
\]
Lemma \ref{L10lem} allows us to bound the first factor by
\[
|e^{-\hat{V}_{\Delta'}(\phi)}|_{\partial\phi,h_\ast}\le 2\ .
\]

The quantity $U_{\Delta'}(\phi)$ is estimated in the same as 
$\delta V_{\Delta'}(\lambda,\phi)$ in Lemma \ref{L36lem} with $\lambda=1$.
except that the sum over $k$ goes from $1$ to $4$ instead of from $0$ to $4$, and $\lambda$ is now set equal to $1$.
Hence
\[
|U_{\Delta'}(\lambda,\phi)|_{\partial\phi,h_\ast}\le 400 \times L^9
\times\left[
4\mathcal{O}_{25}+5\mathcal{O}_{26}+5\mathcal{O}_{24}
\right]\times \bar{g}^{\frac{11}{12}-\frac{1}{3}\eta_R}
\]
or equivalently
\[
|U_{\Delta'}(\lambda,\phi)|_{\partial\phi,h_\ast}\le \log\left(\frac{3}{2}\right)
\mathcal{O}_{28}
\times L^9
\bar{g}^{\frac{11}{12}-\frac{1}{3}\eta_R}\ .
\]
Now the standard hypothesis (\ref{extrastdeq}) implies 
\[
|U_{\Delta'}(\lambda,\phi)|_{\partial\phi,h_\ast}\le \log\left(\frac{3}{2}\right)
\]
and thus
\[
e^{|U_{\Delta'}(\lambda,\phi)|_{\partial\phi,h_\ast}}\le \frac{3}{2}\ .
\]
Finally we bound $Q'$ as before by writing
\[
|Q'_{\Delta'}(\phi)|_{\partial\phi,h_\ast}\le 
412\times
\max\left[|W'_{5,\Delta'}|,|W'_{6,\Delta'}|\right]\times \max_{0\le a\le 6} |\phi^a|_{\partial\phi,h_\ast}\ ,
\]
but for $0\le a\le 6$, 
\[
|\phi^a|_{\partial\phi,h_\ast}=h_{\ast}^{a}\le h_{\ast}^{6}=\left(2^{\frac{3}{4}}L^{\frac{3+\epsilon}{4}}\right)^6
\le 2^{\frac{9}{2}} L^6\ .
\]
This, together with the $W'$ bounds from the proof of the previous lemma,
provides us with the estimate
\[
|Q'_{\Delta'}(\phi)|_{\partial\phi,h_\ast}\le 
412\times 2^{\frac{9}{2}} L^{\frac{17}{2}} \bar{g}^{2-2\eta}\times 11\times (5+\sqrt{2})\ .
\]
Altogether we collect the bound
\[
|\xi_{R,\Delta'}^{{\rm shift}}(\vec{V})(\phi)|_{\partial\phi,h_\ast}\le 2\times \frac{3}{2}
\times \log\left(\frac{3}{2}\right)
\mathcal{O}_{28}
\times L^9
\bar{g}^{\frac{11}{12}-\frac{1}{3}\eta_R}
\]
\[
\times 412\times 2^{\frac{9}{2}} L^{\frac{17}{2}} \bar{g}^{2-2\eta}\times 11\times (5+\sqrt{2})
\]
which after cleaning up becomes the desired bound.
\qed

We now combine the last two lemmas into a single more convenient result.

\begin{Lemma}\label{L48lem}
For all unit cube $\Delta'$ we have
\[
|||\xi_{R,\Delta'}^{{\rm shift}}|||_{\bar{g}}\le\mathcal{O}_{41} L^{\frac{35}{2}} \bar{g}^{\frac{11}{4}-3\eta}
\]
where $\mathcal{O}_{41}=\max(\mathcal{O}_{39},\mathcal{O}_{40})$.
\end{Lemma}

\noindent{\bf Proof:}
From Lemmas \ref{L46lem} and \ref{L47lem} we immediately obtain
\[
|||\xi_{R,\Delta'}^{{\rm shift}}|||_{\bar{g}}\le
\max\left[
\mathcal{O}_{40} L^{\frac{35}{2}}\bar{g}^{\frac{35}{12}-2\eta-\frac{1}{3}\eta_R}\ ,\ 
\mathcal{O}_{39} L^{\frac{15}{2}}\bar{g}^{\frac{11}{4}-3\eta}
\right]\ .
\]
We have
\[
\left(\frac{35}{12}-2\eta-\frac{1}{3}\eta_R\right)-
\left(\frac{11}{4}-3\eta\right)=\frac{1}{6}+\eta-\frac{1}{3}\eta_R
\]
\[
\ge \frac{1}{6}+\eta-\frac{1}{3}\left(\frac{3}{16}+\frac{9}{4}\eta\right)
=\frac{5}{48}+\frac{1}{4}\eta>0
\]
by the standard hypotheses (\ref{etar2stdeq}) and (\ref{etastdeq}). Therefore
\[
\bar{g}^{\frac{35}{12}-2\eta-\frac{1}{3}\eta_R}\le
\bar{g}^{\frac{11}{4}-3\eta}
\]
and the result follows.
\qed

\begin{Lemma}\label{L49lem}
For all $\Delta'\in\mathbb{L}$ and $\Delta\in[L^{-1}\Delta']$ we have
\[
|J_{\Delta',\Delta_1}(\phi)|_{\partial\phi,L^{[\phi]}h_\ast}\le\mathcal{O}_{42} |||R_{\Delta_1}|||_{\bar{g}}
\]
where $\mathcal{O}_{42}=\mathcal{O}_{23}+250\mathcal{O}_{24}$. 
\end{Lemma}

\noindent{\bf Proof:}
By definition
\[
J_{\Delta',\Delta_1}(\phi)=J_+(\phi)-J_-(\phi)
\]
where
\[
J_+(\phi)=e^{-\frac{1}{2}(f,\Gamma f)_{L^{-1}\Delta'}}\times \int {\rm d}\mu_{\Gamma}(\zeta)
e^{\int_{L^{-1}\Delta'} f \zeta}
R_{\Delta_1}(\phi_1+\zeta)
\]
and
\[
J_-(\phi)=\left(\sum_{k=0}^{4}
\delta\beta_{k,3,\Delta',\Delta_1} :\phi^k:_{C_0}
\right)\times
e^{-\tilde{V}_{\Delta_1}(\phi_1)}\ .
\]
By Lemma \ref{L31lem} we have
\[
|J_+(\phi)|_{\partial\phi,L^{[\phi]}h_\ast}\le\mathcal{O}_{23} |||R_{\Delta_1}|||_{\bar{g}}\ .
\]
By Lemma \ref{L32lem} we also have
\[
|\delta\beta_{k,3,\Delta',\Delta_1}|\le \mathcal{O}_{24} L^{-k[\phi]}h_{\ast}^{-k} |||R_{\Delta_1}|||_{\bar{g}}\ .
\] 
We again use
\[
|:\phi^k:_{C_0}|_{\partial\phi,L^{[\phi]}h_\ast}\le 25\max_{0\le a\le k} 
|\phi^a|_{\partial\phi,L^{[\phi]}h_\ast}\le 25 L^{k[\phi]}h_{\ast}^{k}
\]
since $|\phi^a|_{\partial\phi,L^{[\phi]}h_\ast}=(L^{[\phi]}h_{\ast})^a$ and $L^{[\phi]}h_{\ast}\ge 1$.
Finally, by the chain rule
\[
|e^{-\tilde{V}_{\Delta_1}(\phi_1)}|_{\partial\phi,L^{[\phi]}h_\ast}=
|e^{-\tilde{V}_{\Delta_1}(\psi)}|_{\partial\psi,h_\ast}\le 2
\]
by Lemma \ref{L8lem}.
As result we easily arrive at
\[
|J_-(\phi)|_{\partial\phi,L^{[\phi]}h_\ast}\le 250\mathcal{O}_{24} |||R_{\Delta_1}|||_{\bar{g}}\ .
\]
The latter as well as the previous inequality for  $J_+$ imply the desired estimate. 
\qed

\begin{Lemma}\label{L50lem}
For all $\Delta'\in\mathbb{L}$, $\Delta\in[L^{-1}\Delta']$ and $\phi\in\mathbb{R}$ we have
\[
||J_{\Delta',\Delta_1}(\phi)||_{\partial\phi,\phi,L^{[\phi]}h}\le\mathcal{O}_{43} \bar{g}^{-2} |||R_{\Delta_1}|||_{\bar{g}}
\]
where 
\[
\mathcal{O}_{43}=\exp\left(\frac{\sqrt{2}}{2}\right)+155\mathcal{O}_{24}\ .
\]
\end{Lemma}

\noindent{\bf Proof:}
Clearly, we have
\[
||J_{\Delta',\Delta_1}(\phi)||_{\partial\phi,\phi,L^{[\phi]}h}\le
||J_{+}(\phi)||_{\partial\phi,\phi,L^{[\phi]}h}+
||J_{-}(\phi)||_{\partial\phi,\phi,L^{[\phi]}h}
\]
and both terms will be bounded as follows.
We first write
\[
||J_{+}(\phi)||_{\partial\phi,\phi,L^{[\phi]}h}\le
e^{-\frac{1}{2}\Re(f,\Gamma f)_{L^{-1}\Delta'}}\times
\int {\rm d}\mu_{\Gamma}(\zeta)\ e^{\int_{L^{-1}\Delta'}(\Re f)\zeta}
||R_{\Delta_1}(\phi_1+\zeta)||_{\partial\phi,\phi,L^{[\phi]}h}
\]
and then use the chain rule as well as the definition of the $|||\cdot|||_{\bar{g}}$ norm in order to derive
\[
||R_{\Delta_1}(\phi_1+\zeta)||_{\partial\phi,\phi,L^{[\phi]}h}=
||R_{\Delta_1}(\psi+\zeta)||_{\partial\psi,\phi_1,h}=
||R_{\Delta_1}(\psi)||_{\partial\psi,\phi_1+\zeta,h}
\le \bar{g}^{-2} |||R_{\Delta_1}|||_{\bar{g}}\ .
\]
Besides, as shown before $|(f,\Gamma f)_{L^{-1}\Delta'}|\le \frac{1}{\sqrt{2}}$.
Hence
\[
||J_{+}(\phi)||_{\partial\phi,\phi,L^{[\phi]}h}\le
\exp[2^{-\frac{3}{2}}] \bar{g}^{-2} |||R_{\Delta_1}|||_{\bar{g}}
\int {\rm d}\mu_{\Gamma}(\zeta)\ e^{\int_{L^{-1}\Delta'}(\Re f)\zeta}
\]
\[
\le \exp[2^{-\frac{3}{2}}] \bar{g}^{-2} |||R_{\Delta_1}|||_{\bar{g}}
e^{\frac{1}{2}(\Re f,\Gamma\Re f)_{L^{-1}\Delta'}}
\]
by Lemma \ref{L3lem} with $X=\emptyset$ or simply exact computation.
Again one easily gets that $|(\Re f,\Gamma\Re f)_{L^{-1}\Delta'}|\le \frac{1}{\sqrt{2}}$
which results in
\[
||J_{+}(\phi)||_{\partial\phi,\phi,L^{[\phi]}h}\le
\exp[2^{-\frac{1}{2}}] \bar{g}^{-2} |||R_{\Delta_1}|||_{\bar{g}}\ .
\]
From the definition of $J_-(\phi)$
we immediately get
\[
||J_{-}(\phi)||_{\partial\phi,\phi,L^{[\phi]}h}\le
\sum_{k=0}^{4}|\delta\beta_{k,3,\Delta',\Delta_1}|\times ||:\phi^k:_{C_0}||_{\partial\phi,\phi,L^{[\phi]}h}
||e^{-\tilde{V}_{\Delta_1}(\phi_1)}||_{\partial\phi,\phi,L^{[\phi]}h}\ .
\]
By the chain rule and Lemma \ref{L8lem}
\[
||e^{-\tilde{V}_{\Delta_1}(\phi_1)}||_{\partial\phi,\phi,L^{[\phi]}h}=
||e^{-\tilde{V}_{\Delta_1}(\psi)}||_{\partial\psi,\phi_1,h}\le
2 e^{-\frac{1}{2}(\Re\beta_{4,\Delta_1})\phi_1^4}\ .
\]
Again by undoing the Wick ordering we have
\[
||:\phi^k:_{C_0}||_{\partial\phi,\phi,L^{[\phi]}h}\le 25\max_{0\le a\le k}
||\phi^a||_{\partial\phi,\phi,L^{[\phi]}h}\ .
\]
But
\[
||\phi^a||_{\partial\phi,\phi,L^{[\phi]}h}=(L^{[\phi]}h+|\phi|)^a\le (L^{[\phi]}h+|\phi|)^k
\]
since $L^{[\phi]}h\ge 1$ as follows from $h\ge h_\ast$, i.e., from (\ref{c1c2stdeq}).
Now
\[
||:\phi^k:_{C_0}||_{\partial\phi,\phi,L^{[\phi]}h}\le
25 (L^{[\phi]}h+|\phi|)^k
=25\sum_{n=0}^{k}\left(
\begin{array}{c} k \\ n \end{array}
\right)\ \left(L^{[\phi]}c_1\bar{g}^{-\frac{1}{4}}\right)^{k-n} |\phi|^n
\]
\[
\le 25\sum_{n=0}^{k}\left(
\begin{array}{c} k \\ n \end{array}
\right)\ \left(L^{[\phi]}c_1\bar{g}^{-\frac{1}{4}}\right)^{k-n}
\left(\frac{n}{2e}\right)^{\frac{n}{4}}(\gamma\bar{g})^{-\frac{n}{4}}
e^{\gamma(\Re\beta_{4,\Delta_1})\phi^4}
\]
by Lemma \ref{L4lem} and for any $\gamma>0$. Here we choose $\gamma=\frac{1}{2}L^{-4[\phi]}$
which entails
\[
||:\phi^k:_{C_0}||_{\partial\phi,\phi,L^{[\phi]}h}\le
25\times\left(\max_{0\le n\le 4} \left(\frac{n}{2e}\right)^{\frac{n}{4}} \right)
\times e^{\frac{1}{2}(\Re\beta_{4,\Delta_1})\phi_1^4}
\]
\[
\times \sum_{n=0}^{k}\left(
\begin{array}{c} k \\ n \end{array}
\right)\ \left(L^{[\phi]}c_1\bar{g}^{-\frac{1}{4}}\right)^{k-n}
\left(\frac{1}{2}L^{-4[\phi]}\bar{g}\right)^{-\frac{n}{4}}\ .
\]
As a result of the previous considerations we arrive at
\[
||J_{-}(\phi)||_{\partial\phi,\phi,L^{[\phi]}h}\le 50\times
\left(\max_{0\le n\le 4} \left(\frac{n}{2e}\right)^{\frac{n}{4}} \right)
\]
\[
\times
\sum_{k=0}^{4}|\delta\beta_{k,3,\Delta',\Delta_1}|\times 
L^{k[\phi]}\bar{g}^{-\frac{k}{4}}
\left(\sum_{n=0}^{k}\left(
\begin{array}{c} k \\ n \end{array}
\right) c_1^{k-n} 2^{\frac{n}{4}}
\right)\ .
\]
Since $n\le 4$ we simply bound $\frac{n}{2e}$ by $1$. We also
use Lemma \ref{L32lem} in order to write
\[
||J_{-}(\phi)||_{\partial\phi,\phi,L^{[\phi]}h}\le 50\times
\mathcal{O}_{24} |||R_{\Delta_1}|||_{\bar{g}}
\times \sum_{n=0}^{k} \bar{g}^{-\frac{k}{4}} h_{\ast}^{-\frac{k}{4}}
(c_1+2^{\frac{1}{4}})^k\ .
\]
Now we bound $h_{\ast}^{-1}$ by $1$, $\bar{g}^{-\frac{k}{4}}$ by the worst case scenario $\bar{g}^{-1}\le \bar{g}^{-2}$
and finally $c_1+2^{\frac{1}{4}}$ by $2$. Since $1+2+\cdots+2^4=31$
we then obtain
\[
||J_{-}(\phi)||_{\partial\phi,\phi,L^{[\phi]}h}\le 50\times 31
\times
\mathcal{O}_{24} \bar{g}^{-2} |||R_{\Delta_1}|||_{\bar{g}} \ .
\]
The latter inequality, combined with the previous one for $J_+$, gives us the desired result.
\qed

\begin{Lemma}\label{L51lem}
For all unit cube $\Delta'\in\mathbb{L}$ we have
\[
|||\mathcal{L}_{\Delta'}^{(\vec{\beta},f)}(R)|||_{\bar{g}}\le\mathcal{O}_{44}
\times L^{3-5[\phi]}\times
\max_{\Delta_1\in[L^{-1}\Delta']} |||R_{\Delta_1}|||_{\bar{g}} 
\]
where
\[
\mathcal{O}_{44}=2^{10}\times \mathcal{O}_{14}
\times (\mathcal{O}_{42}+\mathcal{O}_{43})\ .
\]
\end{Lemma}

\noindent{\bf Proof:}
Recall that
\[
\mathcal{L}_{\Delta'}^{(\vec{\beta},f)}(R)=
\sum_{\Delta_1\in[L^{-1}\Delta']} \left(
\prod_{\substack{\Delta\in[L^{-1}\Delta'] \\ \Delta\neq\Delta_1}} e^{-\tilde{V}_{\Delta}(\phi_1)}
\right)\times 
J_{\Delta',\Delta_1}(\phi)\ .
\]
Hence
\[
|\mathcal{L}_{\Delta'}^{(\vec{\beta},f)}(R)|_{\partial\phi,h_\ast}\le
\sum_{\Delta_1\in[L^{-1}\Delta']} \left|
\prod_{\substack{\Delta\in[L^{-1}\Delta'] \\ \Delta\neq\Delta_1}} e^{-\tilde{V}_{\Delta}(\phi_1)}
\right|_{\partial\phi,h_\ast}\times 
|J_{\Delta',\Delta_1}(\phi)|_{\partial\phi,h_\ast}\ .
\]
Now by Lemma \ref{L9lem} with $Y_0=[L^{-1}\Delta']\backslash\{\Delta_1\}$ we have
\[
\left|
\prod_{\substack{\Delta\in[L^{-1}\Delta'] \\ \Delta\neq\Delta_1}} e^{-\tilde{V}_{\Delta}(\phi_1)}
\right|_{\partial\phi,h_\ast}\le 2\ .
\]
By definition of the seminorms
\[
|J_{\Delta',\Delta_1}(\phi)|_{\partial\phi,h_\ast}=
\sum_{n=0}^{9}\frac{h_{\ast}^{n}}{n!}|J_{\Delta',\Delta_1}^{(n)}(0)|\ .
\]
However, by construction in \S\ref{algdefsec},
the derivatives $J_{\Delta',\Delta_1}^{(n)}(0)$ vanish when $0\le n\le 4$.
As a result
\[
|J_{\Delta',\Delta_1}(\phi)|_{\partial\phi,h_\ast}=
\sum_{n=5}^{9}\frac{h_{\ast}^{n}}{n!}|J_{\Delta',\Delta_1}^{(n)}(0)|
=
\sum_{n=5}^{9} L^{-n[\phi]}
\frac{(h_{\ast}L^{[\phi]})^{n}}{n!}|J_{\Delta',\Delta_1}^{(n)}(0)|
\]
\[
\le L^{-5[\phi]}
\sum_{n=5}^{9} 
\frac{(h_{\ast}L^{[\phi]})^{n}}{n!}|J_{\Delta',\Delta_1}^{(n)}(0)|
=L^{-5[\phi]}|J_{\Delta',\Delta_1}(\phi)|_{\partial\phi,L^{[\phi]}h_\ast}
\]
and thus by Lemma \ref{L49lem} we have
\begin{eqnarray*}
|\mathcal{L}_{\Delta'}^{(\vec{\beta},f)}(R)|_{\partial\phi,h_\ast} & \le & 2  L^{-5[\phi]}
\sum_{\Delta_1\in[L^{-1}\Delta']}\mathcal{O}_{42} |||R_{\Delta_1}|||_{\bar{g}} \\
 & \le & 2 \mathcal{O}_{42}  L^{3-5[\phi]}
\max_{\Delta_1\in[L^{-1}\Delta']} |||R_{\Delta_1}|||_{\bar{g}}\ .
\end{eqnarray*}
Likewise, we have
\[
||\mathcal{L}_{\Delta'}^{(\vec{\beta},f)}(R)||_{\partial\phi,\phi,h}\le
\sum_{\Delta_1\in[L^{-1}\Delta']} \left|\left|
\prod_{\substack{\Delta\in[L^{-1}\Delta'] \\ \Delta\neq\Delta_1}} e^{-\tilde{V}_{\Delta}(\phi_1)}
\right|\right|_{\partial\phi,\phi,h}\times 
||J_{\Delta',\Delta_1}(\phi)||_{\partial\phi,\phi,h}\ .
\]
If we let $Y_0=[L^{-1}\Delta']\backslash\{\Delta_1\}$, then $|Y_0|\ge \frac{L^3}{2}$ and by Lemma \ref{L9lem}
we have
\[
\left|\left|
\prod_{\substack{\Delta\in[L^{-1}\Delta'] \\ \Delta\neq\Delta_1}} e^{-\tilde{V}_{\Delta}(\phi_1)}
\right|\right|_{\partial\phi,\phi,h}\le 2
e^{-\frac{\bar{g}}{16}\phi^4}\ .
\]
By Lemma \ref{L24lem} with $\beta_4=\bar{g}$ and $\gamma=\frac{1}{16}$ one has
\[
||J_{\Delta',\Delta_1}(\phi)||_{\partial\phi,\phi,h}\le \mathcal{O}_{14}
16^{\frac{9}{4}} e^{\frac{\bar{g}}{16}\phi^4}
\left[|J_{\Delta',\Delta_1}(\phi)|_{\partial\phi,h}
+L^{-9[\phi]}\sup_{\psi\in\mathbb{R}} |J_{\Delta',\Delta_1}(\psi)|_{\partial\psi,\psi,L^{[\phi]}h}
\right]\ .
\] 
By the same argument utilizing the vanishing of the first few derivatives at the origin as before, 
with $h$ instead of $h_\ast$, we get
\[
|J_{\Delta',\Delta_1}(\phi)|_{\partial\phi,h}\le 
L^{-5[\phi]} |J_{\Delta',\Delta_1}(\phi)|_{\partial\phi,L^{[\phi]}h}\ .
\]
Now by Lemma \ref{L49lem}
\[
|J_{\Delta',\Delta_1}(\phi)|_{\partial\phi,L^{[\phi]}h}\le\mathcal{O}_{42}
|||R_{\Delta_1}|||_{\bar{g}}
\]
whereas, by Lemma \ref{L50lem}, one has
\[
\sup_{\psi\in\mathbb{R}} |J_{\Delta',\Delta_1}(\psi)|_{\partial\psi,\psi,L^{[\phi]}h}
\le \mathcal{O}_{43}\bar{g}^{-2}
|||R_{\Delta_1}|||_{\bar{g}}\ .
\]
We then arrive at the estimate
\[
||J_{\Delta',\Delta_1}(\phi)||_{\partial\phi,\phi,h}\le \mathcal{O}_{14}\times 2^9\times
e^{\frac{\bar{g}}{16}\phi^4} |||R_{\Delta_1}|||_{\bar{g}}\times
\left[L^{-5[\phi]}\mathcal{O}_{42}+L^{-9[\phi]}\bar{g}^{-2}\mathcal{O}_{42}
\right]\ .
\]
Using $L^{-5[\phi]}\bar{g}^{-2}$ as a common bound of $L^{-9[\phi]}\bar{g}^{-2}$ and $L^{-5[\phi]}$ 
we immediately get
\[
||\mathcal{L}_{\Delta'}^{(\vec{\beta},f)}(R)(\phi)||_{\partial\phi,\phi,h}\le
\sum_{\Delta_1\in[L^{-1}\Delta']} 2^{10}\mathcal{O}_{14}(\mathcal{O}_{42}+\mathcal{O}_{43})
L^{-5[\phi]}\bar{g}^{-2}|||R_{\Delta_1}|||_{\bar{g}}
\]
and hence
\[
\bar{g}^2\times ||\mathcal{L}_{\Delta'}^{(\vec{\beta},f)}(R)(\phi)||_{\partial\phi,\phi,h}\le
2^{10}\mathcal{O}_{14}(\mathcal{O}_{42}+\mathcal{O}_{43})
L^{3-5[\phi]}\max_{\Delta_1\in[L^{-1}\Delta']} |||R_{\Delta_1}|||_{\bar{g}}\ .
\]
The latter inequality, combined with the previous one for the $|\cdot|_{\partial\phi,h_\ast}$
seminorm, give
\[
|||\mathcal{L}_{\Delta'}^{(\vec{\beta},f)}(R)|||_{\bar{g}}\le 
L^{3-5[\phi]}\left(\max_{\Delta_1\in[L^{-1}\Delta']} |||R_{\Delta_1}|||_{\bar{g}}\right)
\times
\max\left[2\mathcal{O}_{42},2^{10}\mathcal{O}_{14}(\mathcal{O}_{42}+\mathcal{O}_{43})\right]\ .
\]
Since clearly $\mathcal{O}_{14}>1$, the last maximum reduces to the second term, i.e., the given value
of $\mathcal{O}_{44}$.
\qed
 
\begin{Lemma}\label{L52lem}
For all unit cube $\Delta'$ we have
\[
|||\xi_{R,\Delta'}(\vec{V})|||_{\bar{g}}\le\mathcal{O}_{45} L^{60} \bar{g}^{\frac{11}{4}-3\eta}
\]
where $\mathcal{O}_{45}=\mathcal{O}_{37}+\mathcal{O}_{38}+\mathcal{O}_{41}$.
\end{Lemma}

\noindent{\bf Proof:}
Since (\ref{xideceq}) implies
\[
|||\xi_{R,\Delta'}(\vec{V})|||_{\bar{g}}\le
|||\xi_{R,\Delta'}^{{\rm main}}(\vec{V})|||_{\bar{g}}
+|||\xi_{R,\Delta'}^{{\rm higher}}(\vec{V})|||_{\bar{g}}
+|||\xi_{R,\Delta'}^{{\rm shift}}(\vec{V})|||_{\bar{g}}\ ,
\]
Lemmas \ref{L44lem}, \ref{L45lem} and \ref{L48lem} immediately imply the desired result.
\qed

\subsection{Conclusion of the proof of the main estimates for $RG_{\rm ex}$}\label{concthmsec}

We choose $B_{R\mathcal{L}}=\mathcal{O}_{44}$ defined in Lemma \ref{L51lem}.
Then if one fixes $l\ge 1$ or equivalently $L$, we take $B_k=\mathcal{O}_{24} L^{3-k[\phi]}$, for
$0\le k\le 4$. We also set $B_{R\xi}=\mathcal{O}_{45} L^{60}$.
The hypotheses on $\eta$ and $\eta_R$ in the statement of Theorem \ref{mainestthm} easily
imply the properties of $\eta$ and $\eta_R$ mentioned in the standard hypotheses
(\ref{etastdeq}), (\ref{etarstdeq}) and (\ref{etar2stdeq}). 
Now once $A_{\bar{g}}$ has been chosen, the calibrator $\bar{g}$ can be made as small as desired
by taking $\epsilon$ small enough. It is a simple matter of going through the inequalities
in \S\ref{stdhypsec} in order to check that all the standard hypotheses are satisfied for small $\epsilon$
and therefore all the lemma in the previous section hold. In particular, for $0\le k\le 4$,
\[
|\xi_{k,\Delta'}(\vec{V})|\le B_k \max_{\Delta\in[L^{-1}\Delta']} |||R_\Delta|||_{\bar{g}}\qquad
\]
follows from the definition $\xi_{k,\Delta'}(\vec{V})=-\delta\beta_{k,3,\Delta'}$ and Lemma \ref{L32lem}.
Likewise,
\[
|||\mathcal{L}_{\Delta'}^{\vec{\beta},f}(R)|||_{\bar{g}}\le
B_{R\mathcal{L}} L^{3-5[\phi]}
\max_{\Delta\in[L^{-1}\Delta']} |||R_\Delta|||_{\bar{g}}
\]
is the result of Lemma \ref{L51lem}.
Finally, 
\[
|||\xi_{R,\Delta'}(\vec{V})|||_{\bar{g}}\le B_{R\xi} \bar{g}^{\frac{11}{4}-3\eta}
\]
has been established in Lemma \ref{L52lem}.

The statement in Theorem \ref{mainestthm} about sending real data to real data is obvious from the definition the RG map
in \S\ref{algdefsec}. So is the one about analyticity now that the previous bounds on the outcome have been proved.
The proof of Theorem \ref{mainestthm} is now complete.

%% file: BulkRG.tex
\section{The bulk RG}\label{bulksec}

In \S\ref{funcspsec} we defined the complex Banach spaces
$\mathcal{E}$, $\mathcal{E}_{\rm 1B}$, $\mathcal{E}_{\rm bk}$ and $\mathcal{E}_{\rm ex}$.
The transformation $RG_{\rm ex}$ defined in \S\ref{algdefsec} is an analytic map from a domain
in $\mathcal{E}_{\rm ex}$ (given in the hypotheses of Theorem \ref{mainestthm}) into $\mathcal{E}_{\rm ex}$.
In this section we will show that the subspace $\mathcal{E}_{\rm bk}$ is stable by this transformation
and similarly for $\mathcal{E}=\mathbb{C}^2\times C_{\rm bd,ev}^{9}(\mathbb{R},\mathbb{C})$ seen as a subspace
of $\mathcal{E}_{\rm bk}$ and therefore of $\mathcal{E}_{\rm ex}$ too. We will also derive simpler formulas
for the transformation restricted to $\mathcal{E}$.

\begin{Proposition}
The space $\mathcal{E}_{\rm bk}$ is invariant by $RG_{\rm ex}$.
\end{Proposition}

\noindent{Proof:}
This is a trivial consequence of the translation covariance of the definition of $RG_{\rm ex}$ in \S\ref{algdefsec}.
\qed

Let $(g,\mu,R)\in\mathcal{E}$. This corresponds to an element
\[
\vec{V}=(\beta_{4,\Delta},\beta_{3,\Delta},\beta_{2,\Delta},\beta_{1,\Delta},
W_{5,\Delta},W_{6,\Delta},f_{\Delta},R_{\Delta})_{\Delta\in\mathbb{L}}
\]
in $\mathcal{E}_{\rm ex}$ via the specifications
\begin{eqnarray*}
\beta_{4,\Delta}  & = & g \\
\beta_{3,\Delta}  & = & 0 \\
\beta_{2,\Delta}  & = & \mu \\
\beta_{1,\Delta}  & = & 0 \\
W_{5,\Delta}  & = & 0 \\
W_{6,\Delta}  & = & 0 \\
f_{\Delta}  & = & 0 \\
R_{\Delta}  & = & R
\end{eqnarray*}
for all unit cubes $\Delta$.
We introduce the notations
\begin{eqnarray*}
\xi_4(g,\mu,R) & = & \xi_{4,\Delta(0)}(\vec{V}) \\
\xi_2(g,\mu,R) & = & \xi_{2,\Delta(0)}(\vec{V}) \\
\xi_0(g,\mu,R) & = & \xi_{0,\Delta(0)}(\vec{V}) \\
\xi_R(g,\mu,R) & = & \xi_{R,\Delta(0)}(\vec{V}) \\
\mathcal{L}^{(g,\mu)}(R) & = & \mathcal{L}_{\Delta(0)}^{(\vec{\beta},f)}(R) \\
\delta b(g,\mu,R) & = & \delta b_{\Delta(0)}(\vec{V}) 
\end{eqnarray*}
in terms of the previous vector $\vec{V}$. Note that we could have used any box $\Delta'$
instead of $\Delta(0)$, the one containing the origin.

\begin{Proposition}\label{L53prop}
The space $\mathcal{E}$ is invariant by the map $RG_{\rm ex}$. The restricted transformation
\[
\begin{array}{cccc}
RG: & \mathcal{E} & \longrightarrow & \mathcal{E} \\
 & (g,\mu,R) & \longmapsto & (g',\mu',R') 
\end{array}
\]
which we call the bulk RG is, more explicitly, given
by
\[
\left\{
\begin{array}{ccl}
g' & = & L^{\epsilon} g - A_1 g^2+\xi_4(g,\mu,R) \\
\mu' & = & L^{\frac{3+\epsilon}{2}} \mu -A_2 g^2 -A_3 g\mu +\xi_2(g,\mu,R) \\
R' & = & \mathcal{L}^{(g,\mu)}(R) + \xi_R(g,\mu,R)
\end{array}
\right.
\]
where
\begin{eqnarray*}
A_1 & = & 36 L^{3-4[\phi]}\int_{\mathbb{Q}_p^3} \Gamma(x)^2\ {\rm d}^3 x \\
A_2 & = & 48 L^{3-2[\phi]}\left(\int_{\mathbb{Q}_p^3} \Gamma(x)^3\ {\rm d}^3 x\right)
+144 L^{3-4[\phi]}C_0(0)\left(\int_{\mathbb{Q}_p^3} \Gamma(x)^2\ {\rm d}^3 x\right)\\
A_3 & = & 12 L^{3-2[\phi]}\int_{\mathbb{Q}_p^3} \Gamma(x)^2\ {\rm d}^3 x \ .
\end{eqnarray*}
In addition, the vacuum counter-term $\delta b=\delta b(g,\mu,R)$ is given by 
\[
\delta b = A_4 g^2+ A_5 \mu^2+\xi_0(g,\mu,R)
\]
where
\begin{eqnarray*}
A_4 & = & 12 L^{3}\left(\int_{\mathbb{Q}_p^3} \Gamma(x)^4\ {\rm d}^3 x\right)
+ 48 L^{3-2[\phi]}C_0(0)\left(\int_{\mathbb{Q}_p^3} \Gamma(x)^3\ {\rm d}^3 x\right)
+72 L^{3-4[\phi]}C_0(0)^2\left(\int_{\mathbb{Q}_p^3} \Gamma(x)^2\ {\rm d}^3 x\right)   \\
A_5 & = & L^{3}\int_{\mathbb{Q}_p^3} \Gamma(x)^2\ {\rm d}^3 x \ .
\end{eqnarray*}
\end{Proposition}

\noindent{\bf Proof:}
We compute the specialization of the map $\vec{V}\mapsto\vec{V}'$ defined in \S\ref{algdefsec}
to the present situation.
Clearly since $f_\Delta=0$, the new $f'_{\Delta'}$'s defined in (\ref{newfdefeq}) are identically zero.
Likewise, since the $W_6$ are zero the equation for the new one reduces to
\[
W'_{6,\Delta'}=8L^{-6[\phi]}\int_{(L^{-1}\Delta')^2} {\rm d}^3 x\ {\rm d}^3 y\ 
\beta_4(x)\ \Gamma(x-y)\ \beta_4(y)
=8L^{-6[\phi]} g^2 \int_{(L^{-1}\Delta')^2} {\rm d}^3 x\ {\rm d}^3 y\ \Gamma(x-y)\ .
\]
But for $x\in L^{-1}\Delta'$, by a simple change of variables $z=x-y$,
\[
\int_{L^{-1}\Delta'} {\rm d}^3 y\ \Gamma(x-y)=\int_{L^{-1}\Delta(0)} {\rm d}^3 z\ \Gamma(z)
=\int_{\mathbb{Q}_p^3} {\rm d}^3 z\ \Gamma(z)=\widehat{\Gamma}(0)=0
\]
because of the finite range property and the vanishing property at zero momentum.
Therefore $W'_{6,\Delta'}$ vanishes identically and so does $W'_{5,\Delta'}$ for similar reasons.
Now one easily sees from the definition and specification of the input $\vec{V}$ that
\begin{eqnarray*}
\hat{\beta}_{4,\Delta'} & = & L^{\epsilon} g \\
\hat{\beta}_{3,\Delta'} & = & 0 \\
\hat{\beta}_{2,\Delta'} & = & L^{\frac{3+\epsilon}{2}} \mu \\
\hat{\beta}_{1,\Delta'} & = & 0\ .
\end{eqnarray*}
Then we consider the first corrections $\delta\beta_{k,1,\Delta'}$.
These are all zero since in their defining equation (\ref{deltb1eq}) the constraint $b\ge 1$ implies that
at least one $f$ is present. However $f=0$ identically and therefore the Feynman diagram in (\ref{deltb1eq})
vanishes. In sum, $\delta\beta_{k,1,\Delta'}=0$ for all $k$ such that $0\le k\le 4$ and all unit cube $\Delta'$.

We now move on to the computation of the second order corrections $\delta\beta_{k,2,\Delta'}$.
Again since $f$ and the $W$'s are zero, the defining equation (\ref{deltb2eq}) for these quantities reduces to
\[
\delta\beta_{k,2,\Delta'}=
\sum_{a_1,a_2,m}
\bbone\left\{
\begin{array}{c}
a_i+m\le 4 \\
a_i\ge 0\ ,\ m\ge 1
\end{array}
\right\}
\frac{(a_1+m)!\ (a_2+m)!}{a_1!\ a_2!\ m!}
\]
\[
\times \frac{1}{2} C(a_1,a_2|k)\times
L^{-(a_1+a_2)[\phi]}\times C_0(0)^{\frac{a_1+a_2-k}{2}}\times
\beta_{a_1+m}\beta_{a_2+m}\times
\int_{(L^{-1}\Delta')^{2}} {\rm d}^3 x_1\ {\rm d}^3 x_2
\ \Gamma(x_1-x_2)^m\ .
\]
Indeed, nonvanishing imposes the absence of $f$ external vertices and thus $b_1=b_2=m$.
Note that since the $\beta_\nu(x)$ are constants with respect to the location $x$ we pulled them
out of the integral and suppressed the $x$ dependence in the notation.
Now one can rule out the value $m=1$ which gives a vanishing contribution 
for the same reason as explained above when computing $W'_6$.
Another simplification is that $\delta\beta_{k,2,\Delta'}$ vanishes if $k$ is odd.
Indeed, if $k$ is odd and if the connection coefficient $C(a_1,a_2|k)$ is nonzero, then $a_1+a_2$ must 
be odd too and thus also $a_1+m+a_2+m$. Since this forces $a_i+m$ to be odd for one $i=1,2$ then the
contribution in the above sum vanishes. This is because $\beta_{n}$ is nonzero only for even values of $n$,
namely $2$ and $4$.
We now only have three cases to consider: $k=4$ , $2$ and $0$.

\noindent{\bf 1st Case:} Let $k=4$. Then the connection coefficients force $a_1+a_2\ge 4$.
Also $m\ge 2$ and $a_i+m\le 4$ imply $0\le a_i\le 2$ so the only possibility is $(a_1,a_2)=(2,2)$ and $m=2$.
We also have $\beta_{a_1+m}=\beta_{a_1+m}=\beta_{4}=g$.
It is easy to see that the formula reduces to $\delta\beta_{4,2,\Delta'}=A_1 g^2$.

\noindent{\bf 2nd Case:} Let $k=2$. Now the constraints $a_1+a_2\ge 2$, $m\ge 2$, $a_i+m\le 4$, $a_i\ge 0$
and $a_i+m\in\{2,4\}$, without which the contribution would vanish, imply that the only possibilities for the
triple $(a_1,a_2,m)$ are $(2,2,2)$, $(2,0,2)$, $(0,2,2)$ and $(1,1,3)$.
The second and third give the same contribution by symmetry.
A quick computation shows that $(2,2,2)$
contributes
\[
144L^{3-4[\phi]}C_0(0) g^2 \int \Gamma^2
\]
where we used the shorthand
\[
 \int \Gamma^m= \int_{(L^{-1}\Delta')^{2}} {\rm d}^3 x_1\ {\rm d}^3 x_2
\ \Gamma(x_1-x_2)^m\ .
\]
Likewise $(2,0,2)$ and $(0,2,2)$ together contribute
\[
12 L^{3-2[\phi]} g\mu \int \Gamma^2\ .
\]
Finally, $(1,1,3)$ contributes
\[
48 L^{3-2[\phi]} g^2 \int \Gamma^3\ .
\]
Hence
\[
\delta\beta_{2,2,\Delta'}=A_2 g^2+A_3 g\mu\ .
\]

\noindent{\bf 3rd Case:} Let $k=0$.
Note that the connection coefficients also impose $0=k\ge|a_1-a_2|$ and thus the restriction $a_1=a_2$.
Considerations similar to those of the two previous cases show that the only possibilities
for the triple $(a_1,a_2,m)$ are $(0,0,2)$, $(0,0,4)$, $(1,1,3)$ and $(2,2,2)$.
Again a quick computation shows that $(0,0,2)$ contributes
\[
L^3\mu^2\int\Gamma^2\ .
\]
The triple $(0,0,4)$ contributes
\[
12L^3 g^2\int\Gamma^4\ .
\]
The triple $(1,1,3)$ contributes
\[
48 L^{3-2[\phi]}C_0(0) g^2\int\Gamma^3\ .
\]
Finally, the triple $(2,2,2)$ contributes
\[
72L^{3-4[\phi]} C_0(0)^2 g^2\int\Gamma^2\ .
\]
Hence,
\[
\delta\beta_{0,2,\Delta'}=A_4 g^2+A_5 \mu^2\ .
\]

We now show that $\delta\beta_{3,3,\Delta'}$ and $\delta\beta_{1,3,\Delta'}$ are zero because the function
$R$ is even.
First note that the formula (\ref{aieq}) for the $a_i$ reduces to
\[
a_i  =  \exp\left[
-C_0(0) L^{-2[\phi]}\mu+3C_0(0)^2 L^{-4[\phi]}g
\right]  \times L^{-i[\phi]}\times \int {\rm d}\mu_{\Gamma}(\zeta)
\ R^{(i)}(\zeta)\ .
\]
It is easy to see that $a_i=0$ if $i$ is odd.
Indeed, the Gaussian measure ${\rm d}\mu_{\Gamma}$ is centered and therefore one can change $\zeta$
into $-\zeta$ without changing the integral. However, for odd $i$ we have
$R^{(i)}(-\zeta)=-R^{(i)}(\zeta)$ for the $i$-th derivative of the even function $R$.
Thus $a_i=-a_i$ and the stated vanishing property holds.
Now the $M_{k,i}$ in (\ref{mkieq}) are zero unless $k$ and $i$ have the same parity.
Indeed, the $\beta_{l_{\nu}}$'s are nonzero only if $l_{\nu}=4$ or $2$.
This together with the constaint (\ref{mkiconsteq}) imply the desired property.
Therefore the $\delta\beta_{k,3,\Delta',\Delta_1}$ and thus also the $\beta'_{k,\Delta'}$
vanish for $k=1$ and $3$.

Finally in order to complete the proof, all that is needed is to show that $R'$ is an even function.
First notice that Wick powers only involve lower or equal ordinary powers of the same parity.
Since the inputs $\beta_{k}$ are zero for 
$k=1$ and $3$, the $\tilde{V}_\Delta(\phi_1)$
are even functions of $\phi$.
The defining formula (\ref{jddeq}) for $J_{\Delta',\Delta_1}(\phi)$ reduces, in the present situation, to
\[
J_{\Delta',\Delta_1}(\phi)=\left\{\int {\rm d}\mu_{\Gamma} (\zeta)
\ R(\phi_1+\zeta)\right\}
-\left(
\delta\beta_{4,3,\Delta',\Delta_1} :\phi^4:_{C_0} +
\delta\beta_{2,3,\Delta',\Delta_1} :\phi^2:_{C_0} +
\delta\beta_{0,3,\Delta',\Delta_1}  
\right)\times e^{-\tilde{V}_{\Delta_1}(\phi_1)}
\]
which is easily seen to be even thanks to the change of variable $\zeta\rightarrow -\zeta$
and the hypothesis that the input $R$ is even.
Therefore the quantity denoted by $\mathsf{C}'_1$
in \S\ref{algdefsec}, namely, $\mathcal{L}_{\Delta'}^{(\vec{\beta},f)}(R)(\phi)$
is even. Clearly $V_\Delta(\phi)$ is even which results in the invariance of $p_\Delta(\phi_1,\zeta)$ 
with respect to changing the sign of both $\phi$ and $\zeta$. Hence also $P_\Delta(\lambda,\phi_1,\zeta)$
has the same invariance property. Note that since the $W$'s are zero, $Q$ vanishes
and thus $K_\Delta(\lambda,\phi_1,\zeta)=\lambda^3 R(\phi_1+\zeta)$ has that invariance too.
It follows using the change of variable $\zeta\rightarrow -\zeta$ that
$\hat{K}_{\Delta'}(\lambda,\phi)$ is even.
Since $\delta V_{\Delta'}(\lambda,\phi)$ contains no $:\phi^3:_{C_0}$ nor $:\phi^1:_{C_0}$
it is even and as a result $K'_{\Delta'}(\lambda,\phi)$ is also even.
Since the $W'$ have been shown to vanish, the function $Q'_{\Delta'}$ also vanishes.
As a consequence one can easily see that $\xi_{R,\Delta'}(\vec{V})(-\phi)=\xi_{R,\Delta'}(\vec{V})(\phi)$.
Finally, $R'_{\Delta'}=R'_{\Delta(0)}=R'$ must be an even function of the field $\phi$. 
\qed

%% file: DynsysI.tex
\section{The infrared fixed point and local analysis of the bulk RG}\label{dyn1sec}

\subsection{Preparation}\label{prepsubsec}

In this section we make some choices related to the particular application of Theorem \ref{mainestthm}
which will be needed in the sequel.
Note that by Lemma \ref{gamL2lem}, the quantity $A_1$ defined in Proposition \ref{L53prop} is given by
\[
A_1=36L^{\epsilon}\times \frac{(1-p^{-3})(L^{\epsilon}-1)}{p^{\epsilon}-1}>0\ .
\]
If one ignores the $\xi_4$ term in the bulk RG evolution equation for the $\phi^4$ coupling $g$
then the fixed point equation becomes $g=L^{\epsilon}g-A_1g^2$. In addition to the trivial solution $g=0$,
this equation has another solution $g=(L^{\epsilon}-1)/A_1$. This is the approximate value of the $g$ coordinate
of the nontrivial infrared RG fixed point.
We will in the remainder of this article choose the calibrator defining the norms to be this approximate fixed point
value, namely, we set
\[
\bar{g}=\frac{L^{\epsilon}-1}{A_1}=\frac{p^{\epsilon}-1}{36L^{\epsilon}(1-p^{-3})}\ .
\] 
In other words, $\bar{g}=\bar{g}_{\ast}$ where $\bar{g}_{\ast}$ has been defined in \S\ref{formalstatsec}.
Clearly,
\[
\frac{\bar{g}}{\epsilon}\longrightarrow \frac{\log p}{36(1-p^{-3})}
\]
when $\epsilon\rightarrow 0$ with $L$ fixed.
This motivates making the choice
\[
A_{\bar{g}}=\frac{\log p}{36(1-p^{-3})} +1
\]
when applying Theorem \ref{mainestthm}.
This ensures that when $\epsilon$ is made small with $L$ fixed, our choice of $\bar{g}$ will satisfy the
requirement $\bar{g}\in (0,A_{\bar{g}}]$ in Theorem \ref{mainestthm}.
We now also choose $L$, once and for all, so that 
\begin{equation}
B_{R\mathcal{L}} L^{3-5[\phi]}\le \frac{1}{2}
\label{Lchoiceeq}
\end{equation}
holds.
Note that $3-5[\phi]=-\frac{3}{4}+\frac{5}{4}\epsilon$. If we add the harmless condition $\epsilon\le \frac{1}{5}$ which we now assume,
then $3-5[\phi]\le -\frac{1}{2}$. Now we pick $L$ large enough so that $B_{R\mathcal{L}}L^{-\frac{1}{2}}\le \frac{1}{2}$ and therefore
(\ref{Lchoiceeq}) holds.
Note that, contrary to the other $B$ quantities in Theorem \ref{mainestthm},
$B_{R\mathcal{L}}$ is independent of $L$ and indeed is a purely numerical constant. This fact is of course crucial to the previous considerations.
The choices for the parameters $\eta$, $\eta_R$ as well as the exponents $e$ in the definition of the Banach space norm
of $\mathcal{E}_{\rm ex}$ in \S\ref{funcspsec} will be specified later. Once these choice are made, the only free parameter in the problem
is the bifurcation parameter $\epsilon$. All the following results will be established in the regime when this $\epsilon$ is made sufficiently small.

We now apply Theorem \ref{mainestthm} with the choices just mentioned and in concert with Proposition \ref{L53prop}
to obtain that, provided $\epsilon$ is small enough, the bulk RG transformation is well-defined and analytic
on the domain
\[
|g-\bar{g}|<\frac{1}{2}\ ,\ |\mu|<\bar{g}^{1-\eta}\ ,\ |||R|||_{\bar{g}}<\bar{g}^{\frac{11}{4}-\eta_R}
\]
and therein satisfies
\begin{eqnarray*}
|\xi_4(g,\mu,R)| & \le & B_4  |||R|||_{\bar{g}} \\
|\xi_2(g,\mu,R)| & \le & B_2  |||R|||_{\bar{g}} \\
|||\xi_R(g,\mu,R)|||_{\bar{g}} & \le & B_{R\xi} \bar{g}^{\frac{11}{4}-3\eta} \\
|||\mathcal{L}^{(g,\mu)}|||_{\bar{g}} & \le & \frac{1}{2}
\end{eqnarray*}
where $|||\mathcal{L}^{(g,\mu)}|||_{\bar{g}}$ is the operator norm of the linear operator $\mathcal{L}^{(g,\mu)}$
(with respect to the $R$ variable)
corresponding to the norm $|||\cdot|||_{\bar{g}}$.
Note that the statement on analyticity applies not only to the full map $RG$ but also to the constituent pieces such as $\xi_4$, $\xi_2$, $\xi_R$
and $\mathcal{L}^{(g,\mu)}(R)$.

In order to analyze the bulk RG transformation we slightly change our coordinate system
from $(g,\mu,R)$ to $(\delta g,\mu, R)$ where $\delta g=g-\bar{g}$. In this new coordinate system, the bulk RG transformation, still denoted
by $RG$ for simplicity, becomes $(\delta g,\mu,R)\longmapsto RG(\delta g,\mu,R)=(\delta g',\mu',R')$
with
\[
\left\{
\begin{array}{ccl}
\delta g & = & (2-L^{\epsilon})\delta g+\tilde{\xi}_{4}(\delta g,\mu,R) \\
\mu' & = & L^{\frac{3+\epsilon}{2}}\mu+\tilde{\xi}_{2}(\delta g,\mu,R) \\
R' & = & \tilde{\mathcal{L}}^{(\delta g,\mu)(R)} +\tilde{\xi}_{R}(\delta g,\mu,R)
\end{array}
\right.
\]
where
\begin{eqnarray*}
\tilde{\xi}_{4}(\delta g,\mu,R) & = & -A_1 \delta g^2+\xi_{4}(\bar{g}+\delta g,\mu,R) \\
\tilde{\xi}_{2}(\delta g,\mu,R) & = & -A_2 (\bar{g}+\delta g)-A_3(\bar{g}+\delta g)\mu+\xi_{2}(\bar{g}+\delta g,\mu,R) \\
\tilde{\xi}_{R}(\delta g,\mu,R) & = & \xi_{R}(\bar{g}+\delta g,\mu,R) \\
\tilde{\mathcal{L}}^{(\delta g,\mu)}(R) & = & \mathcal{L}^{(\bar{g}+\delta g,\mu)(R)}
\end{eqnarray*}
as follows from an easy computation using the relation $A_1\bar{g}=L^{\epsilon}-1$.
We will commit a similar abuse of notation for the function $\delta b$. Namely, we will write $\delta b(\delta g,\mu,R)$
for what in fact is $\delta b(\bar{g}+\delta g,\mu,R)$.
Note that the norm we will use on such elements $v=(\delta g,\mu,R)\in\mathcal{E}$ is the one induced by the norm of the larger space
$\mathcal{E}_{\rm ex}$ defined in \S\ref{funcspsec}, namely,
\[
||v||=\max\left\{|\delta g|\bar{g}^{-e_4},|\mu|\bar{g}^{-e_2},|||R|||_{\bar{g}}\bar{g}^{-e_R}\right\}\ .
\]

We will assume the following constraints on the exponents defining the norms as well as the parameters $\eta$ and $\eta_R$:
\begin{eqnarray}
e_4 & \ge & 1 \label{econst1eq}\\
e_2 & \ge & 1-\eta \label{econst2eq}\\
e_R & \ge & \frac{11}{4}-\eta_R \label{econst3eq}\\
e_R & > & e_4+1 \label{econst4eq}\\
\frac{11}{4}-3\eta & > & e_R \label{econst5eq}\\
e_2 & < & 2\ . \label{econst6eq}
\end{eqnarray}

The following lemma provides Lipschitz estimates that will be needed in the sequel.

\begin{Lemma}\label{Lipxilem}
For $\epsilon$ small enough we have for all $v=(\delta g,\mu,R)$, $v'=(\delta g',\mu',R')$ in $\mathcal{E}$
such that $||v||$, $||v'||\le \frac{1}{8}$,
\[
|\xi_4(\bar{g}+\delta g,\mu,R)-\xi_4(\bar{g}+\delta g',\mu',R')|
\le 2B_4 \bar{g}^{e_R} ||v-v'||\ ,
\]
\[
|\xi_2(\bar{g}+\delta g,\mu,R)-\xi_4(\bar{g}+\delta g',\mu',R')|
\le 2B_2 \bar{g}^{e_R} ||v-v'||\ ,
\]
\[
|||{\mathcal{L}}^{(\bar{g}+\delta g,\mu)}(R)-
{\mathcal{L}}^{(\bar{g}+\delta g',\mu')}(R')|||_{\bar{g}}
\le \frac{3}{4}\bar{g}^{e_R} ||v-v'||
\]
and
\[
|||\xi_R(\bar{g}+\delta g,\mu,R)-\xi_R(\bar{g}+\delta g',\mu',R')|||_{\bar{g}}\le
3B_{R\xi}\bar{g}^{\frac{11}{4}-3\eta} ||v-v'||\ .
\]
\end{Lemma}

\noindent{\bf Proof:}
If $||v||<\frac{1}{2}$,
then since $\bar{g}\le 1$ for $\epsilon$ small and because of (\ref{econst1eq}),
(\ref{econst2eq}) and (\ref{econst3eq}), we have
\[
\begin{array}{ccccc}
|\delta g| & < & \frac{1}{2}\bar{g}^{e_4} & \le & \frac{1}{2}\bar{g} \\
|\mu| & < & \frac{1}{2}\bar{g}^{e_2} & \le & \frac{1}{2}\bar{g}^{1-\eta} \\
|||R|||_{\bar{g}} & < & \frac{1}{2}\bar{g}^{e_R} & \le & \frac{1}{2}\bar{g}^{\frac{11}{4}-\eta_R}\ .
\end{array}
\]
Hence, by Theorem \ref{mainestthm}
\[
|\xi_4(\bar{g}+\delta g, \mu,R)|\le B_4 |||R|||_{\bar{g}}\le \frac{1}{2}B_4\bar{g}^{e_R}\ .
\]
Therefore the analytic map $v\mapsto \xi_4(\bar{g}+\delta g, \mu,R)$
satisfies the hypotheses of Lemma \ref{Lipschitzlem} with $r_1=\frac{1}{2}$ and $r_2=\frac{1}{2}B_4\bar{g}^{e_R}$.
We pick $\nu=\frac{1}{4}$ which results in 
\[
\frac{r_2(1-\nu)}{r_1(1-2\nu)}=\frac{3}{2} B_4\bar{g}^{e_R}\ .
\] 
With these choices, Lemma \ref{Lipschitzlem} implies the desired Lipschitz estimate where we replaced
the numerical factor $\frac{3}{2}$ by $2$ for a simpler looking formula. The proof of the Lipschitz estimate
for $\xi_2$ is exactly the same  apart from changing $\xi_4$, $B_4$ to $\xi_2$, $B_2$ respectively.

We now do the same for the analytic map $v\mapsto {\mathcal{L}}^{(\bar{g}+\delta g,\mu)}(R)$.
For $||v||<\frac{1}{2}=r_1$ we obtain, as before from Theorem \ref{mainestthm} and from the choice we made when fixing $L$,
\[
|||{\mathcal{L}}^{(\bar{g}+\delta g,\mu)}(R)|||_{\bar{g}}
\le \frac{1}{2} |||R|||_{\bar{g}}\le \frac{1}{2}||v|| \bar{g}^{e_R}\le
r_2
\]
with $r_2=\frac{1}{4}\bar{g}^{e_R}$. Lemma \ref{Lipschitzlem} with $\nu=\frac{1}{4}$
now immediately implies the wanted estimate. Remark that we do not bound the numerical factor $\frac{3}{4}$
by the nearest integer here since it is important that this factor be less than $1$.

Finally, for $\xi_R$ we again note that $||v||<\frac{1}{2}=r_1$
implies
\[
|||\xi_R(\bar{g}+\delta g,\mu,R)|||_{\bar{g}}\le r_2
\]
with $r_2=B_{R\xi}\bar{g}^{\frac{11}{4}-\eta}$. Again, Lemma \ref{Lipschitzlem} with $\nu=\frac{1}{4}$
does the rest.
\qed

\subsection{The local stable manifold}\label{stabmansec}

In order to construct the nontrivial infrared fixed point we first construct its local stable manifold, then show that the RG transformation
is contractive on it and finally obtain the fixed point using the Banach Fixed Point Theorem.
We now proceed with the first step which is the construction of the stable manifold also using the Banach Fixed Point Theorem
in a space of one-sided sequences, in the spirit of Irwin's method~\cite{Irwin2}.
Let $\mathcal{B}_+$ be the Banach space of sequences
\[
\vec{u}=(\mu_0,(\delta g_1,\mu_1,R_1),(\delta g_2,\mu_2,R_2),\ldots)\in \mathbb{C}\times\prod_{n\ge 1}\left[
\mathbb{C}^2\times C_{\rm bd,ev}^{9}(\mathbb{R},\mathbb{C})
\right]
\]
which have finite norm given by
\[
||\vec{u}||=\sup\left\{
|\delta g_j|\bar{g}^{-e_4}\ {\rm for}\ j\ge 1;
|\mu_j|\bar{g}^{-e_2}\ {\rm for}\ j\ge 0;
|||R_j|||_{\bar{g}}\bar{g}^{-e_R}\ {\rm for}\ j\ge 1
\right\}\ .
\]
We will define a map $\mathfrak{m}$ on this space of sequences which depends on parameters $\delta g_0, R_0$ serving as boundary conditions.
Given $\delta g_0$ and $R_0$, the image $\vec{u}'=\mathfrak{m}(u)$ is defined as follows. For $n\ge 1$, we let
\[
\delta g'_n=(2-L^{\epsilon})^n\delta g_0+\sum_{j=0}^{n-1}
(2-L^{\epsilon})^{n-1-j} \tilde{\xi}_4(\delta g_j,\mu_j,R_j)
\]
and
\begin{eqnarray*}
R'_n & = & \tilde{\mathcal{L}}^{(\delta g_{n-1},\mu_{n-1})}\circ\cdots\circ 
\tilde{\mathcal{L}}^{(\delta g_{0},\mu_{0})}(R_0) \\
 & & +\sum_{j=0}^{n-1} \tilde{\mathcal{L}}^{(\delta g_{n-1},\mu_{n-1})}\circ\cdots\circ 
\tilde{\mathcal{L}}^{(\delta g_{j+1},\mu_{j+1})}\left(\tilde{\xi}_R(\delta g_j,\mu_j,R_j)\right)\ .
\end{eqnarray*}
For $n\ge 0$, we let
\[
\mu'_{n}=-\sum_{j=n}^{\infty}L^{-(j-n+1)\left(\frac{3+\epsilon}{2}\right)} \tilde{\xi}_\mu(\delta g_j,\mu_j,R_j)\ .
\]
Given a sufficiently small $\rho>0$ we now show that this map is well defined and analytic on the open ball $B(\vec{0},\rho)\in\mathcal{B}_+$
in the regime of small $\epsilon$ (made small after fixing $\rho$).

\begin{Proposition}\label{stabmanprop} 
If $0<\rho<\frac{1}{12}$, $|\delta g_0|<\frac{\rho}{12}\bar{g}^{e_4}$ and $|||R_0|||_{\bar{g}}
<\frac{\rho}{8}\bar{g}^{e_R}$ then the map $\mathfrak{m}$ is well defined, analytic on $B(\vec{0},\rho)$
and takes its values in the closed ball $\bar{B}(\vec{0},\frac{\rho}{4})$,
provided $\epsilon$ is made sufficiently small after fixing $\rho$.
Moreover, $\mathfrak{m}$ is jointly analytic in $\vec{u}$ and the implicit variables $\delta g_0$ and $R_0$.
\end{Proposition}

\noindent{\bf Proof:}
Recall the choice of constraints (\ref{econst1eq}), (\ref{econst2eq}), (\ref{econst3eq}).
Their purpose is to ensure that the hypothesis $||\vec{u}||<\rho$ guarantees that all triples $(\delta g_j,\mu_j,R_j)$
featuring in the definition of $\mathfrak{m}(\vec{u})$ are in the domain of definition and analyticity of $\tilde{\mathcal{L}}$
and the $\tilde{\xi}$
specified in Theorem \ref{mainestthm}.
Indeed, since $\epsilon$ which controls the size of $\bar{g}$ will be made small, one may assume $\bar{g}\le 1$ and thus
$\bar{g}^{e_4}\le\bar{g}$, $\bar{g}^{e_2}\le\bar{g}^{1-\eta}$ and $\bar{g}^{e_R}\le\bar{g}^{\frac{11}{4}-\eta_R}$.
Note that for $\epsilon>0$ small we have $0<2-L^{\epsilon}<1$.
Hence for all $n\ge 1$,
\[
|\delta g'| \le (2-L^{\epsilon})^n|\delta g_0|+\sum_{j=0}^{n-1}
(2-L^{\epsilon})^{n-1-j} |\tilde{\xi}_4(\delta g_j,\mu_j,R_j)|\ .
\]
From the definition of $\tilde{\xi}_4$, the hypothesis and the bounds provided by Theorem \ref{mainestthm} we have
\[
|\tilde{\xi}_4(\delta g_j,\mu_j,R_j)|\le A_1|\delta g_j|^2+ |\xi_4(\bar{g}+\delta g_j,\mu_j,R_j)|
\le A_1\rho^2\bar{g}^{2e_4}+B_4\rho \bar{g}^{e_R}
\]
and consequently, using $(2-L^{\epsilon})^n\le 1$ for the first term,
\begin{eqnarray*}
\bar{g}^{-e_4} |\delta g'| & \le & (2-L^{\epsilon})^n |\delta g_0|\bar{g}^{-e_4}+\sum_{j=0}^{n-1}
(2-L^{\epsilon})^{n-1-j} \left[A_1\rho^2\bar{g}^{e_4}+B_4\rho \bar{g}^{e_R-e_4}\right] \\
 & \le &  |\delta g_0|\bar{g}^{-e_4}+
\left[A_1\rho^2\bar{g}^{e_4}+B_4\rho \bar{g}^{e_R-e_4}\right]\times\frac{1-(2-L^{\epsilon})^n}{1-(2-L^{\epsilon})} \\
 & \le & |\delta g_0|\bar{g}^{-e_4}+
\left[A_1\rho^2\bar{g}^{e_4}+B_4\rho \bar{g}^{e_R-e_4}\right]\times \frac{1}{L^{\epsilon}-1}\\
 & \le & |\delta g_0|\bar{g}^{-e_4}+
\rho^2\bar{g}^{e_4-1}+A_1^{-1}B_4\rho \bar{g}^{e_R-e_4-1}
\end{eqnarray*}
where in the last line we invoqued the relation $L^{\epsilon}-1=A_1\bar{g}$.
From the hypothesis on $\delta g_0$ we then have
\[
\bar{g}^{-e_4} |\delta g'|  \le \frac{\rho}{12}+
\rho^2\bar{g}^{e_4-1}+A_1^{-1}B_4\rho \bar{g}^{e_R-e_4-1}\ .
\]
Using $\bar{g}\le 1$, $e_4\ge 1$ and $\rho<\frac{\rho}{12}$ we get
\[
\rho^2\bar{g}^{e_4-1}\le \frac{1}{12}\ .
\]
We have $\lim_{\epsilon\rightarrow 0} A_1=36(1-p^{-3})l$, with $l\ge 1$. Since $\epsilon$ will be made
as small as necessary we may assume, e.g., $A_{1,{\rm min}}\le A_1\le A_{1,{\rm max}}$
where $A_{1,{\rm min}}=35(1-p^{-3})l$ and $A_{1,{\rm max}}=37(1-p^{-3})l$.
Then 
\[
A_1^{-1}B_4\rho \bar{g}^{e_R-e_4-1}\le A_{1,{\rm min}}^{-1}B_4\rho \bar{g}^{e_R-e_4-1}<\frac{\rho}{12}
\]
for $\epsilon$ or equivalently $\bar{g}$ small enough because of the requirement (\ref{econst4eq}). 
As a result $\bar{g}^{-e_4} |\delta g'|  \le \frac{\rho}{4}$.

We now bound $R'_n$ using the property that the operator norms of the $\tilde{\mathcal{L}}$ is at most $\frac{1}{2}$.
Indeed, from $\xi_R$ bound provided by Theorem \ref{mainestthm},
\[
|||R'_n|||_{\bar{g}}\le 2^{-n}|||R_0|||_{\bar{g}}+\sum_{j=0}^{n-1} 2^{-(n-1-j)}
B_{R\xi}\bar{g}^{\frac{11}{4}-3\eta}
\]
and therefore, bounding $2^{-n}$ simply by $1$,
\[
\bar{g}^{-e_R}|||R'_n|||_{\bar{g}}\le \bar{g}^{-e_R}|||R_0|||_{\bar{g}}
+2B_{R\xi}\bar{g}^{\frac{11}{4}-3\eta-e_R}\ . 
\]
Now by hypothesis $\bar{g}^{-e_R}|||R_0|||_{\bar{g}}<\frac{\rho}{8}$ and from (\ref{econst5eq})
we see that
\[
2B_{R\xi}\bar{g}^{\frac{11}{4}-3\eta-e_R}\le \frac{\rho}{8}
\]
when $\epsilon$ is small enough. Thus we also get $\bar{g}^{-e_R}|||R'_n|||_{\bar{g}}\le\frac{\rho}{4}$.

Finally, we bound $\mu'_n$ noting that
\[
|\mu'_{n}|\le\sum_{j=n}^{\infty}L^{-(j-n+1)\left(\frac{3+\epsilon}{2}\right)}
|\tilde{\xi}_\mu(\delta g_j,\mu_j,R_j)|\ . 
\]
But we have, using $|\delta g_j|<\frac{1}{2}\bar{g}$ and the $\xi_\mu$ bound from Theorem \ref{mainestthm},
\[
|\tilde{\xi}_\mu(\delta g_j,\mu_j,R_j)|\le |A_2|\frac{9}{4}\bar{g}^2+
+|A_3|\frac{3}{2}\bar{g}\times\rho \bar{g}^{e_2}+B_2\rho\bar{g}^{e_R}\ .
\]
Bounding $L^{-\left(\frac{3+\epsilon}{2}\right)}$ by $L^{-\frac{3}{2}}$ we immediately obtain
\[
\bar{g}^{-e_2}|\mu'_{n}|\le \frac{L^{-\frac{3}{2}}}{1-L^{-\frac{3}{2}}}\times
\left[
\frac{9}{4}|A_2|\bar{g}^{2-e_2}+
+\frac{3}{2}|A_3|\rho \bar{g}+B_2\rho\bar{g}^{e_R-e_2}
\right]\ .
\]
However, from the definition of $A_2$, we have
\begin{eqnarray*}
|A_2| & \le & 4 L^{\frac{3+\epsilon}{2}}||\Gamma||_{L^3}+144 L^{\epsilon} C_0(0)||\Gamma||_{L^2} \\
 & \le & \left[4 L^{2}||\Gamma||_{L^\infty}+144 L^{\epsilon}\times 2\right]\times  ||\Gamma||_{L^2} \\
 & \le & \left[4 L^{2}||\Gamma||_{L^\infty}+144 L^{\epsilon}\times 2\right]\times \frac{1}{36}L^{-\epsilon} A_1 \\
 & \le & A_{2,{\rm max}}
\end{eqnarray*}
with
\[
A_{2,{\rm max}}=2\times [4+144]\times \frac{1}{36}\times A_{1,{\rm max}}\ .
\]
Note that we used our previous bounds on $C_0(0)$ and $||\Gamma||_{L^\infty}$ by $2$ as well as $0<\epsilon\le 1$
in order to eliminate $\epsilon$ from the exponents of $L$.
We also have from the definition of $A_3$ that
\[
|A_3|\le 12 L^{2}||\Gamma||_{L^2}\le A_{3,{\rm max}}
\]
with
\[
A_{3,{\rm max}}= 12 L^{2}\times \frac{1}{36}\times A_{1,{\rm max}}\ .
\]
We now get
\[
\bar{g}^{-e_2}|\mu'_{n}|\le \frac{L^{-\frac{3}{2}}}{1-L^{-\frac{3}{2}}}\times
\left[
\frac{9}{4}A_{2,\rm max}\bar{g}^{2-e_2}+
+\frac{3}{2}A_{3,\rm max}\rho \bar{g}+B_2\rho\bar{g}^{e_R-e_2}
\right]\le \frac{\rho}{4}
\]
for $\epsilon$ small because $e_2<2<e_R$ as follows from (\ref{econst1eq}), (\ref{econst4eq}) and (\ref{econst6eq}).

When showing the absolute convergence of the series for the $\mu'_n$ we proved that the map $\mathfrak{m}$
is well defined. Analyticity follows easily from uniform absolute convergence. The previous estimates
show that $||\vec{u}||<\rho$ implies $||\mathfrak{m}(\vec{u})||\le\frac{\rho}{4}$.
\qed

Using Lemma \ref{Lipschitzlem} with $r_1=\rho$, $r_2=\frac{\rho}{4}$ and $\nu=\frac{1}{3}$
so that
\[
\frac{r_2(1-\nu)}{r_1(1-2\nu)}=\frac{1}{2}
\]
we immediately see that, under the hypotheses of Proposition \ref{stabmanprop},
the closed ball $\bar{B}\left(\vec{0},\frac{\rho}{3}\right)$ is stable by $\mathfrak{m}$
and is a contraction. More precisely, for any $\vec{u}_1$ and $\vec{u}_2$ in that ball, we have
\[
||\mathfrak{m}(\vec{u}_1)-\mathfrak{m}(\vec{u}_2)||\le\frac{1}{2}||\vec{u}_1-\vec{u}_2 ||\ .
\] 
By the Banach Fixed Point Theorem we then have the existence of a unique fixed point denoted by $\vec{u}_\ast$
for the map $\mathfrak{m}$ in the ball $\bar{B}\left(\vec{0},\frac{\rho}{3}\right)$.
Using the representation of this fixed point as
\[
\vec{u}_\ast=\sum_{n=0}^{\infty}\left[\mathfrak{m}^{n+1}(\vec{0})-\mathfrak{m}^{n}(\vec{0})\right]
\]
and by uniform absolute convergence, it is easy to see that $\vec{u}_\ast$ is analytic in the implicit
data $(\delta g_0, R_0)$. In particular the $\mu_0$ component of the sequence $\vec{u}_\ast$
which we will denote by $\mu_{\rm s}(\delta g_0,R_0)$ is analytic on the domain
given by $|\delta g_0|<\frac{\rho}{12}\bar{g}^{e_4}$ and $|||R_0|||_{\bar{g}}<\frac{\rho}{8}\bar{g}^{e_R}$.

We will now show that, for elements $v=(\delta g,\mu,R)\in\mathcal{E}$, the equation
$\mu=\mu_{\rm s}(\delta g,R)$ characterizes those on the stable manifold of the sought for fixed point.
We now define a set $W^{\rm s,loc}$ which will be our candidate for this local stable manifold.
It will be defined in terms the radius $\rho$ which is supposed to satisfy the hypothesis of Proposition \ref{stabmanprop}.
We let
\[
W^{\rm s,loc}=\left\{
(\delta g,\mu,R)\in\mathcal{E}|\ |\delta g|\le\frac{\rho}{13}\bar{g}^{e_4}, |||R|||_{\bar{g}}\le\frac{\rho}{13}\bar{g}^{e_R}, 
\mu=\mu_{\rm s}(\delta g,R)
\right\}\ .
\]
We will also need the subset
\[
W_{\rm int}^{\rm s,loc}=\left\{
(\delta g,\mu,R)\in\mathcal{E}|\ |\delta g|<\frac{\rho}{13}\bar{g}^{e_4}, |||R|||_{\bar{g}}<\frac{\rho}{13}\bar{g}^{e_R}, 
\mu=\mu_{\rm s}(\delta g,R)
\right\}\ .
\]

\begin{Proposition}\label{L55prop}
For fixed $\rho\in(0,\frac{1}{12})$ and for $\epsilon$ small enough, an equivalent description of $W^{\rm s,loc}$ is as the set
of triples $(\delta g,\mu,R)\in\mathcal{E}$ that satisfy all of the following properties: 
\begin{itemize}
\item
$|\delta g|\le\frac{\rho}{13}\bar{g}^{e_4}$,
\item
$|||R|||_{\bar{g}}\le\frac{\rho}{13}\bar{g}^{e_R}$,
\item
there exists a sequence $(\delta g_n,\mu_n,R_n)_{n\ge 0}$ in $\mathcal{E}$ such that $\delta g_0=\delta g$, $\mu_0=\mu$, $R_0=R$,
$\forall n\ge 1$, $|\delta g_n|\le\frac{\rho}{3}\bar{g}^{e_4}$ and $|||R_n|||_{\bar{g}}\le\frac{\rho}{3}\bar{g}^{e_R}$,
$\forall n\ge 0$, $|\mu_n|\le\frac{\rho}{3}\bar{g}^{e_2}$, and
$\forall n\ge 0$, $(\delta g_{n+1},\mu_{n+1},R_{n+1})=RG(\delta g_n,\mu_n,R_n)$.
\end{itemize} 
\end{Proposition}

\noindent{\bf Proof:}
Suppose $(\delta g,\mu,R)\in W^{\rm s,loc}$. We let $\delta g_0=\delta g$ and $R_0=R$ and consider the fixed point $\vec{u}_\ast$
for the map $\mathfrak{m}$ associated to the data $(\delta g_0,R_0)$ given by Proposition \ref{stabmanprop}.
This makes sense since $\frac{\rho}{13}$ is smaller than $\frac{\rho}{12}$ and $\frac{\rho}{8}$.
Since
\[
\vec{u}_\ast=\left(\mu_0,(\delta g_1,\mu_1,R_1),(\delta g_2,\mu_2,R_2),\ldots\right)\in
\bar{B}\left(\vec{0},\frac{\rho}{3}\right)\ ,
\]
the $(\delta g_n,\mu_n,R_n)$, $n\ge 0$, are well-defined, belong to the domain of definition of the map $RG$
and satisfy the wanted bounds. We just need to check that this sequence forms a trajectory for $RG$.
From $\vec{u}_{\ast}=\mathfrak{m}(\vec{u}_\ast)$ we get, for all $n\ge 1$,
\begin{equation}
\delta g_n=(2-L^{\epsilon})^n\delta g_0+\sum_{j=0}^{n-1}
(2-L^{\epsilon})^{n-1-j} \tilde{\xi}_4(\delta g_j,\mu_j,R_j)\ .
\label{ufpforgeq}
\end{equation}
If $n=1$, (\ref{ufpforgeq}) reduces to the wanted equation, namely,
\[
\delta g_1=(2-L^{\epsilon})\delta g_0+ \tilde{\xi}_4(\delta g_0,\mu_0,R_0)\ .
\]
If $n\ge 2$, (\ref{ufpforgeq}) can be rewritten as 
\[
\delta g_n=(2-L^{\epsilon})\left[
(2-L^{\epsilon})^{n-1}\delta g_0+\sum_{j=0}^{n-2}
(2-L^{\epsilon})^{n-2-j} \tilde{\xi}_4(\delta g_j,\mu_j,R_j)
\right]+\tilde{\xi}_4(\delta g_{n-1},\mu_{n-1},R_{n-1})
\]
\[
=(2-L^{\epsilon})\delta g_{n-1}+\tilde{\xi}_4(\delta g_{n-1},\mu_{n-1},R_{n-1})
\]
by (\ref{ufpforgeq}) for $n-1$ instead of $n$.
Likewise, the $R$ projections of the sequence fixed point equation $\vec{u}_{\ast}=\mathfrak{m}(\vec{u}_\ast)$
imply by similar manipulations that, for all $n\ge 1$,
\[
R_n  =  \tilde{\mathcal{L}}^{(\delta g_{n-1},\mu_{n-1})(R_{n-1})} +\tilde{\xi}_{R}(\delta g_{n-1},\mu_{n-1},R_{n-1})\ .
\]
Now for the $\mu$'s we first write, for all $n\ge 0$,
\begin{equation}
\mu_n=-\sum_{j=n}^{\infty}L^{-(j-n+1)\left(\frac{3+\epsilon}{2}\right)} \tilde{\xi}_\mu(\delta g_j,\mu_j,R_j)
\label{ufpformueq}
\end{equation}
as results from $\vec{u}_{\ast}=\mathfrak{m}(\vec{u}_\ast)$.
Hence
\[
\mu_n=-L^{-\left(\frac{3+\epsilon}{2}\right)}\tilde{\xi}_\mu(\delta g_n,\mu_n,R_n)
-L^{-\left(\frac{3+\epsilon}{2}\right)}
\sum_{j=n+1}^{\infty}L^{-(j-(n+1)+1)\left(\frac{3+\epsilon}{2}\right)} \tilde{\xi}_\mu(\delta g_j,\mu_j,R_j)
\]
\[
=L^{-\left(\frac{3+\epsilon}{2}\right)}\tilde{\xi}_\mu(\delta g_n,\mu_n,R_n)
-L^{-\left(\frac{3+\epsilon}{2}\right)}\mu_{n+1}
\]
by (\ref{ufpformueq}) for $n+1$ instead of $n$.
Thus
\[
\mu_{n+1}=L^{\frac{3+\epsilon}{2}}\mu_n+\tilde{\xi}_\mu(\delta g_n,\mu_n,R_n)\ .
\]
We therefore proved that for all $n\ge 0$, $(\delta g_{n+1},\mu_{n+1},R_{n+1})=RG(\delta g_n,\mu_n,R_n)$
and consequently all the requirements in the statement of the proposition are satisfied.

We now prove the converse and assume that $(\delta g,\mu,R)$ satisfies the listed properties.
We then define $\vec{u}$ using the given $RG$ trajectory $(\delta g_n,\mu_n,R_n)_{n\ge 0}$, simply by setting
\[
\vec{u}=\left(\mu_0,(\delta g_1,\mu_1,R_1),(\delta g_2,\mu_2,R_2),\ldots\right)\ .
\] 
By hypothesis, one clearly has $\vec{u}\in\bar{B}\left(\vec{0},\frac{\rho}{3}\right)$.
For any $n\ge 1$, we have
\[
\delta g_n=(2-L^{\epsilon})\delta g_{n-1}+\tilde{\xi}_4(\delta g_{n-1},\mu_{n-1},R_{n-1})\ .
\]
We apply this to $n-1$ instead of $n$ and substitute in the first term of the previous equation only. We do the same for $n-2$
in the resuting equation and continue this backwards iteration. This immediately establishes (\ref{ufpforgeq}).
The same argument also shows the $R$ parts of the sequence fixed point equation.
As for the $\mu's$, we have for all $n\ge 0$
\[
\mu_{n+1}=L^{\frac{3+\epsilon}{2}}\mu_n+\tilde{\xi}_\mu(\delta g_n,\mu_n,R_n)
\]
which can be rewritten as
\[
\mu_n=-L^{-\left(\frac{3+\epsilon}{2}\right)}\mu_{n+1}+
L^{-\left(\frac{3+\epsilon}{2}\right)}\tilde{\xi}_2(\delta g_n,\mu_n,R_n)\ .
\]
We apply this to $n+1$ instead of $n$ and substitute in the first term of the previous equation.
Iterating this procedure forward $k$ times we get
\[
\mu_n=-\sum_{j=n}^{n+k-1} L^{-(j-n+1)\left(\frac{3+\epsilon}{2}\right)}\tilde{\xi}_2(\delta g_j,\mu_j,R_j)
+L^{-k\left(\frac{3+\epsilon}{2}\right)}\mu_{n+k}\ .
\]
Since by hypothesis the $\mu_j$'s are bounded by $\frac{\rho}{3}\bar{g}^{e_2}$,
\[
\lim_{k\rightarrow\infty} L^{-k\left(\frac{3+\epsilon}{2}\right)}\mu_{n+k}=0
\]
and the $\mu$ part of the sequence fixed point equation holds.
We therefore proved $\vec{u}=\mathfrak{m}(\vec{u})$. By the uniqueness part of the Banach Fixed Point Theorem,
$\vec{u}$ and $\vec{u}_\ast$ are equal and therefore so are their $\mu_0$ components.
Given the previous definition, this establishes $\mu=\mu_{\rm s}(\delta g,R)$ and finally $(\delta g,\mu,R)
\in W^{\rm s,loc}$ as wanted.
\qed

\begin{Proposition}\label{L56prop}
For fixed $\rho\in(0,\frac{1}{12})$ and for $\epsilon$ small enough,
$W^{\rm s,loc}$ is stable by $RG$. In fact one has the stronger statement
$RG\left(W^{\rm s,loc}\right)\subset W_{\rm int}^{\rm s,loc}$.
\end{Proposition}

\noindent{\bf Proof:}
We use the characterization provided by Proposition \ref{L55prop} both ways.
Let $(\delta g,\mu,R)\in W^{\rm s,loc}$ and let $(\delta g_n,\mu_n,R_n)_{n\ge 0}$
be the trajectory provided by Proposition \ref{L55prop}
such that $(\delta g_0,\mu_0,R_0)=(\delta g,\mu,R)$.
We will apply the reverse direction of Proposition \ref{L55prop}
to $(\delta g_1,\mu_1,R_1)$ together with $(\delta g_{n+1},\mu_{n+1},R_{n+1})_{n\ge 0}$,
the shifted trajectory,
in order to show $(\delta g_1,\mu_1,R_1)=RG(\delta g,\mu,R)\in W^{\rm s,loc}$.
All we need is to show the more restrictive inequalities
$|\delta g_1|\le \frac{\rho}{13}\bar{g}^{e_4}$ and $|||R_1|||_{\bar{g}}\le\frac{\rho}{13}\bar{g}^{e_R}$.
In fact we will prove the stronger estimates $|\delta g_1|< \frac{\rho}{13}\bar{g}^{e_4}$ and
$|||R_1|||_{\bar{g}}<\frac{\rho}{13}\bar{g}^{e_R}$ which will show $(\delta g_1,\mu_1,R_1)$ belongs
to  $W_{\rm int}^{\rm s,loc}$ by definition of the latter.

From
\[
R_1=\tilde{\mathcal{L}}^{(\delta g _0,\mu_0)}(R_0)+\tilde{\xi}_{R}(\delta g _0,\mu_0,R_0)
\]
we get
\[
|||R_1|||_{\bar{g}}\le \frac{1}{2}|||R_0|||_{\bar{g}}+B_{R\xi}\bar{g}^{\frac{11}{4}-3\eta}\ .
\]
Hence
\[
\bar{g}^{-e_R}|||R_1|||_{\bar{g}}\le \frac{\rho}{26}+B_{R\xi}\bar{g}^{\frac{11}{4}-3\eta-e_R}<\frac{\rho}{13}
\]
when $\epsilon$ and therefore $\bar{g}$ are small enough, because of the hypothesis (\ref{econst5eq}).

From
\[
\delta g_1=(2-L^{\epsilon})\delta g_0+\tilde{\xi}_4(\delta g _0,\mu_0,R_0)\ ,
\]
i.e.,
\[
\delta g_1=(2-L^{\epsilon})\delta g_0-A_1\delta g_0^2+\xi_4(\bar{g}+\delta g _0,\mu_0,R_0)\ ,
\]
we obtain, using $A_1>0$,
\[
|\delta g_1|=(2-L^{\epsilon})|\delta g_0|+A_1|\delta g_0|^2+|\xi_4(\bar{g}+\delta g _0,\mu_0,R_0)|
\]
\[
\le (2-L^{\epsilon})\left(\frac{\rho}{13}\bar{g}^{e_4}\right)
+A_1\left(\frac{\rho}{13}\bar{g}^{e_4}\right)^2+B_4\times\frac{\rho}{13}\bar{g}^{e_R}
=\frac{\rho}{13}\bar{g}^{e_4}+\Omega
\]
with
\[
\Omega=-(L^{\epsilon}-1)\frac{\rho}{13}\bar{g}^{e_4}+
A_1\left(\frac{\rho}{13}\bar{g}^{e_4}\right)^2+B_4\times\frac{\rho}{13}\bar{g}^{e_R}\ .
\]
Using the relation $A_1\bar{g}=L^{\epsilon}-1$, we have
\[
\Omega=-A_1\frac{\rho}{13}\bar{g}^{e_4+1}+
A_1\left(\frac{\rho}{13}\right)^2\bar{g}^{2e_4}+B_4\times\frac{\rho}{13}\bar{g}^{e_R}\ .
\]
From (\ref{econst1eq}) we have $\bar{g}^{2e_4}\le\bar{g}^{e_4+1}$ and therefore
\[
\Omega\le -A_1\frac{\rho}{13}\bar{g}^{e_4+1}+
A_1\left(\frac{\rho}{13}\right)^2\bar{g}^{e_4+1}+B_4\frac{\rho}{13}\bar{g}^{e_R}
=\frac{\rho}{13}\bar{g}^{e_4+1}\left[-A_1\left(1-\frac{\rho}{13}\right)+B_4\bar{g}^{e_R-e_4-1}\right]\ .
\]
Since $A_1\ge A_{1,{\rm min}}$ defined in the proof or Proposition \ref{stabmanprop} and since $\rho<\frac{1}{12}$
we have
\[
A_1\left(1-\frac{\rho}{13}\right)\ge \frac{155}{156}A_{1,{\rm min}}\ .
\]
As a result
\[
|\delta g_1|\le \frac{\rho}{13}\bar{g}^{e_4}+
\frac{\rho}{13}\bar{g}^{e_4+1}\left[-\frac{155}{156}A_{1,{\rm min}}+B_4\bar{g}^{e_R-e_4-1}\right]\ .
\]
Because of the assumption (\ref{econst1eq}) we have
\[
-\frac{155}{156}A_{1,{\rm min}}+B_4\bar{g}^{e_R-e_4-1}<0
\]
for $\epsilon$ small enough
and thus $|\delta g_1|<\frac{\rho}{13}\bar{g}^{e_4}$ as desired.
\qed

\subsection{A dichotomy lemma}

We now prove an important lemma which gives quantitative growth or decay estimates
which provide a separation between expanding and contracting directions.
We first introduce some notation. Clearly $\mathcal{E}=\mathcal{E}_1\oplus \mathcal{E}_2$
where
\[
\mathcal{E}_1=\{(\delta g,0,R)| \delta g\in\mathbb{C}, R\in C_{\rm bd,ev}^{9}(\mathbb{R},\mathbb{C})\}
\]
and
\[
\mathcal{E}_2=\{(0,\mu,0)| \mu\in\mathbb{C}\}\ . 
\]
We denote by $v_1$ and $v_2$ the pieces of the unique decomposition $v=v_1+v_2$ of an element $v\in\mathcal{E}$.
Note that we will commit a slight abuse of notation by writing
$v_1=(\delta g,R)$ and $v_2=\mu$ if $v=(\delta g,\mu,R)$
or, in other words, by making use of the identifications
$\mathcal{E}_1\simeq\mathbb{C}\times C_{\rm bd,ev}^{9}(\mathbb{R},\mathbb{C})$ and $\mathcal{E}_2\simeq \mathbb{C}$
as Banach spaces. In particular the norms we will be using all come from that of $\mathcal{E}$ and thus ultimately from that
of $\mathcal{E}_{\rm ex}$ in \S\ref{funcspsec}.
For example, following up on the previous set-up we have
\[
||v_1||=\max\left[|\delta g|\bar{g}^{-e_4},|||R|||_{\bar{g}}\bar{g}^{-e_R}\right]\ \ {\rm and}
\ \ ||v_2||=|\mu|\bar{g}^{-e_2}\ .
\]
Finally if $v$ is in the domain of definition for the map $RG$ we write $RG_1(v)=[RG(v)]_1$ and $RG_2(v)=[RG(v)]_2$
for better readability.
Our dichotomy lemma, in the spirit of \cite[Lemma 2.2]{Irwin1} is the following result.

\begin{Lemma}\label{L57lem}
There exists $\epsilon_{0}>0$ and functions $c_1(\epsilon)$, $c_2(\epsilon)$, $c_3(\epsilon)$, $c_4(\epsilon)$,
on $(0,\epsilon_0)$
which satisfy $0<c_1(\epsilon)<1$, $L^{\frac{3}{4}}\ge c_2(\epsilon)>1$,
$2L^{\frac{3}{2}}\ge c_3(\epsilon)\ge L^{\frac{3}{2}}$ and $0<c_4(\epsilon)<1$
(in fact $\lim_{\epsilon\rightarrow 0} c_4(\epsilon)=0$) on that interval such that for all 
$v,v'\in \bar{B}\left(0,\frac{1}{8}\right)\subset\mathcal{E}$ the following statements hold:
\begin{enumerate}
\item
unconditionally, $||RG_1(v)-RG_1(v')||\le c_1(\epsilon)||v-v'||$;
\item
if $L^{\frac{3}{4}}||v_2-v'_2||\ge||v_1-v'_1||$ then
$||RG_2(v)-RG_2(v')||\ge c_2(\epsilon)||v-v'||$;
\item
unconditionally, $||RG_2(v)-RG_2(v')||\le c_3(\epsilon)||v-v'||$;
\item
unconditionally,
\[
||RG_2(v)-RG_2(v')-L^{\frac{3+\epsilon}{2}}(v_2-v'_2)||\le c_4(\epsilon)||v-v'||\ .
\]
\end{enumerate}
More explicitly, the $c(\epsilon)$ functions
are given by the formulas
\begin{eqnarray*}
c_1(\epsilon) & = & \max\left[
1-\frac{3}{4}(L^{\epsilon}-1)+2B_4\bar{g}^{e_R-e_4},
\frac{3}{4}+3B_{R\xi}\bar{g}^{\frac{11}{4}-3\eta-\eta_R}
\right] \\
c_2(\epsilon) & = & L^{\frac{3}{4}}-\frac{9}{4}A_{2,{\rm max}}\bar{g}^{e_4-e_2+1}
-\frac{5}{4}A_{3,{\rm max}}\bar{g}-2B_2\bar{g}^{e_R-e_2} \\
c_3(\epsilon)  & = & L^{\frac{3+\epsilon}{2}}+\frac{9}{4}A_{2,{\rm max}}\bar{g}^{e_4-e_2+1}
+\frac{5}{4}A_{3,{\rm max}}\bar{g}+2B_2\bar{g}^{e_R-e_2} \\
c_4(\epsilon)  & = & \frac{9}{4}A_{2,{\rm max}}\bar{g}^{e_4-e_2+1}
+\frac{5}{4}A_{3,{\rm max}}\bar{g}+2B_2\bar{g}^{e_R-e_2}\ . 
\end{eqnarray*}
\end{Lemma}

\noindent{\bf Proof:}
Since $\frac{1}{8}<\frac{1}{2}$, $\bar{g}\le 1$ for $\epsilon$ small, and because of (\ref{econst1eq}),
(\ref{econst2eq}) and (\ref{econst3eq}), $v$ and $v'$ are in the domain of definition of $RG$ as provided by
Theorem \ref{mainestthm}. Let $v=(\delta g,\mu,R)$, $v'=(\delta g',\mu',R')$,
$RG(v)=(\widehat{\delta g},\widehat{\mu},\widehat{R})$ and
$RG(v')=(\widehat{\delta g}',\widehat{\mu}',\widehat{R}')$.
From the formulas defining the bulk RG transformation we have
\[
\widehat{\delta g}-\widehat{\delta g}'=(2-L^{\epsilon})(\delta g-\delta g')
-A_1(\delta g^2-{\delta g'}^2)+\xi_4(\bar{g}+\delta g,\mu,R)-\xi_4(\bar{g}+\delta g',\mu',R')
\] 
and thus
\[
|\widehat{\delta g}-\widehat{\delta g}'|\le (2-L^{\epsilon})|\delta g-\delta g'|
+A_1|\delta g-\delta g'|(|\delta g|+|\delta g'|)
+|\xi_4(\bar{g}+\delta g,\mu,R)-\xi_4(\bar{g}+\delta g',\mu',R')|\ .
\]
Using Lemma \ref{Lipxilem} we therefore obtain
\[
|\widehat{\delta g}-\widehat{\delta g}'|\le (2-L^{\epsilon})|\delta g-\delta g'|
+A_1|\delta g-\delta g'|\times 2\times \frac{1}{8}\bar{g}^{e_4}
+2B_4\bar{g}^{e_R}||v-v'||\ .
\] 
By definition of the norm on $\mathcal{E}$ we then get
\[
\bar{g}^{-e_4}|\widehat{\delta g}-\widehat{\delta g}'|\le
||v-v'||\times\left\{
2-L^{\epsilon}+\frac{1}{4}A_1 \bar{g}^{e_4}+2 B_4\bar{g}^{e_R-e_4}
\right\}
\]
\[
=||v-v'||\times\left\{1-(L^{\epsilon}-1)\left(1-\frac{1}{4}\bar{g}^{e_4-1}\right)+2 B_4\bar{g}^{e_R-e_4}
\right\}
\]
where we used $A_1\bar{g}=L^{\epsilon}-1$. Also using (\ref{econst1eq}) we get the simpler bound
\[
\bar{g}^{-e_4}|\widehat{\delta g}-\widehat{\delta g}'|
\le ||v-v'||\times\left\{1-\frac{3}{4}(L^{\epsilon}-1)+2 B_4\bar{g}^{e_R-e_4}
\right\}\ .
\]

We now turn to the $R$ components and deduce from formulas for the bulk RG and Lemma \ref{Lipxilem}
that
\[
|||\widehat{R}-\widehat{R}'|||_{\bar{g}}
\le 
|||{\mathcal{L}}^{(\bar{g}+\delta g,\mu)}(R)-
{\mathcal{L}}^{(\bar{g}+\delta g',\mu')}(R')|||_{\bar{g}}+
|||\xi_R(\bar{g}+\delta g,\mu,R)-\xi_R(\bar{g}+\delta g',\mu',R')|||_{\bar{g}}
\]
\[
\le
\left(\frac{3}{4}\bar{g}^{e_R}+3B_{R\xi}\bar{g}^{\frac{11}{4}-3\eta}\right) ||v-v'||\ .
\]
This together with the previous estimate on the $\delta g$ part and the definition of the $||\cdot||$
on $\mathcal{E}_1$ immediately implies Part 1) of the lemma with the given $c_1(\epsilon)$.
What remains is to show that this quantity is in the interval $(0,1)$ for $\epsilon$ sufficiently small.
Note that $\bar{g}$ is of the same order as $\epsilon$ in this regime. 
Therefore $L^{\epsilon}-1\sim \epsilon \log L$ as well as the assumptions (\ref{econst4eq})
and (\ref{econst5eq}) ensure that $c_1(\epsilon)$ has the wanted property.

We now tackle Part 2) and assume $L^{\frac{3}{4}}||v_2-v'_2||\ge ||v_1-v'_1||$.
For the formulas for $RG$ we get
\[
\widehat{\mu}-\widehat{\mu}'=
L^{\frac{3+\epsilon}{2}}(\mu-\mu')
-A_2\left[(\bar{g}+\delta g)^2-(\bar{g}+\delta g)^2\right]
\]
\begin{equation}
-A_3\left[(\bar{g}+\delta g)\mu-(\bar{g}+\delta g')\mu'\right]
+\xi_2(\bar{g}+\delta g,\mu,R)-\xi_2(\bar{g}+\delta g',\mu',R')
\label{mudiffeq}
\end{equation}
and therefore
\[
|\widehat{\mu}-\widehat{\mu}'|\ge
L^{\frac{3+\epsilon}{2}}|\mu-\mu'|
-A_{2,{\rm max}}|\delta g-\delta g'|(2\bar{g}+|\delta g|+|\delta g'|)
\]
\[
-A_{3,{\rm max}}\left|(\bar{g}+\delta g)\mu-(\bar{g}+\delta g')\mu'\right|
-|\xi_2(\bar{g}+\delta g,\mu,R)-\xi_2(\bar{g}+\delta g',\mu',R')|
\]
where $A_{2,{\rm max}}$ and $A_{3,{\rm max}}$ have been defined in the proof of Proposition \ref{stabmanprop}.
We use (\ref{econst1eq}) and the hypothesis on $v$ and $v'$ in order to write
\[
2\bar{g}+|\delta g|+|\delta g'| \le  2\bar{g}+\bar{g}^{e_4}||v||+\bar{g}^{e_4}||v'||\le \bar{g}
\times\left[2+\frac{1}{8}+\frac{1}{8}\right]\ .
\]
We also have, for the same reasons,
\begin{eqnarray*}
\left|(\bar{g}+\delta g)\mu-(\bar{g}+\delta g')\mu'\right| & = & \left|
\bar{g}(\mu-\mu')+\delta g(\mu-\mu')+(\delta g-\delta g')\mu'\right| \\
 & \le & \bar{g}\times \bar{g}^{e_2}||v-v'||+\bar{g}^{e_4}||v||\times  \bar{g}^{e_2}||v-v'||
 + \bar{g}^{e_4}||v-v'||\times \bar{g}^{e_2}||v|| \\
  & \le & \bar{g}^{1+e_2}||v-v'||\times\left[1+\frac{1}{8}+\frac{1}{8}\right]\ .
\end{eqnarray*}
The last ingredients are the use of the $\xi_2$ Lispchitz estimate in Lemma \ref{Lipxilem}
and the simplification $L^{\frac{3+\epsilon}{2}}\ge L^{\frac{3}{2}}$.
Altogether, this gives
\[
||RG_2(v)-RG_2(v')||=
\bar{g}^{-e_2}|\widehat{\mu}-\widehat{\mu}'|\ge
L^{\frac{3}{2}}\bar{g}{-e_2}|\mu-\mu'|
-\frac{9}{4}A_{2,{\rm max}}\bar{g}^{1+e_4-e_2}||v-v'||
\]
\begin{equation}
-\frac{5}{4}A_{3,{\rm max}}\bar{g}||v-v'||
-2 B_2\bar{g}^{e_R-e_2}||v-v'||\ .
\label{mulowbdeq}
\end{equation}
From the hypothesis $L^{\frac{3}{4}}||v_2-v'_2||\ge ||v_1-v'_1||$
we get
\[
||v-v'||=\max\left[||v_1-v'_1||,||v_2-v'_2||\right]\le L^{\frac{3}{4}}||v_2-v'_2||
=L^{\frac{3}{4}}\bar{g}^{-e_2}|\mu-\mu'|\ .
\]
Therefore, losing a fraction $L^{\frac{3}{4}}$ of the initial $L^{\frac{3}{2}}$ factor,
(\ref{mulowbdeq}) implies Part 2) of the lemma with the given function $c_2(\epsilon)$.
From (\ref{econst1eq}) and (\ref{econst6eq}) which also imply $1+e_4-e_2>0$ we see that 
$c_2(\epsilon)\rightarrow L^{\frac{3}{4}}>1$ when $\epsilon\rightarrow 0$ as wanted. 

For Part 3) we use (\ref{mudiffeq})to write
\[
|\widehat{\mu}-\widehat{\mu}'|\le
L^{\frac{3+\epsilon}{2}}|\mu-\mu'|
+A_{2,{\rm max}}|\delta g-\delta g'|(2\bar{g}+|\delta g|+|\delta g'|)
\]
\[
+A_{3,{\rm max}}\left|(\bar{g}+\delta g)\mu-(\bar{g}+\delta g')\mu'\right|
+|\xi_2(\bar{g}+\delta g,\mu,R)-\xi_2(\bar{g}+\delta g',\mu',R')|\ .
\]
Then the same bounds as before on the last three terms give
\[
||RG_2(v)-RG_2(v')||=
\bar{g}^{-e_2}|\widehat{\mu}-\widehat{\mu}'|\le
L^{\frac{3+\epsilon}{2}}\bar{g}{-e_2}|\mu-\mu'|
+\frac{9}{4}A_{2,{\rm max}}\bar{g}^{1+e_4-e_2}||v-v'||
\]
\[
+\frac{5}{4}A_{3,{\rm max}}\bar{g}||v-v'||
+2 B_2\bar{g}^{e_R-e_2}||v-v'||\ .
\]
Since $\bar{g}^{-e_2}|\mu-\mu'|\le ||v-v'||$ holds by definition of the norm, the estimate in Part 3) follows.
The bounds on $c_3(\epsilon)$ in the small $\epsilon$ regime are also immediate.

For Part 4), starting from (\ref{mudiffeq}) we transfer $L^{\frac{3+\epsilon}{2}}(\mu-\mu')$ to the left-hand side
and use the same bounds to arrive at
\[
\bar{g}^{-e_2}|\widehat{\mu}-\widehat{\mu}'-L^{\frac{3+\epsilon}{2}}(\mu-\mu')|\le
\left[\frac{9}{4}A_{2,{\rm max}}\bar{g}^{1+e_4-e_2}
+\frac{5}{4}A_{3,{\rm max}}\bar{g}
+2 B_2\bar{g}^{e_R-e_2}\right]\times ||v-v'||
\]
which is the desired result.
\qed

\subsection{The infrared RG fixed point}

\begin{Lemma}\label{L58lem}
If $v\neq v'$ belong to $W^{\rm s,loc}$ then $||v_1-v'_1||> L^{\frac{3}{4}}||v_2-v'_2||$.
\end{Lemma}

\noindent{\bf Proof:}
Note that by the prevailing assumptions we have $\rho<\frac{1}{12}<\frac{1}{8}$ and thus Lemma \ref{L57lem} is applicable
to all elements of $W^{\rm s,loc}$ and their $RG$ iterates by stability of that set.
We proceed by contradiction and
suppose that $||v_1-v'_1||\le L^{\frac{3}{4}}||v_2-v'_2||$.
Then by Lemma \ref{L57lem} Part 1) and 2)
\[
||RG_1(v)-RG_1(v')||\le c_1(\epsilon)||v-v'||\le c_1(\epsilon)c_2(\epsilon)^{-1}||RG_2(v)-RG_2(v')||\ .
\]
From the bounds we have on $c_1(\epsilon)$ and $c_2(\epsilon)$
we trivially get $c_1(\epsilon)c_2(\epsilon)^{-1}<L^{\frac{3}{4}}$ and therefore
\[
||RG_1(v)-RG_1(v')||\le L^{\frac{3}{4}} ||RG_2(v)-RG_2(v')||\ ,
\]
i.e., the first iterates $RG(v)$ and $RG(v')$ satisfy the same hypothesis as $v$ and $v'$.
By an easy induction we then have
\[
\forall n\ge 0,\ \ ||RG_1^n(v)-RG_1^n(v')||\le L^{\frac{3}{4}} ||RG_2^n(v)-RG_2^n(v')||
\]
for the higher iterates where $RG_1^n(\cdot)$ means $(RG^n(\cdot))_1$ and likewise for the second components.
By Lemma \ref{L57lem} Part 2) we obtain, for all $n\ge 0$,
\[
||RG_2^{n+1}(v)-RG_2^{n+1}(v')||\ge c_2(\epsilon)||RG^n(v)-RG^n(v')||\ge c_2(\epsilon)||RG_2^n(v)-RG_2^n(v')||\ .
\]
Again by a trivial induction we get, for all $n\ge 0$,
\[
||RG_2^n(v)-RG_2^n(v')||\ge c_2(\epsilon)^n ||v_2-v'_2||\ .
\]
But $c_2(\epsilon)>1$, so if $||v_2-v'_2||>0$ we have
\[
\lim_{n\rightarrow\infty} ||RG_2^n(v)-RG_2^n(v')||=\infty
\]
which contradicts the stability and boundedness of the set $W^{\rm s,loc}$.
Therefore $||v_2-v'_2||=0$ which also entails $||v_1-v'_1||=0$ by the assumtion made at the beginning of this proof.
This therefore leads to $v=v'$ which is the desired contradition.
\qed

\begin{Lemma}\label{L59lem}
For all $v,v'\in W^{\rm s,loc}$ we have $||RG(v)-RG(v')||\le c_1(\epsilon)||v-v'||$.
\end{Lemma}

\noindent{\bf Proof:}
By the previous lemma and the stability of $W^{\rm s,loc}$ we have
\[
||RG_2(v)-RG_2(v')||\le L^{-\frac{3}{4}} ||RG_1(v)-RG_1(v')||\le ||RG_1(v)-RG_1(v')||
\]
and therefore
\[
||RG(v)-RG(v')||=||RG_1(v)-RG_1(v')||\ .
\]
As a result, the desired conclusion follows from Lemma \ref{L57lem} Part 1).
\qed

\begin{Proposition}\label{L60prop}
The map $RG$ is a contraction when restricted to $W^{\rm s,loc}$ and thus has a unique fixed point
$v_\ast=(\delta g_\ast,\mu_\ast,R_\ast)$ in that set.
In fact $v_\ast$ belongs to the smaller set $W_{\rm int}^{\rm s,loc}$.
\end{Proposition}

\noindent{\bf Proof:}
Note that $W^{\rm s,loc}$ is a closed subset of the Banach space $\mathcal{E}$.
Indeed, $\mu_{\rm s}$ is analytic and thus continuous on an open domain containing that given by the condition
$||(\delta g,R)||\le\frac{\rho}{13}$. Since $W^{\rm s,loc}$ is therefore a complete metric space for the distance
coming from the $||\cdot||$ norm, and since $RG$ restricted to this set is a contraction
as follows form Lemma \ref{L59lem} and $c_1(\epsilon)<1$, the Banach Fixed Point Theorem establishes the present lemma.
The fixed point is in $W_{\rm int}^{\rm s,loc}$ since $v_\ast$ is its own image by application of the stronger conclusion
of Proposition \ref{L57lem}. 
\qed

\subsection{The unstable manifold}

We now construct the local unstable manifold following a procedure similar to that of \S\ref{stabmansec}.
Let $\mathcal{B}_-$ be the Banach space of sequences
\[
\vec{u}=(\ldots,(\delta g_{-2},\mu_{-2},R_{-2}),(\delta g_{-1},\mu_{-1},R_{-1}),\delta g_0,R_0)\in \prod_{n\le -1}\left[
\mathbb{C}^2\times C_{\rm bd,ev}^{9}(\mathbb{R},\mathbb{C})
\right]\times \mathbb{C}\times C_{\rm bd,ev}^{9}(\mathbb{R},\mathbb{C})
\]
which have finite norm given by
\[
||\vec{u}||=\sup\left\{
|\delta g_j|\bar{g}^{-e_4}\ {\rm for}\ j\le 0;
|\mu_j|\bar{g}^{-e_2}\ {\rm for}\ j\le -1;
|||R_j|||\bar{g}^{-e_R}\ {\rm for}\ j\le 0
\right\}\ .
\]
We will define a map $\mathfrak{n}$ on this space of sequences which depends on the parameter $\mu_0$ serving as boundary conditions.
Given $\mu_0$, the image $\vec{u}'=\mathfrak{n}(u)$ is defined as follows. For $n\le 0 $, we let
\[
\delta g'_n=\sum_{j\le n-1}
(2-L^{\epsilon})^{n-1-j} \tilde{\xi}_4(\delta g_j,\mu_j,R_j)
\]
and
\[
R'_n =\sum_{j\le n-1} \tilde{\mathcal{L}}^{(\delta g_{n-1},\mu_{n-1})}\circ\cdots\circ 
\tilde{\mathcal{L}}^{(\delta g_{j+1},\mu_{j+1})}\left(\tilde{\xi}_R(\delta g_j,\mu_j,R_j)\right)\ .
\]
For $n\le -1$, we let
\[
\mu'_{n}=L^{n\left(\frac{3+\epsilon}{2}\right)}\mu_0
-\sum_{j=n}^{-1}L^{-(j-n+1)\left(\frac{3+\epsilon}{2}\right)} \tilde{\xi}_\mu(\delta g_j,\mu_j,R_j)\ .
\]
Given a sufficiently small $\rho'>0$ we will show that this map is well defined and analytic on the open ball $B(\vec{0},\rho')\in\mathcal{B}_-$
in the regime of small $\epsilon$ (made small after fixing $\rho'$).

\begin{Proposition}\label{unstabmanprop} 
If $0<\rho'\le \frac{1}{8}$, $|\mu_0|<\frac{\rho'}{8}\bar{g}^{e_2}$
then the map $\mathfrak{n}$ is well defined, analytic on $B(\vec{0},\rho')$
and takes its values in the closed ball $\bar{B}(\vec{0},\frac{\rho'}{4})$,
provided $\epsilon$ is made sufficiently small after fixing $\rho'$.
Moreover, $\mathfrak{n}$ is jointly analytic in $\vec{u}$ and the implicit variable $\mu_0$.
\end{Proposition}

\noindent{\bf Proof:}
Again the choice of constraints (\ref{econst1eq}), (\ref{econst2eq}), (\ref{econst3eq}) and
the hypothesis $||\vec{u}||<\rho'<\frac{1}{2}$ guarantees that all triples $(\delta g_j,\mu_j,R_j)$
featuring in the definition of $\mathfrak{n}(\vec{u})$ are in the domain of definition and analyticity of $\tilde{\mathcal{L}}$
and the $\tilde{\xi}$ coming from Theorem \ref{mainestthm}.

Hence for all $n\le 0$,
\[
|\delta g'| \le \sum_{j\le n-1}^{n-1}
(2-L^{\epsilon})^{n-1-j} |\tilde{\xi}_4(\delta g_j,\mu_j,R_j)|\ .
\]
As in the proof of Proposition \ref{stabmanprop}
\[
|\tilde{\xi}_4(\delta g_j,\mu_j,R_j)|
\le A_1{\rho'}^2\bar{g}^{2e_4}+B_4\rho' \bar{g}^{e_R}
\]
and consequently,
\begin{eqnarray*}
\bar{g}^{-e_4} |\delta g'| & \le & \sum_{j\le n-1}
(2-L^{\epsilon})^{n-1-j} \left[A_1{\rho'}^2\bar{g}^{e_4}+B_4\rho' \bar{g}^{e_R-e_4}\right] \\
 & \le &  
\left[A_1{\rho'}^2\bar{g}^{e_4}+B_4\rho' \bar{g}^{e_R-e_4}\right]\times\frac{1}{L^{\epsilon}-1} \\
 & \le & 
\left[A_1{\rho'}^2\bar{g}^{e_4}+B_4\rho' \bar{g}^{e_R-e_4}\right]\times \frac{1}{A_1\bar{g}} \\
 & \le & {\rho'}^2+A_{1,{\rm min}}^{-1}B_4\rho' \bar{g}^{e_R-e_4-1}
\end{eqnarray*}
where we used the relation $L^{\epsilon}-1=A_1\bar{g}$ and later the simple bound $\bar{g}^{e_4-1}\le 1$ due to (\ref{econst1eq}).
The hypothesis on $\rho'$ implies ${\rho'}^2\le \frac{\rho'}{8}$ whereas (\ref{econst4eq})
ensures that $A_{1,{\rm min}}^{-1}B_4\rho' \bar{g}^{e_R-e_4-1}\le \frac{\rho'}{8}$ for $\epsilon$ small.
Therefore, the previous estimates
show that the series defining $\delta g'$ converges and that the latter satisfies the bound $\bar{g}^{-e_4} |\delta g'|\le \frac{\rho'}{4}$
in the small $\epsilon$ regime.

We now bound $R'_n$ using the property that the operator norms of the $\tilde{\mathcal{L}}$ is at most $\frac{1}{2}$.
Indeed, similarly to the proof of Proposition \ref{stabmanprop},
\[
|||R'_n|||_{\bar{g}}\le \sum_{j\le n-1} 2^{-(n-1-j)}
B_{R\xi}\bar{g}^{\frac{11}{4}-3\eta} = 2B_{R\xi}\bar{g}^{\frac{11}{4}-3\eta} \ .
\]
Hence
\[
\bar{g}^{-e_R}|||R'_n|||_{\bar{g}}\le 
2B_{R\xi}\bar{g}^{\frac{11}{4}-3\eta-e_R}\le \frac{\rho'}{4} 
\]
for $\epsilon$ small because of  (\ref{econst5eq}).
We also showed by the same token the convergence in the Banach space $\mathcal{E}$ of the series defining
$R'_n$.

Finally, we bound $\mu'_n$ as in the proof of Proposition \ref{stabmanprop} by writing
\[
|\mu'_{n}|\le L^{n\left(\frac{3+\epsilon}{2}\right)}|\mu_0|+
\sum_{j=n}^{-1 }L^{-(j-n+1)\left(\frac{3+\epsilon}{2}\right)}
|\tilde{\xi}_2(\delta g_j,\mu_j,R_j)|
\]
\[
\le L^{n\left(\frac{3+\epsilon}{2}\right)}|\mu_0|+
\sum_{j=n}^{-1 }L^{-(j-n+1)\left(\frac{3+\epsilon}{2}\right)}
\times \left[
\frac{9}{4} A_{2,{\rm max}}\bar{g}^2+
+\frac{3}{2} A_{3,{\rm max}}\rho' \bar{g}^{e_2+1}+B_2\rho'\bar{g}^{e_R}\right]\ .
\]
Using the simple bound $L^{n\left(\frac{3+\epsilon}{2}\right)}\le 1$, since $n\le -1$,
as well as
\[
\sum_{j=n}^{-1 }L^{-(j-n+1)\left(\frac{3+\epsilon}{2}\right)}\le
\sum_{k=1}^{\infty}L^{-k\left(\frac{3+\epsilon}{2}\right)}\le \sum_{k=1}^{\infty}L^{-k\left(\frac{3}{2}\right)}
\]
we obtain
\[
\bar{g}^{-e_2}|\mu'_{n}|\le \bar{g}^{-e_2} |\mu_0|+
\frac{L^{-\frac{3}{2}}}{1-L^{-\frac{3}{2}}}\times
\times \left[
\frac{9}{4} A_{2,{\rm max}}\bar{g}^{2-e_2}+
+\frac{3}{2} A_{3,{\rm max}}\rho' \bar{g}+B_2\rho'\bar{g}^{e_R-e_2}\right]\ .
\]
The first term is bounded by $\frac{\rho'}{8}$ by hypothesis.
Besides, the second also satisfies the same bound when $\epsilon$ is small enough
because of  $e_2<2<e_R$ which follows from (\ref{econst1eq}), (\ref{econst4eq}) and (\ref{econst6eq}).
Hence $\bar{g}^{-e_2}|\mu'_{n}|\le \frac{\rho'}{4}$.

Again, when showing the uniform absolute convergence of the series for the $\delta g'_n$ and $R'_n$ we proved that the map $\mathfrak{n}$
is well defined, analytic and satisfies the bound $||\mathfrak{n}(\vec{u})||\le\frac{\rho'}{4}$ when $||\vec{u}||<\rho'$.
\qed

Again using Lemma \ref{Lipschitzlem} with $r_1=\rho'$, $r_2=\frac{\rho'}{4}$ and $\nu=\frac{1}{3}$
so that
\[
\frac{r_2(1-\nu)}{r_1(1-2\nu)}=\frac{1}{2}
\]
we see that, under the hypotheses of Proposition \ref{unstabmanprop},
the closed ball $\bar{B}\left(\vec{0},\frac{\rho'}{3}\right)$ is stable by $\mathfrak{n}$
and is a contraction. More precisely, for any $\vec{u}_1$ and $\vec{u}_2$ in that ball, we have
\[
||\mathfrak{n}(\vec{u}_1)-\mathfrak{n}(\vec{u}_2)||\le\frac{1}{2}||\vec{u}_1-\vec{u}_2 ||\ .
\] 
By the Banach Fixed Point Theorem we have the existence of a unique fixed point which we again denote by $\vec{u}_\ast$
for the map $\mathfrak{n}$ in the ball $\bar{B}\left(\vec{0},\frac{\rho'}{3}\right)$.
Using the representation of this fixed point as
\[
\vec{u}_\ast=\sum_{n=0}^{\infty}\left[\mathfrak{n}^{n+1}(\vec{0})-\mathfrak{n}^{n}(\vec{0})\right]
\]
and by uniform absulote convergence, we see that $\vec{u}_\ast$ is analytic in the implicit
data $\mu_0$. In particular the $\delta g_0$, $R_0$ components of the sequence $\vec{u}_\ast$
which we will denote by $\delta g_{\rm u}(\mu_0)$, $R_{\rm u}(\mu_0)$ respectively 
are analytic on the domain
given by $|\mu_0|<\frac{\rho'}{8}\bar{g}^{e_2}$.

As in section \S\ref{stabmansec}, the next step will be to show that, for elements $v=(\delta g,\mu,R)\in\mathcal{E}$, the equation
$(\delta g,R)=(\delta g_{\rm u}(\mu),R_{\rm u}(\mu))$ 
characterizes those on the unstable manifold of the bulk RG fixed point $v_\ast$.
We now define a set $W^{\rm u,loc}$ which will be our candidate for this local unstable manifold.
It will be defined in terms the radius $\rho'$ which is supposed to satisfy the hypothesis of Proposition \ref{unstabmanprop}.
We let
\[
W^{\rm u,loc}=\left\{
(\delta g,\mu,R)\in\mathcal{E}|\ |\mu|<\frac{\rho'}{8}\bar{g}^{e_2}, \delta g=\delta g_{\rm u}(\mu), 
R=R_{\rm u}(\mu)
\right\}\ .
\]

\begin{Proposition}\label{L62prop}
For fixed $\rho'\in(0,\frac{1}{8}]$ and for $\epsilon$ small enough, an equivalent description of $W^{\rm u,loc}$ is as the set
of triples $(\delta g,\mu,R)\in\mathcal{E}$ that satisfy all of the following properties: 
\begin{itemize}
\item
$|\mu|<\frac{\rho'}{8}\bar{g}^{e_2}$,
\item
there exists a sequence $(\delta g_n,\mu_n,R_n)_{n\le 0}$ in $\mathcal{E}$ such that $\delta g_0=\delta g$, $\mu_0=\mu$, $R_0=R$,
$\forall n\le 0$, $|\delta g_n|\le\frac{\rho'}{3}\bar{g}^{e_4}$ and $|||R_n|||_{\bar{g}}\le\frac{\rho'}{3}\bar{g}^{e_R}$,
$\forall n\le -1$, $|\mu_n|\le\frac{\rho'}{3}\bar{g}^{e_2}$, and
$\forall n\le -1$, $(\delta g_{n+1},\mu_{n+1},R_{n+1})=RG(\delta g_n,\mu_n,R_n)$.
\end{itemize} 
\end{Proposition}

\noindent{\bf Proof:}
Suppose $(\delta g,\mu,R)\in W^{\rm u,loc}$. We let $\mu_0=\mu$ and consider the fixed point $\vec{u}_\ast$
for the map $\mathfrak{n}$ associated to the data $\mu_0$ given by Proposition \ref{unstabmanprop}.
We write
\[
\vec{u}_\ast=\left(\ldots,(\delta g_{-2},\mu_{-2},R_{-2}),(\delta g_{-1},\mu_{-1},R_{-1}),\delta g_0,R_0\right)
\in\bar{B}\left(\vec{0},\frac{\rho'}{3}\right)\ ,
\]
and note that the $(\delta g_n,\mu_n,R_n)$, $n\le -1$, are well-defined, belong to the domain of definition of the map $RG$
and satisfy the wanted bounds. 
We need to check that this sequence, also including the $n=0$ term, forms a trajectory for $RG$.
Form $\vec{u}_{\ast}=\mathfrak{n}(\vec{u}_\ast)$ we get, for all $n\le -1$,
\[
L^{\left(\frac{3+\epsilon}{2}\right)}\mu_n+\tilde{\xi}_2(\delta g_n,\mu_n,R_n)
=L^{\left(\frac{3+\epsilon}{2}\right)}\left[
L^{n\left(\frac{3+\epsilon}{2}\right)}\mu_0-\sum_{j=n}^{-1}
L^{-(j-n+1)\left(\frac{3+\epsilon}{2}\right)}\tilde{\xi}_2(\delta g_j,\mu_j,R_j)
\right]+\tilde{\xi}_2(\delta g_n,\mu_n,R_n)
\]
\[
=L^{(n+1)\left(\frac{3+\epsilon}{2}\right)}\mu_0-\sum_{j=n+1}^{-1}
L^{-(j-n)\left(\frac{3+\epsilon}{2}\right)}\tilde{\xi}_2(\delta g_j,\mu_j,R_j)\ .
\]
The last quantity is equal to $\mu_{n+1}$ if $n\le -2$ by the fixed point equation for the sequence $\vec{u}_\ast$.
Otherwise if $n=-1$ the same quantity reduces to $\mu_0=\mu_{n+1}$ because the sum is empty.

Likewise and still for $n\le -1$ we have
\[
(2-L^{\epsilon})\delta g_n+\tilde{\xi}_4(\delta g_n,\mu_n,R_n)=
(2-L^{\epsilon})\left[
\sum_{j\le n-1}
(2-L^{\epsilon})^{n-1-j} \tilde{\xi}_4(\delta g_j,\mu_j,R_j)
\right]+\tilde{\xi}_4(\delta g_n,\mu_n,R_n)
\]
\[
=\sum_{j\le n}
(2-L^{\epsilon})^{n-j} \tilde{\xi}_4(\delta g_j,\mu_j,R_j)=\delta g_{n+1}\ .
\]
Similary, the $R$ projections of the sequence fixed point equation $\vec{u}_{\ast}=\mathfrak{n}(\vec{u}_\ast)$
imply by analogous manipulations that, for all $n\le -1$,
\[
\tilde{\mathcal{L}}^{(\delta g_{n},\mu_{n})(R_{n})} +\tilde{\xi}_{R}(\delta g_{n},\mu_{n},R_{n})=R_{n+1}\ .
\]
We therefore proved that for all $n\le -1$, $(\delta g_{n+1},\mu_{n+1},R_{n+1})=RG(\delta g_n,\mu_n,R_n)$
and consequently all the requirements in the statement of the proposition are satisfied.

We now prove the converse and assume that $(\delta g,\mu,R)$ satisfies the listed properties.
We then define $\vec{u}$ using the given $RG$ trajectory $(\delta g_n,\mu_n,R_n)_{n\le 0}$, simply by setting
\[
\vec{u}=\left(\ldots,(\delta g_{-2},\mu_{-2},R_{-2}),(\delta g_{-1},\mu_{-1},R_{-1}),\delta g_0,R_0\right)\ .
\] 
By hypothesis, we clearly have $\vec{u}\in\bar{B}\left(\vec{0},\frac{\rho'}{3}\right)$.
For any $n\le 0$, we have
\[
\delta g_n=(2-L^{\epsilon})\delta g_{n-1}+\tilde{\xi}_4(\delta g_{n-1},\mu_{n-1},R_{n-1})\ .
\]
We apply this to $n-1$ instead of $n$ and substitute in the first term of the previous equation only. We do the same for $n-2$
in the resulting equation and continue this backwards iteration. 
We thus obtain for any $k\ge 1$,
\[
\delta g_n=(2-L^{\epsilon})^{k}\delta g_{n-k}+\sum_{j=n-k}^{n-1}(2-L^{\epsilon})^{n-1-j}
\tilde{\xi}_4(\delta g_{j},\mu_{j},R_{j})\ .
\]
But $0<2-L^{\epsilon}<1$ and the sequence of $\delta g$'s is bounded and therefore $(2-L^{\epsilon})^{k}\delta g_{n-k}\rightarrow 0$
when $k\rightarrow \infty$. Hence
\[
\delta g_n=\sum_{j\le n-1}(2-L^{\epsilon})^{n-1-j}
\tilde{\xi}_4(\delta g_{j},\mu_{j},R_{j})\ .
\]
A similar argument for the $R$'s gives for all $n\le 0$ and all $k\ge 1$
\[
R_n=\tilde{\mathcal{L}}^{(\delta g_{n-1},\mu_{n-1})}\circ\cdots\circ 
\tilde{\mathcal{L}}^{(\delta g_{n-k},\mu_{n-k})}(R_{n-k})+
\sum_{j=n-k}^{n-1} \tilde{\mathcal{L}}^{(\delta g_{n-1},\mu_{n-1})}\circ\cdots\circ 
\tilde{\mathcal{L}}^{(\delta g_{j+1},\mu_{j+1})}\left(\tilde{\xi}_R(\delta g_j,\mu_j,R_j)\right)\ .
\]
However,
\[
|||\tilde{\mathcal{L}}^{(\delta g_{n-1},\mu_{n-1})}\circ\cdots\circ 
\tilde{\mathcal{L}}^{(\delta g_{n-k},\mu_{n-k})}(R_{n-k})|||_{\bar{g}}\le
\left(\frac{1}{2}\right)^{k}|||R_{n-k}|||_{\bar{g}}\le 2^{-k}\times\frac{\rho'}{3}\bar{g}^{e_R}
\]
and thus this boundary term disappears when $k\rightarrow\infty$ and we then get
\[
R_n=
\sum_{j\le n-1} \tilde{\mathcal{L}}^{(\delta g_{n-1},\mu_{n-1})}\circ\cdots\circ 
\tilde{\mathcal{L}}^{(\delta g_{j+1},\mu_{j+1})}\left(\tilde{\xi}_R(\delta g_j,\mu_j,R_j)\right)\ .
\]

As for the $\mu's$, we have for all $n\le -1$
\[
\mu_{n+1}=L^{\frac{3+\epsilon}{2}}\mu_n+\tilde{\xi}_\mu(\delta g_n,\mu_n,R_n)
\]
or equivalently
\[
\mu_n=-L^{-\left(\frac{3+\epsilon}{2}\right)}\mu_{n+1}+
L^{-\left(\frac{3+\epsilon}{2}\right)}\tilde{\xi}_2(\delta g_n,\mu_n,R_n)\ .
\]
Provided $n+1\le -1$, we apply this to $n+1$ instead of $n$ and substitute in the first of the previous equation.
Iterating this procedure forward until one hits the boundary term $\mu_0$ gives
\[
\mu_n=L^{n\left(\frac{3+\epsilon}{2}\right)}\mu_0
-\sum_{j=n}^{-1} L^{-(j-n+1)\left(\frac{3+\epsilon}{2}\right)}\tilde{\xi}_2(\delta g_j,\mu_j,R_j)\ .
\]
We therefore proved $\vec{u}=\mathfrak{n}(\vec{u})$. By the uniqueness part of the Banach Fixed Point Theorem,
$\vec{u}$ and $\vec{u}_\ast$ are equal and therefore so are their $\delta g_0$ and $R_0$ components.
This establishes $(\delta g,R)=(\delta g_{\rm u}(\mu),R_{\rm u}(\mu))$ and finally $(\delta g,\mu,R)
\in W^{\rm u,loc}$ as wanted.
\qed

\begin{Lemma}\label{L63lem}
Provided $\rho$ and $\rho'$ are chosen so that $\rho<\frac{3}{8}\rho'$,
we have $v_\ast\in W^{\rm u,loc}$ as well as the equations
\[
\mu_\ast=\mu_{\rm s}(\delta g_\ast,R_\ast)\ ,\ \delta g_\ast=\delta g_{\rm u}(\mu_\ast)
\ ,\ R_\ast=R_{\rm u}(\mu_\ast)\ .
\]
\end{Lemma}

\noindent{\bf Proof:}
From $v_\ast\in W^{\rm s,loc}$ we get
\begin{equation}
|\delta g_\ast|\le \frac{\rho}{13}\bar{g}^{e_4}\ ,\ 
|\mu_\ast|\le \frac{\rho}{3}\bar{g}^{e_2}\ ,\ 
|||R_\ast|||_{\bar{g}}\le \frac{\rho}{13}\bar{g}^{e_R}\ . 
\label{fromunsteq}
\end{equation}
We also know that $RG(v_\ast)=v_\ast$.
Define $(\delta g_n,\mu_n,R_n)=v_\ast$ for all $n\le 0$.
Since this is an $RG$ trajectory, all we need in order to prove
$v_\ast=(\delta_0,\mu_0,R_0)\in W^{\rm u,loc}$ via Lemma \ref{L62prop} are the inequalities
\[
|\delta g_\ast|\le \frac{\rho'}{3}\bar{g}^{e_4}\ ,\ 
|\mu_\ast|< \frac{\rho'}{8}\bar{g}^{e_2}\ ,\ 
|||R_\ast|||_{\bar{g}}\le \frac{\rho'}{3}\bar{g}^{e_R}\ . 
\]
The latter easily follow from (\ref{fromunsteq}) and the hypothesis $\rho<\frac{3}{8}\rho'$.
Finally the three equations satisfied by $v_\ast$ are tautological.
\qed

\begin{Lemma}\label{L64lem}
If $v\neq v'$ belong to $W^{\rm u,loc}$ then
$||v_2-v'_2||>||v_1-v'_1||$.
\end{Lemma}

\noindent{\bf Proof:}
By Proposition \ref{L62prop} there exists sequences $(w_n)_{n\le 0}$ and $(w'_n)_{n\le 0}$ in $\mathcal{E}$ such that
$w_0=v$, $w'_0=v'$, $w_{n+1}=RG(w_{n})$, $w'_{n+1}=RG(w'_{n})$  for all $n\le -1$ and 
such that the inequalities in Proposition \ref{L62prop} hold for both sequences.
The latter imply that the corresponding points all are in the domain of application of Lemma \ref{L57lem} because
$\frac{\rho'}{3}<\rho'\le\frac{1}{8}$ by assumption.
For all $n\le 0$, $v\neq v'$ can be rewritten $RG^{(-n)}(w_n)\neq RG^{(-n)}(w'_n)$ and thus $w_n\neq w'_n$.
We proceed by contradiction and suppose $||v_1-v'_1||\ge||v_2-v'_2||$.
This provides the $n=0$ instance of the property $\forall n\le 0,\ ||w_{n,1}-w'_{n,1}||\ge||w_{n,2}-w'_{n,2}||$
which we prove by descending induction.
Suppose the inequality is true for $n$.
Then by Lemma \ref{L57lem} 1)
\begin{equation}
||w_{n,1}-w'_{n,1}||=||RG_1(w_{n-1})-RG_1(w'_{n-1})||\le c_1(\epsilon) ||w_{n-1}-w'_{n-1}||\ .
\label{dichw1eq}
\end{equation}
We now examine two possible cases.

\noindent{\bf 1st Case:} Suppose $L^{\frac{3}{4}}||w_{n-1,2}-w'_{n-1,2}||\ge||w_{n-1,1}-w'_{n-1,1}||$.
By Part 2) of Lemma \ref{L57lem}
\begin{equation}
||w_{n,2}-w'_{n,2}||=||RG_2(w_{n-1})-RG_2(w'_{n-1})||\ge c_2(\epsilon) ||w_{n-1}-w'_{n-1}||\ .
\label{dichw2eq}
\end{equation}
Combining (\ref{dichw1eq}) and (\ref{dichw2eq}) we obtain
\begin{equation}
||w_{n,1}-w'_{n,1}||\le c_1(\epsilon) c_2(\epsilon)^{-1} ||w_{n,2}-w'_{n,2}||
\label{dichw3eq}
\end{equation}
However $c_1(\epsilon) c_2(\epsilon)^{-1}<1$ makes (\ref{dichw3eq}) incompatible with the induction hypothesis
unless $||w_{n,2}-w'_{n,2}||=0$. The latter implies, via (\ref{dichw3eq}), that $||w_{n,1}-w'_{n,1}||=0$
and therefore $w_n=w'_n$ which has been shown to be impossible.
In fact, this 1st Case does not occur.

\noindent{\bf 2nd Case:} Suppose $L^{\frac{3}{4}}||w_{n-1,2}-w'_{n-1,2}||<||w_{n-1,1}-w'_{n-1,1}||$.
Since $L^{\frac{3}{4}}>1$, this immediately implies the induction hypothesis for $n-1$.

From the inequalities we just proved by induction and the definition of the norms we have, for all $n\le 0$,
\[
||w_{n}-w'_{n}||=||w_{n,1}-w'_{n,1}|| \ .
\]
Thus (\ref{dichw1eq}) becomes
\[
||w_{n-1,1}-w'_{n-1,1}||\ge c_1(\epsilon)^{-1} ||w_{n,1}-w'_{n,1}||
\]
for all $n\le 0$. Trivial iteration gives
\[
||w_{n,1}-w'_{n,1}||\ge \left[c_1(\epsilon)^{-1}\right]^{-n} ||w_{0,1}-w'_{0,1}||
\]
which contradicts the boundedness of the $w$ and $w'$ sequences in the $n\rightarrow -\infty$ limit because $c_1(\epsilon)<1$,
unless $||w_{0,1}-w'_{0,1}||=0$. Hence $||v_1-v'_1||=0$ which also implies  $||v_2-v'_2||=0$ by the assumption made at the beginning.
We then arrive at $v=v'$ which is impossible.
\qed

\begin{Corollary}\label{stabunstabcor}
Under the hypotheses of Lemma \ref{L63lem} $W^{\rm s,loc}\cap W^{\rm u,loc}=\{v_\ast\}$.
\end{Corollary}

\noindent{\bf Proof:}
We already know that the fixed point $v_\ast$ belongs to the intersection.
Suppose $v\neq v_\ast$ does too.
Since $v$ and $v_\ast$ are distinct and belong to $W^{\rm s,loc}$, then Lemma \ref{L58lem} and the fact $L^{\frac{3}{4}}\ge 1$
imply
\[
||v_1-v_{\ast,1}||>||v_2-v_{\ast,2}||\times L^{\frac{3}{4}}\ge ||v_2-v_{\ast,2}||
\]
However since $v$ and $v_\ast$ are distinct and belong to $W^{\rm u,loc}$, Lemma \ref{L64lem} implies
\[
||v_2-v_{\ast,2}||>||v_1-v_{\ast,1}||
\]
which gives a contradiction.
\qed

We conclude this section by giving an analogue of Proposition \ref{L56prop} for the unstable manifold.
Since this corresponds to an expanding direction for the bulk RG map, one cannot hope for the stability of $W^{\rm u,loc}$.
However we will consider a smaller set $W_{\rm tiny}^{\rm u,loc}$ and show, under suitable additional hypotheses, that
$RG(W_{\rm tiny}^{\rm u,loc})\subset W^{\rm u,loc}$.
For $\rho''>0$ to be suitably adjusted we first define
\[
W_{\rm small}^{\rm u,loc}=\left\{
(\delta g_{\rm u}(\mu),\mu,R_{\rm u}(\mu))|\ |\mu-\mu_\ast|<\rho'' \bar{g}^{e_2}
\right\}
\]
According to the prevailing hypotheses, as in the statement of Lemma \ref{L63lem}, we have $\frac{\rho}{3}<\frac{\rho'}{8}$.
This leaves the possibility of adding the new constraint $\rho''<\frac{\rho'}{8}-\frac{\rho}{3}$ on the new parameter $\rho''$.
From the proof of Lemma \ref{L63lem} we get $|\mu_{\ast}|\le \frac{\rho}{3}\bar{g}^{e_2}$ and therefore
$|\mu-\mu_\ast|<\rho'' \bar{g}^{e_2}$ implies
\[
|\mu|\le |\mu-\mu_\ast|+|\mu_\ast|<\rho''\bar{g}^{e_2}+\frac{\rho}{3}\bar{g}^{e_2}<\frac{\rho'}{8}\bar{g}^{e_2}\ .
\]
This garantees $W_{\rm small}^{\rm u,loc}\subset W^{\rm u,loc}$.
Also note that $\rho'\le \frac{1}{8}$ implies that $W^{\rm u,loc}$
and therefore $W_{\rm small}^{\rm u,loc}$ are contained in a domain where $RG$ is well-defined and analytic.
Therefore, the set $RG(W)$ is also well-defined for any subset $W$ of $W_{\rm small}^{\rm u,loc}$.
We now define
\[
W_{\rm tiny}^{\rm u,loc}=\left\{
(\delta g_{\rm u}(\mu),\mu,R_{\rm u}(\mu))|\ |\mu-\mu_\ast|<\rho''' \bar{g}^{e_2}
\right\}
\]
where
\[
\rho'''=\min\left\{
\rho'',\frac{3}{16} c_3(\epsilon)^{-1}\left(\frac{\rho'}{8}-\frac{\rho}{3}\right)\right\}
\]
where $c_3(\epsilon)$ has been defined in Lemma \ref{L57lem}.
Our working assumptions are now that
\[
0<\rho<\frac{\rho'}{8}\ \ {\rm and}\ \ 0<\rho''<\frac{1}{3}\left(
\frac{\rho'}{8}-\rho
\right)
\]
which are stronger
than the previous ones and garantee $0<\rho'''\le \rho''$ so $W_{\rm tiny}^{\rm u,loc}$ is indeed a subset
of $W_{\rm small}^{\rm u,loc}$.

\begin{Proposition}\label{L68prop}
In the small $\epsilon$ regime,
$W_{\rm tiny}^{\rm u,loc}$ satisfies
\[
RG(W_{\rm tiny}^{\rm u,loc})\subset W^{\rm u,loc}\ .
\]
\end{Proposition}

\noindent{\bf Proof:}
Let $(\delta g,\mu,R)\in W_{\rm small}^{\rm u,loc}$ and consider the associated backwards trajectory $(\delta g_n,\mu_n,R_n)_{n\le 0}$
produced by Proposition \ref{L62prop}. Let $(\delta g_1,\mu_1,R_1)=RG(\delta g,\mu,R)$. We will show that the extended
sequence $(\delta g_n,\mu_n,R_n)_{n\le 1}$ satisfies the conditions stated in Proposition \ref{L62prop} (with suitable and obvious shift
in indexation). For $n\le -1$, the bounds we need are the ones we already have.
For $n=0$ the bounds we have are stronger than the ones we need. Indeed, $|\mu_0|<\frac{\rho'}{8}\bar{g}^{e_2}$
trivially implies $|\mu_0|\le\frac{\rho'}{3}\bar{g}^{e_2}$.
We now focus on the $n=1$ case.
By Lemma \ref{L57lem}
\[
||RG_1(\delta g_0,\mu_0,R_0)-RG_1(v_\ast)||\le c_1(\epsilon) ||(\delta g_0,\mu_0,R_0)-v_{\ast}||\ ,
\]
namely,
\[
||(\delta g_1-\delta g_{\ast},R_1-R_\ast)||\le  c_1(\epsilon)\
||(\delta g_{\rm u}(\mu)-\delta g_{\rm u}(\mu_{\ast}),\mu-\mu_{\ast}, R_{\rm u}(\mu)-R_{\rm u}(\mu_{\ast}))||\ .
\]
Note that, by construction in Proposition \ref{unstabmanprop},
the analytic function $\delta g_{\rm u}$ satisfies the hypotheses of Lemma \ref{Lipschitzlem}
with $r_1=\frac{\rho'}{8}\bar{g}^{e_2}$ and $r_2=\frac{\rho'}{3}\bar{g}^{e_4}$.
If we choose $\nu=\frac{1}{3}$ then resulting Lipschitz estimate will give us
\[
|\delta g_{\rm u}(\mu)-\delta g_{\rm u}(\mu_{\ast})|\le \frac{16}{3}\bar{g}^{e_4-e_2}|\mu-\mu_\ast|
\]
provided both $\mu$ and $\mu_\ast$ are bounded by $\frac{\rho'}{24}\bar{g}^{e_2}$.
However these two requirements are garanteed by the hypotheses $\rho<\frac{\rho'}{8}$
and $\rho'''\le\frac{1}{3}\left(\frac{\rho'}{8}-\rho\right)$ together with
\[
|\mu_\ast|\le\frac{\rho}{3} \bar{g}^{e_2}\ \ {\rm and}\ \ 
|\mu|\le |\mu_{\ast}|+|\mu-\mu_{\ast}|< \frac{\rho}{3}\bar{g}^{e_2}+\rho'''\bar{g}^{e_2}\ .
\]
We therefore have
\[
\bar{g}^{-e_4} |\delta g_{\rm u}(\mu)-\delta g_{\rm u}(\mu_{\ast})|\le \frac{16}{3}\bar{g}^{-e_2}|\mu-\mu_\ast|\ .
\]
By the same reasoning and use of Lemma \ref{Lipschitzlem} for the function $R_{\rm u}$ we also have
\[
\bar{g}^{-e_R} |R_{\rm u}(\mu)-R_{\rm u}(\mu_{\ast}|\le \frac{16}{3}\bar{g}^{-e_2}|\mu-\mu_\ast|\ .
\]
As a result
\[
||(\delta g_{\rm u}(\mu)-\delta g_{\rm u}(\mu_{\ast}),\mu-\mu_{\ast}, R_{\rm u}(\mu)-R_{\rm u}(\mu_{\ast}))||
\le \frac{16}{3}\bar{g}^{-e_2}|\mu-\mu_\ast|<\frac{16}{3}\rho'''
\]
and thus
\[
||(\delta g_1-\delta g_{\ast},R_1-R_\ast)||<  \frac{16}{3} c_1(\epsilon)\rho'''<\frac{16}{3}\rho'''\ .
\]
In view of $|\delta g_\ast|\le\frac{\rho}{13}\bar{g}^{e_4}$ and
$|||R_\ast|||_{\bar{g}}\le\frac{\rho}{13}\bar{g}^{e_R}$
provided by the fact $v_\ast\in W^{\rm s,loc}$ and by simple triangle inequalities we obtain
\begin{eqnarray*}
|\delta g_1| & < & \left(\frac{16}{3}\rho'''+\frac{\rho}{13}\right)\bar{g}^{e_4}\ , \\
|||R_1|||_{\bar{g}} & < & \left(\frac{16}{3}\rho'''+\frac{\rho}{13}\right)\bar{g}^{e_R}\ .
\end{eqnarray*}
Since $c_3(\epsilon)\ge L^{\frac{3}{4}}>1$, the definition of $\rho'''$ implies
\[
\frac{16}{3}\rho'''+\frac{\rho}{13}< \frac{16}{3}c_3(\epsilon)\rho'''+\frac{\rho}{3}\le \frac{\rho'}{8}<\frac{\rho'}{3}\ .
\]
Hence, $|\delta g_1|< \frac{\rho'}{3}\bar{g}^{e_4}$
and $|||R_1|||_{\bar{g}}<\frac{\rho'}{3}\bar{g}^{e_R}$.

By Lemma \ref{L57lem}
\[
||RG_2(\delta g_0,\mu_0,R_0)-RG_2(v_\ast)||\le c_3(\epsilon) ||(\delta g_0,\mu_0,R_0)-v_{\ast}||\ .
\]
By the same bound on the right-hand side as before we thus get
\[
||RG_2(\delta g_0,\mu_0,R_0)-RG_2(v_\ast)||\le c_3(\epsilon)\times \frac{16}{3}\bar{g}^{-e_2}|\mu-\mu_\ast|<
\frac{16}{3} c_3(\epsilon)\rho'''\ .
\]
Hence
\[
|\mu_1-\mu_\ast|<\frac{16}{3} c_3(\epsilon)\rho'''\bar{g}^{e_2}
\]
which together with $|\mu_{\ast}|\le\frac{\rho}{3}\bar{g}^{e_2}$
and the hypothesis on $\rho'''$ implies $|\mu_1|<\frac{\rho'}{8}$.

We therefore proved the required shifted bounds on the sequence $(\delta g_n,\mu_n,R_n)_{n\le 1}$
in order to conclude by reverse use of Proposition \ref{L62prop} that $(\delta g_1,\mu_1,R_1)\in W^{\rm u,loc}$.
\qed

\subsection{Study of the differential of the RG map at the fixed point and quantitative transversality}

We now study the differential $D_{v_\ast}RG$ of the map $RG$ at the fixed point $v_{\ast}$ in relation
to the invariant linear subspaces $\mathcal{E}^{\rm s}$ and $\mathcal{E}^{\rm u}$ corresponding to the tangent spaces
to the stable and unstable manifolds at the fixed point respectively.
We first define $\mathcal{E}^{\rm s}$ as the kernel of the $\mathbb{C}$-linear form
\[
(\delta g,\mu,R)\mapsto \mu-D_{v_{\ast,1}}\mu_{\rm s}[\delta g,R]
\]
where $D_{v_{\ast,1}}\mu_{\rm s}$ is the differential of $\mu_s$ at $v_{\ast,1}=(\delta g_{\ast},R_{\ast})$.
This linear form is clearly nonzero. It is also continuous by analyticity of $\mu_s$.
Therefore $\mathcal{E}^{\rm s}$ is a closed complex hyperplane in $\mathcal{E}$.

We likewise define $\mathcal{E}^{\rm u}$ as the kernel of the $\mathbb{C}$-linear map
\[
\left\{
\begin{array}{ccl}
\mathcal{E} & \longrightarrow & \mathcal{E}_1 \\
(\delta g,\mu,R) & \longmapsto & (\delta g-D_{v_{\ast,2}}\delta g_{\rm u}[\mu], R-D_{v_{\ast,2}}R_{\rm u}[\mu])
\end{array}
\right.
\]
in terms of the differentials at $v_{\ast,2}=\mu_{\ast}$ of the analytic maps $\delta g_{\rm u}$ and $R_{\rm u}$.
Again, $\mathcal{E}^{\rm u}$ is a closed subspace of $\mathcal{E}$.
In fact, it is easy to see that $\mathcal{E}^{\rm u}$ is equal to the complex line $\mathbb{C}e_{\rm u}$
with
\[
e_{\rm u}=(D_{v_{\ast,2}}\delta g_{\rm u}[1],1,D_{v_{\ast,2}}R_{\rm u}[1])\ .
\]

\begin{Lemma}\label{L65lem}
For all $v\in \mathcal{E}^{\rm u}$ we have $||v_1||\le ||v_2||$.
For all $v\in \mathcal{E}^{\rm s}$ we have $L^{\frac{3}{4}}||v_2||\le ||v_1||$.
As a consequence we have the direct sum decomposition $\mathcal{E}=\mathcal{E}^{\rm s}\oplus\mathcal{E}^{\rm u}$.
\end{Lemma}

\noindent{\bf Proof:}
Define the complex curve parametrized by 
$\gamma(\mu)=(\delta g_{\rm u}(\mu_\ast+\mu),\mu_\ast+\mu,R_{\rm u}(\mu_\ast+\mu))$ for $\mu\in\mathbb{C}$ small.
By Proposition \ref{L62prop}, $\gamma(\mu)\in W^{\rm u,loc}$ for $\mu$ small. Since we also
have $v_{\ast}\in W^{\rm u,loc}$, then Lemma \ref{L64lem} gives us
\[
||\gamma(\mu)_{1}-v_{\ast,1}||\le||\gamma(\mu)_{2}-v_{\ast,2}|| 
\]
for $\mu$ small.
However by analyticity and therefore differentiability we have
$\gamma(\mu)=v_{\ast}+\mu e_{\rm u}+\mu\omega(\mu)$ where $\omega(\mu)\rightarrow 0$ when $\mu\rightarrow 0$.
The previous inequality becomes
\[
||\mu e_{{\rm u},1}+\mu \omega(\mu)_{1}||\le
||\mu e_{{\rm u},2}+\mu \omega(\mu)_{2}||\ .
\]
For $\mu\neq 0$ we divide by $|\mu|$ and then let $\mu$ go to $0$ which gives $||e_{{\rm u},1}||\le||e_{{\rm u},2}||$
and therefore $||v_1||\le ||v_2||$ for all $v\in \mathcal{E}^{\rm u}=\mathbb{C} e_{\rm u}$.

Now for $v=(\delta g,\mu,R)\in \mathcal{E}^{\rm s}$ we this time let, for $t\in\mathbb{C}$ small,
\[
\gamma(t)=(\delta g_{\ast}+t\delta g,\mu_{\rm s}(\delta g_{\ast}+t\delta g,R_{\ast}+tR),R_{\ast}+tR)\ .
\]
Since $v_\ast\in W_{\rm int}^{\rm s,loc}$ and $\mu_s$ is analytic and therefore continuous we have $\gamma(t)\in
W_{\rm int}^{\rm s,loc}\subset W^{\rm s,loc}$ for $t$ small. Lemma \ref{L58lem} thus gives the inequality
\[
L^{\frac{3}{4}}||\gamma(t)_{2}-v_{\ast,2}||\le||\gamma(t)_{1}-v_{\ast,1}|| \ .
\]
Again one can write
\[
\gamma(t)=v_{\ast}+t(\delta g, D_{v_{\ast},1}\mu_{\rm s}[v_1],R)+t\omega(t)
\]
where the new function $\omega$ satisfies
$\omega(t)\rightarrow 0$ when $t\rightarrow 0$.
The previous inequality becomes
\[
L^{\frac{3}{4}}||t D_{v_{\ast},1}\mu_{\rm s}[v_1] +t \omega(t)_{2}||\le
||t v_1 +t \omega(t)_{1}||\ .
\]
Again dividing by $|t|$ for $t\neq 0$ and then letting $t\rightarrow 0$ we get
\[
L^{\frac{3}{4}}||D_{v_{\ast},1}\mu_{\rm s}[v_1]||\le ||v_1||\ ,
\]
i.e., $L^{\frac{3}{4}}||v_2||\le ||v_1||$ by the defining equation 
$v_2=\mu=D_{v_{\ast},1}\mu_{\rm s}[v_1]$
of $\mathcal{E}^{\rm s}$.

A vector $v$ which satisfies both inequalities  $L^{\frac{3}{4}}||v_2||\le ||v_1||$
and $||v_1||\le ||v_2||$ must clearly satisfy $||v_1||=||v_2||=0$ because $L^{\frac{3}{4}}>1$.
Namely, $v$ must vanish. This shows $\mathcal{E}^{\rm s}\cap \mathcal{E}^{\rm u}=\{0\}$.
Since $e_{\rm u}\neq 0$ is in $\mathcal{E}^{\rm u}$ we get $e_{\rm u}\notin\mathcal{E}^{\rm s}$.
This proves the direct sum property since we are considering a complex line spanned by $e_{\rm u}$ and a 
complex hyperplane.
\qed

\begin{Lemma}\label{L66lem}
The subspace $\mathcal{E}^{\rm s}$ is invariant by $D_{v_\ast}RG$.
\end{Lemma}

\noindent{\bf Proof:}
For $v=(\delta g,\mu,R)\in \mathcal{E}^{\rm s}$ we again use the curve
\[
\gamma(t)=(\delta g_{\ast}+t\delta g,\mu_{\rm s}(\delta g_{\ast}+t\delta g,R_{\ast}+tR),R_{\ast}+tR)\ .
\]
from the proof of the previous lemma and which satisfies
$\gamma'(0)=v$.
For $t$ small $\gamma(t)$ is well-defined, takes values in $W_{\rm int}^{\rm s,loc}\subset W^{\rm s,loc}$,
is analytic in $t$ and belongs to an open set where $RG$ is well-defined and analytic.
Thus, $t\mapsto RG(\gamma(t))$ is analytic near $t=0$.
By Proposition \ref{L56prop}, $RG(\gamma(t))\in W_{\rm int}^{\rm s,loc}$ and therefore
$RG(\gamma(t))=\mu_{\rm s}(RG_1(\gamma(t)))$. We differentiate this at $t=0$ using the chain rule and obtain
\[
\left(D_{v_\ast}RG[\gamma'(0)]\right)_2=D_{v_{\ast},1}\mu_{\rm s}\left[
\left(D_{v_\ast}RG[\gamma'(0)]\right)_1
\right]\ ,
\]
i.e.,
\[
\left(D_{v_\ast}RG[v]\right)_2=D_{v_{\ast},1}\mu_{\rm s}\left[
\left(D_{v_\ast}RG[v]\right)_1
\right]\ .
\]
Hence $D_{v_\ast}RG[v]\in\mathcal{E}^{\rm s}$ by definition of $\mathcal{E}^{\rm s}$.
\qed

\begin{Lemma}\label{L69lem}
The subspace $\mathcal{E}^{\rm u}$ is invariant by $D_{v_\ast}RG$.
\end{Lemma}

\noindent{\bf Proof:}
We reuse the curve
\[
\gamma(\mu)=(\delta g_{\rm u}(\mu_\ast+\mu),\mu_\ast+\mu,R_{\rm u}(\mu_\ast+\mu))
\]
from the proof of Lemma \ref{L65lem}. Clearly, $\gamma(\mu)$ lies in
$W_{\rm tiny}^{\rm u,loc}$ when $\mu$ is small.
Therefore $RG(\gamma(\mu))$ is analytic and lies in $W^{\rm u,loc}$ when $\mu$ is small because of
Proposition \ref{L68prop}. By definition of $W^{\rm u,loc}$ we thus have
\[
RG_1(\gamma(\mu))=\left(\delta g_{\rm u}(RG_2(\gamma(\mu))), R_{\rm u}(RG_2(\gamma(\mu)))\right)\ .
\]
We differentiate at $\mu=0$ using $\gamma'(0)=e_{\rm u}$ and the chain rule.
This gives
\[
\left(D_{v_\ast}RG[e_{\rm u}]\right)_1=
\left(
D_{v_{\ast,2}}\delta g_{\rm u}\left[\left(D_{v_\ast}RG[e_{\rm u}]\right)_2\right],
D_{v_{\ast,2}}R_{\rm u}\left[\left(D_{v_\ast}RG[e_{\rm u}]\right)_2\right]
\right)\ .
\]
In other words $D_{v_\ast}RG[e_{\rm u}]$ satisfies the defining equation of $\mathcal{E}^{\rm u}=\mathbb{C} e_{\rm u}$
which therefore is invariant by the differential of $RG$ at the fixed point $v_\ast$.
\qed

\begin{Lemma}\label{alphabdlem}
The restriction $\left.D_{v_\ast}RG\right|_{\mathcal{E}^{\rm u}}$ is the multiplication by an eigenvalue
$\alpha_{\rm u}$ which is real and greater than $1$. One also has the more precise estimate
\[
|\alpha_{\rm u}-L^{\frac{3+\epsilon}{2}}|\le c_4(\epsilon)
\]
where $c_4(\epsilon)$ has been defined in Lemma \ref{L57lem}.
\end{Lemma}

\noindent{\bf Proof:}
By Lemma \ref{L69lem} and the unidimensional property of $\mathcal{E}^{\rm u}$ we have
$D_{v_\ast}RG[e_{\rm u}]=\alpha_{\rm u} e_{\rm u}$ for some possibly complex $\alpha_{\rm u}$.
However, by Theorem \ref{mainestthm}, the map $RG$ sends real data to real data. Thus, so does the map
on sequences $\mathfrak{m}$ used in Proposition \ref{stabmanprop}.
Therefore the corresponding fixed point $\vec{u}_\ast$ in the space of sequences obtained by iteration starting from
the null sequence $\vec{0}$ which is real is also real provided the implicit data $(\delta g_0,R_0)$ is too.
As a result the map $\mu_{\rm s}$ sends real data to real data. In other words, if $\delta g\in\mathbb{R}$ and
if $R$ is a real-valued even function then $\mu_{\rm s}(\delta g, R)\in\mathbb{R}$.
Similar statements also hold for the functions $\delta g_{\rm u}$ and $R_{\rm u}$ used for the parametrization of the local
unstable manifold $W^{\rm u,loc}$. It is also easy to see that the fixed point $v_{\ast}$ is real.
Finally, the eigenvalue $\alpha_{\rm u}$ which coincides with the second or $\mu$-component
of $D_{v_\ast}RG[e_{\rm u}]$ is easily seen to be a real number.

Now again consider the curve $\gamma(\mu)$ as in the proof of Lemma \ref{L69lem}.
For $\mu$ small we have by Lemma \ref{L57lem}
\[
||RG_2(\gamma(\mu))-RG_2(\gamma(0))-L^{\frac{3+\epsilon}{2}}(\gamma(\mu)-\gamma(0))||\le c_4(\epsilon)
||\gamma(\mu)-\gamma(0)||\ .
\]
We divide by $|\mu|\neq 0$ and then take the limit when $\mu$ goes to $0$. This results in
\[
||\left(D_{v_\ast}RG[e_{\rm u}]\right)_2-L^{\frac{3+\epsilon}{2}} e_{\rm u}||\le c_4(\epsilon)||e_{\rm u}||\ ,
\]
i.e.,
\[
|\alpha_{\rm u}-L^{\frac{3+\epsilon}{2}}|\times ||e_{{\rm u},2}||\le c_4(\epsilon)||e_{\rm u}||\ .
\]
Since $e_{\rm u}$ belongs to $\mathcal{E}^{\rm u}$ we have by Lemma \ref{L65lem} the equality $||e_{{\rm u},2}||=||e_{\rm u}||$.
Since $e_{\rm u}$ is nonzero we can simplify by $||e_{\rm u}||$ and we end up with the desired estimate.
Finally, in the small $\epsilon$ regime, $c_4(\epsilon)$ goes to zero which readily implies $\alpha_{\rm u}>1$.
\qed

\begin{Lemma}\label{L72lem}
The restriction $\left.D_{v_\ast}RG\right|_{\mathcal{E}^{\rm s}}$ is a contraction on the subspace
$\mathcal{E}^{\rm s}$. More precisely, for every $v\in\mathcal{E}^{\rm s}$, we have $D_{v_\ast}RG[v]\in \mathcal{E}^{\rm s}$ and
\[
||D_{v_\ast}RG[v]||\le c_1(\epsilon)||v||
\]
where $c_1(\epsilon)\in (0,1)$ has been defined in Lemma \ref{L57lem}.
\end{Lemma}

\noindent{\bf Proof:}
For $v=(\delta g,\mu,R)\in \mathcal{E}^{\rm s}$ we again use the curve
$\gamma(t)$ as in the proof of 
Lemma \ref{L66lem}. For small $t$ we have $\gamma(t)\in W^{\rm s,loc}$. We can thus derive from
Lemma \ref{L59lem} the inequality
\[
||RG(\gamma(t))-RG(\gamma(0))||\le c_1(\epsilon)||\gamma(t)-\gamma(0)||\ .
\]
We divide by $|t|\neq 0$ and take the $t\rightarrow 0$ limit in order to obtain
\[
||D_{v_\ast}RG[\gamma'(0)]||\le  c_1(\epsilon)||\gamma'(0)||\ ,
\]
i.e.,
\[
||D_{v_\ast}RG[v]||\le  c_1(\epsilon)||v||
\]
since, as one can easily see, $\gamma'(0)=v$.
Stability has already been shown in Lemma \ref{L66lem}.
\qed

\subsection{Explicit equivalence of norms}

For the needs of \S\ref{Koenigssec} we introduce another norm $||\cdot||_{\substack{\ \\ \Diamond}}$ on $\mathcal{E}$.
Recall that the latter is the direct sum $\mathcal{E}_1\oplus\mathcal{E}_2$ and the original
norm $||\cdot||$ behaves well with respect to this decomposition.
Indeed, if $v=v_1+v_2$ is the decomposition of a vector according to this direct sum, we
have
\[
||v||=\max(||v_1||,||v_2||)\ .
\]
The $||\cdot||_{\substack{\ \\ \Diamond}}$ is designed in order to satisfy a similar property with respect to
the direct sum $\mathcal{E}=\mathcal{E}^{\rm u}\oplus\mathcal{E}^{\rm s}$. Using the notations
$P_{\rm u}$ and $P_{\rm s}$ for the corresponding projections on the two subspaces $\mathcal{E}^{\rm u}$
and $\mathcal{E}^{\rm s}$ respectively,
we let by definition
\[
||v||_{\substack{\ \\ \Diamond}}=\max(||P_{\rm u}(v)||,||P_{\rm s}(v)||)\ .
\]

\begin{Lemma}\label{eqnormlem}
We have the explicit equivalence of norms
\[
\frac{1}{2}||v||\le ||v||_{\substack{\ \\ \Diamond}}\le 5||v||
\]
for all $v\in\mathcal{E}$.
\end{Lemma}

\noindent{\bf Proof:}
For such a $v$ let us write for simplicity $v^{\rm u}=P_{\rm u}(v)$
and $v^{\rm s}=P_{\rm s}(v)$. We decompose all three vectors $v$, $v^{\rm u}$ and $v^{\rm s}$ according
to the old direct sum $\mathcal{E}_1\oplus\mathcal{E}_2$ as
\begin{eqnarray*}
v & = & v_1+v_2 \\
v^{\rm u} & = & v_1^{\rm u}+v_2^{\rm u} \\
v^{\rm s} & = & v_1^{\rm s}+v_2^{\rm s}
\end{eqnarray*}
noting that we must then have the relations $v_1= v_1^{\rm u}+v_1^{\rm s}$ and $v_2= v_2^{\rm u}+v_2^{\rm s}$.
Armed with this observation and the inequalities in Lemma \ref{L65lem}
one easily checks that
\begin{eqnarray*}
||v^{\rm u}|| & = & \max(||v_1^{\rm u}||,||v_2^{\rm u}||) \\
 & = & ||v_2^{\rm u}|| \\
 & = & ||v_2-v_2^{\rm s}|| \\
 & \le & ||v_2||+||v_2^{\rm s}|| \\
 & \le & ||v_2||+L^{-\frac{3}{4}}||v_1^{\rm s}|| \\
 & = & ||v_2||+L^{-\frac{3}{4}}||v_1-v_1^{\rm u}|| \\
 & \le & ||v_2||+L^{-\frac{3}{4}}||v_1||+L^{-\frac{3}{4}}||v_1^{\rm u}|| \\
 & \le & (1+L^{-\frac{3}{4}}) ||v||+ L^{-\frac{3}{4}}||v^{\rm u}||
\end{eqnarray*}
which results in
\[
||v^{\rm u}||\le \frac{1+L^{-\frac{3}{4}}}{1-L^{-\frac{3}{4}}} ||v||\ .
\]
Similarly we have
\begin{eqnarray*}
||v^{\rm s}|| & = & \max(||v_1^{\rm s}||,||v_2^{\rm s}||) \\
 & = & ||v_1^{\rm s}|| \\
 & = & ||v_1-v_1^{\rm u}|| \\
 & \le & ||v_1||+||v_1^{\rm u}|| \\
 & \le & ||v_1||+||v_2^{\rm u}|| \\
 & = & ||v_1||+||v_2-v_2^{\rm s}|| \\
 & \le & ||v_1||+||v_2||+||v_2^{\rm s}|| \\
 & \le &  ||v_1||+||v_2||+L^{-\frac{3}{4}}||v_1^{\rm s}|| \\
 & \le & 2 ||v||+L^{-\frac{3}{4}} ||v^{\rm s}||
\end{eqnarray*}
which entails
\[
||v^{\rm s}||\le \frac{2}{1-L^{-\frac{3}{4}}} ||v||\ .
\]
Since $L^{-\frac{3}{4}}<1$ we get
\[
||v||_{\substack{\ \\ \Diamond}}=\max(||v^{\rm u}||,||v^{\rm s}||)\le \frac{2}{1-L^{-\frac{3}{4}}} ||v||
\le \frac{2}{1-2^{-\frac{3}{4}}} ||v||\le 5 ||v||
\]
where we used the simplification $\frac{2}{1-2^{-\frac{3}{4}}}\simeq 4.933\ldots<5$.

The other inequality is much simpler. Indeed,
\[
||v_1||=||v_1^{\rm u}+v_1^{\rm s}||\le ||v_1^{\rm u}|| +||v_1^{\rm s}||\le  ||v^{\rm u}||+||v^{\rm s}||\le 
2||v||_{\substack{\ \\ \Diamond}}
\]
and
\[
||v_2||=||v_2^{\rm u}+v_2^{\rm s}||\le ||v_2^{\rm u}|| +||v_2^{\rm s}||\le  ||v^{\rm u}||+||v^{\rm s}||\le 
2||v||_{\substack{\ \\ \Diamond}}\ .
\]
Hence
\[
||v||=\max(||v_1||,||v_2||)\le  2||v||_{\substack{\ \\ \Diamond}}
\]
as desired.
\qed

%% file: DynsysII.tex
\section{Partial analytic linearization}\label{Koenigssec}

The crucial ingredient for the proof of existence of anomalous dimension is
an infinite-dimensional generalization of the K{\oe}nigs Linearization Theorem in one-dimensional
holomorphic dynamics. This is the object of Theorem \ref{unnamedtheorem1} below.
As a preliminary step towards establishing this theorem, we prove some lemmas which give us
some explicit control on the second differential of $RG$.

\begin{Lemma}\label{DRG2lem}
In the small $\epsilon$ regime we have, for all $v$ such that
$\displaystyle ||v||<\frac{1}{4}$,

\[
||D^{2}_{v}RG||
\le 
17\ .
\]

\end{Lemma}

\noindent{\bf Proof:}
In $(\delta g, \mu, R)$ coordinates we have:

\[
RG[g,\mu,R] 
= 
RG^{\rm explicit}(\delta g , \mu, R)
+ 
RG^{\rm implicit}(\delta g, \mu, R)
\]

where

\[
RG^{\rm explicit}(\delta g, \mu, R)
=
\begin{pmatrix}
(2-L^{\epsilon})\delta g - A_{1} \delta g^2\\
L^{\frac{3+\epsilon}{2}} \mu - A_{2} (\bar{g} + \delta g)^2 - A_{3}(
\bar{g} + \delta g)\mu\\
0\\
\end{pmatrix}^{T}
\]

and

\[
RG^{\rm implicit}(\delta g, \mu, R)
=
\begin{pmatrix}
\xi_{4}(g + \delta g, \mu,R)\\
\xi_{2}(\bar{g} + \delta g ,\mu, R)\\
\mathcal{L}^{(\bar{g}+\delta g, \mu)}(R) + \xi_{R}(\bar{g} + \delta g,
\mu, R)\\
\end{pmatrix}^{T}
\]

$T$ meaning transpose.\\

An easy computation shows

\[
D^{2}_{v}RG^{\rm explicit}[ v', v''] =
\begin{pmatrix}
-2A_{1}\delta g' \delta g''\\
-2A_{2} \delta g' \delta g'' - A_{3} \delta g' \mu'' - A_{3} \mu'
\delta g''\\
0\\
\end{pmatrix}^{T}
\]

where $v = (\delta g, \mu, R)$, $v ' = (\delta g', \mu', R')$, $v'' =
(\delta g'', \mu'', R'')$.\\

Here $D^{2}_{v}$, the second differential at $v$, is seen as a
bilinear map acting on pairs of vectors $(v',v'')$.\\

It is immediate from the definition of the norm $|| \cdot ||$ that

\[
|| D^2_{v} RG^{\rm explicit}[v',v'']|| \le 2 ||v'|| \times ||v''|| \times
\max
\left[ A_{1,{\rm max}} \bar{g}^{e_{4}}, A_{2,{\rm max}} \bar{g}^{2e_{4} -
  e_{2}}+A_{3,{\rm max}} \bar{g}^{e_{3}} \right]\ .
\]\\

On the other hand if $\displaystyle ||v|| < \frac{1}{4}$ and $v',v''$ are nonzero,
we can use Cauchy's formula to write

\[
D^{2}_{v} RG^{\rm implicit}[ v',v'']
=
\frac{1}{(2\pi i)^2}
\oint \frac{d\lambda_{1}}{\lambda_{1}^2}
\oint \frac{d\lambda_{2}}{\lambda_{2}^2}
RG^{\rm implicit}(v + \lambda_{1} v' + \lambda_{2} v'')
\]

where the contours of integration are the positively oriented circles
given by $\displaystyle |\lambda_{1}| = \frac{1}{8||v'||}$,
$\displaystyle |\lambda_{2}| =
\frac{1}{8||v''||}$.\\

Since clearly $\displaystyle ||v + \lambda_{1} v' + \lambda_{2} v''|| <\frac{1}{2}$
we are in the domain of analyticity specified by the specialization
of Theorem \ref{mainestthm} in \S\ref{prepsubsec}.
Thus we have

\begin{equation*}
\begin{split}
||RG^{\rm implicit}(v+\lambda_{1}v'+\lambda_{2}v'')||
\le&
\max
\left[
B_{4} \bar{g}^{e_{R}-e_{4}} ||v+\lambda_{1}v'+\lambda_{2}v''||,
B_{2} \bar{g}^{e_{R} - e_{2}} ||v + \lambda_{1}v'+\lambda_{2}v''||, \right.\\
&\ \ \ \ \ \  \left. \frac{1}{2}||v+\lambda_{1}v'+\lambda_{2}v''|| + B_{R\xi}
\bar{g}^{\frac{11}{4}-3\eta-e_{R}} 
\right]\\
\le&
\max
\left[
\frac{1}{2}B_{4}\bar{g}^{e_{R} - e_{4}}, 
\frac{1}{2} B_{2} \bar{g}^{e_{R}-e_{2}},
\frac{1}{4} + B_{R\xi} \bar{g}^{\frac{11}{4} - 3 \eta-e_{R}}
\right]
\end{split}
\end{equation*}

and therefore

\[
||D^{2}_{v} RG^{\rm implicit} [v', v'']||
\le
\left(
8
||v'||
\right)
\left(
8
||v''||
\right)
\times
\max
\left[
\frac{1}{2}B_{4}\bar{g}^{e_{R} - e_{4}}, 
\frac{1}{2} B_{2} \bar{g}^{e_{R}-e_{2}},
\frac{1}{4} + B_{R\xi} \bar{g}^{\frac{11}{4} - 3 \eta-e_{R}}
\right]\ .
\]

In terms of the norm on bilinear forms induced by the vector space
norm  $||\cdot||$ we have:

\begin{equation*}
\begin{split}
||D^{2}_{v}RG|| 
\le&
 2 \max 
\left[A_{1,{\rm max}}\bar{g}^{e_{4}}, 
A_{2,{\rm max}} \bar{g}^{2e_{4}-e_{2}} 
+ A_{3,{\rm max}}\bar{g}^{e_{4}}
\right]\\
&\ \ + 64 \max
\left[ \frac{1}{2} B_{4} \bar{g}^{e_{R}-e_{4}},
\frac{1}{2}B_{2}\bar{g}^{e_{R}-e_{2}},
\frac{1}{4}+B_{R\xi} \bar{g}^{\frac{11}{4}-3\eta-e_{R}} \right]
\le 1 + \frac{64}{4} \le 17\ .
\end{split}
\end{equation*}

In going to the last line we used the assumption on $\epsilon$ being sufficiently
small and the inequalities for exponents indicated in \S\ref{prepsubsec}.
\qed

For $v \in W^{\rm s,loc}$ and $n \ge 0$ we define the continuous linear
map

\[
T_{n}(v)
=
\alpha_{\rm u}^{-n} D_{v} RG^{n}
=
\alpha_{\rm u}^{-n}
D_{RG^{n-1}(v)}RG
\circ
\cdots
\circ
D_{RG(v)}RG
\circ
D_{v}RG\ .
\]

It is well defined by the stability of $W^{\rm s,loc}$ which lies in the domain
of analyticity of $RG$.\\

On the same domain as $RG$ (e.g., the domain $||v|| < \frac{1}{2}$) we
define the map $H$ by 

\[
H(v)
=
RG(v)
-
D_{v_{*}}RG[v]
\]
which is also analytic.

Note that $D_{v_{*}} RG = \alpha_{\rm u}P_{\rm u}+A_{\rm s}P_{\rm s}$ where $A_{\rm s} =
D_{v_{*}} RG \big|_{\mathcal{E}^{\rm s}}$. Thus we have that $RG =
\alpha_{\rm u} P_{\rm u} + A_{\rm s}P_{\rm s} + H$. 

\begin{Lemma}\label{DHDRGlem}
In the small $\epsilon$ regime:
\begin{enumerate}

\item If $\displaystyle ||v|| < \frac{1}{4}$ then $||D_{v}H||\le 17 ||v-v_{*}||$.

\item  If $v$ and $w$ satisfy $\displaystyle ||v|| < \frac{1}{4}$
  and $\displaystyle ||v+w|| <  \frac{1}{4}$ then

\[
||RG(v+w) - RG(v) - D_{v}RG[w]|| \le \frac{17}{2}||w||^2\ .
\]

\end{enumerate}

\end{Lemma}

\noindent{\bf Proof:}
Note that by $||D_{v}H||$ we refer to the norm induced on linear
operators on $\mathcal{E}$ by the norm $||\cdot||$ on vectors in $\mathcal{E}$.
Remark that
$v_{*}\in W^{\rm s,loc}$ implies
$||v_{*}||\le\frac{\rho}{3}<\frac{1}{4}$
since  $\rho < \frac{1}{12}$.

Since the ball of vectors $v$ with $||v|| < \frac{1}{4}$ is convex
we can use the mean value theorem to deduce

\[
||D_{v}H - D_{v_{*}}H|| 
\le ||v-v_{*}|| 
\times
\sup_{0 \le t \le 1} ||D^{2}_{v_{*}+t(v-v_{*})}H||\ .
\]

However by construction $D_{v_{*}}H=0$ and $D^{2}H=D^{2}RG$ so Lemma
\ref{DRG2lem} implies

\[
||D_{v}H||
\le
17
||v-v_{*}||.
\]

By the mean value theorem, or Taylor's formula with integral
remainder, we have 

\[
||RG(v+w) - RG(v) - D_{v}RG[w]||
\le
\frac{1}{2}
||w||^{2}
\times
\sup_{0\le t \le 1} ||D^{2}_{v+tw}RG||
\]

and the desired inequality in Part 2) follows by
Lemma \ref{DRG2lem}.\qed\\

We now give a lemma that shows boundedness for $||T_{n}(v)||_{\substack{\ \\ \Diamond}}$, namely,
the operator norm induced by the norm $||\cdot||_{\substack{\ \\ \Diamond}}$ on vectors.

\begin{Lemma}\label{Tbdlem}
For all $v \in W^{\rm s,loc}$ and all $n \ge 0$ we have

\[
||T_{n}(v)||_{\substack{\ \\ \Diamond}}
\le
\mathcal{C}_{1}(\epsilon),
\]\\

where $\displaystyle \mathcal{C}_{1}(\epsilon) = \exp \left[
  \frac{85}{\alpha_{\rm u}(1-c_{1}(\epsilon))} \right]$.

\end{Lemma}

\noindent{\bf Proof:}
From $\displaystyle
T_{n}(v)
=
\alpha_{\rm u}^{-n} D_{v} RG^{n}
=
\alpha_{\rm u}^{-n}
D_{RG^{n-1}(v)}RG
\circ
\cdots
\circ
D_{RG(v)}RG
\circ
D_{v}RG
$
we immediately get

\[
||T_{n}(v)||_{\substack{\ \\ \Diamond}}
\le
\prod_{k=0}^{n-1}
|| \alpha_{\rm u}^{-1} D_{RG^{k}(v)}RG ||_{\substack{\ \\ \Diamond}}\ .
\]

But

\[
D_{RG^{k}(v)}RG
=
\alpha_{\rm u} P_{\rm u}
+
A_{\rm s} P_{\rm s}
+
D_{RG^{k}(v)}H
\]
results in
\[
|| \alpha_{\rm u}^{-1} D_{RG^{k}(v)}RG ||_{\substack{\ \\ \Diamond}}
\le
||P_{\rm u} + \alpha_{\rm u}^{-1} A_{\rm s} P_{\rm s} ||_{\substack{\ \\ \Diamond}}
+
||\alpha_{\rm u}^{-1} D_{RG^{k}(v)} H ||_{\substack{\ \\ \Diamond}}\ .
\]

From the definition of $||\cdot ||_{\substack{\ \\ \Diamond}}$,
$P_{\rm u}$, $P_{\rm s}$, $\alpha_{\rm u} > 1$ and the contraction $A_{\rm s}$ we get
$||P_{\rm u} + \alpha_{\rm u}^{-1} A_{\rm s} P_{\rm s} ||_{\substack{\ \\ \Diamond}}
\le 1$. Note that this would not work if we had used the
$||\cdot||$ operator norm instead.

Hence

\[
||T_{n}(v)||_{\substack{\ \\ \Diamond}}
\le
\prod_{k=0}^{n-1}
\left(
1 + \alpha_{\rm u}^{-1}||D_{RG^{k}(v)}RG||_{\substack{\ \\ \Diamond}}
\right)\ .
\]

Now from Lemmas \ref{eqnormlem},\ref{DHDRGlem} and \ref{L59lem} we derive

\begin{equation*}
\begin{split}
||D_{RG^{k}(v)}H||_{\substack{\ \\ \Diamond}}
=&
\sup_{w \not = 0}
\frac{||D_{RG^{k}(v)}H[w]||_{\substack{\ \\ \Diamond}}}
{||w||_{\substack{\ \\ \Diamond}}}
\le
\frac{5}{\frac{1}{2}}
\sup_{w \not = 0}
\frac{||D_{RG^{k}(v)}H[w]||}
{||w||}\\
\le &
170
||RG^{k}(v) - v_{*}|| \\
\le & 
170 
c_{1}(\epsilon)^{k} 
||v-v_{*}||\ .
\end{split}
\end{equation*}

Since $\displaystyle ||v - v_{*}|| \le ||v|| + ||v_{*}|| \le \frac{1}{4} +
\frac{1}{4} \le \frac{1}{2}$ we have:

\[
||T_{n}(v)||_{\substack{\ \\ \Diamond}}
\le
\exp 
\left( 
\sum_{k=0}^{n-1} 
\alpha_{\rm u}^{-1} 
\times 
170 
\times
\frac{1}{2}
 c_{1}(\epsilon)^{k}
 \right)
\le
\mathcal{C}_{1}(\epsilon)
\]
as wanted.
\qed\\

We now extract geometric decay in $n$ from $||P_{s}T_{n}(v)||_{\substack{\ \\
    \Diamond}}$.\\

\begin{Lemma}\label{pstnlem}
For all $v \in W^{\rm s,loc}$ and $n \ge 0$

\[
||
P_{\rm s}T_{n}(v)
||_{\substack{\ \\ \Diamond}}
\le 
\mathcal{C}_{2}(\epsilon)
c_{1}(\epsilon)^{\frac{n}{2}}
\]\\

where $\displaystyle \mathcal{C}_{2}(\epsilon) = \mathcal{C}_{1}(\epsilon)
\left[ 1 + \frac{85
    \mathcal{C}_{1}(\epsilon)}{c_{1}(\epsilon)(1-c_{1}(\epsilon))} \right]$

\end{Lemma}

\noindent{\bf Proof:}
We write

\[
P_{\rm s}T_{n}(v)
=
P_{\rm s}(M+N_{n-1})
 \circ \cdots \circ
(M + N_{0})
\]

\[
\textnormal{where  }
M=P_{\rm u} + \alpha_{\rm u}^{-1}A_{\rm s}P_{\rm s}
\]
\[
\textnormal{  and  }
N_{k}=\alpha_{\rm u}^{-1}D_{RG^{k}(v)}H\ .
\]

Let $m$ be such that $ 0 \le m \le n$, then

\[
P_{\rm s}
T_{n}(v)
=
P_{\rm s}
(M+N_{n-1})
\cdots
(M+N_{m})
T_{m}(v)\ .
\]

Then by Lemma \ref{Tbdlem}

\[
||P_{\rm s}T_{n}(v)||_{\substack{\ \\ \Diamond}}
\le
\mathcal{C}_{1}(\epsilon)
|| P_{\rm s}
(M+N_{n-1})
\circ \cdots \circ
(M+N_{m})
||_{\substack{\ \\ \Diamond}}\ .
\]

Now

\begin{equation*}
\begin{split}
|| P_{\rm s}
(M+N_{n-1})
\circ \cdots \circ
(M+N_{m})
||_{\substack{\ \\ \Diamond}}
\le&
|| P_{\rm s}
(M+N_{n-1})
\circ \cdots \circ
(M+N_{m})
-
P_{\rm s}
M^{n-m}
||_{\substack{\ \\ \Diamond}}
+ 
||P_{\rm s} M^{n-m}||_{\substack{\ \\ \Diamond}}\\
\le&
||(M+N_{n-1})
\circ \cdots \circ
(M+N_{m})
-
M^{n-m}||_{\substack{\ \\ \Diamond}} 
+ 
||P_{\rm s} M^{n-m}||_{\substack{\ \\ \Diamond}}\ .
\end{split}
\end{equation*}

In going to the last line we used that $||P_{\rm s}||_{\substack{\ \\ \Diamond}} \le 1$
(in fact one has $||P_{\rm s}||_{\substack{\ \\ \Diamond}} =1$).\\

It is easy to see that

\begin{equation*}
\begin{split}
P_{\rm s} M^{n-m}
=&
P_{\rm s} 
\left( 
P_{\rm u} 
+
\alpha_{\rm u}^{-(n-m)}
A_{\rm s}^{n-m}
P_{\rm s} 
\right)\\
=&
\alpha_{\rm u}^{-(n-m)} 
A_{\rm s}^{n-m}
P_{\rm s}\ .
\end{split}
\end{equation*}

Note that $||w|| = ||w||_{\substack{\ \\ \Diamond}}$ if $w \in
\mathcal{E}^{\rm s}$. Then one easily gets from Lemma \ref{L72lem}

\[
||P_{\rm s}
M^{n-m}
||_{\substack{\ \\ \Diamond}} 
\le
\alpha_{\rm u}^{-(n-m)}c_{1}(\epsilon)^{n-m}\ .
\]

On the other hand, we can write

\begin{equation*}
\begin{split}
(M+N_{n-1})
\circ \cdots \circ
(M+N_{m})
-
 M^{n-m}
=&\ \ \ \ 
N_{n-1}\circ
(M+N_{n-2})
\circ \cdots \circ
(M+N_{m})\\
& +
M\circ
N_{n-2}\circ
(M+N_{n-3})
\circ \cdots \circ
(M+N_{m})\\
&\ \ \vdots\\
& +
M^{n-m-2} 
\circ N_{m+1}
\circ
(M + N_{m})\\
& +
M^{n-m-1} 
\circ
N_{m}\\
=&\ \ \ \  
N_{n-1} 
\circ
T_{n-m-1} (RG^{m}(v))\\
& + 
M 
\circ
N_{n-2} 
\circ
T_{n-m-2}( RG^{m}(v) )\\
&\ \ \vdots\\
& + 
M^{n-m-2}
\circ
N_{m+1}
\circ
T_{1}(RG^{m}(v))\\
& +
M^{n-m-1}
\circ
N_{m}
\circ
T_{0}(RG^{m}(v))\ .
\end{split}
\end{equation*}

In the first equality we are expanding the product and ordering terms
with respect to the leftmost factor of $N_{\bullet}$ that appears.\\

Remembering that $||M||_{\substack{\ \\ \Diamond}}  \le 1$ and using
Lemma \ref{Tbdlem} one gets

\[
||
(M+N_{n-1})
\circ \cdots \circ
(M+N_{m})
-
 M^{n-m}
||_{\substack{\ \\ \Diamond}} 
\le
\mathcal{C}_{1}(\epsilon)
\left[
||N_{n-1}||_{\substack{\ \\ \Diamond}} 
+
\cdots
+
||N_m||_{\substack{\ \\ \Diamond}} 
\right].
\]

Thus one has

\[
|| P_{\rm s} T_{n}(v) ||_{\substack{\ \\ \Diamond}} 
\le
\mathcal{C}_{1}(\epsilon)
\alpha_{\rm u}^{-(n-m)}
c_{1}(\epsilon)^{n-m}
+
\mathcal{C}_{1}(\epsilon)^{2}
\left[
||N_{n-1}||_{\substack{\ \\ \Diamond}} 
+
\cdots
+
||N_m||_{\substack{\ \\ \Diamond}} 
\right].
\]

We now note that the proof of Lemma \ref{Tbdlem} tells us that:

\begin{equation}
||
N_{k}
||_{\substack{\ \\ \Diamond}}
\le
170
c_{1}(\epsilon)^{k}
||v-v_{*}||
\le
85
c_{1}(\epsilon)^{k}.
\label{Ndiameq}
\end{equation}\\

Using this in the previous inequality gives the bound

\[
||
P_{\rm s}T_{n}(v)
||_{\substack{\ \\ \Diamond}}
\le 
\mathcal{C}_{1}(\epsilon)
[\alpha_{\rm u}^{-1}c_{1}(\epsilon)]^{n-m}
+
85
\mathcal{C}_{1}(\epsilon)^{2}
\frac{c_{1}(\epsilon)^{m}}{1-c_{1}(\epsilon)}.
\]

Now take $\displaystyle m =
\left\lfloor \frac{n}{2} \right\rfloor$. Then
$\displaystyle m > \frac{n}{2}-1$ and one has

\[ 
c_{1}(\epsilon)^{m} 
\le
c_{1}(\epsilon)^{\frac{n}{2}-1}.
\]

One also has $\displaystyle n - m \ge \frac{n}{2}$ so that

\[
[\alpha_{\rm u}^{-1}c_{1}(\epsilon)]^{n-m}
\le
[\alpha_{\rm u}^{-1}c_{1}(\epsilon)]^{\frac{n}{2}}
\le
c_{1}(\epsilon)^{\frac{n}{2}}
\]

Note that we used the fact that
$\alpha_{\rm u} > 1$. Now inserting these two bounds into our last bound
for $|| P_{\rm s}T_{n}(v)||_{\substack{\ \\ \Diamond}}$ gives

\[
||
P_{\rm s}T_{n}(v)
||_{\substack{\ \\ \Diamond}}
\le 
c_{1}(\epsilon)^{\frac{n}{2}}
\left[
\mathcal{C}_{1}(\epsilon)
+
\frac{85\mathcal{C}_{1}(\epsilon)^{2}}{c_{1}(\epsilon)(1-c_{1}(\epsilon))}
\right].
\]
\qed 

Now we bound differences of the form $||T_{n+1}(v) - T_{n}(v) ||_{\substack{\ \\
    \Diamond}}$. 
    
\begin{Lemma}\label{Tcvlem}
For $v \in W^{\rm s,loc}$ and $n \ge 0$

\[
||
T_{n+1}(v)
-
T_{n}(v)
||_{\substack{\ \\ \Diamond}} 
\le
\mathcal{C}_{3}(\epsilon)
c_{1}(\epsilon)^{\frac{n}{2}}
\]

where $\displaystyle \mathcal{C}_{3}(\epsilon) = 85
\mathcal{C}_{1}(\epsilon) +
(1+\alpha_{\rm u}^{-1}c_{1}(\epsilon))\mathcal{C}_{2}(\epsilon)$.

\end{Lemma}   
    
\noindent{\bf Proof:}
Using the same notation as earlier we have

\begin{equation*}
\begin{split}
T_{n+1}(v) - T_{n}(v)
=&
(M+N_{n} - I)
\circ
T_{n}(v)\\
=&
N_{n}
\circ
T_{n}(v)
-
(I-M)
\circ
T_{n}(v)\ ,
\end{split}
\end{equation*}

but $M = P_{\rm u} + \alpha_{\rm u}^{-1} A_{\rm s} P_{\rm s}$ and $I=P_{\rm u}+P_{\rm s}$ so
we have

\[
I
-
M
=
(I - \alpha_{\rm u}^{-1}A_{\rm s})
\circ
P_{\rm s}\ .
\]

Hence
\begin{equation*}
\begin{split}
||
(I-M)
\circ
T_{n}(v)
||_{\substack{\ \\ \Diamond}}
=&
||
(I - \alpha_{\rm u}^{-1}A_{\rm s}P_{\rm s})
P_{\rm s}
T_{n}(v)
||_{\substack{\ \\ \Diamond}}\\
\le&
||
I
-
\alpha_{\rm u}^{-1}
A_{\rm s}
P_{\rm s}
||_{\substack{\ \\ \Diamond}}
\times
||
P_{\rm s}
T_{n}(v)
||_{\substack{\ \\ \Diamond}}\\
\le&
(1 + \alpha_{\rm u}^{-1})
\times
( \mathcal{C}_{2}(\epsilon) c_{1}(\epsilon)^{\frac{n}{2}})\ .\\
\end{split}
\end{equation*}

In going to last line we used Lemma \ref{pstnlem} to bound
$||P_{\rm s}T_{n}(v)||_{\substack{\ \\ \Diamond}}$. Now by
the estimate in (\ref{Ndiameq}) and Lemma \ref{Tbdlem} we have:

\[
||
N_{n}
T_{n}(v)
||_{\substack{\ \\ \Diamond}} 
\le
85c_{1}(\epsilon)^{n}
\times 
\mathcal{C}_{1}(\epsilon)\ .
\]

Thus we have the bound

\[
||
T_{n+1}(v)
-
T_{n}(v)
||_{\substack{\ \\ \Diamond}}
\le
( 1 + \alpha_{\rm u}^{-1}c_{1}(\epsilon))
\mathcal{C}_{2}(\epsilon)
c_{1}(\epsilon)^{\frac{n}{2}}
+
85
\mathcal{C}_{1}(\epsilon)
c_{1}(\epsilon)^{n}
\]
and the lemma follows.
\qed\\

The last lemma implies that $\displaystyle T_{\infty}(v) = \lim_{n
  \rightarrow \infty} T_{n}(v)$ exists and is a continuous linear
operator on $\mathcal{E}$. We also have the following as consequences
of Lemma \ref{Tcvlem} and Lemma \ref{Tbdlem}:

\[ 
P_{\rm s}T_{\infty}(v) = 0 
\textnormal{ and }
||T_{\infty}(v) ||_{\substack{\ \\ \Diamond}} 
\le
\mathcal{C}_{1}(\epsilon)\ .
\]\\

At the heart of the proof of Theorem \ref{unnamedtheorem1}
is a somewhat involved telescopic sum argument which will reappear many times in the remainder of this section
in slightly different forms. It first features in the following lemma which merely ensures that quantities of interest
are well defined.

\begin{Lemma}\label{RGaRGblem}
The following holds in the small $\epsilon$ regime. For all $v \in
W^{\rm s,loc}$ and $w\in\mathcal{E}$ such that

\[
||w||
\le
\frac{1}{240\mathcal{C}_{1}(\epsilon)}
\]

we have that, for all integers $n,a,b \ge 0$ such that $a+b \le n$, the
following expression is well defined:

\[
RG^{a}
\left( 
RG^{b}(v) 
+ 
D_{v}
RG^{b}[\alpha_{\rm u}^{-n}w]
\right).
\]\\

We also have the bound

\[
||
RG^{a}
\left( 
RG^{b}(v) 
+ 
D_{v}
RG^{b}[\alpha_{\rm u}^{-n}w]
\right)
||
\le
\frac{1}{8}.
\]

\end{Lemma}

\noindent{\bf Proof:}
Note that since $v \in W^{\rm s,loc}$ and $W^{\rm s,loc}$ is stable under
$RG$ we have that $RG^{k}(v)$ is well defined for all $k \ge 0$. We also have that $RG^{k}(v) \in W^{\rm s,loc}$ so we by
definition get the bound

\[
||RG^{k}(v)||
\le
\frac{\rho}{3} 
\le
\frac{1}{36}
\]
since
$\rho < \frac{1}{12}$. \\

For $0 \le k \le n$ we have that $\displaystyle D_{v}
RG^{k}[\alpha_{\rm u}^{-n}w]$ is well defined and in fact equal to
$\displaystyle \alpha_{\rm u}^{-(n-k)} T_{k}(v)[w]$. Noting that
$\alpha_{\rm u} > 1$ gives the estimate

\begin{equation*}
\begin{split}
||
D_{v}
RG^{k}[\alpha_{\rm u}^{-n}w]
||
\le&
||T_{k}(v)[w]||
\le
2
||
T_{k}(v)[w]
||_{\substack{\ \\ \Diamond}}\\
\le&
2
||T_{k}(v) ||_{\substack{\ \\ \Diamond}} 
||w ||_{\substack{\ \\ \Diamond}}\\
\le&
2
\mathcal{C}_{1}(\epsilon)
\times
5
||w||\ .
\end{split}
\end{equation*}

Thus we have that

\[
||
D_{v}RG^{k}[\alpha_{\rm u}^{-n}w]
||
\le
10
\mathcal{C}_{1}(\epsilon)
||w||
\le
\frac{1}{24}\ .
\]

We prove the assertion of the lemma by looking at various cases while
applying induction on $a+b$.\\

For our first case assume $a = 0$. We are then looking at

\[
RG^{a}
\left( 
RG^{b}(v)
+
D_{v}RG^{b}[\alpha_{\rm u}^{-n}w] 
\right)
=
RG^{b}(v)
+
D_{v}RG^{b}[\alpha_{\rm u}^{-n}w].
\]

The right hand side is well defined by the previous remarks on the two
pieces composing it. We also have the bound

\[
||
RG^{b}(v)
+
D_{v}RG^{b}[\alpha_{\rm u}^{-n}w]
||
\le
\frac{1}{24}
+
\frac{1}{24}
=
\frac{1}{12}\ .
\]

This proves the claims of our lemma whenever $a=0$ and also covers the
induction base case $a+b=0$.\\

For the second case assume $a+b >0$. Note that if
$a=0$ then we are again under the previous case for which
the assertions has been proved. Therefore we assume $a>0$. By our
induction hypothesis we have that 

\[
RG^{a-1}
\left(
RG^{b}(v)
+
D_{v}
RG^{b}[\alpha_{\rm u}^{-n}w]
\right)
\textnormal{  is well defined}
\]

and we also have the bound

\[
||
RG^{a-1}
\left(
RG^{b}(v)
+
D_{v}
RG^{b}[\alpha_{\rm u}^{-n}w]
\right)
||
\le
\frac{1}{8}.
\]

This places it within the domain of $RG$ (which is defined on
vectors of norm less than $\frac{1}{2}$).\\

Thus
\[ 
RG^{a}
\left(
RG^{b}(v)
+
D_{v}
RG^{b}[\alpha_{\rm u}^{-n}w]
\right)
\textnormal{ is well defined.}
\]

By the same argument the following quantities are also well defined:

\begin{equation*}
\begin{split}
RG^{a-1}
\left(
RG^{b+1}(v)
+
D_{v}
RG^{b+1}[\alpha_{\rm u}^{-n}w]
\right),\\
RG^{a-2}
\left(
RG^{b+2}(v)
+
D_{v}
RG^{b+2}[\alpha_{\rm u}^{-n}w]
\right),\\
\vdots \hspace{3cm} \\
RG
\left(
RG^{b+a-1}(v)
+
D_{v}
RG^{b+a-1}[\alpha_{\rm u}^{-n}w]
\right).
\end{split}
\end{equation*}

We write the telescopic sum

\begin{multline}\label{abtelescopesum}
RG^{a}
\left( 
RG^{b}(v)
+ 
D_{v}RG^{b}[\alpha_{\rm u}^{-1}w]
\right)
=
RG^{a+b}(v)
+
D_{v}RG^{a+b}[\alpha_{\rm u}^{-n}w]\\
+
\sum_{j=0}^{a-1}
\left\{
RG^{j+1}
\left(
RG^{b+a-j-1}(v)
+
D_{v}
RG_{v}^{b+a-j-1}[\alpha_{\rm u}^{-n}w]
\right) 
-
RG^{j}
\left(
RG^{b+a-j}(v)
+
D_{v}
RG^{b+a-j}[\alpha_{\rm u}^{-n}w]
\right)
\right\}\ .
\end{multline}

Note that Part 1) and Part 3) of Lemma \ref{L57lem} can be combined into
a single Lipschitz estimate

\[
||RG(w')
-
RG(w'')
||
\le
c_{3}(\epsilon)
||w'-w''||
\]
for all $w'$, $w''$ in $\mathcal{E}$ such that
$||w'||, ||w''||
\le
\frac{1}{8}$.

For $j \ge 1$ within the telescoping sum our induction hypothesis
tells us that we can repeatedly use our Lipschitz estimate $j$ times.
At every step the
arguments of the map $RG$ will be within $\displaystyle
\bar{B} \left( 0, \frac{1}{8} \right) \subset \mathcal{E}$. Thus we have:

\begin{gather*}
||
RG^{j+1}
\left(
RG^{b+a-j-1}(v)
+
D_{v}RG^{b+a-j}[\alpha_{\rm u}^{-n}w]
\right)
-
RG^{j}
\left(
RG^{b+a-j-1}(v)
+
D_{v}RG^{b+a-j-1}[\alpha_{\rm u}^{-n}w]
\right)
||\\
\le
c_{3}(\epsilon)^{j}
||
RG
\left(
RG^{b+a-j}(v)
+
D_{v}
RG^{b+a-j-1}[\alpha_{\rm u}^{-n}w]
\right)
-
\left(
RG^{b+a-j}(v)
+
D_{v}
RG^{b+a-j}[\alpha_{\rm u}^{-n}w]
\right)
||\ .
\end{gather*}

Note that the bound above holds for $j=0$ as well so we can apply this
estimate to all the terms of the telescoping sum:

\begin{equation*}
\begin{split}
&\left| \left|
\sum_{j=0}^{a-1}
\left\{
RG^{j+1}
\left(
RG^{b+a-j-1}(v)
+
D_{v}
RG_{v}^{b+a-j-1}[\alpha_{\rm u}^{-n}w]
\right) 
-
RG^{j}
\left(
RG^{b+a-j}(v)
+
D_{v}
RG^{b+a-j}[\alpha_{\rm u}^{-n}w]
\right)
\right\}
\right| \right|\\
&\le
\sum_{j=0}^{a-1}
c_{3}(\epsilon)^{j}
\left| \left|
RG
\left(
RG^{b+a-j-1}(v)
+
D_{v}
RG^{b+a-j-1}[\alpha_{\rm u}^{-n}w]
\right)
\right. \right.\\
&
\hspace{2cm}
\left. \left.
-RG \left( RG^{b+a-j-1}(v) \right)
-D_{RG^{b+a-j-1}(v)}RG
\left[
D_{v}
RG^{b+a-j-1}[\alpha_{\rm u}^{-n}w]
\right]
\right| \right|.
\end{split}
\end{equation*}

Above we used the chain rule for Frechet differentials.\\

 By  the earlier
remarks we know that $\displaystyle ||RG^{b+a-j-1}(v)|| \le \frac{1}{24}$
and $\displaystyle ||D_{v} RG_{b+a-j-1}[\alpha_{\rm u}^{-n}w]|| \le
\frac{1}{24}$. Thus the quantities appearing in the sum above can be
estimated using Lemma \ref{DHDRGlem} which tells us that

\begin{equation*}
\begin{split}
&\left| \left|
RG
\left(
RG^{b+a-j-1}(v)
+
D_{v}
RG^{b+a-j-1}[\alpha_{\rm u}^{-n}w]
\right)
\right. \right.\\
&
\hspace{.5cm}
\left. \left.
-RG \left( RG^{b+a-j-1}(v) \right)
-D_{RG^{b+a-j-1}(v)}RG
\left[
D_{v}
RG^{b+a-j-1}[\alpha_{\rm u}^{-n}w]
\right]
\right| \right|\\
\le&
\frac{17}{2}
||
D_{v}
RG^{b+a-j-1}[\alpha_{\rm u}^{-n}w]
||^{2}\\
=&
\frac{17}{2}
 \left[ 
\alpha_{\rm u}^{-n+(b+a-j-1)}
||T_{b+a-j-1}(v)[w]||
\right]^{2}\\
\le&
\frac{17}{2}
\left[
\alpha_{\rm u}^{-n+(b+a-j-1)}
\times
10
\times
||T_{b+a-j-1}(v)||_{\substack{\ \\ \Diamond}}
||w||
\right]^{2}\\
\le&
50
\times
17
\alpha_{\rm u}^{-2(n-a-b+j+1)}
\mathcal{C}_{1}(\epsilon)^2
||w||^{2}\ .
\end{split}
\end{equation*}

Inserting all of our bounds into \eqref{abtelescopesum} yields the inequality

\[
||
RG^{a} 
\left(
RG^{b}(v)
+
D_{v}RG^{b}[\alpha_{\rm u}^{-n}w]
\right)
||
\le
\frac{1}{24}
+
\frac{1}{24}
+
\sum_{j=0}^{a-1} 
\left[
c_{3}(\epsilon)^j 
850 
\alpha_{\rm u}^{-2(j+1)}
\mathcal{C}_{1}(\epsilon)^2
||w||^2
\right]\ .
\]

Above we used that $\alpha_{\rm u} > 1$ and $n \ge a + b$ which means
$\alpha_{\rm u}^{-2(n-a-b)} \le 1$.\\

 We would like to sum the geometric
series but for this we need to show that $\alpha_{\rm u}^{-2}c_{3}(\epsilon) < 1$
which we now do.\\

We have that  $\displaystyle \alpha_{\rm u} \ge L^{\frac{3+\epsilon}{2}} - c_{4}(\epsilon)
> 0$ for $\epsilon$ small by Lemma \ref{alphabdlem} and $\displaystyle c_{3}(\epsilon) =
L^{\frac{3+\epsilon}{2}} + c_{4}(\epsilon) $ by the definitions given in Lemma \ref{L57lem}. We also know
that $\displaystyle \lim_{\epsilon \rightarrow 0} c_{4}(\epsilon) =0$
so

\[
\alpha_{\rm u}^2 - c_{3}(\epsilon)
\ge
(L^{\frac{3+\epsilon}{2}} - c_{4}(\epsilon))^{2}
- 
(L^{\frac{3+\epsilon}{2}} + c_{4}(\epsilon))
\rightarrow L^{3} - L^{\frac{3}{2}}
\]
when $\epsilon \rightarrow 0$. 

We note that for $L \ge 2$ one has $\displaystyle
L^{3}-L^{\frac{3}{2}} > \frac{1}{2}L^{3}$. It then follows that for
$\epsilon$ sufficiently small one has

\begin{equation}
\alpha_{\rm u}^2 
- 
c_{3}(\epsilon)
\ge
\frac{1}{2}L^{3}.
\label{auc3eq}
\end{equation}

We have shown $\alpha_{\rm u}^{-2}c_{3}(\epsilon) < 1$. Thus

\begin{equation*}
\begin{split}
||
RG^{a} 
\left(
RG^{b}(v)
+
D_{v}RG^{b}[\alpha_{\rm u}^{-n}w]
\right)
||
\le&
\frac{1}{24}
+
\frac{1}{24}
+
850
\mathcal{C}_{1}(\epsilon)^2
||w||^2
\frac{\alpha_{\rm u}^{-2}}{1-\alpha_{\rm u}^{-2}c_{3}(\epsilon)}
\\
\le&
\frac{1}{24}
+
\frac{1}{24}
+
850 \times
\left(
\frac{1}{240}
\right)^{2}\times 
\frac{2}{L^{3}}
<
\frac{1}{24}
+ 
\frac{1}{24}
+
\frac{1}{24}
=
\frac{1}{8}.
\end{split}
\end{equation*}

In going to the last line we have used the lemma's assumption that
$||w|| \le \frac{1}{240\mathcal{C}_{1}(\epsilon)}$ along with the fact
that $\displaystyle
\frac{\alpha_{\rm u}^{-2}}{1-\alpha_{\rm u}^{-2}c_{3}(\epsilon)} =
\frac{1}{\alpha_{\rm u}^{2} - c_{3}(\epsilon)} \le \frac{2}{L^{3}}$. This
proves the bound asserted by the lemma. \qed\\

Note that the constant $\mathcal{C}_1(\epsilon)$ featuring in the domain definition for $w$ is a
very bad one. Indeed it essentially blows up as $\exp(\epsilon^{-1})$ when $\epsilon\rightarrow 0$.
This is because the previous lemma uniformly covers all starting points $v$ in the local
stable manifold $W^{\rm s,loc}$ on which the convergence to the fixed point $v_{\ast}$ is very slow.
In the special case $v=v_{\ast}$ and $w\in\mathcal{E}^{\rm u}$
one can obtain significantly better estimates which is what we do next.

\begin{Lemma}\label{RGaRGbbislem}
The following holds in the small $\epsilon$ regime. For all $w\in\mathcal{E}^{\rm u}$ such that

\[
||w||
\le
\frac{1}{24}
\]

we have that, for all integers $n,a,b \ge 0$ such that $a+b \le n$, the
following expression is well defined:

\[
RG^{a}
\left( 
RG^{b}(v_{\ast}) 
+ 
D_{v_{\ast}}
RG^{b}[\alpha_{\rm u}^{-n}w]
\right).
\]\\

We also have the bound

\[
||
RG^{a}
\left( 
RG^{b}(v_{\ast}) 
+ 
D_{v_{\ast}}
RG^{b}[\alpha_{\rm u}^{-n}w]
\right)
||
\le
\frac{1}{8}.
\]

\end{Lemma}

\noindent{\bf Proof:}
The proof is exaclty the same as that of the Lemma \ref{RGaRGblem} except for the following modifications.
When estimating $D_{v} RG^k[\alpha_{\rm u}^{-n}w]$ we now have the tremendous simplification
\[
D_{v_{\ast}} RG^k[\alpha_{\rm u}^{-n}w]=\alpha_{\rm u}^{-n}(D_{v_{\ast}}RG)^k[w]=\alpha_{\rm u}^{k-n} w
\]
by Lemma \ref{alphabdlem} and the hypothesis $w\in\mathcal{E}^{\rm u}$.
When we estimated the quantity $||T_{b+a-j-1}(v)[w]||$ we had to pay a factor of $10\mathcal{C}_{1}(\epsilon)$.
Now we simply note that $T_{b+a-j-1}(v)[w]=w$ and therefore we obtain the same bounds without this bad factor.
\qed

We now attack the main linearization theorem.\\

We proceed by showing that for $v \in W^{\rm s,loc}$ and $w$ sufficiently
small the following sum converges:

\[
\sum_{n=0}^{\infty}
||RG^{n+1}(v + \alpha_{\rm u}^{-(n+1)}w)
-
RG^{n}(v+\alpha_{\rm u}^{-n}w)||
\]
as results for the next lemma.

\begin{Lemma}\label{unnamedlemma1} 
In the small $\epsilon$ regime, for all $v \in W^{\rm s,loc}$ and all $w$
with $\displaystyle ||w|| \le \frac{1}{240 \mathcal{C}_{1}(\epsilon)}$\\

we have, for all $n \ge 0$,

\[
||
RG^{n+1}(v+\alpha_{\rm u}^{-(n+1)}w)
-
RG^{n}(v+\alpha_{\rm u}^{-n}w)
||
\le
\mathcal{C}_{4}(\epsilon)
c_{1}(\epsilon)^{\frac{n}{4}}
\]\\

where $ \displaystyle \mathcal{C}_{4}(\epsilon) =
\frac{\alpha_{\rm u}^{2}}{12c_{3}(\epsilon)}+\frac{1}{2} +
\frac{\mathcal{C}_{3}(\epsilon)}{24\mathcal{C}_{1}(\epsilon)}$.

\end{Lemma}

\noindent{\bf Proof:}
By Lemma \ref{RGaRGblem} with $b=0$ and $a=n$ the quantities involved
in the sum are all well defined if $\displaystyle ||w|| \le
\frac{1}{240\mathcal{C}_{1}(\epsilon)}$.\\

We now proceed similarly to the proof of Lemma \ref{RGaRGblem}.
For $n \ge 1$ and $0 \le k \le n$ we write the telescoping sum

\begin{multline*}
RG^{n}(v+\alpha_{\rm u}^{-n}w)
-
RG^{k}(RG^{n-k}(v) + D_{v}RG^{n-k}[\alpha_{\rm u}^{-n}w])\\
=
\sum_{j=k}^{n-1} 
\left\{
RG^{j+1}
\left(
RG^{n-j-1}(v)
+
D_{v}
RG^{n-j-1}[\alpha_{\rm u}^{-n}w]
\right)
-
RG^{j}
\left(
RG^{n-j}(v)
+
D_{v}
RG^{n-j}[\alpha_{\rm u}^{-n}w]
\right)
\right\}\ .
\end{multline*}

Since $k \le j \le n-1$ the arguments of $RG$ remain small enough to
allow for repeated use of
Lipschitz estimates giving the bound:

\begin{equation*}
\begin{split}
&||
RG^{j+1}
\left(
RG^{n-j-1}(v)
+
D_{v}
RG^{n-j-1}[\alpha_{\rm u}^{-n}w]
\right)
-
RG^{j}
\left(
RG^{n-j}(v)
+
D_{v}
RG^{n-j}[\alpha_{\rm u}^{-n}w]
\right)
||\\
\le&
c_{3}(\epsilon)^{j}
||
RG
\left(
RG^{n-j-1}(v)
+
D_{v}
RG^{n-j-1}[\alpha_{\rm u}^{-n}w]
\right)
-
RG^{n-j}(v)
-
D_{v}
RG^{n-j}[\alpha_{\rm u}^{-n}w]
||\\
=&
c_{3}(\epsilon)^{j}
||
RG
\left(
RG^{n-j-1}(v)
+
D_{v}
RG^{n-j-1}[\alpha_{\rm u}^{-n}w]
\right)
-
RG^{n-j}(v)
-
D_{RG^{n-j-1}(v)}
RG
\left[
D_{v}
RG^{n-j-1}[\alpha_{\rm u}^{-n}w]
\right]
||\\
\le&
c_{3}(\epsilon)^{j}
\frac{17}{2}
||
D_{v}
RG^{n-j-1}[\alpha_{\rm u}^{-n}w]
||^{2}\\
\le&
c_{3}(\epsilon)^{j}
\frac{17}{2}
\left[
\alpha_{\rm u}^{-(j+1)}
||
T_{n-j-1}(v)[w]
||
\right]^{2}\\
\le& 
850
\mathcal{C}_{1}(\epsilon)^{2}
||w||^{2}
\times
\alpha_{\rm u}^{-2}
\times
\left[
c_{3}(\epsilon)
\alpha_{\rm u}^{-2}
\right]^{j}\ .
\end{split}
\end{equation*}

Using this estimate in the earlier telescoping sum expression gives the bound

\begin{equation*}
\begin{split}
||
RG^{n}(v+\alpha_{\rm u}^{-n}w)
-
RG^{k}(RG^{n-k}(v) +& D_{v}RG^{n-k}[\alpha_{\rm u}^{-n}w])
||\\
\le&
850
\mathcal{C}_{1}(\epsilon)^{2}
||w||^{2}
\alpha_{u}^{-2}
\sum_{j=k}^{n-1}
\left( 
c_{3}(\epsilon)
\alpha_{\rm u}^{-2}
\right)^{j}\\
\le& 
850
\mathcal{C}_{1}(\epsilon)^{2}
||w||^{2} 
\left( 
c_{3}(\epsilon)
\alpha_{\rm u}^{-2}
\right)^{k}
\frac{\alpha_{\rm u}^{-2}}{1-\alpha_{\rm u}^{-2}c_{3}(\epsilon)}
\le
\frac{1}{24}
\left( 
c_{3}(\epsilon)
\alpha_{\rm u}^{-2}
\right)^{k}\ .
\end{split}
\end{equation*}

For the last inequality we proceeded
just as we did at the end of the proof of Lemma \ref{RGaRGblem}. Note
that this bound holds for any $n,k \ge 0$ with $k \le n$ so we get a
valid estimate if we replace $n$ with $n+1$:

\[
||
RG^{n+1}(v+\alpha_{\rm u}^{-(n+1)}w)
-
RG^{k}(RG^{n+1-k}(v) + D_{v}RG^{n+1-k}[\alpha_{\rm u}^{-(n+1)}w])
||
\le
\frac{1}{24}
\left( 
c_{3}(\epsilon)
\alpha_{\rm u}^{-2}
\right)^{k}\ .
\]

As a result, the triangle inequality gives

\begin{multline*}
||
RG^{n+1}
(v+\alpha_{\rm u}^{-(n+1)}w)
-
RG^{n}(v+\alpha_{\rm u}^{-n}w)
||\\
\le
\frac{1}{12}
\left(
c_{3}(\epsilon)
\alpha_{\rm u}^{-2}
\right)^{k}\\
\ \ \ +
\left|\left|
RG^{k} 
\left( 
RG^{n+1-k}(v)
+
D_{v}RG^{n+1-k}[\alpha_{\rm u}^{-(n+1)}w]
\right)
-
RG^{k} 
\left( 
RG^{n-k}(v)
+
D_{v}RG^{n-k}[\alpha_{\rm u}^{-n}w]
\right)
\right|\right|\ .
\end{multline*}
We again repeatedly used the Lipschitz estimate to bound the second term
on the bottom line. Indeed, Lemma \ref{RGaRGblem} guarantees that the
arguments of the outermost $RG$'s remain in the domain of validity of
our Lipschitz estimate. This gives the bound

\begin{equation*}
\begin{split}
&\left|\left|
RG^{k} 
\left( 
RG^{n+1-k}(v)
+
D_{v}RG^{n+1-k}[\alpha_{\rm u}^{-(n+1)}w]
\right)
-
RG^{k} 
\left( 
RG^{n-k}(v)
+
D_{v}RG^{n-k}[\alpha_{\rm u}^{-n}w]
\right)
\right|\right|\\[1.5ex]
\le&
c_{3}(\epsilon)^{k}
||
RG^{n+1-k}(v)
+
D_{v}RG^{n+1-k}[\alpha_{\rm u}^{-(n+1)}w]
-
RG^{n-k}(v)
-
D_{v}RG^{n-k}[\alpha_{\rm u}^{-n}w]
||\\[1.5ex]
\le&
c_{3}(\epsilon)^{k}
\left[
||
RG^{n+1-k}(v)
-
RG^{n-k}(v)
||
+
||
D_{v}
RG^{n+1-k}[\alpha_{\rm u}^{-(n+1)}w]
-
D_{v}
RG^{n-k}[\alpha_{\rm u}^{-n}w]
||
\right]\\[1.5ex]
=&
c_{3}(\epsilon)^{k}
\left[
||
RG^{n-k} 
\left(
RG(v)
\right)
-
RG^{n-k}(v)
||
+
\alpha_{\rm u}^{-k}
||
T_{n-k+1}(v)[w]
-
T_{n-k}(v)[w]
||
\right]\\[1.5ex]
\le&
c_{3}(\epsilon)^{k}
\left[
c_{1}(\epsilon)^{n-k}
||
RG(v)
-
v
||
+
\alpha_{\rm u}^{-k}
\times
10
\times
||
T_{n-k+1}(v)
-
T_{n-k}(v)
||_{\substack{\ \\ \Diamond}}
||w||
\right]\\[1.5ex]
\le& 
\frac{1}{2}
c_{3}(\epsilon)^{k}
c_{1}(\epsilon)^{n-k}
+
c_{3}(\epsilon)^{k}
\alpha_{\rm u}^{-k}
\times
10
\times
\mathcal{C}_{3}(\epsilon)
c_{1}(\epsilon)^{\frac{n-k}{2}}
\times
\frac{1}{240 \mathcal{C}_{1}(\epsilon)}\ .
\end{split}
\end{equation*}

In going to the fifth line from the fourth line we used Lemma
\ref{L59lem} since both $v$ and $RG(v)$ are in $W^{\rm s,loc}$. In going from the
fifth line to the last line we used $\displaystyle ||RG(v) - v|| \le
||RG(v)|| + ||v|| \le \frac{1}{4} + \frac{1}{4} = \frac{1}{2}$ when
bounding the first term. For the second term we used Lemma
\ref{Tcvlem} and our working assumption on the size of $||w||$.\\ 

We now arrive at

\begin{equation*}
\begin{split}
||
RG^{n+1}
(v+\alpha_{\rm u}^{-(n+1)}w)
-
RG^{n}(v+\alpha_{\rm u}^{-n}w)
||
\le&
\frac{1}{2}
\left( c_{3}(\epsilon)\alpha_{\rm u}^{-2} \right)^{k}
+
\frac{1}{2}
c_{3}(\epsilon)^{k}
c_{1}(\epsilon)^{n-k}
+
\frac{\mathcal{C}_{3}(\epsilon)}{24\mathcal{C}_{1}(\epsilon)}
c_{3}(\epsilon)^{k}
\alpha_{u}^{-k}
c_{1}(\epsilon)^{\frac{n-k}{2}}\ .
\end{split}
\end{equation*}

We now choose $k$ adequately as a function of $n$. We will set 

\[
k
=
\lfloor 
\sigma
n
\rfloor
\textnormal{  for suitable }
\sigma
\in
[0,1]\ .
\]

We then have that $0 \le k \le n$. From previous arguments we know
that $0 < c_{3}(\epsilon) \alpha_{\rm u}^{-2} < 1$. Then since $ k > \sigma n
-1$ we have

\[
\left( c_{3}(\epsilon) \alpha_{\rm u}^{-2} \right)^{k}
\le
\left( c_{3}(\epsilon) \alpha_{\rm u}^{-2} \right)^{\sigma n - 1}.
\]

Since $c_{3}(\epsilon) \ge 1$, $c_{1}(\epsilon) < 1$, and $k \le \sigma
n$ we have

\[
c_{3}(\epsilon)^{k}
c_{1}(\epsilon)^{n-k}
\le
c_{3}(\epsilon)^{\sigma n} 
c_{1}(\epsilon)^{n - \sigma n}.
\]

Since $ \displaystyle c_{3}(\epsilon)\alpha_{\rm u}^{-1} \ge 1$ which is
a consequence of $\alpha_{\rm u} \le L^{\frac{3+\epsilon}{2}} + c_{4}(\epsilon) =
c_{3}(\epsilon)$ and since  $\sqrt{c_{1}(\epsilon)} \le 1$, the inequality $k \le
\sigma n$ implies

\[
c_{3}(\epsilon)^{k}
\alpha_{\rm u}^{-k}
c_{1}(\epsilon)^{\frac{n-k}{2}}
\le
c_{3}(\epsilon)^{\sigma n}
\alpha_{\rm u}^{-\sigma n}
c_{1}(\epsilon)^{\frac{ n - \sigma n}{2}}.
\]

Using these three statements in our previous bound we see that

\[
||
RG^{n+1}(v+\alpha_{\rm u}^{-(n+1)} w) 
- 
RG^{n}(v+\alpha_{\rm u}^{-n} w)
||
\le
\left[
\frac{\alpha_{\rm u}^{2}}{12c_{3}(\epsilon)}
+
\frac{1}{2}
+
\frac{\mathcal{C}_{3}(\epsilon)}{24\mathcal{C}_{1}(\epsilon)}
\right]
\times
\gamma^{n}
\]

\[
\textnormal{ with }
\gamma
=
\max
\left[
c_{3}(\epsilon)^{\sigma} \alpha_{\rm u}^{-2\sigma},
c_{3}(\epsilon)^{\sigma} c_{1}(\epsilon)^{1-\sigma},
c_{3}(\epsilon)^{\sigma} \alpha_{\rm u}^{-\sigma} c_{1}(\epsilon)^{\frac{1 - \sigma}{2}}
\right]\ .
\]

Recall that for $\epsilon$ small we have

\[
\alpha_{\rm u}^{2}
>
c_{3}(\epsilon)
\ge
\alpha_{\rm u}
>
1
>
c_{1}(\epsilon)
>
0\ .
\]

Therefore when $\sigma \in [0,1]$ one has

\[
c_{3}(\epsilon)^{\sigma}
\alpha_{\rm u}^{-\sigma}
c_{1}(\epsilon)^{\frac{1 - \sigma}{2}}
\le
\left[
c_{3}(\epsilon)^{\sigma}
c_{1}(\epsilon)^{1-\sigma}
\right]^{\frac{1}{2}}.
\]

Hence

\[
\gamma
\le
\max
\left[
c_{3}(\epsilon)^{\sigma} \alpha_{\rm u}^{-2\sigma},
c_{3}(\epsilon)^{\sigma} c_{1}(\epsilon)^{1-\sigma},
\left(
c_{3}(\epsilon)^{\sigma}
c_{1}(\epsilon)^{1-\sigma}
\right)^{\frac{1}{2}}
\right]
<
1\ .
\]

The last inequality holds provided
$c_{3}(\epsilon)^{\sigma}\alpha_{\rm u}^{-2\sigma} < 1$ and
$c_{3}(\epsilon)^{\sigma}c_{1}(\epsilon)^{1-\sigma} < 1$ which can be guaranteed by
choosing $\sigma$ so that

\[
0
<
\sigma
<
\frac{-\log\ c_{1}(\epsilon)}{\log\ c_{3}(\epsilon) - \log\
  c_{1}(\epsilon)}\ .
\]

For simplicity,
we pick
\[
\sigma = \frac{1}{2} \times \frac{-\log\
  c_{1}(\epsilon)}{\log\ c_{3}(\epsilon) - \log\ c_{1}(\epsilon)} \in
\left( 0, \frac{1}{2} \right)\ . 
\]

We then have

\[
c_{3}(\epsilon)^{\sigma}
c_{1}(\epsilon)^{1-\sigma}
=
c_{1}(\epsilon)^{\frac{1}{2}}.
\]

So

\[
\gamma
\le
\max
\left[
c_{3}(\epsilon)^{\sigma}\alpha_{\rm u}^{-2\sigma},
c_{1}(\epsilon)^{\frac{1}{2}},
c_{1}(\epsilon)^{\frac{1}{4}}
\right]
=
c_{1}(\epsilon)^{\frac{1}{4}}
\textnormal{ in the small } \epsilon \textnormal{ regime.}\]

Indeed the strict inequality

\[
c_{3}(\epsilon)^{\sigma}
\alpha_{\rm u}^{-3\sigma} 
<
c_{1}(\epsilon)^{\frac{1}{4}} 
= 
c_{3}(\epsilon)^{\frac{\sigma}{2}}
c_{1}(\epsilon)^{\frac{1-\sigma}{2}}
\textnormal{ is successively equivalent to}
\]

\begin{equation*}
\begin{split}
\alpha_{\rm u}^{-2} 
&< 
c_{3}(\epsilon)^{-\frac{\sigma}{\frac{1}{2}-\sigma}
\times
\frac{1 - \sigma}{2}
-
\frac{\sigma}{2}}\ ,\\
\alpha_{\rm u}^{2}
&>
c_{3}(\epsilon)^{\frac{1-\sigma}{1-2\sigma}+\frac{1}{2}}\ ,\\
2 \log\ \alpha_{\rm u}
&>
\left(
1
+
\frac{1}{2(1-2\sigma)}
\right)
\log\ c_{3}(\epsilon)\ .
\end{split}
\end{equation*}

However, Lemma \ref{alphabdlem} gives

\[
2 \log\ \alpha_{\rm u}
\rightarrow
3\log\ L\ 
\textnormal{  when  }
\epsilon
\rightarrow
0
\]

while

\[
\left(
1
+
\frac{1}{2(1-2\sigma)}
\log\ c_{3}(\epsilon)
\right)
=
\frac{3}{2}
\log\ c_{3}(\epsilon) 
-
\frac{1}{2}
\log\ c_{1}(\epsilon)
\rightarrow
\frac{9}{4} 
\log\ L
<
3 \log\ L
\textnormal{  when  }
\epsilon
\rightarrow
0
\]
as is readily checked from the definitions in Lemma \ref{L57lem}.

Therefore $c_{3}(\epsilon)\alpha_{\rm u}^{-2\sigma} <
c_{1}(\epsilon)^{\frac{1}{4}}$ in the small $\epsilon$ regime and the result is proved.
\qed\\

The last lemma also has an improved version in the special case $v=v_{\ast}$, $w\in\mathcal{E}^{\rm u}$.

\begin{Lemma}\label{unnamedlemma1bis} 
In the small $\epsilon$ regime, for all $w\in\mathcal{E}^{\rm u}$
with $\displaystyle ||w|| \le \frac{1}{24}$\\

we have, for all $n \ge 0$,

\[
||
RG^{n+1}(v_{\ast}+\alpha_{\rm u}^{-(n+1)}w)
-
RG^{n}(v_{\ast}+\alpha_{\rm u}^{-n}w)
||
\le
\mathcal{C}'_{4}(\epsilon)
c_{1}(\epsilon)^{\frac{n}{4}}
\]\\

where $ \displaystyle \mathcal{C}'_{4}(\epsilon) =
\frac{\alpha_{\rm u}^{2}}{12c_{3}(\epsilon)}+\frac{1}{2}$.

\end{Lemma}

\noindent{\bf Proof:}
The proof is the same as that of Lemma \ref{unnamedlemma1} except that we use Lemma \ref{RGaRGbbislem}
instead of Lemma \ref{RGaRGblem}.
We bound $||T_{n-j-1}(v_{\ast})[w]||$ simply by $||w||\le \frac{1}{24}<\frac{2}{17}$ and do not pay a bad factor
$10\mathcal{C}_{1}(\epsilon)$. Also $T_{n-k+1}(v_{\ast})[w]-T_{n-k}(v_{\ast})[w]=0$ so the new constant $\mathcal{C}'_{4}(\epsilon)$
does not have the third term of $\mathcal{C}_{4}(\epsilon)$.
\qed\\

We now are in a position to state and prove our partial linearization theorem.\\

\begin{Theorem}\label{unnamedtheorem1}
For $v \in W^{\rm s,loc}$ and $\displaystyle ||w|| \le
\frac{1}{240\mathcal{C}_{1}(\epsilon)}$ the quantity

\[
\Psi(v,w) 
= 
\lim_{n \rightarrow \infty}
RG^{n}(v+\alpha_{\rm u}^{-n}w)
\textnormal{  exists in  }
\mathcal{E}
\]

and defines a function of $(v,w)$ with the following properties:

\begin{enumerate}

\item $\Psi$ is continuous in the domain $v \in W^{\rm s,loc}$ and $||w|| \le \frac{1}{240\mathcal{C}_{1}(\epsilon)}$. Over this set one has
the uniform bound $\displaystyle ||\Psi(v,w)|| \le \frac{1}{8}$.

\item $\Psi$ is jointly analytic in $v_{1}$ and $w$ in the domain
$\displaystyle ||v_{1}|| < \frac{\rho}{13}$, $\displaystyle ||w||
< \frac{1}{240\mathcal{C}_{1}(\epsilon)}$ where we have implied the
use of the parameterization

\[
v_{1} 
\mapsto 
v 
= 
(v_{1},v_{2}) 
= 
(v_{1},\mu_{\rm s}(v_{1}))
\textnormal{  of  }
W^{\rm s,loc}_{\rm int}\ .
\]\\

\item For all $v \in W^{\rm s,loc}$, $w$ such that $\displaystyle ||w|| \le
  \frac{1}{240\mathcal{C}_{1}(\epsilon)\alpha_{\rm u}}$ we have the intertwining relation

\[
RG(\Psi(v,w))
=
\Psi(v,\alpha_{\rm u}w).
\]

\item  For all $v \in W^{\rm s,loc}$, $w$ such that $\displaystyle ||w|| \le
  \frac{1}{2400\mathcal{C}_{1}(\epsilon)^2}$, and all integers $q \ge
  0$, we have

\[
\Psi(v,w)
=
\Psi(RG^{q}(v),T_{q}(v)[w]).
\] \\

\item For all $v \in W^{\rm s,loc}$ and $w$ such that $\displaystyle ||w|| \le
  \frac{1}{2400\mathcal{C}_{1}(\epsilon)^2}$, we have

\[
\Psi(v,w)
=
\Psi(v_{*},T_{\infty}(v)[w]).
\]

\end{enumerate}

\end{Theorem}

\noindent{\bf Proof:}
Parts 1) and 2) are immediate consequences of the $\displaystyle
\frac{1}{8}$ bound in Lemma \ref{RGaRGblem} and the uniform absolute
convergence proved in Lemma \ref{unnamedlemma1}.\\

For Part 3) note 

\[
\Psi(v,\alpha_{\rm u}w)
=
\lim_{n\rightarrow \infty}
RG
\left(
RG^{n-1}(v+\alpha_{\rm u}^{-(n-1)}w)
\right)
\]

and the continuity of $RG$ in the ball of radius $\displaystyle
\frac{1}{8}$.\\

For Part 4) we use the $c_{3}(\epsilon)$ Lipschitz estimate and Lemma
\ref{DHDRGlem} which are justified by Lemma \ref{RGaRGblem} in order to write for
fixed $q$ and $n \ge 0$:

\begin{equation*}
\begin{split}
&|| 
RG^{n+q}(v+\alpha_{\rm u}^{-(n+q)}w) 
- 
RG^{n}
(RG^{q}(v) 
+ 
D_{v} RG^{q}[\alpha_{\rm u}^{-(n+q)}w])
||
\\
\le&
c_{3}(\epsilon)^{n}
||
RG^{q}(v+\alpha_{\rm u}^{-(n+q)}w)
-
RG^{q}(v)
-
D_{v} RG^{q}[\alpha_{\rm u}^{-(n+q)}w]
||\\
\le&
c_{3}(\epsilon)^{n}
\sum_{j=0}^{q-1}
||
RG^{j+1}
\left(
RG^{q-(j+1)}(v)
+
D_{v}RG^{q-(j+1)}[\alpha_{\rm u}^{-(n+q)}w]
\right)
-
RG^{j}
\left(
RG^{q-j}(v)
+
D_{v}RG^{q-j}[\alpha_{\rm u}^{-(n+q)}w]
\right)
||\\
\le&
c_{3}(\epsilon)^{n}
\sum_{j=0}^{q-1}
c_{3}(\epsilon)^{j}
||
RG
(RG^{q-(j+1)}(v)
+
D_{v}
RG^{q-(j+1)}[\alpha_{\rm u}^{-(n+q)}w])
-
RG^{q-j}(v)
-
D_{v}RG^{q-j}[\alpha_{\rm u}^{-(n+q)}w]
||
\\
= &
c_{3}(\epsilon)^{n}
\sum_{j=0}^{q-1}
c_{3}(\epsilon)^{j}
||
RG
(RG^{q-(j+1)}(v)
+
D_{v}
RG^{q-(j+1)}[\alpha_{\rm u}^{-(n+q)}w])
-
RG^{q-j}(v)\\
 & -
D_{RG^{q-(j+1)}(v)}RG
\left[
D_{v}RG^{q-(j+1)}[\alpha_{\rm u}^{-(n+q)}w]
\right]
||\\
\le&
c_{3}(\epsilon)^{n}
\sum_{j=0}^{q-1}
c_{3}(\epsilon)^{j}
\times
\frac{17}{2}
\times
||D_{v}RG^{q-(j+1)}[\alpha_{\rm u}^{-(n+q)}w]||^2\\
\le&
c_{3}(\epsilon)^{n}
\sum_{j=0}^{q-1}
c_{3}(\epsilon)^{j}
\times
\frac{17}{2}
\times
\left[ 
\alpha_{\rm u}^{-(n+j+1)} 
\times 
10 
\times
 ||T_{q-(j+1)}(v)||_{\substack{\ \\ \Diamond}} 
\times 
||w|| 
\right]^{2}\ .
\end{split}
\end{equation*}

We now note that we can extract a factor of
$\left(c_{3}(\epsilon)\alpha_{\rm u}^{-2}\right)^{n}$ which will drive the
    expression to $0$ as  $n \rightarrow \infty$. Thus

\begin{equation*}
\begin{split}
\Psi(v,w)
=&
\lim_{n \rightarrow \infty}
RG^{n+q}[v+\alpha_{\rm u}^{-(n+q)}w]\\
=&
\lim_{n \rightarrow \infty}
RG^{n} 
\left[ 
RG^{q}(v) 
+ 
\alpha_{\rm u}^{-n} 
\left( 
\alpha_{\rm u}^{-q}
D_{v}RG^{q}[w]
\right)
\right]\\
=&
\Psi (RG^{q}(v), T_{q}(v)[w])
\end{split}
\end{equation*}

since $\alpha_{\rm u}^{-q}D_{v}RG^{q}[w] = T_{q}(v)[w]$ has norm bounded by

\[
||
T_{q}(v)
||
\times
||w||
\le
10
||
T_{q}(v)
||_{\substack{\ \\ \Diamond}} 
\times
\frac{1}{2400\mathcal{C}_{1}(\epsilon)^2}
\le
\frac{1}{240\mathcal{C}_{1}(\epsilon)}
\]

from Lemma \ref{Tbdlem}.\\

Part 5) follows from Part (4), Lemmas \ref{L59lem} and \ref{Tcvlem} when taking the $q\rightarrow\infty$ limit.
\qed\\

At this point it could seem possible that we went through all this trouble
in order to define a conjugation $\Psi$ which in fact is identically zero.
Our next theorem will rule this out thanks to the consideration of the special case
$v=v_{\ast}$ and $w\in\mathcal{E}^{\rm u}$.

\begin{Theorem}\label{unnamedtheorem1bis}
On the domain $||w||<\frac{1}{24}$ of the one-dimensional space $\mathcal{E}^{\rm u}$
the limit
\[
\lim_{n \rightarrow \infty}
RG^{n}(v_{\ast}+\alpha_{\rm u}^{-n}w)
\]
exists and defines an analytic function of $w$ which will be denoted by $\Psi(v_{\ast},w)$ since it coincides
with the previous one on the common domain of definition.
On the domain $B\left(0,\frac{1}{24}\right)\cap\mathcal{E}^{\rm u}$, this function satisfies the bound
\[
||\Psi(v_{\ast},w)||\le\frac{1}{8}
\]
as well as
\[
||\Psi(v_{\ast},w)-v_{\ast}-w||\le\frac{17}{8}||w||^{2}\ .
\]
In particular, the differential with respect to $w$ at $w=0$ is the identity on $\mathcal{E}^{\rm u}$.
On the domain $B\left(0,\frac{1}{24}\right)\cap\mathcal{E}^{\rm u}$ we also have the intertwining relation
\[
RG(\Psi(v_{\ast},\alpha_{\rm u}^{-1}w))=\Psi(v_{\ast},w)\ .
\]
For $w$ small enough in $\mathcal{E}^{\rm u}$ we have $\Psi(v_{\ast,w})\in W_{\rm tiny}^{\rm u,loc}$
\end{Theorem}

\noindent{\bf Proof:}
Lemma \ref{RGaRGbbislem} garantees that the quantities $RG^{n}(v_{\ast}+\alpha_{\rm u}^{-n}w)$
are well defined and bounded in norm by $\frac{1}{8}$. Lemma \ref{unnamedlemma1bis} shows the limit exists and is analytic
in $w$. Finally the same telescopic sum argument as in the proof of Lemma \ref{unnamedlemma1}, with $k=0$,
gives the estimate
\[
||RG^{n}(v_{\ast}+\alpha_{\rm u}^{-n}w)-(RG^n(v_{\ast})+D_{v_{\ast}}RG[\alpha_{\rm u}^{-n}w])||
\le \sum_{j=0}^{n-1} c_3(\epsilon)^j\times\frac{17}{2}\left[
\alpha_{\rm u}^{-(j+1)} ||T_{n-j-1}(v_{\ast})[w]||
\right]^{2}
\]
which in the present situation simply boils down to
\[
||RG^{n}(v_{\ast}+\alpha_{\rm u}^{-n}w)-v_{\ast}-w||
\le \sum_{j=0}^{n-1} c_3(\epsilon)^j\times\frac{17}{2}\left[
\alpha_{\rm u}^{-(j+1)} ||w||
\right]^{2}
\]
from which the wanted estimate follows easily.
The intertwining relation follows as in Part 3) of Theorem \ref{unnamedtheorem1}.
Using this intertwining relation to construct the backwards $RG$ trajectory $\Psi(v_{\ast},\alpha_{\rm u}^{n}w)$, for $n\le 0$,
and thanks to the criterion in Proposition \ref{L62prop} one easily see that $\Psi(v_{\ast},w)$ is in the local unstable manifold
$W^{\rm u,loc}$. By continuity of $\Psi(v_{\ast},\bullet)$ at zero one also gets the stronger conclusion that
$\Psi(v_{\ast},w)\in W_{\rm tiny}^{\rm u,loc}$. 
\qed

Before concluding this section we state a lemma which is on the same theme as Lemmas \ref{RGaRGblem}
and \ref{unnamedlemma1} and which will be needed in the sequel.

\begin{Lemma}\label{devseedlem}
In the small $\epsilon$ regime for all $v\in W^{\rm s,loc}$ and all $w$ with
$||w||\le\frac{1}{240\mathcal{C}_1(\epsilon)}$ we have
\[
||RG^{n}(v+\alpha_{\rm u}^{-n}w)-RG^{n}(v)||\le 11\mathcal{C}_1(\epsilon) ||w||\ .
\]
\end{Lemma}

\noindent{\bf Proof:}
By the same telescopic sum argument as in the beginning of the proof of Lemma \ref{unnamedlemma1}, with $k=0$, we get
\[
\begin{split}
||RG^{n}(v+\alpha_{\rm u}^{-n}w)-RG^{n}(v)-D_{v}RG^{n}[\alpha_{\rm u}^{-n}w]|| & \le
\sum_{j=0}^{n-1} 850\mathcal{C}_{1}(\epsilon)^2||w||^2\times \alpha_{\rm u}^{-2}\times(c_3(\epsilon)\alpha_{\rm u}^{-2})^j\\
 & \le  850\mathcal{C}_{1}(\epsilon)^2||w||^2\times \frac{1}{4}
\end{split}
\]
where we used (\ref{auc3eq}) and $L\ge 2$.
As a result we have
\[
\begin{split}
||RG^{n}(v+\alpha_{\rm u}^{-n}w)-RG^{n}(v)|| & \le ||T_n(v)[w]||+\frac{425}{2} \mathcal{C}_{1}(\epsilon)^2||w||^2\\
 & \le  10\mathcal{C}_1(\epsilon) ||w||+\frac{425}{2} \mathcal{C}_{1}(\epsilon)||w||\times\frac{1}{240}
\end{split}
\]
because of our hypothesis on $||w||$. Since $2\times 240>425$ the lemma follows.
\qed

%% file: Deviations.tex
\section{Control of the deviation from the bulk}\label{controldevsec}

\subsection{Algebraic considerations}\label{algconsec}

We now pick up the thread from \S\ref{algdefsec} where we consider for test
functions $\tilde{f},\tilde{j} \in
S_{q_{-},q_{+}}(\mathbb{Q}_{p}^{3},\mathbb{C})$ the quantity

\[
\mathcal{S}_{r,s}(\tilde{f},\tilde{j})
=
\frac{\mathcal{Z}_{r,s}(\tilde{f},\tilde{j})}
{\mathcal{Z}_{r,s}(0,0)}
\]

which is the moment generating function with UV and IR cutoffs $r$ and
$s$ respectively.\\

Introduce

\begin{equation*}
\begin{split}
\mathcal{S}^{\rm T}_{r,s}(\tilde{f},\tilde{j})
=&
-Y_{0}
Z_{0}^{r}
 \int_{\mathbb{Q}_{p}^{3}} \tilde{j}(x)\ d^3x
+
\frac{1}{2}
\sum_{r \le q < s}
\left(
f^{(r,q)},
\Gamma
f^{(r,q)}
\right)_{\Lambda_{s-q}}\\
&+ 
\sum_{r \le q < s}
\sum_{\substack{\Delta \in \mathbb{L} \\ \Delta \subset
    \Lambda_{s-q-1}}}
\left(
\delta b_{\Delta}
\left[
\vec{V}^{(r,q)} (\tilde{f},\tilde{j})
\right]
-
\delta b_{\Delta}
\left[
\vec{V}^{(r,q)} (0,0)
\right]
\right)\\
&+ 
{\rm Log}
\left(
\frac{\partial \mathcal{Z}_{r,s}(\tilde{f},\tilde{j})}{\partial
  \mathcal{Z}_{r,s}(0,0)}
\right)
\end{split}
\end{equation*}

where ${\rm Log}$ is the principal logarithm with argument in $(-\pi,\pi]$.\\

We will show that it is indeed a well defined quantity which boils
down to making sure  all the $RG$ iterates $\vec{V}^{(r,q)}$ are in
the domain of definition and analyticity for $RG_{\rm ex}$ provided by
Theorem \ref{mainestthm}.  One also needs to check that $\displaystyle \frac{\partial
  \mathcal{Z}_{r,s}(\tilde{f},\tilde{j})} {\partial
  \mathcal{Z}_{r,s}(0,0)}$ is well defined and nonzero.\\

Once this is verified then it immediately follows from the
considerations in \S\ref{algdefsec} that

\[
\mathcal{S}_{r,s}(\tilde{f},\tilde{j})
=
\exp
\left(
\mathcal{S}_{r,s}^{\rm T}(\tilde{f},\tilde{j})
\right)\ .
\]

The brunt of the remaining work is controlling the $r \rightarrow
-\infty$ and $s \rightarrow \infty$ limits of the log-moment
generating function $\mathcal{S}^{\rm T}_{r,s}(\tilde{f},\tilde{j})$.\\

Recall that for the denominator, i.e. when $\tilde{f},\tilde{j}=0$,
the initial condition for the $RG_{\rm ex}$ iterations is 

\begin{gather*}
\vec{V}^{(r,r)}(0,0) 
=
(g,0,\mu_{\rm c}(g),0,0,0,0,0)\\
\textnormal{  with  }
\mu_{\rm c}(g) = \mu_{\rm s}(g - \bar{g},0)
\textnormal{  by definition.  }
\end{gather*}

If $\iota$ is the affine isometric injection $\mathcal{E} \rightarrow
\mathcal{E}_{\rm ex}$ which sends $(\delta g, \mu, R)$ to the vector 

\[
\vec{V}
=
\left(
\beta_{4,\Delta},
\beta_{3,\Delta},
\beta_{2,\Delta},
\beta_{1,\Delta},
W_{5,\Delta},
W_{6,\Delta},
f_{\Delta},
R_{\Delta}
\right)_{\Delta\in\mathbb{L}}
\]

where for all $\Delta \in \mathbb{L}$

\begin{equation*}
\begin{split}
\beta_{4,\Delta}&=\bar{g}+\delta g\\
\beta_{3,\Delta}&=0\\
\beta_{2,\Delta}&=\mu\\
\beta_{1,\Delta}&=0\\
W_{5,\Delta}&=0\\
W_{6,\Delta}&=0\\
f_{\Delta}&=0\\
R_{\Delta}&=R
\end{split}
\end{equation*}

then $\vec{V}^{r,r}(0,0) = \iota(v)$ with $v = (\delta g,\mu_{\rm s}(\delta g,0),0)$ where $\delta g = g -
\bar{g}$.\\

By construction $v \in W^{\rm s,loc}$ and therefore all of its iterates
are well defined and we have

\[
\vec{V}^{(r,q)}(0,0)
=
\iota
\left(
RG^{q-r}
(v)
\right)
\longrightarrow
\iota(v_{*})
\textnormal{  where  }
r
\rightarrow
-\infty
\textnormal{  with  }
q
\textnormal{  fixed.  }
\]

The purpose of this section is to derive estimates which control the
deviations from this bulk trajectory due to the test functions
$\tilde{f}$ and $\tilde{j}$. We will break up the log-moment
generating function into five pieces which will be analyzed
separately.\\

Namely, we write

\begin{equation*}
\begin{split}
\mathcal{S}^{\rm T}_{r,s}(\tilde{f},\tilde{j})
=&
\mathcal{S}^{\rm T,FR}_{r,s}(\tilde{f},\tilde{j})
+
\mathcal{S}^{\rm T,UV}_{r,s}(\tilde{f},\tilde{j})
+
\mathcal{S}^{\rm T,MD}_{r,s}(\tilde{f},\tilde{j})\\
&+
\mathcal{S}^{\rm T,IR}_{r,s}(\tilde{f},\tilde{j})
+
\mathcal{S}^{\rm T,BD}_{r,s}(\tilde{f},\tilde{j})
\end{split}
\end{equation*}

where

\begin{equation*}
\begin{split}
\mathcal{S}^{\rm T,FR}_{r,s}(\tilde{f},\tilde{j})
=&
\frac{1}{2}
\sum_{r \le q < s}
\left(
f^{(r,q)},\Gamma f^{(r,q)}
\right)_{\Lambda_{s-q}}\\[1.5ex]
\mathcal{S}^{\rm T,UV}_{r,s}(\tilde{f},\tilde{j})
=&
-Y_{0}Z_{0}^{r}
\int_{\mathbb{Q}_{p}^{3}}
\tilde{j}(x)\ d^{3}x
+
\sum_{r \le q < q_{-}}
\sum_{\substack{\Delta \in \mathbb{L} \\ \Delta \subset
    \Lambda_{s-q-1}}}
\left(
\delta b_{\Delta}
[
\vec{V}^{(r,q)}(\tilde{f},\tilde{j})
]
-
\delta b_{\Delta}
[
\vec{V}^{(r,q)}(0,0)
]
\right)\\[1.5ex]
\mathcal{S}_{r,s}^{\rm T,MD}(\tilde{f},\tilde{j})
=&
\sum_{q_{-} \le q < q_{+}}
\sum_{\substack{\Delta \in \mathbb{L} \\ \Delta \subset
    \Lambda_{s-q-1}}}
\left(
\delta b_{\Delta}
[
\vec{V}^{(r,q)}(\tilde{f},\tilde{j})
]
-
\delta b_{\Delta}
[
\vec{V}^{(r,q)}(0,0)
]
\right)\\[1.5ex]
\mathcal{S}_{r,s}^{\rm T,IR}(\tilde{f},\tilde{j})
=&
\sum_{q_{+} \le q < s}
\sum_{\substack{\Delta \in \mathbb{L} \\ \Delta \subset
    \Lambda_{s-q-1}}}
\left(
\delta b_{\Delta}
[
\vec{V}^{(r,q)}(\tilde{f},\tilde{j})
]
-
\delta b_{\Delta}
[
\vec{V}^{(r,q)}(0,0)
]
\right)\\
\textnormal{and}\\
\mathcal{S}^{\rm T,BD}_{r,s}(\tilde{f},\tilde{j})
=&
{\rm Log}
\left(
\frac{\partial \mathcal{Z}_{r,s}(\tilde{f},\tilde{j})}{\partial
  \mathcal{Z}_{r,s}(0,0)}
\right).
\end{split}
\end{equation*}

The subscript ``${\rm FR}$'' stands for the free contribution. Indeed, an
easy exercise shows that

\[
\lim_{\substack{ r \rightarrow -\infty \\ s \rightarrow \infty}}
\mathcal{S}^{\rm T,FR}_{r,s}(\tilde{f},\tilde{j})
=
\frac{1}{2}
\left(
\tilde{f},
C_{-\infty} \tilde{f}
\right)
\]
which corresponds to the free massless measure without cut-offs, i.e.,
the Gaussian measure with covariance $C_{-\infty}$.\\

The quantity $\mathcal{S}^{\rm T,UV}_{r,s}(\tilde{f},\tilde{j})$ collects
the ultraviolet contributions while
$\mathcal{S}^{\rm T,IR}_{r,s}(\tilde{f},\tilde{j})$ contains the infrared
contributions. Most of the influence of the test functions is felt in
the middle regime $q_{-} \le q < q_{+}$, hence the abbreviation
``${\rm MD}$''. Finally $\mathcal{S}^{\rm T,BD}_{r,s}(\tilde{f},\tilde{j})$
corresponds to the a boundary term left after the RG iterations have
shrunk the confining volume $\Lambda$ down to a single unit cube.\\

The analysis will make use of the following observations with are of
an algebraic or combinatorial nature. Since the RG runs from UV
scales to IR scales we will first have a closer look at the terms
featuring in $\mathcal{S}^{\rm T,UV}_{r,s}(\tilde{f},\tilde{j})$.\\

From the definition of $RG_{\rm ex}$ in \S\ref{algdefsec} one sees that this map is given
by a collection of independent operations performed locally.\\

Indeed the output $\displaystyle \left(
  \beta'_{4,\Delta'},\dots,\beta'_{1,\Delta'},W'_{5,\Delta'},W'_{6,\Delta'},f'_{\Delta'},
  R_{\Delta'} \right)$ as well as the output $\displaystyle \delta
b_{\Delta'}$ produced for a cube $\Delta'$ only involves the data
$\displaystyle \left( \beta_{4,\Delta},\dots, \beta_{1,\Delta},
  W_{5,\Delta}, W_{6,\Delta}, f_{\Delta}, R_{\Delta} \right)_{\Delta
  \in [L^{-1}\Delta']}$.\\

In other words, $RG_{\rm ex}$ is made up of independent copies of a map
$(\mathcal{E}_{1{\rm B}})^{\times L^{3}} \longrightarrow
\mathcal{E}_{1{\rm B}}$.\\

Let $\widetilde{\Delta} \in \mathbb{L}_{q_{-}}$ so that $\tilde{f}$ and
$\tilde{j}$ are constant on $\widetilde{\Delta}$ taking the values
$\tilde{f}_{\widetilde{\Delta}}$ and $\tilde{j}_{\widetilde{\Delta}}$
respectively. If $\displaystyle \widetilde{\Delta} \not \in \Lambda_{q_{+}}$ then
$\tilde{f}_{\widetilde{\Delta}} = \tilde{j}_{\widetilde{\Delta}}=0$.\\

First let us see what happens for the first iteration, i.e., $q=r$.

If a unit cube $\Delta$ is in $\displaystyle \Lambda_{s-r} \setminus
\Lambda_{q_{+}-r}$ then the $\Delta$ component of
$\vec{V}^{(r,r)}(\tilde{f},\tilde{j})$ of
$\vec{V}^{(r,r)}(\tilde{f},\tilde{j})$ is exactly the same as that of
the bulk $\vec{V}^{(r,r)}(0,0) = \iota( \delta g, \mu, 0)$ with $\mu =
\mu_{\rm s}(\delta g, 0)$.\\

If $\displaystyle \Delta \in \Lambda_{q_{+}-r}$ then there is a unique
$\widetilde{\Delta} \in \mathbb{L}_{q_{-}}$, $\widetilde{\Delta} \subset
\Lambda_{q_{+}}$ such that $\Delta\subset L^r \widetilde{\Delta}$

In this case:

\[
\vec{V}^{(r,r)}(\tilde{f},\tilde{j})
=
(
g,
0,
\mu
-
Y_{2}Z_{2}^{r}
L^{(3-2[\phi])r}
\tilde{j}_{\widetilde{\Delta}},
0,
0,
0,
L^{(3-[\phi])r} \tilde{f}_{\widetilde{\Delta}},
0
).
\]

Now we choose $Z_{2}$ so that $Z_{2} = \alpha_{\rm u}L^{-(3-2[\phi])}$ and
thus

\[
\vec{V}^{(r,r)}_{\Delta}(\tilde{f},\tilde{j})
=
(
g,
0,
\mu-Y_{2}\alpha_{\rm u}^{r}\tilde{j}_{\widetilde{\Delta}},
0,
0,
0,
L^{(3-[\phi])r}
\tilde{f}_{\widetilde{\Delta}},
0
).
\]

If $q=r < q_{-}$ then all immediate neighbors $\Delta$ carry the same
data. Here by neighbors we mean the $L^{3}-1$ other unit cubes
contained in the same $L$-block $L^{-1}\Delta'$ as $\Delta$. \\

Therefore the computation producing $\displaystyle \delta b_{\Delta'}
[ \vec{V}^{(r,r)}(\tilde{f},\tilde{j})]$ as well as
$\vec{V}^{(r,r+1)}_{\Delta'}(\tilde{f},\tilde{j})$ is the same as the
RG acting on the space $\mathcal{E}_{\rm bk}$. In fact the computation
reduces to the map $RG$ on the even smaller subspace $\mathcal{E}$,
except for the presence of the $f$-component
$L^{(3-[\phi])r}\tilde{f}_{\widetilde{\Delta}}$.\\

The key observation is that this component evolves by
averaging without influencing or being influenced by the other
variables. \\

This again results from the property that $\displaystyle
\int_{L^{-1}\Delta'} \Gamma(x-y) {\rm d}^{3}y = 0$ for all $x \in
L^{-1}\Delta'$, as in the proof of Proposition \ref{L53prop}. \\

Indeed for the explicit diagrams in the RG transformation the possible
effect of $f$ is through legs attached to $f$-vertices of valence 1 which
precisely contribute a factor of the type $\displaystyle
\int_{L^{-1}\Delta'} \Gamma(x-y) {\rm d}^{3}y = 0$ because $f$ is constant
over the $L$-block $L^{-1}\Delta'$.\\

For the other $\mathcal{L}$ or $\xi$ terms, observe that one has $\displaystyle
e^{\int_{L^{-1}\Delta'} f \zeta} = 1$ because $f$ is constant on
$L^{-1}\Delta'$ and $\displaystyle \int_{L^{-1}\Delta'} \zeta = 0$
almost surely by the property of the fluctuation covariance $\Gamma$.\\

As a result

\[
\vec{V}^{(r,r+1)}_{\Delta'}(\tilde{f},\tilde{j})
=
(
g',
0,
\mu',
0,
0,
0,
L^{(3-[\phi])(r+1)}\tilde{f}_{\widetilde{\Delta}},
R'
)
\]

where

\[
(
g'-\bar{g},
\mu'
,R'
)
=
RG
(
g-\bar{g},
\mu
-
\alpha_{\rm u}^{r}Y_{2} \tilde{j}_{\widetilde{\Delta}},
0
)
\]

and also

\[
\delta b_{\Delta'} [ \vec{V}^{(r,r)}(\tilde{f},\tilde{j})]
=
\delta b
( 
g - \bar{g},
\mu
-
\alpha_{\rm u}^{r}Y_{2} \tilde{j}_{\widetilde{\Delta}},
0
).
\]\\

The same decoupling applies to subsequent iterates
$\vec{V}^{(r,q+1)}(\tilde{f},\tilde{j}) = RG_{\rm ex}[
\vec{V}^{(r,q)}(\tilde{f},\tilde{j})]$ as long as $q < q_{-}$, i.e. as
long as $f^{(r,q)}$ is constant over each individual $L$-block.\\

Hence, in the quantity

\[
\sum_{r \le q < q_{-}}
\sum_{\substack{\Delta \in \mathbb{L} \\ \Delta \subset
    \Lambda_{s-q-1}}}
\left(
\delta b_{\Delta}
[
\vec{V}^{(r,q)}(\tilde{f},\tilde{j})
]
-
\delta b_{\Delta}
[
\vec{V}^{(r,q)}(0,0)
]
\right)
\]
appearing in $\mathcal{S}^{\rm T,UV}_{r,s}(\tilde{f},\tilde{j})$, only boxes
$\displaystyle \Delta \subset \Lambda_{ q_{+}-r-1}$ will contribute and these can be organized
according to $\displaystyle \widetilde{\Delta} \in \mathbb{L}_{q_{-}}$,
$\displaystyle \widetilde{\Delta} \subset \Lambda_{q_{+}}$ such that
$L^{q+1}\widetilde{\Delta}$ contains $\Delta$. All $\displaystyle L^{3(q_{-}-q-1)}$
boxes $\Delta$ which satisfy that condition for given $\widetilde{\Delta}$
produce the same contribution.\\

In other words, the previous expression can be rewritten as

\begin{gather*}
\sum_{ \substack{\widetilde{\Delta}\in\mathbb{L}_{q_{-}} \\ \widetilde{\Delta}
      \subset \Lambda_{q_{+}}}}
\sum_{r \le q < q_{-}}
L^{3(q_{-}-q-1)}
\left(
\delta b
\left[
RG^{q-r}
\left(
v 
- 
\alpha_{\rm u}^{r} Y_{2}\tilde{j}_{\widetilde{\Delta}}
e_{\phi^{2}}
\right)
\right]
-
\delta b
\left[
RG^{q-r}
(
v
)
\right]
\right)\\
\textnormal{  where  }
v
=
(
\delta g,
\mu_{\rm s}(\delta g, 0),
0
)
\textnormal{  with  }
\delta g
=
g
-
\bar{g}\\
\textnormal{  and  }
e_{\phi^{2}}
=
(0,1,0)
\in
\mathcal{E}.
\end{gather*}

Here $e_{\phi^{2}}$ gives the direction of pure $: \phi^{2} :$
  perturbations in the bulk. \\

We are thus reduced to $\displaystyle L^{3(q_{+}-q_{-})}$ separate and
independent bulk RG trajectories as considered in \S\ref{bulksec}, one
for each $\widetilde{\Delta}$. Also note that the effect of $\tilde{f}$ is
completely absent form the UV regime contribution.\\

By also organizing the explicit extra linear term in $\tilde{j}$
according to boxes $\widetilde{\Delta}$ of size $\displaystyle L^{q_{-}}$
we can write

\[
\mathcal{S}^{\rm T,UV}_{r,s}(\tilde{f},\tilde{j})
=
\sum_{\substack{\widetilde{\Delta} \in \mathbb{L}_{q_{-}} \\
    \widetilde{\Delta} \subset \Lambda_{q_{+}}}}
\mathcal{K}_{\widetilde{\Delta}}
\]

with

\[
\mathcal{K}_{\widetilde{\Delta}}
=
-Y_{0}
Z_{0}^{r}
L^{3q_{-}} 
\tilde{j}_{\widetilde{\Delta}}
+
\sum_{r \le q < q_{-}}
L^{3(q_{-}-q-1)}
\left(
\delta b
\left[
RG^{q-r}
\left(
v 
- 
\alpha_{\rm u}^{r} Y_{2}\tilde{j}_{\widetilde{\Delta}}
e_{\phi^{2}}
\right)
\right]
-
\delta b
\left[
RG^{q-r}
(
v
)
\right]
\right)\ .
\]\\

We now look at the middle regime and note that

\[
\mathcal{S}_{r,s}^{\rm T,MD}(\tilde{f},\tilde{j})
=
\sum_{q_{-} \le q < q_{+}}
\sum_{\substack{\Delta \in \mathbb{L} \\ \Delta \subset
    \Lambda_{q_{+}-q-1}}}
\left(
\delta b_{\Delta}
[
\vec{V}^{(r,q)}(\tilde{f},\tilde{j})
]
-
\delta b_{\Delta}
[
\vec{V}^{(r,q)}(0,0)
]
\right).
\]

Here we replaced the $s$ that appeared earlier with $q_{+}$ when
describing the summation over boxes $\Delta$. Indeed, if $\Delta
\subset \Lambda_{s-q-1}$ is outside the rescaling $\displaystyle
\Lambda_{q_{+}-q - 1}$ of the set $\displaystyle \Lambda_{q_{+}}$
containing the supports of the $\tilde{f}$ and $\tilde{j}$, then the
effect of $\Delta$ is nil.\\

What we need here is a more precise description of the vector
$\vec{V}^{(r,q_{-})}(\tilde{f},\tilde{j})$ delivered by the RG
evolution in the UV regime. This involves a fusion of
$\displaystyle L^{3(q_{+}-q_{-})}$ data living in $\mathcal{E}$ into a
single vector in $\mathcal{E}_{\rm ex}$. \\

For $m \ge 0$ we introduce the reinjection map

\begin{gather*}
\mathcal{J}_{m}: 
S_{0,m}
(\mathbb{Q}_{p}^{3},\mathbb{C})
\times
\left( 
\prod_{
\substack{ \Delta \in \mathbb{L} \\ \Delta \subset \Lambda_{m}}}
\mathcal{E}
\right)
\longrightarrow
\mathcal{E}_{\rm ex}\\
\left(F, 
(\delta g_{\Delta},
\mu_{\Delta}, 
R_{\Delta})_{ 
\substack{\Delta \in \mathbb{L} \\ \Delta \subset \Lambda_{m}}},
(\delta g,
\mu,
R)
\right)
\mapsto
\vec{V}'
=
\left(
 \beta'_{4,\Delta},
\dots,
\beta'_{1,\Delta},
W'_{5,\Delta},
W'_{6,\Delta},
f'_{\Delta'},
R'_{\Delta'}
\right)_{\Delta \in \mathbb{L}}
\end{gather*}

defined as follows.\\

We let

\begin{equation*}
\begin{split}
\beta'_{4,\Delta} 
&=
\begin{cases}
\bar{g} + \delta g_{\Delta} & \textnormal{  if  } \Delta \subset
\Lambda_{m}\\
\bar{g} + \delta g  & \textnormal{  if  } \Delta \not \subset
\Lambda_{m}
\end{cases}\\[1.5ex]
\beta'_{2,\Delta}
&=
\begin{cases}
\mu_{\Delta} & \textnormal{  if  } \Delta \subset
\Lambda_{m}\\
\mu  & \textnormal{  if  } \Delta \not \subset
\Lambda_{m}
\end{cases}\\[1.5ex]
R'_{\Delta}
&=
\begin{cases}
R_{\Delta} & \textnormal{  if  } \Delta \subset
\Lambda_{m}\\
R  & \textnormal{  if  } \Delta \not \subset
\Lambda_{m}
\end{cases}\\[1.5ex]
\beta'_{3,\Delta}
=
\beta'_{1,\Delta}
&=
W'_{5,\Delta}
=
W'_{6,\Delta}
=0
\end{split}
\end{equation*}

and finally $f'_{\Delta}$ is defined by

\[
f'_{\Delta(x)} 
= 
F(x)
\textnormal{  for all }
x
\in
\mathbb{Q}_{p}^{3}.
\]

Namely, via the correspondance between $\mathbb{L}$-indexed vectors
and functions that are constant on unit cubes, $f'=F$. Recall indeed
that $F$ is assumed constant on unit cubes and with support contained
in $\Lambda_{m}$.\\

Now it is easy to see from the previous considerations that

\[
\vec{V}^{(r,q_{-})}(\tilde{f},\tilde{j})
=
\mathcal{J}_{q_{+}-q_{-}}
\left(
\tilde{f}_{\rightarrow (-q_{-})},
\left(
RG^{q_{-}-r}
\left(
v 
- 
\alpha_{\rm u}^{r}
Y_{2}
\tilde{j}_{L^{-q_{-}}\Delta}
e_{\phi^{2}}
\right)
\right)_{
\substack{\Delta \in \mathbb{L} \\ \Delta \subset \Lambda_{q_{+}-q_{-}}}}
,
RG^{q_{-}-r}(v)
\right).
\]

Note in particular that $\displaystyle f^{(r,q_{-})} = \left(
  \tilde{f}_{\rightarrow (-r)}
\right)_{\rightarrow(q_{-}-r)} = \tilde{f}_{\rightarrow (-q_{-})}$.\\

We also have the special case

\begin{equation*}
\begin{split}
\vec{V}^{(r,q_{-})}(0,0)
=&
\mathcal{J}_{q_{+}-q_{-}}
\left(
0,
\left(
RG^{q_{-}-r}(v)
\right)_{
\substack{\Delta \in \mathbb{L} \\ \Delta \subset
  \Lambda_{q_{+}-q_{-}}}},
RG^{q_{-}-r}(v)
\right)\\
=&
\iota
\left(
RG^{q_{-}-r}(v)
\right)
\in
\mathcal{E}_{\rm bk}\subset\mathcal{E}.
\end{split}
\end{equation*}

Finally, in the infrared regime

\begin{gather*}
\mathcal{S}^{\rm T,IR}_{r,s}(\tilde{f},\tilde{j})
=
\sum_{q_{+} \le q < s}
\sum_{ \substack{ \Delta \in \mathbb{L} \\ \Delta \in
    \Lambda_{s-q-q}}}
\left(
\delta b_{\Delta}
\left[
RG^{q-q_{+}}
\left( 
\vec{V}^{(r,q_{+})}(\tilde{f},\tilde{j})
\right)
\right]
-
\delta b_{\Delta}
\left[
RG^{q-q_{+}} 
\left(
 \vec{V}^{(r,q_{+})}(0,0)
\right)
\right]
\right)\\
\textnormal{  where  }
\vec{V}^{(r,q_{+})}(\tilde{f},\tilde{j})
=
RG^{q_{+}-q_{-}}
\left(
\vec{V}^{(r,q_{-})}(\tilde{f},\tilde{j})
\right)\ .
\end{gather*}

Since $\displaystyle \vec{V}^{(r,q_{-})}(\tilde{f},\tilde{j})$ agrees
with $\displaystyle \vec{V}^{(r,q_{-})}(0,0)$ on all unit cubes
$\displaystyle \Delta \not \subset \Lambda_{q_{+}-q_{-}}$, it is easy
to see that

\begin{gather*}
RG^{q_{+}-q_{-}}
\left(
\vec{V}^{(r,q_{-})}(\tilde{f},\tilde{j})
\right)
\textnormal{  agrees with  }
RG^{q_{+}-q_{-}}
\left(
\vec{V}^{(r,q_{-})}(0,0)
\right)\\ 
\textnormal{  on all unit cubes  }
\Delta \not \subset \Lambda_{0} = \Delta(0)
\textnormal{  the unit cube containing the origin}.
\end{gather*}

Thus

\[
\vec{V}^{(r,q_{+})}(\tilde{f},\tilde{j})
-
\vec{V}^{(r,q_{+})}(0,0)
\in
\mathcal{E}_{\rm pt}
\]

or

\[
\vec{V}^{(r,q_{+})}(\tilde{f},\tilde{j})
\in
\iota(\mathcal{E})
\oplus
\mathcal{E}_{\rm pt}
\subset
\mathcal{E}_{\rm bk}
\oplus
\mathcal{E}_{\rm pt}\ .
\]

This property remains true for the next iterates since the only
difference with the bulk now only happens in $\Delta(0)$.\\

Therefore no summation over $\Delta$ is needed in the formula for
$\mathcal{S}^{\rm T,IR}_{r,s}(\tilde{f},\tilde{j})$ which thus reduces to

\[
\mathcal{S}^{\rm T,IR}_{r,s}(\tilde{f},\tilde{j})
=
\sum_{q_{+} \le q < s}
\left(
\delta b_{\Delta(0)}
\left[
RG^{q-q_{+}}
\left( 
\vec{V}^{(r,q_{+})}(\tilde{f},\tilde{j})
\right)
\right]
-
\delta b_{\Delta(0)}
\left[
RG^{q-q_{+}} 
\left(
 \vec{V}^{(r,q_{+})}(0,0)
\right)
\right]
\right).
\]

After these prepatory steps we can now address the estimates needed in
order to take the $r \rightarrow - \infty$ and $s \rightarrow \infty$
limits.\\

\subsection{The ultraviolet regime}

We first need an analogue of Lemma \ref{DRG2lem} for the function $\delta b$.\\

\begin{Lemma}\label{d2bbound}
For $\epsilon$ small, for all $v$ such that $\displaystyle ||v|| <
\frac{1}{4}$ we have

\[
||
D^{2}_{v}
\delta b
||
=
\sup_{v',v'' \not = 0}
\frac{
\left|
D^{2}_{v}
\delta b
[v',v'']
\right|}
{
||v'||
\times
||v''||
}
\le
2.
\]
\end{Lemma}

\noindent{\bf Proof:}
Recall that

\[
\delta b ( \delta g, \mu, R )
=
\delta b^{\rm explicit}( \delta g , \mu, R) 
+ 
\delta b^{\rm implicit} (\delta g, \mu, R)
\]

where

\[
\delta b^{\rm explicit}( \delta g , \mu, R) 
=
A_{4}
(
\bar{g}
+
\delta g
)^{2}
+
A_{5}
\mu^{2}
\]

and

\[
\delta b^{\rm implicit} (\delta g, \mu, R)
=
\xi_{0}(\bar{g}+\delta g,\mu,R)\ .
\]

using the same notations as in Lemma \ref{DRG2lem} we have

\[
D^{2}_{v} \delta b^{\rm explicit} [ v', v'']
=
2
A_{4}
\delta g'
\delta g''
+
2
A_{5}
\mu'
\mu''.
\]

Now, using $C_0(0),||\Gamma||_{L^{\infty}}\le 2$, we have

\begin{equation*}
\begin{split}
|A_{4}| 
\le&
12
L^{3}
||\Gamma||_{L^{\infty}}^{2}
||\Gamma||_{L^{2}}^{2}
+
48
L^{\frac{3+\epsilon}{2}}
\times
2
\times
||\Gamma||_{L^{\infty}}
||\Gamma||_{L^{2}}^{2}
+
72
L^{\epsilon}
\times
4
\times
||\Gamma||_{L^{2}}^{2}\\
\le&
\left[
48L^{3}
+
192
L^{\frac{3+\epsilon}{2}}
+
288
L^{\epsilon}
\right]
||\Gamma||_{L^{2}}^{2}\\
= &
\left(
48L^{3}
+
192
L^{\frac{3+\epsilon}{2}}
+
288
L^{\epsilon}
\right)
\frac{1}{36}
L^{-\epsilon}
A_{1}
\le
A_{4,{\rm max}}
\end{split}
\end{equation*}

with

\[
A_{4,{\rm max}}
=
\left(
48L^{3}
+
192
L^{\frac{3}{2}}
+
288
\right)
\times
\frac{1}{36}
A_{1,\max}\ .
\]

Likewise $|A_{5}| \le A_{5,{\rm max}}$ with

\[
A_{5,{\rm max}} 
= 
L^{3} 
\times 
\frac{1}{36} A_{1,{\rm max}}.
\]

Thus

\[
\left|
D^{2}_{v}
\delta b^{\rm explicit}
[v',v'']
\right|
\le
2
||v'||
\times
||v''||
\left[
2
A_{4,{\rm max}}
\bar{g}^{2e_{4}}
+
2
A_{5,{\rm max}}
\bar{g}^{2e_{2}}
\right].
\]

Now if $\displaystyle ||v|| < \frac{1}{4}$ then we can use Cauchy's
formula

\[
D^{2}_{v}
\delta b^{\rm implicit}[v',v'']
=
\frac{1}{(2 i \pi)^{2}}
\oint
\frac{d \lambda_{1}}{\lambda_{1}^{2}}
\oint
\frac{d \lambda_{2}}{\lambda_{2}^{2}}\ 
\delta b^{\rm implicit}[v+\lambda_{1} v' + \lambda_{2} v'']
\]

where the contours are given by $\displaystyle
|\lambda_{1}|=\frac{1}{8||v'||}$, $\displaystyle
|\lambda_{2}|=\frac{1}{8||v''||}$. Then we get from Theorem \ref{mainestthm}

\begin{equation*}
\begin{split}
\left|
D^{2}_{v}
\delta b^{\rm implicit}[v',v'']
\right|
\le&
8
||v'||
\times
8
||v''||
\times
B_{0} \bar{g}^{e_{R}}
\times
\sup_{\lambda_{1},\lambda_{2}}
||
v
+
\lambda_{1}v'
+
\lambda_{2}v''
||\\
\le&
64
||v'||
\times
||v''||
\times
B_{0}
\bar{g}^{e_{R}}
\times
\frac{1}{2}\ .
\end{split}
\end{equation*}

Combining both bounds we obtain

\[
\left|
D^{2}_{v}
\delta b[v',v'']
\right|
\le
||v'||
\times
||v''||
\left[
4
A_{4,{\rm max}}
\bar{g}^{2e_{4}}
+
4
A_{5,{\rm max}}
\bar{g}^{2e_{2}}
+
32
B_{0}
\bar{g}^{e_{R}}
\right].
\]

Since $e_{4},e_{2},e_{R} > 0$ the lemma follows by making $\bar{g}$, i.e., $\epsilon$
small enough.
\qed\\

\begin{Lemma}\label{ddblem}
For $\epsilon$ small and for all $v$ such that $\displaystyle ||v|| <
\frac{1}{4}$

\[
||
D_{v} \delta b
||=\sup_{v' \not = 0} \frac{ \left|D_{v} \delta b[v']
  \right|}{||v'||}
\le
1.
\]
\end{Lemma}

\noindent{\bf Proof:}
Now for $\displaystyle ||v|| \le \frac{1}{2}$

\[
D_{v}
\delta b^{\rm explicit}
[v']
=
2 A_{4}
(\bar{g}+\delta g)
\delta g'
+
2
A_{5}
\mu
\mu'
\]

\[
\left|
D_{v}
\delta b^{\rm explicit}[v']
\right|
\le
2 A_{4,{\rm max}}
\times
\frac{3}{2} \bar{g}
\times
\bar{g}^{e_{4}}
||v'||
+
2 A_{5,{\rm max}}
\times
\frac{1}{2}
\bar{g}^{e_{2}}
\times
\bar{g}^{e_{2}}
||v'||.
\]

If furthermore $\displaystyle ||v|| < \frac{1}{4}$ then we can  write Cauchy's formula

\begin{gather*}
D_{v}
\delta b^{\rm implicit}
[v']
=
\frac{1}{2i \pi }
\oint
\frac{d \lambda}
{\lambda^{2}}
\ 
\delta b^{\rm implicit}(v+\lambda v')\\
\textnormal{  on the contour  }
|\lambda| 
= 
\frac{1}{4 ||v'||}
\end{gather*}

and deduce

\[
\left|
D_{v}
\delta b^{\rm implicit}
[v']
\right|
\le
4
||v'||
\times
B_{0}
\bar{g}^{e_{R}}
\times
\frac{1}{2}\ .
\]

Thus for $\displaystyle ||v|| < \frac{1}{4}$, 

\[
||
D_{v} \delta b
||
\le
3
A_{4,{\rm max}}
\bar{g}^{e_{4}+1}
+
A_{5,{\rm max}} 
\bar{g}^{2e_{2}}
+
2 B_{0} \bar{g}^{e_{R}}.
\]
Again the last expression can be made as small as we want provided $\epsilon$ is small enough.
\qed\\

Using the mean value theorem or Taylor's formula with integral remainder we
immediately obtain as before the following lemma.

\begin{Lemma}\label{dbbound}
For $\epsilon$ small

\begin{enumerate}

\item For all $v$ with $\displaystyle ||v|| < \frac{1}{4}$

\[
||
D_{v}
\delta b 
- 
D_{v_{*}}
\delta b
|| 
\le 
2 
||v-v_{*}||\ .
\]

\item For all $v,w$ such that $\displaystyle ||v||, ||w|| <
  \frac{1}{4}$

\[
||
\delta b(v+w) 
- 
\delta b(v)
||
\le
||w||
\]
and
\[
||
\delta b(v+w)
-
\delta b(v)
-
D_{v}
\delta b[w]
||
\le
||w||^{2}.
\]

\end{enumerate}

\end{Lemma}

\qed\\

We now resume the analysis of the expression for $\mathcal{S}_{r,s}^{T,UV}(\tilde{f},\tilde{j})$
derived in the last section.

Adding and subtracting terms linear in
$\tilde{j}_{\widetilde{\Delta}}$ we write

\begin{gather*}
\mathcal{K}_{\widetilde{\Delta}}
=
\tilde{j}_{\widetilde{\Delta}}
\left[
-Y_{0}Z_{0}^{r}L^{3q_{-}}
+
\sum_{r \le q < q_{-}}
L^{3(q_{-}-q-1)}
D_{v}
\left(
\delta b 
\circ
RG^{q-r}
\right)
\left[
-\alpha_{\rm u}^{r} Y_{2} e_{\phi^{2}}
\right]
\right]
+
\sum_{r \le q < q_{-}}  L^{3(q_{-}-q-1)}
\mathcal{K}_{\widetilde{\Delta},q}\\
\textnormal{  where  }
\mathcal{K}_{\widetilde{\Delta},q}
=
\delta b 
\left[
RG^{q-r}
\left(
v 
- 
\alpha_{\rm u}^{r}Y_{2}\tilde{j}_{\widetilde{\Delta}} e_{\phi^{2}}
\right)
\right]
-
\delta
b
\left[
RG^{q-r}
(v)
\right]
+
D_{v}
\left(
\delta b
\circ
RG^{q-r}
\right)
\left[
\alpha_{\rm u}^{r}Y_{2}\tilde{j}_{\widetilde{\Delta}} e_{\phi^{2}}
\right]\ .
\end{gather*}\\

Now $\mathcal{K}_{\widetilde{\Delta},q} =
\mathcal{K}'_{\widetilde{\Delta},q} +
\mathcal{K}''_{\widetilde{\Delta},q}$ where

\[
\mathcal{K}'_{\widetilde{\Delta},q}
=
\delta b 
\left[
RG^{q-r}
\left(
v 
- 
\alpha_{\rm u}^{r}Y_{2}\tilde{j}_{\widetilde{\Delta}} e_{\phi^{2}}
\right)
\right]
-
\delta
b
\left[
RG^{q-r}
(v)
\right]
-
D_{RG^{q-r}(v)}
\delta b
\left[
RG^{q-r}
\left(
v 
- 
\alpha_{\rm u}^{r}Y_{2}\tilde{j}_{\widetilde{\Delta}} e_{\phi^{2}}
\right)
-
RG^{q-r}
(v)
\right]
\]
and
\[
\mathcal{K}''_{\widetilde{\Delta},q}
=
D_{RG^{q-r}(v)}
\delta b
\left[
RG^{q-r}
\left(
v 
- 
\alpha_{\rm u}^{r}Y_{2}\tilde{j}_{\widetilde{\Delta}} e_{\phi^{2}}
\right)
-
RG^{q-r}
(v)
+
D_{v}
RG^{q-r} 
\left[
\alpha_{\rm u}^{r}Y_{2}\tilde{j}_{\widetilde{\Delta}} e_{\phi^{2}}
\right]
\right]\ .
\]

We already know $\displaystyle ||RG^{q-r}(v)|| <\frac{1}{4}$.

If the same is true for $\displaystyle RG^{q-r} \left( v -
  \alpha_{\rm u}^{r}Y_{2}\tilde{j}_{\widetilde{\Delta}} e_{\phi^{2}} \right)$
then Lemmas \ref{dbbound} and \ref{ddblem} imply

\[
||
\mathcal{K}'_{\widetilde{\Delta},q}
||\le
||
RG^{q-r}
\left(
v 
- 
\alpha_{\rm u}^{r}Y_{2}\tilde{j}_{\widetilde{\Delta}} e_{\phi^{2}}
\right)
-
RG^{q-r}(v)
||^{2}
\]

and

\[
||
\mathcal{K}''_{\widetilde{\Delta},q}
||
\le
||
RG^{q-r}
\left(
v 
- 
\alpha_{\rm u}^{r}Y_{2}\tilde{j}_{\widetilde{\Delta}} e_{\phi^{2}}
\right)
-
RG^{q-r}(v)
+
D_{v}
RG^{q-r}
\left[
\alpha_{\rm u}^{r}Y_{2}\tilde{j}_{\widetilde{\Delta}} e_{\phi^{2}}
\right]
||.
\]\\

We assume $\displaystyle || \alpha_{\rm u}^{q_{-}-1} Y_{2}
\tilde{j}_{\widetilde{\Delta}} e_{\phi^{2}} || \le \frac{1}{240
  \mathcal{C}_{1}(\epsilon)}$ which implies $ \displaystyle ||\alpha_{\rm u}^{q} Y_{2}
  \tilde{j}_{\widetilde{\Delta}} e_{\phi^{2}}|| \le \frac{1}{240
    \mathcal{C}_{1}(\epsilon)}$ for all $q < q_{-}$.\\

Lemma \ref{RGaRGblem} guarantees that 

\[
RG^{q-r}
\left(
v 
- 
\alpha_{\rm u}^{r}
Y_{2}
\tilde{j}_{\widetilde{\Delta}}
e_{\phi^{2}}
\right)
=
RG^{q-r}
\left(
v
+
\alpha_{\rm u}^{-(q-r)}
\left(
-\alpha_{\rm u}^{q}
Y_{2} \tilde{j}_{\widetilde{\Delta}} e_{\phi^{2}}
\right)
\right)
\]

is well defined and has norm at most $\displaystyle \frac{1}{8}$. 

The telescoping sum argument at the beginning of the proof of lemma \ref{unnamedlemma1} with
$n=q-1$, $k=0$, and $\displaystyle w = -
\alpha_{\rm u}^{q}Y_{2}\tilde{j}_{\widetilde{\Delta}}e_{\phi^{2}}$ gives

\begin{equation*}
\begin{split}
&||
RG^{q-r}
\left(
v
-
\alpha_{\rm u}^{r}
Y_{2}
\tilde{j}_{\widetilde{\Delta}}
e_{\phi^{2}}
\right)
-
RG^{q-r}
(v)
+
D_{v}
RG^{q-r}
\left[
\alpha_{\rm u}^{r}
Y_{2}
\tilde{j}_{\widetilde{\Delta}}
e_{\phi^{2}}
\right]
||\\
= &
||
RG^{q-r}
\left(
v
+
\alpha_{\rm u}^{-(q-r)}
w
\right)
-
RG^{q-r}
(v)
-
D_{v}
RG^{q-r}
\left[
\alpha_{\rm u}^{-(q-r)}
w
\right]
||\\
\le&
\sum_{i=0}^{q-r}
850
\mathcal{C}_{1}(\epsilon)^{2}
||w||^{2}
\alpha_{\rm u}^{-2}
\left(  c_{3}(\epsilon) \alpha_{\rm u}^{-2} \right)^{i}\\
\le&
850
\mathcal{C}_{1}(\epsilon)^{2}
||w||^{2}
\frac{1}{\alpha_{\rm u}^{2}-c_{3}(\epsilon)}\\
\le&
\frac{1700}{L^{3}} 
\mathcal{C}_{1}(\epsilon)^{2}
||w||^{2}\\
 = &
1700
L^{-3}
\mathcal{C}_{1}(\epsilon)^{2}
\alpha_{\rm u}^{-2(q_{-}-q-1)}
||
\alpha_{\rm u}^{q_{-}-1}
Y_{2}
\tilde{j}_{\widetilde{\Delta}}
e_{\phi^{2}}
||^{2}\\
\le&
\frac{1700}{240^{2} \times 8}
\alpha_{\rm u}^{-2(q_{-}-q-1)}
<
\alpha_{\rm u}^{-2(q_{-}-q-1)}
\end{split}
\end{equation*}
where we used the bound in (\ref{auc3eq}) to go from the fourth
to the fifth line as well as $L\ge 2$ in the last line.

On the other hand,

\[
\begin{split}
RG^{q-r}
\left(
v
-
\alpha_{\rm u}^{r}
Y_{2}
\tilde{j}_{\widetilde{\Delta}}
e_{\phi^{2}}
\right)
-
RG^{q-r}
(v)
= &
RG^{q-r}
\left(
v
-
\alpha_{\rm u}^{r}
Y_{2}
\tilde{j}_{\widetilde{\Delta}}
e_{\phi^{2}}
\right)
-
RG^{q-r}
(v)
+
D_{v}
RG^{q-r}
\left[
\alpha_{\rm u}^{r}
Y_{2}
\tilde{j}_{\widetilde{\Delta}}
e_{\phi^{2}}
\right]\\
 & +
T_{q-r}(v)[w]\ .
\end{split}
\]

So by the previous bound, Lemma \ref{Tbdlem} and Lemma \ref{eqnormlem} we obtain

\begin{equation*}
\begin{split}
||
RG^{q-r}
\left(
v
-
\alpha_{\rm u}^{r}
Y_{2}
\tilde{j}_{\widetilde{\Delta}}
e_{\phi^{2}}
\right)
-
RG^{q-r}
(v)
||
\le&
\alpha_{\rm u}^{-2(q_{-}-q-1)} 
+ 
10 
\mathcal{C}_{1}(\epsilon)
||w||\\
\le&
\alpha_{u}^{-2(q_{-}-q-1)} 
+
10
\mathcal{C}_{1}(\epsilon)
\times
\alpha_{\rm u}^{-(q_{-}-q-1)}
\frac{1}{240 \mathcal{C}_{1}(\epsilon)}\\
\le&
\frac{25}{24} \alpha_{\rm u}^{-(q_{-}-q-1)}\ .
\end{split}
\end{equation*}

Hence $ \displaystyle ||\mathcal{K}''_{\widetilde{\Delta},q}|| \le
\alpha_{\rm u}^{-2(q_{-}-q-1)}$ and $\displaystyle
||\mathcal{K}'_{\widetilde{\Delta},q}|| \le \left( \frac{25}{24} \right)^{2}
\alpha_{\rm u}^{-2(q_{-} - q -1)}$. With these two bounds in hand we can write the
estimate $\displaystyle ||\mathcal{K}_{\widetilde{\Delta},q}|| \le 3
\alpha_{\rm u}^{-2(q_{-}-q-1)}$ for simplicity.\\

$Y_{2}$ is a strictly positive quantity that will be fixed later and
we have that $ \displaystyle ||e_{\phi^{2}}|| = || (0,1,0) || =
\bar{g}^{-e_{2}}$. So the previous construction and bounds work if

\begin{equation}\label{jbound}
||\tilde{j}||_{L^{\infty}}
\le
\left[ 
240 
\mathcal{C}_{1}(\epsilon)
\alpha_{\rm u}^{q_{-}-1}
Y_{2}
\bar{g}^{-e_{2}}
\right]^{-1}\ .
\end{equation}\\

We will later also show $L^{3} \alpha_{\rm u}^{-2} < 1$ which will imply
that $\displaystyle \sum_{r \le q <q_{-}} L^{3(q_{-}-q-1)}
||\mathcal{K}_{\widetilde{\Delta},q}||$ is summable with uniform
bounds with respect to the UV cut-off $r$. \\

We now analyze the quantity

\[
\Omega_{r}
=
-Y_{0}
Z_{0}^{r}
L^{3q_{-}}
+
\sum_{r \le q < q_{-}}
L^{3(q_{-}-q-1)}
D_{v}
\left(
\delta b
\circ
RG^{q-r}
\right)
\left[
-\alpha_{\rm u}^{r}
Y_{2}
e_{\phi^{2}}
\right].
\]\\

We change the summation index to $n=q-r$ and rewrite the differential
using the chain rule and get

\begin{equation*}
\begin{split}
\Omega_{r}
=&
L^{3q_{-}}
\left(
-Y_{0}
Z_{0}^{r}
-
Y_{2}
\sum_{n=0}^{q_{-}-r-1}
L^{-3(n+r+1)}
\alpha_{\rm u}^{r}
D_{RG^{n}(v)}
\delta b 
\left[
D_{v}
RG^{n}
[
e_{\phi^{2}}
]
\right]
\right)\\[1.5ex]
=&
L^{3q_{-}}
\left(
-Y_{0}
Z_{0}^{r}
-
Y_{2}
\sum_{n=0}^{q_{-}-r-1}
L^{-3(n+r+1)}
\alpha_{\rm u}^{r+n}
D_{RG^{n}(v)}
\delta b 
\left[
T_{n}(v)[e_{\phi^{2}}]
\right]
\right)\\[1.5ex]
=&
L^{3q_{-}}
\left(
-Y_{0}
Z_{0}^{r}
-
Y_{2}
L^{-3}
(L^{-3}\alpha_{\rm u})^{r}
\sum_{n=0}^{q_{-}-r-1}
(L^{-3}\alpha_{\rm u})^{n}
\Xi_{n}
\right)\\[2ex]
&\textnormal{      with  } \Xi_{n}=D_{RG^{n}(v)}\delta b
\left[
T_{n}(v)[e_{\phi^{2}}]
\right]\ .
\end{split}
\end{equation*}

Note that from Lemma \ref{alphabdlem} we have $L^{-3}\alpha_{u} < 1$. But from Lemmas \ref{Tbdlem}, \ref{eqnormlem}
and \ref{ddblem} we have

\begin{equation*}
\begin{split}
\left|
D_{RG^{n}} 
\delta b 
\left[
T_{n}(v)
[e_{\phi^{2}}]
\right]
\right|
\le &
||
D_{RG^{n}(v)}
\delta b||
\times
10
\times
||
T_{n}(v) 
||_{\substack{\ \\ \Diamond}}
\times
||e_{\phi^{2}}||\\
\le&
10
\mathcal{C}_{1}(\epsilon)
\bar{g}^{-e_{2}}\ .
\end{split}
\end{equation*}

We then see that $\Xi_{n}$ is bounded uniformly with respect to $n$. Hence

\[ 
\Upsilon 
= 
\sum_{n=0}^{\infty} 
(L^{-3}\alpha_{\rm u})^{n}
\ \Xi_{n} 
\textnormal{  converges  }
\]

and we can  write

\[
\Omega_{r}
=
L^{3q_{-}}
\left(
-Y_{0}
Z_{0}^{r}
-
Y_{2}
L^{-3}
(L^{-3} \alpha_{\rm u})^{r}
\Upsilon
+
Y_{2}
L^{-3}
(L^{-3}\alpha_{\rm u})^{r}
\sum_{n = q_{-} - r}^{\infty}
(L^{-3}\alpha_{\rm u})^{n}
\ \Xi_{n}
\right)\ .
\]

Since $L^{-3}\alpha_{\rm u} < 1$ and $r \rightarrow -\infty$ we choose
$Y_{0}$, $Y_{2}$, and $Z_{0}$ so that the dangerous first two terms
cancel.\\

Namely, we set:

\begin{equation*}
\begin{split}
Z_{0} =& L^{-3}\alpha_{\rm u}\ ,\\
Y_{0} =& -L^{-3}Y_{2} \Upsilon\ .\\
\end{split}
\end{equation*}

Then

\begin{equation*}
\begin{split}
\Omega_{r}
=&
L^{3q_{-}}
Y_{2}
L^{-3}
(L^{-3}\alpha_{\rm u})^{r}
\sum_{n=q_{-}-r}^{\infty}
(L^{-3}\alpha_{\rm u})^{n}
\ \Xi_{n}\\
=&
Y_{2}
L^{-3}
\alpha_{\rm u}^{q_{-}}
\sum_{k=0}^{\infty}
(L^{-3}\alpha_{\rm u})^{k}
\ \Xi_{k+q_{-}-r}
\end{split}
\end{equation*}
after changing the summation index to $k=n-q_{-}+r$. \\

Provided one shows that $\displaystyle \lim_{n \rightarrow \infty}
\Xi_{n} = \Xi_{\infty}$ exists, the discrete dominated convergence
theorem will immediately imply

\[
\lim_{r \rightarrow -\infty}
\Omega_{r}
=
\frac{ Y_{2} L^{-3}\alpha_{\rm u}^{q_{-}} \Xi_{\infty}}
{1 - L^{-3}\alpha_{\rm u}}\ .
\]\\

Now

\begin{equation*}
\begin{split}
\left|
\Xi_{n}
-
D_{v_{*}}
\delta b
\left[
T_{\infty}
(v)
[e_{\phi^{2}}]
\right]
\right|
\le&
\left|
D_{RG^{n}(v)} 
\delta b 
\left[ 
T_{n}(v) 
[e_{\phi^{2}}]
\right]
-
D_{v_{*}}
\delta b
\left[
T_{n}
(v)
[e_{\phi^{2}}]
\right]
\right|
+
\left|
D_{v_{*}} \delta b
\left[
T_{n}(v)
[e_{\phi^{2}}]
- 
T_{\infty}(v)
[e_{\phi^{2}}]
\right]
\right|\\
\le&
2
||RG^{n}(v) - v_{*}||
\times
10
\mathcal{C}_{1}(\epsilon)
||e_{\phi^{2}}||
+
||T_{n}(v)-T_{\infty}(v)||
\times
||e_{\phi^{2}}||\ .
\end{split}
\end{equation*}

Above we used Lemmas \ref{dbbound}, \ref{ddblem} and \ref{Tbdlem}. 
Finally, Proposition \ref{L60prop} and Lemma \ref{Tcvlem} ensure that the
limit of the $\Xi_n$ exists and is given by $\Xi_{\infty}$.

As a consequence of the previous considerations and Theorem \ref{unnamedtheorem1}
we see that
\[
\lim_{
\substack{r \rightarrow -\infty \\ s \rightarrow \infty}}
\mathcal{S}^{\rm T,UV}_{r,s}(\tilde{f},\tilde{j})
=
\mathcal{S}^{\rm T,UV}(\tilde{f}, \tilde{j})
\textnormal{  with  }
\]

\begin{equation*}
\begin{split}
\mathcal{S}^{\rm T,UV}(\tilde{f}, \tilde{j})
=&
\sum_{
\substack{\widetilde{\Delta} \in \mathbb{L}_{q_{-}} \\ 
\widetilde{\Delta} \subset \Lambda_{q_{+}}}}
\bigg\{
\tilde{j}_{\widetilde{\Delta}}
\frac{Y_{2} \alpha_{\rm u}^{q_{-}}}
{L^{3}-\alpha_{\rm u}}
D_{v_{*}}\delta b 
\left[ 
T_{\infty}(v)[e_{\phi^{2}}]
\right]\\
&+
\sum_{q < q_{-}}
L^{3(q_{-}-q - 1)}
\left(
\delta b
\left(
\Psi
(v,
-\alpha_{\rm u}^{q}Y_{2}\tilde{j}_{\widetilde{\Delta}}e_{\phi^{2}}
)
\right)
-
\delta b (v_{*})
+
\alpha_{\rm u}^{q}
Y_{2}
\tilde{j}_{\widetilde{\Delta}}
D_{v_{*}}
\delta b [ T_{\infty}(v) (e_{\phi^{2}}) ]
\right)
\bigg\}\ .
\end{split}
\end{equation*}

The latter is easily seen to be analytic in $\tilde{j}$ in the domain
$\displaystyle ||\tilde{j}||_{L^{\infty}}
<
\left[ 
240 
\mathcal{C}_{1}(\epsilon)
\alpha_{\rm u}^{q_{-}-1}
Y_{2}
\bar{g}^{-e_{2}}
\right]^{-1}$ of $S_{q_{-},q_{+}}(\mathbb{Q}_{p}^{3},\mathbb{C})$.\\

Note that there is no dependence on $\tilde{f}$ for this piece. In
fact the finite cut-off quantity
$\mathcal{S}^{\rm T,UV}_{r,s}(\tilde{f},\tilde{j})$ does not depend on
$\tilde{f}$ nor $s$.

\subsection{The middle regime}

From here onwards we make additional requirements on the
exponents defining our norms:

\begin{equation}
\label{econstbetaeq}
1 - \eta \le e_{1} \le e_{2} \le e_{3} \le e_{4} < 2 - 2\eta
\end{equation}

\begin{equation}
\label{econst12eq}
e_{1}, e_{2} \le 1
\end{equation}

\begin{equation}
\label{econstWeq}
2 - 2\eta \le e_{W} \le \min( e_{3},1) + e_{4}
\end{equation}

We also introduce the notation $\bar{V}$ for the approximate fixed point in
$\mathcal{E}_{\rm bk}$. Namely we set $\bar{V}_{\Delta} = (\bar{g},0,\dots,0)$
for all $\Delta \in \mathbb{L}$. We note that $RG_{\rm ex}$ is well
defined and analytic on $B( \bar{V}, \frac{1}{2})$. 

We will next establish some very coarse bounds on the expansion of deviations
which will be enough for the control of the middle regime.
The next lemmas all assume that one is in the small $\epsilon$ regime.
The first one is a refinement of Lemma \ref{L33lem}.

\begin{Lemma}
\label{LAJ1}
Suppose that $\vec{V} \in B(\bar{V},\frac{1}{2})$. Then for
$k=1,2,3,4$ one
has the following bound for all $\Delta' \in \mathbb{L}$

\[
\left| \delta \beta_{k,1,\Delta'}[\vec{V}] \right|
\bar{g}^{-e_{k}}
\le
\bbone\{ 1 \le k < 4\} 
\mathbf{O}_1
L^{\frac{5}{2}},
\]

where $\mathbf{O}_1=\frac{27}{4}$ .

\end{Lemma}

\noindent{\bf Proof:}
From the definition we get

\[
\left|\delta \beta_{k,1,\Delta'}[\vec{V}] \right| 
\le 
\sum_{b}
\bbone\left\{ 
\begin{array}{c}
k+b\le 4 \\
b\ge 1
\end{array}
\right\}
\frac{(k+b)!}{k!\ b!}\ L^{-k[\phi]}\
\left| 
\parbox{2.1cm}{
\psfrag{a}{$\beta_{k+b}$}\psfrag{f}{$f$}\psfrag{b}{$b$}
\raisebox{-1ex}{
\includegraphics[width=2.1cm]{Fig5.eps}}
}
\ \ \right|
\]

The fact that $\delta \beta_{k,1,\Delta'}[\vec{V}]$ vanishes for $k=4$
is immediate.\\

We bound the Feynman Diagrams appearing in the formula above:

\begin{align*}
\bbone\left\{ 
\begin{array}{c}
k+b\le 4 \\
b\ge 1
\end{array}
\right\}
\left|
\parbox{2.1cm}{
\psfrag{a}{$\beta_{k+b}$}\psfrag{f}{$f$}\psfrag{b}{$b$}
\raisebox{-1ex}{
\includegraphics[width=2.1cm]{Fig5.eps}}
}
\ \ \right|
&\le
\bbone\left\{ 
\begin{array}{c}
k+b\le 4 \\
b\ge 1
\end{array}
\right\}
||f|_{L^{-1}\Delta'}||_{L^{\infty}}^{b}
\times
||\Gamma||_{L^{1}}^{b}
\times 
L^{3}
\times
\max_{\Delta \in [L^{-1}\Delta']}
|\beta_{k+b,\Delta}|\\
&\le
\bbone\left\{ 
b\ge 1
\right\}
L^{3}
\left( \frac{1}{\sqrt{2}} L^{3-2[\phi]}\times\frac{1}{2} L^{-(3-[\phi])} \right)^{b}
\max_{k+1 \le j \le 4} \left( \max_{\Delta \in [L^{-1}\Delta']} |\beta_{j,\Delta}|\right)\\
&\le
\frac{3}{4} L^{\frac{5}{2}} \bar{g}^{\min(e_{k+1},1)}
\end{align*}

In the last line we used $b \ge 1$ and $\epsilon \le 1$ so $3 - b[\phi]
\le \frac{5}{2}$. We also dropped the factor of
$\frac{1}{\sqrt{2}}$. We used \eqref{econstbetaeq} to bound
$\displaystyle \max_{k+1
  \le j \le 4} \left( \max_{\Delta \in [L^{-1}\Delta']}
  |\beta_{j,\Delta}|\right)$ by $\displaystyle \max\left(\frac{3}{2}\bar{g},
\frac{1}{2}\bar{g}^{e_{k+1}}\right)$.\\

We use this bound on Feynman diagrams to get the following bound
valid for $k=1$ and $2$:

\begin{align*}
\left| \delta \beta_{k,1,\Delta'}[\vec{V}] \right|
\bar{g}^{-e_{k}}
\le &
\bar{g}^{-e_{k}}
\sum_{b}
\bbone\left\{ 
\begin{array}{c}
k+b\le 4 \\
b\ge 1
\end{array}
\right\}
\frac{(k+b)!}{k!\ b!}\ L^{-k[\phi]}\
\left(\frac{3}{4} L^{\frac{5}{2}} \bar{g}^{\min(e_{k+1},1)}\right)\\
\le &
\frac{3}{4} 
L^{\frac{5}{2}}
\sum_{b}
\bbone\left\{ 
\begin{array}{c}
k+b\le 4 \\
b\ge 1
\end{array}
\right\}\frac{(k+b)!}{k!\ b!}
\end{align*}

In going to the second line we used that for $k=1,2$ one has
$\min(e_{k+1},1) \ge e_{k}$, this is a consequence of
\eqref{econstbetaeq} and \eqref{econst12eq}. We also dropped the
factors of $L^{-k[\phi]}$. \\

For $k=3$ which forces $b=1$ we only have one diagram to estimate:

\begin{align*}
\left|
\parbox{1cm}{
\psfrag{w}{$\beta_4$}\psfrag{f}{$f$}
\raisebox{1ex}{
\includegraphics[width=1cm]{Fig2.eps}}
}
\ \ \right|
\le &
\left|
\parbox{1cm}{
\psfrag{w}{$(\beta_4-\bar{g})$}\psfrag{f}{$f$}
\raisebox{1ex}{
\includegraphics[width=1cm]{Fig2.eps}}
}
\ \ \ \right|
+
\left|
\parbox{1cm}{
\psfrag{w}{$\bar{g}$}\psfrag{f}{$f$}
\raisebox{1ex}{
\includegraphics[width=1cm]{Fig2.eps}}
}
\ \ \right|\\
 = &
\left|
\parbox{1cm}{
\psfrag{w}{$(\beta_4-\bar{g})$}\psfrag{f}{$f$}
\raisebox{1ex}{
\includegraphics[width=1cm]{Fig2.eps}}
}
\ \ \ \right|\\
\le &
||f|_{L^{-1}\Delta'}||_{L^{\infty}}
\times
||\Gamma||_{L^{1}}
\times
L^{3}
\times
\max_{\Delta \in [L^{-1}\Delta']}
| \beta_{4,\Delta} - \bar{g} |\\
\le &
\left( \frac{1}{\sqrt{2}} L^{3-2[\phi]}\times\frac{1}{2} L^{-(3-[\phi])} \right)
L^{3}\times\frac{1}{2}
\bar{g}^{e_{4}}\\
\le & \frac{1}{4}L^{\frac{5}{2}} \bar{g}^{e_{4}} 
\end{align*}

In going to the second line we used that  
\[
\parbox{1cm}{
\psfrag{w}{$\bar{g}$}\psfrag{f}{$f$}
\raisebox{1ex}{
\includegraphics[width=1cm]{Fig2.eps}}
}
\ \ =0
\]

since $\Gamma$
integrates to $0$. \\

Therefore we have:

\begin{align*}
|\delta \beta_{3,1,\Delta'}[\vec{V}] |
 \bar{g}^{-e_{3}}
\le&
\bar{g}^{-e_{3}}
\sum_{b}
\bbone\left\{ 
\begin{array}{c}
k+b\le 4 \\
b\ge 1
\end{array}
\right\}
\frac{(k+b)!}{k!\ b!}\ L^{-k[\phi]}\
L^{\frac{5}{2}} \bar{g}^{e_{4}}\\
\le &
\frac{1}{4}L^{\frac{5}{2}}
\sum_{b}
\bbone\left\{ 
\begin{array}{c}
k+b\le 4 \\
b\ge 1
\end{array}
\right\}\frac{(k+b)!}{k!\ b!}
\end{align*}

In going to the second line we used that $e_{4} \ge e_{3}$ which is a
consequence of \eqref{econstbetaeq}.\\

We now observe that:

\[
\sum_{b}
\bbone\left\{ 
\begin{array}{c}
k+b\le 4 \\
b\ge 1
\end{array}
\right\}\frac{(k+b)!}{k!\ b!}
\le
9
\]

This proves the lemma.

\qed

\begin{Lemma}
\label{LAJ2}
Suppose that $\vec{V} \in B(\bar{V},\frac{1}{2})$. Then one
has the following bounds for the $W_{5}'$ and $W_{6}'$ components of
$\vec{V}' = RG_{\rm ex}[\vec{V}]$.\\

For all $\Delta' \in \mathbb{L}$ and for $k=5$ or $6$

\[
\left|W_{k,\Delta'}\right| \bar{g}^{-e_W}
\le
\mathbf{O}_2 L^{\frac{5}{2}},
\]

where $\mathbf{O}_2 =14$.

\end{Lemma}

\noindent{\bf Proof:}
For $k = 5$ we have:

\[
\left|W'_{5,\Delta'}\right|
\le
L^{3-5[\phi]} \max_{\Delta \in [L^{-1}\Delta']} \left| W_{5,\Delta'}
\right| 
+ 
6 L^{-5[\phi]} \left|
\parbox{1cm}{
\psfrag{w}{$W_6$}\psfrag{f}{$f$}
\raisebox{1ex}{
\includegraphics[width=1cm]{Fig2.eps}}
}
\ \ \right|
+
12 L^{-5[\phi]} \left|
\parbox{1.5cm}{
\psfrag{a}{$\beta_4$}\psfrag{b}{$\beta_3$}
\raisebox{-5ex}{
\includegraphics[width=1.6cm]{Fig3.eps}}
}
\ \ \ \right|
+
48 L^{-5[\phi]} \left|
\parbox{2.1cm}{
\psfrag{b}{$\beta_4$}\psfrag{f}{$f$}
\raisebox{1ex}{
\includegraphics[width=2.1cm]{Fig4.eps}}
}
\ \ \right|.
\]

We bound each of the diagrams: 

\begin{align*}
6 L^{-5[\phi]} \left|
\parbox{1cm}{
\psfrag{w}{$W_6$}\psfrag{f}{$f$}
\raisebox{1ex}{
\includegraphics[width=1cm]{Fig2.eps}}
}
\ \ \right| 
\le&
6 L^{-5[\phi]} 
|| f |_{L^{-1}\Delta'} ||_{L^{\infty}}
\times
||\Gamma||_{L^{1}}
\times L^{3}
\times
\max_{\Delta \in [L^{-1}\Delta']} |W_{6,\Delta}|\\
\le&
6 L^{-5[\phi]}\times\frac{1}{2} L^{-(3-[\phi])} \left( \frac{1}{\sqrt{2}} L^{3-2[\phi]} \right) L^{3} \times\frac{1}{2}\bar{g}^{e_{W}}\\
\le&
\frac{3}{2} L^{3-6[\phi]} \bar{g}^{e_{W}}\ .
\end{align*}

In going to the last line we dropped the factor of
$\frac{1}{\sqrt{2}}$. We continue to bound the other two
diagrams:

\begin{align*}
12 L^{-5[\phi]} \left|
\parbox{1.5cm}{
\psfrag{a}{$\beta_4$}\psfrag{b}{$\beta_3$}
\raisebox{-5ex}{
\includegraphics[width=1.6cm]{Fig3.eps}}
}
\ \ \ \right|
= &
12 L^{-5[\phi]}  \left|
\parbox{1.5cm}{
\psfrag{a}{$(\beta_4-\bar{g})$}\psfrag{b}{$\beta_3$}
\raisebox{-5ex}{
\includegraphics[width=1.6cm]{Fig3.eps}}
}\ \ \  
+  
\parbox{1.5cm}{
\psfrag{a}{$\bar{g}$}\psfrag{b}{$\beta_3$}
\raisebox{-5ex}{
\includegraphics[width=1.6cm]{Fig3.eps}}
}\ \ \ 
\right|\\
= &
12 L^{-5[\phi]}  \left|
\parbox{1.5cm}{
\psfrag{a}{$(\beta_4-\bar{g})$}\psfrag{b}{$\beta_3$}
\raisebox{-5ex}{
\includegraphics[width=1.6cm]{Fig3.eps}}
}\ \ \ 
\right|\\
= & 
 12 L^{-5[\phi]}  
\left| \int_{(L^{-1}\Delta')^2} {\rm d}^3x\ {\rm d}^3y\
( \beta_{4}(x) - \bar{g}) \Gamma(x-y) \beta_{3}(y) \right|\\
\le& 12 L^{-5[\phi]}  
\times 
\max_{\Delta \in [L^{-1}\Delta']} \left| \beta_{4,\Delta} - \bar{g}
\right| 
\times
||\Gamma||_{L^{1}}
\times 
\max_{\Delta \in [L^{-1}\Delta']} \left| \beta_{3,\Delta}\right|\times L^3\\
\le&
 12 L^{-5[\phi]} 
\left( \frac{1}{\sqrt{2}} L^{3-2[\phi]} \right)
L^{3} \times\frac{1}{4}
\bar{g}^{e_{3}+e_{4}} 
\le 
3 L^{6-7[\phi]} \bar{g}^{e_{3}+e_{4}}
\le 
3 L^{\frac{5}{2}} \bar{g}^{e_{3}+e_{4}} \ .
\end{align*}

In going to the second line we used the fact that
$\Gamma$ integrates to zero. In the last line we used $\epsilon \le 1$ so
that $6-7[\phi] \le \frac{5}{2}$. We now move to the third diagram.

\begin{align*}
48 L^{-5[\phi]} \left|
\parbox{2.1cm}{
\psfrag{b}{$\beta_4$}\psfrag{f}{$f$}
\raisebox{1ex}{
\includegraphics[width=2.1cm]{Fig4.eps}}
}\ \ 
\right|
 = &
48 L^{-5[\phi]} \left| 
\parbox{3cm}{
\psfrag{a}{$(\beta_4-\bar{g})$}\psfrag{b}{$\beta_4$}\psfrag{f}{$f$}
\raisebox{1ex}{
\includegraphics[width=3cm]{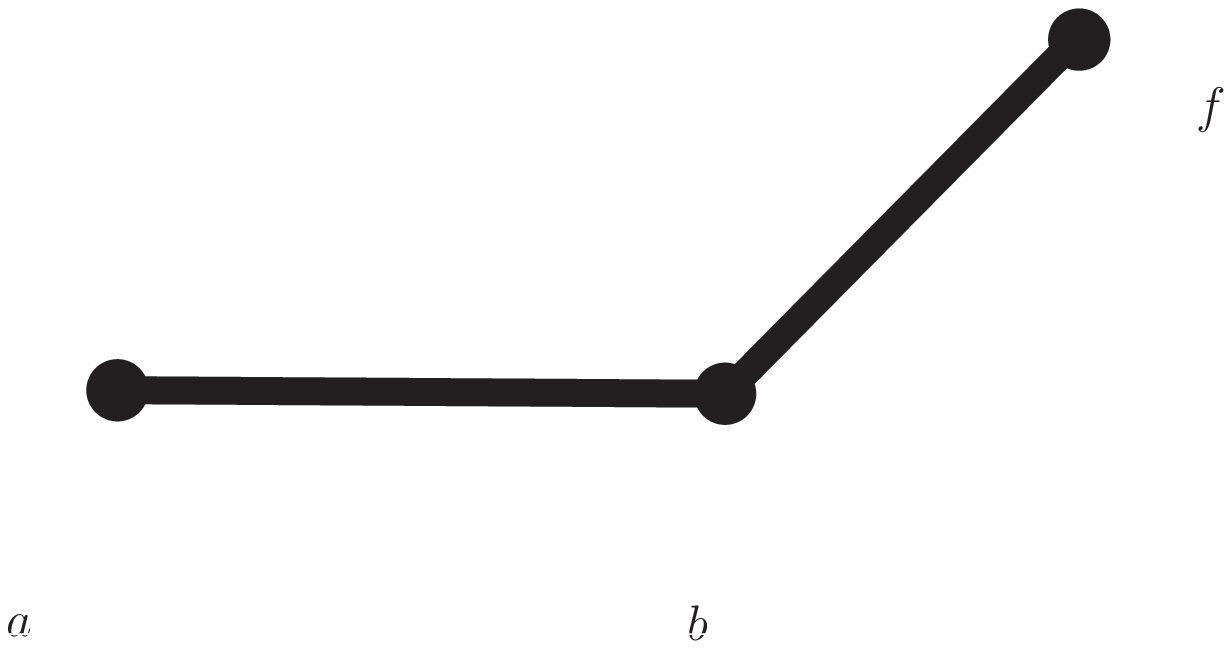}}
}\ \ \ 
+
\parbox{2.1cm}{
\psfrag{a}{$\bar{g}$}\psfrag{b}{$\beta_4$}\psfrag{f}{$f$}
\raisebox{1ex}{
\includegraphics[width=2.1cm]{AJfig1.eps}}
}\ \ 
\right|\\
= &
48 L^{-5[\phi]} \left|
\parbox{3cm}{
\psfrag{a}{$(\beta_4-\bar{g})$}\psfrag{b}{$\beta_4$}\psfrag{f}{$f$}
\raisebox{1ex}{
\includegraphics[width=3cm]{AJfig1.eps}}
}\ \ 
\right|\\
\le& 
48 L^{-5[\phi]}
\times
|| f |_{L^{-1}\Delta'} ||_{L^{\infty}} 
\times
||\Gamma||_{L^1}\\
 & \times
 \max_{\Delta \in [L^{-1}\Delta']}| \beta_{4,\Delta}|
\times
 ||\Gamma||_{L^1} 
\times
L^3
\times
 \max_{\Delta \in [L^{-1}\Delta']} |\beta_{4,\Delta} -\bar{g}|\\
\le&
48 L^{-5[\phi]}
\left( \frac{1}{\sqrt{2}} L^{3-2[\phi]} \right)^2
\times\frac{1}{2}L^{-(3-[\phi])}\times L^{3}
\left( \frac{3}{4} \bar{g}^{1+e_{4}} \right)
\le
9 L^{2\epsilon}\bar{g}^{1+e_{4}}\ .
\end{align*}

Above we used that $6-8[\phi]=2\epsilon$. Putting this together with our assumption on the size of the
unprimed $W_{5}$ gives us: 

\[
\left| W_{5,\Delta'} \right|
\le
\frac{1}{2} \bar{g}^{e_{W}}
+
\frac{3}{2}\bar{g}^{e_{W}}
+
12 L^{\frac{5}{2}} \bar{g}^{\min(1,e_{3})+e_{4}}
\]
where we
bounded $L^{3-5[\phi]}$ and $L^{3-6[\phi]}$ by $1$ as well as $L^{2\epsilon}$ by $L^{\frac{5}{2}}$.
We also
recall that $\min(1,e_{3})+e_{4} \ge e_{W}$ (assumed in \eqref{econstWeq}) to end up with estimate:

\[
\left|W_{5,\Delta'} \right|
\bar{g}^{-e_{W}}
\le
14
L^{\frac{5}{2}}\ .
\]

This proves the lemma for the case $k=5$. For $k=6$ we have:

\[
\left|W'_{6,\Delta'}\right|
\le
L^{-6[\phi]} \sum_{\Delta \in [L^{-1}\Delta']}
\left|W_{6,\Delta}\right|
+
8L^{-6[\phi]}
\left|
\parbox{1.5cm}{
\psfrag{b}{$\beta_4$}
\raisebox{-5ex}{
\includegraphics[width=1.6cm]{Fig1.eps}}
}\ \ \ 
\right|
\]

We bound the diagram above:

\begin{align*}
8L^{-6[\phi]}
\left|
\parbox{1.5cm}{
\psfrag{b}{$\beta_4$}
\raisebox{-5ex}{
\includegraphics[width=1.6cm]{Fig1.eps}}
}\ \ \ 
\right|
= &
8L^{-6[\phi]}
\left|
\parbox{1.5cm}{
\psfrag{b}{$\bar{g}$}
\raisebox{-5ex}{
\includegraphics[width=1.6cm]{Fig1.eps}}
}\ \ \ 
+ 
\parbox{1.7cm}{
\psfrag{b}{$(\beta_4-\bar{g})$}
\raisebox{-5ex}{
\includegraphics[width=1.7cm]{Fig1.eps}}
}\ \ \ \qquad
+
2\ \ 
\parbox{1.7cm}{
\psfrag{a}{$(\beta_4-\bar{g})$}\psfrag{b}{$\bar{g}$}
\raisebox{-5ex}{
\includegraphics[width=1.7cm]{Fig3.eps}}
}\ \qquad
\right|\\
= &
8L^{-6[\phi]}
\left|
\parbox{1.7cm}{
\psfrag{b}{$(\beta_4-\bar{g})$}
\raisebox{-5ex}{
\includegraphics[width=1.7cm]{Fig1.eps}}
}\ \ \ \ \qquad
\right|\\
\le &
8L^{-6[\phi]}
\times
||\Gamma||_{L^{1}}
\times
L^{3}
\times
\left( 
\max_{\Delta \in [L^{-1}\Delta']}
 \left| \beta_{4,\Delta} -\bar{g} \right| 
\right)^2\\
\le&
2 L^{2\epsilon}
\bar{g}^{2e_{4}}\ .
\end{align*}

We plug this back into our
earlier estimate for $\left| W_{6,\Delta'} \right|$ to get:

\[
\left| W_{6,\Delta'} \right| 
\le
\frac{1}{2}L^{3-6[\phi]} \bar{g}^{e_{W}} 
+ 
2 L^{2\epsilon}\bar{g}^{2e_{4}}\ .
\]

We again bound $L^{3-6[\phi]}$ by $1$. We also bound $L^{2\epsilon}$ by
$2$ in the small $\epsilon$ regime.
We also note that
$2e_{4} \ge e_{W}$ (this is a consequence of \eqref{econstWeq} and
\eqref{econstbetaeq}). This leaves us with the bound:

\[
\left| W_{6,\Delta'} \right| 
\bar{g}^{-e_{W}}
\le\frac{9}{2}
\]

This finishes the proof of the lemma.\qed\\

\begin{Lemma}\label{LAJ3}

Suppose that $\vec{V}$ in
$B(\bar{V},\frac{1}{2})$. Then one has the following bound
for the $R'$ component of $\vec{V}'=RG_{\rm ex}[\vec{V}]$: for all
$\Delta' \in \mathbb{L}$

\[
|||R'_{\Delta'}|||_{\bar{g}}
\bar{g}^{-e_{R}}
\le
\frac{3}{8}\ .
\]

\end{Lemma}

\noindent{\bf Proof:}
We use estimates from \eqref{mainestthm}. 

\begin{equation}\label{Rbound}
\begin{split}
|||R'_{\Delta}|||_{\bar{g}}\bar{g}^{-e_{R}} \le&
\left[ |||\mathcal{L}^{\vec{\beta},f}_{\Delta'}(R)|||_{\bar{g}} +
|||\xi_{R,\Delta'}(\vec{V})|||_{\bar{g}} \right] \bar{g}^{-e_{R}}\\
\le & \left[ \mathcal{B}_{R \mathcal{L}} L^{3-5[\phi]} \max_{\Delta \in
  [L^{-1}\Delta']} |||R_{\Delta}|||_{\bar{g}} +
|||\xi_{R,\Delta'}(\vec{V})|||_{\bar{g}} \right] \bar{g}^{-e_{R}} \\
\le & \left[\frac{1}{4} \bar{g}^{e_{R}} + B_{R\xi}
  \bar{g}^{\frac{11}{4} - 3 \eta} \right] \bar{g}^{-e_{R}}\\
\le & \frac{1}{4} + B_{R\xi} \bar{g}^{\frac{11}{4} - 3 \eta -e_{R}} \le  \frac{3}{8}
\end{split}
\end{equation}

In going to the third line we used that $L$ has been fixed to guarantee
$\mathcal{B}_{R\mathcal{L}} L^{3-5[\phi]} \le \frac{1}{2}$. This was done in
\eqref{Lchoiceeq}.\\
In the last line we used that $\frac{11}{4} - 3 \eta -e_{R} > 0$ which
was assumed in \eqref{econst5eq}. Thus by requiring $\epsilon$ is
sufficiently small we can guarantee  $B_{R\xi} \bar{g}^{\frac{11}{4} -
  3 \eta -e_{R}} \le \frac{1}{8}$. 
\qed\\

\begin{Lemma} \label{LAJ4}
Suppose that $\vec{V}$ in
$B(\bar{V},\frac{1}{2})$. Then one has the following bound
for the $\beta'_{4}$ component of $\vec{V}'=RG_{\rm ex}[\vec{V}]$: For all
$\Delta' \in \mathbb{L}$

\[
\left| \beta'_{4,\Delta'} - \bar{g} \right|\bar{g}^{-e_{4}} 
\le
\mathbf{O}_{3}
\]

where $\mathbf{O}_{3}=434 + \mathcal{O}_{26}$ with $\mathcal{O}_{26}$ defined in the statement of Lemma \ref{L34lem}.
\end{Lemma}

\noindent{\bf Proof:}
Due to our assumption on $\vec{V}$ for any $\Delta \in
\mathbb{L}$ we can write $\beta_{4,\Delta} = \bar{g} + \delta
g_{\Delta}$ where $|\delta g_{\Delta}| <
\frac{1}{2}\bar{g}^{e_{4}}$. We
substitute this into the flow equation to get the following:

\begin{equation}\label{beta4bound}
\begin{split}
\beta'_{4,\Delta'}& = L^{3-4[\phi]} \bar{g} + L^{-4[\phi]}\sum_{\Delta \in [L^{-1}\Delta']} \delta g_{\Delta} - \delta \beta_{4,2,\Delta'}[\vec{V}] + \xi_{4,\Delta'}[\vec{V}]\\
&  = L^{3-4[\phi]} \bar{g} - 36L^{-4[\phi]}
\ \ \parbox{1.9cm}{
\psfrag{a}{$\bar{g}+\delta g$}
\psfrag{b}{$\bar{g}+\delta g$}
\raisebox{-1ex}{
\includegraphics[width=1.9cm]{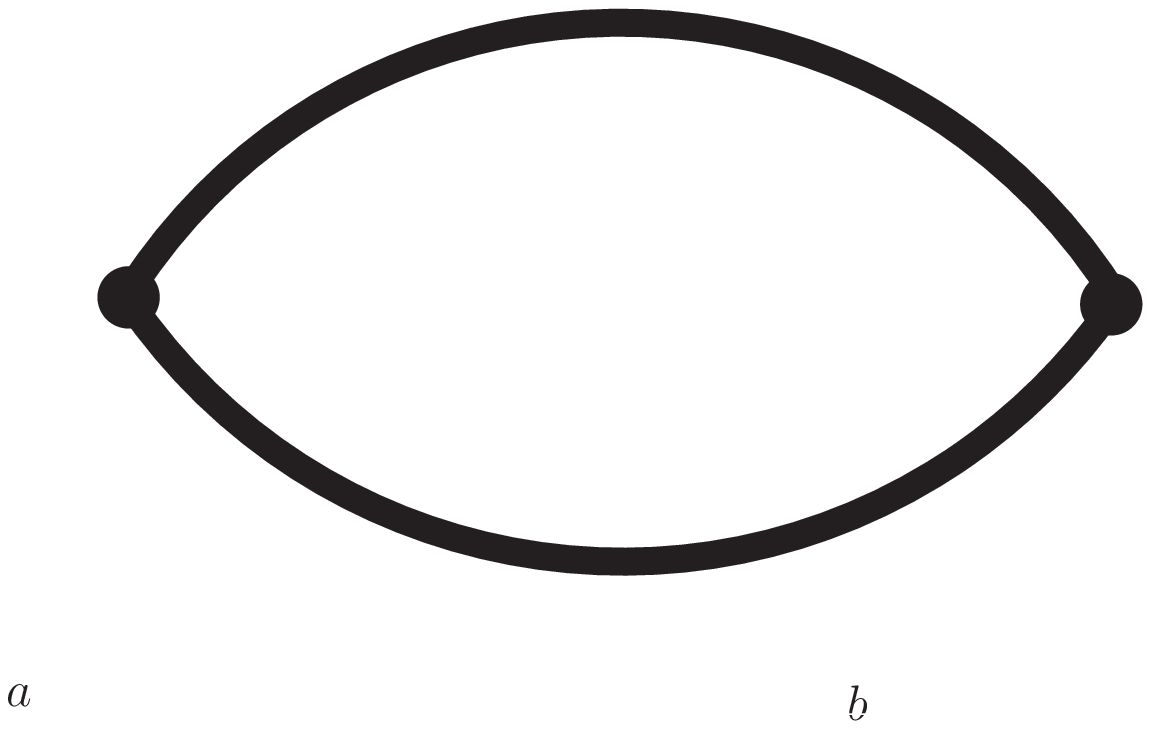}}
}\qquad
- \widetilde{\delta
    \beta}_{4,2,\Delta'}[\vec{V}] + \xi_{4,\Delta'}[\vec{V}] +  L^{-4[\phi]}\sum_{\Delta \in [L^{-1}\Delta']} \delta g_{\Delta}\ .
\end{split}
\end{equation}

We have used the fact that $\delta
\beta_{4,1,\Delta}[\vec{V}]=0$. In the formula above $\widetilde{\delta
  \beta}_{4,2,\Delta'}[\vec{V}]$ is defined to be
$\delta\beta_{4,2,\Delta'}[\vec{V}]$ with the graph that we have made
explicit removed:

\[
\widetilde{\delta\beta}_{4,2,\Delta'} \left[ \vec{V} \right] :=\sum_{a_1,a_2,b_1,b_2,m}
\bbone\left\{
\begin{array}{c}
a_i+b_i\le 4 \\
a_i\ge 0\ ,\ b_i\ge 1\\
m = 1
\end{array}
\right\}
\frac{(a_1+b_1)!\ (a_2+b_2)!}{a_1!\ a_2!\ m!\ (b_1-m)!\ (b_2-m)!}
\]
\[
\times \frac{1}{2} C(a_1,a_2|4)\times
L^{-(a_1+a_2)[\phi]}\times C_0(0)^{\frac{a_1+a_2-4}{2}}\times
\parbox{4cm}{
\psfrag{a}{$\beta_{a_1+b_1}$}
\psfrag{b}{$\beta_{a_2+b_2}$}
\psfrag{c}{$\scriptstyle{b_1-m}$}
\psfrag{d}{$\scriptstyle{b_2-m}$}
\psfrag{m}{$m$}
\psfrag{f}{$f$}
\raisebox{-1ex}{
\includegraphics[width=4cm]{Fig6.eps}}
}
\]
\[
+\sum_b
\bbone\left\{
\begin{array}{c}
4+b=5\ {\rm or}\ 6 \\
b\ge 0
\end{array}
\right\}
\frac{(k+b)!}{k!\ b!}\  L^{-k[\phi]}\ 
\parbox{2.1cm}{
\psfrag{a}{$W_{k+b}$}\psfrag{f}{$f$}\psfrag{b}{$b$}
\raisebox{-1ex}{
\includegraphics[width=2.1cm]{Fig5.eps}}
}
\]
\[
= \delta \beta_{4,2,\Delta'}[\vec{V}] -  \frac{1}{2} L^{-4[\phi]} \frac{4!4!}{2!2!2!} 
\ \ \parbox{1.9cm}{
\psfrag{a}{$\bar{g}+\delta g$}
\psfrag{b}{$\bar{g}+\delta g$}
\raisebox{-1ex}{
\includegraphics[width=1.9cm]{AJfig2.eps}}
}\qquad
\]
Indeed, first note that there is no graph with $m=3$. This is because this would imply $a_1,a_2\le 1$
which contradicts $a_1+a_2\ge 4$ imposed by the nonvanishing of the connection coefficient $C(a_1,a_2|4)$.
Also the removed graph is the only one with $m=2$. This is because $b_1,b_2\ge 2$ implies $a_1,a_2\le 4-2=2$, but
the connection coefficient requires $a_1+a_2\ge 4$ so we are forced to have $a_1=a_2=2$ which implies $b_1,b_2\le 2$
and therefore $b_1=b_2=2$. 

We note that we can decompose the graph above as follows:

\[
\parbox{1.9cm}{
\psfrag{a}{$\bar{g}+\delta g$}
\psfrag{b}{$\bar{g}+\delta g$}
\raisebox{-1ex}{
\includegraphics[width=1.9cm]{AJfig2.eps}}
}\qquad
 =
\ \ \parbox{1.9cm}{
\psfrag{a}{$\bar{g}$}
\psfrag{b}{$\bar{g}$}
\raisebox{-1ex}{
\includegraphics[width=1.9cm]{AJfig2.eps}}
}\qquad
+ 2 
\parbox{1.9cm}{
\psfrag{a}{$\bar{g}$}
\psfrag{b}{$\delta g$}
\raisebox{-1ex}{
\includegraphics[width=1.9cm]{AJfig2.eps}}
}\qquad
+
\parbox{1.9cm}{
\psfrag{a}{$\delta g$}
\psfrag{b}{$\delta g$}
\raisebox{-1ex}{
\includegraphics[width=1.9cm]{AJfig2.eps}}
}\qquad
\]

We now use the fact that $\bar{g}$ is an approximate fixed point:

\[
\bar{g} = L^{\epsilon}\bar{g}-A_1 \bar{g}^2=
 L^{3-4[\phi]} \bar{g} - 36 L^{-4[\phi]}
\ \ \parbox{1.9cm}{
\psfrag{a}{$\bar{g}$}
\psfrag{b}{$\bar{g}$}
\raisebox{-1ex}{
\includegraphics[width=1.9cm]{AJfig2.eps}}
}\qquad
\]

Using this we can write:

\begin{equation}\label{beta4expansionwork}
\begin{split}
\beta'_{4,\Delta'} =& \bar{g} + L^{-4[\phi]} \sum_{\Delta \in
  [L^{-1}\Delta']} \delta g_{\Delta}\\
& - 36 L^{-4[\phi]}  \left(2
\ \ \parbox{1.9cm}{
\psfrag{a}{$\bar{g}$}
\psfrag{b}{$\delta g$}
\raisebox{-1ex}{
\includegraphics[width=1.9cm]{AJfig2.eps}}
}\qquad
+
\ \ \parbox{1.9cm}{
\psfrag{a}{$\delta g$}
\psfrag{b}{$\delta g$}
\raisebox{-1ex}{
\includegraphics[width=1.9cm]{AJfig2.eps}}
}\qquad
\right)\\
& - \widetilde{\delta \beta}_{4,2,\Delta'} [ \vec{V} ] +
\xi_{4,\Delta'}[\vec{V}]\ .
\end{split}
\end{equation}

We now describe how to bound the second and third lines of \eqref{beta4expansionwork}. By the same
arguments as used in Lemma \ref{L34lem}  the contribution of the  two
graphs on the second line can
each be bounded by $4 L^{5}\bar{g}^{2-2\eta}$ as follows from the very coarse bounds $\bar{g}\le \bar{g}^{1-\eta}$
and $|\delta g|\le \bar{g}^{1-\eta}$.
This gives us:

\begin{align*}
\left| 
36
L^{-4[\phi]} 
\left(2\ \  
\parbox{1.9cm}{
\psfrag{a}{$\bar{g}$}
\psfrag{b}{$\delta g$}
\raisebox{-1ex}{
\includegraphics[width=1.9cm]{AJfig2.eps}}
}\ \ 
+
\ \ \parbox{1.9cm}{
\psfrag{a}{$\delta g$}
\psfrag{b}{$\delta g$}
\raisebox{-1ex}{
\includegraphics[width=1.9cm]{AJfig2.eps}}
}\ \ 
\right)
\right|
\le &
36
\left[
2 
\left|\ \ 
\parbox{1.9cm}{
\psfrag{a}{$\bar{g}$}
\psfrag{b}{$\delta g$}
\raisebox{-1ex}{
\includegraphics[width=1.9cm]{AJfig2.eps}}
}\ \ 
\right|
+
\left|
\ \ \parbox{1.9cm}{
\psfrag{a}{$\delta g$}
\psfrag{b}{$\delta g$}
\raisebox{-1ex}{
\includegraphics[width=1.9cm]{AJfig2.eps}}
}\ \ 
\right|
\right]\\
\ & \ \\
\le&
36\times 3\times 4
L^{5} 
\bar{g}^{2-2\eta}
= 432 L^{5}  \bar{g}^{2-2\eta}\ .
\end{align*}

Note that in the first line we dropped the factor of $L^{-4[\phi]}$. The quantity
$\widetilde{\delta \beta}_{4,2,\Delta'} [ \vec{V} ]$ on the third line
of  \eqref{beta4expansionwork} can be
bounded by $\mathcal{O}_{26}L^{5}\bar{g}^{2-2\eta}$ as in Lemma
\ref{L34lem} (we are overestimating since we are summing over fewer
graphs). We combine this with the estimate on
$\xi_{4,\Delta'}[\vec{V}]$ from Theorem \ref{mainestthm} to get

\begin{equation*}
\begin{split}
 \left| \beta'_{4,\Delta'} - \bar{g} \right|\bar{g}^{-e_{4}} \le&
\frac{1}{2}  L^{3-4[\phi]} + \left(432 + \mathcal{O}_{26}
   \right)L^{5}\bar{g}^{2-2\eta-e_{4}} + B_{4}\bar{g}^{e_{R}-e_{4}} \\
\le& \frac{1}{2}  L^{3-4[\phi]} + \left(432 + \mathcal{O}_{26}\right) + 1\\
\le& 1+ \left(432 + \mathcal{O}_{26}\right) + 1
\end{split}
\end{equation*}

Note that in going to the second line that we used $2 - 2\eta - e_{4} > 0$, this is a consequence of
\eqref{econstbetaeq}. Indeed this allows to have $L^{5}\bar{g}^{2-2\eta-e_{4}}\le 1$ in the small $\epsilon$ regime.
We also used that $e_{4} < e_{R}$ (a consequence
of \eqref{econst4eq}), thus we can guarantee
$B_{4}\bar{g}^{e_{R}-e_{4}} \le 1$ for $\epsilon$ sufficiently
small. In going to the third
we used the bound $L^{3-4[\phi]} = L^{\epsilon} \le 2$ for $\epsilon$ small.
\qed\\

\begin{Lemma}\label{LAJ5}
Suppose that $\vec{V}$ in $B(\bar{V},\frac{1}{2})$. Then one has the following bound
for the $\beta'_{k}$ components of $\vec{V}'=RG_{\rm ex}[\vec{V}]$ when $k=1,2,3$: for all
$\Delta' \in \mathbb{L}$

\[
|\beta'_{k,\Delta'}| \bar{g}^{-e_{k}}
\le
\mathbf{O}_{4}
L^{\frac{5}{2}}
\]

where $\mathbf{O}_{4} = \mathbf{O}_{1}+2 $.
 
\end{Lemma}

\noindent{\bf Proof:}
From the flow equations one has:

\begin{equation*}
\begin{split}
|\beta'_{k,\Delta'}| &\le \left| L^{-k[\phi]} \sum_{\Delta \in
    [L^{-1}\Delta']} \beta_{k,\Delta} \right| + \left|\delta  \beta_{k,1,\Delta'}[\vec{V}] \right| +\left|\delta b_{k,2,\Delta'}[\vec{V}] \right| + \left| \xi_{k,\Delta'}[\vec{V}]
\right| \\
& \le L^{3-k[\phi]} \bar{g}^{e_{k}} +\mathbf{O}_{1} L^{5/2}
\bar{g}^{e_{k}} + \mathcal{O}_{26}L^{5}\bar{g}^{2-2\eta} +
\frac{1}{2}B_{k}\bar{g}^{ e_{R} }\ .
\end{split}
\end{equation*}

The bound on the third term on the right hand side of the first line
is from Lemma \ref{L34lem} and the bound on the last term of the first
line is from Theorem \eqref{mainestthm}. We used Lemma \ref{LAJ1} to bound $ \left|\delta  \beta_{k,1,\Delta'}[\vec{V}]
\right|$. Then for $k=1,2,3$ we have:

\begin{equation*}
\begin{split}
|\beta'_{k,\Delta'}| \bar{g}^{-e_{k}} \le & L^{3-k[\phi]} +
\mathbf{O}_{1} L^{5/2} \bar{g}^{e_{k+1}-e_{k}} +
\mathcal{O}_{26}L^{5}\bar{g}^{2-2\eta - e_{k}} + \frac{1}{2}B_{k}\bar{g}^{e_{R}-e_{k}}\\
\le &  L^{\frac{5}{2}} + \mathbf{O}_{1} L^{5/2} + 1\\
\end{split}
\end{equation*}

In going to the second line  we used that $e_{k+1} \ge e_{k}$ which is a
consequence of \eqref{econstbetaeq}. We also used that $e_{R} > e_{k}$
and $2-2\eta > e_{k}$ which come from \eqref{econstbetaeq} and
\eqref{econst4eq}. Thus the sum of the last two terms on the first line can be
made smaller than $1$ by requiring that $\epsilon$ be sufficiently small. We also used that for $\epsilon \le 1$ and $k \ge 1$ one has $3-k[\phi]
\le \frac{5}{2}$.
\qed\\

\begin{Lemma}\label{LAJ6}
$RG_{\rm ex}$ is well defined and analytic on
$B(\bar{V},\frac{1}{2})$. Additionally one has the following
uniform bound for $\vec{V} \in B(\bar{V},\frac{1}{2})$:

\begin{equation}\label{expansionuniformbound}
||RG_{\rm ex}[\vec{V}]-\bar{V}|| \le \mathbf{O}_{5} L^{\frac{5}{2}}
\end{equation}

where $\mathbf{O}_{5} = \max \left( \mathbf{O}_{2},\mathbf{O}_{3},\mathbf{O}_{4} \right)$.

\end{Lemma}

\noindent{\bf Proof:}

The fact that the map is well defined and is analytic comes from the
Theorem \ref{mainestthm} and inspection of the formulas for $\delta
\beta_{k,j,\Delta}$ for $j=1,2$ and $k=1,2,3,4$. We now
establish the uniform bound.\\

Let $\vec{V}'=RG_{\rm ex}[\vec{V}]$. We have the sufficient estimates on
$\beta'_{k}$ for $k=1,2,3,4$ from Lemmas \ref{LAJ4} and \ref{LAJ5}. We
have sufficient estimates on $W'_{k}$ for $k=5$ and $6$ from Lemma
\ref{LAJ2}. A sufficient estimate on $R'$ comes from Lemma \ref{LAJ3}. All that is left is
estimating $f'$.\\

 Note that for any $\Delta' \in \mathbb{L}$ we have

\begin{align*}
\left| f'_{\Delta'} \right| L^{3-[\phi]}
\le&
L^{3-[\phi]}
L^{3-[\phi]}
\max_{\Delta \in [L^{-1}\Delta']} |f_{\Delta}|\\
\le&
L^{3-[\phi]}
L^{3-[\phi]}
\left( \frac{1}{2} L^{-(3-[\phi])} \right)
= \frac{1}{2} L^{3-[\phi]} \le \frac{1}{2}L^{\frac{5}{2}}\ .
\end{align*}

On the last line we used the assumption $\epsilon \le 1$ which
implies $3-[\phi] \le \frac{5}{2}$. Finally note that $\frac{3}{8}<\frac{1}{2}<14=\mathbf{O}_2$
to get the formula for the constant $\mathbf{O}_5$.
\qed\\

\begin{Proposition}\label{expansionbound}
For any $\vec{V}^1, \vec{V}^2 \in \bar{B}(\bar{V},\frac{1}{6})$
one has:

\[
||RG_{\rm ex}[\vec{V}^{1}] 
- RG_{\rm ex}[\vec{V}^{2}]||
\le
\mathbf{O}_{6}
L^{\frac{5}{2}}
||\vec{V}^{1} - \vec{V}^{2}||,
\]

where $\mathbf{O}_{6}=4\mathbf{O}_{5}$.
\end{Proposition}

\noindent{\bf Proof:}

By Lemma \ref{LAJ6} we know that $RG_{\rm ex}$ is an analytic map taking
$B(\bar{V},\frac{1}{2})$ into
$\bar{B}(\bar{V},\mathbf{O}_{5}L^{\frac{5}{2}})$. We get the
desired inequality by applying Lemma
\ref{Lipschitzlem} with the choice $\nu=\frac{1}{3}$.\qed\\

After the previous estimates we now return to the analysis of the $r\rightarrow -\infty$
and $s\rightarrow \infty$ limits of $\mathcal{S}_{r,s}^{\rm T,MD}(\tilde{f},\tilde{j})$
which in fact does not depend on $s$ such that $s\ge q_{+}$.
Since the summation range $q_{-}\le q<q_{+}$ is fixed and finite, all we need is to show that $RG_{\rm ex}$
remain in the domains of definition and analyticity, despite the temporary expansion with rate controlled by
Lemma \ref{LAJ6} and Proposition \ref{expansionbound}.

The quantity of interest, as delivered by \S\ref{algconsec}, is
\[
\mathcal{S}_{r,s}^{\rm T,MD}(\tilde{f},\tilde{j})
=
\sum_{q_{-} \le q < q_{+}}
\sum_{\substack{\Delta \in \mathbb{L} \\ \Delta \subset
    \Lambda_{q_{+}-q-1}}}
\left(
\delta b_{\Delta}
[
\vec{V}^{(r,q)}(\tilde{f},\tilde{j})
]
-
\delta b_{\Delta}
[
\vec{V}^{(r,q)}(0,0)
]
\right)
\]
where
\[
\vec{V}^{(r,q)}(\tilde{f},\tilde{j})=RG_{\rm ex}^{q-q_{-}}\left(
\vec{V}^{(r,q_{-})}(\tilde{f},\tilde{j})\right)
\]
with
\[
\vec{V}^{(r,q_{-})}(\tilde{f},\tilde{j})
=
\mathcal{J}_{q_{+}-q_{-}}
\left(
\tilde{f}_{\rightarrow (-q_{-})},
\left(
RG^{q_{-}-r}
\left(
v 
- 
\alpha_{\rm u}^{r}
Y_{2}
\tilde{j}_{L^{-q_{-}}\Delta}
e_{\phi^{2}}
\right)
\right)_{
\substack{\Delta \in \mathbb{L} \\ \Delta \subset \Lambda_{q_{+}-q_{-}}}}
,
RG^{q_{-}-r}(v)
\right).
\]
It follows from our definitions for the norms and the reinjection map $\mathcal{J}$ that
\[
||\vec{V}^{(r,q_{-})}(\tilde{f},\tilde{j})-
\vec{V}^{(r,q_{-})}(0,0)||\qquad\qquad\qquad\qquad\qquad
\]
\[
\qquad\qquad\qquad
=\max\left\{
||\tilde{f}_{\rightarrow (-q_{-})}||_{L^{\infty}},
\max\limits_{\substack{\Delta \in \mathbb{L} \\ \Delta \subset \Lambda_{q_{+}-q_{-}}}}
||RG^{q_{-}-r}(v-\alpha_{\rm u}^{r}Y_2\tilde{j}_{L^{-q_{-}}\Delta}e_{\phi^2})-RG^{q_{-}-r}(v)||
\right\}\ .
\]
We also have
\[
||\tilde{f}_{\rightarrow (-q_{-})}||_{L^{\infty}}=L^{(3-[\phi])q_{-}}
||\tilde{f}||_{L^{\infty}}\ .
\]
We slightly strengthen the requirement in (\ref{jbound}) by imposing
\[
||\tilde{j}||_{L^{\infty}}\le
[240\mathcal{C}_{1}(\epsilon)\alpha_{\rm u}^{q_{-}}Y_2\bar{g}^{-e_2}]^{-1}
\]
which implies
\[
||-\alpha_{\rm u}^{q_{-}}Y_2\tilde{j}_{L^{-q_{-}}\Delta}e_{\phi^2}||\le\frac{1}{240\mathcal{C}_{1}(\epsilon)}
\]
for all $\Delta\in\mathbb{L}$ such that $\Delta\subset\Lambda_{q_{+}-q_{-}}$.
Thus by Lemma \ref{devseedlem}
\[
\begin{split}
||RG^{q_{-}-r}(v-\alpha_{\rm u}^{r}Y_2\tilde{j}_{L^{-q_{-}}\Delta}e_{\phi^2})-RG^{q_{-}-r}(v)||
 & \le 11\mathcal{C}_{1}(\epsilon)||-\alpha_{\rm u}^{q_{-}}Y_2\tilde{j}_{L^{-q_{-}}\Delta}e_{\phi^2}||\\
 & \le 11\mathcal{C}_{1}(\epsilon)\alpha_{\rm u}^{q_{-}}Y_2 \bar{g}^{-e_2}\times ||\tilde{j}||_{L^{\infty}}
\end{split}
\]
and therefore
\[
||\vec{V}^{(r,q_{-})}(\tilde{f},\tilde{j})-
\vec{V}^{(r,q_{-})}(0,0)||\le\max
\left\{
L^{(3-[\phi])q_{-}}
||\tilde{f}||_{L^{\infty}},
11\mathcal{C}_{1}(\epsilon)\alpha_{\rm u}^{q_{-}}Y_2 \bar{g}^{-e_2}\times ||\tilde{j}||_{L^{\infty}}
\right\}\ .
\]

On the other hand, minding the $\bar{g}$ shift for $\beta_4$ components only, we easily see that
\[
||\vec{V}^{(r,q_{-})}(0,0)-\bar{V}||=||\iota(RG^{q_{-}-r}(v))-\bar{V}||=||RG^{q_{-}-r}(v)||
\]
where the latter quantity can be computed as in section \S\ref{dyn1sec}, i.e., via the norm inherited by $\mathcal{E}$
from $\mathcal{E}_{\rm ex}$
and expressed in $(\delta g,\mu,R)$ coordinates. 

By construction of $W^{\rm s,loc}$, $||RG^{q_{-}-r}(v)||\le \frac{\rho}{3}$
with $\rho\in\left(0,\frac{1}{12}\right)$ as yet unspecified.
We thus have
\[
||\vec{V}^{(r,q_{-})}(0,0)-\bar{V}||\le \frac{1}{12}\ .
\]
Provided we also have
\[
\left(\mathbf{O}_{6}L^{\frac{5}{2}}\right)^{q_{+}-q_{-}}\times
\max
\left\{
L^{(3-[\phi])q_{-}}
||\tilde{f}||_{L^{\infty}},
11\mathcal{C}_{1}(\epsilon)\alpha_{\rm u}^{q_{-}}Y_2 \bar{g}^{-e_2}\times ||\tilde{j}||_{L^{\infty}}
\right\}
\le \frac{1}{12}
\]
then a trivial inductive application of Proposition \ref{expansionbound}
will garantee that for all $q$, $q_{-}\le q\le q_{+}$,
\[
||\vec{V}^{(r,q_{-})}(\tilde{f},\tilde{j})-\bar{V}||\le \frac{1}{12}
\]
so one remains, throughout the iterations, in the domain of definition and analyticity of $RG_{\rm ex}$
as well as the $\delta b$ functions.

As a result of Theorem \ref{unnamedtheorem1} we then immediately obtain, regardless of the order of limits,
\[
\lim_{\substack{r\rightarrow -\infty \\ s\rightarrow \infty}}
\mathcal{S}_{r,s}^{T,MD}(\tilde{f},\tilde{j})=\mathcal{S}^{T,MD}(\tilde{f},\tilde{j})
\]
where
\[
\mathcal{S}^{T,MD}(\tilde{f},\tilde{j})
=
\sum_{q_{-} \le q < q_{+}}
\sum_{\substack{\Delta \in \mathbb{L} \\ \Delta \subset
    \Lambda_{q_{+}-q-1}}}
\left(
\delta b_{\Delta}
[
\vec{V}^{(-\infty,q)}(\tilde{f},\tilde{j})
]
-
\delta b_{\Delta}
[
\iota(v_{\ast})
]
\right)
\]
with
\[
\vec{V}^{(-\infty,q)}(\tilde{f},\tilde{j})=RG_{\rm ex}^{q-q_{-}}\left(
\vec{V}^{(-\infty,q_{-})}(\tilde{f},\tilde{j})\right)
\]
for
\begin{equation}
\vec{V}^{(-\infty,q_{-})}(\tilde{f},\tilde{j})
=
\mathcal{J}_{q_{+}-q_{-}}
\left(
\tilde{f}_{\rightarrow (-q_{-})},
\left(
\Psi
v, 
- \alpha_{\rm u}^{q_{-}}
Y_{2}
\tilde{j}_{L^{-q_{-}}\Delta}
e_{\phi^{2}}
\right)_{
\substack{\Delta \in \mathbb{L} \\ \Delta \subset \Lambda_{q_{+}-q_{-}}}}
,
v_{\ast}
\right).
\label{Vqminuseq}
\end{equation}
Analyticity of $\mathcal{S}^{T,MD}(\tilde{f},\tilde{j})$ is also immediate.

For the purposes of the next section we also note that $\vec{V}^{(r,q_{+})}(\tilde{f},\tilde{j})$
satisfies the bound
\begin{equation}
||\vec{V}^{(r,q_{+})}(\tilde{f},\tilde{j})-
\vec{V}^{(r,q_{+})}(0,0)||\le
\left(\mathbf{O}_{6}L^{\frac{5}{2}}\right)^{q_{+}-q_{-}}\times
\max
\left\{
L^{(3-[\phi])q_{-}}
||\tilde{f}||_{L^{\infty}},
11\mathcal{C}_{1}(\epsilon)\alpha_{\rm u}^{q_{-}}Y_2 \bar{g}^{-e_2}\times ||\tilde{j}||_{L^{\infty}}
\right\}\ .
\label{seedforIReq}
\end{equation}

\subsection{The infrared regime}\label{infregsec}

In this section we are concerned with showing that essentially the differential of
$RG_{\rm ex}$ at any suitable $\vec{V}_{\rm bk} \in \mathcal{E}_{\rm bk}$ in any
direction $\dot{V} \in \mathcal{E}_{\rm pt}$ is a
contraction. \\

We will introduce new notation to facilitate the lemmas below. For
$\vec{V}_{\rm bk}\in \mathcal{E}_{\rm bk}$ we write:

\[
\vec{V}_{\rm bk}
=
\left\{
V_{\rm bk}
\right\}_{\Delta \in \mathbb{L}}
=
\left\{
(
\beta_{4,{\rm bk}},
\dots,
\beta_{1,{\rm bk}}, 
W_{5,{\rm bk}},
W_{6,{\rm bk}},
f_{\rm bk},
R_{\rm bk}
)
\right\}_{\Delta \in \mathbb{L}}\ .
\]
Note that we do need to burden the notation with $\Delta$ subscripts since
the quantities above are independent of the box $\Delta$ by definition of being in $\mathcal{E}_{\rm bk}$.

Similarly for $\dot{V} \in \mathcal{E}_{\rm pt}$ we write:

\[
\dot{V}
=
\left\{
\dot{V}_{\Delta}
\right\}_{\Delta \in \mathbb{L}}
=
\left\{
(
\dot{\beta}_{4,\Delta},
\dots,
\dot{\beta}_{1,\Delta}, 
\dot{W}_{5,\Delta},
\dot{W}_{6,\Delta},
\dot{f}_{\Delta},
\dot{R}_{\Delta}
)
\right\}_{\Delta \in \mathbb{L}}\ .
\]

Note that $\dot{V}_{\Delta} = 0$ for $\Delta \not = \Delta(0)$. We
also recall that $RG_{\rm ex}[\vec{V}_{\rm bk} + \dot{V} ] - RG_{\rm ex}
[\vec{V}_{\rm bk}] \in \mathcal{E}_{\rm pt}$ and so in our estimates we are
only concerned with the $\Delta(0)$ component of $RG_{\rm ex}[\vec{V}_{\rm bk} + \dot{V} ] - RG_{\rm ex}
[\vec{V}_{\rm bk}] \in \mathcal{E}_{\rm pt}$.

\begin{Lemma}\label{LAJ8}
Let $\vec{V}_{\rm bk} \in B(\bar{V},\frac{1}{4}) \cap \mathcal{E}_{\rm bk}$ and $\dot{V}
\in B(0,\frac{1}{4}) \cap \mathcal{E}_{\rm pt}$. Then one has
the bound for $k=1,2,3,4$ and all $\Delta' \in \mathbb{L}$:

\[
\left|
\delta \beta_{k,1,\Delta(0)}
\left[ \vec{V}_{\rm bk} + \dot{V} \right]
-
\delta \beta_{k,1,\Delta(0)}
\left[ \vec{V}_{\rm bk}\right]
\right|
\bar{g}^{-e_{k}}
\le
\bbone \{  1 \le k \le 3\}
\mathbf{O}_{7}
L^{-\frac{9}{4}}
||\dot{V}||^{2},
\]
where $\mathbf{O}_{7} = \left(6 + 21 \times 2^{\frac{3}{2}} \right)$.
\end{Lemma}

\noindent{\bf Proof:}
We again note that the vanishing for $k=4$ follows by inspection of
the definition of $\delta \beta_{k,1,\Delta(0)}$. We now observe that $\delta \beta_{k,1,\Delta(0)} [\vec{V}_{\rm bk} ]$
vanishes. Indeed, by definition we have

\[
\delta \beta_{k,1,\Delta(0)}
[\vec{V}_{\rm bk}]
=
- \sum_{b}
\bbone
\left\{ 
\begin{array}{c}
k+b\le 4 \\
b\ge 1
\end{array}
\right\}
\frac{(k+b)!}{k!\ b!}\  L^{-k[\phi]}
\ \ \ \parbox{2.8cm}{
\psfrag{a}{$\hspace{-.2cm}\beta_{k+b,{\rm bk}}$}
\psfrag{f}{$\hspace{-.2cm}f_{\rm bk}$}
\psfrag{g}{$f_{\rm bk}$}
\psfrag{b}{$b$}
\raisebox{-1ex}{
\includegraphics[width=2.8cm]{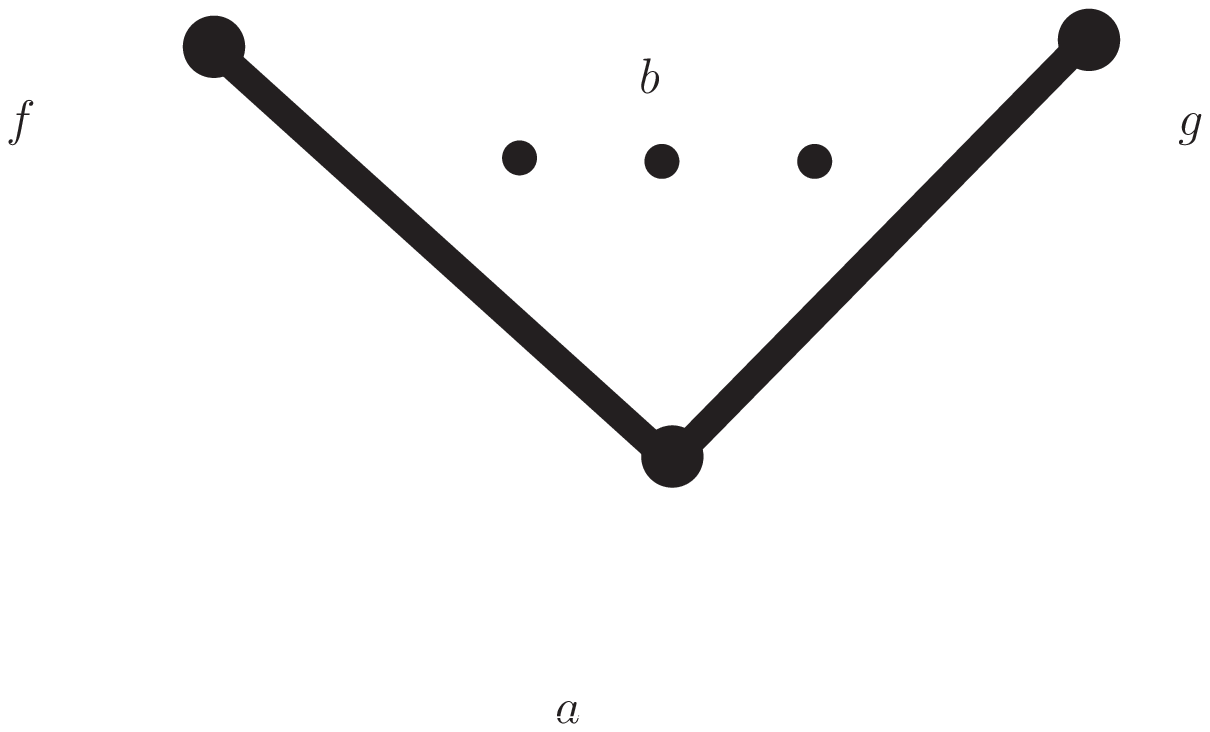}}
}
\]

However one has that 

\[
\parbox{2.8cm}{
\psfrag{a}{$\hspace{-.2cm}\beta_{k+b,{\rm bk}}$}
\psfrag{f}{$\hspace{-.2cm}f_{\rm bk}$}
\psfrag{g}{$f_{\rm bk}$}
\psfrag{b}{$b$}
\raisebox{-1ex}{
\includegraphics[width=2.8cm]{AJfig4.eps}}
}\qquad
=
0\ .
\]

This is because we have at least one integration vertex of degree $1$
which has been assigned a coupling $f_{\rm bk}$ which is constant over the
integration region $L^{-1}\Delta(0)$. Using ultrametricity and the fact that $\Gamma$
integrates to $0$ allows one to show that after integrating any of the
$f_{\rm bk}$ vertices the entire integral vanishes. So \\

\[
\delta \beta_{k,1,\Delta(0)}
[\vec{V}_{\rm bk}]
=0\ .
\]

We now turn to $\delta \beta_{k,1,\Delta(0)} [\vec{V}_{\rm bk} + \dot{V}]$.
From the definition we have:

\[
\delta \beta_{1,k,\Delta(0)}
[\vec{V}_{\rm bk} + \dot{V}]
=
-\sum_{b}
\bbone
\left\{ 
\begin{array}{c}
k+b\le 4 \\
b\ge 1
\end{array}
\right\}
\frac{(k+b)!}{k!\ b!}\  L^{-k[\phi]}
\qquad\qquad\parbox{2.8cm}{
\psfrag{a}{$\hspace{-1cm}\beta_{k+b,{\rm bk}}+\dot{\beta}_{k+b}$}
\psfrag{f}{$\hspace{-1cm}f_{\rm bk}+\dot{f}$}
\psfrag{g}{$f_{\rm bk}+\dot{f}$}
\psfrag{b}{$b$}
\raisebox{-1ex}{
\includegraphics[width=2.8cm]{AJfig4.eps}}
}
\]

Under the assumption that $b \ge 1$ we have:

\[
\parbox{2.8cm}{
\psfrag{a}{$\hspace{-.9cm}\beta_{k+b,{\rm bk}}+\dot{\beta}_{k+b}$}
\psfrag{f}{$\hspace{-.9cm}f_{\rm bk}+\dot{f}$}
\psfrag{g}{$f_{\rm bk}+\dot{f}$}
\psfrag{b}{$b$}
\raisebox{-1ex}{
\includegraphics[width=2.8cm]{AJfig4.eps}}
}\qquad\qquad
= 
\sum_{j=0}^{b}
\binom{b}{j}
\qquad\parbox{3cm}{
\psfrag{a}{$\hspace{-.2cm}\begin{array}{l}\ \\ \beta_{k+b,{\rm bk}} \\ +\dot{\beta}_{k+b}\end{array}$}
\psfrag{c}{$\hspace{.5cm}j$}
\psfrag{d}{$b-j$}
\psfrag{f}{$\hspace{-.2cm}f_{\rm bk}$}
\psfrag{g}{$\dot{f}$}
\raisebox{-1ex}{
\includegraphics[width=3cm]{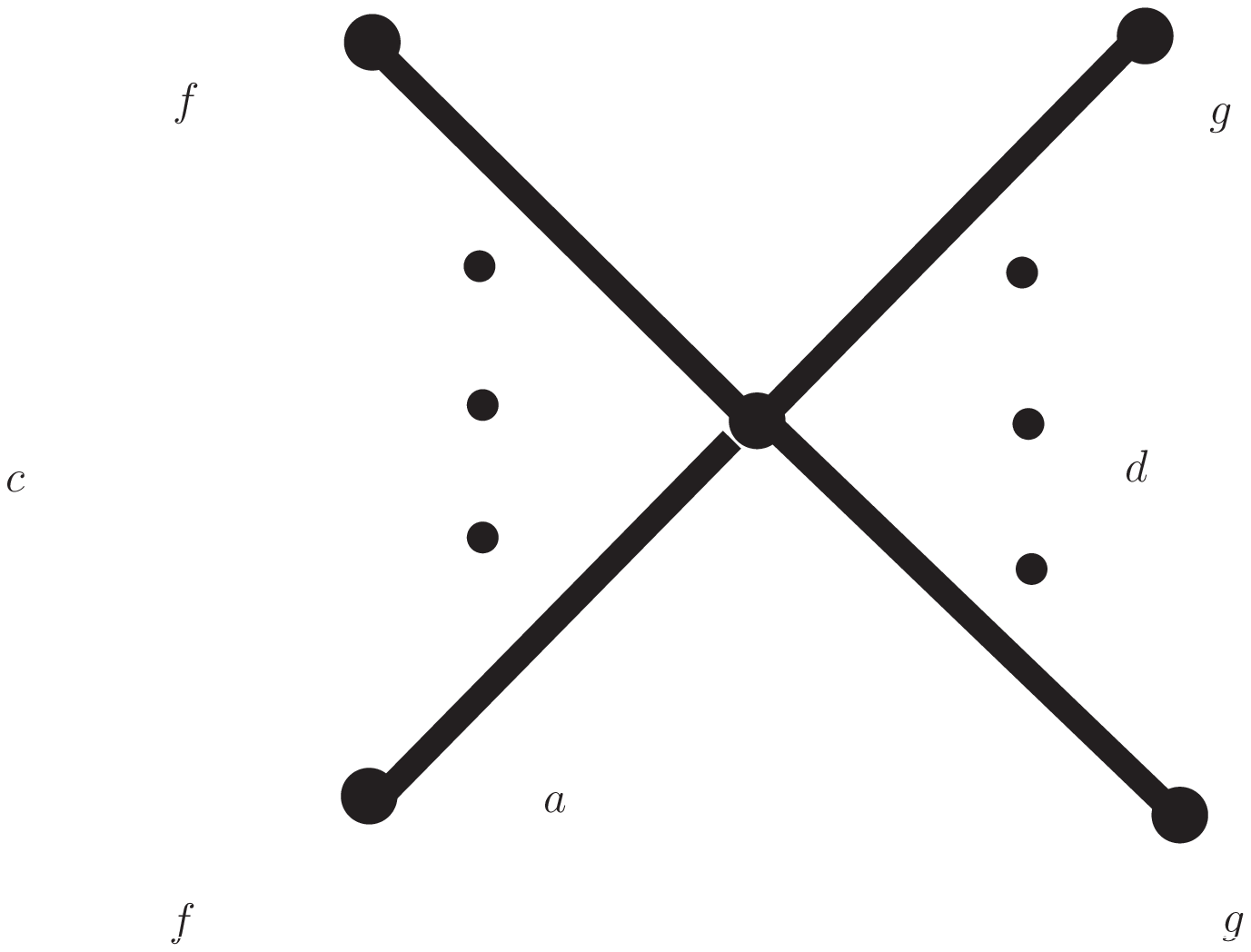}}
}
\]

In the sum above only the $j=0$ term can be non-vanishing, all other
diagrams will have at least one integration vertex of degree $1$ with
a bulk variable assigned to it. We substitute this back into our
formula for $\delta \beta_{k,1,\Delta(0)}$ and perform more manipulations:\\

\begin{align}\label{onevertexdeviationi}
\delta \beta_{k,1,\Delta(0)}
[\vec{V}_{\rm bk} + \dot{V}]
=&
-\sum_{b}
\bbone
\left\{ 
\begin{array}{c}
k+b\le 4 \\
b\ge 1
\end{array}
\right\}
\frac{(k+b)!}{k!\ b!}\  L^{-k[\phi]}
\qquad
\parbox{2.8cm}{
\psfrag{a}{$\hspace{-.9cm}\beta_{k+b,{\rm bk}}+\dot{\beta}_{k+b}$}
\psfrag{f}{$\hspace{-.1cm}\dot{f}$}
\psfrag{g}{$\dot{f}$}
\psfrag{b}{$b$}
\raisebox{-1ex}{
\includegraphics[width=2.8cm]{AJfig4.eps}}
}
\\
=& 
-(k+1) L^{-k[\phi]}
\ \ \parbox{1cm}{
\psfrag{w}{$\hspace{-.3cm}\beta_{k+1,{\rm bk}}+\dot{\beta}_{k+1}$}
\psfrag{f}{$\dot{f}$}
\raisebox{1ex}{
\includegraphics[width=1cm]{Fig2.eps}}
}\nonumber \\
 & - \sum_{b}
\bbone
\left\{ 
\begin{array}{c}
k+b\le 4 \\
b\ge 2
\end{array}
\right\}
\frac{(k+b)!}{k!\ b!}\  L^{-k[\phi]}
\qquad
\parbox{2.8cm}{
\psfrag{a}{$\hspace{-.9cm}\beta_{k+b,{\rm bk}}+\dot{\beta}_{k+b}$}
\psfrag{f}{$\hspace{-.1cm}\dot{f}$}
\psfrag{g}{$\dot{f}$}
\psfrag{b}{$b$}
\raisebox{-1ex}{
\includegraphics[width=2.8cm]{AJfig4.eps}}
}
\end{align}
where we have isolated the $b=1$ term.
Note that for $k=3$ the sum on the last line is empty. We now bound the diagrams appearing above:

\begin{equation}\label{onevertexdeviationii}
\begin{split}
\left|\ \ 
\parbox{1cm}{
\psfrag{w}{$\hspace{-.3cm}\beta_{k+1,{\rm bk}}+\dot{\beta}_{k+1}$}
\psfrag{f}{$\dot{f}$}
\raisebox{1ex}{
\includegraphics[width=1cm]{Fig2.eps}}
}\qquad\qquad
\right|
\le&
\left|\ \  
\parbox{1cm}{
\psfrag{w}{$\hspace{-.2cm}\dot{\beta}_{k+1}$}
\psfrag{f}{$\dot{f}$}
\raisebox{1ex}{
\includegraphics[width=1cm]{Fig2.eps}}
}\ \ \ 
\right|
+
\left| 
\parbox{1cm}{
\psfrag{w}{$\beta_{k+1,{\rm bk}}$}
\psfrag{f}{$\dot{f}$}
\raisebox{1ex}{
\includegraphics[width=1cm]{Fig2.eps}}
}\ \ \ 
\right|\\
=&
\left|\  
\parbox{1cm}{
\psfrag{w}{$\hspace{-.2cm}\dot{\beta}_{k+1}$}
\psfrag{f}{$\dot{f}$}
\raisebox{1ex}{
\includegraphics[width=1cm]{Fig2.eps}}
}\ \ 
\right|\\
= &
\left| 
\dot{f}(0) \times \Gamma(0) \times \dot{\beta}_{k+1,\Delta(0)}
\right|\\
\le&
 2 \left( L^{-(3-[\phi])}||\dot{V}|| \right) \left( ||\dot{V}|| \bar{g}^{e_{k+1}} \right)\\
\le& 2 L^{-\frac{9}{4}}  ||\dot{V}||^2 \bar{g}^{e_{k+1}}\ .
\end{split}
\end{equation}

In going to the third to last line we used local constancy at unit scale and the fact that all the
couplings were supported at $\Delta(0)$ so we did not really do any integration.
In going to the second to last line we used the bound $|\Gamma(0)| \le
2$ which comes from Lemma \ref{gamL00lem}. In going to the last line
we used the bound $-(3-[\phi]) \le - \frac{9}{4}$.\\

For $k=3$ we immediately have the bound:

\begin{align*}
\left| \delta \beta_{3,1,\Delta(0)} \right|
\bar{g}^{-e_{3}}
\le&
 4L^{-3[\phi]} \times 2 L^{-\frac{9}{4}}  ||\dot{V}||^2 \bar{g}^{e_{4}}
\bar{g}^{-e_{3}}\\
\le &
8 L^{-\frac{9}{4}}  ||\dot{V}||^2
\end{align*}

Note that in going to the last line we dropped the factor of $L^{-3[\phi]}$
and used $e_4\ge e_3$. This proves the lemma for the case $k=3$. We now bound the remaining
diagrams to prove the lemma for the cases $k=1$ and $k=2$. 
Before note that in these two cases $k+b=3$ or $4$ because we also assume $b\ge 2$.

If $k+b=4$ then, because of the domain hypotheses for our lemma and noting the $\bar{g}$ shift for the $\beta_4$ component
of the bulk, we must have
\[
|\beta_{k+b,{\rm bk}}|+|\dot{\beta}_{b+k}|\le
\bar{g}+\frac{1}{4}\bar{g}^{e_4} +\frac{1}{4}\bar{g}^{e_4} 
\le \frac{3}{2}\bar{g}\le \frac{3}{2}\bar{g}^{e_k}\ .
\]
This is because of our assumptions $e_1,e_2\le 1\le e_4$.

If $k+b=3$ then
\[
|\beta_{k+b,{\rm bk}}|+|\dot{\beta}_{b+k}|\le
\frac{1}{4}\bar{g}^{e_3} +\frac{1}{4}\bar{g}^{e_3}\le
\frac{3}{2}\bar{g}^{e_k}
\]
because of the assumption $e_1\le e_2\le e_3$.
So in all relevent cases we can use $\frac{3}{2}\bar{g}^{e_k}$ as a bound, as we do next.

\begin{equation}\label{onevertexdeviationiii}
\begin{split}
\bbone
\left\{ 
\begin{array}{c}
k+b\le 4 \\
b\ge 2
\end{array}
\right\}
\left|\qquad
\parbox{2.8cm}{
\psfrag{a}{$\hspace{-.9cm}\beta_{k+b,{\rm bk}}+\dot{\beta}_{k+b}$}
\psfrag{f}{$\hspace{-.1cm}\dot{f}$}
\psfrag{g}{$\dot{f}$}
\psfrag{b}{$b$}
\raisebox{-1ex}{
\includegraphics[width=2.8cm]{AJfig4.eps}}
}
\qquad
\right|
\le &
\bbone
\left\{ 
\begin{array}{c}
k+b\le 4 \\
b\ge 2
\end{array}
\right\}
\times
\left| \dot{f}(0) \right|^{b}\\
 & \times
\left(
\left|\beta_{b+k,{\rm bk}}\right|
+ 
\left|\dot{\beta}_{b+k} \right|
\right)
\times
\int_{\mathbb{Q}_{p}^{3}} {\rm d}^3x\ \left| \Gamma(x) \right|^b\\
\le &
\left( L^{-(3-[\phi])} ||\dot{V}|| \right)^2 
\times 
\frac{3}{2}\bar{g}^{e_k}
\times
2^{5/2} L^{3-2[\phi]}\\
\le&
3\times2^{3/2}L^{-3} || \dot{V} ||^2 \bar{g}^{e_k}\ .
\end{split}
\end{equation}

For the bound on the first line we used the fact that all the
$\dot{f}$ vertices are pinned to the origin and the only integration
occurs at the $\beta_{b+k,{\rm bk}} + \dot{\beta}_{b+k}$ vertex which has
been left with $b$ copies of the fluctuation covariance. \\

In going to the second to last line we used the bound $|\dot{f}(0)|^b \le |\dot{f}(0)|^2$
since $b \ge 2$ and $|\dot{f}(0)|\le 1$. For that same line also used the following
bound which is valid for $2 \le b \le 4$:

\begin{align*}
\int_{\mathbb{Q}_{p}^{3}} {\rm d}^3x\ |\Gamma(x)|^{b}
\le& 
||\Gamma||_{L^{1}} ||\Gamma||_{L^{\infty}}^{b-1}\\
\le&
\left( \frac{1}{\sqrt{2}} L^{3-2[\phi]} \right) 2^{b-1}\\
\le& 2^{5/2} L^{3-2[\phi]}\ .
\end{align*}

Note that we have used fluctuation covariance bounds of Corollary
\ref{gamL1cor} and Lemma \ref{gamL00lem}. Thus we
can use \eqref{onevertexdeviationii} to get the following bound
for $k=1$ and $k=2$:

\begin{equation}\label{onevertexdeviationiv}
\begin{split}
\left|
\sum_{b}
\bbone
\left\{ 
\begin{array}{c}
k+b\le 4 \\
b\ge 2
\end{array}
\right\}
\frac{(k+b)!}{k!\ b!}\  L^{-k[\phi]}
\qquad
\parbox{2.8cm}{
\psfrag{a}{$\hspace{-.9cm}\beta_{k+b,{\rm bk}}+\dot{\beta}_{k+b}$}
\psfrag{f}{$\hspace{-.1cm}\dot{f}$}
\psfrag{g}{$\dot{f}$}
\psfrag{b}{$b$}
\raisebox{-1ex}{
\includegraphics[width=2.8cm]{AJfig4.eps}}
}\qquad
\right|
\le&
\sum_{b}
\bbone
\left\{ 
\begin{array}{c}
k+b\le 4 \\
b\ge 2
\end{array}
\right\}
\frac{(k+b)!}{k!\ b!} \\
 & \times
3\times2^{3/2}L^{-3} || \dot{V} ||^2 \bar{g}^{e_k}\\
\le&
21 \times 2^{3/2}L^{-3} || \dot{V} ||^2 \bar{g}^{e_k}\ .
\end{split}
\end{equation}

Note in going to the last line we dropped the factors of
$L^{-k[\phi]}$ and used that

\[
\max_{k=1,2} \sum_{b}
\bbone
\left\{ 
\begin{array}{c}
k+b\le 4 \\
b\ge 2
\end{array}
\right\}\ \frac{(k+b)!}{k!\ b!}
=
7\ .
\]

Finally by inserting the bound \eqref{onevertexdeviationii} and
\eqref{onevertexdeviationiv} into \eqref{onevertexdeviationi} we get
the following bound for $k=1$ and $k=2$:

\begin{equation*}
\begin{split}
\left|
\delta \beta_{k,1,\Delta(0)}[\vec{V}_{\rm bk}+\dot{V}]
\right|
\bar{g}^{-e_{k}}
\le& (k+1)\times
2 L^{-\frac{9}{4}}  ||\dot{V}||^2 +
21\times2^{3/2}L^{-3} || \dot{V} ||^2 \\
\le& \left(6 + 21 \times 2^{\frac{3}{2}} \right)L^{-\frac{9}{4}} ||\dot{V} ||^2\ .
\end{split}
\end{equation*}

In going to the last line we simply bounded $L^{-3}$ by
$L^{-\frac{9}{4}}$.  This proves the lemma
  for $k=1$ and $k=2$ which finishes the proof. \qed\\

\begin{Lemma} \label{LAJ9}
For $k=1,2,3,4$ and for all $\Delta' \in \mathbb{L}$ one has that
$\delta b_{k,2,\Delta'}[\bullet]$ and $\xi_{k,\Delta'}[\bullet]$ are analytic functions on
$B( \bar{V}, \frac{1}{2})$ taking values in $\mathbb{C}$. In
particular one has the following bounds for any $\vec{V}^{1},
\vec{V}^{2} \in \bar{B}( \bar{V} ,\frac{1}{6} )$

\begin{equation}\label{LAJ9i}
\left|
\delta \beta_{k,2,\Delta'} \left[ \vec{V}^{1} \right] 
- 
\delta \beta_{k,2,\Delta'} \left[ \vec{V}^{2} \right] 
\right|
\bar{g}^{-e_{k}}
\le
\frac{1}{100}
||\vec{V}^{1} - \vec{V}^{2}||,
\end{equation}

\begin{equation}\label{LAJ9ii}
\left|
\xi_{k,\Delta'} \left[ \vec{V}^{1} \right] 
- 
\xi_{k,\Delta'} \left[ \vec{V}^{2} \right] 
\right|
\bar{g}^{-e_{k}}
\le
\frac{1}{100}
||\vec{V}^{1} - \vec{V}^{2}||.
\end{equation}

\end{Lemma}

\noindent{\bf Proof:}
The statement of analyticity for $\delta \beta_{k,2,\Delta'}$ immediate from the formulas that
define $\delta \beta_{k,2,\Delta'}$. To establish the bound
\eqref{LAJ9i} we first use Lemma \ref{L34lem} which gives us the
following uniform bound for all $\vec{V} \in B(\bar{V},\frac{1}{2})$:

\[
\left|
\delta \beta_{k,2,\Delta'} \left[ \vec{V} \right] 
\right|
\le
\mathcal{O}_{26}L^5
\bar{g}^{2-2\eta}\ .
\]

Thus by applying Lemma \ref{Lipschitzlem} with $\nu =
\frac{1}{3}$ we have:

\[
\left|
\delta \beta_{k,2,\Delta'} \left[ \vec{V}^{1} \right] 
-
\delta \beta_{k,2,\Delta'} \left[ \vec{V}^{2} \right] 
\right|
\bar{g}^{-e_{k}}
\le
4
\mathcal{O}_{26}L^5
\bar{g}^{2-2\eta - e_{k}}
||\vec{V}^{1} - \vec{V}^{2}||\ .
\]

Note that we have $2-2\eta - e_{k} > 0$ as a consequence of
\eqref{econstbetaeq}, thus by making epsilon sufficiently small we can
guarantee $2 \mathcal{O}_{26}L^{5} \bar{g}^{2-2\eta - e_{k}} \le
\frac{1}{100}$ which proves \eqref{LAJ9i}.\\

For $\xi_{k,\Delta'}$ we have both analyticity and the following
uniform bound for $\vec{V} \in B(\bar{V},\frac{1}{2})$  as consequences of Theorem \ref{mainestthm}:

\[
\left|
\xi_{k,\Delta'}[\vec{V}]
\right|
\le
B_{k}
\times \frac{1}{2}
\bar{g}^{e_{R}}\ .
\]

We again use Lemma \ref{Lipschitzlem} with $\nu =
\frac{1}{3}$ to get:

\[
\left|
\xi_{k,\Delta'}[\vec{V}^{1}]
-
\xi_{k,\Delta'}[\vec{V}^{1}]
\right|\bar{g}^{-e_{k}}
\le
2 B_{k}
\bar{g}^{e_{R} - e_{k}}
||\vec{V}^{1} - \vec{V}^{2}||\ .
\]

Note that $e_{R} > e_{k}$ because of the assumptions \eqref{econst4eq}
and \eqref{econstbetaeq}, thus by requiring  that $\epsilon$ be sufficiently small
we can guarantee $2 B_{k} \bar{g}^{e_{R} - e_{k}} \le \frac{1}{100}$
which proves \eqref{LAJ9ii}. \qed\\

Given $\vec{V}_{\rm bk} \in B(\bar{V},\frac{1}{4}) \cap \mathcal{E}_{\rm bk}$ and $\dot{V}
\in B(0,\frac{1}{4}) \cap \mathcal{E}_{\rm pt}$ we define:

\[
RG_{\rm dv} [\vec{V}_{\rm bk}, \dot{V} ]
=
RG_{\rm ex}[ \vec{V}_{\rm bk} + \dot{V}] 
- 
RG_{\rm ex}[\vec{V}_{\rm bk}].
\]

Note that, as a subspace of $\mathcal{E_{\rm ex}}$, the space
$\mathcal{E}_{\rm bk} \oplus \mathcal{E}_{\rm pt}$ is invariant by
$RG_{\rm ex}$. Since $\vec{V}_{\rm bk} + \dot{V} \in \mathcal{E}_{\rm bk} \oplus
\mathcal{E}_{\rm pt}$ one has a unique decomposition $RG_{\rm ex} [
\vec{V}_{\rm bk} + \dot{V}] = \vec{V}_{\rm bk}'+\dot{V}'$ with $\vec{V}_{\rm bk}'
\in \mathcal{E}_{\rm bk}$ and $\dot{V}' \in \mathcal{E}_{\rm pt}$. Using the
locality of $RG_{\rm ex}$ it is not hard to see that
$\vec{V}_{\rm bk}' = RG_{\rm ex} [\vec{V}_{\rm bk}]$ and $\dot{V}'= RG_{\rm dv} [\vec{V}_{\rm bk},
\dot{V} ]$. In particular $RG_{\rm ex}[\bullet, \bullet]$ takes values in
$\mathcal{E}_{\rm pt}$. \\

\begin{Lemma}\label{LAJ10}
Suppose that $\vec{V}_{\rm bk} \in B(\bar{V},\frac{1}{12}) \cap \mathcal{E}_{\rm bk}$ and $\dot{V}
\in B(0,\frac{1}{12}) \cap \mathcal{E}_{\rm pt}$. Define
$\dot{V}' =  RG_{\rm dv} [\vec{V}_{\rm bk}, \dot{V} ]$ and for $k=1,2,3,4$ let
$\dot{\beta}'_{k}$ be the corresponding components of $\dot{V}'$. \\

We then have the following bound for $k=1,2,3,4$

\[
\left| \dot{\beta}'_{k,\Delta(0)} \right|
\bar{g}^{-e_{k}}
\le
\frac{4}{5}
|| \dot{V} ||
+
\mathbf{O}_{7}
||\dot{V}||^2
\]
where $\mathbf{O}_{7}$ has been defined in Lemma \ref{LAJ8}.
\end{Lemma}

\noindent{\bf Proof:}
By definition we have:

\begin{equation} \label{deviationrelevant1}
\begin{split}
\dot{\beta}'_{k,\Delta(0)} =& L^{-k[\phi]} \sum_{\Delta \in[L^{-1}\Delta(0)]} \dot{\beta}_{k,\Delta}  \\
&\ + \left( \delta \beta_{k,1,\Delta(0)} \left[ \vec{V}_{\rm bk} \right]-
  \delta \beta_{k,1,\Delta(0)} \left[ \vec{V}_{\rm bk} + \dot{V} \right] \right)\\
&\ + \left( \delta \beta_{k,2,\Delta(0)} \left[ \vec{V}_{\rm bk} \right] -
  \delta \beta_{k,2,\Delta(0)} \left[ \vec{V}_{\rm bk} + \dot{V} \right] \right)\\
&\ + \left( \xi_{k,\Delta(0)} \left[ \vec{V}_{\rm bk} + \dot{V} \right] - \xi_{k,\Delta(0)} \left[ \vec{V}_{\rm bk}\right] \right)\ .
\end{split}
\end{equation}\\

Since $\dot{\beta}_{k,\Delta}$ is supported on $\Delta = \Delta(0)$
the sum on the first line only has one non-zero term. We then
have the bound:

\begin{align*}
\left| 
L^{-k[\phi]} \sum_{\Delta \in[L^{-1}\Delta(0)]} \dot{\beta}_{k,\Delta}
\right|
\le& 
L^{-k[\phi]} ||\dot{V}||\bar{g}^{e_{k}}\\
\le& 2^{-\frac{1}{2}} ||\dot{V}||\bar{g}^{e_{k}}\ .
\end{align*}

In going to the last line we are assuming that $L \ge 2$ and $\epsilon
\le 1$ so that we have $L^{-k[\phi]} \le L^{-[\phi]} \le
  2^{-\frac{1}{2}}$.  We then combine the estimate above
  with the lemmas \ref{LAJ8}, \ref{LAJ9},
and \ref{LAJ10} to get the bound:

\[
\left| \dot{\beta}'_{k,\Delta(0)} \right| \bar{g}^{-e_{k}}
\le  
2^{-\frac{1}{2}} ||\dot{V}|| 
+\mathbf{O}_{7} \bbone \{  1 \le k \le 3\} L^{-\frac{9}{4}} ||\dot{V}||^2
+\frac{1}{100} ||\dot{V}||
+\frac{1}{100} ||\dot{V}||\ .
\]

We get the bound of this lemma by dropping the factor of
$L^{-\frac{9}{4}}$ and the indiciator function in the fourth term
while also observing that
  $2^{-\frac{1}{2}} + \frac{1}{100} + \frac{1}{100} < \frac{4}{5}$.
  \qed\\

\begin{Lemma}\label{LAJ11}

Suppose that $\vec{V}_{\rm bk} \in B(\bar{V},\frac{1}{4}) \cap \mathcal{E}_{\rm bk}$ and $\dot{V}
\in B(0,\frac{1}{4}) \cap \mathcal{E}_{\rm pt}$. Let
$\dot{V}' =  RG_{\rm dv} [\vec{V}_{\rm bk}, \dot{V} ]$ and for $k=5,6$ let
$\dot{W}'_{k}$ be the corresponding components of $\dot{V}'$. \\

We then have the following bound for $k=5,6$

\[
\left|\dot{W}'_{k,\Delta(0)} \right| \le 2^{-\frac{5}{2}} ||\dot{V}||
+ \mathbf{O}_{8} ||\dot{V}||^2,
\]
where $\mathbf{O}_{8} = \left(18+\frac{9}{\sqrt{2}} \right)$.

\end{Lemma}

\noindent{\bf Proof:}
For $k=5$ we have:

\begin{equation}\label{W5deviationflow}
\begin{split}
 \dot{W}'_{5,\Delta(0)} =& L^{-5[\phi]}  \sum_{\Delta \in[L^{-1}\Delta(0)]}
 \dot{W}_{5,\Delta}\\
&\ + 48 L^{-5[\phi]} \left(
\qquad\parbox{2.1cm}{
\psfrag{a}{$\hspace{-.7cm}\beta_{4,{\rm bk}}+\dot{\beta}_4$}
\psfrag{b}{$\beta_{4,{\rm bk}}+\dot{\beta}_4$}
\psfrag{f}{$f_{\rm bk}+\dot{f}$}
\raisebox{1ex}{
\includegraphics[width=2.1cm]{AJfig1.eps}}
}\qquad\qquad
-
\parbox{2.1cm}{
\psfrag{a}{$\beta_{4,{\rm bk}}$}
\psfrag{b}{$\beta_{4,{\rm bk}}$}
\psfrag{f}{$f_{\rm bk}$}
\raisebox{1ex}{
\includegraphics[width=2.1cm]{AJfig1.eps}}
}\qquad
\right)\\
\ &\ \\
 & \ + 6 L^{-5[\phi]} \left(
\qquad\parbox{1cm}{
\psfrag{w}{$\hspace{-.7cm}W_{6,{\rm bk}}+\dot{W}_6$}\psfrag{f}{$f_{\rm bk}+\dot{f}$}
\raisebox{1ex}{
\includegraphics[width=1cm]{Fig2.eps}}
}\qquad\qquad
-
\qquad\parbox{1cm}{
\psfrag{w}{$W_{6,\rm bk}$}
\psfrag{f}{$f_{\rm bk}$}
\raisebox{1ex}{
\includegraphics[width=1cm]{Fig2.eps}}
}\qquad\qquad
\right)\\
\ &\ \\
&\ + 12 L^{-5[\phi]} \left(
\qquad\parbox{1cm}{
\psfrag{w}{$\hspace{-.7cm}\beta_{4,{\rm bk}}+\dot{\beta}_4$}
\psfrag{f}{$\beta_{3,{\rm bk}}+\dot{\beta}_3$}
\raisebox{1ex}{
\includegraphics[width=1cm]{Fig2.eps}}
}\qquad\qquad
-
\parbox{1cm}{
\psfrag{w}{$\beta_{4,{\rm bk}}$}
\psfrag{f}{$\beta_{3,{\rm bk}}$}
\raisebox{1ex}{
\includegraphics[width=1cm]{Fig2.eps}}
}\qquad
\right)\ .
\end{split}
\end{equation}

As before using that $\dot{W}_{5,\Delta}$ is supported on $\Delta =
\Delta(0)$ gives us the bound:

\[
\left|L^{-5[\phi]}  \sum_{\Delta \in[L^{-1}\Delta(0)]}
  \dot{W}_{5,\Delta} \right|
\le
L^{-5[\phi]} \bar{g}^{e_{W}}||\dot{V}||\ .
\]

We now bound the various graphs appearing in \eqref{W5deviationflow}. We again note that when a graph has an integration vertex of degree
one that has been assigned a bulk variable the graph will vanish. This
tells us that:
\[
\parbox{2.1cm}{
\psfrag{a}{$\beta_{4,{\rm bk}}$}
\psfrag{b}{$\beta_{4,{\rm bk}}$}
\psfrag{f}{$f_{\rm bk}$}
\raisebox{1ex}{
\includegraphics[width=2.1cm]{AJfig1.eps}}
}\qquad
=\qquad
\parbox{1cm}{
\psfrag{w}{$\hspace{-.5cm}\beta_{4,{\rm bk}}$}
\psfrag{f}{$\beta_{3,{\rm bk}}$}
\raisebox{1ex}{
\includegraphics[width=1cm]{Fig2.eps}}
}\qquad\qquad
=\ 
\parbox{1cm}{
\psfrag{w}{$W_{6,{\rm bk}}$}
\psfrag{f}{$f_{\rm bk}$}
\raisebox{1ex}{
\includegraphics[width=1cm]{Fig2.eps}}
}\qquad\qquad
=0.
\]

We use this same observation to break up the non-vanishing graphs and show
that their contribution is second order in $||\dot{V}||$. For example:

\begin{equation*}
\begin{split}
\parbox{2.1cm}{
\psfrag{a}{$\hspace{-.7cm}\beta_{4,{\rm bk}}+\dot{\beta}_4$}
\psfrag{b}{$\beta_{4,{\rm bk}}+\dot{\beta}_4$}
\psfrag{f}{$f_{\rm bk}+\dot{f}$}
\raisebox{1ex}{
\includegraphics[width=2.1cm]{AJfig1.eps}}
}\qquad\qquad
= &
\ \parbox{2.1cm}{
\psfrag{a}{$\hspace{-.2cm}\dot{\beta}_4$}
\psfrag{b}{$\beta_{4,{\rm bk}}+\dot{\beta}_4$}
\psfrag{f}{$\dot{f}$}
\raisebox{1ex}{
\includegraphics[width=2.1cm]{AJfig1.eps}}
}\qquad\qquad
+
\parbox{2.1cm}{
\psfrag{a}{$\beta_{4,{\rm bk}}$}
\psfrag{b}{$\beta_{4,{\rm bk}}+\dot{\beta}_4$}
\psfrag{f}{$f_{\rm bk}$}
\raisebox{1ex}{
\includegraphics[width=2.1cm]{AJfig1.eps}}
} \\
\ & \ \\
\ & +
\parbox{2.1cm}{
\psfrag{a}{$\beta_{4,{\rm bk}}$}
\psfrag{b}{$\beta_{4,{\rm bk}}+\dot{\beta}_4$}
\psfrag{f}{$\dot{f}$}
\raisebox{1ex}{
\includegraphics[width=2.1cm]{AJfig1.eps}}
}\qquad\qquad
+
\parbox{2.1cm}{
\psfrag{a}{$\hspace{-.2cm}\dot{\beta}_4$}
\psfrag{b}{$\beta_{4,{\rm bk}}+\dot{\beta}_4$}
\psfrag{f}{$f_{\rm bk}$}
\raisebox{1ex}{
\includegraphics[width=2.1cm]{AJfig1.eps}}
}
\\
\ & \ \\
=&
\parbox{2.1cm}{
\psfrag{a}{$\hspace{-.2cm}\dot{\beta}_4$}
\psfrag{b}{$\beta_{4,{\rm bk}}+\dot{\beta}_4$}
\psfrag{f}{$\dot{f}$}
\raisebox{1ex}{
\includegraphics[width=2.1cm]{AJfig1.eps}}
}
\end{split}
\end{equation*}
after expanding the two outer vertices of valence one.

We then have

\begin{equation}\label{W5deviationgraphboundi}
\begin{split}
\left|\qquad
\parbox{2.1cm}{
\psfrag{a}{$\hspace{-.7cm}\beta_{4,{\rm bk}}+\dot{\beta}_4$}
\psfrag{b}{$\beta_{4,{\rm bk}}+\dot{\beta}_4$}
\psfrag{f}{$f_{\rm bk}+\dot{f}$}
\raisebox{1ex}{
\includegraphics[width=2.1cm]{AJfig1.eps}}
}\qquad
\right| = &
\left|\ \ 
\parbox{2.1cm}{
\psfrag{a}{$\hspace{-.2cm}\dot{\beta}_4$}
\psfrag{b}{$\beta_{4,{\rm bk}}+\dot{\beta}_4$}
\psfrag{f}{$\dot{f}$}
\raisebox{1ex}{
\includegraphics[width=2.1cm]{AJfig1.eps}}
}\qquad\qquad
\right|\\
\le& \left|\dot{f}(0) \right| \times \left| \dot{\beta}_{4,\Delta(0)} \right|
\times \left( \left|\beta_{4,{\rm bk}} \right| + \left| \dot{\beta}_{4,\Delta(0)}
\right| \right) \times \int_{\mathbb{Q}_{p}^{3}} {\rm d}^3x\ \left|
\Gamma(x)\right|^{2}\\
\le&L^{-(3-[\phi])} ||\dot{V}|| \times
||\dot{V}|| \bar{g}^{e_{4}} \times \frac{3}{2}\bar{g} \times \left(
  2^{\frac{1}{2}} L^{3-2[\phi]} \right)\\
\le& 3\times 2^{-\frac{1}{2}}  L^{-[\phi]} \bar{g}^{1+e_{4}} ||\dot{V}||^2\\
\le& 3\times 2^{-\frac{1}{2}} \bar{g}^{1+e_{4}} ||\dot{V}||^2\ .
\end{split}
\end{equation}

Note that in going to the second to last line we again used the bound:

\[
\int_{\mathbb{Q}_{p}^{3}} {\rm d}^3x \left| \Gamma(x) \right|^{n} \le 2^{n
- \frac{3}{2}} L^{3-2[\phi]}\ .
\]

Proceeding similarly for the other graphs we have:

\begin{equation}\label{W5deviationgraphboundii}
\begin{split}
\left|\qquad
\parbox{1cm}{
\psfrag{w}{$\hspace{-.7cm}W_{6,{\rm bk}}+\dot{W}_6$}\psfrag{f}{$f_{\rm bk}+\dot{f}$}
\raisebox{1ex}{
\includegraphics[width=1cm]{Fig2.eps}}
}\qquad\qquad
\right|
=&
\left|\ \ \   
\parbox{1cm}{
\psfrag{w}{$\hspace{-.3cm}\dot{W}_6$}\psfrag{f}{$\dot{f}$}
\raisebox{1ex}{
\includegraphics[width=1cm]{Fig2.eps}}
}\ \ 
\right|
= \left|\dot{f}(0) \right| \times
\left|\dot{W}_{6,\Delta(0)} \right| \times \left| \Gamma(0) \right|\\
\le & L^{-(3-[\phi])} \bar{g}^{e_{W}} \left| \Gamma(0) \right| \times
||\dot{V}||^2 \\
\le &2 \bar{g}^{e_{W}} ||\dot{V}||^2\ .
\end{split}
\end{equation}

In going to the last line we used the bound $|\Gamma(0)| \le 2$ which
is a consequence of Corollary \ref{gamL1cor}. We also dropped the
factor of $L^{-(3-[\phi])}\le L^{-\frac{9}{4}}\le 1$. We continue to the last graph we need to
bound for $\dot{W}_{5}'$:

\begin{equation}\label{W5deviationgraphboundiii}
\begin{split}
\left|\qquad
\parbox{1cm}{
\psfrag{w}{$\hspace{-.7cm}\beta_{4,{\rm bk}}+\dot{\beta}_4$}
\psfrag{f}{$\beta_{3,{\rm bk}}+\dot{\beta}_3$}
\raisebox{1ex}{
\includegraphics[width=1cm]{Fig2.eps}}
}\qquad\qquad\ \ \ 
\right|
\le&
\left|\ \ 
\parbox{1cm}{
\psfrag{w}{$\hspace{-.2cm}\dot{\beta}_4$}
\psfrag{f}{$\dot{\beta}_3$}
\raisebox{1ex}{
\includegraphics[width=1cm]{Fig2.eps}}
}\ \ \ \ 
\right|\\
= &
\left| \dot{\beta}_{4,\Delta(0)} \right|
\times
\left|\dot{\beta}_{3,\Delta(0)}\right|
\times
\left| \Gamma(0) \right|\\
\le&
2||\dot{V} ||^2 \bar{g}^{e_{4}+e_{3}}\ .
\end{split}
\end{equation}

Using the
bounds \eqref{W5deviationgraphboundi}, \eqref{W5deviationgraphboundii},
and \eqref{W5deviationgraphboundiii} in \eqref{W5deviationflow} gives
us the bound:

\begin{equation*}
\begin{split}
\left|\dot{W}_{5,\Delta(0)} \right| \bar{g}^{-e_{W}}
\le&
L^{-5[\phi]} ||\dot{V}|| + L^{-5[\phi]} \left[48\times 3\times 2^{-\frac{1}{2}}  \bar{g}^{1+e_{4} -
  e_{W}} + 6 \times 2 + 12
\times 2 \bar{g}^{e_{4}+e_{3} - e_{W}} \right]||\dot{V}||^2\\
\le&
2^{-\frac{5}{2}} ||\dot{V}|| + 2^{-\frac{5}{2}}  \left[ 48\times 3\times 2^{-\frac{1}{2}}+ 6 \times 2 + 12
\times 2 \right] ||\dot{V}||^2\\
=& 2^{-\frac{5}{2}} ||\dot{V}|| + \left(18+\frac{9}{\sqrt{2}} \right) ||\dot{V}||^2\ .
\end{split}
\end{equation*}

In going to the second line we used the fact that $\epsilon \le 1$ and
$L \ge 2$ to bound $L^{-5[\phi]} \le 2^{-\frac{5}{2}}$. We also used that
$1+e_{4} - e_{W} \ge 0$ and $1 + e_{4}-e_{W} \ge 0$, these are both
consequences of \eqref{econstWeq}. This proves the lemma for $k=5$.

For $k=6$ we have:

\begin{equation}\label{W6deviationflow}
\begin{split}
\dot{W}'_{6,\Delta(0)}
=&
L^{-6[\phi]} \sum_{\Delta \in [L^{-1}\Delta(0)]} \dot{W}_{6,\Delta}\\
&\ +
8L^{-6[\phi]} \left(
\qquad\parbox{1cm}{
\psfrag{w}{$\hspace{-.7cm}\beta_{4,{\rm bk}}+\dot{\beta}_4$}
\psfrag{f}{$\beta_{4,{\rm bk}}+\dot{\beta}_4$}
\raisebox{1ex}{
\includegraphics[width=1cm]{Fig2.eps}}
}\qquad\qquad
-\ \ 
\parbox{1cm}{
\psfrag{w}{$\hspace{-.2cm}\beta_{4,{\rm bk}}$}
\psfrag{f}{$\beta_{4,{\rm bk}}$}
\raisebox{1ex}{
\includegraphics[width=1cm]{Fig2.eps}}
}\qquad\qquad
\right)\ .
\end{split}
\end{equation}

Proceeding as last time we see:

\[
\left| 
L^{-6[\phi]} \sum_{\Delta \in [L^{-1}\Delta(0)]} \dot{W}_{6,\Delta}
\right|
\le
L^{-6[\phi]}
||\dot{V}|| 
\bar{g}^{e_{W}}
\]
and

\[ 
\parbox{1cm}{
\psfrag{w}{$\beta_{4,{\rm bk}}$}
\psfrag{f}{$\beta_{4,{\rm bk}}$}
\raisebox{1ex}{
\includegraphics[width=1cm]{Fig2.eps}}
}\qquad
= 0,
\qquad\qquad
\parbox{1cm}{
\psfrag{w}{$\hspace{-.7cm}\beta_{4,{\rm bk}}+\dot{\beta}_4$}
\psfrag{f}{$\beta_{4,{\rm bk}}+\dot{\beta}_4$}
\raisebox{1ex}{
\includegraphics[width=1cm]{Fig2.eps}}
}\qquad\qquad
=
\ \parbox{1cm}{
\psfrag{w}{$\hspace{-.2cm}\dot{\beta}_4$}
\psfrag{f}{$\dot{\beta}_4$}
\raisebox{1ex}{
\includegraphics[width=1cm]{Fig2.eps}}
}
\]
which simplifies the right-hand side of (\ref{W6deviationflow}).
We now bound the contributing graph:

\begin{equation}
\begin{split}
\label{W6deviationgraphi}
\left|\ \ 
\parbox{1cm}{
\psfrag{w}{$\hspace{-.2cm}\dot{\beta}_4$}
\psfrag{f}{$\dot{\beta}_4$}
\raisebox{1ex}{
\includegraphics[width=1cm]{Fig2.eps}}
}\ \ \ 
\right|
= &
|\dot{\beta}_{4,\Delta(0)}|^2 \times |\Gamma(0)|\\
\le& 2 ||\dot{V}||^2 \bar{g}^{2e_{4}}\ .
\end{split}
\end{equation}

Inserting \eqref{W6deviationgraphi} along with the our earlier bound
into \eqref{W6deviationflow} gives us:

\begin{equation*}
\begin{split}
\left| \dot{W}'_{6,\Delta(0)} \right|e^{-e_{W}}
\le&
L^{-6[\phi]} ||\dot{V}|| + 8 L^{-6[\phi]} \times 2 ||\dot{V}||^2
\bar{g}^{2e_{4}-e_{W}}\\
\le& 2^{-3} ||\dot{V}|| + 2 ||\dot{V}||^2\ . 
\end{split}
\end{equation*}

In going to the last line we used our assumption that $\epsilon \le 1$
and $L \ge 2 $ to bound $L^{-6[\phi]} \le 2^{-3}$. We also used that
$2e_{4}-e_{W} \ge 0$ which is a consequence of \eqref{econstWeq} and
\eqref{econstbetaeq}. This proves the bound of our lemma for the case $k=6$ which finishes
the proof. \qed

\begin{Lemma}\label{LAJ12}
For any $\Delta' \in \mathbb{L}$ let $R'_{\Delta'}[\vec{V}]$ be the
corresponding component of $RG_{\rm ex}[\vec{V}]$.\\

 Let $\vec{V}^{1}, \vec{V}^{2} \in
\bar{B}(\bar{V},\frac{1}{20})$. Then one has the following bound:

\[
|||R'_{\Delta'}[\vec{V}^{1}]  - R'_{\Delta}[\vec{V}^{2}]|||_{\bar{g}}\ \bar{g}^{-e_{R}}
\le
\frac{27}{32}
||\vec{V}^{1} - \vec{V}^{2}||\ .
\]
\end{Lemma}

\noindent{\bf Proof:}
By Theorem \ref{mainestthm} and Lemma \ref{LAJ3}
we have that $R'_{\Delta'}[ \bullet ]$ is
an analytic function from $B(\bar{V},\frac{1}{2}) \subset
\mathcal{E}_{\rm ex}$ into $\bar{B}(0,\frac{3}{8} \bar{g}^{e_{R}})
\subset C^{9}_{\rm bd}(\mathbb{R},\mathbb{C})$ where we are using the norm
$|||\bullet|||_{\bar{g}}$ on $C^{9}_{\rm bd}(\mathbb{R},\mathbb{C})$. One can
then use Lemma \ref{Lipschitzlem} with $\nu = \frac{1}{10}$ to get the
bound:

\begin{equation*}
\begin{split}
|||R'_{\Delta'}[\vec{V}^{1}]  - R'_{\Delta}[\vec{V}^{2}]|||_{\bar{g}}
\le&
\frac{\frac{3}{8}\bar{g}^{e_{R}}
\left(1-\frac{1}{10}\right)}{\frac{1}{2}
\left(1-\frac{2}{10}\right)} ||\vec{V}^{1} - \vec{V}^{2}||\\
=&\frac{27}{32} \bar{g}^{e_{R}}||\vec{V}^{1} - \vec{V}^{2}||\ .
\end{split}
\end{equation*}

This proves the lemma. \qed

\begin{Lemma}\label{LAJ13}

Suppose that $\vec{V}_{\rm bk} \in \bar{B}(\bar{V},\frac{1}{40}) \cap \mathcal{E}_{\rm bk}$ and $\dot{V}
\in \bar{B}(0,\frac{1}{40}) \cap \mathcal{E}_{\rm pt}$. Let
$\dot{V}' =  RG_{\rm dv} [\vec{V}_{\rm bk}, \dot{V} ]$. Then one has the
following bound:

\[
||\dot{V}'|| \le \frac{27}{32} ||\dot{V}|| + \mathbf{O}_{9}
||\dot{V}||^2
\]
where $\mathbf{O}_{9}=\max(\mathbf{O}_{7},\mathbf{O}_{8})$.
\end{Lemma}

\noindent{\bf Proof:}
Note that $\dot{V}'_{\Delta'}$ is supported on $\Delta'=\Delta(0)$. The necessary estimates for the $\dot{\beta}'$ and $\dot{W}'$
components of $\dot{V}'$ come from Lemmas \ref{LAJ10} and
\ref{LAJ11}. For the $R$ bound we note that by Lemma \ref{LAJ12} one
has:

\begin{equation*}
\begin{split}
|||\dot{R}'_{\Delta(0)}|||_{\bar{g}}\bar{g}^{-e_{R}}
 & =
|||R'_{\Delta(0)}[\vec{V}_{\rm bk} + \dot{V}] -
R'_{\Delta(0)}[\vec{V}_{\rm bk}]|||_{\bar{g}}\bar{g}^{-e_{R}}\\
 & \le \frac{27}{32} ||\dot{V}||\ .
\end{split}
\end{equation*}

The last component we must estimate is $\dot{f}'$ which can be done
easily.

\begin{equation*}
\begin{split}
|\dot{f}'_{\Delta(0)}|L^{-(3-[\phi])}
=&
L^{-[\phi]} \left|\ \sum_{\Delta \in [L^{-1}\Delta(0)]} \dot{f}_{\Delta}\ \right|
L^{-(3-[\phi])}\\
=& L^{-[\phi]} |\dot{f}_{\Delta(0)}|L^{-(3-[\phi])}\\
\le & L^{-[\phi]} ||\dot{V}||\\
\le & \frac{27}{32} ||\dot{V}||\ .
\end{split}
\end{equation*}

In going to the second line we used the fact that $\dot{f}_{\Delta}$ is
supported on $\Delta = \Delta(0)$. In going to the last line our
assumptions that $\epsilon \le 1$ and $L \ge 2$ give us that
$L^{-[\phi]} \le 2^{-\frac{1}{2}} < \frac{27}{32}$. The lemma is then
proved.

\qed

\begin{Proposition}\label{contractionbound}

Suppose that $\vec{V}_{\rm bk} \in \bar{B}(\bar{V},\frac{1}{40}) \cap \mathcal{E}_{\rm bk}$ and $\dot{V}
\in \bar{B}(0,\mathbf{O}_{10}) \cap \mathcal{E}_{\rm pt}$ where
  $\mathbf{O}_{10} = \min( \frac{1}{40}, \frac{3}{32}
  \mathbf{O}_{9}^{-1})$. Let $\dot{V}' = RG_{\rm dv}[\vec{V}_{\rm bk},\dot{V}]$ Then one has the following bound:

\[
||\dot{V}'||
\le
\frac{15}{16}
||\dot{V}||\ .
\]

\end{Proposition}

\noindent{\bf Proof:}
This proposition is a direct
consequence of Lemma \ref{LAJ13}.

\qed

For the control of the infrared contributions to the log-moment generating function
we will finally need a very coarse Lipschitz estimate on the $\delta b$ functions.

\begin{Lemma}\label{dbforIRlem}
For all $\vec{V}^{1}$, $\vec{V}^{2}$ in $\bar{B}\left(\bar{V},\frac{1}{6}\right)$
we have
\[
|\delta b_{\Delta(0)}[\vec{V}^{1}]-\delta b_{\Delta(0)}[\vec{V}^{2}]|
\le 4 ||\vec{V}^{1}-\vec{V}^{2}||\ .
\]
\end{Lemma}

\noindent{\bf Proof:}
By our assumptions on exponents, $||\vec{V}-\bar{V}||<\frac{1}{2}$
implies one is in the domain of applicability of Theorem \ref{mainestthm}
as well as all the lemmas that led to its proof.
In particular Lemma \ref{L38lem} with $\lambda =1$
gives us the bound $|\delta b_{\Delta(0)}[\vec{V}]|\le 1$
provided  $\mathcal{O}_{30}L^{5}\bar{g}^{\frac{11}{12}-\frac{1}{3}\eta_{R}}\le 1$.
However we can take the latter for granted since we are in the small $\epsilon$ regime
and $\frac{11}{12}-\frac{1}{3}\eta_{R}>0$.
Now Lemma \ref{Lipschitzlem} with $\nu=\frac{1}{3}$ immediately produces the desired estimate.
\qed

Now recall from \S\ref{algconsec}
that
\[
\mathcal{S}^{\rm T,IR}_{r,s}(\tilde{f},\tilde{j})
=
\sum_{q_{+} \le q < s}
\left(
\delta b_{\Delta(0)}
\left[
\vec{V}^{(r,q)}(\tilde{f},\tilde{j})
\right]
-
\delta b_{\Delta(0)}
\left[
\vec{V}^{(r,q)}(0,0)
\right]
\right)
\]
where
\[
\vec{V}^{(r,q)}(\tilde{f},\tilde{j})=RG_{\rm ex}^{q-q_{+}}\left(
\vec{V}^{(r,q_{+})}(\tilde{f},\tilde{j})
\right)\ .
\]
With a view to lighten the notation we write
\[
\vec{V}^{(r,q)}(\tilde{f},\tilde{j})=\vec{V}_{\rm bk}^{(r,q)}+\dot{V}^{(r,q)}
\]
where
\[
\vec{V}_{\rm bk}^{(r,q)}=\vec{V}^{(r,q)}(0,0)=\iota(RG^{q-r}(v))\in\mathcal{E}_{\rm bk}
\]
and
\[
\dot{V}^{(r,q)}=\vec{V}^{(r,q)}(\tilde{f},\tilde{j})-\vec{V}^{(r,q)}(0,0)\in\mathcal{E}_{\rm pt}\ .
\]
We will control the latter via Proposition \ref{contractionbound}.

First note that
\[
||\vec{V}_{\rm bk}^{((r,q))}-\bar{V}||=||RG^{q-r}(v)||\le\frac{\rho}{3}\ .
\]
To make this at most $\frac{1}{40}$ we add the new requirement on $\rho$:
\[
\rho\le\frac{3}{40}\ .
\]
If we can ensure that $||\dot{V}^{(r,q_{+})}||\le \mathbf{O}_{10}$
then a trivial inductive use of Proposition
\ref{contractionbound} will imply that
\[
||\dot{V}^{(r,q)}||\le \mathbf{O}_{10}\times\left(\frac{15}{16}\right)^{q-q_{+}}
\]
for all $q$, such that $q_{+}\le q\le s$. We again include the value $s$ although
it does not belong to what we called the infrared regime in order to pass the baton
to the next section about controlling the boundary term.
In view of (\ref{seedforIReq}), we now impose the new domain condition
\begin{equation}
\left(\mathbf{O}_{6}L^{\frac{5}{2}}\right)^{q_{+}-q_{-}}\times
\max\left\{
L^{(3-[\phi])q_{-}}
||\tilde{f}||_{L^{\infty}},
11\mathcal{C}_{1}(\epsilon)\alpha_{\rm u}^{q_{-}}Y_2 \bar{g}^{-e_2}\times ||\tilde{j}||_{L^{\infty}}
\right\}\le \mathbf{O}_{10}\ .
\label{finaldomaineq}
\end{equation}
Now Proposition \ref{contractionbound} followed by Lemma \ref{dbforIRlem}
imply that for any $q$ with $q_{+}\le q<s$ we have
\[
\left|
\delta b_{\Delta(0)}\left[\vec{V}^{(r,q)}(\tilde{f},\tilde{j})\right]
-\delta b_{\Delta(0)}\left[\vec{V}^{(r,q)}(0,0)\right]
\right|\le
4\mathbf{O}_{10}\times\left(\frac{15}{16}\right)^{q-q_{+}}\ .
\]
Hence we get the uniform absolute convergence
of the sum over $q$ needed to say
\[
\lim_{\substack{r\rightarrow -\infty \\ s\rightarrow \infty}}
\mathcal{S}_{r,s}^{T,IR}(\tilde{f},\tilde{j})=\mathcal{S}^{T,IR}(\tilde{f},\tilde{j})
\]
with
\[
\mathcal{S}^{T,IR}(\tilde{f},\tilde{j})
=
\sum_{q =q_{+}}^{\infty}
\left(
\delta b_{\Delta(0)}
[
\vec{V}^{(-\infty,q)}(\tilde{f},\tilde{j})
]
-
\delta b_{\Delta(0)}
[
\iota(v_{\ast})
]
\right)
\]
where
\[
\vec{V}^{(-\infty,q)}(\tilde{f},\tilde{j})=RG_{\rm ex}^{q-q_{-}}\left(
\vec{V}^{(-\infty,q_{-})}(\tilde{f},\tilde{j})\right)
\]
and $\vec{V}^{(-\infty,q_{-})}(\tilde{f},\tilde{j})$ has been defined in (\ref{Vqminuseq}).
The limit $\mathcal{S}^{T,IR}(\tilde{f},\tilde{j})$
is analytic and the order of the $r\rightarrow -\infty$, $s\rightarrow\infty$ limits is immaterial.

\subsection{The boundary term}

Let $\vec{V}\in\mathcal{E}_{\rm ex}$
and simply denote by
\[
(\beta_4,\beta_3,\beta_2,\beta_1, W_5, W_6,f,R)
\in\mathbb{C}^{7}\times C_{\rm bd}^{9}(\mathbb{R},\mathbb{C})
\]
its component at $\Delta=\Delta(0)$.
We let
\[
\begin{split}
\partial\mathcal{Z}[\vec{V}] & =
\int {\rm d}\mu_{C_0}(\phi)\ e^{f\phi}\times
\left\{
\exp\left(
-\beta_4 :\phi^4:_{C_0}
-\beta_3 :\phi^3:_{C_0}
-\beta_2 :\phi^2:_{C_0}
-\beta_1 :\phi:_{C_0}
\right)
\right.\\
 & \left.\times(1+W_5:\phi^5:_{C_0}+W_6:\phi^6:_{C_0})+R(\phi)\right\}
\end{split}
\]
which reduces to an integral over a single real variable still denoted by $\phi$.
Let $\partial\mathcal{Z}_{\ast}=\partial\mathcal{Z}[\iota(v_{\ast})]$
which is the value at the infrared fixed point.
We have
\[
\partial\mathcal{Z}_{\ast}=
\int {\rm d}\mu_{C_0}(\phi)\ 
\left\{\exp\left(
-g_{\ast} :\phi^4:_{C_0}
-\mu_{\ast} :\phi^2:_{C_0}
\right)
+R_{\ast}(\phi)\right\}
\]
with $g_{\ast}=\bar{g}+\delta g_{\ast}$.
Recall that $g_{\ast}$, $\mu_{\ast}$, $R_{\ast}$ are real.
Note that by Jensen's inequality and the basic properties of Wick ordering
on has the lower bound
\[
\int {\rm d}\mu_{C_0}(\phi)\ 
\exp\left(
-g_{\ast} :\phi^4:_{C_0}
-\mu_{\ast} :\phi^2:_{C_0}
\right)\ge
\exp\left(
-\int {\rm d}\mu_{C_0}(\phi)\ \left(
g_{\ast} :\phi^4:_{C_0}
+\mu_{\ast} :\phi^2:_{C_0}
\right)\right)=1\ .
\]
Besides
\[
\begin{split}
\left|
\int {\rm d}\mu_{C_0}(\phi)\ R_{\ast}(\phi)
\right| & \le
\sup_{\phi\in\mathbb{R}} |R_{\ast}(\phi)|\le
\sup_{\phi\in\mathbb{R}} ||R_{\ast}(\phi)||_{\partial\phi,\phi,h} \\
 & \le \bar{g}^{-2}|||R_{\ast}|||_{\bar{g}}\le \bar{g}^{e_{R}-2}\frac{\rho}{13}\ .
\end{split}
\]
Since $e_R>e_4+1\ge 2$, $\bar{g}\le 1$ and $\rho<\frac{3}{40}$, we clearly have
$\partial\mathcal{Z}_{\ast}\ge \frac{1}{2}$.

Now if $||\vec{V}-\bar{V}||<\frac{1}{2}$
it is easy to see that $|\partial\mathcal{Z}[\vec{V}]|\le \mathcal{C}_{5}(\epsilon)$ with
\[
\begin{split}
\mathcal{C}_{5}(\epsilon) = &
\int {\rm d}\mu_{C_0}(\phi)\ e^{\frac{1}{2}L^{3-[\phi]}|\phi|}\times \\
 & \left\{
\exp\left[
-\frac{1}{2}\bar{g}\phi^4+\frac{3}{4}\bar{g}^{1-\eta}
\left(|\phi|^3+13\phi^2+7|\phi|+14\right)
\right] 
\right. \\
 & \times\left(
1+\frac{1}{2}\bar{g}^{2-2\eta}\left(|\phi|^5+20|\phi|^3+60|\phi|\right)
 +\frac{1}{2}\bar{g}^{2-2\eta}\left(\phi^6+30\phi^4+180\phi^2+120\right)
\right) \\
 & \left.+\frac{1}{2}\bar{g}^{e_R-2}\right\}\ .
\end{split}
\]
Indeed, by undoing the Wick ordering
\[
-\Re\left[
\beta_4 :\phi^4:_{C_0}
+\beta_3 :\phi^3:_{C_0}
+\beta_2 :\phi^2:_{C_0}
+\beta_1 :\phi:_{C_0}
\right]=-\bar{g}\phi^4-Y(\phi)
\]
with
\[
\begin{split}
Y(\phi)= & \Re(\beta_4-\bar{g})\phi^4 \\
 & +(\Re\beta_3)\phi^3 \\
 & +(\Re\beta_2-6C_0(0)\Re\beta_4)\phi^2\\ 
 & +(\Re\beta_1-3C_0(0)\Re\beta_3)\phi\\
 & +(-C_0(0)\Re\beta_2+3C_0(0)^2\Re\beta_4)\ .  
\end{split}
\]
Using $|\Re(\beta_4-\bar{g})|<\frac{1}{2}\bar{g}^{e_4}\le\frac{1}{2}\bar{g}$ for the fourth degree monomial
and $|\Re\beta_k|\le\frac{3}{2}\bar{g}^{1-\eta}$ for $k=1,2,3,4$ when bounding the lower degree monomials, 
and finally using $C_0(0)\le 2$
we obtain
\[
|Y(\phi)|\le  \frac{1}{2}\bar{g}\phi^4+\frac{3}{4}\bar{g}^{1-\eta}
\left(|\phi|^3+13\phi^2+7|\phi|+14\right)\ .
\]
The bounds on $W_k:\phi^k:_{C_0}$, for $k=5,6$ are similar and use the explicit Wick ordering
formulas given in the proof of Lemma \ref{L17lem}.

Since $\partial\mathcal{Z}[\vec{V}]$ is clearly analytic in the domain
$||\vec{V}-\bar{V}||<\frac{1}{2}$, Lemma \ref{Lipschitzlem} with $\nu=\frac{1}{3}$
tell us that for all $\vec{V}^{1}$, $\vec{V}^{2}$ in $\bar{B}\left(\bar{V},\frac{1}{6}\right)$
one has the Lipschitz estimate
\[
|\partial\mathcal{Z}[\vec{V}^{1}]-\partial\mathcal{Z}[\vec{V}^{2}]|\le
4\mathcal{C}_{5}(\epsilon)||\vec{V}^{1}-\vec{V}^{2}||\ .
\]
We now have, using the outcome of the discussion for the infrared regime
\begin{eqnarray*}
|\partial\mathcal{Z}_{r,s}(\tilde{f},\tilde{j})-\partial\mathcal{Z}_{\ast}|
 & = & |\partial\mathcal{Z}[\vec{V}^{(r,s)}(\tilde{f},\tilde{j})]-\partial\mathcal{Z}[\iota(v_{\ast})]| \\
 & \le & 4\mathcal{C}_{5}(\epsilon)\times
 \left[
 ||\vec{V}^{(r,s)}(\tilde{f},\tilde{j})-\vec{V}^{(r,s)}(0,0)||+||\vec{V}^{(r,s)}(0,0)-\iota(v_{\ast})||
 \right] \\
 & \le & 4\mathcal{C}_{5}(\epsilon)\times
 \left[ ||\dot{V}^{(r,s)}||+||RG^{s-r}(v)-v_{\ast}||
 \right] \\
 & \le & 4\mathcal{C}_{5}(\epsilon)\times
 \left[ \mathbf{O}_{10}\times\left(\frac{15}{16}\right)^{s-q_{+}}+ c_1(\epsilon)^{s-r}||v-v_{\ast}||
 \right]\ .
\end{eqnarray*}
One of course has a similar and simpler estimate for the quantity $\partial\mathcal{Z}_{r,s}(0,0)$
appearing in the denominator of the boundary ratio. Namely, the $\mathbf{O}_{10}$ term is absent.
Bounding $c_1(\epsilon)^{s-r}$ by $c_1(\epsilon)^{s-q_{+}}$ and
using the previous lower bound $\partial\mathcal{Z}_{\ast}\ge\frac{1}{2}$
we see that
\[
\frac{\partial\mathcal{Z}_{r,s}(\tilde{f},\tilde{j})}{\partial\mathcal{Z}_{r,s}(0,0)}
\longrightarrow 1
\]
when $s\rightarrow\infty$, uniformly in $r\le q_{-}$.
Therefore the boundary term $\mathcal{S}^{\rm T,BD}$
disappears when $r\rightarrow -\infty$, $s\rightarrow\infty$ regardless of the order of limits.

%% file: Constructmeas.tex
\section{Construction of the limit measures and invariance properties}

As a consequence of what we have shown in the previous section we see that
\[
\mathcal{S}_{r,s}(\tilde{f},\tilde{j})=\exp\left(
\mathcal{S}_{r,s}^{\rm T}(\tilde{f},\tilde{j})
\right)
\]
converges uniformly to the analytic function
\[
\mathcal{S}(\tilde{f},\tilde{j})=\exp\left(
\mathcal{S}^{\rm T}(\tilde{f},\tilde{j})
\right)
\]
in a suitable neighborhood of $\tilde{f}=\tilde{j}=0$
in $S_{q_{-},q_{+}}(\mathbb{Q}_{p}^{3},\mathbb{C})$,
when $r\rightarrow -\infty$ and $s\rightarrow \infty$.
Using the multivariate Cauchy formula it is immediate that the cut-off correlators
\[
\left\langle
\tilde{\phi}(\tilde{f}_1)\cdots\tilde{\phi}(\tilde{f}_n)\ 
N_{r}[\tilde{\phi}^2](\tilde{j}_1)\cdots
N_{r}[\tilde{\phi}^2](\tilde{j}_m)
\right\rangle_{r,s}=
\]
\[
\frac{1}{(2i\pi)^{n+m}}\oint\cdots\oint
\prod_{j=1}^{n}\frac{{\rm d}z_j}{z_j^2}\ 
\prod_{k=1}^{m}\frac{{\rm d}u_k}{u_k^2}\ 
\mathcal{S}_{r,s}(z_1\tilde{f}_1+\cdots+z_n\tilde{f}_n,u_1\tilde{j}_1+\cdots+u_m\tilde{j}_m)
\]
converge to the similar integrals with $\mathcal{S}$ instead of $\mathcal{S}_{r,s}$.
The contours of integration are governed by the domain condition (\ref{finaldomaineq}).
We define our mixed correlators by
\[
\left\langle
\tilde{\phi}(\tilde{f}_1)\cdots\tilde{\phi}(\tilde{f}_n)\ 
N[\tilde{\phi}^2](\tilde{j}_1)\cdots
N[\tilde{\phi}^2](\tilde{j}_m)
\right\rangle=
\]
\[
\frac{1}{(2i\pi)^{n+m}}\oint\cdots\oint
\prod_{j=1}^{n}\frac{{\rm d}z_j}{z_j^2}\ 
\prod_{k=1}^{m}\frac{{\rm d}u_k}{u_k^2}\ 
\mathcal{S}(z_1\tilde{f}_1+\cdots+z_n\tilde{f}_n,u_1\tilde{j}_1+\cdots+u_m\tilde{j}_m)
\]
which are multilinear in the $\tilde{f}$'s and $\tilde{j}$'s.
Because of the uniform bounds on $\mathcal{S}_{r,s}^{\rm T}$, and therefore on $\mathcal{S}^{\rm T}$,
proved in the last section and thanks to Cauchy's formula, it is immediate that the
pure $\tilde{\phi}$ or $N[\tilde{\phi^2}]$ correlators will satisfy Condition 4) in Theorem \ref{Hamburgerthm}.
The other conditions are satisfied by the cut-off correlators $\langle\cdots\rangle_{r,s}$
as joint moments of random variables obtained from the probability measures $\nu_{r,s}$.
As these properties are preserved in the limit $r\rightarrow -\infty$ and $s\rightarrow \infty$ we can use
Theorem \ref{Hamburgerthm} to affirm the existence and uniqueness of the measures ${\nu}_{\phi}$ and $\nu_{\phi^2}$
mentioned in Theorem \ref{themainthm}.
By the uniqueness part of Theorem \ref{Hamburgerthm}, the invariance properties of the measures ${\nu}_{\phi}$ and $\nu_{\phi^2}$
follow from those of the moments. Hence
it is enough to show Parts 1) and 2) of Theorem \ref{themainthm}.
These are easier to prove from the functional integral definitions of the cut-off correlators.

Indeed, one can trivially check that for $M\in GL_3(\mathbb{Z}_{p})$
one has
\[
\left\langle
\tilde{\phi}(M\cdot\tilde{f}_1)\cdots\tilde{\phi}(M\cdot\tilde{f}_n)\ 
N_{r}[\tilde{\phi}^2](M\cdot\tilde{j}_1)\cdots
N_{r}[\tilde{\phi}^2](M\cdot\tilde{j}_m)
\right\rangle_{r,s}=
\]
\[
\left\langle
\tilde{\phi}(\tilde{f}_1)\cdots\tilde{\phi}(\tilde{f}_n)\ 
N_{r}[\tilde{\phi}^2](\tilde{j}_1)\cdots
N_{r}[\tilde{\phi}^2](\tilde{j}_m)
\right\rangle_{r,s}
\]
because ${\rm d}\mu_{C_r}$ is invariant by rotation and $M\cdot\Lambda_s=\Lambda_s$.

Also if $y\in\mathbb{Q}_{p}^{3}$ with $|y|\le L^s$ then
\[
\left\langle
\tilde{\phi}(\tau_y\tilde{f}_1)\cdots\tilde{\phi}(\tau_y\tilde{f}_n)\ 
N_{r}[\tilde{\phi}^2](\tau_y\tilde{j}_1)\cdots
N_{r}[\tilde{\phi}^2](\tau_y\tilde{j}_m)
\right\rangle_{r,s}=
\]
\[
\left\langle
\tilde{\phi}(\tilde{f}_1)\cdots\tilde{\phi}(\tilde{f}_n)\ 
N_{r}[\tilde{\phi}^2](\tilde{j}_1)\cdots
N_{r}[\tilde{\phi}^2](\tilde{j}_m)
\right\rangle_{r,s}
\]
because $\Lambda_s$ is unchanged by this translation as results from ultrametricity.

Finally, by changing variables from $\tilde{\phi}$ to $\tilde{\phi}_{\leadsto 1}$, one has
\[
\left\langle
\tilde{\phi}(L\cdot\tilde{f}_1)\cdots\tilde{\phi}(L\cdot\tilde{f}_n)\ 
N_{r}[\tilde{\phi}^2](L\cdot\tilde{j}_1)\cdots
N_{r}[\tilde{\phi}^2](L\cdot\tilde{j}_m)
\right\rangle_{r,s}=
\]
\[
\left\langle
\tilde{\phi}(\tilde{f}_1)\cdots\tilde{\phi}(\tilde{f}_n)\ 
N_{r}[\tilde{\phi}^2](\tilde{j}_1)\cdots
N_{r}[\tilde{\phi}^2](\tilde{j}_m)
\right\rangle_{r+1,s+1}\times
\left[L^{-(3-[\phi])}
\right]^n\times
\left[L^{-(3-2[\phi])}Z_{2}^{-1}\right]^m\ .
\]
Noting that $|L|=L^{-1}$ and $Z_2=L^{-\frac{1}{2}\eta_{\phi^2}}$ by definition of $\eta_{\phi^2}$, and
from the existence of the $r\rightarrow -\infty$, $s\rightarrow\infty$ limits, we see
that the property in Part 3) of Theorem \ref{themainthm} holds for $\lambda=L$.
Thus it holds for the subgroup $L^{\mathbb{Z}}$ it generates.

A trivial consequence of these invariance properties is that
\[
\langle N[\tilde{\phi}^2](\tilde{j})\rangle =0
\] 
identically. Namely, there is no one-point function.
Indeed, it is enough to show this for $\tilde{j}=\bbone_{\mathbb{Z}_{p}^{3}}$.
In that case, by translation invariance followed by scale invariance
\begin{eqnarray*}
\langle N[\tilde{\phi}^2](\bbone_{\mathbb{Z}_{p}^{3}})\rangle & = & 
L^3 \langle N[\tilde{\phi}^2](\bbone_{(L\mathbb{Z}_{p})^{3}})\rangle \\
 & = & L^3\times L^{-3+2[\phi]+\frac{1}{2}\eta_{\phi^2}}\times
 \langle N[\tilde{\phi}^2](\bbone_{\mathbb{Z}_{p}^{3}})\rangle \\
 & = & L^3\alpha_{\rm u}^{-1}\times \langle N[\tilde{\phi}^2](\bbone_{\mathbb{Z}_{p}^{3}})\rangle\ .
\end{eqnarray*}
By Lemma \ref{alphabdlem} it is clear that $L^3\alpha_{\rm u}^{-1}>1$ for $\epsilon$ small and the vanishing follows.
We leave it as an exercise to show this same fact explicitly using
the $\mathcal{S}^{\rm T,UV}+\mathcal{S}^{\rm T,IR}$ expression for $\langle N[\tilde{\phi}^2](\bbone_{\mathbb{Z}_{p}^{3}})\rangle$.
This hinges on showing that the vector in $\mathcal{E}_{\rm pt}$ corresponding
to an $e_{\rm u}$ perturbation in the box $\Delta(0)$ only is an eigenvector of $D_{\iota(v_{\ast})}RG_{\rm ex}$
with eigenvalue $L^{-3}\alpha_{\rm u}$.
One has a similar statement for the evaluation of $D_{\iota(v_{\ast})}\delta b_{\Delta(0)}$ on that vector.

%% file: Nontriv.tex
\section{Nontriviality and proof of existence of anomalous dimension}

\subsection{The two-point and four-point functions of the elementary field}

We have constructed the generalized random field
$\tilde{\phi}$ via constructing and proving the analyticity of
$\mathcal{S}^{\rm T}(\tilde{f}, 0)$, the
cumulant generating function. We now show that the process $\tilde{\phi}$ is not
Gaussian. In particular we show that in the small $\epsilon$ regime one has

\[
\frac{d^4} {dz^4} \Big|_{z=0}
\mathcal{S}^{\rm T}(z\bbone_{\mathbb{Z}_p^3}, 0)
=
\langle
\tilde{\phi} ( \bbone_{\mathbb{Z}_{p}^3} )^4 
\rangle
- 3 \langle
 \tilde{\phi} ( \bbone_{\mathbb{Z}_{p}^3})^2
\rangle < 0\ .
\]

 We establish the inequality above by expanding
$\mathcal{S}^{\rm T}(z\bbone_{\mathbb{Z}_{p}^3},0)$ and isolating a
part that explicitly contains first order pertubation theory. We will
calculate the derivative by hand for this explicit part and use Cauchy
bounds to estimate the contribution of the remainder. From now on we
will drop the tildes from the notation for the fields $\tilde{\phi}$ and $N[\tilde{\phi}^2]$
but we will still use tildes for test functions if needed.\\

Since $z \bbone_{\mathbb{Z}_{p}^{3}} \in
S_{0,0}(\mathbb{Q}_{p}^{3}, \mathbb{C})$ we can set $q_{-} = q_{+} =
0$. From section \S\ref{controldevsec} and in particular the domain condition (\ref{finaldomaineq})
we know that $\mathcal{S}^{T}(z\bbone_{\mathbb{Z}_{p}^{3}},0)$ is an analytic function for $z$ such
that $|z| < \mathbf{O}_{10}$. This condition is assumed throughout this section.
We will repeatedly make use of the fact that for
$z$ in this domain $|z|
\le 1$ which follows from $\mathbf{O}_{10}\le \frac{1}{40}$. In particular for $z$ in that domain we have

\begin{equation*}
\begin{split}
\mathcal{S}^{\rm T}(z \bbone_{\mathbb{Z}_{p}^{3}},0)
=&
\mathcal{S}^{\rm T,FR}(z \bbone_{\mathbb{Z}_{p}^{3}},0)
+
\mathcal{S}^{\rm T,UV}(z \bbone_{\mathbb{Z}_{p}^{3}},0)
+
\mathcal{S}^{\rm T,MD}(z \bbone_{\mathbb{Z}_{p}^{3}},0)
+
\mathcal{S}^{\rm T,IR}(z \bbone_{\mathbb{Z}_{p}^{3}},0)\ .
\end{split}
\end{equation*}

For our choice of test function we have:

\begin{equation*}
\begin{split}
\mathcal{S}^{\rm T,FR}(z \bbone_{\mathbb{Z}_{p}^{3}},0) =&\ \frac{1}{2}
z^{2} \left( \bbone_{\mathbb{Z}_{p}^{3}}, C_{-\infty}
  \bbone_{\mathbb{Z}_{p}^{3}} \right)\\[1.5ex]
\mathcal{S}^{\rm T,UV}(z \bbone_{\mathbb{Z}_{p}^{3}},0) =&\ 0
\textnormal{   since  } \tilde{j}=0\\[1.5ex]
\mathcal{S}^{\rm T,MD}(z \bbone_{\mathbb{Z}_{p}^{3}},0) =&\ 0
\textnormal{   since  } q_{-} = q_{+} = 0\\[1.5ex]
\mathcal{S}^{\rm T,IR}(z \bbone_{\mathbb{Z}_{p}^{3}},0) =&
\ \sum_{q=0}^{\infty} 
\left( 
\delta b_{\Delta(0)} 
\left[ 
\vec{V}^{(-\infty,q)} 
(z \bbone_{\mathbb{Z}_{p}^{3}}, 0) 
\right]
-
\delta b_{\Delta(0)} [ \vec{V}_{*} ]
\right)
\end{split}
\end{equation*}

where $\vec{V}_{*} = \iota( v_{*}) = \vec{V}^{(-\infty,q)}(0,0)$. \\

By previous considerations we
know that up to scale $q_{-}=0$ the test function $\tilde{f} = z
\bbone_{\mathbb{Z}_{p}}$ does not influence the evolution of the other
parameters, thus for scales $q \le q_{-} = 0$ all components of
$\vec{V}^{(-\infty,q)}(z \bbone_{\mathbb{Z}_{p}^{3}}, 0)$
other than the $f$ component take their fixed point value. Additionally we
know that for scales $q \ge q_{+}=0$  the vector
$\vec{V}^{(-\infty,q)}$ deviates from $\vec{V}_{*}$ only at $\Delta =
\Delta(0)$. 

We write

\[ 
\vec{V}^{(-\infty,q)}(z \bbone_{\mathbb{Z}_{p}^{3}}, 0)
=
\left(
(\beta^{(q)}_{4,\Delta},\dots,
\beta^{(q)}_{1,\Delta},W^{(q)}_{5,\Delta}, W^{(q)}_{6,\Delta},
f^{(q)}_{\Delta}, R^{(q)}_{\Delta})
\right)_{\Delta \in \mathbb{L}}\ .
\]

Keeping our previous observations in mind for $k=1,2,3,4$ we decompose $\beta^{(q)}_{k,\Delta}$ as follows:

\begin{equation*}
\begin{split}
\beta^{(q)}_{4,\Delta}
=&
\begin{cases}
g_{*} + \beta^{(q,{\rm exp})}_{4} +  \beta^{(q,{\rm imp})}_{4} &
\textnormal{  if  } \Delta = \Delta(0)\\
g_{*} &
\textnormal{  if  } \Delta \not = \Delta(0)
\end{cases}\\[1.5ex]
\beta^{(q)}_{3,\Delta}
=&
\begin{cases}
\beta^{(q,{\rm exp})}_{3} + \beta^{(q,{\rm imp})}_{3} &
\textnormal{  if  } \Delta = \Delta(0)\\
0 &
\textnormal{  if  } \Delta \not = \Delta(0)
\end{cases}\\[1.5ex]
\beta^{(q)}_{2,\Delta}
=&
\begin{cases}
\mu_{*} + \beta^{(q,{\rm exp})}_{2} + \beta^{(q,{\rm imp})}_{2} &
\textnormal{  if  } \Delta = \Delta(0)\\
\mu_{*} &
\textnormal{  if  } \Delta \not = \Delta(0)
\end{cases}\\[1.5ex]
\beta^{(q)}_{1,\Delta}
=&
\begin{cases}
\beta^{(q,{\rm exp})}_{1} + \beta^{(q,{\rm imp})}_{1} &
\textnormal{  if  } \Delta = \Delta(0)\\
0 &
\textnormal{  if  } \Delta \not = \Delta(0)
\end{cases}
\end{split}
\end{equation*}

Here ``exp'' and ``imp'' are abbreviations for explicit and
implicit. The quantities $\beta^{(q,{\rm exp})}_{k}$ and
$\beta^{(q,{\rm imp})}_{k}$ will be defined inductively starting from
$q=0$. We start with the following intital condition:

\[
\textnormal{for  }
k=1,2,3,4
\textnormal{  we set  }
\beta^{(0,{\rm exp})}_{k} = \beta^{(0,{\rm imp})}_{k} = 0.
\]

Now we prepare to give the inductive part of the definition. Recall that for
$k=1,2,3,4$ the evolution
of our couplings is given by 

\begin{equation*}
\begin{split}
\beta^{(q+1)}_{k,\Delta(0)} 
= &
L^{-k[\phi]} 
\left( 
\sum_{\Delta \in [L^{-1}\Delta(0)]} \beta^{(q)}_{k,\Delta} 
\right) 
- 
\delta \beta_{k,1,\Delta(0)} 
\left[ \vec{V}^{(-\infty,q)}(z\bbone_{\mathbb{Z}_{p}^{3}},0)
\right]\\
&- 
\delta \beta_{k,2,\Delta(0)} 
\left[
\vec{V}^{(-\infty,q)}(z\bbone_{\mathbb{Z}_{p}^{3}},0)
\right]
+ 
\xi_{k,\Delta(0)} 
\left[
\vec{V}^{(-\infty,q)}(z\bbone_{\mathbb{Z}_{p}^{3}},0)
\right]\ .
\end{split}
\end{equation*}

We introduce some more short hand. For $k=1,2,3,4$ we define $\beta_{k}^{*}$ to be the
corresponding component of $\vec{V}_{*} \in \mathcal{E}_{\rm bk}$. In
particular $\beta_{4}^{*} = g_{*}$, $\beta_{3}^{*} = 0$,
$\beta_{2}^{*} = \mu_{*}$, and $\beta_{1}^{*}=0$. These are also seen as constant vectors in $\mathbb{C}^{\mathbb{L}}$.\\

We now use the fact that $\vec{V}_{*}$ is a fixed point of $RG_{\rm ex}$
to arrive at the following formula:

\begin{equation}\label{4pointflow}
\begin{split}
\beta^{(q+1)}_{k,\Delta(0)}  =& \beta_{k}^{*} + L^{-k[\phi]} \left(
  \beta^{(q,{\rm exp})}_{k} + \beta^{(q,{\rm imp})}_{k} \right)\\
&- \delta \beta_{k,1,\Delta(0)} 
\left[
  \vec{V}^{(-\infty,q)}(z\bbone_{\mathbb{Z}_{p}^{3}},0)
\right]\\
&+ 
\left( 
\delta \beta_{k,2,\Delta(0)}
\left[
\vec{V}_{*}
\right]
-
\delta \beta_{k,2,\Delta(0)}
\left[
\vec{V}^{(-\infty,q)}(z\bbone_{\mathbb{Z}_{p}^{3}},0)
\right] 
\right)\\
&-  
\left( 
\xi_{k,\Delta(0)} 
\left[
\vec{V}_{*}
\right] 
- 
\xi_{k,\Delta(0)} 
\left[
\vec{V}^{(-\infty,q)}(z\bbone_{\mathbb{Z}_{p}^{3}},0)
\right]
\right).
\end{split}
\end{equation}

Above we have used the fact that $\delta
b_{k,1,\Delta}\left[\vec{V}_{*}\right]=0$. We now decompose $\delta
\beta_{k,1,\Delta(0)} \left[ \vec{V}^{(-\infty,q)}(z\bbone_{\mathbb{Z}_{p}^{3}},0) \right]$.  For $0 \le k< l \le 4$ and $\beta, f \in
\mathbb{C}^{\mathbb{L}}$ define

\[
F_{k,l}
\left[ \beta ,f \right] 
= 
\binom{l}{k} 
\int_{\left( L^{-1} \Delta(0) \right)^{l-k}} {\rm d}^3 a\ {\rm d}^3 b_1\cdots {\rm d}^3 b_{l-k}\ 
\beta(a)\times \prod_{i=1}^{l-k}
\left[\Gamma(a-b_i) f(b_{i})\right].
\]

With this notation we have:

\[
\delta \beta_{k,1,\Delta(0)} 
\left[ \vec{V}^{(-\infty,q)}(z\bbone_{\mathbb{Z}_{p}^{3}},0) \right]
= 
- 
\sum_{l=k+1}^{4}
L^{-k[\phi]} 
F_{k,l}\left[ \beta^{(q)}_{l},f^{(q)} \right].
\]

With this notation we define the evolution for
$\beta^{(q,{\rm exp})}_{k}$ and $\beta^{(q,{\rm imp})}_{k}$ as follows:

\begin{equation}\label{4pointexp}
\beta_{k}^{(q+1),{\rm exp}} = L^{-k[\phi]} \beta^{(q,{\rm exp})}_{k} + \sum_{l=k+1}^{4} L^{-k[\phi]} F_{k,l}\left[ \beta_{l}^{*}
  + \beta^{(q,{\rm exp})}_{l} \bbone_{\Delta(0)}, f^{(q)}
\right]
\end{equation}

\begin{equation}\label{4pointerr}
\begin{split}
\beta_{k}^{(q+1),{\rm imp}} &= L^{-k[\phi]} \beta^{(q,{\rm imp})}_{k} + \sum_{l=k+1}^{4} L^{-k[\phi]} F_{k,l}\left[ \beta^{(q,{\rm imp})}_{l}\bbone_{\Delta(0)},f^{(q)}
\right] \\
&+ \left( \delta \beta_{k,2,\Delta(0)}[\vec{V}_{*}]-\delta \beta_{k,2,\Delta(0)}[
  \vec{V}^{(-\infty,q)}(z\bbone_{\mathbb{Z}_{p}^{3}},0)] \right) + \left( \xi_{k,\Delta(0)} [\vec{V}_{*}]- \xi_{k,\Delta(0)} [
  \vec{V}^{(-\infty,q)}(z\bbone_{\mathbb{Z}_{p}^{3}},0)] \right)\ .
\end{split}
\end{equation}\\

Here we have designated $\bbone_{\Delta(0)} : \mathbb{L} \rightarrow
\mathbb{C} $ as the indicator function
of $\{\Delta(0)\}$.\\

We also impose a splitting of the difference of vacuum
renormalizations at $\Delta(0)$. For $q \ge 0$ we have:

\[
\delta b_{\Delta(0)} \left[ \vec{V}^{(-\infty,q)}(z\bbone_{\mathbb{Z}_{p}^{3}},0) \right] 
-
\delta b_{\Delta(0)} [ \vec{V}_{*} ]
= 
\delta b^{(q,{\rm exp})} + \delta b^{(q,{\rm imp})}\ .
\]

We define

\begin{equation}\label{4pointdeltabexp}
\delta b^{(q,{\rm exp})} = -\sum_{l=1}^{4} F_{0,l}\left[
  \beta_{l}^{*} + \beta^{(q,{\rm exp})}_{l}\bbone_{\Delta(0)},f^{(q)} \right]\ ,
\end{equation}

\begin{equation}\label{4pointdeltaberr}
\begin{split}
\delta b^{(q,{\rm imp})} =& -\sum_{l=1}^{4} F_{0,l}\left[
  \beta^{(q,{\rm imp})}_{l}\bbone_{\Delta(0)},f^{(q)} \right]\\
 &+ \left( \delta
   \beta_{0,2,\Delta(0)}[\vec{V}^{(-\infty,q)}(z\bbone_{\mathbb{Z}_{p}^{3}},0)] - \delta
   \beta_{0,2,\Delta(0)}[\vec{V}_{*}] \right)\\
&+ \left( \xi_{0,\Delta(0)}[\vec{V}^{(-\infty,q)}(z\bbone_{\mathbb{Z}_{p}^{3}},0)] -
  \xi_{0,\Delta(0)}[\vec{V}_{*}] \right)\ .
\end{split}
\end{equation}

We now derive explicit formulas for $\beta^{(q,{\rm exp})}_{k}$ and $\delta b^{(q,{\rm exp})}$.\\

\begin{Lemma}\label{betaexp}
Given the previous inductive definitions for $\beta^{(q,{\rm exp})}_{k}$ for
$q \ge 0$ and $k=1,2,3,4$ we have the following explicit formulas:

\begin{equation*}
\begin{split}
\beta^{(q,{\rm exp})}_{4} = & 0 \\[2ex]
\beta^{(q,{\rm exp})}_{3} = & 0 \\[2ex]
\beta^{(q,{\rm exp})}_{2} = & 6qL^{-2q[\phi]}z^{2}g_{*} ||\Gamma||_{L^{2}}^{2} \\[2ex]
\beta^{(q,{\rm exp})}_{1} = & 
z^3 
g_{*}
L^{-q[\phi]} 
\left[ 
4
\frac{1-L^{-2q[\phi]}}{1-L^{-2[\phi]}}
\left( 
\int_{\mathbb{Q}_{p}^{3}}{\rm d}^3x\ \Gamma(x)^3 \right) 
+ 
12 
\left(
\sum_{n=0}^{q-1} n L^{-2n[\phi]} \right) ||\Gamma||_{L^{2}}^{2}
\times
\Gamma(0) 
\right]\ .
\end{split}
\end{equation*}\\

For $q \ge 0$ we also have

\begin{equation*}
\begin{split}
\delta b^{(q,{\rm exp})}
=&
-z^{4} g_{*} \left[
 L^{-4q[\phi]} \left( \int_{\mathbb{Q}_{p}^{3}} {\rm d}^{3}x\ \Gamma(x)^{4} \right)
+ 
6 L^{-4q[\phi]}q ||\Gamma||_{L^{2}}^{2} \Gamma(0)^{2} 
+12 L^{-2q[\phi]}\left(
\sum_{n=0}^{q-1} n L^{-2n[\phi]} \right) ||\Gamma||_{L^{2}}^{2}
\Gamma(0)^2 
\right.\\
 & + \left.
4L^{-2q[\phi]}\frac{1-L^{-2q[\phi]}}{1-L^{-2[\phi]}}
\Gamma(0) \left( \int_{\mathbb{Q}_{p}^{3}} {\rm d}^{3}x\ \Gamma(x)^{3}
\right) \right] 
-z^{2} \mu_{*} L^{-2[\phi]q} ||\Gamma||_{L^{2}}^{2}\ .
\end{split}
\end{equation*}

\end{Lemma}

\noindent{\bf Proof:} We first note that below one often sees expressions of the form
$\displaystyle \int_{L^{-1}\Delta(0)} \Gamma(x)^{n}$. In the statement
  of the theorem we extended the integration to all of
  $\mathbb{Q}_{p}^{3}$, we can do this since $\Gamma$ is supported on
  $L^{-1}\Delta(0)$.\\

 For $\beta^{(q,{\rm exp})}_{4}$ the result is immediate after recalling
 that $\beta^{(0,{\rm exp})}_{4} = 0$ and noticing
 the evolution for this parameter reduces to multiplication by $L^{-4[\phi]}$.\\

For $\beta^{(q,{\rm exp})}_{3}$ we have 

\begin{equation*}
\begin{split}
\beta^{(q,{\rm exp})}_{3}
=&
\sum_{n=0}^{q-1} L^{-3[\phi](q-n)}
F_{3,4}[ \beta_{4}^{*} + \beta^{(n,{\rm exp})}_{4} \bbone_{\Delta(0)}, f^{(n)} ]\\
=& 
\sum_{n=0}^{q-1} L^{-3[\phi](q-n)} F_{3,4}[ g_{*}, f^{(n)} ]\\
=& \sum_{n=0}^{q-1} 0\ .
\end{split}
\end{equation*}

The last line follows from ultrametricity and the fact that $\Gamma$
integrates to $0$. In particular $F_{j,j+1} \left[ \beta_{j}^{*},
f^{(\cdot)} \right]$ will always vanish.\\

For $\beta^{(q,{\rm exp})}_{2}$ we have 

\begin{equation*}
\begin{split}
\beta^{(q,{\rm exp})}_{2} 
=&
\sum_{n=0}^{q-1} L^{-2(q-n)[\phi]} \left(  
F_{2,4}
\left[\beta_{4}^{*} + \beta^{(n,{\rm exp})}_{4} \bbone_{\Delta(0)},
  f^{(n)} \right] 
+
F_{2,3}
\left[\beta_{3}^{*} + \beta^{(n,{\rm exp})}_{3} \bbone_{\Delta(0)},
  f^{(n)} \right]
\right)\\
=&
\sum_{n=0}^{q-1} L^{-2(q-n)[\phi]} 
F_{2,4}
\left[ g_{*}, f^{(n)} \right]\\
=&
\sum_{n=0}^{q-1} L^{-2(q-n)[\phi]} 6
\left(
\int_{(L^{-1}\Delta(0))^{3}}{\rm d}^{3}a\ {\rm d}^{3}b_{1}\ {\rm d}^{3}b_{2}\ 
g_{*} \prod_{i=1,2} \left[ \Gamma(a-b_{i}) L^{-n[\phi]}
  z \bbone_{\mathbb{Z}_{p}}(b_{i}) \right] \right)\\
=& \sum_{n=0}^{q-1} L^{-2(q-n)[\phi]} 6 z^{2} g_{*} L^{-2n[\phi]} \left(
    \int_{L^{-1}\Delta(0)}{\rm d}^{3}a\ \Gamma(a)^{2} \right)\\
=& \sum_{n=0}^{q-1} L^{-2q[\phi]} 6 z^{2} g_{*} ||\Gamma||_{L^{2}}^{2}
\end{split}
\end{equation*}
from which the formula for $\beta_{2}^{(q,{\rm exp})}$ follows.
Note that above we used the fact that $f^{(n)} = L^{-n[\phi]}
\bbone_{\Delta(0)}$ as a vector in $\mathbb{C}^{\mathbb{L}}$
or $L^{-n[\phi]}
\bbone_{\mathbb{Z}_{p}^{3}}$ as function on $\mathbb{Q}_{p}^{3}$.\\

For $\beta^{(q,{\rm exp})}_{1}$ we have 

\begin{equation*}
\begin{split}
\beta^{(q,{\rm exp})}_{1} =& 
\sum_{n=0}^{q-1} L^{-(q-n)[\phi]} \left(  
F_{1,4}
\left[\beta_{4}^{*} + \beta^{(n,{\rm exp})}_{4} \bbone_{\Delta(0)},
  f^{(n)} \right] 
+
F_{1,3}
\left[\beta_{3}^{*} + \beta^{(n,{\rm exp})}_{3} \bbone_{\Delta(0)},
  f^{(n)} \right]\right.\\
& \qquad\qquad
+\left.
F_{1,2}
\left[\beta_{2}^{*} + \beta^{(n,{\rm exp})}_{2} \bbone_{\Delta(0)},
  f^{(n)} \right]
\right)\\
=& 
\sum_{n=0}^{q-1} L^{-(q-n)[\phi]} \left(  
F_{1,4}
\left[g_{*}, f^{(n)} \right] 
+
F_{1,2} \left[ \mu_{*} + \beta^{(n,{\rm exp})}_{2} \bbone_{\Delta(0)},
  f^{(n)} \right]
\right).
\end{split}
\end{equation*}

Looking at the terms involved one sees

\[
F_{1,4}
\left[g_{*}, f^{(n)} \right] 
=
4
g_{*}
z^3
L^{-3n[\phi]}
\left( \int_{L^{-1}\Delta(0)} {\rm d}^3x\ \Gamma(x)^{3} \right)
\]

and

\begin{equation*}
\begin{split}
F_{1,2} \left[ \mu_{*} + \beta^{(n,{\rm exp})}_{2} \bbone_{\Delta(0)},
  f^{(n)} \right] =& F_{1,2} \left[ \mu_{*}, f^{(n)} \right] + F_{1,2} \left[ \beta^{(n,{\rm exp})}_{2} \bbone_{\Delta(0)},
  f^{(n)} \right]\\
=& F_{1,2} \left[ \beta^{(n,{\rm exp})}_{2} \bbone_{\Delta(0)},
  f^{(n)} \right]\\
=& 2 L^{-n[\phi]}z \Gamma(0) \times \left(  6nL^{-2n[\phi]}z^{2}g_{*}
  ||\Gamma||_{L^{2}}^{2} \right).
\end{split}
\end{equation*}

The formula for $\beta_{1}^{(q,{\rm exp})}$ then follows.\\

We now move on to $\delta b^{(q,{\rm exp})}$. To keep things lighter we have
left out terms with a vanishing contribution:

\[
\delta b^{(q,{\rm exp})} 
= 
- F_{0,4} \left[ g_{*}, f^{(q)} \right]
- F_{0,2} \left[ \beta^{(q,{\rm exp})}_{2} \bbone_{\Delta(0)} , f^{(q)} \right]
- F_{0,2} \left[\mu_{*}, f^{(q)} \right] 
- F_{0,1} \left[ \beta^{(q,{\rm exp})}_{1} \bbone_{\Delta(0)}  , f^{(q)} \right]\ .
\]

We calculate each of the terms appearing above:

\begin{equation*}
\begin{split}
F_{0,4} \left[ g_{*}, f^{(q)} \right] =& z^4 g_{*} L^{-4q[\phi]} \left(
  \int_{L^{-1}\Delta(0)}{\rm d}^3x\ \Gamma(x)^4 \right)\\[1.5ex]
F_{0,2} \left[ \beta^{(q,{\rm exp})}_{2} \bbone_{\Delta(0)} , f^{(q)}
\right]=& z^{2} L^{-2q[\phi]} \Gamma(0)^{2} \times \left[
  6qL^{-2q[\phi]}z^{2}g_{*} ||\Gamma||_{L^{2}}^{2} \right]\\[1.5ex]
F_{0,2} \left[\mu_{*}, f^{(q)} \right] 
=& z^{2} L^{-2q[\phi]} ||\Gamma||_{L^{2}}^{2} \mu_{*} \\[1.5ex]
F_{0,1} \left[ \beta^{(q,{\rm exp})}_{1} \bbone_{\Delta(0)}  , f^{(q)} \right]
=& z L^{-q[\phi]} \Gamma(0) \times
\left\{
z^3 g_{*}L^{-q[\phi]} 
\left[ 
4\frac{1-L^{-2q[\phi]}}{1-L^{-2[\phi]}}
\left( \int_{\mathbb{Q}_{p}^{3}}{\rm d}^3x\ \Gamma(x)^{3} \right)
\right.\right. \\
 & \left.\left. + 12 
\left(\sum_{n=0}^{q-1} n L^{-2[\phi]n} \right) ||\Gamma||_{L^{2}}^{2}\times\Gamma(0) 
\right]
\right\}\ .
\end{split}
\end{equation*}

This proves the formula for $\delta b^{(q,{\rm exp})}$. 

\qed\\

We now calculate running bounds for the $\beta^{(q,{\rm imp})}_{k}$.

\begin{Lemma}\label{betaimp}
In the small $\epsilon$ regime one has the following bounds for $q \ge
0$

\begin{equation*}
\begin{split}
|\beta_{4}^{(q,{\rm imp})}| \le& 
\mathbf{O}_{11} 
\times
q
\times
L^{8} 
\bar{g}^{2-2\eta} 
\left( \frac{15}{16} \right)^{q}\\[1.5ex]
|\beta_{3}^{(q,{\rm imp})}| \le& 
17
\times
\mathbf{O}_{11} 
\times
q
\times
L^{8}
\bar{g}^{2-2\eta}
\left( \frac{15}{16} \right)^{q}\\[1.5ex]
|\beta_{2}^{(q,{\rm imp})}| \le& 
253
\times
\mathbf{O}_{11} 
\times
q
\times
L^{8}
\bar{g}^{2-2\eta}
\left( \frac{15}{16} \right)^{q}\\[1.5ex]
|\beta_{1}^{(q,{\rm imp})}| \le& 
2497
\times
\mathbf{O}_{11} 
\times
q
\times
L^{8}
\bar{g}^{2-2\eta}
\left( \frac{15}{16} \right)^{q}\\[1.5ex]
|\delta b^{(q,{\rm imp})}|
\le&
\mathbf{O}_{12} 
\times
L^{8}
\times
\bar{g}^{2-2\eta}
\left( \frac{15}{16} \right)^{q}
\end{split}
\end{equation*}

where $\mathbf{O}_{11} = \left( 4 \mathcal{O}_{26} + 1\right)$ and
$\mathbf{O}_{12} = 319617 \times \mathbf{O}_{11}$.
\end{Lemma}

\noindent{\bf Proof:}
We note that for all $q \ge 0$ one has $\vec{V}^{(-\infty,q)}(z\bbone_{\mathbb{Z}_{p}^{3}},0),
\vec{V}_{*} \in \bar{B}(0,\frac{1}{6})$. Thus by the proof of
Lemma \ref{LAJ9} we have the following bounds for all $q \ge 0$ and
for $k=0,1,2,3,4$. 

\begin{equation*}
\begin{split}
\left|
\delta \beta_{k,2,\Delta(0)}
\left[
\vec{V}_{*}
\right]
-
\delta \beta_{k,2,\Delta(0)}
\left[
\vec{V}^{(-\infty,q)}(z\bbone_{\mathbb{Z}_{p}^{3}},0)
\right] 
\right| \le& 4 \mathcal{O}_{26} L^{5} \bar{g}^{2-2\eta}
||\vec{V}^{(-\infty,q)}(z\bbone_{\mathbb{Z}_{p}^{3}},0) - \vec{V}_{*} ||\\
\left|
\xi_{k,\Delta(0)} 
\left[
\vec{V}_{*}
\right] 
- 
\xi_{k,\Delta(0)} 
\left[
\vec{V}^{(-\infty,q)}(z\bbone_{\mathbb{Z}_{p}^{3}},0)
\right]
\right| \le& 2B_{k} \bar{g}^{e_{R}} ||\vec{V}^{(-\infty,q)}(z\bbone_{\mathbb{Z}_{p}^{3}},0) - \vec{V}_{*} ||\ .
\end{split}
\end{equation*}

We also note that by applying the bound of Proposition \ref{contractionbound} $q$-times one has:

\begin{equation*}
\begin{split}
||\vec{V}^{(-\infty,q)}(z\bbone_{\mathbb{Z}_{p}^{3}},0) - \vec{V}_{*} ||
=&
||\vec{V}^{(-\infty,q)}(z\bbone_{\mathbb{Z}_{p}^{3}},0) - \vec{V}^{(-\infty,q)}(0,0) ||\\ 
\le&
\left( \frac{15}{16} \right)^{q} ||\vec{V}^{(-\infty,0)}(z\bbone_{\mathbb{Z}_{p}^{3}},0)
- \vec{V}^{(-\infty,0)}(0,0) ||\\
\le& \left( \frac{15}{16} \right)^{q}\ .
\end{split}
\end{equation*}

Now we note that by \eqref{econst4eq} and \eqref{econst1eq} one has
$e_{R} > e_{4} + 1 \ge 2$. Thus $e_{R} > 2-2\eta$ so in the $\epsilon$
small regime one has:

\begin{gather*}
\left|
\delta \beta_{k,2,\Delta(0)}
[
\vec{V}_{*}
]
-
\delta \beta_{k,2,\Delta(0)}
[
\vec{V}^{(-\infty,q)}(z\bbone_{\mathbb{Z}_{p}^{3}},0)
] 
\right|
+
\left|
\xi_{k,\Delta(0)} 
[
\vec{V}_{*}
] 
- 
\xi_{k,\Delta(0)} 
[
\vec{V}^{(-\infty,q)}(z\bbone_{\mathbb{Z}_{p}^{3}},0)
]
\right|\\
\le 
\left( 4 \mathcal{O}_{26} + 1 \right)
L^{5}
\bar{g}^{2-2\eta}
||\vec{V}^{(-\infty,q)}(z\bbone_{\mathbb{Z}_{p}^{3}},0) - \vec{V}_{*} ||\\
=
\mathbf{O}_{11}
L^{5}
\bar{g}^{2-2\eta}
||\vec{V}^{(-\infty,q)}(z\bbone_{\mathbb{Z}_{p}^{3}},0) - \vec{V}_{*} ||\\
\le
\mathbf{O}_{11}
L^{5}
\bar{g}^{2-2\eta}
\left( \frac{15}{16} \right)^{q}\ .
\end{gather*}

We start with estimating $\beta^{(q,{\rm imp})}_{4}$:

\begin{equation*}
\begin{split}
| \beta^{(q,{\rm imp})}_{4} |
\le&
L^{4[\phi]} 
\sum_{n=0}^{q-1}
L^{-4(q-n)[\phi]} 
\left( 
\left|
\delta \beta_{4,2,\Delta(0)}
[
\vec{V}_{*}
]
-
\delta \beta_{4,2,\Delta(0)}
[
\vec{V}^{(-\infty,n)}(z\bbone_{\mathbb{Z}_{p}^{3}},0)
] 
\right|
\right.\\
& \ \ \ \ \ \ \ \ \ \ \ +
\left.
\left|
\xi_{4,\Delta(0)} 
[
\vec{V}_{*}
] 
- 
\xi_{4,\Delta(0)} 
[
\vec{V}^{(-\infty,n)}(z\bbone_{\mathbb{Z}_{p}^{3}},0)
]
\right|
\right)\\
\le&
L^{4[\phi]} 
\mathbf{O}_{11} 
\times
L^{5}
\times
\bar{g}^{2-2\eta}
\sum_{n=0}^{q-1} L^{-4(q-n)[\phi]}  \left(
  \frac{15}{16} \right)^{n} \\
\le&
L^{4[\phi]} 
\mathbf{O}_{11}
\times 
L^{5}
\bar{g}^{2-2\eta} 
\sum_{n=0}^{q-1} 
\left( 
\frac{15}{16} 
\right)^{(q-n)} 
\left(
\frac{15}{16} 
\right)^{n} \\
\le&
\mathbf{O}_{11} 
\times
q 
L^{5+4[\phi]}  
\bar{g}^{2-2\eta}
\left( \frac{15}{16} \right)^{q}\ .
\end{split}
\end{equation*}

In going to the second to last line we used the fact that for $L \ge
2$ and $\epsilon \le 1$ we have the following inequality : $\displaystyle L^{-4[\phi]} \le
L^{-[\phi]} \le 2^{-\frac{1}{2}}<\left( \frac{15}{16} \right)$. Then by bounding
$L^{5+4[\phi]} \le L^{8}$ we get the desired bound for $\beta^{(q,{\rm imp})}_{4}$.

For $\beta^{(q,{\rm imp})}_{3}$ we have

\begin{equation*}
\begin{split}
|\beta^{(q,{\rm imp})}_{3}| 
\le&
L^{3[\phi]} 
\bigg[
\sum_{n=0}^{q-1}
L^{-3(q-n)[\phi]} 
\left( 
\left|
\delta \beta_{3,2,\Delta(0)}
[
\vec{V}_{*}
]
-
\delta \beta_{3,2,\Delta(0)}
[
\vec{V}^{(-\infty,n)}(z\bbone_{\mathbb{Z}_{p}^{3}},0)
] 
\right|
\right.\\
& \ \ \ \ \ \ \ \ \ \ \ +
\left.
\left|
\xi_{3,\Delta(0)} 
[
\vec{V}_{*}
] 
- 
\xi_{3,\Delta(0)} 
[
\vec{V}^{(-\infty,n)}(z\bbone_{\mathbb{Z}_{p}^{3}},0)
]
\right|
\right)
\bigg]\\
&\ \ \ + 
\left[
\sum_{n=0}^{q-1} 
L^{-3(q-n)[\phi]}  
\left| F_{3,4} 
\left[
\beta^{(n,{\rm imp})}_{4}\bbone_{\Delta(0)} , f^{(n)} 
\right]
\right|
\right] \\
\le& 
L^{3[\phi]} 
\mathbf{O}_{11} 
L^{5} 
\bar{g}^{2-2\eta}
\left[ 
\sum_{n=0}^{q-1} 
L^{-3(q-n)[\phi]} 
\left( 
\frac{15}{16} 
\right)^{n}
\right]
+
\left[
\sum_{n=0}^{q-1} 
L^{-3(q-n)[\phi]}  
\left| F_{3,4} 
\left[
\beta^{(n,{\rm imp})}_{4}\bbone_{\Delta(0)} , f^{(n)} 
\right]
\right|
\right]
\\
\le&
\mathbf{O}_{11} 
q 
L^{5+3[\phi]} 
\bar{g}^{2-2\eta}
\left( \frac{15}{16} \right)^{q} 
+ 
\left[
\sum_{n=0}^{q-1} 
L^{-3(q-n)[\phi]}  
\left| F_{3,4} 
\left[
\beta^{(n,{\rm imp})}_{4}\bbone_{\Delta(0)} , f^{(n)} 
\right]
\right|
\right]\ .
\end{split}
\end{equation*} 

In the above expressions the first term was bounded just as it was for
$\beta^{(q,{\rm imp})}_{4}$. We now try to estimate the summands appearing inside of
the second term. We will use $|z|\le 1$.

\begin{equation*}
\begin{split}
\left| F_{3,4} 
\left[
\beta^{(n,{\rm imp})}_{4}\bbone_{\Delta(0)} , f^{(n)} 
\right] 
\right|
\le&
4 L^{-n[\phi]} |\beta^{(n,{\rm imp})}_{4}| \times |\Gamma(0)|\\
\le& 
8
L^{-n[\phi]} 
\mathbf{O}_{11} 
n
L^{8} \bar{g}^{2-2\eta}
\left( \frac{15}{16} \right)^{n}\\
\le&
16 
\mathbf{O}_{11} 
L^{8} 
\bar{g}^{2-2\eta}
\left( \frac{15}{16} \right)^{n}\ .
\end{split}
\end{equation*}

In going to the second to last line we used the bound $|\Gamma(0)| \le
||\Gamma||_{L^{\infty}} \le 2$. In going to the last line note that
for $\epsilon \le 1$ and $L \ge 2$ one has $\displaystyle nL^{-n[\phi]}
\le n2^{-\frac{n}{2}} \le \frac{2}{e\times\log(2)} \le 2$. Inserting this
into our previous inequality gives us

\begin{equation*}
\begin{split}
|\beta^{(q,{\rm imp})}_{3}| 
\le&
\mathbf{O}_{11} 
q 
L^{5+3[\phi]} 
\bar{g}^{2-2\eta}
\left( \frac{15}{16} \right)^{q} 
+
16 
\mathbf{O}_{11} 
L^{8} 
\bar{g}^{2-2\eta}
\sum_{n=0}^{q-1} 
\left[
L^{-3(q-n)[\phi]} 
\left( \frac{15}{16} \right)^{n}
\right]\\
\le&
\mathbf{O}_{11} 
q 
L^{5+3[\phi]} 
\bar{g}^{2-2\eta}
\left( \frac{15}{16} \right)^{q}
+
16
\mathbf{O}_{11} 
q
L^{8}
\bar{g}^{2-2\eta}
\left( \frac{15}{16} \right)^{q}\\
\le&
17
\mathbf{O}_{11}
q
L^{8}
\bar{g}^{2-2\eta}
\left( \frac{15}{16} \right)^{q}\ .
\end{split}
\end{equation*} 

Not that in going to the second line we used the bound
$L^{-3(q-n)[\phi]} \le \left(\frac{15}{16}\right)^{(q-n)}$.\\

We start on $\beta_{2}^{(q,{\rm imp})}$ by making the following estimates:

\begin{equation*}
\begin{split}
\left| F_{2,4} \left[ \beta^{(n,{\rm imp})}_{4}\bbone_{\Delta(0)}, f^{(n)} \right] \right|
\le&
6
\times
\left| \beta^{(n,{\rm imp})}_{4} \right| 
\times 
\Gamma(0)^{2}
\times
L^{-2n[\phi]}\\
\le&
24
\times
\mathbf{O}_{11} 
n
L^{8} 
\bar{g}^{2-2\eta} 
\left( \frac{15}{16} \right)^{n}
L^{-2n[\phi]}\\
\le&
48
\times
\mathbf{O}_{11} 
L^{8} 
\bar{g}^{2-2\eta} 
\left( \frac{15}{16} \right)^{n}\ .
\end{split}
\end{equation*}

Similarly one gets the bound

\[
\left| F_{2,3} \left[ \beta^{(n,{\rm imp})}_{3}\bbone_{\Delta(0)}, f^{(n)} \right] \right| 
\le
204
\times
\mathbf{O}_{11}
L^{6}
\bar{g}^{2-2\eta}
\left( \frac{15}{16} \right)^{n}\ .
\]

The bound for $\beta^{(q,{\rm imp})}_{2}$ then proceeds along familiar
lines. One uses the same arguments to prove the estimate for
$\beta^{(q,{\rm imp})}_{1}$. In particular

\begin{equation*}
\begin{split}
\left|
F_{1,4} \left[
\beta^{(n,{\rm imp})}_{4}\bbone_{\Delta(0)}, f^{(n)}
\right]
\right|
\le&
64
\times
\mathbf{O}_{11}
L^{8}
\bar{g}^{2-2\eta}
\left( \frac{15}{16} \right)^{n}\\
\left|
F_{1,3} \left[
\beta^{(n,{\rm imp})}_{3}\bbone_{\Delta(0)}, f^{(n)}
\right]
\right|
\le&
408
\times
\mathbf{O}_{11}
L^{8}
\bar{g}^{2-2\eta}
\left( \frac{15}{16} \right)^{n}\\
\left|
F_{1,2} \left[
\beta^{(n,{\rm imp})}_{2}\bbone_{\Delta(0)}, f^{(n)}
\right]
\right|
\le&
2024
\times
\mathbf{O}_{11}
L^{8}
\bar{g}^{2-2\eta}
\left( \frac{15}{16} \right)^{n}\ .
\end{split}
\end{equation*}

To bound $\delta b^{(q,{\rm imp})}$ we first make the following
estimate. For $k = 1,2,3,4$ one has:

\begin{equation*}
\begin{split}
\left|
F_{0,k} \left[
\beta^{(q,{\rm imp})}_{k}\bbone_{\Delta(0)}, f^{(q)}
\right]
\right|
\le&
L^{-kq[\phi]}
\times
\Gamma(0)^{k}
\times
|\beta_{k}^{(q,{\rm imp})}|\\
\le&
L^{-q[\phi]}
\times
2^{4}
\times
2497
\times
\mathbf{O}_{11}
\times
q
L^{8}
\bar{g}^{2-2\eta}
\left( \frac{15}{16} \right)^{q}\\
\le&
79904
\times
\mathbf{O}_{11}
L^{8}
\bar{g}^{2-2\eta}
\left( \frac{15}{16} \right)^{q}\ .
\end{split}
\end{equation*}

We then have

\begin{equation*}
\begin{split}
|\delta b^{(q,{\rm imp})}|
\le&
\left|
\delta \beta_{0,2,\Delta(0)}
[
V_{*}
]
-
\delta \beta_{0,2,\Delta(0)}
[
\vec{V}^{(-\infty,q)}(z\bbone_{\mathbb{Z}_{p}^{3}},0)
] 
\right|\\
& +
\left|
\xi_{0,\Delta(0)} 
[
V_{*}
] 
- 
\xi_{0,\Delta(0)} 
[
\vec{V}^{(-\infty,q)}(z\bbone_{\mathbb{Z}_{p}^{3}},0)
]
\right|
\\
&+ 
\left[
\sum_{k=1}^{4} 
\left| F_{0,k} 
\left[
\beta^{(q,{\rm imp})}_{k}\bbone_{\Delta(0)} , f^{(q)} 
\right]
\right|
\right]\\
\le& 
\mathbf{O}_{11}
\times
L^{5}
\bar{g}^{2-2\eta}
\left( \frac{15}{16} \right)^{q}
+
4
\times
79904
\times
\mathbf{O}_{11}
\times
L^{8}
\bar{g}^{2-2\eta}
\left( \frac{15}{16} \right)^{q}\\
\le &
319617
\times
\mathbf{O}_{11}
\times
L^{8}
\bar{g}^{2-2\eta}
\left( \frac{15}{16} \right)^{q}\ .
\end{split}
\end{equation*}

This gives the desired bound. \qed \\

\begin{Lemma}
In the $\epsilon$ small regime
and on the domain $\{ z \in \mathbb{C}\ |\ |z| < \mathbf{O}_{10} \}$
one has the decomposition
\[
\mathcal{S}^{\rm T}(z\bbone_{\mathbb{Z}_p^3}, 0)
= S^{\rm T,exp}(z) + S^{\rm T,imp}(z)\ .
\]

All three of the above functions are analytic on the above
domain. Additionally, over this domain one has the following explicit formula

\begin{equation*}
\begin{split}
\mathcal{S}^{\rm T,exp}(z)
=& 
-\sum_{q =0}^{\infty}
\Bigg\{ z^{4} g_{*}
\left[ L^{-4q[\phi]}
\left( \int_{\mathbb{Q}_{p}^{3}} {\rm d}^{3}x\ \Gamma(x)^{4} \right)
+ 6 L^{-4q[\phi]}q ||\Gamma||_{L^{2}}^{2} \times \Gamma(0)^{2}
\right. \\
 & +12L^{-2q[\phi]}\left(\sum_{n=0}^{q-1} n L^{-2n[\phi]}
 \right) ||\Gamma||_{L^{2}}^{2} \times \Gamma(0)^{2} \\
 & + 
\left.
4L^{-2q[\phi]}\frac{1-L^{-2q[\phi]}}{1-L^{-2[\phi]}}
\Gamma(0) \left( \int_{\mathbb{Q}_{p}^{3}} {\rm d}^{3}x\ \Gamma(x)^{3}\right)
\right] + z^{2} \mu_{*} L^{-2q[\phi]} ||\Gamma||_{L^{2}}^{2}
\Bigg\} \\
 & +\frac{z^2}{2}
\left(\bbone_{\mathbb{Z}_{p}^{3}}, C_{-\infty} \bbone_{\mathbb{Z}_{p}^{3}} \right)
\end{split}
\end{equation*} 

and the following uniform bound

\[
|\mathcal{S}^{\rm T,imp}(z)|
\le
\mathbf{O}_{13}
L^{8}
\bar{g}^{2-2\eta}.
\]

where $\mathbf{O}_{13}=16 \times \mathbf{O}_{12}$.
\end{Lemma}

\noindent{\bf Proof:}
From earlier definitions we have that

\[
\mathcal{S}^{\rm T}(z\bbone_{\mathbb{Z}_p^3}, 0)
=
\frac{z^2}{2} \left( \bbone_{\mathbb{Z}_{p}^{3}},
  C_{-\infty} \bbone_{\mathbb{Z}_{p}^{3}} \right)
+
\sum_{q = 0}^{\infty}
\left( \delta b^{(q,{\rm exp})} + \delta b^{(q,{\rm imp})} \right)\ .
\]

We define

\[
\mathcal{S}^{\rm T,exp}(z)
=
\frac{z^2}{2} \left( \bbone_{\mathbb{Z}_{p}^{3}},
  C_{-\infty} \bbone_{\mathbb{Z}_{p}^{3}} \right)
+
\sum_{q = 0}^{\infty} \delta b^{(q,{\rm exp})} \ .
\]

\[
\mathcal{S}^{\rm T,imp}(z)
=
\sum_{q = 0}^{\infty} \delta b^{(q,{\rm imp})}\ .
\]
 
The explicit formula given for $\mathcal{S}^{\rm T,exp}(z)$ comes from
substitution of the explicit formula for the $\delta b^{(q,{\rm exp})}$ from
Lemma \ref{betaexp}. Since $[\phi] > 0$ for $\epsilon \in (0,1]$
it is not hard to see that infinite sum in the expression for
$\mathcal{S}^{\rm T,exp}(z)$ is uniformly absolutely summable on our
domain. Analyticity follows from the explicit formula.\\

On the other hand we have

\begin{equation*}
\begin{split}
|\mathcal{S}^{\rm T,imp}(z)|
\le&
\sum_{q=0}^{\infty} | \delta b^{(q,{\rm imp})}|\\
\le& 
\mathbf{O}_{12} 
\times
L^{6}
\times
\bar{g}^{2-2\eta}
\sum_{q=0}^{\infty}
\left( \frac{15}{16} \right)^{q}\\
\le&
16
\times
\mathbf{O}_{12} 
\times
L^{8}
\times
\bar{g}^{2-2\eta}\ .
\end{split}
\end{equation*}

We have then proved the desired uniform bound and we have uniform
absolute convergence yielding analyticity
as well. \qed

\begin{Lemma}
In the small $\epsilon$ regime one has

\[
\left|
\frac{d^2} {dz^2} \Big|_{z=0}
\mathcal{S}^{\rm T}(z\bbone_{\mathbb{Z}_p^3}, 0)
- U_{2}
\right|
\le \mathbf{O}_{14} L^{8} \bar{g}^{2-2\eta}
\]

where 

\[ 
U_{2} 
=
\left( \bbone_{\mathbb{Z}_{p}^{3}}, C_{-\infty}
  \bbone_{\mathbb{Z}_{p}^{3}} \right)
-
2 ||\Gamma||_{L^{2}}^{2} \times \frac{1}{1-L^{-2[\phi]}} \times
\mu_{*}
\]

and

\[
\left|
\frac{d^4} {dz^4} \Big|_{z=0}
\mathcal{S}^{\rm T}(z\bbone_{\mathbb{Z}_p^3}, 0)
- U_{4}
\right|
\le \mathbf{O}_{15}L^{8}\bar{g}^{2-2\eta}
\]

where

\[
\begin{split}
U_{4}
= & 
-24 g_{*} \sum_{q=0}^{\infty} \left[
 L^{-4q[\phi]} \left( \int_{\mathbb{Q}_{p}^{3}} {\rm d}^{3}x\ \Gamma(x)^{4} \right)
+ 
6 L^{-4q[\phi]}q ||\Gamma||_{L^{2}}^{2} \Gamma(0)^{2} 
+12L^{-2q[\phi]}\left(\sum_{n=0}^{q-1} n L^{-2n[\phi]}\right)||\Gamma||_{L^{2}}^{2} \Gamma(0)^{2} 
\right. \\
 & + \left.
4L^{-2q[\phi]}\frac{1-L^{-2q[\phi]}}{1-L^{-2[\phi]}}
\Gamma(0) \left( \int_{\mathbb{Q}_{p}^{3}} {\rm d}^{3}x\ \Gamma(x)^{3}
\right) \right].
\end{split}
\]

Here we have used the following numerical constants: $\mathbf{O}_{14}
= 8 \times \mathbf{O}_{10}^{-2} \mathbf{O}_{13}$ and
$\mathbf{O}_{15} = 384 \times \mathbf{O}_{10}^{-4} \mathbf{O}_{13}$.

\end{Lemma}

\noindent{\bf Proof:}
We note that for $j = 2,4$ we have that $\displaystyle U_{j} =
\frac{d^j} {dz^j} \Big|_{z=0} \mathcal{S}^{\rm T,exp}(z)$.\\

 By the previous
lemma the bounds above will follow if we have the necessary bounds on
$\displaystyle \left| \frac{d^j} {dz^j} \Big|_{z=0}  \mathcal{S}^{\rm T,imp}(z)
\right|$. By Cauchy's formula we have

\[
\frac{d^j} {dz^j} \Big|_{z=0} \mathcal{S}^{\rm T,imp}(z) 
= 
\frac{j!}{2 i \pi} 
\oint \frac{d \lambda}{\lambda^{j+1}} \ \mathcal{S}^{\rm T,imp}(\lambda)
\]

Here we are integrating around the contour $\displaystyle |\lambda| =
\frac{1}{2} \mathbf{O}_{10}$. Utilizing the uniform bound on
$\mathcal{S}^{\rm T,imp}(z)$ from the previous lemma we get the estimate:

\[
\left|
\frac{d^j} {dz^j} \Big|_{z=0} \mathcal{S}^{\rm T,imp}(z)
\right|
\le
j!\times 2^{j}  \mathbf{O}_{10}^{-j} \times
\mathbf{O}_{13} \times L^{8} \times \bar{g}^{2-2\eta}\ .
\]

This proves the lemma. \qed\\

\begin{Proposition}\label{nontriviality}
In the small $\epsilon$ regime

\[
\frac{d^4} {dz^4} \Big|_{z=0}
\mathcal{S}^{\rm T}(z\bbone_{\mathbb{Z}_p^3}, 0) \le -\frac{1}{4} \bar{g} < 0\ .
\]

\end{Proposition}

\noindent{\bf Proof:}
We observe that since $\hat{\Gamma}(k) \ge 0$ one has

\begin{gather*}
\Gamma(0) = \int_{\mathbb{Q}_{p}^{3}}{\rm d}^{3}k\ \hat{\Gamma}(k) \ge 0 \\
\int_{\mathbb{Q}_{p}^{3}} {\rm d}^{3}x\ \Gamma(x)^{3} = \  \left(
  \hat{\Gamma} \ast \hat{\Gamma} \ast \hat{\Gamma} \right)(0) \ge 0\ .
\end{gather*}

In the above expression $\ast$ denotes convolution. It then
follows by only keeping the first $q=0$ term that

\begin{equation*}
\begin{split}
U_{4} \le& -24 g_{*} \int_{\mathbb{Q}_{p}^{3}} {\rm d}^{3}x\ \Gamma(x)^{4}\\
\le& -24g_{*}  \int_{\mathbb{Z}_{p}^{3}} {\rm d}^{3}x\ \Gamma(x)^{4}\\
=& -24g_{*} \Gamma(0)^{4}\\
=& -24 g_{*}  \times \left[ \frac{1 - p^{-3}}{1 - p^{-2[\phi]}} \left( 1 - L^{-2[\phi]} \right)\right]^{4}\ .
\end{split}
\end{equation*}

In going to the last line we used Lemma \ref{L1lem}. Now we note that
$p, L \ge 2$ and $\epsilon \le 1$ implies that $-2[\phi] \le -1$

\begin{equation*}
\begin{split}
U_{4} \le& -24 g_{*} \times \left[\frac{ 1 - \frac{1}{2^3}}{1} \times \left(1 -
\frac{1}{2} \right)\right]^{4} \\
=& - 24\left(\frac{7}{16}\right)^4 g_{*}\\
\le& -12\left(\frac{7}{16}\right)^4 \bar{g} \\
\le & -\frac{1}{3} \bar{g}\ .
\end{split}
\end{equation*}

Note that in going to the third line we used that $g_{*} > \frac{1}{2} \bar{g}$. Now using the previous lemma we have:

\begin{equation*}
\begin{split}
\frac{d^4} {dz^4} \Big|_{z=0}
\mathcal{S}^{\rm T}(z\bbone_{\mathbb{Z}_p^3}, 0) \le& U_{4} +
\mathbf{O}_{15}L^{8}\bar{g}^{2-2\eta}\\
\le& -\frac{1}{3} \bar{g} +
\mathbf{O}_{15}L^{8}\bar{g}^{2-2\eta}\ .
\end{split}
\end{equation*}

Since $2-2\eta > e_4\ge 1$
we can take $\epsilon$ sufficiently
small to guarantee that
$\mathbf{O}_{15}L^{8}\bar{g}^{2-2\eta} \le \frac{1}{12}
\bar{g}$. This proves the proposition. \qed

\subsection{The two-point function for the composite field}

We now study the $\phi^{2}$ correlation when smeared with the
characteristic function of $\mathbb{Z}_{p}^{3}$, i.e., the quantity

\begin{equation*}
\begin{split}
\frac{d^2}{dz^2} \Big|_{z=0}
\mathcal{S}^{T}(0,z\bbone_{\mathbb{Z}_{p}^{3}})
=&\ 
\langle N[\phi^{2}](\bbone_{\mathbb{Z}_{p}^{3}})^{2} \rangle
-
\langle N[\phi^{2}](\bbone_{\mathbb{Z}_{p}^{3}})\rangle^{2}\\
=&\ 
\langle N[\phi^{2}](\bbone_{\mathbb{Z}_{p}^{3}})^{2}\rangle
\end{split}
\end{equation*}

since the one-point function is identically zero.\\

Here $q_{-}=q_{+}=0$ so there is no contribution from the middle
regime. Thus

\begin{gather*}
\langle N[\phi^{2}](\bbone_{\mathbb{Z}_{p}^{3}})^{2}\rangle
=
\langle N[\phi^{2}](\bbone_{\mathbb{Z}_{p}^{3}})^{2}\rangle^{\rm{UV}}
+
\langle
N[\phi^{2}](\bbone_{\mathbb{Z}_{p}^{3}})^{2}\rangle^{\rm{IR}}\\[1.2ex]
\textnormal{where  }
\langle
N[\phi^{2}](\bbone_{\mathbb{Z}_{p}^{3}})^{2}\rangle^{\rm{UV}}
=
\frac{d^2}{dz^2} \Big|_{z=0}
\mathcal{S}^{\rm{T},\rm{UV}} ( 0, z
\bbone_{\mathbb{Z}_{p}^{3}})\\[1.2ex]
\textnormal{and  }
 \langle
N[\phi^{2}](\bbone_{\mathbb{Z}_{p}^{3}})^{2}\rangle^{\rm{IR}}
=
\frac{d^2}{dz^2} \Big|_{z=0}
\mathcal{S}^{\rm{T},\rm{IR}} ( 0, z
\bbone_{\mathbb{Z}_{p}^{3}}).
\end{gather*}

Clearly, since we can derive term-by-term in the sum over $q$ and
since the constant and linear parts disappear

\begin{equation*}
\begin{split}
\langle
N[\phi^{2}](\bbone_{\mathbb{Z}_{p}^{3}})^{2}\rangle^{\rm{UV}}
=&
\frac{d^2}{dz^2} \Big|_{z=0}
\sum_{q < 0} L^{-3(q+1)} 
\delta b 
\left[ 
\Psi( v, - \alpha_{\rm{u}}^{q}Y_{2}ze_{\phi^{2}})
\right]\\
=&\ 
Y_{2}^{2}
\times
\left(
\sum_{q < 0} L^{-3(q+1)} 
\alpha_{\rm{u}}^{2q}
\right)
\times
\frac{d^2}{dz^2} \Big|_{z=0}
\delta b 
\left[ 
\Psi( v, ze_{\phi^{2}})
\right]
\end{split}
\end{equation*}

by the chain rule. This also uses $L^{3} \alpha_{\rm{u}}^{-2} < 1$ which will be proved
shortly.\\

We will use the more convenient notation $\displaystyle \Psi_{v}(w)$
instead of $\displaystyle \Psi(v,w)$.\\

Now for $w$ small we have by Theorem \ref{unnamedtheorem1}

\[
\Psi_{v}(w)
=
\Psi_{v_{*}}
\left(
T_{\infty}
(v)
[w]
\right).
\]

By the remark following Lemma \ref{Tcvlem}

\[
P_{s}
T_{\infty}(v)
[
e_{\phi^{2}}
]
=
0
\]

i.e. $\displaystyle T_{\infty}(v)[e_{\phi^{2}}]$ is in
$\mathcal{E}^{\rm{u}}$ and therefore is proportional to
$e_{\rm{u}}$.\\

We define $\displaystyle \varkappa_{\phi^{2}}$ as the proportionality constant,
i.e., by

\[
T_{\infty}
(v)
[
e_{\phi^{2}}
]
=
\varkappa_{\phi^{2}}
e_{\rm{u}}.
\]

Hence

\[
\Psi
(v,ze_{\phi^{2}})
=
\Psi_{v_{*}}
(
z
\varkappa_{\phi^{2}}
e_{\rm{u}}
)
\]

and as a result

\[
\langle 
N[\phi^{2}](\bbone_{\mathbb{Z}_{p}^{3}})^{2}
\rangle^{\rm{UV}}
= 
Y_{2}^{2}
\varkappa^{2}_{\phi^{2}}
\langle
N[\phi^{2}](\bbone_{\mathbb{Z}_{p}^{3}})^2
\rangle^{\rm{UV}}_{\rm{reduced}}
\]

with

\[
\langle
N[\phi^{2}](\bbone_{\mathbb{Z}_{p}^{3}})^2
\rangle^{\rm{UV}}_{\rm{reduced}}
=
\frac{1}{\alpha_{\rm{u}}^{2} - L^{3}}
\times
D_{0}^{2}
(
\delta b
\circ
\Psi_{v_{*}}
)
[e_{\rm{u}}, e_{\rm{u}}].
\]

On the other hand we easily see that

\[
\langle
N[\phi^{2}]
(\bbone_{\mathbb{Z}_{p}^{3}})^2
\rangle^{\rm{IR}}
=
\sum_{q \ge 0}
\frac{d^2}{dz^2} \Big|_{z=0}
\delta b_{\Delta(0)}
\left[
RG^{q}_{\rm{ex}}
\left(
\vec{V}^{(-\infty,0)}(0,z \bbone_{\mathbb{Z}_{p}^{3}})
\right)
\right]
\]

where

\[
\vec{V}^{(-\infty,0)}(0,z \bbone_{\mathbb{Z}_{p}^{3}})
=
\mathcal{J}_{0}
\left(
0,
(
\Psi_{v}
(-Y_{2}ze_{\phi^{2}})
),
v_{*}
\right).
\]\\

We define the affine isometric map $\varpi: \mathcal{E} \rightarrow
\mathcal{E}_{\rm{pt}}$ which sends $v = (\delta g, \mu, R)$ to
$\vec{V}= ( V_{\Delta} )_{\Delta \in \mathbb{L}} = \varpi(v)$ such
that

\[ 
V_{\Delta} 
= 
(
\beta_{4,\Delta}, 
\beta_{3,\Delta}, 
\beta_{2,\Delta},
\beta_{1,\Delta}, 
W_{5,\Delta}, 
W_{6,\Delta}, 
f_{\Delta}, 
R_{\Delta}
)
\]

is zero for $\Delta \not = \Delta(0)$ and equal to 

\[
(
\delta g - \delta g_{*},
0,
\mu - \mu_{*},
0,
0,
0,
0,
R - R_{*}
)
\]

for $\Delta = \Delta(0)$.\\

It easily follows from the definitions that

\begin{equation*}
\begin{split}
\vec{V}^{(-\infty,0)}(0,z\bbone_{\mathbb{Z}_{p}^{3}})
=&\ 
\iota ( v_{*} )
+
\varpi \circ \Psi_{v}( -Y_{2} z e_{\phi^{2}} )\\
=&\ 
\iota(v_{*})
+
\varpi \circ \Psi_{v_{*}}
(
-Y_{2}
\varkappa_{\phi^{2}}
z
e_{\rm{u}}
)
\end{split}
\end{equation*}

for $z$ small.\\

Hence by the chain rule

\[
\langle
N[\phi^{2}]
(\bbone_{\mathbb{Z}_{p}^{3}})^2
\rangle^{\rm{IR}}
=
Y_{2}^{2} \varkappa_{\phi^{2}}^{2}
\langle
N[\phi^{2}]
(\bbone_{\mathbb{Z}_{p}^{3}})^2
\rangle^{\rm{IR}}_{\rm{reduced}}
\]

where

\[
\langle
N[\phi^{2}]
(\bbone_{\mathbb{Z}_{p}^{3}})^2
\rangle^{\rm{IR}}_{\rm{reduced}}
=
\sum_{q \ge 0}
\frac{d^2}{dz^{2}} \Big|_{z=0}
\delta b_{\Delta(0)}
\left[
\iota(v_{*})
+
RG^{q}_{{\rm dv},\iota(v_{*})}
\circ
\varpi
\circ
\Psi_{v_{*}}
( z e_{\rm{u}})
\right]
\]

where we introduced the more convenient notation $\displaystyle
RG_{{\rm dv},\vec{V}_{\rm{bk}}}[\dot{V}]$ for $\displaystyle
RG_{\rm{dv}}[\vec{V}_{\rm{bk}}, \dot{V} ]$ of section \S\ref{infregsec}.\\

In what follows we will show that when $\epsilon \rightarrow 0$, $\langle
N[\phi^{2}]
(\bbone_{\mathbb{Z}_{p}^{3}})^2
\rangle^{\rm{IR}}_{\rm{reduced}}$ remains bounded while $\langle
N[\phi^{2}]
(\bbone_{\mathbb{Z}_{p}^{3}})^2
\rangle^{\rm{UV}}$ blows up.\\

This will need new constraints on the norm exponents which we list
redundantly as:

\begin{equation*}
\begin{split}
e_{4} &> 1\\
e_{4} &> e_{2}\\
e_{4} &> 2 e_{2} - 1\\
e_{R} &> e_{2} + 1\\
e_{R} &> 2e_{2}\ .
\end{split}
\end{equation*}

We first introduce the subspace $\mathcal{E}_{\rm{ex},\rm{ev}}$ of
$\mathcal{E}_{\rm{ex}}$.\\

It is the space of vectors

\[
(
\beta_{4,\Delta}, 
\beta_{3,\Delta}, 
\beta_{2,\Delta},
\beta_{1,\Delta}, 
W_{5,\Delta}, 
W_{6,\Delta}, 
f_{\Delta}, 
R_{\Delta}
)_{\Delta \in \mathbb{L}}
\]

such that for all $\Delta \in \mathbb{L}$,

\begin{gather*}
\beta_{3,\Delta}
=
\beta_{1,\Delta}
=
W_{5,\Delta}
=
W_{6,\Delta}
=
f_{\Delta}
=
0\\
\textnormal{and  }
R_{\Delta} 
\in
C^{9}_{\rm{bd},\rm{ev}}(\mathbb{R}, \mathbb{C}).
\end{gather*}

Using the same line of reasoning as in the proof of Proposition
\ref{L53prop} or in \S \ref{algconsec} it is easy to see that
$\mathcal{E}_{\rm{ex},\rm{ev}}$ is invariant by $RG_{\rm{ex}}$.

\begin{Lemma}

In the small $\epsilon$ regime and for $\vec{V} \in B(\bar{V},
\frac{1}{2}) \cap \mathcal{E}_{\rm{ex},\rm{dv}}$\\
we have for all $\Delta' \in \mathbb{L}$

\[
| \delta b_{\Delta'} [ \vec{V} ] |
\le
\mathbf{O}_{16} 
L^{5} 
\bar{g}^{2e_{2}}
\]

where

\[
\mathbf{O}_{16}
=
1
+
9
\sum_{a_{1},a_{2},m}
\bbone
\left\{
\begin{array}{c}
a_{i}+m = 2 \textnormal{ or } 4\\
a_{i} \ge 0,\  m \ge 1
\end{array}
\right\}
\times
C(a_{1},a_{2}|0)
\times
2^{\frac{a_{1}+a_{2}}{2}}.
\]

\end{Lemma}

\noindent{\bf Proof: }
Recall that

\[
\delta
b_{\Delta'}
[\vec{V}]
=
\delta
\beta_{0,1,\Delta'}
+
\delta\beta_{0,2,\Delta'}
+
\xi_{0,\Delta'}(\vec{V}).
\]
Since there are no $f$'s we have $\beta_{0,1,\Delta'} = 0$. Similarly
the $\delta\beta_{0,2,\Delta'}$ contirbution reduces to

\begin{gather*}
\delta\beta_{0,2,\Delta'}
=
\sum_{a_{1},a_{2},m}
\bbone
\left\{
\begin{array}{c}
a_{i}+m = 2 \textnormal{ or } 4\\
a_{i} \ge 0,\  m \ge 1
\end{array}
\right\}
\frac{(a_{1}+m)!(a_{2}+m)!}{a_{1}!a_{2}!m!}
\times
\frac{1}{2}
C(a_{1},a_{2}|0)\\
\times
L^{-(a_{1}+a_{2})[\phi]}
C_{0}(0)^{\frac{a_{1}+a_{2}}{2}}
\times
\int_{(L^{-1}\Delta')^{2}}
{\rm d}^{3}x_{1}
{\rm d}^{3}x_{2}\ 
\beta_{a_{1}+m}(x_{1})
\beta_{a_{2}+m}(x_{2})\ 
\Gamma(x_{1}-x_{2})^{m}.
\end{gather*}

We use the bound

\[
\begin{split}
\left|
\int_{(L^{-1}\Delta')^{2}}
{\rm d}^{3}x_{1}
{\rm d}^{3}x_{2}\ 
\beta_{a_{1}+m}(x_{1})
\beta_{a_{2}+m}(x_{2})\ 
\Gamma(x_{1}-x_{2})^{m}
\right|
\le &
L^{3}
||\Gamma||_{L^{\infty}}^{m-1}
\times
||\Gamma||_{L^{1}} \\
 & \times
\sup_{x \in L^{-1}\Delta'}
|\beta_{a_{1}+m}(x)|
\times
\sup_{x \in L^{-1}\Delta'}
|\beta_{a_{2}+m}(x)|.
\end{split}
\]
We bound the supremums by noting that $\beta_{a_{1}+m}$ can only be
$\beta_{2}$ or $\beta_{4}$.
Since $\bar{V}$ has no $\beta_{2}$ component

\[
|\beta_{2}(x)|
\le
||\vec{V} - \bar{V}||
\bar{g}^{e_{2}}
\le
\frac{1}{2}
\bar{g}^{e_{2}}.
\]
On the other hand
\[
|\beta_{4}(x)| 
\le 
\bar{g} 
+
||\vec{V} - \bar{V}|| \bar{g}^{e_{4}}
\le
\frac{3}{2}
\bar{g}
\le
\frac{3}{2}
\bar{g}^{e_{2}}
\]
since $e_{4} > 1 \ge e_{2}$.
As a result the previous integral is bounded by

\[
L^{3}
||\Gamma||_{L^{\infty}}^{m-1}
\times
||\Gamma||_{L^{1}}
\times
\frac{9}{4}
\bar{g}^{2e_{2}}
\le
L^{3}
\times
2^{m-1}
\times
\frac{1}{\sqrt{2}}
L^{3-2[\phi]}
\times
\frac{9}{4}
\bar{g}^{2e_{2}}
\le
18L^{5}\bar{g}^{2e_{2}}
\]
where we used $\epsilon \le 1$ so $3-2[\phi] \le 2$, and $m \le 4$
while dropping $\sqrt{2}$.
Finally $|\xi_{0,\Delta'}(\vec{V})| \le \frac{1}{2} B_{4}
\bar{g}^{e_{R}}$ by Theorem \ref{mainestthm}. Noting that
$\frac{1}{2}B_{4}\bar{g}^{e_{R}-2e_{2}} \le 1$ for $\epsilon$
small the lemma follows. \qed\\

\begin{Lemma}
For $\vec{V}^{1},\vec{V}^{2} \in \bar{B}(\bar{V},\frac{1}{6}) \cap
\mathcal{E}_{\rm{ex},\rm{dv}}$,
we have the Lipschitz estimate

\[
\left|
\delta b_{\Delta(0)} [ \vec{V}^{1} ]
-
\delta b_{\Delta(0)} [ \vec{V}^{2}]
\right|
\le
4
\mathbf{O}_{16}
L^{5}
\bar{g}^{2e_{2}}
|| \vec{V}^{1} - \vec{V}^{2}||\ .
\]

\end{Lemma}

\noindent{\bf Proof: }
This is an immediate consequence  of the previous lemma and Lemma
\ref{Lipschitzlem} with $\nu = \frac{1}{2}$. \qed\\

Since we are computing second derivatives there is no harm in writing

\[
\langle
N[\phi^{2}]
(\bbone_{\mathbb{Z}_{p}^{3}})^2
\rangle^{\rm{IR}}_{\rm{reduced}}
=
\sum_{q \ge 0}
\frac{d^2}{dz^2} \Big|_{z=0}
\left\{
\delta b_{\Delta(0)}
\left[
\iota(v_{*})
+
RG^{q}_{\rm{dv},\iota(v_{*})}
\circ
\varpi
\circ
\Psi_{v_{*}}
(
z
e_{\rm{u}}
)
\right]
-
\delta b_{\Delta(0)}
\left[
\iota(v_{*})
\right]
\right\}.
\]

If $z$ is small enough so that

\[
||
\Psi(
z e_{\rm{u}}
)
-
v_{*}
||
\le
\mathbf{O}_{10}
\]
which is the same as saying that $|| \varpi \circ
\Psi_{v_{*}}(ze_{\rm{u}})|| \le \mathbf{O}_{10}$, then Proposition
\ref{contractionbound} along with the last lemma will imply

\[
\left|
\delta b_{\Delta(0)}
\left[
\iota(v_{*})
+
RG^{q}_{\rm{dv},\iota(v_{*})}
\circ
\varpi
\circ
\Psi_{v_{*}}
(
z
e_{\rm{u}}
)
\right]
-
\delta b_{\Delta(0)}
\left[
\iota(v_{*})
\right]
\right|
\le
4
\mathbf{O}_{16}
L^{5}
\bar{g}^{2e_{2}}
\left(
\frac{15}{16}
\right)^q
\times
\mathbf{O}_{10}.
\]

Let $z_{\rm{max}} > 0$ be such that $|z| \le z_{\rm{max}}$ implies
$||\Psi_{v_{*}}(ze_{\rm{u}}) - v_{*}|| \le \mathbf{O}_{10}$.
Then by extracting the derivatives with Cauchy's formula we easily
arrive at the bound

\[
\left|
\langle
N[\phi^{2}]
(\bbone_{\mathbb{Z}_{p}^{3}})^2
\rangle^{\rm{IR}}_{\rm{reduced}}
\right|
\le
4
\mathbf{O}_{10}
\mathbf{O}_{16}
L^{5}
\bar{g}^{2e_{2}}
\times
\frac{1}{1-\frac{15}{16}}
\times
2!\times
z_{\rm{max}}^{-2}\ .
\]

Now from Theorem \ref{unnamedtheorem1bis} $||z e_{\rm{u}}|| < \frac{1}{24}$ implies 

\begin{equation*}
\begin{split}
||\Psi_{v_{*}}(ze_{\rm{u}}) - v_{*}|| 
\le&
||ze_{\rm{u}}||
\left(
1
+
\frac{17}{18}
\times
\frac{1}{24}
\right)\\
\le&
2
||ze_{\rm{u}}||
\end{split}
\end{equation*}
for simplicity. So $z_{\rm{max}} = \frac{1}{2} \mathbf{O}_{10}
||e_{\rm{u}}||^{-1}$ works because $\frac{1}{2}\mathbf{O}_{10} \le
\frac{1}{80} < \frac{1}{24}$. Also by lemma \ref{L65lem},
$||e_{\rm{u}}|| = \bar{g}^{-e_{2}}$.
Hence in the small $\epsilon$ regime we have the bound

\[
\left|
\langle
N[\phi^{2}]
(\bbone_{\mathbb{Z}_{p}^{3}})^2
\rangle^{\rm{IR}}_{\rm{reduced}}
\right|
\le
512
\mathbf{O}_{10}^{-1}
\times
\mathbf{O}_{16}
\times
L^{5}\ .
\]
Namely, the infrared contribution remains finite when $\epsilon
\rightarrow 0$. 

We now examine the ultraviolet contribution more closely. From Theorem
\ref{unnamedtheorem1bis} the small $z$ expansion of $\Psi_{v_{*}}(z e_{\rm{u}})$ is of the
form

\begin{equation}\label{Psiexpansion}
\Psi_{v_{*}}
(
z
e_{\rm{u}}
)
=
v_{*}
+
z
e_{\rm{u}}
+
z^{2}
\Theta
+
O(z^{3})
\end{equation}
for some vector $\Theta$ to be determined shortly. 
Using the decomposition in Lemma \ref{d2bbound}

\[
D_{0}^{2}
(
\delta b
\circ
\Psi_{v_{*}}
)
[e_{\rm{u}},e_{\rm{u}}]
=
\frac{d^{2}}{dz^2} \Big|_{z=0}
\delta b^{\rm{explicit}} 
\left(
\Psi_{v_{*}}(z e_{\rm{u}})
\right)
+
\frac{d^{2}}{dz^2} \Big|_{z=0}
\delta b^{\rm{implicit}} 
\left(
\Psi_{v_{*}}(z e_{\rm{u}})
\right).
\]

If $\displaystyle |z| \le \frac{1}{30}\bar{g}^{e_{2}}$ then as before we get 

\begin{equation*}
\begin{split}
|| \Psi_{v_{*}}(z e_{\rm{u}}) ||
\le&
||v_{*}||
+
2
||z e_{\rm{u}}||\\
\le&
\frac{1}{40} + \frac{1}{15} < \frac{1}{2}.
\end{split}
\end{equation*}
So by Theorem \ref{mainestthm}

\[
\left|
\delta b^{\rm{implicit}} 
\left(
\Psi_{v_{*}}(z e_{\rm{u}})
\right)
\right|
\le
\frac{1}{2}B_{0}\bar{g}^{e_{R}}.
\]
Cauchy's formula then immediately implies

\[
\left|
\frac{d^2}{dz^2} \Big|_{z=0}
\delta b^{\rm{implicit}}
\left(
\Psi_{v_{*}}(z e_{\rm{u}})
\right)
\right|
\le
2!
\left(
\frac{1}{30}
\bar{g}^{e_{2}}
\right)^{-2}
\times
\frac{1}{2}
B_{0}
\bar{g}^{e_{R}}.
\]
Since $e_{R} > 2e_{2}$, we must have
\[
\lim_{\epsilon \rightarrow 0}
\frac{d^2}{dz^2} \Big|_{z=0}
\delta b^{\rm{implicit}}
\left(
\Psi_{v_{*}}(z e_{\rm{u}})
\right)
=
0.
\]

Now recall that

\begin{gather*}
\delta b^{\rm{explicit}} 
=
A_{4} \bar{g}^{2}
+
\delta b^{\rm{explicit}}_{\rm\MakeUppercase{\romannumeral 1}}
( \delta g, \mu, R)
+
\delta b^{\rm{explicit}}_{\rm\MakeUppercase{\romannumeral 2}}
(\delta g, \mu, R)\\[1.2ex]
\textnormal{  where  }
\delta b^{\rm{explicit}}_{\rm\MakeUppercase{\romannumeral 1}}
( \delta g, \mu, R)
=
2A_{4}\bar{g}\delta g
+
A_{4} \delta g^{2}\\[1.2ex]
\textnormal{  and  }
\delta b^{\rm{explicit}}_{\rm\MakeUppercase{\romannumeral 2}}
(\delta g, \mu, R)
=
A_{5} \mu^{2}\ .
\end{gather*}

Note that the $A_{4} \bar{g}^{2}$ term disappears in the computation
of derivatives while $\delta
b^{\rm{explicit}}_{\rm\MakeUppercase{\romannumeral 1}}$ can be treated
as we treated $\delta b^{\rm{implicit}}$. Indeed by Cauchy's formula
and Theorem \ref{mainestthm}

\[
\left|
\frac{d^2}{dz^2} \Big|_{z=0}
\delta b^{\rm{explicit}}_{\rm\MakeUppercase{\romannumeral 1}}
\left(
\Psi_{v_{*}}(z e_{\rm{u}})
\right)
\right|
\le
2!
\left(
\frac{1}{30}
\bar{g}^{e_{2}}
\right)^{-2}
\times
A_{4,\rm{max}}
\left[
2
\times
\bar{g}
\times
\frac{1}{2}
\bar{g}^{e_{4}}
+
\left(
\frac{1}{2}
\bar{g}^{e_{4}}
\right)^{2}
\right].
\]

Since $e_{4} > 2e_{2}-1$ and $e_{4} > e_{2}$ we have

\[
\lim_{\epsilon \rightarrow 0}
\frac{d^2}{dz^2} \Big|_{z=0}
\delta b^{\rm{explicit}}_{\rm\MakeUppercase{\romannumeral 1}}
\left(
\Psi_{v_{*}}(z e_{\rm{u}})
\right)
=
0.
\]

As a result of the formula $e_{\rm{u}} = \left( \delta
g'_{\rm{u}}(\mu_{*}), 1, R'_{\rm{u}}(\mu_{*}) \right)$ and the
expansion \eqref{Psiexpansion} we easily compute

\[
\frac{d^2}{dz^2} \Big|_{z=0}
\delta b^{\rm{explicit}}_{\rm\MakeUppercase{\romannumeral 2}}
\left(
\Psi_{v_{*}}(z e_{\rm{u}})
\right)
=
2
A_{5}
\left(
1
+
2
\mu_{*}
\Theta_{\mu}
\right)
\]

where $\Theta_{\mu}$ is the $\mu$ component of $\Theta \in
\mathcal{E}$.

We determine the latter using the intertwining relation in Theorem
\ref{unnamedtheorem1bis} for small $z$. 

We have by an easy calculation using \eqref{Psiexpansion}

\[
RG
\left(
\Psi_{v_{*}}(z e_{\rm{u}})
\right)
=
v_{*}
+
D_{v_{*}}
RG
[e_{\rm{u}}]
+
z^{2}
\left(
D_{v_{*}}
RG
[
\Theta
]
+
\frac{1}{2}
D^{2}_{v_{*}}
RG
[e_{\rm{u}},e_{\rm{u}}]
\right)
+
O(z^{3})\ .
\]

But this is the same as

\[
\Psi_{v_{*}}(\alpha_{\rm{u}}z e_{\rm{u}})
=
v_{*}
+
z
\alpha_{\rm{u}}
e_{\rm{u}}
+
z^{2}
\alpha_{\rm{u}}^{2}
\Theta
+
O(z^{3}).
\]

Thus

\begin{equation}\label{Thetarelation}
\alpha_{\rm{u}}^{2}
\Theta
=
D_{v_{*}}
RG
[
\Theta
]
+
\frac{1}{2}
D^{2}_{v_{*}}
RG
[e_{\rm{u}},e_{\rm{u}}].
\end{equation}

On the other hand $\displaystyle \Psi_{v_{*}} \in W^{\rm{u},\rm{loc}}$
for $z$ small and therefore

\[
\left[
\Psi_{v_{*}}(z e_{\rm{u}})
\right]_{\delta g}
=
\delta g_{\rm{u}}
\left(
\left[
\Psi_{v_{*}}(z e_{\rm{u}})
\right]_{ \mu }
\right)
\]

and

\[
\left[
\Psi_{v_{*}}(z e_{\rm{u}})
\right]_{R}
=
R_{\rm{u}}
\left(
\left[
\Psi_{v_{*}}(z e_{\rm{u}})
\right]_{\mu}
\right)
\]

where $[ \cdots ]_{\delta g}$, $[\cdots]_{\mu}$, and $[\cdots]_{R}$
refer to the $\delta g$, $\mu$, and $R$ components respectively.\\

Expanding these relations up to second order imply

\begin{gather*}
\Theta
=
(
\Theta_{\delta g},
\Theta_{\mu},
\Theta_{R}
)
=
\Theta_{\mu}
e_{\rm{u}}
+
\frac{1}{2}
c_{\rm{u}}\\
\textnormal{  where  }
c_{\rm{u}}
=
(
\delta g_{\rm{u}}''(\mu_{*}),
0,
R_{\rm{u}}''(\mu_{*})
).
\end{gather*}

Taking the $\mu$ component of \eqref{Thetarelation} we see that

\[
\alpha_{\rm{u}}^{2}
\Theta_{\mu}
=
\Theta_{\mu}
\alpha_{\rm{u}}
+
\frac{1}{2}
\left[
D_{v_{*}}
RG
[c_{\rm{u}}]
\right]_{\mu}
+
\frac{1}{2}
\left[
D^{2}_{v_{*}}
RG
[
e_{\rm{u}},e_{\rm{u}}
]
\right]_{\mu}
\]
where we have used $[e_{\rm{u}}]_{\mu}=1$ and
$D_{v_{*}}RG[e_{\rm{u}}]=\alpha_{\rm{u}}e_{\rm{u}}$.\\
Since $\alpha_{\rm u}$ we have
\[
\Theta_{\mu}=\frac{1}{2\alpha_{\rm u}(\alpha_{\rm u}-1)}
\left\{
\left[
D_{v_{*}}
RG
[c_{\rm{u}}]
\right]_{\mu}
+
\left[
D^{2}_{v_{*}}
RG
[
e_{\rm{u}},e_{\rm{u}}
]
\right]_{\mu}
\right\}\ .
\]

Now $\displaystyle |\mu - \mu_{*}| < \rho''\bar{g}^{e_{2}}$ implies
$\displaystyle |\delta g_{\rm{u}}(\mu)| \le
\frac{\rho'}{3}\bar{g}^{e_{4}}$ and $\displaystyle
|||R_{\rm{u}}(\mu)|||_{\bar{g}} \le \frac{\rho'}{3} \bar{g}^{e_{R}}$.
Using $\displaystyle |\mu - \mu_{*}| = \frac{1}{2} \rho''
\bar{g}^{e_{2}}$ as a contour of integration, Cauchy's formula implies
the following estimates:

\begin{equation*}
\begin{split}
|\delta g'_{\rm{u}}(\mu_{*})| 
&\le 
\frac{2 \rho'}{3 \rho''}
\bar{g}^{e_{4}-e_{2}}\\[1.2ex]
|\delta g''_{\rm{u}}(\mu_{*})| 
&\le 
\frac{8 \rho'}{3 (\rho'')^2}
\bar{g}^{e_{4}-2e_{2}}\\[1.2ex]
|||R_{\rm{u}}'(\mu_{*})|||_{\bar{g}}
\le&
\frac{2 \rho'}{3 \rho''}
\bar{g}^{e_{R}-e_{2}}\\[1.2ex]
|||R_{\rm{u}}''(\mu_{*})|||_{\bar{g}}
\le&
\frac{8 \rho'}{3 (\rho'')^2}
\bar{g}^{e_{R}-2e_{2}}.
\end{split}
\end{equation*}

As a result

\[
||c_{\rm{u}}||
=
\max
\left\{
|\delta g''_{\rm{u}}(\mu_{*})| \bar{g}^{-e_{4}},
|||R_{\rm{u}}''(\mu_{*})|||_{\bar{g}}\bar{g}^{-e_{R}}
\right\}
\le
\frac{8 \rho'}{3 (\rho'')^2} 
\bar{g}^{-2e_{2}}\ .
\]

From the explicit formulas in the proof of Lemma \ref{DRG2lem} and
following the same notation

\begin{equation}\label{DRGmu}
\left[
D_{v}
RG
[v']
\right]_{\mu}
=
L^{\frac{3+\epsilon}{2}}
\mu'
-
2
A_{2}
(\bar{g}+\delta g)
\delta g'
-
A_{3}
(\bar{g} + \delta g)
\mu'
-
A_{3}
\mu \ 
\delta g'
+
\left[
D_{v}
RG^{\rm{implicit}}
[v']
\right]_{\mu}.
\end{equation}

For $v=v_{*}$ and $v' = c_{\rm{u}}$ this gives

\[
\left[
D_{v_{*}}
RG
[v']
\right]_{\mu}
=
-
2
A_{2}
(\bar{g}+\delta g_{*})
\delta g_{\rm{u}}''(\mu_{*})
-
A_{3}
\mu_{*}
\delta g_{\rm{u}}''(\mu_{*})
+
\left[
D_{v_{*}}
RG^{\rm{implicit}}
[c_{\rm{u}}]
\right]_{\mu}.
\]

The infinitesimal version of the $\xi_{2}$ Lipschitz estimate in Lemma
\ref{Lipxilem} immediately implies

\[
||
\left[
D_{v_{*}}RG^{\rm{implicit}}
\right]_{\mu}
||
\le
2
B_{2}
\bar{g}^{e_{R}}
\]

for the operator norm induced on linear maps from $\mathcal{E}$ to
$\mathbb{C}$ by the norm $|| \cdot ||$ on $\mathcal{E}$ and the
modulus on $\mathbb{C}$. As a result we have

\[
\left|
\left[
D_{v_{*}}
RG
[c_{\rm{u}}]
\right]_{\mu}
\right|
\le
\left(
2
A_{2,\rm{max}}
\times
\frac{3}{2}
\bar{g}
+
A_{3,\rm{max}}
\times
\frac{1}{2}
\bar{g}^{e_{2}}
\right)
\times
\frac{8 \rho'}{3 (\rho'')^{2}}
\bar{g}^{e_{4}-2e_{2}}
+
2
B_{2}
\bar{g}^{e_{R}}
\times
\frac{8 \rho'}{3 (\rho'')^{2}}
\bar{g}^{-2e_{2}}\ .
\]

Since $e_{4} > e_{2}$, $e_{4} > 2e_{2}-1$, and $e_{R} > 2e_{2}$ we
have

\[
\lim_{\epsilon \rightarrow 0}
\left[
D_{v_{*}}
RG
[c_{\rm{u}}]
\right]_{\mu}
=
0\ .
\]

Also from the formulas in Lemma \ref{DRG2lem}

\[
\left[
D_{v_{*}}^{2}
RG
[
e_{\rm{u}},
e_{\rm{u}}
]
\right]_{\mu}
=
-
2
A_{2}
\delta g'_{\rm{u}}(\mu_{*})^2
-
2
A_{3}
\delta g'_{\rm{u}}(\mu_{*})
+
\left[
D_{v_{*}}^{2}
RG^{\rm{implicit}}
[
e_{\rm{u}},
e_{\rm{u}}
]
\right]_{\mu}
\]

and therefore

\[
\left|
\left[
D_{v_{*}}^{2}
RG
[
e_{\rm{u}},
e_{\rm{u}}
]
\right]_{\mu}
\right|
\le
2
A_{2,\rm{max}}
\left(
\frac{2 \rho'}{3 \rho''}
\bar{g}^{e_{4}-e_{2}}
\right)^{2}
+
2
A_{3,\rm{max}}
\left(
\frac{2 \rho'}{3 \rho''}
\bar{g}^{e_{4}-e_{2}}
\right)
+
\left|
\left[
D_{v_{*}}^{2}
RG^{\rm{implicit}}
[
e_{\rm{u}},
e_{\rm{u}}
]
\right]_{\mu}
\right|.
\] 

The argument of Lemma \ref{DRG2lem} only applied to the $\mu$
component of $RG^{\rm{implicit}}$ gives

\[
||
D_{v_{*}}^{2}
RG^{\rm{implicit}}_{\mu}
||
\le
32
B_{2}
\bar{g}^{e_{R}-e_{2}}
\]

for the norm of the second differential.\\

Since $\displaystyle ||e_{\rm{u}}|| = \bar{g}^{-e_{2}}$ we obtain 

\[
\left|
\left[
D_{v_{*}}^{2}
RG^{\rm{implicit}}
[
e_{\rm{u}},
e_{\rm{u}}
]
\right]_{\mu}
\right|
\le
32
B_{2}
\bar{g}^{e_{R}-3e_{2}}.
\]

Since $\displaystyle |\mu_{*}| \le \frac{1}{2} \bar{g}^{e_{2}}$ and
$e_{R} > 2e_{2}$ we have

\[
\lim_{\epsilon \rightarrow 0}
\mu_{*}
\left[
D_{v_{*}}^{2}
RG^{\rm{implicit}}
[
e_{\rm{u}},
e_{\rm{u}}
]
\right]_{\mu}
=
0.
\]

Since also $\displaystyle \lim_{\epsilon \rightarrow 0} \mu_{*} = 0$
and $\displaystyle \lim_{\epsilon \rightarrow 0} \alpha_{\rm{u}} =
L^{\frac{3}{2}} > 1$ we have enough to affirm

\[
\lim_{\epsilon \rightarrow 0}
\mu_{*}
\Theta_{\mu}
=
0.
\]

Thus

\[
\lim_{\epsilon \rightarrow 0}
\frac{d^2}{dz^2}\Big|_{z=0}
D_{0}^{2}( \delta b \circ \Psi_{v_{*}})
[
e_{\rm{u}},
e_{\rm{u}}
]
=
2
\lim_{\epsilon \rightarrow 0}
A_{5}
=
2 L^{3} ( 1- p^{-3}) \times l > 0 \textnormal{  by Lemma
  \ref{gamL2lem}.}
\]\\

We now study the $\epsilon \rightarrow 0$ asymptotics of
$\alpha_{\rm{u}}$ more closely. One way to get a precise hold on this
eigenvalue is to note that

\[
\alpha_{\rm{u}}
=
\left[
D_{v_{*}}
RG[
e_{\rm{u}}
]
\right]_{\mu}.
\]

Then by the formula in \eqref{DRGmu} we have

\[
\alpha_{\rm{u}}
=
L^{\frac{3+\epsilon}{2}}
-
2
A_{2}
(\bar{g} + \delta g_{*})
\delta g'_{\rm{u}}(\mu_{*})
-
A_{3}
(\bar{g} + \delta g_{*})
-
A_{3}
\mu_{*}
\delta g'_{\rm{u}}(\mu_{*})
+
\left[
D_{v_{*}}
RG^{\rm{implicit}}
[e_{\rm{u}}]
\right]_{\mu}
\]
since $e_{\rm{u}} = \left( \delta g'_{\rm{u}}(\mu_{*}), 1,
  R'_{\rm{u}}(\mu_{*}) \right)$. \\

As before

\begin{equation*}
\begin{split}
\left|
\left[
D_{v_{*}}
RG^{\rm{implicit}}
[ e_{\rm{u}} ]
\right]_{\mu}
\right|
\le&
||
D_{v_{*}} RG^{\rm{implicit}}
||
\times
||e_{\rm{u}}||\\
\le&
2 
B_{2} 
\bar{g}^{e_{R}}
\times
\bar{g}^{-e_{2}}.
\end{split}
\end{equation*}

But $e_{R} > e_{2}+1$ and $\bar{g}$ is of order $\epsilon$ so

\[
\left[
D_{v_{*}}
RG^{\rm{implicit}}
[e_{\rm{u}}]
\right]_{\mu}
=
o(\epsilon).
\]

We have

\[
\left|
-2
A_{2}
(\bar{g} + \delta g_{*})
\delta g'_{\rm{u}}(\mu_{*})
\right|
\le
2
A_{2,\rm{max}}
\times
\frac{3}{2}
\bar{g}
\times
\frac{2 \rho'}{3 \rho''}
\bar{g}^{e_{4}-e_{2}}
\]
so this is also an $o(\epsilon)$ term because $e_{4} > e_{2}$. \\

Likewise

\[
\left|
-A_{3}
\delta g_{*}
\right|
\le
A_{3,\rm{max}}
\times
\frac{1}{2}\bar{g}^{e_{4}}
\]
so this is $o(\epsilon)$ because $e_{4} > 1$. \\

Finally,

\[
\left|
-A_{3}
\mu_{*}
\delta g'_{\rm{u}}(\mu_{*})
\right|
\le
A_{3,\rm{max}}
\times
\frac{1}{2}
\bar{g}^{e_{2}}
\times
\frac{ 2 \rho'}{3 \rho''}
\bar{g}^{e_{4}-e_{2}}
\]

so this is an $o(\epsilon)$ term too.\\

As a result we have

\begin{equation*}\begin{split}
\alpha_{\rm{u}}
=&
L^{\frac{3+\epsilon}{2}}
-
A_{3}
\bar{g}
+
o(\epsilon)\\
=&
L^{\frac{3+\epsilon}{2}}
-
12
\times
L^{\frac{3+\epsilon}{2}}
\times
\frac{A_{1}}{36 L^{\epsilon}}
\bar{g}
+
o(\epsilon)\\
=&
L^{\frac{3+\epsilon}{2}}
\left(
1
-
\frac{1}{3}
\left(
\frac{L^{\epsilon}-1}{L^{\epsilon}}
\right)
\right)
+
o(\epsilon)
\end{split}
\end{equation*}

from the relations between $A_{3}$, $A_{1}$, and $\bar{g}$.

It is now a simple calculus exercise to derive

\[
\eta_{\phi^{2}} 
= 
\frac{2}{3} \epsilon 
+ 
o(\epsilon)
\]

where $\eta_{\phi^{2}}$ is defined by

\[
L^{\frac{1}{2} \eta_{\phi^{2}}}
=
Z_{2}^{-1}
=
L^{\frac{3+\epsilon}{2}}
\alpha_{\rm{u}}^{-1}.
\]

We also easily get

\[
L^{3}
\alpha_{\rm{u}}^{-2}
=
1
-
\frac{\log(L)}{3}\epsilon
+
o(\epsilon)
\]

which proves the earlier statement that

\[
L^{3} \alpha_{\rm{u}}^{-2} 
< 
1
\]

in the small $\epsilon$ regime which was crucial for the convergence
and analyticity in the ultraviolet regime.\\

Another byproduct is 

\[
\frac{1}{\alpha_{\rm{u}}^{2} - L^{3}}
\sim
\frac{ 3 L^{-3}}{\log(L)} 
\times 
\frac{1}{\epsilon}
\]
and therefore

\[
\langle 
\langle
N[\phi^{2}]
(\bbone_{\mathbb{Z}_{p}^{3}})^2
\rangle^{\rm{UV}}_{\rm{reduced}}
\sim
\frac{6(1-p^{-3})}{\log(p)}
\times
\frac{1}{\epsilon}
\]
when $\epsilon \rightarrow 0$.\\

Since $\langle
N[\phi^{2}]
(\bbone_{\mathbb{Z}_{p}^{3}})^2
\rangle^{\rm{IR}}_{\rm{reduced}}$ remains bounded, the quantity

\[
\langle
N[\phi^{2}]
(\bbone_{\mathbb{Z}_{p}^{3}})^2
\rangle_{\rm{reduced}}
=
\langle
N[\phi^{2}]
(\bbone_{\mathbb{Z}_{p}^{3}})^2
\rangle^{\rm{UV}}_{\rm{reduced}}
+
\langle
N[\phi^{2}]
(\bbone_{\mathbb{Z}_{p}^{3}})^2
\rangle^{\rm{IR}}_{\rm{reduced}}
\]

is strictly positive for $\epsilon$ small enough.\\

Provided $\varkappa_{\phi^{2}} \not = 0$ we can then impose by definition

\[
Y_{2}
=
|\varkappa_{\phi^{2}} |^{-1}
\times
\left\{
\langle
N[\phi^{2}]
(\bbone_{\mathbb{Z}_{p}^{3}})^2
\rangle_{\rm{reduced}}
\right\}^{-\frac{1}{2}}
\]
and thus force the normalization

\[
\langle
N[\phi^{2}]
(\bbone_{\mathbb{Z}_{p}^{3}})^2
\rangle
=
1.
\]\\

We now address the issue of showing $\varkappa_{\phi^{2}} \not = 0$. While
most of the proof so far relied on quantitative estimates, here we had
to use a more qualitative approach. This is because of the slow
convergence to the fixed point on the stable manifold and the fact
that we do not have much freedom of choice for our starting point $v$. The latter
has to be on the $R=0$ bare surface and therefore we cannot choose it
as close to $v_{*}$ as we would like to.\\

Recall that $W^{\rm{s},\rm{loc}}_{\rm{int}}$ is parametrized as

\[
v_{1}
\mapsto
( v_{1}, \mu_{\rm{s}}(v_{1}))
\]
for $\displaystyle ||v_{1}|| < \frac{\rho}{13}$ in $\mathcal{E}_{1}$.
For $v \in W^{\rm{s},\rm{loc}}_{\rm{int}}$ we consider the tangent
space $T_{v}W^{\rm{s}}$ defined as the kernel of the linear form

\[
( w_{1}, w_{2})
\mapsto
w_{2}
-
D_{v_{1}}\mu_{\rm{s}}[w_{1}]
\]
via the identification $\mathcal{E}_{2} \simeq \mathbb{C}$.

This linear form is continuous and does not vanish identically, so
$T_{v}W^{\rm{s}}$ is a closed complex hyperplane in $\mathcal{E}$. If
$w \in \mathcal{E}$ satisfies $w \not \in T_{v}W^{\rm{s}}$ then we
have a direct sum decomposition $\mathcal{E} = \mathbb{C} \oplus
T_{v}W^{\rm{s}}$.

We have the following infinitesimal version of Parts 1) and Parts 2) of
Lemma \ref{L57lem} and Lemma \ref{L58lem}.

\begin{Lemma}\label{infteslem}

For all $v \in W_{\rm{int}}^{\rm s,loc}$ we have:\\

1) for all $w \in \mathcal{E}$, 

\[
|| ( D_{v} RG[w] )_{1} || 
\le 
c_{1}(\epsilon)||w||
\]\\

2) for all $w \in \mathcal{E}$, such that $L^{\frac{3}{4}}||w_{2}||
\ge ||w_{1}||$,

\[
||
(D_{v}RG[w])_{2}
||
\ge
c_{2}(\epsilon)
||w||
\]

3) for all $w \in T_{v} W^{\rm{s}}$,

\[
||w_{1}||
\ge
L^{\frac{3}{4}}
||w_{2}||\ .
\]

\end{Lemma}

\noindent{\bf Proof: }
Consider the complex curve $\gamma(t) = v + t w$ for $t$ small which
ensures that $\Gamma(t) \in \bar{B}\left(0,\frac{1}{8}\right)$.
Lemma \ref{L57lem} Part 1) implies

\[
||
RG_{1}(\gamma(t))
-
RG_{1}(\gamma(0))
||
\le
c_{1}(\epsilon)
||tw||.
\]
Dividing by $|t|$ and taking $t \rightarrow 0$ we immediately get
$||(D_{v}(RG[w])_{1}|| \le c_{1}(\epsilon)||w||$.

Now if $L^{\frac{3}{4}} ||w_{2}|| \ge ||w_{1}||$ then we have

\[
L^{\frac{3}{4}} 
||\gamma(t)_{2}
- \gamma(0)_{2}||
\ge
||\gamma(t)_{1}
- \gamma(0)_{1}||
\]
and thus

\[
||
RG_{2}(\gamma(t))
-
RG_{2}(\gamma(0))
||
\ge
c_{2}(\epsilon)
||tw||
\]
by Lemma \ref{L57lem} Part 2). Taking the $t \rightarrow 0$ limit as
before we obtain

\[
||
(D_{v}RG[w])_{2}
||
\ge
c_{2}(\epsilon)
||w||.
\]

For the third part we use Lemma \ref{L58lem} to write

\[
||(v_{1}+tw_{1})-v_{1}||
\ge
L^{\frac{3}{4}}
||
\mu_{\rm{s}}(v_{1}+tw_{1})
-
\mu_{\rm{s}}(v_{1})
||
\]

for t small. Dividing by $|t|$ and taking $t \rightarrow 0$ gives

\[
||w_{1}||
\ge
L^{\frac{3}{4}}
||D_{v_{1}}\mu_{\rm{s}}[w_{1}]||
=
L^{\frac{3}{4}}
||w_{2}||
\]

since $w \in T_{v}W^{\rm{s}}$. \qed\\

\begin{Lemma}\label{nontanglem}

For all $v \in W^{\rm{s},\rm{loc}}_{\rm{int}}$ and $w \in \mathcal{E}$
we have the implication

\[
L^{\frac{3}{4}}
||w_{2}||
>
||w_{1}||
\Rightarrow
D_{v}
RG[w]
\not
\in
T_{RG(v)}
W^{\rm{s}}\ .
\]

\end{Lemma}

\noindent{\bf Proof:}
We proceed by contradiction. Suppose

\[
L^{\frac{3}{4}}
||w_{2}||
>
||w_{1}||
\textnormal{  and  }
D_{v}
RG[w]
\in
T_{RG(v)}
W^{\rm{s}}.
\]
Then by Lemma \ref{infteslem} Parts 1), 2), 3) we have

\[
c_{1}(\epsilon)
||w||
\ge
||
(D_{v}RG[w])_{1}
||,
\]

\[
||
(D_{v}RG[w])_{2}
||
\ge
c_{2}(\epsilon)
||w||
\]
and

\[
||
(D_{v}RG[w])_{1}
||
\ge
L^{\frac{3}{4}}
||
(D_{v}RG[w])_{2}
||
\]
respectively. As a result 

\[
c_{1}(\epsilon)
||w||
\ge
L^{\frac{3}{4}}
c_{2}(\epsilon)
||w||.
\]

But $c_{1}(\epsilon) < 1 < L^{\frac{3}{4}} c_{2}(\epsilon)$ so
$||w||=0$ which contradicts the strict inequality
$L^{\frac{3}{4}}||w_{2}|| > ||w_{1}||$. \qed\\

\begin{Lemma}

For all $v \in W^{\rm{s},\rm{loc}}_{\rm{int}}$ and $w \in
T_{v}W^{\rm{s}}$

\[
T_{1}(v)[w]
\in
T_{RG(v)}W^{\rm{s}}
\]
and

\[
T_{\infty}(v)[w] 
\in 
T_{v_{*}}W^{\rm{s}}\ .
\]

\end{Lemma}

\noindent{\bf Proof:}
Consider the curve $t \mapsto \left( v_{1}+tw_{1},
  \mu_{\rm{s}}(v_{1}+tw_{1})\right)$ in
$W^{\rm{s},\rm{loc}}_{\rm{int}}$ for $t$ small.
By Proposition \ref{L56prop} and the parametrization of
$W^{\rm{s},\rm{loc}}_{\rm{int}}$ we have

\[
RG_{2}(v_{1}+tw_{1},
\mu_{\rm{s}}(v_{1}+tw_{1}))
=
\mu_{\rm{s}}
\left(
RG_{1}(v_{1}+tw_{2},
\mu_{\rm{s}}(v_{1}+tw_{1}))
\right).
\]
Differentiating this at $t=0$ gives

\[
\left(
D_{v}
RG
[
(w_{1},
D_{v_{1}} \mu_{\rm{s}}[w_{1}]
)
]
\right)_{2}
=
D_{RG_{1}(v)}
\mu_{\rm{s}}
\left[
\left(
D_{v}
RG
[(w_{1},D_{v_{1}} \mu_{\rm{s}}[w_{1}])]
\right)_{1}
\right],
\]
i.e.,

\[
\left(
D_{v}
RG[w]
\right)_{2}
=
D_{RG_{1}(v)}
\mu_{\rm{s}}
\left[
\left(
D_{v}
RG[w]
\right)_{1}
\right].
\]

Hence $D_{v}RG[w]$ belongs to $T_{RG(v)}W^{\rm{s}}$ and so does
$T_{1}(v)[w]=\alpha_{\rm{u}}^{-1}D_{v}RG[w]$.

By iteration this immediately implies

\[
T_{n}(v)[w]
\in
T_{RG^{n}(v)}
W^{\rm{s}}
\]

for all integer $n \ge 0$.

Namely, we have

\[
\left(
T_{n}(v)[w] 
\right)_{2}
=
D_{ \left(RG^{n}(v)\right)_{1}}
\mu_{\rm{s}}
\left[
\left(
T_{n}(v)[w]
\right)_{1}
\right].
\]

Using continuity, the remark following Lemma \ref{Tcvlem}, and the fact
that $RG^{n}(v) \rightarrow v_{*}$, we can take the $n \rightarrow
\infty$ limit in the previous equality and obtain

\[
\left(
T_{\infty}
(v)[w] \right)_{2}
=
D_{v_{*,1}}
\mu_{\rm{s}}
\left[
\left(
T_{\infty}(v)[w]
\right)_{1}
\right].
\]

This proves $ T_{\infty}(v)[w] \in T_{v_{*}}W^{\rm{s}} =
\mathcal{E}^{\rm{s}}$ by definition of $\mathcal{E}^{\rm{s}}$. \qed\\

\begin{Lemma}\label{tanzerolem}
For all $v \in W^{\rm{s},\rm{loc}}_{\rm{int}}$ and $w \in
T_{v}W^{\rm{s}}$

\[
D_{0}
\Psi_{v}[w]
=
0,
\]

where the differential is with respect to the $w$ variable at $w=0$
for the function $\Psi_{v}(\bullet)=\Psi(v,\bullet)$.

\end{Lemma}

\noindent{\bf Proof: }
By Theorem \ref{unnamedtheorem1} Part 5)

\[
\Psi_{v}
=
\Psi_{v_{*}}
\circ
T_{\infty}(v)
\]
and thus by the chain rule

\[
D_{0}
\Psi_{v}
[w]
=
D_{0}
\Psi_{v_{*}}
\left[
T_{\infty}(v)
[w]
\right].
\]

However by the previous lemma $ T_{\infty}(v)[w] \in
\mathcal{E}^{\rm{s}}$ so $P_{\rm{u}}T_{\infty}(v)[w]=0$. But we also
have $P_{s}T_{\infty}(v)[w] = 0$ as a follow up to Lemma \ref{Tcvlem}.

As a result, $T_{\infty}(v)[w] = 0$ and consequently $D_{0}
\Psi_{v}[w]=0$. \qed\\

\begin{Lemma}

For all $v \in W_{\rm{int}}^{\rm{s},\rm{loc}}$, 
if $D_{0} \Psi_{v} = 0$ then $D_{0} \Psi_{RG(v)} = 0$.

\end{Lemma}

\noindent{\bf Proof:}
By Theorem \ref{unnamedtheorem1} Part 4)

\[
\Psi_{v} = \Psi_{RG(v)} \circ T_{1}(v)
\]

near $w=0$. Differentiating at zero gives

\begin{equation}\label{PsiTrelation}
D_{0}\Psi_{v}
=
D_{0}\Psi_{RG(v)}
\circ
T_{1}(v).
\end{equation}

Pick some vector $u \in \mathcal{E}$ satisfying the hypothesis of
lemma \ref{nontanglem}.
For instance $e_{\phi^{2}}$ works since $L^{\frac{3}{4}}
||e_{\phi^{2},2}|| = L^{\frac{3}{4}} \bar{g}^{-e_{2}} >
||e_{\phi^{2},1}|| = 0$.
By the same lemma $\displaystyle T_{1}(v)[u] \not \in T_{RG(v)}W^{\rm{s}}$ and
therefore $\displaystyle \mathcal{E} = \mathbb{C} T_{1}(v)[u] \oplus
T_{RG(v)}W^{\rm{s}}$. 

Let $w \in \mathcal{E}$. We decompose it as $w = \lambda T_{1}(v)[u] +
w'$ with $w' \in T_{RG(v)}W^{\rm{s}}$. 
Then by \eqref{PsiTrelation}:

\begin{equation*}
\begin{split}
D_{0}
\Psi_{RG(v)}
[w]
=&
\lambda
D_{0}
\Psi_{v}[u]
+
D_{0}
\Psi_{RG(v)}[w']\\
=&
0
\end{split}
\end{equation*}
by the hypothesis and the previous lemma for $RG(v)$ instead of $v$. 
Hence the differential $D_{0} \Psi_{RG(v)}$ vanishes. \qed\\

Iterating the last lemma we see that if $D_{*} \Psi_{v} = 0$ then
$D_{0} \Psi_{RG^{n}(v)} = 0$ for all $n \ge 0$. 
By the joint analyticity in Theorem \ref{unnamedtheorem1} we can take the
$n \rightarrow \infty$ limit which gives $D_{0} \Psi_{v_{*}} = 0$ and
therefore

\[
\frac{d}{dz} \Big|_{z=0}
\Psi_{v_{*}} ( z e_{\rm{u}} ) 
= 
0
\]
which contradicts (\ref{Psiexpansion}) and $e_{\rm{u}}\neq 0$.\\

We have proved $D_{0}\Psi_{v} \not = 0$ for all $v \in
W_{\rm{int}}^{\rm{s},\rm{loc}}$. Now since $e_{\phi^{2}}$ satisfies
$L^{\frac{3}{4}} ||e_{\phi^{2},2}|| > ||e_{\phi^{2},1}||$ we know that
$e_{\phi^{2}} \not \in T_{v}W^{\rm{s}}$ by Lemma \ref{infteslem}. Thus
$\mathcal{E} = \mathbb{C} e_{\phi^{2}} \oplus T_{v} W^{\rm{s}}$. 

Recall that $D_{0} \Psi_{v} = D_{0} \Psi_{v_{*}} \circ T_{\infty}(v)$
so $D_{0} \Psi_{v}[e_{\phi^{2}}] = \varkappa_{\phi^{2}} D_{0} \Psi_{v_{*}}
[e_{\rm{u}}]$ by definition of $\varkappa_{\phi^{2}}$. 
If the latter vanishes then $D_{0}\Psi_{v}$ vanishes on
$\mathbb{C}e_{\phi^{2}}$ and therefore on all of $\mathcal{E}$ by
Lemma \ref{tanzerolem}. This contradicts $D_{0}\Psi_{v} \not = 0$. We
have now finally proved $\varkappa_{\phi^{2}} \not = 0$.\\

The remaining items to be settled are the mini-universality result and
the choice of parameters $\eta$, $e$, etc.\\

The mini-universality should be clear at this point: the generating
function $\mathcal{S}^{T}(\tilde{f},\tilde{j})$ does not depend on the
starting point $v = (g - \bar{g}, \mu_{\rm{c}}(g),0) \in
W^{\rm{s},\rm{loc}}_{\rm{int}}$ for the $RG$ iterations.
Indeed using $\Psi_{v} = \Psi_{v_{*}} \circ T_{\infty}(v)$ we see that
the effect of $v$ is entirely in the multiplying factor
$\varkappa_{\phi^{2}}$ which however always comes in the combination $Y_{2}
\varkappa_{\phi^{2}}$. By our choice of normalization, $Y_{2}
\varkappa_{\phi^{2}}$ is defined in terms of the reduced $N[\phi^{2}]$
two-point function which only involves data at the fixed point
$v_{*}$. 

Finally, to complete our very long proof of Theorem \ref{themainthm} we have to
pick a choice of parameters which satisfies all the required
inequalities.
We pick

\begin{equation*}
\begin{split}
\eta &= 0\\
\eta_{R} &= \frac{1}{8}\\
e_{1}=e_{2}&=e_{3}=1\\
e_{4}&=\frac{3}{2}\\
e_{W}&=2\\
e_{R}&=\frac{21}{8}\\
\rho'&=\frac{1}{8}\\
\rho&=\frac{1}{128}\\
\rho''&=\frac{1}{768}\ .
\end{split}
\end{equation*}

We leave it to the reader to check that these choices indeed satisfy
all the previously stated inequalities. Note that $\rho'''$ was already
defined and was only needed in the local analysis at the fixed
point. The proof of Theorem \ref{themainthm} is now complete.